%% file: main.tex
\documentclass[11pt]{article}

\usepackage{iftex}
\ifPDFTeX
  \usepackage[T1]{fontenc}
  \usepackage[utf8]{inputenc}
  \usepackage{lmodern}
\else
  \usepackage{fontspec}
\fi
\usepackage[margin=1in]{geometry}
\usepackage{amsmath,amssymb,amsthm,mathtools}
\usepackage{microtype}
\usepackage{enumitem}
\usepackage{needspace}
\usepackage{xcolor}
\usepackage{booktabs}
\usepackage{graphicx}
\usepackage{tikz}
\usetikzlibrary{arrows.meta,calc,positioning,shapes.misc,shapes.geometric,
  patterns,decorations.pathreplacing,fit,backgrounds,intersections}
\usepackage{pgfplots}
\pgfplotsset{compat=1.18}
\usepackage[font=small,labelfont=bf]{caption}
\usepackage{subcaption}
\usepackage{thmtools,thm-restate}
\usepackage[colorlinks=true,linkcolor=blue!50!black,citecolor=green!35!black,
  urlcolor=blue!50!black,bookmarksdepth=2]{hyperref}

\input{macros}

\input{statements}

\ifdefined\PAsinglesetup\PAsinglesetup\fi

\title{TAP states at upper deviations of the free energy:\\
existence, localization and marginal stability in Ising spin glasses}
\author{Yan Ru Pei}
\date{\vspace{-6ex}}
\hypersetup{pdftitle={TAP states at upper deviations of the free energy: existence,
  localization and marginal stability in Ising spin glasses},pdfauthor={Yan Ru Pei}}

\begin{document}

\ifdefined\PAsinglebody
  \PAsinglebody
\else

\maketitle

\begin{abstract}
We study upper large deviations of the free energy of mixed $p$-spin Ising spin glasses through
the generalized TAP free energy $\FTAP$ of Chen, Panchenko and Subag. For every mixture with
radius of convergence greater than one, $\max\FTAP$ and $\log\ZN$ have the same large deviations
at speed $N$; the proof combines their band theorem with Ramsey's theorem. For convex mixtures
with $\sum_p2^p\beta_p^2<\infty$, upper deviations at level $f$ are carried by generalized TAP
critical points, which exist at the free-energy rate without Boursier's strict Plefka condition.
Outside an event
of smaller exponential order, near-maximizers at level $f$ lie near the contact set of a
constrained Parisi obstacle, follow the Auffinger--Chen field law, and have Hessian bulk near a
reflected Pastur law whose edge is nonpositive and vanishes exactly at contacts where the
obstacle is flat to second order. Where every contact is of this kind, a fixed stability margin
costs rate.
\end{abstract}

\input{sections/sec1-introduction}

\setcounter{tocdepth}{1}
\begingroup
\makeatletter
\renewcommand*\l@section[2]{%
  \ifnum\c@tocdepth>\z@
    \addpenalty\@secpenalty
    \addvspace{0.45em \@plus\p@}%
    \setlength\@tempdima{1.5em}%
    \begingroup
      \parindent\z@ \rightskip\@pnumwidth
      \parfillskip-\@pnumwidth
      \leavevmode\bfseries
      \advance\leftskip\@tempdima
      \hskip-\leftskip
      #1\nobreak\hfil
      \nobreak\hb@xt@\@pnumwidth{\hss #2\kern-\p@\kern\p@}\par
    \endgroup
  \fi}
\makeatother
\tableofcontents
\endgroup

\input{sections/sec2-maximum}
\input{sections/sec3-localization}
\input{sections/sec4-hessian}
\input{sections/sec5-stability}
\input{sections/sec6-equilibrium}

\input{sections/sec7-open}

\appendix
\input{sections/appA-parisi}
\input{sections/appB-glued}
\input{sections/appC-gaussian}
\input{sections/appD-matrix}
\input{sections/appE-kacrice}
\input{sections/appF-nearcritical}
\input{sections/appG-moments}

\fi

\input{refs}
\end{document}

%% file: macros.tex
\theoremstyle{plain}

\newtheorem{proposition}{Proposition}[section]
\newtheorem{lemma}[proposition]{Lemma}
\newtheorem{corollary}[proposition]{Corollary}

\newtheorem*{theorem*}{Theorem}
\theoremstyle{definition}
\newtheorem{definition}[proposition]{Definition}
\newtheorem{remark}[proposition]{Remark}

\newtheorem*{remark*}{Remark}
\counterwithin{equation}{section}

\makeatletter
\@for\PA@env:=proposition,lemma,corollary,claim,definition,remark,example\do{%
  \expandafter\def\csname theH\PA@env\endcsname{\theHsection.\arabic{proposition}}}
\makeatother

\newlist{thmlist}{enumerate}{1}
\setlist[thmlist]{label=\textup{(\alph*)}, ref=(\alph*), leftmargin=*, widest=b,
  topsep=3pt, itemsep=2pt}
\newlist{romanlist}{enumerate}{1}
\setlist[romanlist]{label=\textup{(\roman*)}, ref=(\roman*), leftmargin=*, widest=iii,
  topsep=3pt, itemsep=2pt}
\newlist{itemlist}{enumerate}{2}
\setlist[itemlist]{label=\textup{\arabic*.}, ref=\arabic*, leftmargin=*, widest=3,
  topsep=2pt, itemsep=1pt}
\newlist{steplist}{enumerate}{1}
\setlist[steplist]{label=\textup{(\arabic*)}, ref=(\arabic*), leftmargin=*, widest=3,
  topsep=3pt, itemsep=2pt}

\newcommand{\E}{\mathbb{E}}
\newcommand{\Prob}{\mathbb{P}}
\newcommand{\Pout}{\mathbb{P}^{*}}            %
\newcommand{\Pin}{\mathbb{P}_{*}}            %
\newcommand{\R}{\mathbb{R}}
\newcommand{\N}{\mathbb{N}}
\newcommand{\Cplx}{\mathbb{C}}
\newcommand{\ind}{\mathbf{1}}
\newcommand{\id}{\mathrm{id}}
\DeclareMathOperator{\supp}{supp}
\DeclareMathOperator{\sech}{sech}
\DeclareMathOperator{\Law}{Law}
\DeclareMathOperator{\Shift}{Shift}
\DeclareMathOperator{\dist}{dist}
\DeclareMathOperator{\diag}{diag}
\DeclareMathOperator{\ESD}{ESD}
\DeclareMathOperator{\Var}{Var}
\DeclareMathOperator{\Cov}{Cov}
\DeclareMathOperator{\rank}{rank}
\DeclareMathOperator{\tr}{tr}
\DeclareMathOperator*{\argmin}{arg\,min}

\DeclareMathOperator*{\esssup}{ess\,sup}
\DeclareMathOperator*{\essinf}{ess\,inf}
\newcommand{\dL}{d_{\mathrm{L}}}              %
\newcommand{\Wone}{W_1}                       %
\newcommand{\Pspace}{\mathcal{P}}             %
\newcommand{\GOE}{W_N}                        %

\newcommand{\HN}{H_N}                         %
\newcommand{\ZN}{Z_N}                         %
\newcommand{\cube}{[-1,1]^N}
\newcommand{\ocube}{(-1,1)^N}
\newcommand{\hcube}{\{-1,1\}^N}
\newcommand{\ov}[2]{R(#1,#2)}                 %
\newcommand{\MN}{M_N}                         %
\newcommand{\FTAP}{F_{\mathrm{TAP}}}          %
\newcommand{\TAPcps}{\mathrm{TAP}}            %
\newcommand{\Phistar}{\Phi^{*}}               %
\newcommand{\Qth}{\mathbb{Q}_\theta}          %
\newcommand{\Band}{B_n(m,\varepsilon,\delta)} %
\newcommand{\glue}{\mathrm{glue}}             %
\newcommand{\xc}{\mathsf{x}}                  %
\newcommand{\cdf}[1]{\bar{#1}}                %
\newcommand{\vth}{\vartheta}                  %

\newcommand{\feq}{f_{\mathrm{eq}}}            %
\newcommand{\fone}{f_1}                       %
\newcommand{\fzero}{f_0}                      %
\newcommand{\Lambar}{\overline{\Lambda}}      %
\newcommand{\thf}{\theta_f}                   %
\newcommand{\zf}{\zeta_f}                     %
\newcommand{\zP}{\zeta_{\mathrm{P}}}          %
\newcommand{\thstar}{\theta_*}                %
\newcommand{\Obs}{\mathcal{O}}                %
\newcommand{\Omin}{\mathcal{O}_{\min}}        %
\newcommand{\Cf}{C_f}                         %
\newcommand{\qf}{q_f}                         %
\newcommand{\OP}{\mathcal{O}_{\mathrm{P}}}    %
\newcommand{\CP}{C_{\mathrm{P}}}              %
\newcommand{\Lq}[1]{L_{#1}}                   %
\newcommand{\LPq}[1]{L^{\mathrm{P}}_{#1}}     %
\newcommand{\Gz}[1]{\Gamma_{#1}}              %

\newcommand{\Hess}{\nabla^2}
\newcommand{\Hpm}{\nabla^2_{\pm}}             %
\newcommand{\aQ}{a_Q}                         %
\newcommand{\UQ}{U_Q}                         %
\newcommand{\chiQ}{\chi_Q}                    %
\newcommand{\Plef}{\Delta}                    %
\newcommand{\rhoQ}{\rho_Q}                    %
\newcommand{\nuQ}{\nu_Q}                      %
\newcommand{\Pastur}[2]{\nu[#1,#2]}           %
\newcommand{\lamp}{\lambda_+}                 %

\newcommand{\Eex}[2]{\mathcal{E}_{N,#1}(#2)}                     %
\newcommand{\Emarg}[3]{\mathcal{E}^{#1}_{N,#2}(#3)}               %
\newcommand{\Efix}[4]{\mathcal{E}^{#1,#2}_{N,#3,\mathrm{fix}}(#4)}%

\newcommand{\classC}{\mathcal{C}}             %
\newcommand{\condR}{\hyperref[hyp:R]{\textup{(R)}}}

\newcommand{\condH}{\hyperref[hyp:H]{\textup{(H)}}}
\newcommand{\condFM}{\hyperref[hyp:FM]{\textup{(FM)}}}
\newcommand{\condHmarg}{\hyperref[hyp:Hmarg]{\textup{(H-marg)}}}
\newcommand{\FM}[1]{\mathrm{FM}(\xi,#1)}      %
\newcommand{\FMge}[1]{\mathrm{FM}^{\ge}(\xi,#1)}

%% file: statements.tex
\newcommand{\DisplayCondC}{%
\begin{equation*}\tag{$\mathcal{C}$}\label{hyp:C}
 \sum_{p\ge2}2^p\beta_p^2<\infty .
\end{equation*}}
\newcommand{\DisplayCondH}{%
\begin{equation*}\tag{H}\label{hyp:H}
 \xi\ \text{satisfies \condR\ and is convex on } [-1,1].
\end{equation*}}
\newcommand{\DisplayCondFM}{%
\begin{equation*}\tag{FM}\label{hyp:FM}
 \lim_{N\to\infty}\frac1N\log\E\ZN^{\theta}=\Lambda(\theta)
 \qquad\text{for every } \theta\in(0,1).
\end{equation*}}
\newcommand{\DisplayCondHmarg}{%
\begin{equation*}\tag{H-marg}\label{hyp:Hmarg}
 \Plef(Q)=0\ \text{ for every } Q\in\Cf,
 \quad\text{equivalently}\quad \Obs_{\zf}''=0\ \text{ on } \Cf .
\end{equation*}}

\newcommand{\DisplayEventExist}{%
For $\varepsilon>0$ and $f\in\R$ let $\Eex{\varepsilon}{f}$ be the event that some
$\zeta\in\Pspace([0,1])$ and $x\in\ocube$ satisfy
\begin{equation}\label{eq:event-exist}
 \nabla F_\zeta(x)=0,\qquad
 \bigl|N^{-1}F_\zeta(x)-f\bigr|\le\varepsilon,\qquad
 \bigl|N^{-1}\FTAP(x)-f\bigr|\le\varepsilon .
\end{equation}}
\newcommand{\DisplayEventsStab}{%
For $g,\rho>0$ let $\Emarg{g}{\varepsilon}{f}\subset\Eex{\varepsilon}{f}$ be the event
that some pair $(\zeta,x)$ satisfies \eqref{eq:event-exist} and
$\nabla^2_\sigma F_\zeta(x)\preceq-gI$ for some $\sigma\in\{+,-\}$, and let
$\Efix{g}{\rho}{\varepsilon}{f}\subset\Emarg{g}{\varepsilon}{f}$ be the event that some
$x\in\ocube$ satisfies \eqref{eq:event-exist} with $\zeta=\zf$ and
\begin{equation}\label{eq:event-fix}
 \Hpm F_{\zf}(y)\preceq-gI\qquad\text{for every } y\in\ocube \text{ with } |y-x|\le\rho\sqrt N .
\end{equation}}

\newcommand{\StateThmMax}{%
\begin{restatable}[The TAP maximum carries the free energy]{theorem}{thmmax}
\label{thm:max}
Assume \condR. The maximum $\MN=\max_{\cube}\FTAP$ is attained almost surely, and the
following hold.
\begin{thmlist}
\item\label{thm:max-a} For every $t>0$,
 \[
  \lim_{N\to\infty}\frac1N\log\Prob\bigl(|\MN-\log\ZN|>tN\bigr)=-\infty .
 \]
\item\label{thm:max-b} For every $\theta\in\R$,
 \[
  \lim_{N\to\infty}\frac1N\Bigl[\log\E e^{\theta\MN}-\log\E\ZN^{\theta}\Bigr]=0 .
 \]
\item\label{thm:max-c} For every $f\in\R$, $0<t<\varepsilon$ and $c>0$, for all large $N$,
 \begin{align*}
  \Prob\bigl(|N^{-1}\MN-f|\le\varepsilon\bigr)
  &\ge\Prob\bigl(|N^{-1}\log\ZN-f|\le\varepsilon-t\bigr)-e^{-cN},\\
  \Prob\bigl(|N^{-1}\MN-f|\le\varepsilon\bigr)
  &\le\Prob\bigl(|N^{-1}\log\ZN-f|\le\varepsilon+t\bigr)+e^{-cN}.
 \end{align*}
 Consequently, for every $f\in\R$, the limits as $\varepsilon\downarrow0$ of
 $\liminf_N N^{-1}\log\Prob(|N^{-1}\MN-f|\le\varepsilon)$ and of
 $\limsup_N N^{-1}\log\Prob(|N^{-1}\MN-f|\le\varepsilon)$ equal the corresponding limits
 for $N^{-1}\log\ZN$. Moreover, $N^{-1}\MN$ satisfies a large deviation principle at
 speed $N$ with a good rate function if and only if $N^{-1}\log\ZN$ does, and then the two
 rate functions coincide.
\end{thmlist}
\end{restatable}}

\newcommand{\StateThmExist}{%
\begin{restatable}[Generalized TAP critical points at every upper level]{theorem}{thmexist}
\label{thm:exist}
Assume $\xi\in\classC$.
\begin{thmlist}
\item\label{thm:exist-a} Almost surely, for every $N$, every maximizer $x^*$ of $\FTAP$ over
 $\cube$ lies in $\ocube$. Let $\zeta^*$ be the unique minimizer of $\mu\mapsto F_\mu(x^*)$
 over $\Pspace([q_{x^*},1])$. Then $F_{\zeta^*}$ is continuously differentiable on $\ocube$,
 \[
  \nabla F_{\zeta^*}(x^*)=0
  \qquad\text{and}\qquad
  F_{\zeta^*}(x^*)=\FTAP(x^*)=\MN .
 \]
\item\label{thm:exist-b} The function $\Lambda$ is continuously differentiable on $(0,1)$, with
 $\Lambda'$ nondecreasing and equal to $\feq$ on $(0,\thstar]\cap(0,1)$, where
 $\thstar=\zP(\{0\})$. If $\thstar<1$, then $\Lambda'$ is a strictly increasing
 homeomorphism of $(\thstar,1)$ onto $(\feq,\fone)$; if $\thstar=1$, then $\fone=\feq$.
 Every $f\in(\feq,\fone)$ is regular: the minimizer $\thf$ of
 $\theta\mapsto\Lambar(\theta)-\theta f$ over $[0,1]$ is unique, $\thf\in(\thstar,1)$,
 $\Lambda'(\thf)=f$ and $\Sigma(f)=\Lambda(\thf)-\thf f$.
\item\label{thm:exist-c} If $\FM{\theta}$ holds at some $\theta\in(0,1)$, then, at that
 $\theta$,
 \[
  \lim_{N\to\infty}\frac1N\log\E e^{\theta\MN}=\Lambda(\theta).
 \]
 Moreover, with no assumption on fractional moments,
 $\lim_N N^{-1}\log\E e^{\MN}=\xi(1)/2=\Lambar(1)$.
\item\label{thm:exist-d} Assume \condFM\ and let $f\in(\feq,\fone)$. Then
 \[
  \lim_{\varepsilon\downarrow0}\liminf_{N\to\infty}\frac1N
   \log\Prob\bigl(|N^{-1}\MN-f|\le\varepsilon\bigr)
  =\lim_{\varepsilon\downarrow0}\limsup_{N\to\infty}\frac1N
   \log\Prob\bigl(|N^{-1}\MN-f|\le\varepsilon\bigr)
  =\Sigma(f),
 \]
 and the same holds with the event $\{|N^{-1}\MN-f|\le\varepsilon\}$ replaced by
 $\Eex{\varepsilon}{f}$, for its inner and for its outer probability.
\item\label{thm:exist-e} If $\xi$ is convex on $[-1,1]$, then \condFM\ holds. Hence
 \ref{thm:exist-c} holds at every $\theta\in(0,1)$, and \ref{thm:exist-d} holds at every
 $f\in(\feq,\fone)$.
\end{thmlist}
\end{restatable}}

\newcommand{\StateThmStruct}{%
\begin{restatable}[Structure of the states that carry an upper deviation]{theorem}{thmstruct}
\label{thm:struct}
Assume \condH\ and let $f\in(\feq,\fone)$. Let $\zf$ be the minimizer of $\zeta\mapsto P_\zeta$
over $\{\zeta\in\Pspace([0,1]):\zeta(\{0\})\ge\thf\}$ (it is unique by Lemma~\ref{lem:var-basic}),
and let $\Cf$, $\Lq{q}$, $\nuQ$,
$\lamp$ and $\Plef$ be defined from $\zf$ as in Section~\ref{sec:intro-structure}.
\begin{romanlist}
\item\label{thm:struct-i} \textup{(Localization.)} For every $\kappa>0$ there are
 $\delta,c>0$ such that, for all large $N$,
 \[
  \Pout\bigl(\exists m\in\cube:\ \FTAP(m)\ge N(f-\delta),\ \dist(q_m,\Cf)\ge\kappa\bigr)
  \le e^{N(\Sigma(f)-c)} .
 \]
\item\label{thm:struct-ii} \textup{(Field law.)} For every $\kappa>0$ there are
 $\delta,c>0$ such that, for all large $N$,
 \[
  \Pout\bigl(\exists m\in\cube:\ \FTAP(m)\ge N(f-\delta),\ \Wone(\mu_m,\Lq{q_m})\ge\kappa\bigr)
  \le e^{N(\Sigma(f)-c)} .
 \]
\end{romanlist}
Assume moreover $\xi\in\classC$. Put $\kappa_f=\dist(\Cf,\{0,1\})$,
which is positive.
\begin{romanlist}[resume]
\item\label{thm:struct-iii} \textup{(Hessian bulk.)} For every $\kappa\in(0,\kappa_f)$ there
 are $\varepsilon,c>0$ such that, for all large $N$, outside an event of outer probability at
 most $e^{N(\Sigma(f)-c)}$, every $x\in\ocube$ with $\FTAP(x)\ge N(f-\varepsilon)$ satisfies
 $\dist(q_x,\Cf)<\kappa$ and
 \[
  \dL\bigl(\ESD(\Hpm F_\alpha(x)),\nu_{q_x}\bigr)\le\kappa
 \]
 both for $\alpha=\zf$ and for every $\alpha\in\Pspace([0,1])$ with
 $F_\alpha(x)\le\FTAP(x)+2\varepsilon N$. The same holds for the compressions of these
 Hessians to the orthogonal complement of $x$.
\item\label{thm:struct-iv} \textup{(Edge.)} For every $Q\in(0,1)$, $\lamp(Q)\le0$, with
 equality if and only if $\Plef(Q)=0$. On $\Cf$ one has $\Plef=2\Obs_{\zf}''/\xi''\ge0$, so
 a contact $Q\in\Cf$ has $\lamp(Q)=0$ if and only if $\Obs_{\zf}''(Q)=0$.
\end{romanlist}
\end{restatable}}

\newcommand{\StateCorTitle}{%
\begin{restatable}[Upper deviations of the free energy are carried by TAP states]{maincorollary}{cortitle}
\label{cor:title}
Assume \condH, $\xi\in\classC$ and $f\in(\feq,\fone)$. Let $\zf$ be the minimizer of
$\zeta\mapsto P_\zeta$ over $\{\zeta(\{0\})\ge\thf\}$ (it is unique by
Lemma~\ref{lem:var-basic}), and let $\Cf$, $\Lq{q}$ and $\nuQ$ be defined from $\zf$.
For a maximizer $x^*$ of $\FTAP$ let $\zeta^*$ be as in
Theorem~\ref{thm:exist}\ref{thm:exist-a}. Conditional outer probabilities are
$\Pout(B\mid A)=\Pout(B\cap A)/\Pin(A)$ for events $A$ with $\Pin(A)>0$; for measurable $A$ this is
$\Pout(B\cap A)/\Prob(A)$.
\begin{thmlist}
\item\label{cor:title-a} For every $\kappa>0$ there are $\delta,c>0$ such that for every
 $t>0$ there is $N_0$ with the following property. For $N\ge N_0$,
 $\Prob(N^{-1}\log\ZN\ge f)>0$, and conditionally on $\{N^{-1}\log\ZN\ge f\}$, outside an
 event of conditional outer probability at most $e^{-cN}$:
 \begin{itemlist}
 \item\label{cor:title-1} $|\MN-\log\ZN|\le tN$, and every maximizer $x^*$ of $\FTAP$ lies in
  $\ocube$ and satisfies $\nabla F_{\zeta^*}(x^*)=0$ and
  $F_{\zeta^*}(x^*)=\FTAP(x^*)=\MN$;
 \item\label{cor:title-2} every $m\in\cube$ with $\FTAP(m)\ge N(f-\delta)$ has
  $\dist(q_m,\Cf)<\kappa$ and $\Wone(\mu_m,\Lq{q_m})<\kappa$;
 \item\label{cor:title-3} every such $m$ in $\ocube$ has $q_m\in(0,1)$, and every matrix $H$
  in the following list satisfies $\dL(\ESD(H),\nu_{q_m})\le\kappa$ and has at most
  $\kappa N$ eigenvalues larger than $\kappa$: the one-sided Hessians $\Hpm F_{\zf}(m)$; the
  one-sided Hessians $\Hpm F_\alpha(m)$ for every $\alpha\in\Pspace([0,1])$ with
  $F_\alpha(m)\le\FTAP(m)+2\delta N$; and the compressions of all of these to the orthogonal
  complement of $m$.
 \end{itemlist}
 In particular, if $t\le\delta$, every maximizer $x^*$ of $\FTAP$ satisfies
 items~\ref{cor:title-2} and~\ref{cor:title-3}, the latter with $\alpha=\zeta^*$.
\item\label{cor:title-b} For every $\kappa>0$ there are $\delta,c>0$ and $N_0$ such that,
 for $N\ge N_0$, $\Prob(\MN\ge N(f-\delta))>0$, and conditionally on
 $\{\MN\ge N(f-\delta)\}$, items~\ref{cor:title-2} and~\ref{cor:title-3} hold outside an
 event of conditional outer probability at most $e^{-cN}$.
\end{thmlist}
\end{restatable}}

\newcommand{\StateThmStab}{%
\begin{restatable}[Stability costs rate at marginal levels]{theorem}{thmstab}
\label{thm:stab}
Assume \condH\ and $\xi\in\classC$, and let $f\in(\feq,\fone)$.
\begin{thmlist}
\item\label{thm:stab-a} \textup{(Kac--Rice construction.)} Assume that $\xi$ has at least
 two active degrees. For every $\varepsilon>0$ there are $g,\rho>0$ such that
 \[
  \liminf_{N\to\infty}\frac1N\log\Prob\bigl(\Efix{g}{\rho}{\varepsilon}{f}\bigr)
  \ge\Sigma(f)-\varepsilon .
 \]
 If the first positive contact $\qf=\min\Cf$ is strictly Plefka stable, that is
 $\Plef(\qf)>0$, then $g$ and $\rho$ can be chosen independently of $\varepsilon$, so that at
 a fixed margin
 \[
  \lim_{\varepsilon\downarrow0}\liminf_{N\to\infty}\frac1N
   \log\Prob\bigl(\Efix{g}{\rho}{\varepsilon}{f}\bigr)=\Sigma(f).
 \]
\item\label{thm:stab-b} \textup{(Fixed margins are exponentially rarer.)} Assume
 \condHmarg. For every $g>0$ there are $\varepsilon,c>0$ such that, for all large $N$, the
 outer probability that some $x\in\ocube$ and $\alpha\in\Pspace([0,1])$ satisfy
 \[
  \FTAP(x)\ge N(f-\varepsilon),\qquad
  F_\alpha(x)\le\FTAP(x)+2\varepsilon N,\qquad
  \min_{\pm}\lambda_{\max}\bigl(\Hpm F_\alpha(x)\bigr)\le-g
 \]
 is at most $e^{N(\Sigma(f)-c)}$, and the same holds with $\alpha=\zf$ and without the second
 condition. In particular, for every $\varepsilon'\in(0,\varepsilon]$,
 \[
  \limsup_{N\to\infty}\frac1N\log\Prob\bigl(\Emarg{g}{\varepsilon'}{f}\bigr)\le\Sigma(f)-c .
 \]
\item\label{thm:stab-c} \textup{(A near-critical family.)} There is $\varepsilon_*>0$ such
 that every covariance
 \[
  \xi(t)=bt^2+et^4,\qquad b>\tfrac12,\qquad e>0,\qquad b-\tfrac12+e<\varepsilon_*,
 \]
 admits $\fzero\in(\feq,\fone]$ with the following properties at every
 $f\in(\feq,\fzero)$.
 \begin{itemlist}
 \item\label{thm:stab-c1} $\Cf=\supp\zf\cap(0,1]$ is an interval of positive length, and
  \condHmarg\ holds.
 \item\label{thm:stab-c2} If $g(\varepsilon),\rho(\varepsilon)>0$ satisfy the conclusion
  of \ref{thm:stab-a} for every $\varepsilon>0$, then $g(\varepsilon)\to0$ as
  $\varepsilon\downarrow0$.
 \item\label{thm:stab-c3} For every sufficiently small $\kappa>0$ there are
  $\kappa',\varepsilon,c>0$ such that, for all large $N$, the inner probability $\Pin(\Eex{\varepsilon}{f})$ is positive and, conditionally on
  $\Eex{\varepsilon}{f}$ and outside an event of conditional outer probability
  (Corollary~\ref{cor:title}) at most
  $e^{-cN}$, every pair $(\alpha,x)$ satisfying \eqref{eq:event-exist} has
  $\dist(q_x,\Cf)<\kappa$ and $\Wone(\mu_x,\Lq{q_x})<\kappa$, and each of the matrices
  $\Hpm F_\alpha(x)$, $\Hpm F_{\zf}(x)$ and their compressions to the orthogonal complement
  of $x$ has its $\lceil\kappa'N\rceil$-th largest eigenvalue in $[-\kappa,\kappa]$ and
  fewer than $\kappa'N$ eigenvalues larger than $\kappa$.
 \end{itemlist}
\end{thmlist}
\end{restatable}}

\newcommand{\StateCorEq}{%
\begin{restatable}[The equilibrium level]{maincorollary}{coreq}
\label{cor:eq}
Assume \condH. Let $\zP$ be a minimizer of $\zeta\mapsto P_\zeta$ over $\Pspace([0,1])$, let
$\OP$ be its obstacle, $\CP=\argmin_{[0,1]}\OP$, and let $\LPq{q}$ be its field law. Then
$\CP$ contains $0$ and $\supp\zP$. For every $\kappa>0$ there is $c>0$ such that, for all
large $N$,
\[
 \Pout\Bigl(\exists m\in\cube:\ \FTAP(m)\ge N(\feq-c)\ \text{and}\
 \bigl[\dist(q_m,\CP)\ge\kappa\ \text{or}\ \Wone(\mu_m,\LPq{q_m})\ge\kappa\bigr]\Bigr)
 \le e^{-cN}.
\]
In particular, for every $Q\in[0,1]\setminus\CP$ there is $c(Q)>0$ such that, for all large
$N$, $\sup\{\FTAP(m):q_m=Q\}\le N(\feq-c(Q))$ with probability at least $1-e^{-c(Q)N}$.
\end{restatable}}

%% file: sections/sec1-introduction.tex
\section{Introduction}\label{sec:intro}

Let $\HN$ be the centered Gaussian field on $\cube$ with $\E\HN(x)\HN(y)=N\xi(\ov xy)$, where
$\xi(t)=\sum_{p\ge2}\beta_p^2t^p\not\equiv0$ and $\ov xy=N^{-1}x\cdot y$; on $\hcube$ it is the
mixed $p$-spin Hamiltonian. The free energy $N^{-1}\log\ZN$, $\ZN=2^{-N}\sum_{\sigma}e^{\HN(\sigma)}$,
converges to the Parisi value $\feq$ (for instance when $\sum_p2^p\beta_p^2<\infty$, by the Parisi
formula \cite{PF,Pan14}); its upper deviations are exponentially rare. The maximum of the
generalized TAP free energy $\FTAP$ of Chen, Panchenko and Subag over a sphere $\{\ov mm=q\}$,
$q$ in the support of the Parisi measure, is $\log\ZN+o(N)$ in probability \cite[Theorem~2]{CPS}.
Boursier observed that an upper deviation of the free energy should be realized by an upper deviation of
$\MN=\max_{\cube}\FTAP$, hence by TAP states at the deviating level
\cite[Remark~1.3]{Bou26b}. We say that $\xi$ satisfies \phantomsection\label{hyp:R}\textup{(R)}
if its radius of convergence exceeds one.

\StateThmMax

For $\xi\in\classC$ (Section~\ref{sec:intro-setting}) the upper side of~\ref{thm:max-a} is
essentially in \cite{CPS}, at every rate: the band free energy never exceeds the full one
\cite[(2.16)]{CPS}, so the band theorem \cite[Theorem~1]{CPS} gives $\MN\le\log\ZN+tN$ outside an
event of probability $e^{-cN}$, for every $c$. Theorem~2 of \cite{CPS} and Gaussian concentration
already give $|\MN-\log\ZN|\le tN$ outside an event of probability $e^{-c(t)N}$ for some
$c(t)>0$. The point of Theorem~\ref{thm:max} is that every exponential rate is allowed, which is
what lets~\ref{thm:max-c} transfer large deviations in both directions. To our knowledge, this
superexponential form has not been recorded. The proof (Section~\ref{sec:max}) uses Ramsey's
theorem \cite{Ram30} to find, for every realization of the disorder, an overlap level at which the
Gibbs measure puts mass bounded below independently of $N$ into the band around a barycenter of
replicas; in \cite[Section~5]{CPS} the level is fixed in the support of the Parisi measure and the
band mass is controlled on an event whose probability is not exponentially small
\cite[(5.2)--(5.3)]{CPS}, which suffices in probability. Theorem~\ref{thm:max} does not say where
the maximizers of $\FTAP$ lie, at which self-overlap, or whether they are critical points of a TAP
function: the barycenter in its proof sits at an overlap level that depends on the realization. The other five
results answer these questions.

\subsection{Setting and generalized TAP critical points}\label{sec:intro-setting}

\paragraph{Covariances.} Throughout, $\xi(t)=\sum_{p\ge2}\beta_p^2t^p\not\equiv0$ and there is no
external field. We say that $\xi$ has \emph{at least two active degrees} if $\beta_p\ne0$ for at
least two values of $p$. Besides \condR\ we use the following two conditions.
\DisplayCondC
\DisplayCondH
The class $\classC$ is the model class of \cite{CPS}; it contains every polynomial, and
$\xi\in\classC$ forces radius of convergence at least two, hence \condR. Under \condR, $\xi''$ is
analytic near $[-1,1]$ and not identically zero, so it has finitely many zeros there and a convex
$\xi$ is strictly convex; thus \condH\ is the standing assumption of Boursier
\cite[Assumption~1.1]{Bou26b}. Every even $\xi$ satisfying \condR\ is convex on $[-1,1]$, since
$\xi''$ is then a series of even powers with nonnegative coefficients. Odd degrees are allowed in \condH; a convex $\xi$ with an odd
degree has at least two active degrees, because a single odd monomial is not convex on $[-1,1]$. The convex covariances with one active
degree are the Sherrington--Kirkpatrick (SK) model $\xi(t)=\beta^2t^2$ and the pure even $p$-spin
models.

Theorem~\ref{thm:max} assumes only \condR;
Theorem~\ref{thm:struct}\ref{thm:struct-i}--\ref{thm:struct-ii} and Corollary~\ref{cor:eq}
assume \condH; Theorem~\ref{thm:exist} assumes $\xi\in\classC$; and
Theorem~\ref{thm:struct}\ref{thm:struct-iii}--\ref{thm:struct-iv}, Corollary~\ref{cor:title}
and Theorem~\ref{thm:stab} assume both \condH\ and $\xi\in\classC$. The class $\classC$ enters
only through cited facts on the Parisi equation and the TAP correction (Remark~\ref{rem:rpde});
part~\ref{thm:exist-b} of Theorem~\ref{thm:exist} and the uniqueness of the constrained
minimizers defined below need only \condR\ (Lemma~\ref{lem:var-basic} and
Proposition~\ref{prop:var-regular}).

\paragraph{The Parisi functional and the TAP functions.} We follow the notation and
normalization of Boursier \cite{Bou26b}. For $\zeta\in\Pspace([0,1])$ let $\Phi_\zeta$ solve the
Parisi equation
\begin{equation}\label{eq:intro-parisi-pde}
 \partial_t\Phi_\zeta+\frac{\xi''(t)}2\Bigl(\partial_x^2\Phi_\zeta
 +\zeta([0,t])\,(\partial_x\Phi_\zeta)^2\Bigr)=0
 \quad\text{on }[0,1)\times\R,\qquad \Phi_\zeta(1,x)=\log\cosh x .
\end{equation}
For $\zeta$ with finitely many atoms the solution is given by the Cole--Hopf transformation; in
general it is defined by continuity in $\zeta$ \cite{AC15a} (see
Appendix~\ref{app:pde-toolbox}). Let $\Phistar_\zeta(q,m)=\inf_{b\in\R}\{\Phi_\zeta(q,b)-mb\}$
for $m\in(-1,1)$, extended continuously to $m=\pm1$, and let
\begin{equation}\label{eq:intro-parisi-functional}
 P_\zeta=\Phi_\zeta(0,0)-\frac12\int_0^1t\,\xi''(t)\,\zeta([0,t])\,dt
\end{equation}
be the Parisi functional. For $x\in\cube$ let $q_x=\ov{x}{x}$ be its self-overlap and
$\mu_x=N^{-1}\sum_i\delta_{x_i}$ the empirical law of its coordinates. The TAP function with
order parameter $\zeta$ and the generalized TAP free energy are
\begin{align}
 F_\zeta(x)&=\HN(x)+\sum_{i=1}^N\Phistar_\zeta(q_x,x_i)
 -\frac N2\int_{q_x}^1t\,\xi''(t)\,\zeta([0,t])\,dt,\label{eq:intro-tap-function}\\
 \FTAP(x)&=\inf_{\mu\in\Pspace([q_x,1])}F_\mu(x)=\inf_{\zeta\in\Pspace([0,1])}F_\zeta(x).
 \label{eq:intro-ftap}
\end{align}
The two infima agree because $F_\zeta(x)$ depends on $\zeta$ only through $t\mapsto\zeta([0,t])$
on $[q_x,1]$, so $F_\zeta=F_{\Shift_{q_x}\zeta}$ with
$\Shift_q\zeta=\zeta|_{(q,1]}+\zeta([0,q])\delta_q$. The paper \cite{CPS} uses the terminal value
$\log(2\cosh x)$ and the partition function without the factor $2^{-N}$, which raises $P_\zeta$ by
$\log2$, $F_\zeta$ and $\FTAP$ by $N\log2$, and the free energy by $\log2$; thus
$\FTAP(m)=\HN(m)+N\TAPcps(\mu_m)-N\log2$, with $\TAPcps$ the functional of
\cite[(2.24)--(2.25)]{CPS}. Formulas quoted from \cite{CPS,CPS2} are translated accordingly.

\paragraph{Equilibrium and constrained Parisi problems.} Let $\feq=\inf_\zeta P_\zeta$. For
$\xi\in\classC$, $N^{-1}\log\ZN\to\feq$ almost surely \cite{PF,Pan14}, and the minimizer $\zP$,
the Parisi measure, is unique by the strict convexity of the Parisi functional
\cite[Theorem~2]{AC15b}, \cite[Theorem~2]{JT16}; put $\thstar=\zP(\{0\})$. For $\theta\in[0,1)$
let
\[
 p(\theta)=\inf\bigl\{P_\zeta:\ \zeta\in\Pspace([0,1]),\ \zeta(\{0\})\ge\theta\bigr\},
 \qquad \Lambda(\theta)=\theta\,p(\theta),
\]
let $\Lambar$ be the closed extension of $\Lambda$ (equal to $\Lambda$ on $[0,1)$, to
$\liminf_{u\uparrow1}\Lambda(u)$ at $1$, and to $+\infty$ off $[0,1]$), and put
\begin{equation}\label{eq:intro-rate}
 \Sigma(f)=\inf_{\theta\in[0,1]}\bigl\{\Lambar(\theta)-\theta f\bigr\},\qquad
 \fone=\Lambar'_-(1).
\end{equation}
These are the objects of \cite[Section~1.3]{Bou26b}. For $\xi$ satisfying \condR\ outside
$\classC$ we use $\feq$ and $\Sigma$ as variational quantities;
Theorem~\ref{thm:struct}\ref{thm:struct-i}--\ref{thm:struct-ii} and Corollary~\ref{cor:eq} are
stated relative to them. The interval $(\feq,\fone)$ is nonempty exactly when $\thstar<1$
(Proposition~\ref{prop:var-regular}\ref{prop:var-regular-a}); it is empty when $\zP=\delta_0$. For
$f\ge\fone$ the infimum in~\eqref{eq:intro-rate} is attained at $\theta=1$, so
$\Sigma(f)=\Lambar(1)-f=\xi(1)/2-f$ (Lemma~\ref{lem:var-dual} and
Proposition~\ref{prop:var-regular}\ref{prop:var-regular-c}); at these levels we have only the upper
bound $\Prob(N^{-1}\log\ZN\ge f)\le e^{N\Sigma(f)}$, which holds for every $N$ by Markov's
inequality and $\E\ZN=e^{N\xi(1)/2}$ (for convex $\xi\in\classC$ it is also the case $f\ge\fone$ of
Proposition~\ref{prop:var-tail}\ref{prop:var-tail-e}), and we treat $f\in(\feq,\fone)$. As an infimum of affine functions of $f$, $\Sigma$ is concave, and
$\Sigma\le\Lambar(0)=0$. We call $\Sigma(f)$ the exponent of the level; the rate function is
$-\Sigma$. The link with the free energy is the \emph{fractional-moment formula}: we say that
$\FM{\theta}$ holds if $\lim_NN^{-1}\log\E\ZN^\theta=\Lambda(\theta)$, and that $\FMge{\theta}$
holds if $\liminf_NN^{-1}\log\E\ZN^\theta\ge\Lambda(\theta)$, and we write
\DisplayCondFM
Talagrand computed the fractional moments for even $\xi$ \cite[Section~7]{Tal07}, and the
resulting exponent of the upper deviations is $\Sigma$ \cite[Remark~1.3]{Bou26b}: under \condFM,
Chebyshev's inequality gives $\Prob(N^{-1}\log\ZN\ge f)\le e^{N(\Sigma(f)+o(1))}$ for
$f\in(\feq,\fone)$, and for $\xi\in\classC$ satisfying \condFM, in particular for convex
$\xi\in\classC$ (Theorem~\ref{thm:exist}\ref{thm:exist-e}), Theorems~\ref{thm:max}\ref{thm:max-c}
and~\ref{thm:exist}\ref{thm:exist-d} below show that $\Sigma(f)$ is the exact local rate of
$N^{-1}\log\ZN$ at $f$.

Following \cite[Definition~1.1]{Bou26b}, a level $f\in(\feq,\fone)$ is \emph{regular} if
$\theta\mapsto\Lambar(\theta)-\theta f$ has a unique minimizer $\thf$ on $[0,1]$. Every level of
$(\feq,\fone)$ is regular, with $\thf\in(0,1)$: under \condR\ this is
Proposition~\ref{prop:var-regular} (Theorem~\ref{thm:exist}\ref{thm:exist-b} needs only \condR),
and under \condH\ it is also Lemma~\ref{lem:loc-obstacle}. A \emph{constrained minimizer} at $\theta$ is a minimizer of
$\zeta\mapsto P_\zeta$ over $\{\zeta(\{0\})\ge\theta\}$. It exists, and is unique by
Lemma~\ref{lem:var-basic}. We write $\zf$ for the constrained minimizer at $\thf$; it satisfies
$\zf(\{0\})=\thf$ and $\Sigma(f)=\thf(P_{\zf}-f)$ (Lemma~\ref{lem:loc-obstacle}).

\paragraph{Generalized TAP critical points.}
\DisplayEventExist
This is the event of \cite[Section~1.3]{Bou26b}. It quantifies over all pairs $(\zeta,x)$ and need
not be measurable. For it, and for the events $\Emarg{g}{\varepsilon}{f}$ and
$\Efix{g}{\rho}{\varepsilon}{f}$ of Section~\ref{sec:intro-results}, $\Prob$ stands for the inner
probability $\Pin$ in lower bounds and for the outer probability $\Pout$ in upper bounds, and every
stated equality of limits holds for both (Sections~\ref{sec:max-existence} and~\ref{sec:stab-kr}). Criticality is imposed on $F_\zeta$ with an
order parameter that may depend on the point, and the last condition places the TAP value of $x$
at the same level. We call a pair $(\zeta,x)$ satisfying \eqref{eq:event-exist} a
\emph{witness}, and $x$ a \emph{generalized TAP critical point at level $f$}. Our second result
finds such points at every upper level.

\StateThmExist

\begin{remark}[Which critical points]\label{rem:exist-points}
The points of Theorem~\ref{thm:exist}\ref{thm:exist-a} are critical points of $F_{\zeta^*}$ for an
order parameter $\zeta^*$ that depends on the point. They are not shown to be isolated, to be
strictly stable, or to be critical points of the fixed function $F_{\zf}$. For convex $\xi$ with
at least two active degrees, Theorem~\ref{thm:stab}\ref{thm:stab-a} supplies critical points of $F_{\zf}$
that are locally strongly concave, hence isolated, at rates arbitrarily close to $\Sigma(f)$.
\end{remark}

By~\ref{thm:exist-a}, a maximizer of $\FTAP$ is itself a witness, so
\[
 \bigl\{|N^{-1}\MN-f|\le\varepsilon\bigr\}\subset\Eex{\varepsilon}{f}\subset
 \bigl\{\MN\ge N(f-\varepsilon)\bigr\}
\]
up to a null set, and Theorem~\ref{thm:max}\ref{thm:max-b} with \condFM\ identifies the rates of
the outer events. For convex $\xi\in\classC$, Theorems~\ref{thm:max}\ref{thm:max-c}
and~\ref{thm:exist}\ref{thm:exist-d}--\ref{thm:exist-e} give, at every $f\in(\feq,\fone)$,
\[
 \lim_{\varepsilon\downarrow0}\liminf_{N\to\infty}\frac1N
  \log\Prob\bigl(|N^{-1}\log\ZN-f|\le\varepsilon\bigr)
 =\lim_{\varepsilon\downarrow0}\limsup_{N\to\infty}\frac1N
  \log\Prob\bigl(|N^{-1}\log\ZN-f|\le\varepsilon\bigr)
 =\Sigma(f).
\]
For even $\xi$ this exponent follows, by the tilt of Section~\ref{sec:max-existence}, from
Theorem~\ref{thm:exist}\ref{thm:exist-b} and \cite[Section~7]{Tal07}, the latter extended by
truncation in the degree (Appendix~\ref{app:moments-guerra}); for convex $\xi$ with nonnegative coefficients that is
finite on $\R$, \condFM\ is stated in \cite[Theorem~2.1]{CGKLM}. From these fractional moments
Boursier announced a large deviation principle for $N^{-1}\log\ZN$ above $\feq$, with rate
function $-\Sigma$, for strictly convex $\xi$ \cite[Section~1.4]{Bou26a}. For SK with $\xi''(0)>1$, the order
$e^{-\Theta(Nx^{6/5})}$ of the tail at $\feq+x$ for small $x$ is
\cite[Theorem~1.5]{AL26}. What we add for the free
energy itself is the regularity of every level (Theorem~\ref{thm:exist}\ref{thm:exist-b}) and the
lower half of \condFM\ by the Kac--Rice route for convex $\xi\in\classC$ that is not even
(Proposition~\ref{prop:max-fm}\ref{prop:max-fm-d}).

Boursier \cite[Theorem~1]{Bou26b} proved the conclusion of~\ref{thm:exist-c} at
$\theta=\thf$ and the rate in~\ref{thm:exist-d} for $\Eex{\varepsilon}{f}$, under \condH, at
regular levels whose first positive contact is strictly Plefka stable \cite[(1.8)]{Bou26b}
(defined in Section~\ref{sec:intro-structure}); she noted that the non-strict Plefka inequality always
holds, and expected the strictness and the strict convexity to be removable
\cite[Remark~1.2]{Bou26b} (strict convexity is automatic for convex $\xi$ under \condR, so \condH\
coincides with her assumption). For convex mixtures in $\classC$, Theorem~\ref{thm:exist} removes the
strict Plefka condition and the regularity assumption from \cite[Theorem~1]{Bou26b}: every level
is regular by~\ref{thm:exist-b}, and marginal first contacts are allowed. For the single-degree
SK and pure even $p$-spin models, her theorem with~\ref{thm:exist-b} already applies at every
level whose first contact is strictly Plefka stable, and Theorem~\ref{thm:exist} adds the marginal
levels. Marginal first contacts occur for SK at low temperature at all levels close to $\feq$ (the
companion paper
\cite{PeiB}, not used here) and on an interval of levels in a family of near-critical $2+4$
mixtures (Theorem~\ref{thm:stab}\ref{thm:stab-c}). Boursier also notes that the constrained
minimizer is unique \cite[Section~1.3]{Bou26b} and proves that $\Lambar$ is convex with
$\Lambar'_+(0)=\feq$ and $\fone<\infty$ \cite[Lemma~3.2]{Bou26b}.

Existence is not elementary: through Theorem~\ref{thm:max}\ref{thm:max-b}, the lower bound
of~\ref{thm:exist-d}, asserted at every level of $(\feq,\fone)$, implies the lower half of \condFM\
on $(\thstar,1)$ (on $(0,\thstar]$ it holds under \condR\ alone,
Proposition~\ref{prop:max-fm}\ref{prop:max-fm-b}), and given the upper half the converse is~\ref{thm:exist-d} itself; so for convex $\xi$
existence is exactly as hard as the fractional-moment lower bound. Part~\ref{thm:exist-e}
assembles the inputs. For every convex $\xi$ the upper half of \condFM\ is Guerra's bound
\cite{Gue03} in the form \cite[Theorem~6.1]{Tal07}. For even $\xi$ the lower half is due to
Talagrand \cite[Section~7]{Tal07}, argued there by stated modifications of the proof of the
Parisi formula \cite{PF}, with no Kac--Rice argument; for convex $\xi$ that is not even we derive
it from the Kac--Rice construction behind Theorem~\ref{thm:stab}\ref{thm:stab-a}
(Proposition~\ref{prop:kr-main}) through Theorem~\ref{thm:max}\ref{thm:max-b}
(Proposition~\ref{prop:max-fm}); for pure even $\xi$,
Appendix~\ref{app:moments-pure} gives a second route that does not use \cite[Section~7]{Tal07}.
The formula is stated for every $\xi$ with nonnegative coefficients that is finite on $\R$ in
\cite[Theorem~2.1]{CGKLM}; there the upper bound is sketched for $\xi$ that is not convex on
$[-1,1]$, and the lower bound is presented as a modification of the proof of the Parisi formula.
For convex $\xi$ finite on $\R$ this covers the case that is not even; the Kac--Rice route of
Proposition~\ref{prop:max-fm}\ref{prop:max-fm-d} is an independent proof that also covers
$\xi\in\classC$ with finite radius of convergence. We do not rely on \cite{CGKLM}.

\subsection{The states that carry an upper deviation}\label{sec:intro-structure}

The structure of the states at level $f$ is read off the first variation of the Parisi functional
at $\zf$.

\paragraph{The obstacle, its contact set and the field law.} The Auffinger--Chen process of
$\zeta$ is the solution of
\[
 dX_t=\xi''(t)\,\zeta([0,t])\,\partial_x\Phi_\zeta(t,X_t)\,dt+\sqrt{\xi''(t)}\,dB_t,\qquad X_0=0,
\]
with $B$ a standard Brownian motion; the process $\partial_x\Phi_\zeta(t,X_t)$ is a martingale
\cite[Lemma~2]{AC15b}. Let $\Gz{\zeta}(t)=\E\,\partial_x\Phi_\zeta(t,X_t)^2$ and define the
\emph{obstacle} of $\zeta$,
\begin{equation}\label{eq:intro-obstacle}
 \Obs_\zeta(t)=\frac12\int_t^1\xi''(s)\bigl(\Gz{\zeta}(s)-s\bigr)\,ds,\qquad t\in[0,1].
\end{equation}
It is the first variation of the Parisi functional (\cite[Theorem~2]{Che17} for $\xi\in\classC$,
\cite[Lemma~3.2]{JT17} under \condR): for every
$\mu\in\Pspace([0,1])$,
\[
 \frac{d}{d\lambda}P_{(1-\lambda)\zeta+\lambda\mu}\Big|_{\lambda=0^+}
 =\int\Obs_\zeta\,d(\mu-\zeta),
\]
see Appendix~\ref{app:pde-toolbox}. In particular the derivative toward a point mass $\delta_t$ is
$\Obs_\zeta(t)-\int\Obs_\zeta\,d\zeta$: the minimizers of the obstacle are the cheapest places
to move mass to. For the
constrained minimizer write $\Obs=\Obs_{\zf}$ and $\Omin=\min_{(0,1]}\Obs$, and define the
\emph{contact set} and the \emph{field law}
\[
 \Cf=\argmin_{(0,1]}\Obs,\qquad
 \Lq{q}=\Law\bigl(\partial_x\Phi_{\zf}(q,X_q)\bigr),\quad q\in[0,1],
\]
where $X$ is the Auffinger--Chen process of $\zf$. The contact set is a nonempty compact subset of
$(0,1)$ that contains $\supp\zf\cap(0,1]$, and $\Obs(0)>\Omin$ (Lemma~\ref{lem:loc-obstacle}).
Its points are the positive contact points of \cite[Section~1.3]{Bou26b}, and $\qf=\min\Cf$ is
her first positive contact point. The second moment of $\Lq{q}$ is $\Gz{\zf}(q)$, which equals $q$
at every contact because $\Obs'=-\frac12\xi''(\Gz{\zf}-t)$ vanishes there.

\paragraph{The Plefka factor and the Hessian bulk.} For $Q\in(0,1)$ we write $\aQ=\xi''(Q)$,
$\UQ=\partial_x^2\Phi_{\zf}(Q,X_Q)$, which takes values in $(0,1]$, $\chiQ=\E\UQ$ and
\[
 \Plef(Q)=1-\aQ\,\E\UQ^2 .
\]
For $\xi\in\classC$, at a contact $Q\in\Cf$ the obstacle is twice differentiable with
$\Obs''(Q)=\frac12\aQ\Plef(Q)\ge0$, since $Q$ is an interior minimum
(Lemma~\ref{lem:pde-standing}). So the \emph{Plefka inequality} $\Plef\ge0$ holds on $\Cf$. A
contact $Q$ is \emph{strictly Plefka stable} if $\Plef(Q)>0$ (at $Q=\qf$ this is the strict
Plefka condition \cite[(1.8)]{Bou26b}, named after Plefka's condition for the classical TAP
equations \cite{Ple82}), and \emph{marginal} if $\Plef(Q)=0$, that is, if the obstacle is flat to
second order at $Q$.

For $a\ge0$ and a probability measure $\rho$ on $[0,\infty)$, the \emph{reflected Pastur law}
$\Pastur{a}{\rho}$ is the unique probability measure on $\R$ whose Cauchy transform
$g(z)=\int(z-\lambda)^{-1}\,\Pastur{a}{\rho}(d\lambda)$ satisfies
\begin{equation}\label{eq:intro-pastur}
 g(z)=\int\frac{\rho(dD)}{z+D-a\,g(z)}\qquad\text{for every } z\in\Cplx^+ .
\end{equation}
It is the limit of the empirical spectral distribution of $\sqrt a\,\GOE-D_N$, where $\GOE$ is a
GOE matrix (off-diagonal variance $1/N$) and $D_N$ are diagonal matrices whose empirical
spectral distributions converge to $\rho$ (Lemma~\ref{lem:hess-pastur}). Put
\[
 \rhoQ=\Law\bigl(\UQ^{-1}+\aQ\chiQ\bigr),\qquad \nuQ=\Pastur{\aQ}{\rhoQ},\qquad
 \lamp(Q)=\sup\supp\nuQ .
\]
For $\alpha\in\Pspace([0,1])$ and $x\in\ocube$ with $q_x\in(0,1)$, the function $F_\alpha$ has two
one-sided Hessians $\Hpm F_\alpha(x)$ at $x$, taken from either side of the sphere
$\{q=q_x\}$; they coincide when $\alpha(\{q_x\})=0$ (Lemma~\ref{lem:matrix-hessian}). Eigenvalues
of a symmetric $N\times N$ matrix $A$ are ordered $\lambda_1(A)\ge\dots\ge\lambda_N(A)$, with
$\lambda_{\max}=\lambda_1$ and $\ESD(A)=N^{-1}\sum_i\delta_{\lambda_i(A)}$; $A\preceq B$ means that
$B-A$ is positive semidefinite; $\dL$ is the L\'evy distance on $\Pspace(\R)$ and $\Wone$ the
Wasserstein distance on $\Pspace([-1,1])$. Events that quantify over all points of the cube, as
in Theorems~\ref{thm:struct} and~\ref{thm:stab}\ref{thm:stab-b} and Corollaries~\ref{cor:title}
and~\ref{cor:eq}, need not be measurable; we bound them in outer probability $\Pout$. Our third
result describes every state that carries an upper deviation.

\StateThmStruct

Theorem~\ref{thm:struct} concerns every point at which $\FTAP$ is at least $N(f-\delta)$, not only critical points, and its proof uses no Kac--Rice argument.
Part~\ref{thm:struct-i} is an Ising analogue, at rare levels, of Subag's characterization of the
TAP representation for spherical models at equilibrium \cite[Theorems~5 and~10]{Sub24}.
Part~\ref{thm:struct-ii} proves, for every near-maximizer at level $f$ and uniformly over the cube,
the empirical-law description that Boursier \cite[Section~1.5]{Bou26a} gives heuristically for
typical TAP states and that \cite[(1.29)]{CPS2} recalls informally for ancestor states at
equilibrium, following \cite[Chapter~4]{Pan13} and \cite{AJ19a,AJ19b}. Alone, \ref{thm:struct-ii}
localizes $q_m$ only near stationary points of the obstacle, since
$|\Gz{\zf}(q_m)-q_m|\le2\Wone(\mu_m,\Lq{q_m})$ ($m\mapsto m^2$ is $2$-Lipschitz on $[-1,1]$); \ref{thm:struct-i} also excludes stationary
points that are not minimizers. We do not address the ultrametric organization of the TAP states at
level $f$ conjectured in \cite[Conjecture~3]{Bou26a}.

Parts~\ref{thm:struct-iii}--\ref{thm:struct-iv} extend to generalized TAP functionals at rare
levels the bulk-edge dichotomy proved for SK with external field at Bolthausen's iterates
\cite{Bol14}, strictly below zero where the de~Almeida--Thouless condition \cite{AT78} is strict
and reaching zero on its boundary \cite[Theorem~1.4]{GSS}, and for spherical approximate ground
states \cite[Proposition~1.2 and the remark after it]{Sel25}; its physical origin is the work of
Bray and Moore \cite{BM79,BM80} (see \cite[p.~3]{GSS}) and of Plefka \cite{Ple82} (see also
\cite{ABM04,AM19}). The quadratic
vanishing $\lamp=-\Plef^2/(4\aQ\E\UQ^3)+O(\Plef^3)$ of Lemma~\ref{lem:hess-edge} mirrors the
explicit square in \cite[(1.7)]{Sel25}. The strict side, at a strictly Plefka-stable first
contact and conditionally on criticality, the level and a convergent field law, is
\cite[Proposition~B.1, Remark~B.1]{Bou26b}, and Remark~B.1 there notes that the Plefka condition
places the limiting bulk in $(-\infty,0]$; related observations, for SK under full replica
symmetry breaking and for the classical TAP Hessian at the Plefka condition, are
\cite[Section~2.5]{JSS} and \cite[Proposition~5.1]{BCG}. What~\ref{thm:struct-iii}--\ref{thm:struct-iv} add is generality: generalized TAP functionals,
both $F_{\zf}$ and every near-optimal order parameter of the point itself; every self-overlap
$Q\in(0,1)$ through $\nuQ$; every point at a rare level, uniformly over the cube; every convex
mixture in $\classC$; and the coupling with~\ref{thm:struct-i}--\ref{thm:struct-ii}. We do not
control outlier eigenvalues, a limitation also recorded in \cite[after Theorem~1.4]{GSS} and
\cite[Remark~B.1]{Bou26b}; an index form of~\ref{thm:struct-iii} is
Proposition~\ref{prop:hess-index}. The uniqueness of $\zf$ needs only \condR\
(Lemma~\ref{lem:var-basic}); parts~\ref{thm:struct-i}--\ref{thm:struct-ii} do not use it.

Combined with Theorems~\ref{thm:max} and~\ref{thm:exist}, Theorem~\ref{thm:struct} gives the
statement of the title: conditionally on an upper deviation of the free energy, the deviation is
carried by TAP states of the structure just described.

\StateCorTitle

In~\ref{cor:title-a}, $\delta$ and $c$ depend only on $\kappa$ (and on $f$ and $\xi$); the
closeness $t$ of $\MN$ to $\log\ZN$ is arbitrary and enters only through $N_0$. The conditioning
event has probability $e^{N(\Sigma(f)+o(1))}$, and the exceptional events of
Theorem~\ref{thm:struct} are smaller by a factor $e^{-cN}$. Corollary~\ref{cor:title} does not by
itself assert marginal stability: with Theorem~\ref{thm:struct}\ref{thm:struct-iv}, every listed
Hessian has its bulk within L\'evy distance $\kappa$ of $\nu_{q_m}$, whose edge $\lamp(q_m)$ is
nonpositive, and at a contact this edge vanishes exactly when
$\Obs''=0$ there. Part~\ref{cor:title-b} is the same statement conditioned on the TAP maximum.

\subsection{Stability and the equilibrium level}\label{sec:intro-results}

Theorem~\ref{thm:stab} shows that marginality has a cost when it is forced. The relevant levels
are those at which every contact is marginal:
\DisplayCondHmarg
When $\Cf$ is an interval of positive length, \condHmarg\ is automatic, since $\Gz{\zf}=\id$ on
$\Cf$ gives $\Obs''=\frac12\xi''(1-\Gz{\zf}')=0$ there; whether it holds for SK at the marginal levels of
Section~\ref{sec:intro-setting} is not addressed here (Section~\ref{sec:open}). For the stability statements we use two strengthenings of the existence event.
\DisplayEventsStab
In $\Efix{g}{\rho}{\varepsilon}{f}$ the order parameter is the fixed $\zf$, and the TAP function is
uniformly concave, with margin $g$, on a ball of radius $\rho\sqrt N$ around the critical point.

\StateThmStab

Part~\ref{thm:stab-a} constructs strictly stable, isolated critical points of the fixed function
$F_{\zf}$ (Appendix~\ref{app:kr}), with a fixed margin at a strictly Plefka-stable first contact;
at a marginal first contact the construction first tilts the field law by a small parameter
$\eta$, which opens the Plefka inequality, and its margins are then at most of order $\eta^2$
(Remark~\ref{rem:stab-margin}). Part~\ref{thm:stab-b} shows that this loss is unavoidable when
every contact is marginal: with no criticality assumed, points at the level whose one-sided
Hessian, for $F_{\zf}$ or for a near-optimal order parameter of their own, is $\preceq-gI$ are
exponentially rarer than the level itself, with a rate loss $c=c(g)$ that is not explicit.
Part~\ref{thm:stab-c} shows that~\ref{thm:stab-b} is not vacuous; the family has two active
degrees, so~\ref{thm:stab-a} applies to it. Our last result is the equilibrium counterpart of
Theorem~\ref{thm:struct}\ref{thm:struct-i}--\ref{thm:struct-ii}.

\StateCorEq

Corollary~\ref{cor:eq} is the Ising counterpart of Subag's characterization for spherical models
\cite[Theorems~5 and~10]{Sub24}: off the Parisi contact set the TAP representation fails. There,
for generic spherical mixtures, the representation at radius $q$ holds exactly when $q$ belongs to
a set introduced by Talagrand, which contains the support of the Parisi measure and can be
strictly larger \cite[Proposition~9 and Remark~13]{Sub24}. Here the failure holds for every convex
mixture, with an exponential bound and a field-law clause. For Ising spins and $\xi\in\classC$, \cite[Theorem~2]{CPS}
and \cite[Theorems~2, 3 and~5]{CPS2} concern radii $q\in\supp\zP$ (the survey
\cite[Section~4.4]{Sub26} summarizes these as analogues of the spherical results; no Ising statement
off the support is among them), and at every radius the TAP maximum is at most the free energy
\cite[Section~5]{CPS}. For $\xi\in\classC$, at a fixed radius $Q\notin\CP$ with
$\Gz{\zP}(Q)\ne Q$, and at every $Q\notin\CP$ below $\max\supp\zP$, the fixed-radius clause of
Corollary~\ref{cor:eq} follows in a few lines from the exact formula of \cite[Theorem~36]{CPS},
taken from \cite{CP18} (Remark~\ref{rem:eq-cps}). What Corollary~\ref{cor:eq} adds is the convex
mixtures outside $\classC$, uniformity over the cube, the field law, and the radii $Q\notin\CP$
above $\max\supp\zP$ with $\Gz{\zP}(Q)=Q$, if there are any. For spherical models every point of
Talagrand's set, the spherical analogue of the contact set, is multi-samplable
\cite[Theorem~10]{Sub24} and hence carries the representation \cite[Theorem~5]{Sub24}, also at its
points outside the support. For Ising spins, Boursier indicates that when the equilibrium contact
set is finite, TAP states at level close to $\feq$ can be built at each contact point by an
induction over bands, under a strict stability assumption \cite[abstract and Section~1.3]{Bou26b}; we know of
no written proof, and the general case is open (Section~\ref{sec:open}).

\input{figures/fig1-roadmap}

\subsection{How the proofs work}\label{sec:intro-ideas}

Figure~\ref{fig:roadmap} shows the proofs as three columns, for Theorems~\ref{thm:max}
and~\ref{thm:exist}, for Theorem~\ref{thm:struct} and for Theorem~\ref{thm:stab}, with the two
corollaries and Lemma~\ref{lem:eq-zero} above them and, in two bands below, the techniques that
several columns share: the Parisi equation (Appendix~\ref{app:pde}) and the Kac--Rice
construction (Appendix~\ref{app:kr}). We read the figure one column at a time.

\paragraph{Theorems~\ref*{thm:max} and~\ref*{thm:exist}: Ramsey's theorem and the band theorem.} The
band theorem \cite[Theorem~1]{CPS} identifies $N^{-1}\FTAP(m)$, with probability exponentially
close to one and up to the constant $\log2$ and an error small uniformly in $m$, with the free
energy per replica of $n$ replicas confined to a thin band around $m$, for large $n$. Under
\condR\ we obtain it at every exponential rate by truncating $\xi$ inside the band
(Lemma~\ref{lem:max-truncation}); with \cite[(2.16)]{CPS} this gives the upper bound
(Proposition~\ref{prop:max-upper}). For the lower bound, color each pair of $r_J(2n)$ independent
Gibbs samples, $r_J(2n)$ the Ramsey number for $2n$-cliques in $J$ colors, by which of $J$
intervals of $[-1,1]$ contains their overlap. Some $2n$ samples form a monochromatic clique, and
taking expectations shows that, for some interval, $2n$ independent samples have all mutual
overlaps in it with probability at least $c_{J,2n}>0$, independently of $N$ and of the disorder
(Lemma~\ref{lem:max-ramsey}). By averaging, for some realization of $n$ of them, the
band around their barycenter $m$ carries Gibbs mass at least $c_{J,2n}$ for the remaining $n$
(Lemma~\ref{lem:max-band}); so the band free energy at $m$ falls short of $\log\ZN$
by at most $n^{-1}\log(1/c_{J,2n})$, and $\MN\ge\FTAP(m)\ge\log\ZN-tN$ outside an event of
probability $e^{-cN}$, for every $c$ (Proposition~\ref{prop:max-lower}, Figure~\ref{fig:ramsey}).
Parts~\ref{thm:max-b}--\ref{thm:max-c} follow by H\"older's inequality with Gaussian moment
bounds, and by exponential equivalence \cite[Section~4.2]{DZ}. In
Theorem~\ref{thm:exist}\ref{thm:exist-a}, maximizers are interior since
$\partial_m\Phistar_\zeta(q,m)\to\mp\infty$ as $m\to\pm1$, uniformly in $(q,\zeta)$, while the
gradient of $\HN$ is bounded on the cube, and they are critical points of $F_{\zeta^*}$ by the
uniqueness and continuity of the minimizer in~\eqref{eq:intro-ftap}
\cite[Theorem~12(iii) and Section~11]{CPS} and the
envelope lemma \cite[Lemma~37]{CPS}
(Lemmas~\ref{lem:gauss-interior} and~\ref{lem:gauss-critical});
parts~\ref{thm:exist-c}--\ref{thm:exist-d} follow from \condFM, Theorem~\ref{thm:max} and an
exponential tilt by $e^{\thf\MN}$ (Section~\ref{sec:max-existence}).

\paragraph{Theorem~\ref*{thm:struct}\ref*{thm:struct-i}--\ref*{thm:struct-ii}: the glued comparison.}
By Markov's inequality and $\Sigma(f)=\thf(P_{\zf}-f)$, excluding a sphere $\{q_m=Q\}$ needs a
bound on $\E\exp\{\thf\sup_{q_m=Q}\FTAP(m)\}$ exponentially smaller than $e^{N\thf P_{\zf}}$.
Boursier's comparison \cite[Proposition~2.7]{Bou26b} bounds the exponential moments of the
maximum over the whole cube by Parisi functionals $P_\zeta$ with $\zeta(\{0\})\ge\thf$, hence by
$e^{N\thf P_{\zf}}$ at best, and on a sphere of positive radius a nondecreasing coefficient comes from an
admissible order parameter and gains nothing (Lemma~\ref{lem:loc-monotone}). We therefore run the
comparison on the sphere with a \emph{glued} coefficient $\glue(\xc,\eta)$, equal to a cascade
coefficient $\xc$ on $[0,Q)$ and to the distribution function of a band measure
$\eta\in\Pspace([Q,1])$ on $[Q,1]$, with no relation between the two at $Q$
(Proposition~\ref{prop:glued-comparison}):
\begin{equation}\label{eq:intro-glued}
 \E\exp\Bigl\{\thf\Bigl(\sup_{q_m=Q}\FTAP(m)-N\,P\bigl(\glue(\xc,\eta)\bigr)\Bigr)\Bigr\}\le1 .
\end{equation}
This is Boursier's comparison \cite[Proposition~2.7]{Bou26b} run on a sphere or thin shell. Her
proof uses neither monotonicity of the coefficient across the radius nor the bound by one
(Appendix~\ref{app:glued-map} maps each step), so the coefficient may exceed one below the radius
and may jump down at it. It rests on her shift identity for the Onsager field on a Ruelle cascade
(the first identity of \cite[Remark~2.2]{Bou26b}), which we prove for all positive cascade
parameters (Lemma~\ref{lem:onsager-shift}), and on an outer Poisson level of parameter $\thf$, the
device of \cite{CGKLM}, as in \cite[p.~13]{Bou26b}. Parameters above one appear in Talagrand's
reverse bound for $\E\ZN^a$, $a>1$, whose parameters lie between $1$ and $a$
\cite[Proposition~6.2]{Tal07}, and non-monotone effective coefficients in \cite[Section~6]{HS25}
and \cite{EAMS21}. If $Q\notin\Cf$, pick $r\in\Cf$; when $r<Q$, raising the coefficient by
$\varepsilon$ on $[r,Q)$ changes the Parisi functional by
$\varepsilon(\Obs(r)-\Obs(Q))+O(\varepsilon^2)=-\varepsilon(\Obs(Q)-\Omin)+O(\varepsilon^2)$; a
suitable $\varepsilon$ gives data with $P(\glue(\xc,\eta))\le P_{\zf}-\varphi(\Obs(Q)-\Omin)$, where
$\varphi(d)=\min\{d^2/(4K),\thf d/2\}$ and $K$ depends only on $\xi'(1)$ and $\xi''(1)$
(Proposition~\ref{prop:loc-gap}). When $r>Q$, one lowers the band distribution function on $[Q,r)$
instead. For $Q>0$, either way the gain comes from a downward jump of the glued coefficient at $Q$. A strict
Parisi obstacle gap is already used inside a band comparison in \cite[Proposition~3.9 and
Lemma~3.10, (3.15)--(3.16)]{Bou26b}, for competitors at a fixed relative overlap with a critical
point; we use it for the self-overlaps of all near-maximizers, through a jump at the radius.
Markov's inequality and~\eqref{eq:intro-glued} give
\[
 \Pout\Bigl(\sup_{q_m=Q}\FTAP(m)\ge N(f-\delta)\Bigr)\le e^{N(\Sigma(f)-\thf\varphi+\thf\delta)},
\]
and~\ref{thm:struct-i} follows with thin shells and a union bound (Section~\ref{sec:loc},
Figure~\ref{fig:obstacle}). For~\ref{thm:struct-ii} we add $s\sum_ig(m_i)$ to $\FTAP$, for a
smooth test function $g$; the perturbed Parisi functional is at most
$P_{\zf}+s\,\E g(\Lq{Q})+Cs^2$, so an excess $N^{-1}\sum_ig(m_i)\ge\E g(\Lq{q_m})+\kappa/2$ is
exponentially rare once $s$ is small, and a finite net of test functions with
Kantorovich--Rubinstein duality gives the $\Wone$ bound (Section~\ref{sec:loc-field},
Appendix~\ref{app:var-field}). Perturbing by a parameter and using convexity in it is the device
by which \cite[Theorem~5]{CPS2} computes the energy of TAP states.

\paragraph{Theorem~\ref*{thm:struct}\ref*{thm:struct-iii}--\ref*{thm:struct-iv}: the edge identity.}
At $x\in\ocube$ with $q_x=Q$ one has $\Hpm F_\alpha(x)=\Hess\HN(x)-D_x+E_x$, where $D_x$ is diagonal with
entries $u_i^{-1}+\xi''(Q)N^{-1}\sum_ju_j$, $u_i=\partial_x^2\Phi_\alpha(Q,b_i)$ and
$\partial_x\Phi_\alpha(Q,b_i)=x_i$, and $E_x$ has rank at most two
(Lemma~\ref{lem:matrix-hessian}). With probability at least $1-e^{-MN}$, for any prescribed $M$,
$\ESD(\Hess\HN(x)-D)$ is close to $\Pastur{\xi''(q_x)}{\ESD(D)}$ simultaneously for all $x$ and
all bounded nonnegative diagonal $D$ (Appendix~\ref{app:matrix-uniform}; uniform free-convolution estimates
of this kind appear in \cite[Section~4]{JSS}), so at a localized point the bulk is close to $\nuQ$
(Proposition~\ref{prop:hess-bulk}). An order parameter $\alpha$ of the point itself with
$F_\alpha(x)\le\FTAP(x)+2\varepsilon N$ is close to $\zf$ on $[q_x,1]$, by the uniqueness in
\cite[Theorem~12(iii)]{CPS} (Proposition~\ref{prop:hess-own}), a quantitative form, for rare
states, of the observation \cite[(2.34)]{CPS} that the order parameter of a TAP state follows the
Parisi measure; Proposition~\ref{prop:hess-own} uses
Theorem~\ref{thm:struct}\ref{thm:struct-i}--\ref{thm:struct-ii} at a smaller accuracy. With $U=\UQ$, $a=\aQ$, $\chi=\E U$ and $D=U^{-1}+a\chi$, the function
$z(s)=s+a\,\E(s+D)^{-1}$, $s>-\essinf D$, is strictly convex, its derivative tends to $-\infty$ at
the left end (by a Girsanov bound on the law of $X_Q$), and $\lamp(Q)=\min z$
(Lemma~\ref{lem:hess-edge}). Since
\[
 z(-a\chi)=0,\qquad z'(-a\chi)=1-a\,\E U^2=\Plef(Q),\qquad z''(-a\chi)=2a\,\E U^3>0 ,
\]
$\lamp(Q)\le0$, with equality if and only if $-a\chi$ minimizes $z$, that is, $\Plef(Q)=0$
(Figure~\ref{fig:edge}). This is Theorem~\ref{thm:struct}\ref{thm:struct-iv}, and at contacts
$\Plef=2\Obs''/\xi''$.

\paragraph{Theorem~\ref*{thm:stab}: the stability-rate trade-off.} Part~\ref{thm:stab-b} follows
from the second column: under \condHmarg\ the edge $\lamp$ vanishes on $\Cf$, so $\nuQ$ gives mass
at least $p>0$ to $(-g/4,\infty)$ for all $Q$ near $\Cf$, and by
Theorem~\ref{thm:struct}\ref{thm:struct-iii}, outside the exceptional event, every Hessian in
question has at least $pN/4$ eigenvalues above $-g/2$, so it is not $\preceq-gI$
(Section~\ref{sec:stab-margin}). Part~\ref{thm:stab-a} counts critical points of the fixed function
$F_{\zf}$ with self-overlap near $\qf$ and prescribed empirical field law (Appendix~\ref{app:kr}).
At a marginal first contact the prescribed law is a tilt of $\Lq{\qf}$ of size $\eta$, which
preserves the self-overlap and makes the Plefka factor positive; with two active degrees, a small
increase of the prescribed energy also opens the radial Hessian gap and a gap for competitors at
smaller overlap. The isolation and spherical-code argument of \cite[Section~3]{Bou26b}, which
adapts that of Huang and Sellke \cite{HS}, turns the expected count into an existence probability
at rate $\Sigma(f)-O(\eta)$ (Proposition~\ref{prop:kr-main}). A margin $g$ needs $g\le-\lamp$ for
the tilted law, so the quadratic vanishing of the edge forces $g=O(\eta^2)$
(Remark~\ref{rem:stab-margin}). For~\ref{thm:stab-c}, every point of an interval
$[q_\theta,v_\theta]\subset\supp\zf$ is a contact, and a third-derivative bound with Rolle's
theorem on the plateaus of $\zf$ and a confinement argument, adapting the near-critical
gap-exclusion mechanism of Zhou \cite{Zho26a}, show that
$h=\Gz{\zf}-\id$ is positive on $(0,q_\theta)$ and negative on $(v_\theta,1]$; so
$\Cf=[q_\theta,v_\theta]$ and $\Obs''=\frac12\xi''(1-\Gz{\zf}')=0$ on it
(Lemma~\ref{lem:stab-family}, Appendix~\ref{app:nc}). Part~\ref{thm:stab-a} at a fixed margin
and item~\ref{thm:stab-c3} of~\ref{thm:stab-c} also use Theorem~\ref{thm:exist}\ref{thm:exist-d}--\ref{thm:exist-e}.

\paragraph{The corollaries.} Corollary~\ref{cor:title} composes the first two columns: the
conditioning event costs $e^{N(\Sigma(f)+o(1))}$ by Theorems~\ref{thm:max} and~\ref{thm:exist}
and the continuity of $\Sigma$, and the exceptional events of Theorem~\ref{thm:struct} cost
$e^{N(\Sigma(f)-c)}$ (Section~\ref{sec:hess-cor}). Corollary~\ref{cor:eq} runs the second column
at the equilibrium level, with $\zP$ in place of $\zf$ and a small exponential tilt in place of
$\thf$. Its one additional input is that $0$ is a contact point of the equilibrium obstacle for
every Parisi measure (Lemma~\ref{lem:eq-zero}; for $\xi\in\classC$ it also follows from
$0\in\supp\zP$ \cite[Theorem~1(i)]{AC15a}), which we prove under \condH\ from the convexity of
$\Gz{\zP}$ on the initial gap of $\zP$; it lets the coefficient be raised on $[0,Q)$ at every radius
$Q\notin\CP$ (Section~\ref{sec:eq}).

\subsection{Earlier work and organization}\label{sec:intro-prior}

\paragraph{Boursier's program.} In \cite[Section~1.4]{Bou26a} Boursier observed that the
fractional-moment formula gives, ``with a little more work'', a large deviation principle for the
free energy, whose rate function is, up to sign, the annealed complexity of TAP states that she
conjectured in \cite[Conjecture~1]{Bou26a}; in \cite[Section~1.2]{Bou26b} she revised this,
identifying the Legendre transform with the annealed count of isolated local maxima (the
supersymmetric count) rather than with the ordinary annealed complexity. She read the duality as a
rigorous counterpart of
Monasson's real-replica method \cite{Mon95}, announced it for a forthcoming paper, and suggested
proving Conjecture~1 there through a ``purely probabilistic argument'' that bypasses Kac--Rice
\cite[p.~4]{Bou26a}. In \cite{Bou26b} she proved the existence theorem discussed after
Theorem~\ref{thm:exist}, gave the heuristic of \cite[Remark~1.3]{Bou26b}, and proved the
free-energy lower bound at strictly stable levels \cite[Remark~1.4 and Proposition~3.12]{Bou26b}.
Neither Theorem~\ref{thm:max} nor Theorem~\ref{thm:exist} addresses the ordinary annealed or
quenched complexity of \cite[Conjectures~1 and~2]{Bou26a}.
In physics, upper deviations were studied through the analytic continuation
of the replica number \cite{Kon83,PR}.

\paragraph{TAP free energies.} The TAP functional goes back to \cite{TAP77}; its Hessian, the
complexity of its critical points and its stability were studied in physics
\cite{BM79,BM80,Ple82,ABM04,AM19}, and rigorous work includes \cite{Bol14,CP18,AJ19a,BK19,Gay24}.
The generalized TAP free energy was developed by Subag for spherical models \cite{Sub24} and by
Chen, Panchenko and Subag for Ising spins \cite{CPS,CPS2}, and used by Huang and Sellke \cite{HS};
spin distributions of Gibbs states \cite{AJ19b} are a different object from the empirical laws of
Theorem~\ref{thm:struct}\ref{thm:struct-ii}, and the parametrization and first-variation formula
of the Parisi functional that we use are those of \cite{AC15a,AC15b,Che17,JT16,JT17}.

\paragraph{Partly anticipated statements.} The endpoint identity $\Sigma'(\feq+)=-\thstar$
(Appendix~\ref{app:var-constrained}) is stated in \cite[Remark~1.4]{Bou26a}. Given existence, the
retained isolated states of Appendix~\ref{app:kr} are subexponentially many
(Proposition~\ref{prop:kr-retained}), by the spherical-code bound of \cite[Lemma~3.11 and
Section~3.7]{Bou26b}, used there under the strict Plefka condition. The expected number of
critical points in the supersymmetric window has rate $\Sigma(f)$ at every level, marginal ones
included, for convex $\xi\in\classC$ with two active degrees (Proposition~\ref{prop:stab-susy});
Boursier identifies this exponent under the strict Plefka condition \cite[Proposition~3.6]{Bou26b},
and \cite[Proposition~2.2]{Bou26a} counts a different window, with no Plefka hypothesis stated.

\paragraph{Existing results that we re-prove or extend.} They include the adapted inputs of
Figure~\ref{fig:roadmap}; the shift identity and the parameter-convexity device sit inside
Proposition~\ref{prop:glued-comparison} and Appendix~\ref{app:var-field}, and the facts on the
Parisi equation and the isolation argument are the bands of the figure.
\begin{itemlist}[itemsep=0pt]
\item The band theorem \cite[Theorem~1]{CPS}, which holds at every exponential rate for
 $\xi\in\classC$, extended to \condR\ by
 truncation of $\xi$ inside the band (Lemma~\ref{lem:max-truncation},
 Appendix~\ref{app:gauss-truncation}).
\item The band step of \cite[Section~5]{CPS}, at a barycenter of $n$ replicas, which lies in the
 cube at every overlap level (Lemma~\ref{lem:max-band}).
\item The gradient computation of \cite[Section~11]{CPS} and the envelope lemma
 \cite[Lemma~37]{CPS}, at the maximizers of $\FTAP$ (Lemma~\ref{lem:gauss-critical}).
\item Boursier's comparison \cite[Proposition~2.7]{Bou26b}, on spheres and thin shells, with glued
 coefficients whose cascade part may exceed one and jump down at the radius, and with perturbed
 terminal conditions (Proposition~\ref{prop:glued-comparison}, Appendix~\ref{app:glued}).
\item The shift identity for the Onsager field (the first identity of \cite[Remark~2.2]{Bou26b}),
 stated there for parameters in $(0,1]$, for all positive parameters (Lemma~\ref{lem:onsager-shift}).
\item The parameter-convexity device of \cite[Theorem~5]{CPS2}, applied with a parameter
 multiplying $\sum_ig(m_i)$, for the field law (Appendix~\ref{app:var-field}).
\item The facts on the Parisi equation that we need, for nonnegative, possibly non-monotone
 coefficients (Appendix~\ref{app:pde}); from \cite{AC15a} we take the definition of $\Phi_\zeta$
 for general $\zeta$, and cite its derivative bounds and representation formulas in
 Appendices~\ref{app:pde-standing} and~\ref{app:nc}.
\item The bulk-edge dichotomy of \cite[Theorem~1.4]{GSS} and \cite[Proposition~1.2]{Sel25}, with
 its strict side \cite[Proposition~B.1 and Remark~B.1]{Bou26b}, for the reflected Pastur laws $\nuQ$
 (Lemma~\ref{lem:hess-edge}).
\item Spectral limits at TAP points (\cite[Sections~5--6]{GSS} at Bolthausen's iterates,
 \cite[Appendix~B]{Bou26b} conditionally on criticality), uniform over the cube and over
 near-optimal order parameters, with failure probability $e^{-MN}$ for every $M$
 (Proposition~\ref{prop:hess-bulk}).
\item The isolation and spherical-code argument of \cite[Section~3]{Bou26b}, which adapts
 \cite{HS}, at marginal first contacts after the field-law tilt and the energy shift
 (Proposition~\ref{prop:kr-main}).
\item The near-critical gap-exclusion mechanism of \cite{Zho26a}, for the constrained problem in
 $2+4$ mixtures (Appendix~\ref{app:nc}); Zhou's theorem on SK at equilibrium is not used as a
 statement.
\end{itemlist}

\paragraph{Organization.} Sections~\ref{sec:max}--\ref{sec:eq} follow the columns of
Figure~\ref{fig:roadmap}, each titled by the results it proves; Section~\ref{sec:open} collects
open problems. Appendices~\ref{app:pde}--\ref{app:moments} treat, in order, the Parisi functional
with non-monotone coefficients and the constrained problem; Ruelle cascades, the glued comparison
and the field law; Gaussian bounds and the TAP maximum; random-matrix inputs; the Kac--Rice
construction; the near-critical $2+4$ family; and fractional moments. Theorem, page and equation locators for
\cite{CPS,CPS2,GSS,Sel25,JSS,AC15a,AC15b,JT16,JT17,Che17,Sub24,Sub26,HS25,Zho26a,Bou26a,Bou26b}
refer to the
arXiv versions listed in the bibliography; those for \cite{Tal07} follow the author's manuscript.

%% file: figures/fig1-roadmap.tex
\begin{figure}[tp]
\centering
\newcommand*{\rmloc}[1]{\\[1pt]{\tiny\color{black!75}#1}}%
\begin{NoHyper}%
\begin{tikzpicture}[
  x=1cm, y=1cm, font=\scriptsize,
  execute at begin node={\thickmuskip=5mu\medmuskip=4mu\relax},
  >={Stealth[length=1.5mm,width=1.1mm]},
  thmlook/.style={fill=black!80, draw=black, line width=0.8pt},
  ourslook/.style={rounded corners=1.3mm, fill=black!16, draw=black, line width=0.6pt},
  knownlook/.style={fill=white, draw=black, line width=0.4pt},
  redolook/.style={fill=white, draw=black, line width=0.45pt,
    path picture={\draw[black, line width=0.45pt]
      ([shift={(1.6pt,1.6pt)}]path picture bounding box.south west) rectangle
      ([shift={(-1.6pt,-1.6pt)}]path picture bounding box.north east);}},
  toollook/.style={rounded rectangle, fill=black!9, draw=black, line width=0.6pt},
  box/.style={rectangle, align=center, text width=2.115cm, inner xsep=2pt,
    inner ysep=1.5pt, minimum height=0.9cm},
  thm/.style={rectangle, align=center, inner xsep=1.2pt, inner ysep=2pt,
    minimum height=1.5cm, thmlook, text=white, font=\scriptsize\bfseries},
  ours/.style={box, ourslook},
  known/.style={box, knownlook},
  redo/.style={box, redolook, text width=2.073cm, inner xsep=2.6pt, inner ysep=2.8pt},
  tool/.style={toollook, align=center, inner sep=2.5pt, minimum height=0.5cm},
  sw/.style={rectangle, minimum width=0.5cm, minimum height=0.3cm, inner sep=0pt},
  lg/.style={anchor=west, inner sep=2pt, text height=1.6ex, text depth=0.3ex},
  use/.style={->, line width=0.5pt, draw=black!80},
  rel/.style={use, dashed},
  band/.style={->, line width=0.9pt, draw=black!45},
  col/.style={rounded corners=2.5mm, fill=black!3, draw=black!30, line width=0.4pt},
]
\draw[col] (0.00,1.07) rectangle (5.20,-7.75);
\draw[col] (5.55,1.07) rectangle (10.75,-7.75);
\draw[col] (11.10,1.07) rectangle (16.30,-7.75);
\node[thm, text width=2.235cm] (thmOne) at (1.19,0)
  {Theorem~\ref{thm:max}\\ $\MN\approx\log\ZN$\\ at every rate};
\node[thm, text width=2.235cm] (thmTwo) at (4.01,0)
  {Theorem~\ref{thm:exist}\\ given (FM),\\ critical points\\ at rate $\Sigma(f)$};
\node[thm, text width=5.055cm] (thmThree) at (8.15,0)
  {Theorem~\ref{thm:struct}\\ near-maximizers lie near $\Cf$,\\
   field law $\Lq{q}$; Hessian bulk $\approx\nu_q$\\
   with edge $\lamp\le0$, zero iff $\Plef=0$};
\node[thm, text width=5.055cm] (thmFive) at (13.70,0)
  {Theorem~\ref{thm:stab}\\ with two active degrees, stable\\ critical points of $F_{\zf}$ at
   rates\\ near $\Sigma(f)$. Under (H-marg)\\ every fixed margin costs rate};
\node[thm, text width=7.135cm] (corFour) at (3.64,2.05)
  {Corollary~\ref{cor:title}\\ conditionally on $N^{-1}\log\ZN\ge f$, the deviation\\
   is carried by TAP states as in Theorem~\ref{thm:struct}};
\node[thm, text width=2.985cm] (corSix) at (9.185,2.05)
  {Corollary~\ref{cor:eq}\\ at level $\feq$,\\ near-maximizers\\ lie near $\CP$};
\node[ours, text width=2.215cm] (eqZero) at (12.35,2.05)
  {$0\in\CP$ for\\ every Parisi\\ measure\rmloc{cf.\ \cite{AC15a}; Lem.~\ref{lem:eq-zero}}};
\node[ours]  (colOneA) at (1.35,-2.3) {Ramsey: $2n$\\ replicas in one\\ overlap bin%
  \rmloc{\cite{Ram30};\\ Lem.~\ref{lem:max-ramsey}}};
\node[redo]  (colOneB) at (1.35,-4.5) {Band theorem\\ under (R)%
  \rmloc{\cite[Thm.~1]{CPS}\\ Lem.~\ref{lem:max-truncation}}};
\node[redo]  (colOneC) at (1.35,-6.6) {Band step at\\ a barycenter,\\ at every\\ overlap level%
  \rmloc{\cite[\S5]{CPS}\\ Lem.~\ref{lem:max-band}}};
\node[redo]  (colOneD) at (3.85,-2.3) {Maximizers\\ are critical\\ points of $F_{\zeta^*}$%
  \rmloc{\cite[\S11]{CPS}\\ Lem.~\ref{lem:gauss-critical}}};
\node[known] (colOneE) at (3.85,-4.5) {(FM): even $\xi$;\\ upper half for\\ convex $\xi$%
  \rmloc{\cite{Tal07}; Prop.~\ref{prop:max-fm}(a)--(c)}};
\node[ours]  (colTwoA) at (6.90,-2.3) {Variational gap\\ off $\Cf$, from a\\ downward jump%
  \rmloc{Prop.~\ref{prop:loc-gap}}};
\node[redo]  (colTwoB) at (6.90,-4.3) {Comparison\\ on spheres\\ and shells%
  \rmloc{\cite{Bou26b}\\ Prop.~\ref{prop:glued-comparison}}};
\node[known] (colTwoC) at (6.90,-6.5) {Optimal order\\ parameter of a\\ band problem\\ is unique%
  \rmloc{\cite[Thm.~12(iii)]{CPS}\\ Prop.~\ref{prop:hess-own}}};
\node[redo]  (colTwoD) at (9.40,-2.3) {Hessian bulk\\ near $\nu_{q_x}$ if\\ $\mu_x\approx\Lq{q_x}$%
  \rmloc{\cite{GSS,Bou26b}\\ Prop.~\ref{prop:hess-bulk}}};
\node[redo]  (colTwoE) at (9.40,-4.3) {Edge $\lamp\le0$,\\ zero iff $\Plef=0$%
  \rmloc{\cite{GSS,Sel25}\\ Lem.~\ref{lem:hess-edge}}};
\node[ours]  (colTwoF) at (9.40,-6.5) {Own order\\ parameters stay\\ close to $\zf$%
  \rmloc{Prop.~\ref{prop:hess-own}}};
\node[ours]  (colThreeA) at (12.45,-2.3) {$2+4$ family:\\ every contact\\ is marginal%
  \rmloc{Lem.~\ref{lem:stab-family}, App.~\ref{app:nc}}};
\node[redo]  (colThreeB) at (12.45,-4.5) {Positive part\\ of $\supp\zf$ is\\ an interval%
  \rmloc{\cite{Zho26a}\\ Prop.~\ref{prop:nc-example}}};
\draw[use] (thmOne.east |- 0,0.3) -- (thmTwo.west |- 0,0.3);
\draw[use] (thmThree.east) -- (thmFive.west);
\draw[use] (thmOne.north) -- (thmOne.north |- corFour.south);
\draw[use] (thmTwo.north) -- (thmTwo.north |- corFour.south);
\draw[use] (6.60,0 |- thmThree.north) -- (6.60,0 |- corFour.south);
\draw[use] (eqZero.west) -- (corSix.east);
\draw[use] (colOneA) -- (colOneA |- thmOne.south);
\draw[use] (colOneB.west) -- (0.125,-4.5) -- (0.125,0 |- thmOne.south);
\draw[use] (colOneC.east) -- (2.60,-6.6) -- (2.60,-0.3) -- (thmOne.east |- 0,-0.3);
\draw[use] (colOneD) -- (colOneD |- thmTwo.south);
\draw[use] (colOneE.east) -- (5.075,-4.5) -- (5.075,0 |- thmTwo.south);
\draw[use] (colTwoA) -- (colTwoA |- thmThree.south);
\draw[use] (colTwoB.west) -- (5.675,-4.3) -- (5.675,0 |- thmThree.south);
\draw[use] (colTwoC.east) -- (colTwoF.west);
\draw[use] (colTwoD) -- (colTwoD |- thmThree.south);
\draw[use] (colTwoE.east) -- (10.625,-4.3) -- (10.625,0 |- thmThree.south);
\draw[use] (colTwoF.north) -- (9.40,-5.6) -- (8.15,-5.6) -- (8.15,0 |- thmThree.south);
\draw[use] (colThreeA) -- (colThreeA |- thmFive.south);
\draw[use] (colThreeB) -- (colThreeA);
\draw[rel] (9.185,0 |- thmThree.north) -- (9.185,0 |- corSix.south);
\node[tool, text width=7.9cm] (PDE) at (9.70,-8.35)
  {Parisi equation with nonnegative, non-monotone coefficients;
   Auffinger--Chen process (App.~\ref{app:pde})};
\node[tool, text width=15.5cm] (KR) at (8.15,-9.25)
  {Kac--Rice construction of isolated critical points of $F_{\zf}$ \cite{Bou26b}, after \cite{HS},\\
   extended to a marginal first contact by a field-law tilt (Prop.~\ref{prop:kr-main}, App.~\ref{app:kr})};
\draw[band] (PDE.west) -- (4.60,0 |- PDE.west) -- (4.60,-7.75);
\draw[band] (6.90,0 |- PDE.north) -- (6.90,-7.75);
\draw[band] (12.45,0 |- PDE.north) -- (12.45,-7.75);
\draw[band] (3.85,0 |- KR.north) -- (3.85,-7.75);
\draw[band] (14.95,0 |- KR.north) -- (14.95,-7.75);
\node[sw, thmlook]  (g1) at (0.45,-10.25) {}; \node[lg] at (g1.east) {main result};
\node[sw, ourslook, rounded corners=0.6mm] (g2) at (3.35,-10.25) {};
  \node[lg] at (g2.east) {our ingredient};
\node[sw, toollook, minimum width=0.8cm] (g3) at (6.45,-10.25) {};
  \node[lg] at (g3.east) {shared technique};
\node[sw, knownlook] (g4) at (0.45,-10.75) {}; \node[lg] at (g4.east) {existing result};
\node[sw, redolook]  (g5) at (3.35,-10.75) {};
  \node[lg] at (g5.east) {existing result, adapted: extended or re-proved};
\draw[use]  (0.20,-11.25) -- (0.70,-11.25) node[lg] {used in proof};
\draw[rel]  (3.10,-11.25) -- (3.60,-11.25) node[lg] {other relation};
\draw[band] (6.05,-11.25) -- (6.85,-11.25) node[lg] {shared technique used in the column};
\end{tikzpicture}%
\end{NoHyper}
\caption{How the main results are proved. The three columns lead to Theorems~\ref{thm:max}
and~\ref{thm:exist}, to Theorem~\ref{thm:struct} and to Theorem~\ref{thm:stab}; the two
corollaries sit above them. A solid arrow means that its source is used in the proof of its
target, and the grey arrows from the two bottom bands mark the shared techniques that each column
uses; only the main dependencies are drawn. The arrow from the Kac--Rice band into the first
column is Proposition~\ref{prop:max-fm}\ref{prop:max-fm-d}, the lower half of (FM) for convex
$\xi$ that is not even (ours). The dashed arrow marks the one relation of another
kind: Corollary~\ref{cor:eq} is proved by running the argument of
Theorem~\ref{thm:struct}\ref{thm:struct-i}--\ref{thm:struct-ii} at the equilibrium level. Box
styles give provenance; grey labels give sources and where each item is proved.
Section~\ref{sec:intro-ideas} reads the figure one column at a time and explains the terms in
the boxes, and
Section~\ref{sec:intro-prior} lists the existing results that we re-prove or extend.}
\label{fig:roadmap}
\end{figure}
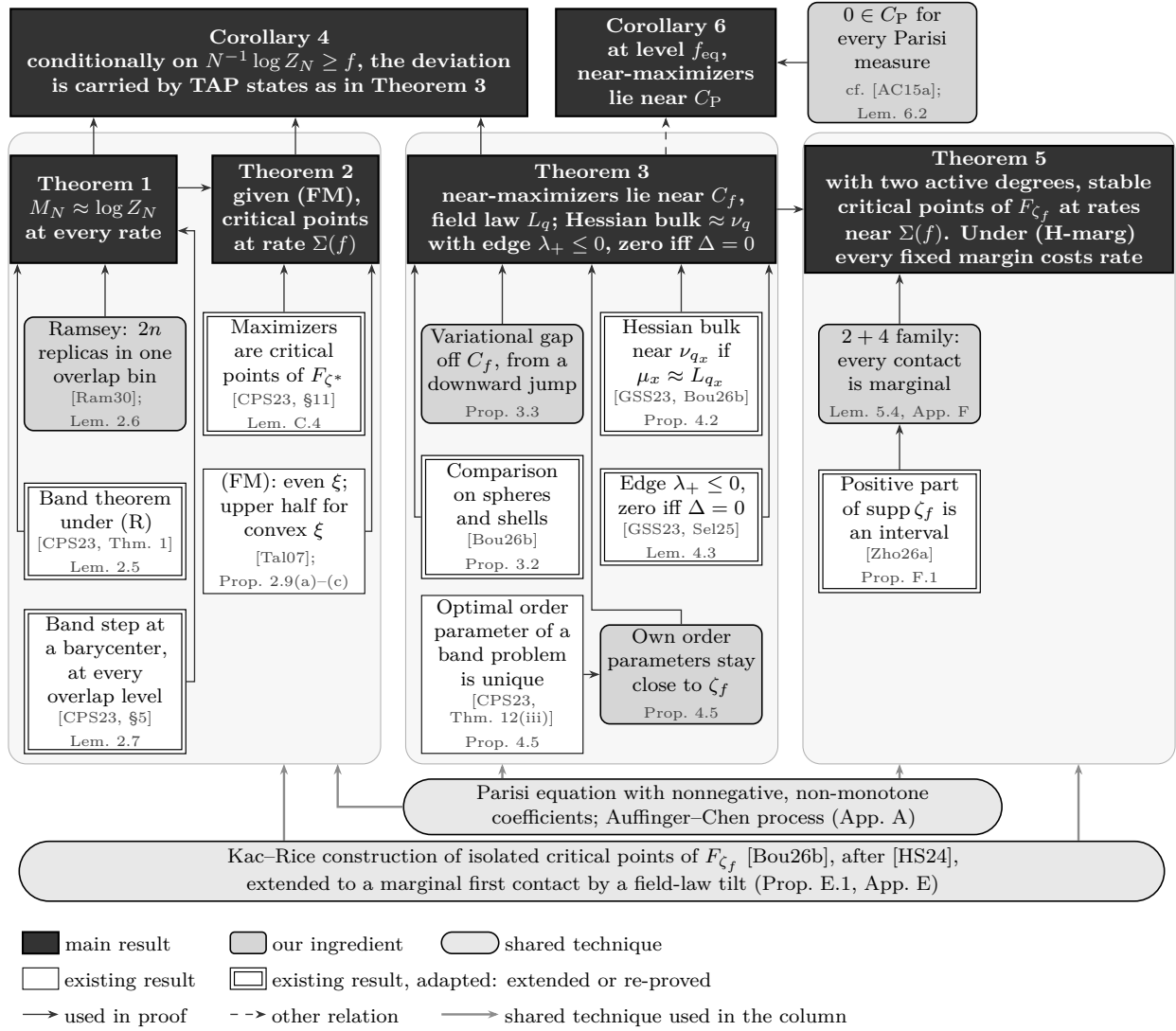

%% file: sections/sec2-maximum.tex
\section{The TAP maximum carries the free energy:
\texorpdfstring{Theorems~\ref*{thm:max} and~\ref*{thm:exist}}{Theorems 1 and 2}}\label{sec:max}

This section proves Theorems~\ref{thm:max} and~\ref{thm:exist}, the first column of
Figure~\ref{fig:roadmap}. Chen, Panchenko and Subag connect $\log\ZN$ and $\MN$ through the free
energy of a thin band of configurations around a point $m$ of the cube: it never exceeds
$\log\ZN$, and it is close to $\FTAP(m)$ uniformly in $m$, outside an event of probability
$e^{-cN}$ for every $c>0$ (\cite[Theorem~1]{CPS}; Lemma~\ref{lem:max-truncation} under \condR). The
upper bound $\MN\le\log\ZN+tN$ follows at once. The lower bound needs a point whose band carries
a share of the Gibbs measure that is not exponentially small; Ramsey's theorem finds one for every realization of the disorder, at a
self-overlap that depends on the realization. Both bounds then hold at every exponential rate,
which is what lets large deviation statements pass between $\MN$ and $\log\ZN$ in both
directions. Theorem~\ref{thm:exist} adds three inputs: maximizers of $\FTAP$ are critical points,
every level is regular, and the fractional-moment formula holds for convex mixtures.

\subsection{Two band estimates and the proof of
\texorpdfstring{Theorem~\ref*{thm:max}}{Theorem 1}}\label{sec:max-bands}

We deduce Theorem~\ref{thm:max} from two band estimates, which hold at every exponential rate.
The upper one is a direct consequence of the band theorem of Chen, Panchenko and Subag; the lower
one is the main step and uses Ramsey's theorem. Both are proved in Section~\ref{sec:max-ramsey}.
Under \condR\ the field $\HN$ is almost surely continuously differentiable on a neighborhood of
the cube, $\FTAP$ is continuous on $\cube$, and $\MN$ is attained
(Lemma~\ref{lem:gauss-interior}); in particular $\MN$ is a random variable.

\begin{proposition}[Lower band]\label{prop:max-lower}
Assume \condR. For all $c,t>0$ and all large $N$,
\[
 \Prob\bigl(\MN<\log\ZN-tN\bigr)\le e^{-cN}.
\]
\end{proposition}

\begin{proposition}[Upper band]\label{prop:max-upper}
Assume \condR. For all $c,t>0$ and all large $N$,
\[
 \Prob\bigl(\MN>\log\ZN+tN\bigr)\le e^{-cN}.
\]
\end{proposition}

Besides the two propositions, the deduction uses only crude a priori bounds: for every
$\theta\ne0$ there is $C=C(\xi,\theta)<\infty$ such that, for all $N$,
\begin{equation}\label{eq:max-crude}
 \E\ZN^{\theta}\ge e^{-CN},\qquad
 \E\ZN^{2\theta}\le e^{CN},\qquad
 \E e^{2\theta\MN}\le e^{CN}
\end{equation}
(Lemma~\ref{lem:gauss-bounds}). They follow from Jensen's inequality
$\log\ZN\ge2^{-N}\sum_\sigma\HN(\sigma)$, a centered Gaussian variable of variance at most
$N\xi(1)$; from Gaussian tails of $\sup_{\cube}\HN$ around its mean, which is of order $N$, with
variance proxy $N\xi(1)$; and from the deterministic bounds
\begin{equation}\label{eq:max-crude-M}
 N\feq=\FTAP(0)\le\MN\le\sup_{\cube}\HN+\tfrac N2\,\xi'(1).
\end{equation}
The identity holds because $\HN(0)=0$ and $\Phi_\zeta(0,\cdot)$ is even and convex, so that
$F_\zeta(0)=NP_\zeta$ for every $\zeta$; the upper bound comes from the trial measure
$\zeta=\delta_1$ in the definition of $\FTAP$.

\begin{proof}[Proof of Theorem~\ref{thm:max}, given Propositions~\ref{prop:max-lower}
and~\ref{prop:max-upper}]
Attainment of $\MN$ was recalled above.

\ref{thm:max-a} By the two propositions, for every $c>0$,
$\Prob(|\MN-\log\ZN|>tN)\le2e^{-cN}$ for all large $N$.

\ref{thm:max-b} The case $\theta=0$ is trivial. Fix $\theta\ne0$ and $t>0$, let $C$ be as in
\eqref{eq:max-crude}, and put $G_t=\{|\MN-\log\ZN|\le tN\}$. By~\ref{thm:max-a},
$\Prob(G_t^c)\le e^{-(3C+2)N}$ for all large $N$. On $G_t$,
$e^{-|\theta|tN}\ZN^\theta\le e^{\theta\MN}\le e^{|\theta|tN}\ZN^\theta$. Off $G_t$, the
Cauchy--Schwarz inequality and \eqref{eq:max-crude} give
\[
 \E\bigl[e^{\theta\MN};G_t^c\bigr]\le\bigl(\E e^{2\theta\MN}\bigr)^{1/2}\Prob(G_t^c)^{1/2}
 \le e^{-(C+1)N}\le e^{-N}\,\E\ZN^\theta,
\]
and the same bound holds for $\E[\ZN^\theta;G_t^c]$. Hence, for all large $N$,
\[
 (1-e^{-N})\,e^{-|\theta|tN}\,\E\ZN^\theta\ \le\ \E e^{\theta\MN}\ \le\
 \bigl(e^{|\theta|tN}+e^{-N}\bigr)\,\E\ZN^\theta .
\]
So $\limsup_N N^{-1}\bigl|\log\E e^{\theta\MN}-\log\E\ZN^\theta\bigr|\le|\theta|t$, and $t>0$ is
arbitrary.

\ref{thm:max-c} For $0<t<\varepsilon$ the inclusions
\[
 \bigl\{|N^{-1}\log\ZN-f|\le\varepsilon-t\bigr\}\cap G_t
 \subset\bigl\{|N^{-1}\MN-f|\le\varepsilon\bigr\}
 \subset\bigl\{|N^{-1}\log\ZN-f|\le\varepsilon+t\bigr\}\cup G_t^c
\]
and~\ref{thm:max-a} give the two displayed bounds. Since $c$ is arbitrary, the additive error
$e^{-cN}$ changes neither $\liminf_NN^{-1}\log$ nor $\limsup_NN^{-1}\log$ of the probabilities
involved. The lower bound with $t=\varepsilon/2$ shows that the local rates of $N^{-1}\MN$ at
window $\varepsilon$ are at least those of $N^{-1}\log\ZN$ at window $\varepsilon/2$. The upper
bound with $(\varepsilon,t)$ replaced by $(\varepsilon/2,\varepsilon/4)$ shows that the local rates of
$N^{-1}\log\ZN$ at window $\varepsilon$ are at least those of $N^{-1}\MN$ at window
$\varepsilon/2$. All of them are monotone in $\varepsilon$, so their limits as
$\varepsilon\downarrow0$ coincide. Finally,~\ref{thm:max-a} says that $N^{-1}\MN$ and
$N^{-1}\log\ZN$ are exponentially equivalent, and exponentially equivalent families satisfy the
same large deviation principles with good rate functions \cite[Section~4.2]{DZ}.
\end{proof}

Tilting the law of the disorder by $\ZN^\theta$ is how upper deviations of the free energy are
produced. For $\theta>0$ let $\Qth(A)=\E[\ZN^\theta\ind_A]/\E\ZN^\theta$; this is well defined
since $\ZN\le e^{\max\HN}$.

\begin{corollary}[Under every tilt]\label{cor:max-tilt}
Assume \condR. For all $\theta>0$ and $t>0$,
\[
 \lim_{N\to\infty}\frac1N\log\Qth\bigl(|\MN-\log\ZN|>tN\bigr)=-\infty .
\]
\end{corollary}

\begin{proof}
By the Cauchy--Schwarz inequality and \eqref{eq:max-crude}, every event $A$ satisfies
\[
 \Qth(A)\le\frac{(\E\ZN^{2\theta})^{1/2}\,\Prob(A)^{1/2}}{\E\ZN^\theta}\le e^{2CN}\Prob(A)^{1/2}.
\]
Apply Theorem~\ref{thm:max}\ref{thm:max-a}.
\end{proof}

\begin{remark}[Every rate versus a fixed rate]\label{rem:max-fixedrate}
Every exponential rate is what transfers large deviations in both directions: in the proofs of
Theorem~\ref{thm:max}\ref{thm:max-b}--\ref{thm:max-c} the error event must be negligible against
moments and probabilities that are themselves exponentially large or small. For $\xi\in\classC$
the upper band at every rate is already in \cite{CPS}, as their Theorem~1 combined with
\eqref{eq:max-band-le} below, and Boursier uses it in this way \cite[Proposition~3.12]{Bou26b}.
The gain is in the lower band.
\end{remark}

\subsection{Ramsey multiplicity and the band estimates}\label{sec:max-ramsey}

We prove Propositions~\ref{prop:max-upper} and~\ref{prop:max-lower}: the upper band from the band
theorem of Chen, Panchenko and Subag, the lower band from Ramsey's theorem and a band step.

\paragraph{Bands.} For $m\in\cube$, $\varepsilon,\delta>0$ and an integer $n\ge1$, Chen, Panchenko
and Subag consider \cite[(2.10), (2.12)]{CPS}
\begin{align*}
 B(m,\varepsilon)&=\bigl\{\sigma\in\hcube:\ |\ov{\sigma}{m}-q_m|<\varepsilon\bigr\},\\
 \Band&=\bigl\{(\sigma^1,\dots,\sigma^n)\in B(m,\varepsilon)^n:\
 |\ov{\sigma^i}{\sigma^j}-q_m|<\delta\ \text{for all } i\ne j\bigr\}.
\end{align*}
Heuristically, if $m$ is an ancestor state of the Parisi tree, the band $B(m,\varepsilon)$ carries
the pure states that descend from $m$, and configurations drawn from distinct descendants have
nearly orthogonal deviations $\sigma^i-m$, which inside the band is the condition defining $\Band$
\cite[Section~2.2]{CPS}. We use the band free energy
\begin{equation}\label{eq:max-band}
 F^{\mathrm{band}}_n(m)=F^{\mathrm{band}}_n(m;\varepsilon,\delta)
 =\frac1n\log\Bigl(2^{-nN}\!\!\!\sum_{(\sigma^1,\dots,\sigma^n)\in\Band}\!\!\!
  \exp\Bigl(\sum_{i\le n}\HN(\sigma^i)\Bigr)\Bigr),
\end{equation}
with $\log0=-\infty$. Writing $G_N(\sigma)=2^{-N}e^{\HN(\sigma)}/\ZN$ for the Gibbs measure on
$\hcube$, the definition reads
\begin{equation}\label{eq:max-band-gibbs}
 F^{\mathrm{band}}_n(m;\varepsilon,\delta)
 =\log\ZN+\frac1n\log G_N^{\otimes n}\bigl(\Band\bigr).
\end{equation}
In particular, for every function $\HN$ on the cube and all $m,\varepsilon,\delta,n$,
\begin{equation}\label{eq:max-band-le}
 F^{\mathrm{band}}_n(m;\varepsilon,\delta)\le\log\ZN .
\end{equation}
In the notation of \cite{CPS},
$F^{\mathrm{band}}_n(m)=\HN(m)+N\,\mathrm{TAP}_{N,n}(m,\varepsilon,\delta)-N\log2$
\cite[(2.14)]{CPS} and $\FTAP(m)=\HN(m)+N\,\TAPcps(\mu_m)-N\log2$, where $\TAPcps$ is the
functional of \cite[(2.24)--(2.25)]{CPS}; the shift by $N\log2$ comes from the factor $2^{-nN}$
in~\eqref{eq:max-band} and, for $\FTAP$, from their terminal condition
$\log(2\cosh x)$ in place of our $\log\cosh x$ (Section~\ref{sec:intro-setting}). Thus
\eqref{eq:max-band-le} is \cite[(2.16)]{CPS}, and their Theorem~1 says that
$F^{\mathrm{band}}_n(m)$ is uniformly close to $\FTAP(m)$. We use it in the following form.

\begin{lemma}[The band theorem at every exponential rate]\label{lem:max-truncation}
Assume \condR. For all $c,t>0$ and $L\ge1$ there are $\varepsilon,\delta>0$ and an integer $n$
with $n\ge L$ and $n\ge L/\min(\varepsilon,\delta)$ such that, for all large $N$,
\begin{equation}\label{eq:max-band-thm}
 \Pout\Bigl(\exists m\in\cube:\
 \bigl|F^{\mathrm{band}}_n(m;\varepsilon,\delta)-\FTAP(m)\bigr|\ge tN\Bigr)\le e^{-cN}.
\end{equation}
\end{lemma}

For $\xi\in\classC$ this is \cite[Theorem~1]{CPS}, read with the order of quantifiers of its proof:
$\varepsilon$ and $\delta$ are fixed small first, and $n$ is taken large afterwards
\cite[Sections~4.1--4.2]{CPS}, so $n$ may be large compared with $1/\min(\varepsilon,\delta)$.
Under \condR\ alone we split $\HN$ by degree into a polynomial part of degree at most $D$ and a
remainder $\HN^{>D}$ with covariance $N\xi_{>D}$, $\xi_{>D}(t)=\sum_{p>D}\beta_p^2t^p$, and apply
\cite[Theorem~1]{CPS} to the polynomial part. The remainder moves $F^{\mathrm{band}}_n$ and $\FTAP$
by at most $\sup_{\cube}|\HN^{>D}|+\frac32N\xi'_{>D}(1)$. Once $D$ is large, this is a small
multiple of $N$ outside an event of probability $e^{-cN}$, because $N^{-1}\E\sup_{\cube}|\HN^{>D}|$
and $\xi_{>D}(1)$, the variance proxy of $\HN^{>D}$ per unit $N$, are small.
The details are in Appendix~\ref{app:gauss-truncation}.

\begin{proof}[Proof of Proposition~\ref{prop:max-upper}]
Take $\varepsilon,\delta,n$ from Lemma~\ref{lem:max-truncation} with the given $c,t$ and $L=1$.
Outside the event in \eqref{eq:max-band-thm}, every $m\in\cube$ satisfies
$\FTAP(m)<F^{\mathrm{band}}_n(m)+tN\le\log\ZN+tN$ by \eqref{eq:max-band-le}. Hence
$\MN<\log\ZN+tN$ there.
\end{proof}

This uses one half of the band theorem, $\FTAP\le F^{\mathrm{band}}_n+tN$. The lower band uses the
other half, and needs in addition a point $m$ whose band is heavy.

\paragraph{The idea for the lower band.} By \eqref{eq:max-band-thm} and \eqref{eq:max-band-gibbs},
it suffices to find, for every realization, a point $m$ with
$G_N^{\otimes n}(\Band)\ge e^{-o(N)}$. A tuple in $\Band$ has its pairwise overlaps, and its
overlaps with $m$, close to $q_m$. Conversely, if $2n$ configurations have all their pairwise
overlaps close to a common value $q$, then the barycenter $m$ of the last $n$ of them has
self-overlap close to $q$, and its band contains the first $n$ (Lemma~\ref{lem:max-band}). It
therefore suffices that, with $G_N^{\otimes 2n}$-probability bounded below, all pairwise overlaps
of $2n$ replicas fall close to a common value. Cut the overlap axis $[-1,1]$ into finitely many
bins. Whatever the measure $G_N$ is, Ramsey's theorem forces some bin to catch all pairwise
overlaps of $2n$ independent replicas with probability bounded below by a constant that depends
only on $n$ and the number of bins (Lemma~\ref{lem:max-ramsey}). Figure~\ref{fig:ramsey} shows the
three steps.

\input{figures/fig2-ramsey}

Recall Ramsey's theorem in its finite multicolor form \cite{Ram30}: for all integers $J\ge1$ and
$s\ge2$ there is $r_J(s)<\infty$ such that, whenever the edges of the complete graph on $r_J(s)$
vertices are colored with $J$ colors, some $s$ vertices span a complete subgraph whose edges all
have the same color.

\begin{lemma}[Ramsey multiplicity]\label{lem:max-ramsey}
Let $J\ge1$ and $s\ge2$, let $[-1,1]=I_1\sqcup\dots\sqcup I_J$ be a partition into Borel sets, and
put $c_{J,s}=\bigl[J\binom{r_J(s)}{s}\bigr]^{-1}$. For every $N$ and every probability measure $G$
on $\hcube$,
\[
 \max_{k\le J}\ G^{\otimes s}\bigl(\ov{\sigma^i}{\sigma^j}\in I_k\ \text{for all }
 1\le i<j\le s\bigr)\ \ge\ c_{J,s}.
\]
\end{lemma}

\begin{proof}
Let $r=r_J(s)$ and sample $\sigma^1,\dots,\sigma^r$ independently from $G$. Color each pair
$\{i,j\}\subset\{1,\dots,r\}$ by the index $k$ with $\ov{\sigma^i}{\sigma^j}\in I_k$; two equal
samples have overlap one and receive the color of the bin that contains one. By Ramsey's theorem,
for every outcome some $s$-element set $S\subset\{1,\dots,r\}$ is monochromatic, so
\[
 \sum_{|S|=s}\ \sum_{k\le J}\ind\bigl\{\text{all pairs in } S \text{ have color } k\bigr\}\ge1
\]
pointwise. Take expectations. For each $S$ the samples $(\sigma^i)_{i\in S}$ are independent with
law $G$, so each of the $\binom rs$ sets contributes
$\sum_kG^{\otimes s}(\ov{\sigma^i}{\sigma^j}\in I_k\ \forall i<j)$, and this sum is therefore at
least $\binom rs^{-1}$. The largest of its $J$ terms is at least the average.
\end{proof}

The bound depends neither on $N$ nor on $G$. The proof uses only exchangeability of independent
replicas and a deterministic combinatorial fact; it uses no property of the Gibbs measure.

For $q\in[-1,1]$, $\tau>0$ and $s\ge2$ let
\[
 A_s(q,\tau)=\bigl\{(\sigma^1,\dots,\sigma^s)\in(\hcube)^s:\
 |\ov{\sigma^i}{\sigma^j}-q|<\tau\ \text{for all } i\ne j\bigr\}.
\]

\begin{lemma}[Band step]\label{lem:max-band}
Let $H$ be a real function on $\cube$, $Z=2^{-N}\sum_{\sigma\in\hcube}e^{H(\sigma)}$ and
$G(\sigma)=2^{-N}e^{H(\sigma)}/Z$, and define $F^{\mathrm{band}}_n$ by \eqref{eq:max-band} with $H$
in place of $\HN$. Let $q\in[-1,1]$, $\varepsilon,\delta>0$, let $n\ge4/\min(\varepsilon,\delta)$ be
an integer and $\tau\in(0,\min(\varepsilon,\delta)/4]$. If $G^{\otimes2n}(A_{2n}(q,\tau))>0$, there
is $m\in\cube$, the barycenter of $n$ points of $\hcube$, such that
$|q_m-q|<(1-q)/n+\tau$ and
\begin{equation}\label{eq:max-band-step}
 F^{\mathrm{band}}_n(m;\varepsilon,\delta)\ \ge\ \log Z+\frac1n\log G^{\otimes2n}\bigl(A_{2n}(q,\tau)\bigr).
\end{equation}
\end{lemma}

\begin{proof}
\emph{An outer block.} Split a $2n$-tuple into an inner block $(\sigma^1,\dots,\sigma^n)$ and an
outer block $\bar\sigma=(\bar\sigma^{n+1},\dots,\bar\sigma^{2n})$. For an outer block whose pairwise
overlaps lie in $(q-\tau,q+\tau)$, let $A_n(\bar\sigma)$ be the set of inner blocks that complete
it to a point of $A_{2n}(q,\tau)$. Then
\[
 G^{\otimes2n}\bigl(A_{2n}(q,\tau)\bigr)=\sum_{\bar\sigma}G^{\otimes n}(\bar\sigma)\,
 G^{\otimes n}\bigl(A_n(\bar\sigma)\bigr),
\]
the sum running over such outer blocks. The weights $G^{\otimes n}(\bar\sigma)$ sum to at most
one, so some outer block $\bar\sigma$ in the sum has
$G^{\otimes n}(A_n(\bar\sigma))\ge G^{\otimes2n}(A_{2n}(q,\tau))$. Fix it.

\emph{Its barycenter.} Let $m=n^{-1}\sum_{k=n+1}^{2n}\bar\sigma^k\in\cube$. Then
\[
 q_m=\frac1n+\frac1{n^2}\sum_{k\ne k'}\ov{\bar\sigma^k}{\bar\sigma^{k'}},\qquad
 q_m-q=\frac{1-q}n+\frac1{n^2}\sum_{k\ne k'}\bigl(\ov{\bar\sigma^k}{\bar\sigma^{k'}}-q\bigr),
\]
and the last sum has $n(n-1)$ terms of modulus less than $\tau$. Hence
$|q_m-q|<(1-q)/n+\tau\le2/n+\tau$. Now let $(\sigma^1,\dots,\sigma^n)\in A_n(\bar\sigma)$. Each
$\ov{\sigma^i}{m}=n^{-1}\sum_k\ov{\sigma^i}{\bar\sigma^k}$ lies in $(q-\tau,q+\tau)$, and so does
each $\ov{\sigma^i}{\sigma^j}$ with $i\ne j$. Therefore $|\ov{\sigma^i}{m}-q_m|$ and
$|\ov{\sigma^i}{\sigma^j}-q_m|$ are less than $2\tau+2/n\le\min(\varepsilon,\delta)$, by the choice
of $n$ and $\tau$. So $A_n(\bar\sigma)\subset\Band$.

\emph{The bound.} By \eqref{eq:max-band-gibbs} with $H$ in place of $\HN$,
\begin{align*}
 F^{\mathrm{band}}_n(m)=\log Z+\frac1n\log G^{\otimes n}\bigl(\Band\bigr)
 &\ge\log Z+\frac1n\log G^{\otimes n}\bigl(A_n(\bar\sigma)\bigr)\\
 &\ge\log Z+\frac1n\log G^{\otimes2n}\bigl(A_{2n}(q,\tau)\bigr).\qedhere
\end{align*}
\end{proof}

The point $m$ is the barycenter itself. It lies in the cube for every level $q\in[-1,1]$, negative
levels included, and its self-overlap need not be prescribed; in \cite[Section~5]{CPS} the
barycenter is rescaled to the sphere of the prescribed self-overlap $q$.

\begin{proof}[Proof of Proposition~\ref{prop:max-lower}]
Fix $c,t>0$. Lemma~\ref{lem:max-truncation} with $(c,t/2,L=4)$ gives $\varepsilon,\delta>0$ and an
integer $n\ge4/\min(\varepsilon,\delta)$ such that the event
\[
 \Omega_N^{\mathrm{bad}}=\bigl\{\exists m\in\cube:\
 |F^{\mathrm{band}}_n(m;\varepsilon,\delta)-\FTAP(m)|\ge tN/2\bigr\}
\]
has outer probability at most $e^{-cN}$ for all large $N$. Put $\tau=\min(\varepsilon,\delta)/4$
and $J=\lceil2/\tau\rceil$, and cut $[-1,1]$ into the intervals
$I_k=[-1+2(k-1)/J,\,-1+2k/J)$ for $k<J$ and $I_J=[1-2/J,1]$, with midpoints $q_k$. An overlap in
$I_k$ lies within $1/J\le\tau/2$ of $q_k$, so the event that all pairwise overlaps lie in $I_k$ is
contained in $A_{2n}(q_k,\tau)$. Lemma~\ref{lem:max-ramsey} with $s=2n$ and $G=G_N$ gives, for
every $N$ and every realization, an index $k$ with
\[
 G_N^{\otimes2n}\bigl(A_{2n}(q_k,\tau)\bigr)\ge c_{J,2n}.
\]
Lemma~\ref{lem:max-band} with $H=\HN$ then gives $m\in\cube$ with
$F^{\mathrm{band}}_n(m)\ge\log\ZN-n^{-1}\log(1/c_{J,2n})$. Outside $\Omega_N^{\mathrm{bad}}$,
$\FTAP(m)>F^{\mathrm{band}}_n(m)-tN/2$, so
\[
 \MN\ge\FTAP(m)>\log\ZN-\frac{tN}2-\frac1n\log\frac1{c_{J,2n}}\ge\log\ZN-tN
\]
as soon as $N\ge2\log(1/c_{J,2n})/(nt)$. Thus $\{\MN<\log\ZN-tN\}\subset\Omega_N^{\mathrm{bad}}$ for
all large $N$, and the claim follows from \eqref{eq:max-band-thm}.
\end{proof}

The constants $\varepsilon,\delta,n,\tau,J$ and $c_{J,2n}$ depend only on $\xi$, $c$ and $t$. The
barycenter $m$ is random and may have coordinates equal to $\pm1$, which is why the band theorem is
needed uniformly over the closed cube.

\begin{remark}[Where the barycenter sits]\label{rem:max-witness}
The barycenter $m$ has self-overlap within $2/n+\tau$ of the bin midpoint $q_k$; the bin depends on
the realization, and nothing in the argument locates it. It may even be negative: positive
semidefiniteness of the overlap matrix of $2n$ replicas only forces $q_k>-1/(2n-1)-\tau$. For convex
$\xi$ the variational side locates the level: outside an event of outer probability at most
$e^{N(\Sigma(f)-c)}$, every $m$ with $\FTAP(m)\ge N(f-\delta)$ sits near the contact set $\Cf$
(Theorem~\ref{thm:struct}\ref{thm:struct-i}), and for convex $\xi\in\classC$ this holds conditionally
on an upper deviation of the free energy (Corollary~\ref{cor:title}). For every mixture,
Proposition~\ref{rem:max-overlap} in
Appendix~\ref{app:var-overlap} gives a cruder lower bound on its self-overlap.
\end{remark}

\subsection{Generalized TAP critical points: proof of
\texorpdfstring{Theorem~\ref*{thm:exist}}{Theorem 2}}\label{sec:max-existence}

Theorem~\ref{thm:exist} turns an upper deviation of $\MN$, which by Theorem~\ref{thm:max} is an
upper deviation of the free energy, into a statement about critical points. Its proof adds three
inputs to Theorem~\ref{thm:max}, none of them, for parts~\ref{thm:exist-a}--\ref{thm:exist-d}, a
construction: every maximizer of $\FTAP$ is an
interior critical point of the TAP function $F_{\zeta^*}$ with its own order parameter; the
duality that defines $\Sigma$ is regular at every level of $(\feq,\fone)$; and the
fractional-moment formula identifies, through Theorem~\ref{thm:max}\ref{thm:max-b}, the
exponential moments of $\MN$, from which the rate follows by the standard tilting argument.
For $\xi$ that is not even, part~\ref{thm:exist-e} rests on the Kac--Rice construction through
Proposition~\ref{prop:max-fm}\ref{prop:max-fm-d}.

\paragraph{Criticality of maximizers.} Part~\ref{thm:exist-a} is
Lemmas~\ref{lem:gauss-interior} and~\ref{lem:gauss-critical}. Maximizers are interior because
moving a coordinate inward from $\pm1$ by $s$ raises the entropy term by at least $s\,\psi(1-s)$,
where $\psi(a)\to\infty$ as $a\uparrow1$ uniformly in $\zeta$ and $q$, while $\HN$ and the
self-overlap change by $O(s)$. At an interior point the directional derivatives of $\FTAP$ and of
the frozen function $F_{\zeta^*}$ coincide, by the envelope argument behind the gradient formula
\cite[Theorem~5 and Section~11]{CPS}, so they vanish at a maximum. The class $\classC$ enters Theorem~\ref{thm:exist} only through part~\ref{thm:exist-a},
which~\ref{thm:exist-d} uses for the existence event, and through Proposition~\ref{prop:max-fm}
in part~\ref{thm:exist-e} (Remark~\ref{rem:rpde}); part~\ref{thm:exist-b} uses only \condR.

\paragraph{Regularity of every level.} Part~\ref{thm:exist-b} is Lemma~\ref{lem:var-basic} and
Proposition~\ref{prop:var-regular} (Appendix~\ref{app:var-constrained}); it uses only \condR, and in
particular not \condFM.

\paragraph{Fractional moments.} The formula \condFM\ is not ours as a statement: it is due to
Talagrand for even $\xi$ \cite[Section~7]{Tal07}, argued there by stated modifications of the
proof of \cite{PF}, with the upper half for convex $\xi$ from Guerra's bound
\cite[Theorem~6.1]{Tal07}, and it is stated for every $\xi$ finite on $\R$ in
\cite[Theorem~2.1]{CGKLM}; for convex $\xi$ that is not even we prove its lower half below
(item~\ref{prop:max-fm-d}). Section~\ref{sec:intro-setting} recalls the sources. We use the
following form, proved in Appendix~\ref{app:moments-guerra}; recall that $\thstar=\zP(\{0\})$.

\begin{proposition}[Fractional moments for convex mixtures]\label{prop:max-fm}
Assume \condH\ and $\xi\in\classC$, and let $\theta\in(0,1)$.
\begin{thmlist}
\item\label{prop:max-fm-a} $\limsup_N N^{-1}\log\E\ZN^\theta\le\Lambda(\theta)$.
\item\label{prop:max-fm-b} If $\theta\le\thstar$, then
 $\liminf_N N^{-1}\log\E\ZN^\theta\ge\Lambda(\theta)=\theta\feq$.
\item\label{prop:max-fm-c} If $\xi$ is even, then $\liminf_N N^{-1}\log\E\ZN^\theta\ge\Lambda(\theta)$.
\item\label{prop:max-fm-d} If $\xi$ is not even, then
 $\liminf_N N^{-1}\log\E\ZN^\theta\ge\Lambda(\theta)$.
\end{thmlist}
In particular \condFM\ holds.
\end{proposition}

Item~\ref{prop:max-fm-a} is Guerra's interpolation bound (Lemma~\ref{lem:moments-guerra}); it uses
only the convexity of $\xi$ on $[-1,1]$. Item~\ref{prop:max-fm-b} holds under \condR\ alone: for
$\theta\le\thstar$ the Parisi measure is admissible in the constrained problem, and
\eqref{eq:max-crude-M} with Proposition~\ref{prop:max-upper} gives
$\ZN^\theta\ge e^{\theta N(\feq-t)}$ with probability tending to one. Item~\ref{prop:max-fm-c} is
\cite[Section~7]{Tal07} for even polynomials, extended to the other even $\xi\in\classC$ by
truncation in the degree (Appendix~\ref{app:moments-guerra}), so for even $\xi$ the proof of
Theorem~\ref{thm:exist} uses no Kac--Rice
argument. For~\ref{prop:max-fm-d}, a non-even $\xi$ convex on $[-1,1]$ has at least two active
degrees, so the Kac--Rice construction (Proposition~\ref{prop:kr-main}) bounds
$\Prob(\Eex{\varepsilon}{f})$ from below by $e^{N(\Sigma(f)-o(1))}$ at every $f\in(\feq,\fone)$; as
$\Eex{\varepsilon}{f}\subset\{\MN\ge N(f-\varepsilon)\}$, Theorem~\ref{thm:max}\ref{thm:max-b}
turns this into the lower half of $\FM{\thf}$, and the levels $\thf$ fill $(\thstar,1)$.
Proposition~\ref{prop:kr-main} uses neither \condFM\ nor
Theorem~\ref{thm:exist}\ref{thm:exist-c}--\ref{thm:exist-e}, so the argument is not circular.

\begin{remark}[Other routes to \condFM]\label{rem:max-routes}
Boursier bounds $\limsup_NN^{-1}\log\E e^{\theta\MN}\le\Lambda(\theta)$ directly, for
$\theta\in(0,1)$ and $\xi$ in her class \cite[Proposition~2.7]{Bou26b}; by
Theorem~\ref{thm:max}\ref{thm:max-b} this is equivalent to
Proposition~\ref{prop:max-fm}\ref{prop:max-fm-a}. For convex $\xi\in\classC$ the lower half of
\condFM\ also has proofs that avoid \cite[Section~7]{Tal07}: with two active degrees the argument
for item~\ref{prop:max-fm-d} applies whether or not $\xi$ is even, and for $\xi(t)=\beta^2t^p$ with
$p$ even Appendix~\ref{app:moments-pure} perturbs $\xi$ by a second degree and transfers $\MN$
back (Proposition~\ref{prop:moments-pure}).
\end{remark}

\begin{proof}[Proof of Theorem~\ref{thm:exist}]
Parts~\ref{thm:exist-a} and~\ref{thm:exist-b} are
Lemmas~\ref{lem:gauss-interior}, \ref{lem:gauss-critical}, \ref{lem:var-basic} and
Proposition~\ref{prop:var-regular}, as discussed above.

\ref{thm:exist-c} If $\FM{\theta}$ holds, Theorem~\ref{thm:max}\ref{thm:max-b} gives
$\lim_NN^{-1}\log\E e^{\theta\MN}=\lim_NN^{-1}\log\E\ZN^\theta=\Lambda(\theta)$. At $\theta=1$ no
assumption is needed. Since $\Var\HN(\sigma)=N\xi(1)$ for $\sigma\in\hcube$, $\E\ZN=e^{N\xi(1)/2}$
exactly, and Theorem~\ref{thm:max}\ref{thm:max-b} gives $\lim_NN^{-1}\log\E e^{\MN}=\xi(1)/2$.
Under \condR, Lemma~\ref{lem:var-dual} gives
$\Lambar(1)=\lim_{\theta\uparrow1}\Lambda(\theta)=P_{\delta_0}=\xi(1)/2$.

\ref{thm:exist-d} Fix $f\in(\feq,\fone)$ and put $\theta=\thf\in(\thstar,1)$. By~\ref{thm:exist-b},
$\Lambda'(\theta)=f$ and $\Sigma(f)=\Lambda(\theta)-\theta f$; by~\ref{thm:exist-c} and \condFM,
$N^{-1}\log\E e^{s\MN}\to\Lambda(s)$ for every $s\in(0,1)$. Let $\varepsilon>0$.

\emph{Upper bound.} Chebyshev's inequality at $\theta$ gives
\[
 \Prob\bigl(|N^{-1}\MN-f|\le\varepsilon\bigr)\le\Prob\bigl(\MN\ge N(f-\varepsilon)\bigr)
 \le e^{-\theta N(f-\varepsilon)}\,\E e^{\theta\MN},
\]
so the $\limsup$ of $N^{-1}\log$ of both probabilities is at most
$\Lambda(\theta)-\theta f+\theta\varepsilon=\Sigma(f)+\theta\varepsilon$.

\emph{Lower bound.} Let $\tilde\Prob=e^{\theta\MN}\,d\Prob/\E e^{\theta\MN}$. Then
\[
 \Prob\bigl(|N^{-1}\MN-f|\le\varepsilon\bigr)
 \ge e^{-\theta N(f+\varepsilon)}\,\E e^{\theta\MN}\,
 \tilde\Prob\bigl(|N^{-1}\MN-f|\le\varepsilon\bigr).
\]
For $0<\eta<\min(\theta,1-\theta)$, Chebyshev's inequality under $\tilde\Prob$ gives
\[
 \tilde\Prob\bigl(\MN>N(f+\varepsilon)\bigr)\le e^{-\eta N(f+\varepsilon)}
 \frac{\E e^{(\theta+\eta)\MN}}{\E e^{\theta\MN}},\qquad
 \tilde\Prob\bigl(\MN<N(f-\varepsilon)\bigr)\le e^{\eta N(f-\varepsilon)}
 \frac{\E e^{(\theta-\eta)\MN}}{\E e^{\theta\MN}} .
\]
After $N^{-1}\log$, the right-hand sides tend to
$\Lambda(\theta\pm\eta)-\Lambda(\theta)\mp\eta f-\eta\varepsilon$, which equals
$-\eta\varepsilon+o(\eta)$ because $\Lambda'(\theta)=f$. Fixing $\eta$ small, we get
$\tilde\Prob(|N^{-1}\MN-f|\le\varepsilon)\to1$, hence
$\liminf_NN^{-1}\log\Prob(|N^{-1}\MN-f|\le\varepsilon)\ge\Lambda(\theta)-\theta(f+\varepsilon)
=\Sigma(f)-\theta\varepsilon$. Letting $\varepsilon\downarrow0$ proves the first claim.

\emph{The existence event.} If $|N^{-1}\MN-f|\le\varepsilon$, then by~\ref{thm:exist-a} every
maximizer $x^*$ and its order parameter $\zeta^*$ satisfy \eqref{eq:event-exist}. Since
$\MN\ge\FTAP(x)$ for every $x$, almost surely
\[
 \bigl\{|N^{-1}\MN-f|\le\varepsilon\bigr\}\subset\Eex{\varepsilon}{f}
 \subset\bigl\{\MN\ge N(f-\varepsilon)\bigr\}.
\]
The lower bound applies to the smaller event and the upper bound to the larger one, so the same
limits hold for $\Eex{\varepsilon}{f}$, for its inner and for its outer probability.

\ref{thm:exist-e} This is Proposition~\ref{prop:max-fm}.
\end{proof}

\paragraph{Existence is as hard as the fractional-moment formula.} Assume \condR\ and
$\xi\in\classC$. If $\liminf_NN^{-1}\log\Prob(\Eex{\varepsilon}{f})\ge\Sigma(f)-o_\varepsilon(1)$ as
$\varepsilon\downarrow0$, for every $f\in(\feq,\fone)$, then the bound
$\E e^{\thf\MN}\ge e^{\thf N(f-\varepsilon)}\Prob(\Eex{\varepsilon}{f})$ and
Theorem~\ref{thm:max}\ref{thm:max-b} give $\FMge{\thf}$, and $\thf$ ranges over $(\thstar,1)$;
for $\theta\le\thstar$ the lower half holds under \condR\ alone
(Proposition~\ref{prop:max-fm}\ref{prop:max-fm-b}).
Conversely, given the upper half of \condFM, which holds for convex $\xi$
(Proposition~\ref{prop:max-fm}\ref{prop:max-fm-a}), its lower half gives the existence lower bound
through Theorem~\ref{thm:exist}\ref{thm:exist-d}. So, given Theorem~\ref{thm:max} and the upper
half, existence at the free-energy rate at every level and the lower half of \condFM\ are
interchangeable. The first direction is how Proposition~\ref{prop:max-fm}\ref{prop:max-fm-d} is
proved; for even $\xi$ it makes the lower bound of Theorem~\ref{thm:exist}\ref{thm:exist-d}
equivalent to the lower half of the formula of \cite[Section~7]{Tal07}.

The consequences for the free energy itself are collected in Remark~\ref{rem:max-subsumed}
(Appendix~\ref{app:var-moments}).

%% file: figures/fig2-ramsey.tex
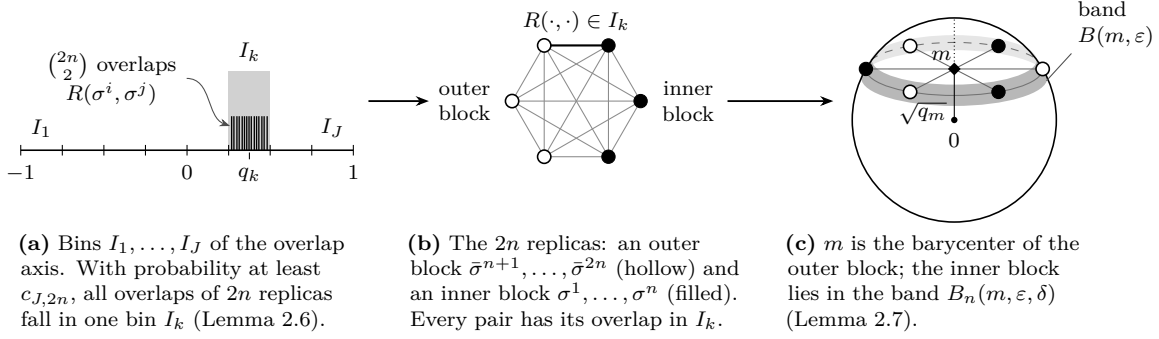
\begin{figure}[tp]
\centering
\begin{NoHyper}
\begin{tikzpicture}[
  font=\scriptsize,
  repout/.style={circle,draw=black,fill=white,line width=0.6pt,inner sep=0pt,minimum size=5pt},
  repin/.style={circle,draw=black,fill=black,inner sep=0pt,minimum size=5pt},
  pairline/.style={draw=black!45,line width=0.4pt},
  flow/.style={-{Stealth[length=5pt]},line width=0.7pt},
  subcap/.style={text width=4.4cm,align=left,anchor=north,inner sep=0pt}]

\begin{scope}[shift={(0,-0.35)}]
  \fill[black!18] (2.75,0) rectangle (3.30,1.05);
  \draw[line width=0.6pt] (0,0) -- (4.4,0);
  \foreach \j in {0,...,8} \draw (0.55*\j,-0.07) -- (0.55*\j,0.07);
  \foreach \x in {2.79,2.82,2.86,2.89,2.93,2.96,2.99,3.02,3.05,3.09,3.12,3.15,3.19,3.22,3.26}
    \draw[line width=0.5pt] (\x,0) -- (\x,0.45);
  \node[above] at (3.025,1.05) {$I_k$};
  \node[above] at (0.275,0.04) {$I_1$};
  \node[above] at (4.125,0.04) {$I_J$};
  \node[below] at (0,-0.07) {$-1$};
  \node[below] at (2.2,-0.07) {$0$};
  \node[below] at (4.4,-0.07) {$1$};
  \draw (3.025,-0.03) -- (3.025,-0.12);
  \node[below] at (3.025,-0.10) {$q_k$};
  \node[align=center] (lab) at (1.2,0.95) {$\binom{2n}{2}$ overlaps\\ $\ov{\sigma^i}{\sigma^j}$};
  \draw[-{Stealth[length=4pt]},black!70] (lab.east) to[out=0,in=160] (2.77,0.32);
\end{scope}
\node[subcap] at (2.2,-1.45) {\textbf{(a)} Bins $I_1,\dots,I_J$ of the overlap axis. With
  probability at least $c_{J,2n}$, all overlaps of $2n$ replicas fall in one bin $I_k$
  (Lemma~\ref{lem:max-ramsey}).};

\draw[flow] (4.6,0.3) -- (5.4,0.3);

\begin{scope}[shift={(7.35,0.3)}]
  \foreach \a in {0,...,4}{%
    \pgfmathtruncatemacro{\s}{\a+1}%
    \foreach \b in {\s,...,5}{\draw[pairline] ({60*\a}:0.85) -- ({60*\b}:0.85);}}
  \draw[line width=0.8pt] (120:0.85) -- (60:0.85)
    node[midway,above=1pt] {$\ov{\cdot}{\cdot}\in I_k$};
  \foreach \a in {120,180,240} \node[repout] at (\a:0.85) {};
  \foreach \a in {0,60,300} \node[repin] at (\a:0.85) {};
  \node[align=center,anchor=east] at (-1.0,0) {outer\\ block};
  \node[align=center,anchor=west] at (1.0,0) {inner\\ block};
\end{scope}
\node[subcap] at (7.35,-1.45) {\textbf{(b)} The $2n$ replicas: an outer block
  $\bar\sigma^{n+1},\dots,\bar\sigma^{2n}$ (hollow) and an inner block $\sigma^1,\dots,\sigma^n$
  (filled). Every pair has its overlap in $I_k$.};

\draw[flow] (9.35,0.3) -- (10.75,0.3);

\begin{scope}[shift={(12.35,0.05)}]
  \fill[black!10] (-1.0911,0.795) arc[start angle=180,end angle=0,x radius=1.0911,y radius=0.3273]
    -- (1.2307,0.555) arc[start angle=0,end angle=180,x radius=1.2307,y radius=0.3692] -- cycle;
  \fill[black!28] (-1.2307,0.555) arc[start angle=180,end angle=360,x radius=1.2307,y radius=0.3692]
    -- (1.0911,0.795) arc[start angle=360,end angle=180,x radius=1.0911,y radius=0.3273] -- cycle;
  \draw[black!55,dashed] (1.1691,0.675) arc[start angle=0,end angle=180,x radius=1.1691,y radius=0.3507];
  \draw[black!55] (-1.1691,0.675) arc[start angle=180,end angle=360,x radius=1.1691,y radius=0.3507];
  \draw[line width=0.6pt] (0,0) circle[radius=1.35];
  \draw[densely dotted] (0,0.675) -- (0,1.35);
  \draw (0,0) -- node[left=-1pt,pos=0.14] {$\sqrt{q_m}$} (0,0.675);
  \fill (0,0) circle[radius=1.2pt] node[below=1pt] {$0$};
  \foreach \p in {0,60,...,300}
    \draw[black!60] (0,0.675) -- ({1.1691*cos(\p)},{0.675+0.3507*sin(\p)});
  \foreach \p in {0,120,240} \node[repout] at ({1.1691*cos(\p)},{0.675+0.3507*sin(\p)}) {};
  \foreach \p in {60,180,300} \node[repin] at ({1.1691*cos(\p)},{0.675+0.3507*sin(\p)}) {};
  \node[diamond,fill=black,inner sep=0pt,minimum size=5pt] at (0,0.675) {};
  \node at (-0.17,0.86) {$m$};
  \draw[black!70] (1.22,0.52) -- (1.55,0.88);
  \node[above right=-2pt,align=left] at (1.55,0.88) {band\\ $B(m,\varepsilon)$};
\end{scope}
\node[subcap] at (12.35,-1.45) {\textbf{(c)} $m$ is the barycenter of the outer block; the
  inner block lies in the band $\Band$ (Lemma~\ref{lem:max-band}).};
\end{tikzpicture}
\end{NoHyper}
\caption{The lower band (Proposition~\ref{prop:max-lower}), drawn with $2n=6$ replicas.
(a)~Whatever the Gibbs measure $G_N$ and whatever $N$, Ramsey's theorem forces some bin $I_k$ of
the overlap axis to catch all pairwise overlaps of $2n$ independent replicas with probability at
least $c_{J,2n}$, which depends only on the number $J$ of bins and on $n$.
(b)~Averaging over the outer block fixes one outer block whose completions to a $2n$-tuple
with all pairwise overlaps near $q_k$ still carry probability at least $c_{J,2n}$.
(c)~Its barycenter $m$ lies in the cube, with $q_m$ close to the bin midpoint $q_k$. After
normalization the configurations lie on the unit sphere. Since
$\ov{\sigma}{m}-q_m=\ov{\sigma-m}{m}$, the band $B(m,\varepsilon)$ is a belt of configurations
whose deviation $\sigma-m$ is nearly orthogonal to $m$, and the pair condition of $\Band$ makes
the deviations of the $n$ replicas nearly orthogonal to each other. Every completing inner block
lies in $\Band$, so the band free energy at $m$ falls short of $\log\ZN$ by at most
$n^{-1}\log(1/c_{J,2n})$. Given $c,t>0$, the band theorem at every exponential rate
(Lemma~\ref{lem:max-truncation}), with $n$ chosen in terms of $c$ and $t$, then gives
$\FTAP(m)\ge\log\ZN-tN$ outside an event of probability at most $e^{-cN}$, for all large $N$.}
\label{fig:ramsey}
\end{figure}

%% file: sections/sec3-localization.tex
\section[Localization on the contact set and the field law:
\texorpdfstring{Theorem~\ref*{thm:struct}\ref*{thm:struct-i}--\ref*{thm:struct-ii}}{Theorem 3(i)-(ii)}]
{Localization on the contact set and the field law:\\
Theorem~\ref*{thm:struct}\ref*{thm:struct-i}--\ref*{thm:struct-ii}}
\label{sec:loc}

This section proves parts~\ref{thm:struct-i} and~\ref{thm:struct-ii} of
Theorem~\ref{thm:struct}. For convex $\xi\in\classC$ the TAP maximum reaches a level
$f\in(\feq,\fone)$ with probability $e^{N\Sigma(f)}$ to leading exponential order
(Theorem~\ref{thm:exist}). Outside an event of probability $e^{N(\Sigma(f)-c)}$, exponentially
smaller, every $m\in\cube$ with $\FTAP(m)\ge N(f-\delta)$, critical or not, has its self-overlap
$q_m$ near the contact set $\Cf$ and the empirical law $\mu_m$ of its coordinates near the field
law $\Lq{q_m}$. The proof uses no existence result, no Kac--Rice formula and no criticality, and
the class $\classC$ plays no role. Throughout, $\xi$ satisfies \condH\ (odd degrees are
allowed), $f\in(\feq,\fone)$ is fixed, and $\zf$ is any minimizer of $\zeta\mapsto P_\zeta$ over
$\{\zeta\in\Pspace([0,1]):\zeta(\{0\})\ge\thf\}$ (it is unique by Lemma~\ref{lem:var-basic}, but
the argument does not use this).

\paragraph{Strategy.}
The proof is an exponential Markov inequality. For a set $S\subset[0,1]$ of radii and
$\theta>0$,
\begin{equation}\label{eq:loc-markov}
 \Pout\bigl(\exists m\in\cube:\ q_m\in S,\ \FTAP(m)\ge N(f-\delta)\bigr)
 \le e^{-N\theta(f-\delta)}\,\E\exp\Bigl(\theta\sup_{q_m\in S}\FTAP(m)\Bigr),
\end{equation}
where the expectation is understood through a measurable majorant of the supremum
(Section~\ref{sec:loc-proof}). Over the whole cube, the comparison of Boursier
\cite[Proposition~2.7]{Bou26b} bounds this expectation by $e^{N(\Lambda(\theta)+o(1))}$ for every
$\theta\in(0,1)$. At $\theta=\thf$ one has $\Lambda(\thf)=\thf P_{\zf}$ and
$\Lambda(\thf)-\thf f=\Sigma(f)$, so~\eqref{eq:loc-markov} returns the rate $\Sigma(f)$ itself
and nothing better. Theorem~\ref{thm:struct}\ref{thm:struct-i} needs more: for every thin
shell $S=[Q,Q+\rho]$ away from $\Cf$, a bound $e^{N\thf(P_{\zf}-\gamma)+o(N)}$ in which
$\gamma>0$ depends on the distance from $\Cf$ but not otherwise on $Q$. Two ingredients provide
it. The \emph{glued comparison} (Proposition~\ref{prop:glued-comparison}) bounds the expectation
in~\eqref{eq:loc-markov} for a shell through the Parisi functional of a coefficient that may
exceed one below $Q$ and jump down at $Q$. The \emph{relaxed variational gap}
(Proposition~\ref{prop:loc-gap}) shows that off $\Cf$ some such coefficient beats $P_{\zf}$ by an
amount that depends only on the height of the obstacle at $Q$ above its minimum; for $Q>0$ a
nondecreasing one gains nothing (Lemma~\ref{lem:loc-monotone}), so the gain comes from the jump
(at $Q=0$ the comparison class has no atom constraint, and its best value is $\feq$).
Markov's inequality and a union bound over finitely many shells then give
Theorem~\ref{thm:struct}\ref{thm:struct-i} (Section~\ref{sec:loc-proof}). The field law
\ref{thm:struct-ii} comes from the same comparison with the terminal condition at the radius
perturbed by a test function of the coordinates (Section~\ref{sec:loc-field}).
Proposition~\ref{prop:glued-comparison} is proved in Appendix~\ref{app:glued}, the remaining
direction of Proposition~\ref{prop:loc-gap} in Appendix~\ref{app:var-gap}, and the calculus of
the Parisi equation with nonnegative, non-monotone coefficients on which both rest in
Appendix~\ref{app:pde-toolbox}.

\subsection{The obstacle and its contact set}\label{sec:loc-obstacle}

\paragraph{Coefficients and the Parisi functional.}
The comparisons below involve Parisi functionals of coefficients that need not be
distribution functions. A \emph{coefficient} is a measurable function $c:[0,1)\to[0,3]$. For a
coefficient $c$ let $\Phi_c$ solve
\[
 \partial_t\Phi_c+\frac{\xi''(t)}2\Bigl(\partial_x^2\Phi_c+c(t)\,(\partial_x\Phi_c)^2\Bigr)=0
 \quad\text{on }[0,1)\times\R,\qquad \Phi_c(1,\cdot)=\log\cosh,
\]
in the sense of Appendix~\ref{app:pde-toolbox}: by the Cole--Hopf formula on each step when
$c$ is a step function, and by continuity in $L^1([0,1))$ in general
(Lemma~\ref{lem:pde-measurable}). Let $X^c$ be the associated Auffinger--Chen process,
$dX_t=\xi''(t)c(t)\,\partial_x\Phi_c(t,X_t)\,dt+\sqrt{\xi''(t)}\,dB_t$ with $X_0=0$, and set
\begin{gather*}
 P(c)=\Phi_c(0,0)-\frac12\int_0^1t\xi''(t)c(t)\,dt,\\
 \Gz{c}(t)=\E\,\partial_x\Phi_c(t,X^c_t)^2,\qquad
 \Obs_c(t)=\frac12\int_t^1\xi''(s)\bigl(\Gz{c}(s)-s\bigr)\,ds .
\end{gather*}
For $c=\cdf{\zeta}$ with $\zeta\in\Pspace([0,1])$ these are the Parisi functional $P_\zeta$,
the function $\Gz{\zeta}$ and the obstacle $\Obs_\zeta$ of
Sections~\ref{sec:intro-setting}--\ref{sec:intro-structure} (Lemma~\ref{lem:pde-measurable}).
Following \cite{Bou26b}, we call $\Obs_\zeta$ the Parisi obstacle of $\zeta$.

Two facts from Appendix~\ref{app:pde-toolbox} carry all the variational arguments of this
section (Corollary~\ref{cor:pde-step} and Lemma~\ref{lem:pde-measurable}). Write
\[
 DP(c)[h]=\frac12\int_0^1\xi''(t)h(t)\bigl(\Gz{c}(t)-t\bigr)\,dt=-\int_0^1h(t)\,\Obs_c'(t)\,dt
\]
for the first variation of $P$ at $c$ in the direction of a bounded measurable $h$. For all
coefficients $c,c'$,
\begin{equation}\label{eq:loc-first-order}
 P(c')-P(c)\ge DP(c)[c'-c],
\end{equation}
and there is a constant $K<\infty$, depending only on $\xi'(1)$ and $\xi''(1)$, such that
\begin{equation}\label{eq:loc-second-order}
 P(c+\varepsilon h)\le P(c)+\varepsilon\,DP(c)[h]+K\varepsilon^2
 \qquad\text{whenever } \|h\|_{L^1}\le1 \text{ and } c+\varepsilon h \text{ is a coefficient.}
\end{equation}
The constant $K$ depends neither on $c$ nor on the direction $h$; this is what will make the
gap below uniform in the radius. Since $\Obs_\zeta\in C^1([0,1])$ and $\Obs_\zeta(1)=0$,
Fubini's theorem turns the first variation towards another distribution function into an
integral of the obstacle:
\begin{equation}\label{eq:loc-moving-mass}
 DP(\cdf{\zeta})\bigl[\cdf{\mu}-\cdf{\zeta}\bigr]=\int\Obs_\zeta\,d(\mu-\zeta),
 \qquad \mu,\zeta\in\Pspace([0,1]).
\end{equation}
So $\Obs_\zeta(t)$ is the first-order cost of moving mass of $\zeta$ to the point $t$. For the
constrained minimizer, mass may be moved only within $\{\mu(\{0\})\ge\thf\}$, and this
determines the shape of the obstacle.

\begin{lemma}[The constrained minimizer and its obstacle]\label{lem:loc-obstacle}
Assume \condH\ and let $f\in(\feq,\fone)$.
\begin{thmlist}[beginpenalty=10000]
\item\label{lem:loc-obstacle-a} $\Lambar$ is convex and lower semicontinuous on $[0,1]$, and
 $\Lambar'_+(0)=\feq$. The function $\theta\mapsto\Lambar(\theta)-\theta f$ has exactly one
 minimizer $\thf$ on $[0,1]$; it lies in $(0,1)$, and $\Sigma(f)=\Lambda(\thf)-\thf f$.
\item\label{lem:loc-obstacle-b} The infimum of $P_\zeta$ over
 $\{\zeta\in\Pspace([0,1]):\zeta(\{0\})\ge\thf\}$ is attained. Every minimizer $\zf$ satisfies
 $\zf(\{0\})=\thf$ and $\Sigma(f)=\thf(P_{\zf}-f)$.
\item\label{lem:loc-obstacle-c} Let $\zf$ be such a minimizer and $\Obs=\Obs_{\zf}$. Then
 $\Obs\in C^1([0,1])$, $\Obs'(t)=-\frac12\xi''(t)(\Gz{\zf}(t)-t)$ and $\Obs(1)=0$. The minimum
 $\Omin$ of $\Obs$ over $(0,1]$ is attained, $\Obs=\Omin$ on $\supp\zf\cap(0,1]$, and
 \[
  \Obs(0)>\Omin,\qquad \Obs(1)>\Omin .
 \]
 In particular $\Cf=\argmin_{(0,1]}\Obs$ is a nonempty compact subset of $(0,1)$ that
 contains $\supp\zf\cap(0,1]$.
\end{thmlist}
\end{lemma}

Except for the uniqueness of $\thf$, the facts in
\ref{lem:loc-obstacle-a}--\ref{lem:loc-obstacle-c} appear in \cite[Lemma~3.2]{Bou26b} and its
proof, under Boursier's standing assumption on $\xi$, which is \condH\
(Section~\ref{sec:intro-setting}). Appendix~\ref{app:var-constrained} gives proofs that use only
the calculus of Appendix~\ref{app:pde-toolbox}; there $\thf$ is shown to be unique by a mixture
argument, so that under \condH\ every level in $(\feq,\fone)$ is regular in the sense of
\cite[Definition~1.1]{Bou26b}.

Part~\ref{lem:loc-obstacle-c} says that $\zf$ puts its mass in $(0,1]$ only where moving mass
is cheapest, while its atom of mass $\thf$ at the origin is held in place by the constraint. The
strict inequality $\Obs(0)>\Omin$ is where $f>\feq$ enters: if it failed, moving a little more
mass to the origin would cost nothing at first order, and $\thf$ would not minimize
$\Lambar(\theta)-\theta f$. At a contact $Q\in\Cf$, an interior minimum of $\Obs$, $\Obs'(Q)=0$
and $\xi''(Q)>0$ give $\Gz{\zf}(Q)=\int t^2\,\Lq{Q}(dt)=Q$: the second moment of the field law
equals the radius.
Theorem~\ref{thm:struct}\ref{thm:struct-i}--\ref{thm:struct-ii} say that near-maximizers at
level $f$ realize this picture. Besides $\supp\zf\cap(0,1]$, the contact set may contain points
that carry no mass. Figure~\ref{fig:obstacle} sketches the obstacle, its contact set and the
glued coefficients of Section~\ref{sec:loc-gap}.

\input{figures/fig3-obstacle}

\subsection{The glued comparison}\label{sec:loc-glued}

Fix a radius $Q\in[0,1)$. For $\eta\in\Pspace([Q,1])$ let $\Phi_\eta$ be the solution on
$[Q,1]\times\R$ of the Parisi equation with coefficient $\cdf{\eta}(t)=\eta([Q,t])$ and terminal
value $\log\cosh$ (Appendix~\ref{app:pde-toolbox}), and let
$\Phistar_\eta(Q,m)=\inf_b\{\Phi_\eta(Q,b)-mb\}$, extended continuously to $m=\pm1$. If
$q_m=Q$, then $\eta$ is admissible in the infimum that defines $\FTAP(m)$, so
\begin{equation}\label{eq:loc-band}
 \FTAP(m)\le F_\eta(m)=\HN(m)+\sum_{i=1}^N\Phistar_\eta(Q,m_i)
 -\frac N2\int_Q^1t\xi''(t)\cdf{\eta}(t)\,dt .
\end{equation}
On a sphere, the band measure $\eta$ thus enters only through a deterministic function of the
coordinates of $m$. The randomness of $\HN$ is compared with a cascade, as in Guerra's
interpolation, and the comparison involves only the overlaps up to $Q$; the parameters of the
cascade are not tied to $\eta$. The next proposition makes this precise. It also allows a shell
of radii $[Q,Q+\rho]$, which a union bound over radii needs, and a perturbation of the terminal
condition at $Q$, which the proof of the field law needs.

\begin{proposition}[Glued comparison on a sphere or a shell]\label{prop:glued-comparison}
Assume \condH. Let $\theta\in(0,1)$, $Q\in[0,1)$ and $\rho\in[0,1-Q]$, with $\rho>0$ if
$Q=0$. Let the data be as follows.
\begin{itemize}[leftmargin=*,topsep=2pt,itemsep=1pt]
\item If $Q>0$: a step function $\xc:[0,Q)\to(\theta,3]$ equal to $\xc_\ell$ on
 $[q_{\ell-1},q_\ell)$ for $\ell=1,\dots,L$, where $0=q_0<q_1<\dots<q_L=Q$ and
 $\theta<\xc_1<\dots<\xc_L$; put $y=\xc_L$. If $Q=0$ there is no $\xc$, and $y=1$.
\item A measure $\eta\in\Pspace([Q,1])$ whose distribution function $\cdf{\eta}$ is a step
 function. No relation between $\xc(Q-)$ and $\cdf{\eta}(Q)$ is assumed.
\item A function $g_0\in C^2([-1,1])$ and $s\in\R$ with $|s|\,\|g_0''\|_\infty\le\frac12$.
\end{itemize}
Let $\Psi_s(z)=\sup_{|m|\le1}\{zm+\Phistar_\eta(Q,m)+s\,g_0(m)\}$ for $z\in\R$; let $W_s$ be
the solution on $[0,Q]\times\R$ of
$\partial_tW+\frac{\xi''}2\bigl(\partial_x^2W+\xc\,(\partial_xW)^2\bigr)=0$ with
$W(Q,\cdot)=\Psi_s$ (and $W_s(0,0)=\Psi_s(0)$ if $Q=0$); and let
\[
 P_s(\xc,\eta)=W_s(0,0)-\frac12\int_0^Qt\xi''(t)\xc(t)\,dt
 -\frac12\int_Q^1t\xi''(t)\cdf{\eta}(t)\,dt .
\]
Then for every $N\ge1$ there is a random variable $\mathcal T$ such that
\begin{equation}\label{eq:loc-glued-pointwise}
 \FTAP(m)+s\sum_{i=1}^Ng_0(m_i)\le N P_s(\xc,\eta)+\mathcal T
 \qquad\text{for every } m\in\cube \text{ with } q_m\in[Q,Q+\rho],
\end{equation}
and
\begin{equation}\label{eq:loc-glued-moment}
 \E e^{\theta\mathcal T}\le
 \exp\Bigl\{\theta\Bigl(NC_R\rho+\ind_{\rho>0}\,\frac{\log(N+1)}{y}\Bigr)\Bigr\},
 \qquad C_R=5\xi''(1).
\end{equation}
\end{proposition}

When $s=0$, one has $\Psi_0=\Phi_\eta(Q,\cdot)$ (Lemma~\ref{lem:pde-fenchel}). Then $W_0$,
continued by $\Phi_\eta$ on $[Q,1]$, is the solution on $[0,1]$ of the Parisi equation with the
\emph{glued coefficient}
\begin{equation}\label{eq:loc-glue}
 \glue(\xc,\eta)(t)=
 \begin{cases}\xc(t), & t\in[0,Q),\\ \cdf{\eta}(t), & t\in[Q,1),\end{cases}
\end{equation}
and $P_0(\xc,\eta)=P(\glue(\xc,\eta))$. Below $Q$ the glued coefficient may exceed one, and at
$Q$ it may jump down. On a sphere ($\rho=0<Q$) the right side of~\eqref{eq:loc-glued-moment} is
$1$ for every $N$.

\paragraph{General data.}
Let $\mathcal A_Q$, for $Q\in(0,1)$, be the set of pairs $(\xc,\eta)$ with $\xc:[0,Q)\to[\theta,2]$
nondecreasing and $\eta\in\Pspace([Q,1])$, and let $\mathcal A_0=\Pspace([0,1])$; at $Q=0$ the
data consist of $\eta$ alone. The quantity $P_s(\xc,\eta)$ is Lipschitz in the data for the
$L^1$ distance of $\xc$ and of $\cdf\eta$, uniformly in $(s,g_0)$, so it extends to
$\mathcal A_Q$, with $P_0(\xc,\eta)=P(\glue(\xc,\eta))$. Approximating the data by step data that
satisfy the hypotheses of Proposition~\ref{prop:glued-comparison}
(Lemma~\ref{lem:pde-approx} and Corollary~\ref{cor:pde-majorant}), one obtains: for every
$(\xc,\eta)\in\mathcal A_Q$, every admissible $(s,g_0)$ and every $\varepsilon_A>0$ there is a
random variable $\mathcal T$ with
\begin{equation}\label{eq:loc-glued-general}
 \sup_{q_m\in[Q,Q+\rho]}\Bigl[\FTAP(m)+s\sum_{i=1}^Ng_0(m_i)\Bigr]
 \le NP_s(\xc,\eta)+N\varepsilon_A+\mathcal T,
 \qquad \E e^{\theta\mathcal T}\le(N+1)\,e^{N\theta C_R\rho},
\end{equation}
uniformly in $Q$ and in the data. For $\zeta\in\Pspace([0,1])$ with $\zeta(\{0\})\ge\theta$, the
pair $(\cdf{\zeta}|_{[0,Q)},\Shift_Q\zeta)$, where
$\Shift_Q\zeta=\zeta|_{(Q,1]}+\zeta([0,Q])\delta_Q$ has distribution function $\cdf{\zeta}$ on
$[Q,1]$, lies in $\mathcal A_Q$ and has glued coefficient $\cdf{\zeta}$. With these data,
\eqref{eq:loc-glued-general} and a union over shells recover Boursier's bound
$\limsup_N N^{-1}\log\E\exp(\theta\MN)\le\theta P_\zeta$ \cite[Proposition~2.7]{Bou26b}.

\paragraph{Relation to Boursier's comparison.}
Proposition~\ref{prop:glued-comparison} is the comparison of Boursier
\cite[Proposition~2.7]{Bou26b}, run on a sphere or a thin shell. Her proof uses neither
monotonicity of the coefficient across the radius nor the bound by one, so the coefficient may
exceed one below the radius and may jump down at it. Appendix~\ref{app:glued-proof} gives the
proof, and Appendix~\ref{app:glued-map} maps each of its steps onto \cite[Section~2]{Bou26b}.
The cascade carries an outer Poisson level of parameter $\theta$, the device of \cite{CGKLM}, as
in \cite[p.~13]{Bou26b}. The part below $Q$ is a Guerra-type bound for a
supremum, and parameters above one are natural there. Parameters larger than one also occur in
Talagrand's reverse bound for $\E\ZN^a$ with $a>1$, whose parameters lie between $1$ and $a$
\cite[Proposition~6.2]{Tal07}, and non-monotone effective coefficients occur in the Guerra-type
bounds of \cite[Section~6]{HS25} and in the extended variational principle of \cite{EAMS21}. The
band $[Q,1]$ enters only through the TAP correction $\Phistar_\eta(Q,\cdot)$, whose coefficient
$\cdf{\eta}$ is at most one.

\subsection{The relaxed variational gap}\label{sec:loc-gap}

From now on $\theta=\thf$ in the definition of $\mathcal A_Q$, $\zf$ is a constrained minimizer
and $\Obs=\Obs_{\zf}$. In the limit of thin shells, the best bound that
\eqref{eq:loc-glued-general} gives at the radius $Q$ with $s=0$ is
\[
 p_Q=\inf\bigl\{P(\glue(\xc,\eta)):(\xc,\eta)\in\mathcal A_Q\bigr\},\qquad Q\in[0,1).
\]
The next proposition compares $p_Q$ with $P_{\zf}=\Lambda(\thf)/\thf$, the value that the
whole cube gives. Put
\begin{equation}\label{eq:loc-gap-function}
 \varphi(d)=\min\Bigl\{\frac{d^2}{4K},\ \frac{\thf d}2\Bigr\},\qquad d\ge0,
\end{equation}
with $K$ the constant of~\eqref{eq:loc-second-order}. Note that $\varphi$ is increasing.

\begin{proposition}[Relaxed variational gap]\label{prop:loc-gap}
Assume \condH\ and let $f\in(\feq,\fone)$.
\begin{thmlist}
\item\label{prop:loc-gap-a} $p_Q\le P_{\zf}$ for every $Q\in[0,1)$.
\item\label{prop:loc-gap-b} $p_Q=P_{\zf}$ for every $Q\in\Cf$.
\item\label{prop:loc-gap-c} For every $Q\in[0,1)\setminus\Cf$ there are data
 $(\xc,\eta)\in\mathcal A_Q$ with
 \[
  P(\glue(\xc,\eta))\le P_{\zf}-\varphi\bigl(\Obs(Q)-\Omin\bigr).
 \]
 In particular $p_Q\le P_{\zf}-\varphi(\Obs(Q)-\Omin)$, a gap that depends on $Q$ only through
 $\Obs(Q)$.
\end{thmlist}
\end{proposition}

By Lemma~\ref{lem:loc-obstacle}\ref{lem:loc-obstacle-c}, $\Obs(Q)-\Omin>0$ for every
$Q\in[0,1)\setminus\Cf$, including $Q=0$, so the gap in~\ref{prop:loc-gap-c} is positive. Part
\ref{prop:loc-gap-b} is a consistency check that the proof of Theorem~\ref{thm:struct} does not
use: at a contact no glued coefficient beats $P_{\zf}$, so the relaxation never contradicts the
existence of states at contacts.

\begin{proof}[Proof of Proposition~\ref{prop:loc-gap}]
We prove~\ref{prop:loc-gap-a}, \ref{prop:loc-gap-b}, and~\ref{prop:loc-gap-c} in full when
$\Cf$ meets $[0,Q)$; the other case of~\ref{prop:loc-gap-c} is the same computation, sketched
here and written out in Appendix~\ref{app:var-gap}.

\ref{prop:loc-gap-a} Since $\cdf{\zf}\ge\zf(\{0\})=\thf$, the pair
$(\cdf{\zf}|_{[0,Q)},\Shift_Q\zf)$ lies in $\mathcal A_Q$, and its glued coefficient is
$\cdf{\zf}$.

\ref{prop:loc-gap-c}, \emph{direction (A)}. Suppose that $\Cf$ contains a point $r<Q$. Let
$h=\ind_{[r,Q)}$ and, for $0<\varepsilon\le\thf$, take $\xc=\cdf{\zf}+\varepsilon h$ on
$[0,Q)$ and $\eta=\Shift_Q\zf$. Then $\xc$ is nondecreasing with values in
$[\thf,1+\thf]\subset[\thf,2]$, so $(\xc,\eta)\in\mathcal A_Q$, and
$\glue(\xc,\eta)=\cdf{\zf}+\varepsilon h$. The first variation is
\[
 DP(\cdf{\zf})[h]=-\int_r^Q\Obs'(t)\,dt=\Obs(r)-\Obs(Q)=-d,\qquad d:=\Obs(Q)-\Omin .
\]
As $\|h\|_{L^1}\le1$, \eqref{eq:loc-second-order} gives
$P(\glue(\xc,\eta))\le P_{\zf}-\varepsilon d+K\varepsilon^2$. If $d\le2K\thf$, the choice
$\varepsilon=d/(2K)$ gives $-\varepsilon d+K\varepsilon^2=-d^2/(4K)$. Otherwise take
$\varepsilon=\thf$; then $K\thf<d/2$ and $-\varepsilon d+K\varepsilon^2\le-\thf d/2$. In both
cases the gain is at least $\varphi(d)$.

\emph{Direction (B).} If $\Cf$ has no point below $Q$ (this includes $Q=0$ and every
$Q<\min\Cf$), it has a point $r\in(Q,1)$ by Lemma~\ref{lem:loc-obstacle}\ref{lem:loc-obstacle-c},
and one lowers the band distribution function by $\varepsilon$ on $[Q,r)$ instead: keep
$\xc=\cdf{\zf}|_{[0,Q)}$ and take $\eta=\Shift_Q\zf-\varepsilon\delta_Q+\varepsilon\delta_r$, a
probability measure for $\varepsilon\le\thf$ because $\cdf{\zf}(Q)\ge\zf(\{0\})=\thf$. The
first variation is again $\Obs(r)-\Obs(Q)$, and the same choice of $\varepsilon$ applies;
Appendix~\ref{app:var-gap} gives the details.

\ref{prop:loc-gap-b} Let $Q\in\Cf$ and $(\xc,\eta)\in\mathcal A_Q$. After a change of $\xc$ on
a countable set, which does not change $P$, $\xc(t)=\alpha([0,t])$ for a finite measure $\alpha$
on $[0,Q)$ with $\alpha(\{0\})\ge\thf$. By Fubini and $\Obs(1)=0$,
\[
 \int_0^1\glue(\xc,\eta)\,(-\Obs')\,dt
 =\int_{[0,Q)}\bigl(\Obs(t)-\Obs(Q)\bigr)\,\alpha(dt)+\int_{[Q,1]}\Obs\,d\eta ,
\]
and the same computation with Lemma~\ref{lem:loc-obstacle} gives
$\int_0^1\cdf{\zf}(-\Obs')\,dt=\int\Obs\,d\zf=\thf\Obs(0)+(1-\thf)\Omin$. Since
$\Obs(Q)=\Omin\le\Obs$ on $(0,1]$, $\Obs(0)>\Omin$ and $\alpha(\{0\})\ge\thf$, the first integral
is at least $\thf(\Obs(0)-\Omin)$, and the second is at least $\Omin$. Hence
\[
 DP(\cdf{\zf})\bigl[\glue(\xc,\eta)-\cdf{\zf}\bigr]
 \ge\thf\bigl(\Obs(0)-\Omin\bigr)+\Omin-\thf\Obs(0)-(1-\thf)\Omin=0,
\]
and~\eqref{eq:loc-first-order} gives $P(\glue(\xc,\eta))\ge P_{\zf}$. With~\ref{prop:loc-gap-a},
$p_Q=P_{\zf}$.
\end{proof}

Why must the glued coefficient jump when $Q>0$? Because otherwise it is the distribution
function of an admissible order parameter.

\begin{lemma}[A monotone glued coefficient gains nothing]\label{lem:loc-monotone}
Assume \condH\ and let $Q\in(0,1)$. If $(\xc,\eta)\in\mathcal A_Q$ and $\glue(\xc,\eta)$ is
nondecreasing on $[0,1)$, then $P(\glue(\xc,\eta))\ge P_{\zf}$.
\end{lemma}

\begin{proof}
Monotonicity across $Q$ forces $\xc\le\cdf{\eta}(Q)\le1$. So the glued coefficient agrees, off a
countable set, with a right-continuous nondecreasing function with values in $[\thf,1]$, that is,
with the distribution function of some $\mu\in\Pspace([0,1])$; and
$\mu(\{0\})=\xc(0+)\ge\thf$. Since $P$ depends on a coefficient only through its class in $L^1$,
$P(\glue(\xc,\eta))=P_\mu\ge\inf\{P_{\mu'}:\mu'(\{0\})\ge\thf\}=P_{\zf}$.
\end{proof}

\begin{remark}[Where the gain comes from]\label{rem:loc-jump}
By Lemma~\ref{lem:loc-monotone}, for $Q>0$ every gain in
Proposition~\ref{prop:loc-gap}\ref{prop:loc-gap-c} comes from a downward jump
$\xc(Q-)>\cdf{\eta}(Q)$ of the glued coefficient. Both directions create such a jump by moving
mass across $Q$ towards a point of $\Cf$. In direction (A) the coefficient is raised by
$\varepsilon$ on $[r,Q)$ and then drops by $\varepsilon$ at $Q$, since $\zf$ has no atom at
$Q\notin\Cf$ (Figure~\ref{fig:obstacle}); in direction (B) the band distribution function is
lowered on $[Q,r)$. The size of the gain is governed by $\Obs(Q)-\Omin$, which compares the
obstacle at $Q$ with its global minimum wherever that minimum is attained. This is why the
radii of near-maximizers localize on the set of global minimizers of $\Obs$, and not merely
near its stationary points (compare Remark~\ref{rem:loc-stationary}). Boursier's band comparison
already uses a strict obstacle gap of this kind, for competitors whose overlap with a given state
is a fixed value below the first contact \cite[Proposition~3.9 and Lemma~3.10,
(3.15)--(3.16)]{Bou26b}; here the gap
controls the self-overlaps of all near-maximizers, through the jump at the radius.
\end{remark}

\subsection{Proof of \texorpdfstring{Theorem~\ref*{thm:struct}\ref*{thm:struct-i}}{Theorem 3(i)}}
\label{sec:loc-proof}

We prove part~\ref{thm:struct-i} of Theorem~\ref{thm:struct} here and part~\ref{thm:struct-ii}
in Section~\ref{sec:loc-field}. Parts~\ref{thm:struct-iii} and~\ref{thm:struct-iv} are proved in
Section~\ref{sec:hess-proof}, the first of them from~\ref{thm:struct-i} and~\ref{thm:struct-ii}.

\paragraph{Measurability.}
The events in~\ref{thm:struct-i} and~\ref{thm:struct-ii} are defined through suprema over the
cube, and we bound their outer probabilities. In each case the event is contained in a finite
union of measurable events $\{\mathcal T_j\ge a_j\}$, with $\mathcal T_j$ random variables as
in~\eqref{eq:loc-glued-general}, and each of these is bounded by Markov's inequality.

\begin{proof}[Proof of Theorem~\ref{thm:struct}\ref{thm:struct-i}]
Fix $\kappa>0$ and write $\theta=\thf$, $\zeta=\zf$, $\Obs=\Obs_{\zf}$. Let
$\mathcal Q_\kappa=\{Q\in[0,1]:\dist(Q,\Cf)\ge\kappa/2\}$. If $\mathcal Q_\kappa$ is empty, no $m$ has
$\dist(q_m,\Cf)\ge\kappa$ and there is nothing to prove. Otherwise
\[
 \gamma_\kappa=\min_{Q\in \mathcal Q_\kappa}\bigl(\Obs(Q)-\Omin\bigr)>0,
\]
because $\Obs$ is continuous on $[0,1]$, $\Obs>\Omin$ on $(0,1]\setminus\Cf$, and
$\Obs(0)>\Omin$ (Lemma~\ref{lem:loc-obstacle}\ref{lem:loc-obstacle-c}). Put
$\varphi_\kappa=\varphi(\gamma_\kappa)$, choose $\rho\in(0,\kappa/2)$ with $C_R\rho\le\varphi_\kappa/8$, and
set $\delta=\varphi_\kappa/8$.

\emph{Shells.} Let $Q_j=j\rho$ and $I_j=[Q_j,\min(Q_j+\rho,1)]$ for the integers $j\ge0$ with
$Q_j<1$. If $\dist(q_m,\Cf)\ge\kappa$, then $q_m\in I_j$ for some such $j$, and
$\dist(Q_j,\Cf)\ge\kappa-\rho>\kappa/2$, so $Q_j\in \mathcal Q_\kappa$. Call such $j$ \emph{bad}.

\emph{One shell.} Let $j$ be bad. Proposition~\ref{prop:loc-gap}\ref{prop:loc-gap-c} gives data
in $\mathcal A_{Q_j}$ whose glued coefficient has
$P\le P_\zeta-\varphi(\Obs(Q_j)-\Omin)\le P_\zeta-\varphi_\kappa$, since $\varphi$ is
increasing. Then~\eqref{eq:loc-glued-general}, with $s=0$, shell width $\min(\rho,1-Q_j)$ and
$\varepsilon_A=\varphi_\kappa/2$, gives a random variable $\mathcal T_j$ with
\[
 \FTAP(m)\le N\Bigl(P_\zeta-\frac{\varphi_\kappa}2\Bigr)+\mathcal T_j\quad\text{for } q_m\in I_j,
 \qquad \E e^{\theta\mathcal T_j}\le(N+1)\,e^{N\theta C_R\rho}.
\]
Hence
$\{\exists m:q_m\in I_j,\ \FTAP(m)\ge N(f-\delta)\}\subset
\{\mathcal T_j\ge N(f-\delta-P_\zeta+\varphi_\kappa/2)\}$, and Markov's inequality, together with
$\Sigma(f)=\theta(P_\zeta-f)$ (Lemma~\ref{lem:loc-obstacle}\ref{lem:loc-obstacle-b}), bounds the
probability of the latter event by
\[
 (N+1)\exp N\Bigl[\theta(P_\zeta-f)-\frac{\theta\varphi_\kappa}2+\theta\delta+\theta C_R\rho\Bigr]
 \le(N+1)\,e^{N(\Sigma(f)-\theta\varphi_\kappa/4)} .
\]

\emph{Union bound.} There are at most $\lceil1/\rho\rceil+1$ shells. The event in
\ref{thm:struct-i} is therefore contained in a finite union of measurable events of total
probability at most $(\lceil1/\rho\rceil+1)(N+1)e^{N(\Sigma(f)-\theta\varphi_\kappa/4)}$, which is at
most $e^{N(\Sigma(f)-\theta\varphi_\kappa/8)}$ for all large $N$. This proves~\ref{thm:struct-i} with
\[
 \delta=\frac{\varphi(\gamma_\kappa)}8,\qquad c=\frac{\thf\,\varphi(\gamma_\kappa)}8 ,
\]
where $\varphi$ is the function~\eqref{eq:loc-gap-function}, whose constant $K$ depends only on
$\xi'(1)$ and $\xi''(1)$, and $\gamma_\kappa$ measures how far the obstacle rises above its
minimum at distance at least $\kappa/2$ from $\Cf$.
\end{proof}

No existence result enters this proof. When $\xi\in\classC$,
Theorem~\ref{thm:exist}\ref{thm:exist-d}--\ref{thm:exist-e} shows that the TAP maximum reaches
level $f-\delta$ with probability at least $e^{N(\Sigma(f)-c/2)}$ for large $N$, and dividing turns
\ref{thm:struct-i} and~\ref{thm:struct-ii} into statements conditional on the deviation:
Corollary~\ref{cor:title}.

\subsection{The field law:
\texorpdfstring{Theorem~\ref*{thm:struct}\ref*{thm:struct-ii}}{Theorem 3(ii)}}
\label{sec:loc-field}

The field law comes from the same comparison, with the terminal condition at the radius
perturbed by $s$ times a test function of the coordinates. We give the argument in outline; the
net of test functions, the approximation of the data and the choice of constants are in
Appendix~\ref{app:var-field}.

\begin{proof}[Outline of the proof of Theorem~\ref{thm:struct}\ref{thm:struct-ii}]
Fix $\kappa>0$ and let $X$ be the Auffinger--Chen process of $\zf$, so that
$\Lq{q}=\Law(\Phi_{\zf,x}(q,X_q))$.

\emph{Finitely many test functions.} By Kantorovich--Rubinstein duality on $[-1,1]$, a net and a
mollification, there are smooth $g_1,\dots,g_J$ with $\|g_j'\|_\infty\le1$ such that, whenever
$\Wone(\mu_m,\Lq{q_m})\ge\kappa$, some $j$ satisfies
\begin{equation}\label{eq:loc-field-bad}
 \frac1N\sum_{i=1}^Ng_j(m_i)\ge\int g_j\,d\Lq{q_m}+\frac\kappa2 .
\end{equation}
The maps $q\mapsto\int g_j\,d\Lq{q}$ are uniformly continuous on $[0,1]$, since
$q\mapsto\Phi_{\zf,x}(q,X_q)$ is almost surely continuous and bounded
(Lemma~\ref{lem:pde-measurable}); let $\omega$ be a common modulus of continuity.

\emph{The tangent bound.} Fix a radius $Q$ and take the base data
$(\cdf{\zf}|_{[0,Q)},\Shift_Q\zf)$, for which $P_0=P_{\zf}$, with $g_0=g_j$. For these data (through step
approximants, Appendix~\ref{app:var-field}) the maximizer in the definition of $\Psi_0(z)$ is $m=\Phi_{\zf,x}(Q,z)$, so
$\Psi_s(z)=\Psi_0(z)+s\,g_j(\Phi_{\zf,x}(Q,z))+O(s^2)$ uniformly in $z$. The terminal sandwich
(Lemma~\ref{lem:pde-terminal}) and the stability of the Auffinger--Chen process under a change of
the terminal condition (Lemma~\ref{lem:pde-stability}) carry this expansion to the initial
value:
\begin{equation}\label{eq:loc-tangent}
 P_s\bigl(\cdf{\zf}|_{[0,Q)},\Shift_Q\zf\bigr)
 \le P_{\zf}+s\int g_j\,d\Lq{Q}+C_2s^2,
 \qquad 0\le s\le\bigl(2\|g_j''\|_\infty\bigr)^{-1},
\end{equation}
with $C_2$ depending only on $\xi'(1)$ and $\xi''(1)$, uniformly in $Q$ and $j$. So, to first
order in $s$, the perturbation $s\sum_ig_j(m_i)$ costs the bound no more than $s$ times the mean
of $g_j$ under the field law at the radius; a near-maximizer whose coordinates give $g_j$ a larger
mean gains more than the bound allows, and Markov's inequality detects this. The perturbation
follows the device by which \cite[Theorem~5]{CPS2} identifies the energy of near-maximizers on a
sphere whose radius lies in the support of the Parisi measure, there by convexity in the inverse
temperature.

\emph{Markov's inequality.} Let $I=[Q,Q+\rho]$ be a shell. If $q_m\in I$,
$\FTAP(m)\ge N(f-\delta)$ and~\eqref{eq:loc-field-bad} holds, then
\[
 \FTAP(m)+s\sum_{i=1}^Ng_j(m_i)\ge N\Bigl[f-\delta+s\Bigl(\int g_j\,d\Lq{Q}-\omega(\rho)
 +\frac\kappa2\Bigr)\Bigr],
\]
while by~\eqref{eq:loc-glued-general} and~\eqref{eq:loc-tangent} the left side is at most
$N[P_{\zf}+s\int g_j\,d\Lq{Q}+C_2s^2+\varepsilon_A]+\mathcal T$. The terms
$s\int g_j\,d\Lq{Q}$ cancel, and Markov's inequality with $\Sigma(f)=\thf(P_{\zf}-f)$ bounds the
probability of this event by
\[
 (N+1)\exp N\Bigl\{\Sigma(f)+\thf\Bigl[C_2s^2+\varepsilon_A+C_R\rho+\delta+s\,\omega(\rho)
 -\frac{s\kappa}2\Bigr]\Bigr\}.
\]
The gain $s\kappa/2$ is linear in $s$ and the cost $C_2s^2$ is quadratic. Appendix~\ref{app:var-field}
takes $s\le\kappa/(8C_2)$, and then $\rho$, $\delta$ and $\varepsilon_A$ small compared with
$s\kappa$, so that the bracket is at most $-\frac5{32}s\kappa$. A union over the finitely many
shells and test functions gives~\ref{thm:struct-ii}, with $c=\thf s\kappa/8$.
\end{proof}

Taking the smaller of the two values of $\delta$, and half the smaller of the two values of $c$,
gives the joint form: for all large $N$ (the union of the two exceptional events has outer
probability at most $2e^{N(\Sigma(f)-2c)}\le e^{N(\Sigma(f)-c)}$), outside an event of outer
probability at most $e^{N(\Sigma(f)-c)}$, every $m$ with
$\FTAP(m)\ge N(f-\delta)$ satisfies both $\dist(q_m,\Cf)<\kappa$ and
$\Wone(\mu_m,\Lq{q_m})<\kappa$. Section~\ref{sec:hess} starts from this statement.

\begin{remark}[The field law and stationarity]\label{rem:loc-stationary}
Since $t\mapsto t^2$ is $2$-Lipschitz on $[-1,1]$, $q_m=\int t^2\,\mu_m(dt)$ and
$\Gz{\zf}(q)=\int t^2\,\Lq{q}(dt)$, one has $|q_m-\Gz{\zf}(q_m)|\le2\,\Wone(\mu_m,\Lq{q_m})$.
So~\ref{thm:struct-ii} alone places $q_m$ near the zero set of $\Gz{\zf}(t)-t$, which consists of
$0$ and the points of $(0,1)$ where $\Obs'$ vanishes. Our proof excludes the stationary points that
are not global minimizers of $\Obs$ through~\ref{thm:struct-i}, that is, through the downward jump
of the glued coefficient (Remark~\ref{rem:loc-jump}).
\end{remark}

%% file: figures/fig3-obstacle.tex
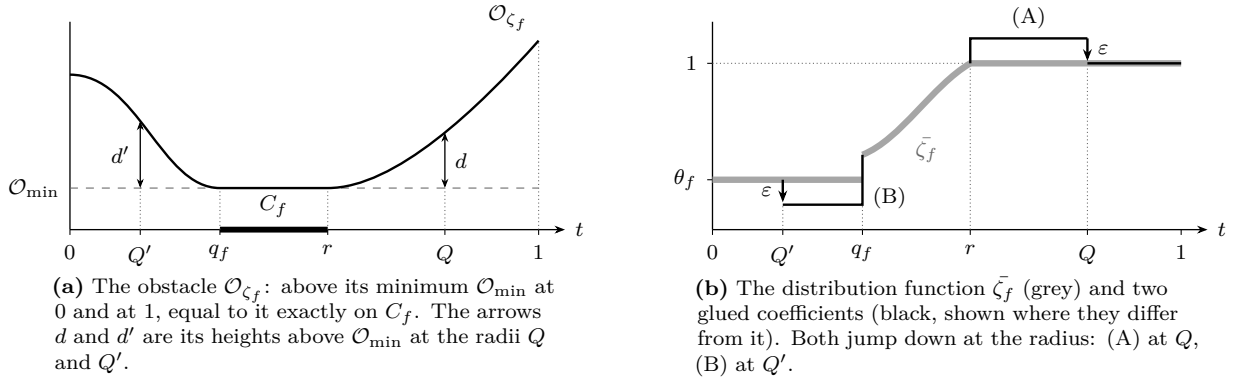
\begin{figure}[bp]
\centering
\begin{NoHyper}
\begin{tikzpicture}[
  font=\scriptsize,
  axis/.style={line width=0.5pt},
  guide/.style={densely dotted,black!60},
  jumparrow/.style={-{Stealth[length=4.5pt]},line width=0.9pt},
  gap/.style={{Stealth[length=3.5pt]}-{Stealth[length=3.5pt]},line width=0.5pt},
  subcap/.style={text width=6.8cm,align=left,anchor=north,inner sep=0pt}]

\begin{scope}
  \draw[axis,-{Stealth[length=4pt]}] (0,0) -- (6.55,0) node[right] {$t$};
  \draw[axis] (0,0) -- (0,2.75);
  \draw[dashed,black!60] (0,0.55) -- (6.2,0.55);
  \node[left] at (0,0.55) {$\Omin$};
  \draw[line width=2.2pt] (1.984,0) -- (3.41,0);
  \node[above] at (2.697,0.03) {$\Cf$};
  \draw[guide] (1.984,0) -- (1.984,0.55);
  \draw[guide] (3.41,0) -- (3.41,0.55);
  \draw[guide] (0.93,0) -- (0.93,0.55);
  \draw[guide] (4.96,0) -- (4.96,0.55);
  \draw[guide] (6.2,0) -- (6.2,2.5);
  \path[name path=obs] (0,2.05) .. controls (0.9,2.05) and (1.2,0.55) .. (1.984,0.55)
    -- (3.41,0.55) .. controls (4.3,0.55) and (5.4,1.6) .. (6.2,2.5);
  \draw[line width=0.9pt] (0,2.05) .. controls (0.9,2.05) and (1.2,0.55) .. (1.984,0.55)
    -- (3.41,0.55) .. controls (4.3,0.55) and (5.4,1.6) .. (6.2,2.5);
  \path[name path=vQ] (4.96,0) -- (4.96,3);
  \path[name path=vQp] (0.93,0) -- (0.93,3);
  \path[name intersections={of=obs and vQ,by=OQ}];
  \path[name intersections={of=obs and vQp,by=OQp}];
  \draw[gap] (4.96,0.55) -- (OQ) node[midway,right] {$d$};
  \draw[gap] (0.93,0.55) -- (OQp) node[midway,left] {$d'$};
  \node[above left] at (6.2,2.5) {$\Obs_{\zf}$};
  \foreach \x/\l in {0/0,0.93/Q',1.984/\qf,3.41/r,4.96/Q,6.2/1}
    \draw (\x,0) -- (\x,-0.07) node[below] {$\l$};
\end{scope}
\node[subcap] at (3.2,-0.62) {\textbf{(a)} The obstacle $\Obs_{\zf}$: above its minimum
  $\Omin$ at $0$ and at $1$, equal to it exactly on $\Cf$. The arrows $d$ and $d'$ are its heights
  above $\Omin$ at the radii $Q$ and $Q'$.};

\begin{scope}[shift={(8.5,0)}]
  \draw[axis,-{Stealth[length=4pt]}] (0,0) -- (6.55,0) node[right] {$t$};
  \draw[axis] (0,0) -- (0,2.75);
  \foreach \y/\l in {0.66/\thf,2.2/1}
    \draw (0,\y) -- (-0.07,\y) node[left] {$\l$};
  \draw[guide] (0,2.2) -- (3.41,2.2);
  \draw[guide] (0.93,0) -- (0.93,0.33);
  \draw[guide] (1.984,0) -- (1.984,0.33);
  \draw[guide] (3.41,0) -- (3.41,2.2);
  \draw[guide] (4.96,0) -- (4.96,2.2);
  \draw[black!35,line width=2.4pt] (0,0.66) -- (1.984,0.66);
  \draw[black!35,line width=0.6pt] (1.984,0.66) -- (1.984,0.99);
  \draw[black!35,line width=2.4pt] (1.984,0.99) .. controls (2.5,1.2) and (2.9,1.9) .. (3.41,2.2)
    -- (6.2,2.2);
  \node[black!60,below right] at (2.55,1.35) {$\cdf{\zf}$};
  \draw[line width=0.9pt] (3.41,2.2) -- (3.41,2.53) -- (4.96,2.53);
  \draw[jumparrow] (4.96,2.53) -- (4.96,2.23) node[midway,right] {$\varepsilon$};
  \draw[line width=0.9pt] (4.96,2.2) -- (6.2,2.2);
  \node[above] at (4.185,2.53) {(A)};
  \draw[jumparrow] (0.93,0.66) -- (0.93,0.36) node[midway,left] {$\varepsilon$};
  \draw[line width=0.9pt] (0.93,0.33) -- (1.984,0.33) -- (1.984,0.99);
  \node[right] at (1.984,0.42) {(B)};
  \foreach \x/\l in {0/0,0.93/Q',1.984/\qf,3.41/r,4.96/Q,6.2/1}
    \draw (\x,0) -- (\x,-0.07) node[below] {$\l$};
\end{scope}
\node[subcap] at (11.7,-0.62) {\textbf{(b)} The distribution function $\cdf{\zf}$ (grey) and two
  glued coefficients (black, shown where they differ from it). Both jump down at the radius:
  (A)~at $Q$, (B)~at $Q'$.};
\end{tikzpicture}
\end{NoHyper}
\caption{The obstacle and the glued coefficient: a schematic, not to scale, for a convex $\xi$
with a hypothetical contact set $\Cf=[\qf,r]$; in general $\Cf$ is only known to be a nonempty
compact subset of $(0,1)$ (Lemma~\ref{lem:loc-obstacle}). The constrained minimizer $\zf$ has an
atom $\thf$ at $0$, puts the rest of its mass on $\Cf$, and may have an atom at $\qf$. To exclude
a radius off $\Cf$, the glued comparison (Proposition~\ref{prop:glued-comparison}) is applied to a coefficient
that moves mass across the radius toward a contact. \emph{Direction (A)}, for a radius $Q$ above
a contact $r$: raise the coefficient by $\varepsilon$ on $[r,Q)$. Here $r=\max\Cf$, so below $Q$ the
cascade coefficient equals $1+\varepsilon>1$, and at $Q$ it drops back to the band
distribution function.
\emph{Direction (B)}, for a radius $Q'$ below every contact: lower the band distribution function
by $\varepsilon$ on $[Q',\qf)$. Either move changes the Parisi functional by
$-\varepsilon d+O(\varepsilon^2)$, where $d$ is the height of the obstacle above its minimum at
the radius, so the bound for the shell at the radius lies at least $\varphi(d)>0$ below
$P_{\zf}$ (Proposition~\ref{prop:loc-gap}). For $Q>0$ a nondecreasing glued coefficient is the
distribution function of an admissible order parameter and gains nothing
(Lemma~\ref{lem:loc-monotone}): the gain comes from the downward jump at the radius.}
\label{fig:obstacle}
\end{figure}

%% file: sections/sec4-hessian.tex
\section[The Hessian at rare TAP states: \texorpdfstring{Theorem~\ref*{thm:struct}\ref*{thm:struct-iii}--\ref*{thm:struct-iv} and Corollary~\ref*{cor:title}}{Theorem 3(iii)-(iv) and Corollary 4}]
  {The Hessian at rare TAP states:\\ Theorem~\ref*{thm:struct}\ref*{thm:struct-iii}--\ref*{thm:struct-iv} and Corollary~\ref*{cor:title}}\label{sec:hess}

By Theorem~\ref{thm:struct}\ref{thm:struct-i}--\ref{thm:struct-ii}, at a level
$f\in(\feq,\fone)$ and outside an event of outer probability exponentially smaller than
$e^{N\Sigma(f)}$, every $m$ with $\FTAP(m)\ge N(f-\delta)$ has $q_m$ near the contact set $\Cf$ and
$\mu_m$ near the field law $\Lq{q_m}$. We read off from these two facts the spectrum of the Hessian
at such states, which proves Theorem~\ref{thm:struct}\ref{thm:struct-iii}--\ref{thm:struct-iv}, and
combine it with Theorems~\ref{thm:max} and~\ref{thm:exist} to prove Corollary~\ref{cor:title}.

\paragraph{Standing assumptions and notation.}
Throughout this section \condH\ holds, $\xi\in\classC$ and $f\in(\feq,\fone)$. Single-degree
covariances, such as the Sherrington--Kirkpatrick covariance $\xi(t)=\beta^2t^2$, are included.
The constrained minimizer $\zf$ is then unique (Lemma~\ref{lem:var-basic}), and $\Cf$,
$\Lq{q}$, $\UQ$, $\chiQ$, $\Plef$, $\rhoQ$, $\nuQ$ and $\lamp$ are defined from $\zf$ as in
Section~\ref{sec:intro-structure}. We write $X$ for the Auffinger--Chen process of $\zf$ and
$M_t=\Phi_{\zf,x}(t,X_t)$, so that $\Lq{Q}=\Law(M_Q)$, and $\UQ=\Phi_{\zf,xx}(Q,X_Q)$ takes
values in $(0,1]$ (Lemma~\ref{lem:pde-class}\ref{lem:pde-class-a} with $G=1$, extended by
Lemma~\ref{lem:pde-measurable}, and Lemma~\ref{lem:pde-standing}\ref{lem:pde-standing-a}). The
contact set is a nonempty compact subset of $(0,1)$
(Lemma~\ref{lem:loc-obstacle}), so $\kappa_f=\dist(\Cf,\{0,1\})$ is positive. We fix
$r_0\in(0,1/4)$ with $2r_0\le\kappa_f$, so that $\Cf\subset[2r_0,1-2r_0]$. Since $\xi\not\equiv0$
has nonnegative coefficients, $\aQ=\xi''(Q)>0$ for every $Q\in(0,1)$.

\paragraph{The argument.}
It has three steps. \emph{Structure and uniformity}
(Sections~\ref{sec:hess-pastur}--\ref{sec:hess-bulk}): the Hessian of $F_\alpha$ at $x\in\ocube$ is
$\Hess\HN(x)$, minus a diagonal matrix whose $i$-th entry depends only on $\alpha$, $q_x$, $x_i$
and the mean susceptibility $\chi_\alpha(x)$, plus a matrix of rank at most
two~\eqref{eq:hess-structure}, and Proposition~\ref{prop:hess-bulk} puts its spectrum near a
reflected Pastur law at all localized points at once (Lemma~\ref{lem:matrix-uniform} gives the
underlying estimate over the whole cube), outside an event of outer probability $e^{-MN}$. \emph{Edge}
(Section~\ref{sec:hess-edge}): the upper edge $\lamp(Q)$ of $\nuQ$ is the minimum of a strictly
convex function that vanishes at an explicit point with slope $\Plef(Q)$
(Lemma~\ref{lem:hess-edge}, Figure~\ref{fig:edge}); this is
Theorem~\ref{thm:struct}\ref{thm:struct-iv}. \emph{Localization}
(Sections~\ref{sec:hess-own}--\ref{sec:hess-proof}): at a near-maximizer $x$,
Proposition~\ref{prop:hess-own}, built on
Theorem~\ref{thm:struct}\ref{thm:struct-i}--\ref{thm:struct-ii}, brings $(q_x,\mu_x)$ near
$(Q,\Lq{Q})$ for a contact $Q$ and every near-optimal order parameter near $\zf$ on $[q_x,1]$; this
gives Theorem~\ref{thm:struct}\ref{thm:struct-iii}. Only this step uses the rare event. In the proof
of Theorem~\ref{thm:struct}\ref{thm:struct-iii}--\ref{thm:struct-iv}, convexity of $\xi$ enters only
through Theorem~\ref{thm:struct}\ref{thm:struct-i}--\ref{thm:struct-ii}; the
class $\classC$ enters through the regularity of the Parisi equation near the contact set
(Lemma~\ref{lem:pde-standing}) and through \cite[Theorem~12(iii)]{CPS}, used in
Proposition~\ref{prop:hess-own}.

\subsection{Reflected Pastur laws}\label{sec:hess-pastur}

Let $\GOE$ be a GOE matrix, with off-diagonal variance $1/N$ and diagonal variance $2/N$; $\ESD$
and the ordering of eigenvalues are as in Section~\ref{sec:intro-structure}. We record the
characterization of reflected Pastur laws and the properties used below.

\begin{lemma}[Reflected Pastur laws]\label{lem:hess-pastur}
Let $a\ge0$ and $\rho\in\Pspace([0,\infty))$.
\begin{thmlist}
\item\label{lem:hess-pastur-a} There is a unique probability measure $\nu=\Pastur{a}{\rho}$ on
 $\R$ whose Cauchy transform $g(w)=\int\nu(d\lambda)/(w-\lambda)$ satisfies
 \begin{equation}\label{eq:hess-pastur}
  g(w)=\int\frac{\rho(dD)}{w+D-a\,g(w)}\qquad\text{for every } w\in\Cplx^+ .
 \end{equation}
 For each $w\in\Cplx^+$, equation~\eqref{eq:hess-pastur} has at most one solution $g(w)$ with
 $\operatorname{Im}g(w)<0$. If $D_N$ are real diagonal matrices with $\ESD(D_N)\to\rho$ weakly,
 then $\ESD(\sqrt a\,\GOE-D_N)\to\nu$ weakly in probability, and the expected empirical spectral
 distribution converges weakly to $\nu$.
\item\label{lem:hess-pastur-b} $\nu((2\sqrt a,\infty))=0$, and
 $\nu((-\infty,-2\sqrt a-s))\le\rho([s,\infty))$ for every $s\ge0$.
\item\label{lem:hess-pastur-c} If $a_n\to a$ and $\rho_n\to\rho$ weakly, then
 $\Pastur{a_n}{\rho_n}\to\Pastur{a}{\rho}$ weakly. Hence, for $A,B>0$, the map
 $(a,\rho)\mapsto\Pastur{a}{\rho}$ is uniformly continuous from $[0,A]\times\Pspace([0,B])$,
 with the metric $|a-a'|+\dL$, to $(\Pspace(\R),\dL)$.
\item\label{lem:hess-pastur-d} Let $D,D'\ge0$ be random variables on one probability space, with
 laws $\rho$ and $\rho'$. If $|D-D'|\le\varepsilon$ almost surely, then
 $\dL(\Pastur{a}{\rho},\Pastur{a}{\rho'})\le\varepsilon$. If $\Prob(D\ne D')\le p$, then
 $\dL(\Pastur{a}{\rho},\Pastur{a}{\rho'})\le p$. Moreover
 $\dL(\Pastur{a}{\rho},\Pastur{a'}{\rho})\le2|\sqrt a-\sqrt{a'}|$.
\end{thmlist}
\end{lemma}

In the language of free probability, $\Pastur{a}{\rho}$ is the free additive convolution of the
semicircle law of variance $a$ with the law of $-D$, $D\sim\rho$; we use only~\eqref{eq:hess-pastur}.
The proof (Appendix~\ref{app:matrix-pastur}) is the standard resolvent argument, with Weyl's
inequalities, interlacing and tightness for~\ref{lem:hess-pastur-b}--\ref{lem:hess-pastur-d}.

\paragraph{Two elementary facts.}
Recall that $\dL(\mu,\nu)$ is the infimum of the $\varepsilon>0$ such that
$\nu((-\infty,t-\varepsilon])-\varepsilon\le\mu((-\infty,t])\le\nu((-\infty,t+\varepsilon])+\varepsilon$
for every $t\in\R$. By right continuity of distribution functions, $\dL(\mu,\nu)\le\eta$ implies,
for every $t\in\R$,
\begin{equation}\label{eq:hess-levy}
 \nu((t+\eta,\infty))-\eta\ \le\ \mu((t,\infty))\ \le\ \nu([t-\eta,\infty))+\eta .
\end{equation}
For a symmetric $N\times N$ matrix $H$ and $x\ne0$, the \emph{compression} $H'$ of $H$ to
$x^\perp$ is the restriction of $PHP$ to $x^\perp$, where $P$ is the orthogonal projection onto
$x^\perp$; its empirical spectral distribution is normalized by $N-1$. Cauchy's interlacing
theorem gives $\lambda_{i+1}(H)\le\lambda_i(H')\le\lambda_i(H)$ for $i\le N-1$. So, for every
$t\in\R$, the number of eigenvalues of $H'$ larger than $t$ lies between that of $H$ minus one and
that of $H$, and comparing distribution functions,
\begin{equation}\label{eq:hess-interlace}
 \dL\bigl(\ESD(H'),\ESD(H)\bigr)\le\sup_{t\in\R}\bigl|\ESD(H')((-\infty,t])-\ESD(H)((-\infty,t])\bigr|
 \le\frac1{N-1}.
\end{equation}

\subsection{The Hessian bulk at localized points}\label{sec:hess-bulk}

\paragraph{Susceptibilities and the Hessian.}
For $\alpha\in\Pspace([0,1])$ and $t\in[0,1]$, the map $b\mapsto\Phi_{\alpha,x}(t,b)$ is an
increasing bijection of $\R$ onto $(-1,1)$. Let $b_\alpha(t,\cdot)$ be its inverse, and put
\begin{equation}\label{eq:hess-u}
 u_\alpha(t,m)=\Phi_{\alpha,xx}\bigl(t,b_\alpha(t,m)\bigr)\quad(|m|<1),\qquad
 u_\alpha(t,\pm1)=0 .
\end{equation}
Thus $u_\alpha(t,m)\in(0,1]$ for $|m|<1$ is the susceptibility $\partial m/\partial b$ of a
coordinate with magnetization $m$ in the effective field $b$. We write $u=u_{\zf}$. Since
$M_Q=\Phi_{\zf,x}(Q,X_Q)$, the definition gives $u(Q,M_Q)=\UQ$, so $u(Q,\cdot)$ maps $\Lq{Q}$ to
the law of $\UQ$.

Let $x\in\ocube$ with $r=q_x\in(0,1)$, and write $u_i=u_\alpha(r,x_i)$ and
$\chi_\alpha(x)=N^{-1}\sum_iu_i$. The self-overlap enters $F_\alpha$ in two places, the time
argument of $\Phistar_\alpha$ and the lower limit of the integral. Differentiating, each produces
a term containing $\cdf{\alpha}(q_x)$, and the two terms cancel in the gradient:
$\nabla F_\alpha(x)=\nabla\HN(x)-(b_\alpha(r,x_i))_{i\le N}-\xi''(r)\chi_\alpha(x)\,x$. The second
derivative still depends on $\cdf{\alpha}(q_x)$. So if $\alpha$ has an atom at $r$, the Hessian of
$F_\alpha$ jumps across the sphere $\{q_y=r\}$, and $\nabla^2_+F_\alpha(x)$ and
$\nabla^2_-F_\alpha(x)$ denote the Hessians computed on its outer and inner side; they coincide
when $\alpha(\{r\})=0$.
Lemma~\ref{lem:matrix-hessian} computes them exactly:
\begin{equation}\label{eq:hess-structure}
 \Hpm F_\alpha(x)=\Hess\HN(x)-\diag\bigl(u_1^{-1},\dots,u_N^{-1}\bigr)-\xi''(r)\chi_\alpha(x)\,I
 +R_\pm(x),
\end{equation}
where $R_\pm(x)$ is a symmetric matrix of rank at most two with range in
$\operatorname{span}\{x,v\}$, $v_i=\Phi_{\alpha,xxx}/\Phi_{\alpha,xx}$ evaluated at
$(r,b_\alpha(r,x_i))$, and $R_+(x)-R_-(x)$ is a multiple of $xx^{\mathsf T}$. In particular the
compressions of $\nabla^2_+F_\alpha(x)$ and $\nabla^2_-F_\alpha(x)$ to $x^\perp$ coincide. At a
critical point of $F_\alpha$ this compression is the Riemannian Hessian of the restriction of
$F_\alpha$ to the sphere $\{q_y=q_x\}$, which is why compressions appear in
Theorem~\ref{thm:struct}.

\paragraph{Heuristic.}
Fix $x$. The matrix $\Hess\HN(x)$ has the law of $\sqrt{\xi''(r)}\,\GOE$ plus a Gaussian matrix of
rank at most two (Appendix~\ref{app:matrix}), so up to rank four~\eqref{eq:hess-structure} is
$\sqrt a\,\GOE-D_N$ with $a=\xi''(r)$ and $D_N=\diag(u_i^{-1}+a\chi_\alpha(x))$, and
Lemma~\ref{lem:hess-pastur} predicts the spectrum $\Pastur{a}{\ESD(D_N)}$. If moreover
$q_x\approx Q$, $\mu_x\approx\Lq{Q}$ and $\alpha\approx\zf$, the diagonal entries are distributed
approximately like $\UQ^{-1}+\aQ\chiQ$, whose law is $\rhoQ$, and the prediction is $\nuQ$.

\paragraph{Uniformity.}
The points of interest depend on the disorder; a separation of scales makes the heuristic uniform.
For fixed $x$ and $D_N$, a Lipschitz spectral statistic $N^{-1}\tr\varphi(\Hess\HN(x)-D_N)$ is
$O(N^{-1})$-Lipschitz in the underlying Gaussian variables, so it deviates from its mean by a fixed
amount only with probability $e^{-cN^2}$, while a net of $e^{CN}$ points and diagonal matrices and
a bound on the third derivative of $\HN$ outside probability $e^{-MN}$ cover the cube.
Lemma~\ref{lem:matrix-uniform} gives, for every $K,\eta,M>0$ and all large $N$, outside an event
of outer probability at most $e^{-MN}$,
\begin{equation}\label{eq:hess-uniform}
 \dL\Bigl(\ESD\bigl(\Hess\HN(x)-D_N\bigr),\Pastur{\xi''(q_x)}{\ESD(D_N)}\Bigr)\le\eta
\end{equation}
for every $x\in\cube$ and every diagonal matrix $D_N$ with entries in $[0,K]$. Jekel, Sandhu and
Shi prove free-convolution estimates simultaneously for all diagonal matrices by the same balance
between concentration at scale $N^2$ and an exponential net \cite[Section~4]{JSS}.

\paragraph{The localized bulk.}
For $Q\in(0,1)$ and $\mu\in\Pspace([-1,1])$ with $\mu(\{-1,1\})=0$, let $\rho_{Q,\mu}$ be the law
of $u(Q,m)^{-1}+\xi''(Q)\int u(Q,\cdot)\,d\mu$ under $m\sim\mu$. Since $u(Q,\cdot)$ maps
$\Lq{Q}$ to the law of $\UQ$, we have $\rho_{Q,\Lq{Q}}=\rhoQ$, and therefore
$\Pastur{\aQ}{\rho_{Q,\Lq{Q}}}=\nuQ$.

\begin{proposition}[The Hessian bulk at localized points]\label{prop:hess-bulk}
Let $\mathcal K$ be a compact subset of $[r_0,1-r_0]\times\Pspace([-1,1])$, with the metric
$|Q-Q'|+\Wone(\mu,\mu')$, such that $\mu(\{-1,1\})=0$ for every $(Q,\mu)\in\mathcal K$. For every
$\kappa>0$ there is $w\in(0,r_0/2]$, and for every $M>0$ there is $N_0$, such that for $N\ge N_0$,
outside an event of outer probability at most $e^{-MN}$, the following holds. For all
$(Q,\mu)\in\mathcal K$, $x\in\ocube$ and $\alpha\in\Pspace([0,1])$ with
\[
 |q_x-Q|\le w,\qquad \Wone(\mu_x,\mu)\le w,\qquad
 \int_{q_x}^1\bigl|\cdf{\alpha}(t)-\cdf{\zf}(t)\bigr|\,dt\le w,
\]
both one-sided Hessians satisfy
$\dL\bigl(\ESD(\Hpm F_\alpha(x)),\Pastur{\xi''(Q)}{\rho_{Q,\mu}}\bigr)\le\kappa$.
\end{proposition}

The proof (Appendix~\ref{app:matrix-uniform}) caps the entries $u_i^{-1}$ at a level $K$, at a cost
in L\'evy distance of the proportion of coordinates with $u_i<1/K$, which is small; applies
\eqref{eq:hess-uniform}; replaces $u_\alpha$ by $u$ (Lemma~\ref{lem:matrix-orderstab}); and moves
the diagonal law from $\mu_x$ to $\mu$ by a coupling and
Lemma~\ref{lem:hess-pastur}\ref{lem:hess-pastur-c}--\ref{lem:hess-pastur-d}. It does not use the
rare event.

\subsection{The edge of the bulk}\label{sec:hess-edge}

We now locate the upper edge of $\nuQ$. Fix $a>0$ and a random variable $U$ with values in
$(0,1]$, and set
\[
 u_{\max}=\esssup U,\qquad \chi=\E U,\qquad D=U^{-1}+a\chi,\qquad \rho=\Law(D),\qquad
 \Delta=1-a\E U^2 .
\]
Let $\nu=\Pastur{a}{\rho}$ and $\lamp=\sup\supp\nu$. The \emph{edge function} is
\begin{equation}\label{eq:hess-z}
 z(s)=s+a\,\E\bigl[(s+D)^{-1}\bigr],\qquad s>s_-:=-\essinf D=-(u_{\max}^{-1}+a\chi).
\end{equation}
To see why $z$ governs the edge, let $g$ be the Cauchy transform of $\nu$. To the right of the
support, \eqref{eq:hess-pastur} reads $g(w)=\E[(s+D)^{-1}]$ with $s=w-ag(w)$, that is,
$w=s+ag(w)=z(s)$. So $z$ inverts the subordination map $w\mapsto w-ag(w)$, and the edge of the
support is where this inversion stops: at the minimum of $z$. The special point is $s=-a\chi$,
where the diagonal resolvent $(s+D)^{-1}$ reduces to $U$ itself. We use the hypothesis
\begin{equation}\label{eq:hess-fhyp}
 \E\bigl[(D-\essinf D)^{-2}\bigr]=\infty,
\end{equation}
which says that $D$ has enough mass near its essential infimum; it ensures that the minimum of $z$
is attained inside $(s_-,\infty)$.

\begin{lemma}[The edge of the bulk]\label{lem:hess-edge}
Assume~\eqref{eq:hess-fhyp}.
\begin{thmlist}
\item\label{lem:hess-edge-a} The function $z$ is smooth and strictly convex on $(s_-,\infty)$,
 $z'(s)\to-\infty$ as $s\downarrow s_-$, and $z'(s)\to1$ as $s\to\infty$. Hence $z$ has a unique
 critical point $s^*$, and
 \[
  \lamp=z(s^*)=\min_{s>s_-}z(s).
 \]
\item\label{lem:hess-edge-b} \textup{(The edge identity.)} $z(-a\chi)=0$, $z'(-a\chi)=\Delta$ and
 $z''(-a\chi)=2a\E U^3$. Hence $\lamp\le0$, with equality if and only if $\Delta=0$.
\item\label{lem:hess-edge-c} \textup{(Quadratic vanishing.)} There are $c_0,C>0$ depending only on
 $a$ and $\E U^3$ such that $|\Delta|\le c_0$ implies
 $\bigl|\lamp+\Delta^2/(4a\E U^3)\bigr|\le C|\Delta|^3$.
\item\label{lem:hess-edge-d} \textup{(Square-root edge.)} Put $C_2=z''(s^*)$,
 $c_2=z''(s^*+1)$ and $t_0=z(s^*+1)-\lamp$. There are $0<c\le C$ and $\kappa_0>0$, continuous
 functions of $(a,C_2,c_2,t_0)\in(0,\infty)^4$, such that
 \[
  c\,\kappa^{3/2}\le\nu([\lamp-\kappa,\lamp])\le C\kappa^{3/2}\qquad(0<\kappa\le\kappa_0).
 \]
\item\label{lem:hess-edge-e} \textup{(Dependence on the radius.)} For every $Q\in(0,1)$, the
 variable $U=\UQ$ with $a=\aQ$ satisfies~\eqref{eq:hess-fhyp}, and then $\nu=\nuQ$,
 $\lamp=\lamp(Q)$ and $\Delta=\Plef(Q)$. The critical point $s^*_Q$, the edge $\lamp(Q)$ and the
 constants of~\ref{lem:hess-edge-d} depend continuously on $Q\in(0,1)$. For $Q\in\Cf$,
 $\Plef(Q)=2\Obs_{\zf}''(Q)/\xi''(Q)\ge0$; hence $\lamp(Q)=0$ if and only if
 $\Obs_{\zf}''(Q)=0$, and
 \[
  \lamp(Q)=-\frac{\Obs_{\zf}''(Q)^2}{\xi''(Q)^3\,\E\UQ^3}+O\bigl(\Obs_{\zf}''(Q)^3\bigr),
 \]
 where the implied constant depends only on $\xi''(Q)$ and $\E\UQ^3$.
\end{thmlist}
\end{lemma}

\input{figures/fig4-edge}

Parts~\ref{lem:hess-edge-a}--\ref{lem:hess-edge-b} carry the mechanism, and we prove them here.
Figure~\ref{fig:edge} shows it: in the coordinate $\sigma=s+a\chi$, every edge function passes
through the origin with slope $\Delta$, and the edge is its minimum.

\begin{proof}[Proof of Lemma~\ref{lem:hess-edge}\ref{lem:hess-edge-a}--\ref{lem:hess-edge-b}]
\emph{Shape of $z$.} For $s>s_-$ we have $s+D\ge s-s_->0$, so $z$ is smooth, with
\[
 z'(s)=1-a\E(s+D)^{-2},\qquad z''(s)=2a\E(s+D)^{-3}>0 .
\]
As $s\downarrow s_-$, $(s+D)^{-2}$ increases to $(D-\essinf D)^{-2}$, so $z'(s)\to-\infty$ by
monotone convergence and~\eqref{eq:hess-fhyp}; as $s\to\infty$, $z'(s)\to1$ by dominated
convergence. Hence $z'$ increases from $-\infty$ to $1$ and vanishes at exactly one point $s^*$,
which minimizes $z$. It remains to show that $\lamp=z(s^*)$; we prove the two inequalities
separately.

\emph{Upper bound: $\lamp\le z(s)$ for every $s>s_-$.} Let $D_N$ be the diagonal matrix whose
entries are the quantiles $d_i=\inf\{y:\rho((-\infty,y])\ge i/(N+1)\}$, $i\le N$. Then
$d_i\ge\essinf D$ and $\ESD(D_N)\to\rho$ weakly. Since $-s<\min_id_i$, the Gaussian comparison
inequality of Lemma~\ref{lem:matrix-comparison} gives
\[
 \E\lambda_{\max}\bigl(\sqrt a\,\GOE-D_N\bigr)\le s+\frac aN\sum_{i\le N}\frac1{d_i+s}
 \xrightarrow[N\to\infty]{}s+a\E(s+D)^{-1}=z(s),
\]
because $y\mapsto(y+s)^{-1}$ is bounded and continuous on $[\essinf D,\infty)$. The largest
eigenvalue is a $\sqrt{2a/N}$-Lipschitz function of the independent standard Gaussian variables
that define $\GOE$, so for every $\eta>0$, $\lambda_{\max}(\sqrt a\,\GOE-D_N)\le z(s)+\eta$ with
probability tending to one. By Lemma~\ref{lem:hess-pastur}\ref{lem:hess-pastur-a},
$\dL(\ESD(\sqrt a\,\GOE-D_N),\nu)\le\eta$ with probability tending to one as well. On both
events, the left inequality of~\eqref{eq:hess-levy} with $t=z(s)+\eta$ gives
$\nu((z(s)+2\eta,\infty))\le\eta$. Letting $\eta\downarrow0$ yields $\nu((z(s),\infty))=0$.

\emph{Lower bound: $\lamp\ge z(s^*)$.} Let $g$ be the Cauchy transform of $\nu$. By
Lemma~\ref{lem:hess-pastur}\ref{lem:hess-pastur-b}, $\lamp\le2\sqrt a$. On $(\lamp,\infty)$ the
function $g$ is real analytic and positive, with $g'(w)=-\int\nu(d\lambda)(w-\lambda)^{-2}<0$. Put
$s(w)=w-ag(w)$. Then $s'(w)=1-ag'(w)\ge1$ and $s(w)\to\infty$ as $w\to\infty$. Let $w_0$ be the
infimum of the $w>\lamp$ such that $s>s_-$ on $(w,\infty)$. The function
$G(\sigma)=\E(\sigma+D)^{-1}$ is analytic off $-\supp\rho\subset(-\infty,s_-]$, so
$F(w)=g(w)-G(w-ag(w))$ is analytic on $\Cplx^+$ and on a complex neighbourhood of $(w_0,\infty)$.
For $w\in\Cplx^+$ we have $\operatorname{Im}(w-ag(w))>0$, and~\eqref{eq:hess-pastur} says that
$F(w)=0$. By the identity theorem $F=0$ on $(w_0,\infty)$, that is, $g(w)=G(s(w))$, and therefore
\[
 z(s(w))=s(w)+aG(s(w))=s(w)+ag(w)=w\qquad(w>w_0).
\]
Differentiating, $z'(s(w))\,s'(w)=1$, so $z'(s(w))>0$ and $s(w)>s^*$ for every $w>w_0$. If
$w_0>\lamp$, continuity of $s$ at $w_0$ would give $s(w_0)=s_-<s^*$, contradicting $s(w)>s^*$
for $w\downarrow w_0$. Hence $w_0=\lamp$, and every $w>\lamp$ satisfies $w=z(s(w))\ge z(s^*)$.
Letting $w\downarrow\lamp$ gives $\lamp\ge z(s^*)$.

\emph{The edge identity.} Since $-a\chi-s_-=u_{\max}^{-1}\ge1$, the point $-a\chi$ lies in
$(s_-,\infty)$, and $(-a\chi+D)^{-1}=U$. Hence
\[
 z(-a\chi)=-a\chi+a\E U=0,\qquad z'(-a\chi)=1-a\E U^2=\Delta,\qquad z''(-a\chi)=2a\E U^3 .
\]
Thus $\lamp=\min z\le z(-a\chi)=0$. If $\Delta=0$, then $-a\chi$ is the critical point $s^*$ and
$\lamp=0$. If $\Delta\ne0$, then $s^*\ne-a\chi$, and strict convexity gives
$\lamp=z(s^*)<z(-a\chi)=0$.
\end{proof}

Parts~\ref{lem:hess-edge-c}--\ref{lem:hess-edge-e} are proved in Appendix~\ref{app:matrix-pastur}:
\ref{lem:hess-edge-c} by a Taylor expansion of $z$ at $-a\chi$; \ref{lem:hess-edge-d} from
$-g'(w)=(1-z'(s(w)))/(a\,z'(s(w)))$ for $w>\lamp$, which turns the quadratic minimum of $z$ into a
square-root edge of $\nu$; \ref{lem:hess-edge-e} from Lemma~\ref{lem:matrix-hyp}, which
gives~\eqref{eq:hess-fhyp} for $\UQ$, and Lemma~\ref{lem:pde-standing}, which gives continuity in
$Q$ and the identity $\Plef=2\Obs_{\zf}''/\xi''\ge0$ on $\Cf$ (a contact minimizes the obstacle);
substituting this identity into~\ref{lem:hess-edge-b}--\ref{lem:hess-edge-c} gives the last
assertions.

\begin{remark}[The edge detects whether $\Plef$ vanishes, not its sign]\label{rem:hess-sign}
By Lemma~\ref{lem:hess-edge}\ref{lem:hess-edge-b}, $\lamp<0$ for both signs of $\Delta\ne0$. At a
contact the sign is fixed, $\Plef\ge0$, because a contact minimizes the obstacle. By the formulas
for $\Obs_{\zf}'$ and $\Obs_{\zf}''$ in Lemma~\ref{lem:pde-standing}, the identity
$\Plef=2\Obs_{\zf}''/\xi''$ holds at every critical point of the obstacle in $(0,1)$, so
$\Plef\le0$ at a local maximum of the obstacle, and the edge is again negative there unless
$\Plef=0$. Since $\lamp(Q)\le0$ for every $Q\in(0,1)$, the Hessian bulk never rules out a radius:
near-maximizers at radii away from $\Cf$ are excluded by the localization of
Theorem~\ref{thm:struct}\ref{thm:struct-i}, not by the spectrum.
\end{remark}

\paragraph{Relation to earlier work.}
Lemma~\ref{lem:hess-edge} is compared with earlier bulk-edge results after Theorem~\ref{thm:struct}
(Section~\ref{sec:intro-structure}). At equilibrium the inequality
$\xi''(q)\,\E\,\partial_x^2\Phi_{\zP}(q,X_q)^2\le1$ on the support of the Parisi measure ($X$ its
Auffinger--Chen process) is \cite[Corollary~3.10]{JT17} (see also \cite[Proposition~1]{Che17} for
$\Gz{\zP}(q)=q$); an identity of the form $\Plef=2\Obs''/\xi''$ appears in \cite[(2.13)]{Zho26b},
stated there for 1-RSB measures, with $q$ a local maximum of the function $f_\mu$ used there.

\subsection{Order parameters of near-optimal states}\label{sec:hess-own}

Proposition~\ref{prop:hess-bulk} applies to $F_\alpha$ as soon as $\alpha$ is close to $\zf$ on
$[q_x,1]$, in the $L^1$ sense of distribution functions. For $\alpha=\zf$ there is nothing to
prove. The order parameters that matter at a near-maximizer $x$ are those that nearly attain the
infimum defining $\FTAP(x)$; examples are the order parameter $\zeta^*$ of a maximizer in
Theorem~\ref{thm:exist}\ref{thm:exist-a} and the order parameter of a witness of
$\Eex{\varepsilon}{f}$. Since $\HN(x)$ does not depend on the order parameter, near-optimality is a
statement about the functional
\begin{equation}\label{eq:hess-V}
 V_\alpha(r,\mu)=\int\Phistar_\alpha(r,m)\,\mu(dm)
 -\frac12\int_r^1t\xi''(t)\,\cdf{\alpha}(t)\,dt,\qquad r\in[0,1],\ \mu\in\Pspace([-1,1]).
\end{equation}
Indeed $F_\alpha(x)=\HN(x)+NV_\alpha(q_x,\mu_x)$, so that
$F_\alpha(x)-\FTAP(x)=N\bigl[V_\alpha(q_x,\mu_x)-\inf_\beta V_\beta(q_x,\mu_x)\bigr]$.

\begin{proposition}[Order parameters of near-optimal states]\label{prop:hess-own}
For every $\kappa>0$ there are $\varepsilon,c>0$ such that, for all large $N$, outside an event of
outer probability at most $e^{N(\Sigma(f)-c)}$, every $x\in\cube$ and $\alpha\in\Pspace([0,1])$
with
\[
 \FTAP(x)\ge N(f-\varepsilon)\qquad\text{and}\qquad F_\alpha(x)\le\FTAP(x)+2\varepsilon N
\]
satisfy
\[
 \dist(q_x,\Cf)<\kappa,\qquad \Wone(\mu_x,\Lq{q_x})<\kappa,\qquad
 \int_{q_x}^1\bigl|\cdf{\alpha}(t)-\cdf{\zf}(t)\bigr|\,dt<\kappa .
\]
\end{proposition}

In particular the conclusion holds for every witness $(\alpha,x)$ of $\Eex{\varepsilon}{f}$, since
there $F_\alpha(x)\le N(f+\varepsilon)\le\FTAP(x)+2\varepsilon N$.

The idea is that the limiting band problem has a unique minimizer. For $Q\in\Cf$, the map
$\beta\mapsto V_\beta(Q,\Lq{Q})$ has the unique minimizer $\Shift_Q\zf$ over $\Pspace([Q,1])$.
Uniqueness is \cite[Theorem~12(iii)]{CPS} (their normalization differs from ours by the constant
$\log 2$, which does not affect minimizers). Started from the law of $X_Q$, the band process of
$\Shift_Q\zf$ follows the dynamics of $X$ on $[Q,1]$, so its first variation is given by
$\Obs_{\zf}$ on $[Q,1]$, which is minimal at the contact $Q$ and on $\Cf\supset\supp\Shift_Q\zf$;
the first-order condition holds, and by convexity it identifies the minimizer. Compactness
(Lemma~\ref{lem:matrix-own}) makes uniqueness quantitative, uniformly in $Q\in\Cf$, and
Theorem~\ref{thm:struct}\ref{thm:struct-i}--\ref{thm:struct-ii} at a smaller accuracy, with the
$\Wone$-continuity of $Q\mapsto\Lq{Q}$, gives the proposition (Appendix~\ref{app:matrix-own}).
Chen, Panchenko and Subag note that the order parameter of a TAP state follows the Parisi measure
on $[q,1]$ and that this can be quantified \cite[Section~2.4, (2.34)]{CPS};
Proposition~\ref{prop:hess-own} does so at a rare level, with $\zf$ in place of the Parisi measure.

\subsection{\texorpdfstring{Proof of Theorem~\ref*{thm:struct}\ref*{thm:struct-iii}--\ref*{thm:struct-iv} and the index form}{Proof of Theorem 3(iii)-(iv) and the index form}}\label{sec:hess-proof}

For $x\in\ocube$ with $q_x\in(0,1)$ and $\varepsilon>0$, let $\mathcal H_\varepsilon(x)$ be the
set of matrices consisting of $\Hpm F_{\zf}(x)$, of $\Hpm F_\alpha(x)$ for every
$\alpha\in\Pspace([0,1])$ with $F_\alpha(x)\le\FTAP(x)+2\varepsilon N$, and of the compressions of
all of these to $x^\perp$. In this notation, Theorem~\ref{thm:struct}\ref{thm:struct-iii} says:
for every $\kappa\in(0,\kappa_f)$ there are $\varepsilon,c>0$ such that, for all large $N$,
outside an event of outer probability at most $e^{N(\Sigma(f)-c)}$, every $x\in\ocube$ with
$\FTAP(x)\ge N(f-\varepsilon)$ has $\dist(q_x,\Cf)<\kappa$ and
$\dL(\ESD(H),\nu_{q_x})\le\kappa$ for every $H\in\mathcal H_\varepsilon(x)$.

\begin{proof}[Proof of Theorem~\ref{thm:struct}\ref{thm:struct-iii}]
Fix $\kappa\in(0,\kappa_f)$. The set $\mathcal K=\{(Q,\Lq{Q}):Q\in[r_0,1-r_0]\}$ is compact,
because $Q\mapsto\Lq{Q}$ is $\Wone$-continuous (Lemma~\ref{lem:pde-standing}), and
$\Lq{Q}(\{-1,1\})=0$ because $M_Q\in(-1,1)$. Let $w\in(0,r_0/2]$ be given by
Proposition~\ref{prop:hess-bulk} for this $\mathcal K$ and accuracy $\kappa/2$, and use that
proposition with the exponent $M'=1-\Sigma(f)>0$ in place of $M$. Put $\kappa_1=\min(\kappa,w,r_0)$, and let $\varepsilon,c>0$ be
given by Proposition~\ref{prop:hess-own} at accuracy $\kappa_1$. For large $N$ the union of the two
exceptional events has outer probability at most
$e^{N(\Sigma(f)-c)}+e^{N(\Sigma(f)-1)}\le e^{N(\Sigma(f)-c')}$, with $c'=\min(c,1)/2$.

Outside this union, let $x\in\ocube$ with $\FTAP(x)\ge N(f-\varepsilon)$. Since $\FTAP(x)$ is an
infimum over order parameters, some $\alpha$ satisfies $F_\alpha(x)\le\FTAP(x)+2\varepsilon N$, and
Proposition~\ref{prop:hess-own} gives $\dist(q_x,\Cf)<\kappa_1$ and
$\Wone(\mu_x,\Lq{q_x})<\kappa_1$. As $\kappa_1\le r_0$ and $\Cf\subset[2r_0,1-2r_0]$, we get
$q_x\in(r_0,1-r_0)$, so $(q_x,\Lq{q_x})\in\mathcal K$. Apply Proposition~\ref{prop:hess-bulk} with
$(Q,\mu)=(q_x,\Lq{q_x})$. The first two window conditions hold because $\kappa_1\le w$. The third
is void for $\alpha=\zf$, and for $\alpha$ with $F_\alpha(x)\le\FTAP(x)+2\varepsilon N$ it holds by
Proposition~\ref{prop:hess-own}, again because $\kappa_1\le w$. Since
$\Pastur{\xi''(q_x)}{\rho_{q_x,\Lq{q_x}}}=\nu_{q_x}$, every one-sided Hessian in
$\mathcal H_\varepsilon(x)$ has $\ESD$ within L\'evy distance $\kappa/2$ of $\nu_{q_x}$. By
\eqref{eq:hess-interlace}, the compressions are within $\kappa/2+1/(N-1)\le\kappa$ for large $N$.
Finally $\dist(q_x,\Cf)<\kappa_1\le\kappa$.
\end{proof}

\begin{proof}[Proof of Theorem~\ref{thm:struct}\ref{thm:struct-iv}]
Let $Q\in(0,1)$. Then $\aQ>0$, the variable $\UQ$ takes values in $(0,1]$, and
Lemma~\ref{lem:hess-edge}\ref{lem:hess-edge-e} identifies $\nuQ$, $\lamp(Q)$ and $\Plef(Q)$ with
the objects of Lemma~\ref{lem:hess-edge} for $U=\UQ$ and $a=\aQ$, for which~\eqref{eq:hess-fhyp}
holds. Lemma~\ref{lem:hess-edge}\ref{lem:hess-edge-b} gives $\lamp(Q)\le0$, with equality if and
only if $\Plef(Q)=0$. On $\Cf$, $\Plef=2\Obs_{\zf}''/\xi''\ge0$
(Lemma~\ref{lem:pde-standing}), so a contact $Q$ has $\lamp(Q)=0$ if and only if
$\Obs_{\zf}''(Q)=0$.
\end{proof}

\paragraph{The index form.}
The L\'evy distance does not see $o(N)$ eigenvalues, so Theorem~\ref{thm:struct}\ref{thm:struct-iii}
says nothing about the largest eigenvalue. It does control the position of the
$\lceil\kappa'N\rceil$-th eigenvalue for $\kappa'$ of order $\kappa^{3/2}$, because the reflected
Pastur law has a square-root edge (Lemma~\ref{lem:hess-edge}\ref{lem:hess-edge-d}). For $r>0$ let
$\Cf^r$ be the closed $r$-neighbourhood of $\Cf$; for $r<\kappa_f$ it is a compact subset of
$(0,1)$. Fix $\kappa_*\in(0,\kappa_f)$, and for $\kappa\in(0,\kappa_*]$ put
\begin{equation}\label{eq:hess-kappaprime}
 \kappa'(\kappa)=\frac12\inf_{Q\in\Cf^\kappa}\nuQ\bigl((\lamp(Q)-\kappa/2,\infty)\bigr).
\end{equation}

\begin{proposition}[Index form of the bulk]\label{prop:hess-index}
\begin{thmlist}
\item\label{prop:hess-index-a} $\kappa'(\kappa)>0$ for every $\kappa\in(0,\kappa_*]$, and there are
 $c_{\mathrm{ix}},C_{\mathrm{ix}},\kappa_{\mathrm{ix}}>0$ such that
 $c_{\mathrm{ix}}\kappa^{3/2}\le\kappa'(\kappa)\le C_{\mathrm{ix}}\kappa^{3/2}$ for
 $\kappa\in(0,\kappa_{\mathrm{ix}}]$.
\item\label{prop:hess-index-b} For every $\kappa\in(0,\kappa_*]$ there are $\varepsilon,c>0$ such
 that, for all large $N$, outside an event of outer probability at most $e^{N(\Sigma(f)-c)}$,
 every $x\in\ocube$ with $\FTAP(x)\ge N(f-\varepsilon)$ has $q_x\in\Cf^\kappa$, and every
 $H\in\mathcal H_\varepsilon(x)$ satisfies, with $\kappa'=\kappa'(\kappa)$,
 \[
  \lambda_{\lceil\kappa'N\rceil}(H)\in\bigl[\lamp(q_x)-\kappa,\lamp(q_x)+\kappa\bigr],\qquad
  \#\bigl\{i:\lambda_i(H)>\lamp(q_x)+\kappa\bigr\}<\kappa'N .
 \]
\end{thmlist}
\end{proposition}

The proof, in Appendix~\ref{app:matrix-index}, combines the square-root edge with
\eqref{eq:hess-levy} and~\eqref{eq:hess-interlace}. In particular the index cannot in general be
taken equal to the window: for the mixtures to which the Kac--Rice construction of
Section~\ref{sec:stab-kr} applies, a single parameter $\kappa$ for both fails at the constructed
states (Remark~\ref{rem:stab-window}).

\begin{remark}[The index window cannot be widened]\label{rem:stab-window}
Assume the hypotheses of Theorem~\ref{thm:stab}\ref{thm:stab-a} and let $f\in(\feq,\fone)$. There is
$\kappa_{\mathrm w}>0$ with the following property. For every $\kappa\in(0,\kappa_{\mathrm w}]$ and every $\delta>0$,
let $\mathcal B_N(\kappa,\delta)$ be the event that some critical point $x\in\ocube$ of $F_{\zf}$ with
$\FTAP(x)\ge N(f-\delta)$ has $\lambda_{\lfloor\kappa N\rfloor}(H)<\lamp(q_x)-\kappa$ for each $H$
among $\Hpm F_{\zf}(x)$ and their compressions to the orthogonal complement of $x$. Then, for
every $c'>0$ and all large $N$, $\mathcal B_N(\kappa,\delta)$ contains an event of probability at
least $e^{N(\Sigma(f)-c')}$. Consequently, for no $\kappa\in(0,\kappa_{\mathrm w}]$ are there $\delta,c>0$
such that, for all large $N$, outside an event of outer probability at most $e^{N(\Sigma(f)-c)}$,
every $m\in\ocube$ with $\FTAP(m)\ge N(f-\delta)$ has
$\lambda_{\lfloor\kappa N\rfloor}(H)\ge\lamp(q_m)-\kappa$ for some $H$ among $\Hpm F_{\zf}(m)$ and
their compressions. The witnesses are the retained points of Section~\ref{sec:stab-kr}. The proof
is in Appendix~\ref{app:matrix-index} (Proposition~\ref{prop:matrix-window}).
\end{remark}

The index form is most useful in terms of a contact rather than of the radius $q_x$.

\begin{corollary}[Bulk marginality at degenerate contacts]\label{cor:hess-marg}
For every $\kappa\in(0,\kappa_*]$ there are $\kappa',\varepsilon,c>0$ such that, for all large $N$,
outside an event of outer probability at most $e^{N(\Sigma(f)-c)}$, every $x\in\ocube$ with
$\FTAP(x)\ge N(f-\varepsilon)$ admits $Q\in\Cf$ such that every $H\in\mathcal H_\varepsilon(x)$
satisfies
\[
 \lambda_{\lceil\kappa'N\rceil}(H)\in\bigl[\lamp(Q)-\kappa,\lamp(Q)+\kappa\bigr],\qquad
 \#\bigl\{i:\lambda_i(H)>\lamp(Q)+\kappa\bigr\}<\kappa'N .
\]
Here $\lamp(Q)=0$ if $\Obs_{\zf}''(Q)=0$; otherwise $\lamp(Q)<0$, and
$\lamp(Q)=-\Obs_{\zf}''(Q)^2/(\xi''(Q)^3\E\UQ^3)+O(\Obs_{\zf}''(Q)^3)$.
\end{corollary}

Its proof (Appendix~\ref{app:matrix-index}) runs that of
Proposition~\ref{prop:hess-index}\ref{prop:hess-index-b} with $\kappa/2$ in place of $\kappa$ and uses the uniform
continuity of $\lamp$ near $\Cf$.

Under~\condHmarg, every contact is degenerate, so $\lamp(Q)=0$ in
Corollary~\ref{cor:hess-marg}: for every matrix in $\mathcal H_\varepsilon(x)$, at least
$\kappa'N$ eigenvalues are at least $-\kappa$ and fewer than $\kappa'N$ exceed $\kappa$, with
$\kappa'$ of order $\kappa^{3/2}$. This is how Section~\ref{sec:stab} uses the corollary.

\begin{remark}[An empirical Plefka condition]\label{rem:hess-plefka}
For $x\in\ocube$ with $q_x\in(0,1)$ let $\Plef_N(x)=1-\xi''(q_x)N^{-1}\sum_iu(q_x,x_i)^2$ be the
Plefka factor computed from the susceptibilities of $x$ itself. For every $\kappa>0$ there are
$\varepsilon,c>0$ such that, for all large $N$, outside an event of outer probability at most
$e^{N(\Sigma(f)-c)}$, every $x\in\ocube$ with $\FTAP(x)\ge N(f-\varepsilon)$ satisfies
$\Plef_N(x)\ge-\kappa$ and $|\Plef_N(x)-\Plef(Q)|\le\kappa$ for some $Q\in\Cf$. Indeed,
$\Plef_N(x)$ and $\Plef(q_x)$ are the integrals of $1-\xi''(q_x)u(q_x,\cdot)^2$ against $\mu_x$ and
$\Lq{q_x}$, because $u(q_x,\cdot)$ maps $\Lq{q_x}$ to the law of $U_{q_x}$; $u^2$ is bounded and
uniformly continuous on $[r_0/2,1-r_0/2]\times[-1,1]$ (Lemma~\ref{lem:pde-standing}), and $\Plef$
is uniformly continuous on $\Cf^{r_0}$. So the joint form of Theorem~\ref{thm:struct}\ref{thm:struct-i}--\ref{thm:struct-ii} (end of Section~\ref{sec:loc-field})
at a small enough accuracy $\kappa_2\le r_0$ gives $Q\in\Cf$ with $|q_x-Q|<\kappa_2$ and
$|\Plef_N(x)-\Plef(q_x)|+|\Plef(q_x)-\Plef(Q)|<\kappa$, by coupling $\mu_x$ with $\Lq{q_x}$; and
$\Plef(Q)\ge0$ gives $\Plef_N(x)\ge-\kappa$. This transfers to every near-maximizer, up to
$\kappa$, the non-strict Plefka inequality that Boursier notes for the first contact
\cite[Remark~1.2]{Bou26b}.
\end{remark}

\subsection{\texorpdfstring{Proof of Corollary~\ref*{cor:title}}{Proof of Corollary 4}}\label{sec:hess-cor}

Corollary~\ref{cor:title} combines the three theorems: Theorem~\ref{thm:struct} controls every
state near level $f$ outside an event of outer probability at most $2e^{N(\Sigma(f)-c_{\mathrm s})}$,
Theorem~\ref{thm:exist} bounds the conditioning event below by $e^{N(\Sigma(f)-c_{\mathrm s}/2)}$, and
Theorem~\ref{thm:max}, at every exponential rate, transfers the conditioning from $\MN$ to
$\log\ZN$. The constants $\delta,c$ depend on $\kappa$ alone; $t$ enters only through $N_0$. Outer
probability is monotone and countably subadditive, and $\Pout(B\cap A)\le\Pout(B)$; in particular
$\Pout(B\mid A)\le\Pout(B)/\Prob(A)$.

\begin{proof}
\emph{Step 0: reduction and the bad events.} Every conclusion of the corollary becomes weaker when
$\kappa$ increases: this is clear for $\dist(q_m,\Cf)<\kappa$, $\Wone(\mu_m,\Lq{q_m})<\kappa$ and
$\dL\le\kappa$, and at most $\kappa N$ eigenvalues are larger than $\kappa$ if at most
$\kappa_2N$ eigenvalues are larger than $\kappa_2\le\kappa$. So we may assume $\kappa<\kappa_f$.
The joint form of Theorem~\ref{thm:struct}\ref{thm:struct-i}--\ref{thm:struct-ii} (end of
Section~\ref{sec:loc-field}) at accuracy $\kappa$ gives $\delta_1,c_1>0$, and
Theorem~\ref{thm:struct}\ref{thm:struct-iii} at accuracy $\kappa$ gives $\varepsilon_2,c_2>0$. Put
$\delta=\min(\delta_1,\varepsilon_2)$ and $c_{\mathrm s}=\min(c_1,c_2)$. Let $B_1$ be the event that some
$m\in\cube$ with $\FTAP(m)\ge N(f-\delta)$ has $\dist(q_m,\Cf)\ge\kappa$ or
$\Wone(\mu_m,\Lq{q_m})\ge\kappa$. Let $B_2$ be the event that some $m\in\ocube$ with
$\FTAP(m)\ge N(f-\delta)$ fails the following: $\dist(q_m,\Cf)<\kappa$, and
$\dL(\ESD(H),\nu_{q_m})\le\kappa$ for every $H\in\mathcal H_\delta(m)$. When
$\dist(q_m,\Cf)<\kappa<\kappa_f$ we have $q_m\in(0,1)$, so every object in this condition is
defined. Since $\delta\le\delta_1$, the event $B_1$ is contained in the exceptional event of that
joint form. Since
$\delta\le\varepsilon_2$ and $\mathcal H_\delta(m)\subset\mathcal H_{\varepsilon_2}(m)$, the event
$B_2$ is contained in the exceptional event of Theorem~\ref{thm:struct}\ref{thm:struct-iii}.
Hence, for $N\ge N_1(\kappa)$,
\[
 \Pout(B_1\cup B_2)\le2e^{N(\Sigma(f)-c_{\mathrm s})} .
\]
Outside $B_2$, every $H\in\mathcal H_\delta(m)$ also has at most $\kappa N$ eigenvalues larger than
$\kappa$. Indeed, $\nu_{q_m}$ is carried by $(-\infty,\lamp(q_m)]\subset(-\infty,0]$ by
Theorem~\ref{thm:struct}\ref{thm:struct-iv}, so the right inequality of~\eqref{eq:hess-levy} gives
$\ESD(H)((\kappa+s,\infty))\le\nu_{q_m}([s,\infty))+\kappa=\kappa$ for every $s>0$; letting
$s\downarrow0$, $\ESD(H)((\kappa,\infty))\le\kappa$. (For a compression the normalization is $N-1$,
and $\kappa(N-1)\le\kappa N$.)

\emph{Step 1(a): the conditioning event in (a).} The rate $\Sigma$ is concave, as an infimum of
affine functions of $f$, and finite on $(\feq,\fone)$ by
Theorem~\ref{thm:exist}\ref{thm:exist-b}; hence it is continuous there. Choose
$\eta\in(0,(\fone-f)/2)$ with $\Sigma(f+2\eta)\ge\Sigma(f)-c_{\mathrm s}/8$, and put
$f'=f+2\eta\in(\feq,\fone)$. By Theorem~\ref{thm:exist}\ref{thm:exist-e}, \condFM\ holds, so
Theorem~\ref{thm:exist}\ref{thm:exist-d} applies at $f'$. Since
$\Prob(|N^{-1}\MN-f'|\le\varepsilon)$ is nondecreasing in $\varepsilon$,
\[
 \liminf_{N\to\infty}\frac1N\log\Prob\bigl(|N^{-1}\MN-f'|\le\eta\bigr)\ge\Sigma(f')
 \ge\Sigma(f)-c_{\mathrm s}/8 .
\]
On $\{|N^{-1}\MN-f'|\le\eta\}\cap\{|\MN-\log\ZN|\le\eta N/2\}$ we have
$N^{-1}\log\ZN\ge f'-3\eta/2=f+\eta/2\ge f$. Theorem~\ref{thm:max}\ref{thm:max-a} at $t=\eta/2$
allows every exponential rate, in particular the rate $c_{\mathrm s}-\Sigma(f)>0$:
$\Prob(|\MN-\log\ZN|>\eta N/2)\le e^{N(\Sigma(f)-c_{\mathrm s})}$ for large $N$. Hence, for
$N\ge N_2(\kappa)$,
\begin{equation}\label{eq:hess-cond}
 \Prob\bigl(N^{-1}\log\ZN\ge f\bigr)\ge e^{N(\Sigma(f)-c_{\mathrm s}/4)}-e^{N(\Sigma(f)-c_{\mathrm s})}
 \ge e^{N(\Sigma(f)-c_{\mathrm s}/2)} .
\end{equation}
In particular this probability is positive.

\emph{Step 1(b): the conditioning event in (b).} By Theorem~\ref{thm:exist}\ref{thm:exist-d} at
$f$ and the same monotonicity,
\[
 \liminf_{N\to\infty}\frac1N\log\Prob\bigl(\MN\ge N(f-\delta)\bigr)
 \ge\liminf_{N\to\infty}\frac1N\log\Prob\bigl(|N^{-1}\MN-f|\le\delta\bigr)\ge\Sigma(f),
\]
so $\Prob(\MN\ge N(f-\delta))\ge e^{N(\Sigma(f)-c_{\mathrm s}/2)}$ for $N\ge N_2'(\kappa)$.

\emph{Step 2: conditional bounds.} Put $c=c_{\mathrm s}/4$ and $A=\{N^{-1}\log\ZN\ge f\}$. By Step 0
and~\eqref{eq:hess-cond},
\[
 \Pout(B_1\cup B_2\mid A)\le\frac{2e^{N(\Sigma(f)-c_{\mathrm s})}}{e^{N(\Sigma(f)-c_{\mathrm s}/2)}}=2e^{-c_{\mathrm s}N/2}.
\]
Given $t>0$, Theorem~\ref{thm:max}\ref{thm:max-a} gives
$\Prob(|\MN-\log\ZN|>tN)\le e^{N(\Sigma(f)-c_{\mathrm s})}$ for $N\ge N_3(t,\kappa)$, so its conditional
probability given $A$ is at most $e^{-c_{\mathrm s}N/2}$. The conclusions of
Theorem~\ref{thm:exist}\ref{thm:exist-a} fail only on a null set. Hence, for
$N\ge N_0(t,\kappa)=\max(N_1,N_2,N_3,\lceil4\log3/c_{\mathrm s}\rceil)$, the conditional outer probability
that one of items~\ref{cor:title-1}--\ref{cor:title-3} fails is at most
$3e^{-c_{\mathrm s}N/2}\le e^{-cN}$. Indeed, outside the three exceptional sets, item~\ref{cor:title-1}
holds by Theorems~\ref{thm:max}\ref{thm:max-a} and~\ref{thm:exist}\ref{thm:exist-a},
item~\ref{cor:title-2} holds off $B_1$, and item~\ref{cor:title-3} holds off $B_2$, with the
eigenvalue count of Step~0.

For the last assertion of (a), let $t\le\delta$. On the good event, every maximizer $x^*$ lies in
$\ocube$ and satisfies $F_{\zeta^*}(x^*)=\FTAP(x^*)=\MN\ge\log\ZN-tN\ge N(f-\delta)$. So
items~\ref{cor:title-2} and~\ref{cor:title-3} apply to $x^*$, and $\alpha=\zeta^*$ is admissible in
item~\ref{cor:title-3} because $F_{\zeta^*}(x^*)=\FTAP(x^*)$.

For (b), Step~1(b) bounds the probability of the conditioning event from below by
$e^{N(\Sigma(f)-c_{\mathrm s}/2)}$. So the conditional outer probability of $B_1\cup B_2$ is at most
\[
 2e^{-c_{\mathrm s}N/2}\le e^{-cN}\qquad\text{for } N\ge\max\bigl(N_1,N_2',\lceil4\log2/c_{\mathrm s}\rceil\bigr),
\]
and outside $B_1\cup B_2$ items~\ref{cor:title-2} and~\ref{cor:title-3} hold.
\end{proof}

The constants $\delta$ and $c$ depend only on $\kappa$, $f$ and $\xi$; the threshold $N_0$ depends
on $\kappa$ and $t$ in (a), and on $\kappa$ in (b).

%% file: figures/fig4-edge.tex
\begin{figure}[tp]
\centering
\begin{NoHyper}
\begin{tikzpicture}
\pgfplotsset{
  edgeaxis/.style={
    width=7.6cm, height=5.6cm, scale only axis=false,
    axis lines=box, tick align=inside,
    tick label style={font=\scriptsize}, label style={font=\footnotesize},
    legend style={font=\scriptsize,draw=black!40,fill=white,inner sep=2pt,
      at={(0.5,1.03)},anchor=south,/tikz/every even column/.append style={column sep=8pt}},
    legend cell align=left,
    scaled ticks=false,
    yticklabel style={/pgf/number format/fixed,/pgf/number format/precision=2},
    xticklabel style={/pgf/number format/fixed,/pgf/number format/precision=1},
    y label style={rotate=-90},
  }}

\begin{axis}[edgeaxis, name=left,
    xmin=-0.4, xmax=0.4, ymin=-0.032, ymax=0.08,
    xtick={-0.4,-0.2,0,0.2,0.4}, ytick={-0.02,0,0.02,0.04,0.06,0.08},
    xlabel={$\sigma=s+a\chi$}, ylabel={$z$}, legend columns=3]
  \addplot[black!35,line width=0.4pt,forget plot] coordinates {(-0.4,0) (0.4,0)};
  \addplot[black,dashed,line width=0.9pt] coordinates {
    (-0.450,0.074715) (-0.440,0.067728) (-0.430,0.061100) (-0.420,0.054818) (-0.410,0.048871) (-0.400,0.043249)
    (-0.390,0.037940) (-0.380,0.032934) (-0.370,0.028222) (-0.360,0.023795) (-0.350,0.019643) (-0.340,0.015759)
    (-0.330,0.012134) (-0.320,0.008761) (-0.310,0.005632) (-0.300,0.002740) (-0.290,0.000079) (-0.280,-0.002359)
    (-0.270,-0.004579) (-0.260,-0.006588) (-0.250,-0.008391) (-0.240,-0.009995) (-0.230,-0.011403) (-0.220,-0.012622)
    (-0.210,-0.013657) (-0.200,-0.014512) (-0.190,-0.015192) (-0.180,-0.015702) (-0.170,-0.016045) (-0.160,-0.016227)
    (-0.150,-0.016250) (-0.140,-0.016120) (-0.130,-0.015838) (-0.120,-0.015411) (-0.110,-0.014839) (-0.100,-0.014129)
    (-0.090,-0.013281) (-0.080,-0.012300) (-0.070,-0.011189) (-0.060,-0.009950) (-0.050,-0.008587) (-0.040,-0.007103)
    (-0.030,-0.005499) (-0.020,-0.003779) (-0.010,-0.001945) (0,0) (0.010,0.002054) (0.020,0.004215)
    (0.030,0.006480) (0.040,0.008847) (0.050,0.011314) (0.060,0.013880) (0.070,0.016541) (0.080,0.019296)
    (0.090,0.022143) (0.100,0.025081) (0.110,0.028107) (0.120,0.031220) (0.130,0.034418) (0.140,0.037699)
    (0.150,0.041063) (0.160,0.044506) (0.170,0.048028) (0.180,0.051627) (0.190,0.055303) (0.200,0.059052)
    (0.210,0.062874) (0.220,0.066769) (0.230,0.070733) (0.240,0.074767) (0.250,0.078868) (0.260,0.083036)
    (0.270,0.087270) (0.280,0.091568) (0.290,0.095929) (0.300,0.100352) (0.310,0.104836) (0.320,0.109380)
    (0.330,0.113983) (0.340,0.118644) (0.350,0.123362) (0.360,0.128136) (0.370,0.132965) (0.380,0.137848)
    (0.390,0.142785) (0.400,0.147774) (0.410,0.152815) (0.420,0.157906) (0.430,0.163048) (0.440,0.168238)
    (0.450,0.173477)
  };
  \addlegendentry{$\Delta=0.2$}
  \addplot[black,line width=1.2pt] coordinates {
    (-0.450,0.205893) (-0.440,0.194660) (-0.430,0.183875) (-0.420,0.173523) (-0.410,0.163589) (-0.400,0.154061)
    (-0.390,0.144925) (-0.380,0.136167) (-0.370,0.127777) (-0.360,0.119743) (-0.350,0.112054) (-0.340,0.104699)
    (-0.330,0.097668) (-0.320,0.090951) (-0.310,0.084540) (-0.300,0.078425) (-0.290,0.072599) (-0.280,0.067051)
    (-0.270,0.061776) (-0.260,0.056765) (-0.250,0.052011) (-0.240,0.047507) (-0.230,0.043246) (-0.220,0.039222)
    (-0.210,0.035429) (-0.200,0.031860) (-0.190,0.028509) (-0.180,0.025372) (-0.170,0.022443) (-0.160,0.019716)
    (-0.150,0.017187) (-0.140,0.014851) (-0.130,0.012702) (-0.120,0.010737) (-0.110,0.008951) (-0.100,0.007339)
    (-0.090,0.005899) (-0.080,0.004625) (-0.070,0.003514) (-0.060,0.002562) (-0.050,0.001766) (-0.040,0.001122)
    (-0.030,0.000626) (-0.020,0.000276) (-0.010,0.000069) (0,0) (0.010,0.000068) (0.020,0.000268)
    (0.030,0.000600) (0.040,0.001059) (0.050,0.001643) (0.060,0.002349) (0.070,0.003176) (0.080,0.004120)
    (0.090,0.005179) (0.100,0.006351) (0.110,0.007634) (0.120,0.009025) (0.130,0.010523) (0.140,0.012124)
    (0.150,0.013828) (0.160,0.015632) (0.170,0.017535) (0.180,0.019534) (0.190,0.021628) (0.200,0.023815)
    (0.210,0.026093) (0.220,0.028461) (0.230,0.030916) (0.240,0.033458) (0.250,0.036085) (0.260,0.038795)
    (0.270,0.041587) (0.280,0.044459) (0.290,0.047411) (0.300,0.050440) (0.310,0.053545) (0.320,0.056725)
    (0.330,0.059979) (0.340,0.063305) (0.350,0.066702) (0.360,0.070170) (0.370,0.073706) (0.380,0.077310)
    (0.390,0.080981) (0.400,0.084718) (0.410,0.088518) (0.420,0.092383) (0.430,0.096309) (0.440,0.100298)
    (0.450,0.104346)
  };
  \addlegendentry{$\Delta=0$}
  \addplot[black!55,dash dot,line width=0.9pt] coordinates {
    (-0.450,0.337072) (-0.440,0.321592) (-0.430,0.306650) (-0.420,0.292227) (-0.410,0.278307) (-0.400,0.264873)
    (-0.390,0.251909) (-0.380,0.239401) (-0.370,0.227333) (-0.360,0.215692) (-0.350,0.204465) (-0.340,0.193638)
    (-0.330,0.183201) (-0.320,0.173142) (-0.310,0.163448) (-0.300,0.154110) (-0.290,0.145118) (-0.280,0.136462)
    (-0.270,0.128131) (-0.260,0.120118) (-0.250,0.112413) (-0.240,0.105008) (-0.230,0.097895) (-0.220,0.091066)
    (-0.210,0.084514) (-0.200,0.078232) (-0.190,0.072211) (-0.180,0.066447) (-0.170,0.060932) (-0.160,0.055660)
    (-0.150,0.050625) (-0.140,0.045821) (-0.130,0.041242) (-0.120,0.036884) (-0.110,0.032741) (-0.100,0.028807)
    (-0.090,0.025078) (-0.080,0.021550) (-0.070,0.018217) (-0.060,0.015074) (-0.050,0.012119) (-0.040,0.009346)
    (-0.030,0.006752) (-0.020,0.004332) (-0.010,0.002082) (0,0) (0.010,-0.001919) (0.020,-0.003678)
    (0.030,-0.005280) (0.040,-0.006729) (0.050,-0.008029) (0.060,-0.009181) (0.070,-0.010189) (0.080,-0.011056)
    (0.090,-0.011785) (0.100,-0.012378) (0.110,-0.012839) (0.120,-0.013170) (0.130,-0.013373) (0.140,-0.013451)
    (0.150,-0.013406) (0.160,-0.013241) (0.170,-0.012958) (0.180,-0.012559) (0.190,-0.012046) (0.200,-0.011422)
    (0.210,-0.010688) (0.220,-0.009847) (0.230,-0.008900) (0.240,-0.007850) (0.250,-0.006698) (0.260,-0.005446)
    (0.270,-0.004095) (0.280,-0.002649) (0.290,-0.001107) (0.300,0.000527) (0.310,0.002254) (0.320,0.004070)
    (0.330,0.005974) (0.340,0.007966) (0.350,0.010043) (0.360,0.012204) (0.370,0.014448) (0.380,0.016773)
    (0.390,0.019178) (0.400,0.021661) (0.410,0.024222) (0.420,0.026859) (0.430,0.029571) (0.440,0.032357)
    (0.450,0.035215)
  };
  \addlegendentry{$\Delta=-0.2$}
  \addplot[black!60,line width=0.4pt,forget plot] coordinates {(-0.08,-0.016) (0.08,0.016)};
  \addplot[black!60,line width=0.4pt,forget plot] coordinates {(-0.08,0.016) (0.08,-0.016)};
  \addplot[only marks,mark=*,mark size=1.6pt,black,forget plot] coordinates
    {(-0.15351,-0.016260) (0,0) (0.14132,-0.013452)};
  \node[font=\scriptsize,anchor=north] at (axis cs:-0.15351,-0.0178) {$\lamp$};
  \node[font=\scriptsize,anchor=north] at (axis cs:0.14132,-0.0150) {$\lamp$};
\end{axis}

\begin{axis}[edgeaxis, at={($(left.east)+(2.0cm,0)$)}, anchor=west,
    xmin=-0.3, xmax=0.3, ymin=-0.05, ymax=0.005,
    xtick={-0.3,-0.2,-0.1,0,0.1,0.2,0.3}, ytick={-0.05,-0.04,-0.03,-0.02,-0.01,0},
    xlabel={$\Delta$}, ylabel={$\lamp$}, legend columns=2]
  \fill[black!10] (axis cs:-0.3,-0.05) rectangle (axis cs:0,0.005);
  \node[font=\scriptsize,align=center] at (axis cs:-0.13,-0.038)
    {$\Delta<0$: not attained\\ at contacts};
  \addplot[black!35,line width=0.4pt,forget plot] coordinates {(-0.3,0) (0.3,0)};
  \addplot[black,line width=1.2pt] coordinates {
    (-0.30,-0.029081) (-0.29,-0.027280) (-0.28,-0.025531) (-0.27,-0.023834) (-0.26,-0.022189) (-0.25,-0.020597)
    (-0.24,-0.019058) (-0.23,-0.017574) (-0.22,-0.016144) (-0.21,-0.014770) (-0.20,-0.013452) (-0.19,-0.012191)
    (-0.18,-0.010987) (-0.17,-0.009841) (-0.16,-0.008754) (-0.15,-0.007727) (-0.14,-0.006759) (-0.13,-0.005853)
    (-0.12,-0.005009) (-0.11,-0.004227) (-0.10,-0.003509) (-0.09,-0.002855) (-0.08,-0.002266) (-0.07,-0.001742)
    (-0.06,-0.001286) (-0.05,-0.000897) (-0.04,-0.000577) (-0.03,-0.000326) (-0.02,-0.000146) (-0.01,-0.000037)
    (0,0) (0.01,-0.000037) (0.02,-0.000148) (0.03,-0.000335) (0.04,-0.000599) (0.05,-0.000940)
    (0.06,-0.001361) (0.07,-0.001861) (0.08,-0.002443) (0.09,-0.003107) (0.10,-0.003856) (0.11,-0.004689)
    (0.12,-0.005609) (0.13,-0.006617) (0.14,-0.007715) (0.15,-0.008903) (0.16,-0.010183) (0.17,-0.011558)
    (0.18,-0.013027) (0.19,-0.014594) (0.20,-0.016260) (0.21,-0.018026) (0.22,-0.019894) (0.23,-0.021867)
    (0.24,-0.023946) (0.25,-0.026132) (0.26,-0.028429) (0.27,-0.030838) (0.28,-0.033361) (0.29,-0.036000)
    (0.30,-0.038758)
  };
  \addlegendentry{$\lamp$}
  \addplot[black,dashed,line width=0.9pt] coordinates {
    (-0.30,-0.025422) (-0.29,-0.023939) (-0.28,-0.022491) (-0.27,-0.021078) (-0.26,-0.019701) (-0.25,-0.018360)
    (-0.24,-0.017057) (-0.23,-0.015793) (-0.22,-0.014568) (-0.21,-0.013383) (-0.20,-0.012240) (-0.19,-0.011139)
    (-0.18,-0.010082) (-0.17,-0.009070) (-0.16,-0.008104) (-0.15,-0.007184) (-0.14,-0.006313) (-0.13,-0.005492)
    (-0.12,-0.004721) (-0.11,-0.004003) (-0.10,-0.003338) (-0.09,-0.002729) (-0.08,-0.002176) (-0.07,-0.001682)
    (-0.06,-0.001247) (-0.05,-0.000874) (-0.04,-0.000565) (-0.03,-0.000321) (-0.02,-0.000144) (-0.01,-0.000036)
    (0,0) (0.01,-0.000037) (0.02,-0.000150) (0.03,-0.000341) (0.04,-0.000612) (0.05,-0.000966)
    (0.06,-0.001406) (0.07,-0.001935) (0.08,-0.002554) (0.09,-0.003268) (0.10,-0.004080) (0.11,-0.004992)
    (0.12,-0.006009) (0.13,-0.007133) (0.14,-0.008369) (0.15,-0.009720) (0.16,-0.011191) (0.17,-0.012786)
    (0.18,-0.014509) (0.19,-0.016365) (0.20,-0.018360) (0.21,-0.020498) (0.22,-0.022785) (0.23,-0.025227)
    (0.24,-0.027830) (0.25,-0.030600) (0.26,-0.033544) (0.27,-0.036670) (0.28,-0.039984) (0.29,-0.043495)
    (0.30,-0.047212)
  };
  \addlegendentry{$-\Delta^2/(4a\,\E U^3)$}
  \addplot[only marks,mark=*,mark size=1.6pt,black,forget plot] coordinates
    {(-0.2,-0.013452) (0,0) (0.2,-0.016260)};
\end{axis}
\end{tikzpicture}
\end{NoHyper}
\caption{The edge of the Hessian bulk (Lemma~\ref{lem:hess-edge}). \emph{Left:} the edge
function $z$ of~\eqref{eq:hess-z} in the coordinate $\sigma=s+a\chi$, for three values of the
Plefka factor $\Delta=1-a\E U^2$. Every curve passes through the origin with slope $\Delta$ (thin
segments) and second derivative $2a\E U^3$, and the upper edge $\lamp$ of the reflected Pastur law
$\Pastur{a}{\rho}$ is its minimum (dots). So $\lamp<0$ when $\Delta\ne0$, and $\lamp=0$ exactly
when $\Delta=0$, where the origin itself is the minimum. \emph{Right:} the edge as a function of
$\Delta$ (solid) and its quadratic approximation (dashed), which differ by $O(|\Delta|^3)$ as
$\Delta\to0$ (Lemma~\ref{lem:hess-edge}\ref{lem:hess-edge-c}). At a contact $Q\in\Cf$ the Plefka factor is
$\Plef(Q)=2\Obs_{\zf}''(Q)/\xi''(Q)\ge0$, so the shaded half-line does not occur there, and the
edge vanishes exactly at the contacts where the obstacle is flat to second order,
$\Obs_{\zf}''(Q)=0$ (Theorem~\ref{thm:struct}\ref{thm:struct-iv}). The curves are computed for an illustrative law,
$U$ uniform on $[0.2,0.9]$, which satisfies~\eqref{eq:hess-fhyp}, with $a=(1-\Delta)/\E U^2$ varying with $\Delta$; they are not computed
from a spin-glass model.}
\label{fig:edge}
\end{figure}
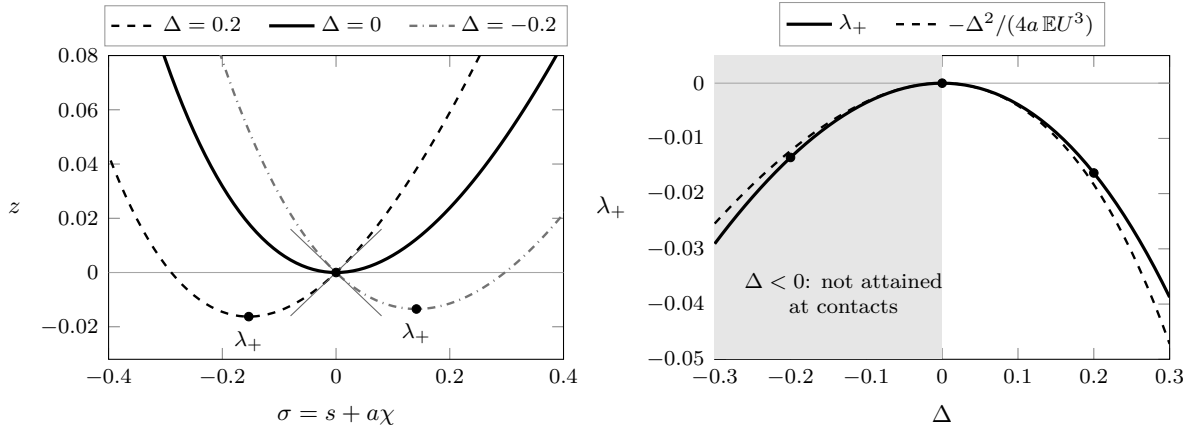

%% file: sections/sec5-stability.tex
\section{\texorpdfstring{Stability costs rate: Theorem~\ref*{thm:stab}}{Stability costs rate:
Theorem 5}}\label{sec:stab}

This section proves Theorem~\ref{thm:stab}. The states of Theorem~\ref{thm:exist} are critical
points of $F_{\zeta^*}$ for a point-dependent order parameter $\zeta^*$, and nothing in their
construction makes them stable (Remark~\ref{rem:exist-points}). By
Theorem~\ref{thm:struct}\ref{thm:struct-iii}--\ref{thm:struct-iv}, the Hessian bulk at every
near-maximizer $x$ at level $f$ follows $\nu_{q_x}$, whose upper edge $\lamp(q_x)$ is nonpositive;
at a contact $Q$, $\lamp(Q)=0$ exactly when the obstacle is flat to second order at $Q$. Can the
rate $\Sigma(f)$ be achieved by critical points of the \emph{fixed} function $F_{\zf}$ at which
$F_{\zf}$ is strongly concave with a fixed margin $g>0$? The Plefka factor on $\Cf$ decides.

Parts~\ref{thm:stab-a} and~\ref{thm:stab-b} form one dichotomy. At a strictly Plefka stable first
contact, part~\ref{thm:stab-a} produces critical points of $F_{\zf}$ at the full rate around which
$F_{\zf}$ is strongly concave with a fixed margin. If every contact is marginal,
part~\ref{thm:stab-b} shows that every fixed margin costs a positive rate $c(g)$, also for states
with their own order parameters, and part~\ref{thm:stab-c} shows that this occurs, on an interval
of levels, for an open set of mixed $2+4$ covariances (open in the parameters $(b,e)$). If the
first contact is marginal but another contact is strictly stable, neither the fixed-margin clause
of part~\ref{thm:stab-a} nor part~\ref{thm:stab-b} applies. The rate $c(g)$ is not explicit: the proof passes through the continuity moduli behind
Theorem~\ref{thm:struct}\ref{thm:struct-iii}, which are not quantitative.

Section~\ref{sec:stab-margin} deduces part~\ref{thm:stab-b} from the edge identity.
Section~\ref{sec:stab-kr} reads part~\ref{thm:stab-a} off the Kac--Rice construction of
Appendix~\ref{app:kr} and states the margin bound and the annealed count of supersymmetric (SUSY)
states, which keeps the exponent $\Sigma(f)$ at marginal levels (proofs in
Appendix~\ref{app:kr-susy}). Section~\ref{sec:stab-family} proves part~\ref{thm:stab-c}.

\subsection{Fixed margins at marginal contacts}\label{sec:stab-margin}

The mechanism behind Theorem~\ref{thm:stab}\ref{thm:stab-b} is a short consequence of the edge
identity. Under \condHmarg, Theorem~\ref{thm:struct}\ref{thm:struct-iv} gives $\lamp(Q)=0$ for
every $Q\in\Cf$, so the bulk law $\nu_Q$ puts positive mass above $-g/4$ for every $g>0$.
Theorem~\ref{thm:struct}\ref{thm:struct-iii} transfers this mass, uniformly over the cube and over
the order parameters, to the empirical spectral distribution of every Hessian in question. A matrix
with a positive fraction of its eigenvalues above $-g/2$ is not bounded above by $-gI$.

\begin{proof}[Proof of Theorem~\ref{thm:stab}\ref{thm:stab-b}]
Fix $g>0$, and recall that $\kappa_f=\dist(\Cf,\{0,1\})>0$.

\emph{Step 1: mass near the edge.} By \condHmarg\ and
Theorem~\ref{thm:struct}\ref{thm:struct-iv}, $\lamp(Q)=\sup\supp\nuQ=0$ for every $Q\in\Cf$. The
open set $(-g/4,\infty)$ contains a point of $\supp\nuQ$, so $\nuQ((-g/4,\infty))>0$ for $Q\in\Cf$.
The map $Q\mapsto\nuQ$ is weakly continuous on $(0,1)$, by the continuity of reflected Pastur laws
in their data (Lemma~\ref{lem:hess-pastur}) and the continuity of $Q\mapsto\Law(\UQ)$
(Lemma~\ref{lem:pde-standing}). By the portmanteau theorem, $Q\mapsto\nuQ((-g/4,\infty))$ is lower
semicontinuous on $(0,1)$. The contact set $\Cf$ is a compact subset of $(0,1)$
(Lemma~\ref{lem:loc-obstacle}), and a lower semicontinuous function attains its minimum on a
compact set, so
\[
 p:=\min_{Q\in\Cf}\nuQ\bigl((-g/4,\infty)\bigr)>0 .
\]
The set $\{Q\in(0,1):\nuQ((-g/4,\infty))>p/2\}$ is open and contains $\Cf$, so it contains
$\{Q:\dist(Q,\Cf)<\kappa_0\}$ for some $\kappa_0\in(0,\kappa_f)$.

\emph{Step 2: transfer to the Hessians.} Put $\kappa:=\min\{\kappa_0,g/4,p/4\}\in(0,\kappa_f)$.
Theorem~\ref{thm:struct}\ref{thm:struct-iii} with this $\kappa$ gives $\varepsilon,c>0$ and, for all
large $N$, an event $B_N$ with $\Pout(B_N)\le e^{N(\Sigma(f)-c)}$ outside of which every
$x\in\ocube$ with $\FTAP(x)\ge N(f-\varepsilon)$ satisfies $\dist(q_x,\Cf)<\kappa$ and
$\dL(\ESD(H),\nu_{q_x})\le\kappa$ for each of the matrices $H=\Hpm F_{\zf}(x)$ and
$H=\Hpm F_\alpha(x)$, $\alpha\in\Pspace([0,1])$ with $F_\alpha(x)\le\FTAP(x)+2\varepsilon N$.

\emph{Step 3: no margin.} Outside $B_N$, let $x$ and $H$ be as in Step~2. If $\dL(\mu,\mu')\le\kappa$,
then $\mu((t,\infty))\ge\mu'((t+\kappa,\infty))-\kappa$ for every $t\in\R$. We apply this with
$t=-g/2$. Since $\kappa\le g/4$, $\dist(q_x,\Cf)<\kappa_0$ and $\kappa\le p/4$,
\[
 \ESD(H)\bigl((-g/2,\infty)\bigr)\ge\nu_{q_x}\bigl((-g/4,\infty)\bigr)-\kappa
 \ge\frac p2-\frac p4=\frac p4 .
\]
So at least $pN/4$ eigenvalues of $H$ exceed $-g/2$, and $\lambda_{\max}(H)>-g$. This holds for
both one-sided Hessians, so $\min_\pm\lambda_{\max}(\Hpm F_\alpha(x))>-g$, both for $\alpha=\zf$ and
for every $\alpha$ with $F_\alpha(x)\le\FTAP(x)+2\varepsilon N$. The two events of
part~\ref{thm:stab-b} are therefore contained in $B_N$, which proves the bound on their outer
probability.

For the last assertion, let $\varepsilon'\in(0,\varepsilon]$ and let $(\zeta,x)$ witness
$\Emarg{g}{\varepsilon'}{f}$, with $\Hess_\sigma F_\zeta(x)\preceq-gI$. Then $x$ lies in $\ocube$,
\begin{gather*}
 \FTAP(x)\ge N(f-\varepsilon'),\qquad
 F_\zeta(x)\le N(f+\varepsilon')\le\FTAP(x)+2\varepsilon' N,\\
 \min_\pm\lambda_{\max}\bigl(\Hpm F_\zeta(x)\bigr)\le\lambda_{\max}\bigl(\Hess_\sigma F_\zeta(x)\bigr)\le-g .
\end{gather*}
Hence $\Emarg{g}{\varepsilon'}{f}\subset B_N$, and
$\Pout(\Emarg{g}{\varepsilon'}{f})\le e^{N(\Sigma(f)-c)}$ for all large $N$.
\end{proof}

No criticality is used, and outlier eigenvalues, which Theorem~\ref{thm:struct}\ref{thm:struct-iii}
does not control, can only increase $\lambda_{\max}$, so they cannot rescue a margin. \condHmarg\ is
used at every contact: at a strictly stable contact $Q$, $\lamp(Q)<0$
(Theorem~\ref{thm:struct}\ref{thm:struct-iv}) leaves room for margins up to $|\lamp(Q)|$. Boursier's
states at a strictly stable first contact are of this kind: they carry a negative Hessian bound on a
neighbourhood, and outliers need a separate argument \cite[Proposition~B.1 and Remark~B.1]{Bou26b}.

\subsection{The Kac--Rice construction}\label{sec:stab-kr}

Theorem~\ref{thm:stab}\ref{thm:stab-a} is a statement about the construction by which
Appendix~\ref{app:kr} proves the existence lower bound of Proposition~\ref{prop:kr-main}. That
construction never changes the function. It counts critical points of the fixed $F_{\zf}$ by a
marked Kac--Rice formula and keeps only those at which a mark holds: $F_{\zf}$ is strongly
concave on a ball of radius $\rho\sqrt N$ around the point, and its value drops by a definite amount at every point of nearly the same
self-overlap whose overlap with the given point lies in $[\delta,\qf-\delta]$. Distinct retained
points then have overlap less than $\delta$, so their number is bounded by a spherical code, and
dividing the expected count by this deterministic bound bounds the existence probability from
below. This isolation and spherical-code argument follows Boursier \cite[Section~3]{Bou26b}, who
adapts the isolated-critical-point argument of Huang and Sellke \cite{HS}.

At a strictly stable first contact the construction runs on the unperturbed data. In that case
Boursier's construction already produces, at the full rate, critical points of $F_{\zf}$ at which
the Hessian is bounded above by a fixed negative multiple of the identity on a neighbourhood
\cite[Definition~3.5, Proposition~B.1 and the proof of Proposition~3.3]{Bou26b}. At a marginal
first contact the unperturbed data have no spectral gap. The construction then tilts the empirical
law of the fields and raises the energy before counting (Appendix~\ref{app:kr-tilt}); both
adjustments act on the data inside the Kac--Rice integral, while $\xi$ and $\zf$ stay fixed. The
energy adjustment is one place where two active degrees enter (they also give the conditional
density of the energy in the count at every tilt, Remark~\ref{rem:kr-summable}): the value of
$\HN$ at a point keeps a positive conditional variance after its gradient is fixed
(Lemma~\ref{lem:kr-energy}).

We record what the construction provides. Write $q=\qf$, $a=\xi''(q)$,
$M(b)=\Phi_{\zf,x}(q,b)$, $U(b)=\Phi_{\zf,xx}(q,b)$ and $\pi_0=\Law(X_q)$, where $X$ is the
Auffinger--Chen process of $\zf$, so that $\Plef(\qf)=1-a\E_{\pi_0}U^2$. If $\Plef(\qf)=0$,
Lemma~\ref{lem:kr-tilt} provides a bounded function $\mathcal R\ne0$ with $\E_{\pi_0}\mathcal R=0$
such that, for $0<\eta\le\eta_0:=(2\|\mathcal R\|_\infty)^{-1}$, the tilted law
$\pi_\eta=(1-\eta\mathcal R)\,\pi_0$ satisfies
\begin{equation}\label{eq:stab-tilt}
 \E_{\pi_\eta}M^2=q,\qquad
 \Plef_\eta:=1-a\E_{\pi_\eta}U^2=a\eta\,\E_{\pi_0}\mathcal R^2>0 .
\end{equation}
Call a tilt $\eta$ \emph{admissible} if $\eta=0$ when $\Plef(\qf)>0$, and if $\eta\in(0,\eta_0]$
when $\Plef(\qf)=0$ (Appendix~\ref{app:kr} uses every such tilt); at $\eta=0$ put
$\pi_\eta=\pi_0$. For an admissible $\eta$, an overlap cutoff $\delta\in(0,q/2)$ and a window
width $w\in(0,w_1]$ as in~\ref{step:kr-par6}, let $\widehat C_N=\widehat C_N(\eta,\delta,w)$ be the
marked count of Appendix~\ref{app:kr-count}, and call the points it counts \emph{retained}. For
$x\in\ocube$ write $\pi_x=N^{-1}\sum_i\delta_{M^{-1}(x_i)}$ for the empirical law of its fields; here
$M$ is an increasing bijection of $\R$ onto $(-1,1)$. The proof of Proposition~\ref{prop:kr-main}
gives the following.
\begin{enumerate}[label=\textup{(K\arabic*)},ref=(K\arabic*),leftmargin=*,widest=3,topsep=3pt,itemsep=2pt]
\item\label{step:stab-k1} \emph{Margin.} For every admissible $\eta$ there are $g_\eta,\rho_\eta>0$,
 not depending on $\delta$ or $w$, such that every retained point $x$ is a critical point of $F_{\zf}$
 in $\ocube$ and $\Hpm F_{\zf}(y)\preceq-g_\eta I$ for every $y\in\ocube$ with
 $|y-x|\le\rho_\eta\sqrt N$ (Proposition~\ref{prop:kr-local}; the mark requires this bound).
\item\label{step:stab-k2} \emph{Windows.} Every retained point satisfies $\Wone(\mu_x,M_\#\pi_\eta)
 \le W_2(\pi_x,\pi_\eta)\le w$, $|q_x-q|\le2w$, and
 \[
  \bigl|N^{-1}F_{\zf}(x)-f\bigr|+\bigl|N^{-1}\FTAP(x)-f\bigr|\le\varepsilon_1(\eta)+\varepsilon_2(w)+2N^{-2},
 \]
 where $\varepsilon_1(0)=0$, $\varepsilon_1(\eta)\to0$ as $\eta\downarrow0$ and $\varepsilon_2(w)\to0$
 as $w\downarrow0$.
\item\label{step:stab-k3} \emph{Rate.} For every admissible $\eta$ and every $\delta\in(0,q/2)$,
 \[
  \liminf_{w\downarrow0}\liminf_{N\to\infty}\frac1N\log\Prob\bigl(\widehat C_N(\eta,\delta,w)\ge1\bigr)
  \ge\Sigma(f)-\varepsilon_3(\eta)+\frac12\log\Bigl(1-\frac\delta q\Bigr),
 \]
 where $\varepsilon_3(0)=0$ and $\varepsilon_3(\eta)\to0$ as $\eta\downarrow0$.
\end{enumerate}

\begin{proof}[Proof of Theorem~\ref{thm:stab}\ref{thm:stab-a}]
Fix $\varepsilon>0$. Choose an admissible $\eta$ with $\varepsilon_1(\eta)\le\varepsilon/3$ and
$\varepsilon_3(\eta)\le\varepsilon/4$; if $\Plef(\qf)>0$ this is $\eta=0$. Put $g=g_\eta$ and
$\rho=\rho_\eta$. Choose $\delta\in(0,q/2)$ with $\frac12\log(1-\delta/q)\ge-\varepsilon/4$. By
\ref{step:stab-k3} there is $w_2>0$ such that
$\liminf_NN^{-1}\log\Prob(\widehat C_N(\eta,\delta,w)\ge1)\ge\Sigma(f)-\varepsilon$ for every
$w\in(0,w_2]$; fix such a $w$ with $\varepsilon_2(w)\le\varepsilon/3$. For $N^{-2}\le\varepsilon/6$,
\ref{step:stab-k1} and \ref{step:stab-k2} show that every retained point $x$ satisfies
\eqref{eq:event-exist} with $\zeta=\zf$, and \eqref{eq:event-fix}. Hence
$\{\widehat C_N\ge1\}\subset\Efix{g}{\rho}{\varepsilon}{f}$ for all large $N$, which gives the first
display of part~\ref{thm:stab-a}.

If $\Plef(\qf)>0$, then $\eta=0$ for every $\varepsilon$, and $(g,\rho)=(g_0,\rho_0)$ does not
depend on $\varepsilon$. The lower bound $\Sigma(f)-\varepsilon$ then holds at this fixed margin for
every $\varepsilon>0$. Since $\Efix{g}{\rho}{\varepsilon}{f}\subset\Eex{\varepsilon}{f}$,
Theorem~\ref{thm:exist}\ref{thm:exist-d}--\ref{thm:exist-e} gives the matching upper bound
$\lim_{\varepsilon\downarrow0}\limsup_NN^{-1}\log\Pout(\Efix{g}{\rho}{\varepsilon}{f})\le\Sigma(f)$.
The lower bound comes from the measurable subevent $\{\widehat C_N\ge1\}$, so the equality holds
for the inner and for the outer probability.
\end{proof}

At a marginal first contact the margin $g_\eta$ depends on the tilt and vanishes with it. The
tilted data are strictly Plefka stable, $\Plef_\eta>0$, but only by $O(\eta)$, and the edge of the
tilted bulk is then only $O(\eta^2)$ below zero. For an admissible $\eta$ write $\chi_\eta=\E_{\pi_\eta}U$ and
\[
 \nu^{(\eta)}=\Pastur{a}{\Law_{b\sim\pi_\eta}\bigl(U(b)^{-1}+a\chi_\eta\bigr)},
 \qquad \lamp^{(\eta)}=\sup\supp\nu^{(\eta)} .
\]
At $\eta=0$ this is the bulk law $\nu_{\qf}$ of Theorem~\ref{thm:struct}, since $U(b)$ under
$b\sim\pi_0$ has the law of $U_{\qf}$.

\begin{remark}[The margins of the construction]\label{rem:stab-margin}
Assume $\Plef(\qf)=0$. For every admissible $\eta>0$, every margin $g_\eta$ as in
\ref{step:stab-k1} satisfies $g_\eta\le-\lamp^{(\eta)}$, and, as $\eta\downarrow0$,
\[
 -\lamp^{(\eta)}=\frac{\Plef_\eta^2}{4a\E_{\pi_\eta}U^3}+O(\eta^3)
 =\frac{a(\E_{\pi_0}\mathcal R^2)^2}{4\E_{\pi_0}U^3}\,\eta^2+O(\eta^3).
\]
So the margins of the construction are $O(\eta^2)$. In the $2+4$ family of
Theorem~\ref{thm:stab}\ref{thm:stab-c} some such vanishing is unavoidable, by
item~\ref{thm:stab-c2} there. The remark is proved in Appendix~\ref{app:kr-susy}.
\end{remark}

The proof uses Proposition~\ref{prop:hess-bulk} at the tilted data and the edge expansion of
Lemma~\ref{lem:hess-edge}.

\paragraph{Counting does not see the margin.}
Part~\ref{thm:stab-b} is a statement about stability, not about criticality. At a marginal first
contact, critical points of $F_{\zf}$ with the self-overlap and field law of the first contact are
still as numerous, in expectation, as the existence probability allows. To state this we use the
windows of Boursier \cite[Definition~3.4]{Bou26b}, where $F_\zeta$ is written $F_{\mathrm{TAP},\zeta}$.

\begin{definition}[SUSY states]\label{def:stab-susy}
Let $f\in(\feq,\fone)$, and let $M$, $\pi_0$ and $\pi_x$ be as above. For $\varepsilon>0$, an
\emph{$\varepsilon$-SUSY state at level $f$} is a critical point $x\in\ocube$ of $F_{\zf}$ such that
\[
 |q_x-\qf|\le\varepsilon,\qquad
 \bigl|N^{-1}F_{\zf}(x)-f\bigr|\le\varepsilon,\qquad
 \bigl|N^{-1}\FTAP(x)-f\bigr|\le\varepsilon,\qquad
 W_2(\pi_x,\pi_0)\le\varepsilon .
\]
Let $\mathsf N_{N,\varepsilon}$ be the number of $\varepsilon$-SUSY states.
\end{definition}

\begin{proposition}[Annealed count of SUSY states]\label{prop:stab-susy}
Assume \condH, $\xi\in\classC$ and that $\xi$ has at least two active degrees, and let
$f\in(\feq,\fone)$.
\begin{romanlist}
\item\label{prop:stab-susy-i} For every $\varepsilon>0$,
 $\liminf_NN^{-1}\log\E\mathsf N_{N,\varepsilon}\ge\Sigma(f)$, and
 \[
  \lim_{\varepsilon\downarrow0}\limsup_{N\to\infty}\frac1N\log\E\mathsf N_{N,\varepsilon}\le\Sigma(f).
 \]
 Consequently both $\lim_{\varepsilon\downarrow0}\liminf_N$ and
 $\lim_{\varepsilon\downarrow0}\limsup_N$ of $N^{-1}\log\E\mathsf N_{N,\varepsilon}$ equal
 \[
  \Sigma(f)=\thf(P_{\zf}-f).
 \]
\item\label{prop:stab-susy-ii} For every $\varepsilon>0$,
 $\liminf_NN^{-1}\log\Prob(\mathsf N_{N,\varepsilon}\ge1)\ge\Sigma(f)$, and
 \[
  \lim_{\varepsilon\downarrow0}\limsup_{N\to\infty}\frac1N
  \log\E\bigl[\mathsf N_{N,\varepsilon}\bigm|\mathsf N_{N,\varepsilon}\ge1\bigr]=0 .
 \]
\end{romanlist}
No strict Plefka condition is imposed: $\Plef(\qf)=0$ is allowed.
\end{proposition}

The limits in $\varepsilon$ exist because $\mathsf N_{N,\varepsilon}$ is nondecreasing in
$\varepsilon$, and $\Lambda(\thf)-\thf f=\thf(P_{\zf}-f)$ because $\zf$ attains the constrained
minimum at $\thf$. The proof is in Appendix~\ref{app:kr-susy}: the lower bounds follow from the
construction through~\ref{step:stab-k1}--\ref{step:stab-k3}, and the upper bound is
Proposition~\ref{prop:kr-susy}, whose only stability input is the Plefka inequality
$\Plef(\qf)\ge0$.

For $\Plef(\qf)>0$, part~\ref{prop:stab-susy-i} is Boursier's \cite[Proposition~3.6]{Bou26b};
Proposition~\ref{prop:stab-susy} adds the marginal case. Her earlier annealed count
\cite[Proposition~2.2]{Bou26a}, with no Plefka hypothesis stated, uses another window
\cite[Definition~2.1]{Bou26a}, which keeps the coordinates in
$[-1+\varepsilon^{1/2},1-\varepsilon^{1/2}]$ and constrains the self-overlap, a scalar functional of
the magnetizations and the optimal order parameter of the point, but not the empirical law of the
fields or $\FTAP$. We derive neither count from the other. Boursier reaches a single active degree
by a perturbation of $\xi$ preserving the strict Plefka inequality \cite[p.~20]{Bou26b}; we do not
pursue that case.

Proposition~\ref{prop:stab-susy} and Theorem~\ref{thm:stab}\ref{thm:stab-b} separate existence from
stability among the same states. Assume in addition \condHmarg, let $g>0$, and let $\varepsilon,c$
be as in Theorem~\ref{thm:stab}\ref{thm:stab-b}. Every $\varepsilon'$-SUSY state $x$ at which
$\Hess_\sigma F_{\zf}(x)\preceq-gI$ for some $\sigma$ witnesses $\Emarg{g}{\varepsilon'}{f}$. For
$\varepsilon'\in(0,\varepsilon]$ and all large $N$, part~\ref{prop:stab-susy-ii} therefore gives
\[
 \Pout\bigl(\Emarg{g}{\varepsilon'}{f}\bigm|\{\mathsf N_{N,\varepsilon'}\ge1\}\bigr)
 \le\frac{\Pout(\Emarg{g}{\varepsilon'}{f})}{\Pin(\mathsf N_{N,\varepsilon'}\ge1)}
 \le e^{-cN/2} .
\]
So, given that SUSY states exist, the conditional outer probability that some witness of
$\Eex{\varepsilon'}{f}$, for instance an $\varepsilon'$-SUSY state, has a one-sided Hessian bounded
above by $-gI$ is at most $e^{-cN/2}$.

\subsection{The near-critical \texorpdfstring{$2+4$}{2+4} family}\label{sec:stab-family}

Theorem~\ref{thm:stab}\ref{thm:stab-b} needs \condHmarg\ at every contact, which is a degenerate
condition. This subsection shows that it holds, on an interval of levels, for an open set of
mixed $2+4$ covariances (open in the parameters $(b,e)$). The input is the following description
of the constrained minimizers near criticality, proved in Appendix~\ref{app:nc} (Proposition~\ref{prop:nc-example}). Its proof
adapts the near-critical gap-exclusion mechanism of Zhou \cite{Zho26a} to the constrained mixed problem.
There is $\varepsilon_*>0$ such that every covariance $\xi(t)=bt^2+et^4$ with $b>\frac12$, $e>0$ and
$b-\frac12+e<\varepsilon_*$ has $\theta_0\in(0,1)$ and $\fzero\in(\feq,\fone]$ with the following
properties. For $\theta\in(0,\theta_0)$ the constrained minimizer $\zeta_\theta$ satisfies
\[
 \zeta_\theta(\{0\})=\theta,\qquad
 \supp\zeta_\theta=\{0\}\cup[q_\theta,v_\theta],\qquad 0<q_\theta<v_\theta<1,
\]
and its first positive contact $q_\theta$ is marginal; every $f\in(\feq,\fzero)$ has
$\thf\in(0,\theta_0)$. Marginality of the first contact alone does not suffice for
part~\ref{thm:stab-b}. The following lemma shows that the whole contact set is marginal.

\begin{lemma}[Every contact is marginal]\label{lem:stab-family}
Let $\xi$ and $\theta_0$ be as above, let $\theta\in(0,\theta_0)$ and $\zeta=\zeta_\theta$, and put
$h(t)=\Gz{\zeta}(t)-t$. Then
\[
 h>0\ \text{ on }(0,q_\theta),\qquad h=0\ \text{ on }[q_\theta,v_\theta],\qquad
 h<0\ \text{ on }(v_\theta,1].
\]
Consequently $\Obs_\zeta$ is strictly decreasing on $[0,q_\theta]$, constant on
$[q_\theta,v_\theta]$ and strictly increasing on $[v_\theta,1]$, so that
$\argmin_{(0,1]}\Obs_\zeta=[q_\theta,v_\theta]=\supp\zeta\cap(0,1]$, and
$\Obs_\zeta''=\frac12\xi''(1-\Gz{\zeta}')=0$ on $[q_\theta,v_\theta]$.
\end{lemma}

The proof (Appendix~\ref{app:nc-marginal}) is a Rolle argument on the two plateaus of $\zeta$:
$h(0)=0$, $h=h'=0$ on the contacts $[q_\theta,v_\theta]$, $h'''\ge2$ on the plateaus below a
threshold $t_{\mathrm c}\in(v_\theta,1)$, and $h<0$ on $[t_{\mathrm c},1]$ (Corollary~\ref{cor:nc-standing}).
So $h$ has no further zero, and $\Obs_\zeta'=-\frac12\xi''h$ and the minimality of $\Obs_\zeta$ on
$[q_\theta,v_\theta]$ fix its signs.

\begin{proof}[Proof of Theorem~\ref{thm:stab}\ref{thm:stab-c}]
Let $\xi$ be in the family and $f\in(\feq,\fzero)$. The covariance is a polynomial with
$\xi''(t)=2b+12et^2>0$, so $\xi$ satisfies \condH, $\xi\in\classC$, and $\xi$ has two active
degrees. Since $\thf\in(0,\theta_0)$, Lemma~\ref{lem:stab-family} applies to $\zf=\zeta_{\thf}$.

\emph{Item~\ref{thm:stab-c1}.} By Lemma~\ref{lem:stab-family},
$\Cf=\supp\zf\cap(0,1]=[q_{\thf},v_{\thf}]$, an interval of positive length, and
$\Obs_{\zf}''=0$ on $\Cf$, which is \condHmarg.

\emph{Item~\ref{thm:stab-c2}.} Suppose that $g(\varepsilon_n)\ge g_0>0$ along a sequence
$\varepsilon_n\downarrow0$. Part~\ref{thm:stab-b} at the margin $g_0$, which applies by
item~\ref{thm:stab-c1}, gives $\varepsilon(g_0),c(g_0)>0$ such that
$\limsup_NN^{-1}\log\Prob(\Emarg{g_0}{\varepsilon'}{f})\le\Sigma(f)-c(g_0)$ for every
$\varepsilon'\in(0,\varepsilon(g_0)]$. A point in $\Efix{g}{\rho}{\varepsilon}{f}$ witnesses
$\Emarg{g}{\varepsilon}{f}$ with $\zeta=\zf$, and $\Emarg{g}{\varepsilon}{f}\subset
\Emarg{g_0}{\varepsilon}{f}$ for $g\ge g_0$. For $n$ with
$\varepsilon_n<\min\{\varepsilon(g_0),c(g_0)\}$, part~\ref{thm:stab-a} therefore gives
\[
 \Sigma(f)-\varepsilon_n\le\liminf_{N\to\infty}\frac1N
 \log\Prob\bigl(\Efix{g(\varepsilon_n)}{\rho(\varepsilon_n)}{\varepsilon_n}{f}\bigr)
 \le\limsup_{N\to\infty}\frac1N\log\Prob\bigl(\Emarg{g_0}{\varepsilon_n}{f}\bigr)
 \le\Sigma(f)-c(g_0),
\]
a contradiction. Hence $g(\varepsilon)\to0$.

\emph{Item~\ref{thm:stab-c3}.} By item~\ref{thm:stab-c1} and
Theorem~\ref{thm:struct}\ref{thm:struct-iv}, $\lamp(Q)=0$ for every $Q\in\Cf$. Let $\kappa>0$ be
small enough for Corollary~\ref{cor:hess-marg}. It gives $\kappa',\varepsilon_4,c_4>0$ such that, for
all large $N$, outside an event of outer probability at most $e^{N(\Sigma(f)-c_4)}$, every
$x\in\ocube$ with $\FTAP(x)\ge N(f-\varepsilon_4)$ admits $Q\in\Cf$ for which each of the matrices
$\Hpm F_{\zf}(x)$, $\Hpm F_\alpha(x)$ with $F_\alpha(x)\le\FTAP(x)+2\varepsilon_4N$, and their
compressions to the orthogonal complement of $x$, has its $\lceil\kappa'N\rceil$-th largest
eigenvalue in $[\lamp(Q)-\kappa,\lamp(Q)+\kappa]=[-\kappa,\kappa]$ and fewer than $\kappa'N$
eigenvalues larger than $\lamp(Q)+\kappa=\kappa$. Theorem~\ref{thm:struct}\ref{thm:struct-i}--\ref{thm:struct-ii}
with this $\kappa$ give $\delta_2,c_2>0$ such that, outside an event of outer probability at most
$2e^{N(\Sigma(f)-c_2)}$, every $m\in\cube$ with $\FTAP(m)\ge N(f-\delta_2)$ satisfies
$\dist(q_m,\Cf)<\kappa$ and $\Wone(\mu_m,\Lq{q_m})<\kappa$. Put
$\varepsilon=\min\{\varepsilon_4,\delta_2\}$ and $c_3=\min\{c_4,c_2\}$, and let $B_N$ be the union of
the two exceptional events, so that $\Pout(B_N)\le3e^{N(\Sigma(f)-c_3)}$.

Let $(\alpha,x)$ satisfy \eqref{eq:event-exist}. Then $x\in\ocube$, $\FTAP(x)\ge N(f-\varepsilon)$ and
$F_\alpha(x)\le N(f+\varepsilon)\le\FTAP(x)+2\varepsilon N$, so outside $B_N$ the pair $(\alpha,x)$
has all the properties listed in item~\ref{thm:stab-c3}. It remains to condition. By
Theorem~\ref{thm:exist}\ref{thm:exist-d}--\ref{thm:exist-e} and the monotonicity of
$\Eex{\varepsilon}{f}$ in $\varepsilon$ (for the inner probability, Section~\ref{sec:max-existence}),
$\liminf_NN^{-1}\log\Pin(\Eex{\varepsilon}{f})\ge\Sigma(f)$. Hence
$\Pin(\Eex{\varepsilon}{f})\ge e^{N(\Sigma(f)-c_3/4)}$ for all large $N$, and
\[
 \Pout\bigl(B_N\bigm|\Eex{\varepsilon}{f}\bigr)
 \le\frac{\Pout(B_N)}{\Pin(\Eex{\varepsilon}{f})}\le3e^{-3c_3N/4}\le e^{-c_3N/2}
\]
for all large $N$. This proves item~\ref{thm:stab-c3} with $c=c_3/2$.
\end{proof}

%% file: sections/sec6-equilibrium.tex
\section{\texorpdfstring{The equilibrium level: Corollary~\ref*{cor:eq}}%
{The equilibrium level: Corollary 6}}\label{sec:eq}

At the equilibrium level the rate vanishes, $\Sigma(\feq)=0$, and the free energy is carried by
typical states. For $\xi\in\classC$ and every $q$ in the support of the Parisi measure,
$N^{-1}\max\{\FTAP(m):q_m=q\}-N^{-1}\log\ZN\to0$ in probability \cite[Theorem~2]{CPS} (their
normalization adds $N\log2$ to both terms), and $N^{-1}\log\ZN\to\feq$ by the Parisi formula
\cite{PF,Pan14}. This section proves Corollary~\ref{cor:eq}, which answers the reverse question:
which radii and magnetization profiles can carry a TAP value within $cN$ of $N\feq$? Only radii
near the contact set $\CP\supset\supp\zP$ of the Parisi obstacle, and only profiles whose
empirical law is close to the Auffinger--Chen field law of $\zP$; away from $\CP$ the maximum of
$\FTAP$ falls short of $N\feq$ by a linear amount, uniformly over all spheres at a fixed distance
from $\CP$, outside an exponentially small event.

The proof is that of Theorem~\ref{thm:struct}\ref{thm:struct-i}--\ref{thm:struct-ii}, run at
the equilibrium level, with two changes. The target probability is $e^{-cN}$ rather than
$e^{N(\Sigma(f)-c)}$, so the tilt $\theta$ of the Markov step is no longer tied to the level; it
only has to stay below the first parameter of the cascade, and we take it small. And the gap of
Proposition~\ref{prop:loc-gap} is replaced by Lemma~\ref{lem:eq-gap}, which uses a single
perturbation: at a radius $Q$ off the contact set, move mass across $Q$ toward the contact
point $0$.

That $0$ is a contact point at equilibrium is the one additional input (Lemma~\ref{lem:eq-zero}),
and Corollary~\ref{cor:eq} would be false without it. Indeed, for every realization,
\begin{equation}\label{eq:eq-tapzero}
 \FTAP(0)=N\feq ,
\end{equation}
as observed by Boursier \cite[Section~2.1]{Bou26b}: $\HN(0)=0$ because $\xi(0)=0$, and
$\Phistar_\zeta(0,0)=\Phi_\zeta(0,0)$ because $\Phi_\zeta(0,\cdot)$ is even and convex, so
$F_\zeta(0)=NP_\zeta$ for every $\zeta$ and $\FTAP(0)=N\inf_\zeta P_\zeta$. If $0$ were not in
$\CP$, Corollary~\ref{cor:eq} with $\kappa<\dist(0,\CP)$ would exclude $m=0$. At an upper level
$f>\feq$ the situation is reversed: $\FTAP(0)=N\feq<N(f-\delta)$ for small $\delta$, and
correspondingly the obstacle of $\zf$ is strictly larger at $0$ than on $\Cf$
(Lemma~\ref{lem:loc-obstacle}).

Section~\ref{sec:eq-obstacle} proves three lemmas on the Parisi obstacle,
Section~\ref{sec:eq-proof} deduces Corollary~\ref{cor:eq} from them and the comparison of
Section~\ref{sec:loc}, and Section~\ref{sec:eq-remarks} relates it to the TAP representation.

\subsection{The Parisi obstacle at equilibrium}\label{sec:eq-obstacle}

Throughout this section $\xi$ satisfies \condH, and $\zP$ is a Parisi measure, that is, a
minimizer of $\zeta\mapsto P_\zeta$ over $\Pspace([0,1])$. Minimizers exist because $P$ is
continuous for the Wasserstein distance $\Wone$ (Corollary~\ref{cor:pde-step} and
Lemma~\ref{lem:pde-measurable}) and $\Pspace([0,1])$ is compact. The minimizer is unique
(Lemma~\ref{lem:var-basic}\ref{lem:var-basic-a}; for $\xi\in\classC$ also \cite[Theorem~2]{AC15b});
we do not use uniqueness. We write
$u(t,x)=\partial_x\Phi_{\zP}(t,x)$, and $X$ for the Auffinger--Chen process of $\zP$, the
solution of
$dX_t=\xi''(t)\cdf{\zP}(t)\,u(t,X_t)\,dt+\sqrt{\xi''(t)}\,dB_t$ with $X_0=0$
(Lemmas~\ref{lem:pde-ac} and~\ref{lem:pde-measurable}). Then
\[
 \Gz{\zP}(t)=\E u(t,X_t)^2,\qquad
 \OP(t)=\frac12\int_t^1\xi''(s)\bigl(\Gz{\zP}(s)-s\bigr)\,ds,\qquad
 \LPq{q}=\Law\bigl(u(q,X_q)\bigr),
\]
and $\CP=\argmin_{[0,1]}\OP$. Unlike $\Cf$, which is an $\argmin$ over $(0,1]$, the set $\CP$ is
an $\argmin$ over the closed interval.

We use three facts from Section~\ref{sec:loc-obstacle}, valid under \condH, on the Parisi
functional $P(c)$ of a coefficient $c:[0,1)\to[0,3]$, so that $P_\zeta=P(\cdf{\zeta})$. Since
$\OP\in C^1([0,1])$ and $\OP(1)=0$, Fubini's theorem gives, for every $\mu\in\Pspace([0,1])$,
\begin{equation}\label{eq:eq-fubini}
 -\int_0^1\cdf{\mu}(t)\,\OP'(t)\,dt=\int_{[0,1]}\OP\,d\mu ;
\end{equation}
the first variation of $P$ at $\cdf{\zP}$ in a direction $h$ is
\begin{equation}\label{eq:eq-firstvar}
 DP(\cdf{\zP})[h]=\frac12\int_0^1\xi''(t)\,h(t)\bigl(\Gz{\zP}(t)-t\bigr)\,dt
 =-\int_0^1h(t)\,\OP'(t)\,dt ;
\end{equation}
and, with the constant $K$ of \eqref{eq:loc-second-order}, which depends only on $\xi'(1)$ and
$\xi''(1)$, and since $P_{\zP}=\feq$,
\begin{equation}\label{eq:eq-second}
 P(\cdf{\zP}+\varepsilon h)\le\feq+\varepsilon\,DP(\cdf{\zP})[h]+K\varepsilon^2
\end{equation}
whenever $\|h\|_{L^1}\le1$, $\varepsilon\ge0$ and $\cdf{\zP}+\varepsilon h$ takes values in
$[0,3]$. The first two give \eqref{eq:loc-moving-mass} at $\zeta=\zP$; the third is
\eqref{eq:loc-second-order} at $c=\cdf{\zP}$.

The first lemma is the first-order optimality condition for Parisi measures, due to Jagannath
and Tobasco \cite[Proposition~1.1]{JT17} (see also \cite{Che17}), together with the argument of
\cite[Lemma~3.8]{JT17} at the right endpoint. We include the short
proofs because they use only the facts above.

\begin{lemma}[First-order conditions at equilibrium]\label{lem:eq-obstacle}
Assume \condH\ and let $\zP$ be a Parisi measure.
\begin{thmlist}
\item\label{lem:eq-obstacle-a} $\int\OP\,d\mu\ge\int\OP\,d\zP$ for every
 $\mu\in\Pspace([0,1])$. Consequently $\int\OP\,d\zP=\min_{[0,1]}\OP$, and
 $\supp\zP\subset\CP$.
\item\label{lem:eq-obstacle-b} $\OP'(1)=\frac12\xi''(1)\bigl(1-\E\tanh^2X_1\bigr)>0$.
 Consequently $\min_{[0,1]}\OP<0=\OP(1)$, and $\CP$ is a nonempty compact subset of $[0,1)$.
\end{thmlist}
\end{lemma}

\begin{proof}
\ref{lem:eq-obstacle-a} Let $\mu\in\Pspace([0,1])$ and $\lambda\in(0,1]$. The measure
$(1-\lambda)\zP+\lambda\mu$ is admissible, and its distribution function is
$\cdf{\zP}+\lambda h$ with $h=\cdf{\mu}-\cdf{\zP}$, which has $\|h\|_{L^1}\le1$. Minimality of
$\zP$ and \eqref{eq:eq-second} give
$0\le P_{(1-\lambda)\zP+\lambda\mu}-\feq\le\lambda\,DP(\cdf{\zP})[h]+K\lambda^2$. Dividing by
$\lambda$ and letting $\lambda\downarrow0$, then applying \eqref{eq:eq-firstvar} and
\eqref{eq:eq-fubini} to $\mu$ and to $\zP$,
\[
 0\le DP(\cdf{\zP})[\cdf{\mu}-\cdf{\zP}]=\int\OP\,d\mu-\int\OP\,d\zP .
\]
With $\mu=\delta_t$ this reads $\OP(t)\ge\int\OP\,d\zP$ for every $t\in[0,1]$, so
$\int\OP\,d\zP=\min\OP$. Hence $\OP=\min\OP$ holds $\zP$-almost everywhere; the set where it
holds is closed, so it contains $\supp\zP$.

\ref{lem:eq-obstacle-b} Since $\Phi_{\zP}(1,x)=\log\cosh x$, we have $u(1,x)=\tanh x$,
so $\Gz{\zP}(1)=\E\tanh^2X_1<1$ and
$\OP'(1)=-\frac12\xi''(1)(\Gz{\zP}(1)-1)>0$, since $\xi''(1)>0$. As $\OP(1)=0$, the
obstacle is negative just below $1$, so $\min\OP<0$ and $1\notin\CP$. The set $\CP$ is closed
in $[0,1]$ by continuity of $\OP$, and nonempty.
\end{proof}

The next lemma is the additional input. If $0\in\supp\zP$ it is part~\ref{lem:eq-obstacle-a} of the
previous lemma. Otherwise $\zP$ leaves an initial gap $[0,\underline q)$ free of mass; there the
coefficient of the Parisi equation vanishes, the Auffinger--Chen process is a time-changed
Brownian motion, and $u(t,X_t)$ is a Doob martingale. Its quadratic variation makes
$\Gz{\zP}$ convex on the gap. A convex function that vanishes at $0$ and meets the diagonal
at $\underline q$ stays below the diagonal, so the obstacle increases across the gap.

\begin{lemma}[The obstacle is minimal at the origin]\label{lem:eq-zero}
Assume \condH. For every Parisi measure $\zP$, $\OP(0)=\min_{[0,1]}\OP$; that is,
$0\in\CP$.
\end{lemma}

\begin{proof}
Let $\underline q=\min\supp\zP$. If $\underline q=0$, then $0\in\CP$ by
Lemma~\ref{lem:eq-obstacle}\ref{lem:eq-obstacle-a}. Suppose $\underline q>0$. By
Lemma~\ref{lem:eq-obstacle}, $\underline q\in\CP$ and $\underline q<1$, so $\underline q$ is an
interior minimum point of the $C^1$ function $\OP$. Hence $\OP'(\underline q)=0$, and since
$\xi''(\underline q)>0$,
\begin{equation}\label{eq:eq-gap-end}
 \Gz{\zP}(\underline q)=\underline q .
\end{equation}

\emph{The gap.} We have $\cdf{\zP}=0$ on $[0,\underline q)$. Approximate $\cdf{\zP}$ in $L^1$ by
step distribution functions that also vanish on $[0,\underline q)$ (for instance, those of the
measures obtained by rounding $\zP$ up to a finite grid). On $[0,\underline q)$ each
approximating Parisi equation has coefficient $0$, so its solution is a Gaussian convolution of
its value at time $\underline q$, and so is its space derivative. Passing to the uniform limit
(Lemma~\ref{lem:pde-measurable}) gives, with $\varphi:=u(\underline q,\cdot)$,
\[
 u(t,x)=\E\,\varphi(x+\sigma_tZ),\qquad
 \sigma_t^2=\xi'(\underline q)-\xi'(t)>0,\qquad 0\le t<\underline q,
\]
where $Z$ is a standard Gaussian variable. By Lemma~\ref{lem:pde-measurable}, $\varphi$ is odd,
nondecreasing and $1$-Lipschitz. The drift of $X$ vanishes on $[0,\underline q)$, so
$X_t=\int_0^t\sqrt{\xi''(r)}\,dB_r$ there: a Gaussian martingale whose increment
$X_{\underline q}-X_t\sim N(0,\sigma_t^2)$ is independent of $\mathcal F_t$, and
$X_{\underline q}$ has a density because $\xi'(\underline q)>0$. Hence, for $t<\underline q$,
\[
 Y_t:=u(t,X_t)=\E\bigl[\varphi(X_{\underline q})\bigm|\mathcal F_t\bigr],
 \qquad
 Y'_t:=\partial_xu(t,X_t)=\E\bigl[\varphi'(X_{\underline q})\bigm|\mathcal F_t\bigr],
\]
where we differentiated under the Gaussian integral ($\varphi$ is Lipschitz, so $\varphi'$ exists
almost everywhere and is bounded by $1$). Both are bounded martingales, so
$t\mapsto\E(Y'_t)^2$ is nondecreasing on $[0,\underline q)$.

\emph{Convexity of $\Gz{\zP}$ on the gap.} On $[0,\underline q)\times\R$ the function $u$ is
smooth and solves $\partial_tu+\frac12\xi''\partial_{xx}u=0$, so It\^o's formula gives
$dY_t=\sqrt{\xi''(t)}\,Y'_t\,dB_t$. Since $Y_0=u(0,0)=0$ by oddness, the It\^o isometry yields
\[
 \Gz{\zP}(t)=\E Y_t^2=\int_0^t\xi''(r)\,\E(Y'_r)^2\,dr,\qquad 0\le t<\underline q .
\]
The integrand is the product of two nonnegative nondecreasing functions, because
$\xi'''\ge0$ on $[0,1]$ (the coefficients of $\xi$ are nonnegative). So $\Gz{\zP}$ is
convex on $[0,\underline q)$, and on $[0,\underline q]$ by continuity, with
$\Gz{\zP}(0)=0$ and $\Gz{\zP}(\underline q)=\underline q$ by \eqref{eq:eq-gap-end}.
Therefore $\Gz{\zP}(t)\le t$ on $[0,\underline q]$, that is,
$\OP'(t)=-\frac12\xi''(t)(\Gz{\zP}(t)-t)\ge0$ there, and
$\OP(0)\le\OP(\underline q)=\min\OP$.
\end{proof}

Lemmas~\ref{lem:eq-obstacle} and~\ref{lem:eq-zero} prove the first assertion of
Corollary~\ref{cor:eq}: $\CP$ contains $0$ and $\supp\zP$. The stronger property
$0\in\supp\zP$ is \cite[Theorem~1(i)]{AC15a}, stated there for $\sum_p2^p\beta_p^2<\infty$.
Lemma~\ref{lem:eq-zero} asks less, holds under \condH, and is exactly what the proof below uses.
The convexity of $\Gz{\zP}$ on an initial gap is also the mechanism of the proof of
\cite[Theorem~1(i)]{AC15a}.

The third lemma, the equilibrium counterpart of Proposition~\ref{prop:loc-gap}, concerns glued
data (Section~\ref{sec:loc-glued}): a cascade coefficient $\xc$ on $[0,Q)$ and a band measure
$\eta\in\Pspace([Q,1])$, with glued coefficient $\glue(\xc,\eta)$ equal to $\xc$ on $[0,Q)$ and
to $\cdf\eta$ on $[Q,1]$.

\begin{lemma}[The variational gap at equilibrium]\label{lem:eq-gap}
Assume \condH. Let $Q\in(0,1)\setminus\CP$ and $\varepsilon\in(0,1]$, and take the glued data
\[
 \xc=\cdf{\zP}+\varepsilon\ \text{ on }[0,Q),\qquad \eta=\Shift_Q\zP\in\Pspace([Q,1]).
\]
Then $\xc$ is nondecreasing with values in $[\varepsilon,2]$, and
\begin{equation}\label{eq:eq-gap}
 P\bigl(\glue(\xc,\eta)\bigr)\le\feq-\varepsilon\bigl(\OP(Q)-\min\OP\bigr)+K\varepsilon^2 .
\end{equation}
\end{lemma}

\begin{proof}
The measure $\Shift_Q\zP=\zP|_{(Q,1]}+\zP([0,Q])\delta_Q$ has distribution function
$\cdf{\zP}$ on $[Q,1]$, so $\glue(\xc,\eta)=\cdf{\zP}+\varepsilon h$ with
$h=\ind_{[0,Q)}$. Here $\|h\|_{L^1}=Q\le1$ and the glued coefficient takes values in
$[0,2]$. By \eqref{eq:eq-second}, \eqref{eq:eq-firstvar} and Lemma~\ref{lem:eq-zero},
\[
 P\bigl(\glue(\xc,\eta)\bigr)\le\feq+\varepsilon\bigl(\OP(0)-\OP(Q)\bigr)+K\varepsilon^2
 =\feq-\varepsilon\bigl(\OP(Q)-\min\OP\bigr)+K\varepsilon^2 .
\]
The claims about $\xc$ are clear.
\end{proof}

Lemma~\ref{lem:eq-gap} is the perturbation used in the proof of Proposition~\ref{prop:loc-gap}
above a contact point, here the contact point $0$, which lies below every radius; the upward
perturbation of direction (B), which uses $\zf(\{0\})=\thf$, is not needed. As
$Q\notin\CP\supset\supp\zP$, $\zP$ has no atom at $Q$, and the glued coefficient jumps down by
$\varepsilon$ at $Q$. The gain comes from this jump: by the argument of
Lemma~\ref{lem:loc-monotone}, a nondecreasing glued coefficient is at most $\cdf\eta(Q)\le1$, hence
a distribution function, with Parisi functional at least $\feq$. The cascade coefficient is at
least $\varepsilon$ whether or not $\zP$ charges $0$, which leaves room for a positive tilt
below it.

\subsection{\texorpdfstring{Proof of Corollary~\ref*{cor:eq}}{Proof of Corollary 6}}%
\label{sec:eq-proof}

\begin{proof}[Proof of Corollary~\ref{cor:eq}]
The containment $\{0\}\cup\supp\zP\subset\CP$ was proved in Section~\ref{sec:eq-obstacle}. The
probability bound has a localization part and a field-law part. We prove them in Steps~1 and~2
by the shell argument of Section~\ref{sec:loc-proof} and Appendix~\ref{app:var-field}, with
Lemma~\ref{lem:eq-gap} as the source of the gain. Step~3 combines them, and Step~4 deduces the
fixed-radius statement. Convexity of $\xi$ on $[-1,1]$ enters only through the Slepian step of
Proposition~\ref{prop:glued-comparison}; the toolbox of Appendix~\ref{app:pde-toolbox} and
the proofs of Lemmas~\ref{lem:eq-obstacle} and~\ref{lem:eq-zero} use only \condR.

\emph{Shells and the comparison.} For $\rho\in(0,1]$ let $Q_j=j\rho$ and
$I_j=[Q_j,\min(Q_j+\rho,1)]$ for the integers $j\ge0$ with $Q_j<1$. There are at most
$\lceil1/\rho\rceil+1$ such shells, and they cover $[0,1]$. We use \eqref{eq:loc-glued-general}
with $s=0$ (Proposition~\ref{prop:glued-comparison} with the approximation of
Corollary~\ref{cor:pde-majorant}): if $\theta\in(0,1)$, $Q\in[0,1)$, $\rho'\in(0,1-Q]$,
$\xc:[0,Q)\to[\theta,2]$ is nondecreasing (there is no cascade part if $Q=0$) and
$\eta\in\Pspace([Q,1])$, then for every $\varepsilon_A>0$ there is a random
variable $\mathcal T$ such that
\begin{equation}\label{eq:eq-majorant}
\begin{aligned}
 \sup\bigl\{\FTAP(m):\ m\in\cube,\ q_m\in[Q,Q+\rho']\bigr\}
 &\le N P\bigl(\glue(\xc,\eta)\bigr)+N\varepsilon_A+\mathcal T,\\
 \E e^{\theta\mathcal T}&\le(N+1)\,e^{N\theta C_R\rho'},
\end{aligned}
\end{equation}
with $C_R=5\xi''(1)$. The events below that involve all $m$ in the cube are contained in finite
unions of events $\{\mathcal T\ge a\}$, which is how their outer probabilities are bounded.

\emph{Step 1: localization.} Fix $\kappa>0$ and let
$A_\kappa=\{Q\in[0,1]:\dist(Q,\CP)\ge\kappa/2\}$. If $A_\kappa$ is empty, no $m$ has $\dist(q_m,\CP)\ge\kappa$, so the event of
\eqref{eq:eq-step1} is empty and there is nothing to prove; set $c_1:=c_2$ in Step~3. Otherwise $A_\kappa$ is compact and disjoint from $\CP$, and
$\OP$ is continuous, so
\[
 \gamma:=\min\bigl\{\OP(Q)-\min\OP:\ Q\in A_\kappa\bigr\}>0 .
\]
Put $\varepsilon=\min\{1,\gamma/(2K)\}$ and $\theta=\varepsilon/2$, choose $\rho\in(0,\kappa/2)$
with $C_R\rho\le\varepsilon\gamma/16$, and put $\varepsilon_A=\varepsilon\gamma/16$ and
$t=\varepsilon\gamma/8$.

If $\dist(q_m,\CP)\ge\kappa$, then $q_m\in I_j$ for some $j$ with
$\dist(Q_j,\CP)\ge\kappa-\rho>\kappa/2$. Such an anchor lies in $A_\kappa$; it satisfies
$Q_j>\kappa/2>0$ because $0\in\CP$ (Lemma~\ref{lem:eq-zero}), and $Q_j<1$. So
Lemma~\ref{lem:eq-gap} applies at $Q=Q_j$, and since $\OP(Q_j)-\min\OP\ge\gamma$ and
$K\varepsilon\le\gamma/2$, its data satisfy
\[
 P\bigl(\glue(\xc,\eta)\bigr)\le\feq-\varepsilon\gamma+K\varepsilon^2\le\feq-\varepsilon\gamma/2 .
\]
The cascade coefficient takes values in $[\varepsilon,2]\subset[\theta,2]$, so
\eqref{eq:eq-majorant} applies on $I_j$ and gives a random variable $\mathcal T_j$ with
$\E e^{\theta\mathcal T_j}\le(N+1)e^{N\theta C_R\rho}$ and
\[
 \bigl\{\exists m:\ q_m\in I_j,\ \FTAP(m)\ge N(\feq-t)\bigr\}
 \subset\bigl\{\mathcal T_j\ge N(\varepsilon\gamma/2-\varepsilon_A-t)\bigr\}
 =\bigl\{\mathcal T_j\ge\tfrac5{16}N\varepsilon\gamma\bigr\}.
\]
By Markov's inequality applied to $e^{\theta\mathcal T_j}$, this event has probability at most
\[
 (N+1)\exp\Bigl\{N\theta\Bigl(C_R\rho-\frac5{16}\varepsilon\gamma\Bigr)\Bigr\}
 \le(N+1)\,e^{-N\theta\varepsilon\gamma/4}.
\]
A union bound over the at most $\lceil1/\rho\rceil+1$ shells gives, for all large $N$,
\begin{equation}\label{eq:eq-step1}
 \Pout\bigl(\exists m\in\cube:\ \FTAP(m)\ge N(\feq-t),\ \dist(q_m,\CP)\ge\kappa\bigr)
 \le e^{-c_1N},
\end{equation}
where $c_1:=\theta\varepsilon\gamma/8=\varepsilon^2\gamma/16\le t$.

\emph{Step 2: field law.} Fix $\kappa>0$. Let $\tilde g_1,\dots,\tilde g_J$ be the smooth test
functions of Lemma~\ref{lem:var-net} for the accuracy $\kappa$: they satisfy
$\|\tilde g_j'\|_\infty\le1$, $\|\tilde g_j''\|_\infty<\infty$, and for probability measures
$\mu,\mu'$ on $[-1,1]$,
\begin{equation}\label{eq:eq-net}
 \Wone(\mu,\mu')\ge\kappa\quad\Longrightarrow\quad
 \int\tilde g_j\,d(\mu-\mu')\ge\kappa/2\ \text{ for some } j .
\end{equation}
The map $q\mapsto\int\tilde g_j\,d\LPq{q}$ is uniformly continuous on $[0,1]$, because
$q\mapsto u(q,X_q)$ is almost surely continuous and bounded by $1$
(Lemma~\ref{lem:pde-measurable}); let $\omega_{\mathrm P}$ be a common modulus of continuity for
$j=1,\dots,J$.

We need the perturbed comparison at a base point other than a minimizer. Let
$\zeta'\in\Pspace([0,1])$ and $\theta\in(0,1)$ with $\zeta'(\{0\})\ge\theta$,
let $L'_Q=\Law(\partial_x\Phi_{\zeta'}(Q,X'_Q))$ be the field law of $\zeta'$ ($X'$ its
Auffinger--Chen process), and let $Q\in[0,1)$, $\rho'\in(0,1-Q]$,
$0<s\le(2\max_j\|\tilde g_j''\|_\infty)^{-1}$ and $\varepsilon_A>0$. Then for each $j$ there is a
random variable $\mathcal T$ with $\E e^{\theta\mathcal T}\le(N+1)e^{N\theta C_R\rho'}$ and
\begin{equation}\label{eq:eq-perturbed}
 \sup\Bigl\{\FTAP(m)+s\sum_{i=1}^N\tilde g_j(m_i):\ q_m\in[Q,Q+\rho']\Bigr\}
 \le N\bigl[P_{\zeta'}+s\int\tilde g_j\,dL'_Q+C_2s^2+\varepsilon_A\bigr]+\mathcal T,
\end{equation}
where $C_2=1+L_\psi$, with $L_\psi$ as in~\eqref{eq:pde-constants}, depends only on $\xi'(1)$
and $\xi''(1)$. This is
Lemma~\ref{lem:var-tangent}\ref{lem:var-tangent-b} with $g=\tilde g_j$: the comparison of
Proposition~\ref{prop:glued-comparison} at the base point $(\cdf{\zeta'}|_{[0,Q)},\Shift_Q\zeta')$,
with the terminal condition perturbed by $s\tilde g_j$. Its proof uses only
$\zeta'(\{0\})\ge\theta$, which makes this base point admissible, and not any minimality of
$\zeta'$; so it applies to the measure $\zeta'$ chosen below.

Let $L_M$ be the Lipschitz constant of the field law in the coefficient
(Corollary~\ref{cor:pde-step}\ref{cor:pde-step-d} and Lemma~\ref{lem:pde-measurable}): under the
coupling by a common Brownian motion,
\[
 \E\bigl|\partial_x\Phi_{\zeta_1}(t,X^{(1)}_t)-\partial_x\Phi_{\zeta_2}(t,X^{(2)}_t)\bigr|\le
 L_M\|\cdf{\zeta_1}-\cdf{\zeta_2}\|_{L^1}.
\]
Choose, in this order,
\[
 s=\min\Bigl\{\frac{\kappa}{8C_2},\ \bigl(2\max_j\|\tilde g_j''\|_\infty\bigr)^{-1}\Bigr\},
 \qquad
 \varepsilon'\in(0,1]\ \text{with}\ K\varepsilon'^2+sL_M\varepsilon'\le\frac{s\kappa}{16},
\]
then $\rho\in(0,1]$ with $C_R\rho\le s\kappa/16$ and $\omega_{\mathrm P}(\rho)\le\kappa/16$, and
put $\varepsilon_A=s\kappa/32$, $t=s\kappa/16$ and $\theta=\varepsilon'/2$. The base point is
\[
 \zeta'=(1-\varepsilon')\zP+\varepsilon'\delta_0,\qquad
 \cdf{\zeta'}=\cdf{\zP}+\varepsilon'(1-\cdf{\zP}),\qquad
 \zeta'(\{0\})\ge\varepsilon'>\theta .
\]
It costs little, because the first variation in its direction vanishes. With
$h=1-\cdf{\zP}=\cdf{\delta_0}-\cdf{\zP}$, \eqref{eq:eq-firstvar}, \eqref{eq:eq-fubini},
Lemma~\ref{lem:eq-zero} and Lemma~\ref{lem:eq-obstacle}\ref{lem:eq-obstacle-a} give
$DP(\cdf{\zP})[h]=\OP(0)-\int\OP\,d\zP=\min\OP-\min\OP=0$, so by \eqref{eq:eq-second}
\begin{equation}\label{eq:eq-basecost}
 P_{\zeta'}\le\feq+K\varepsilon'^2 .
\end{equation}
Since $\|\cdf{\zeta'}-\cdf{\zP}\|_{L^1}\le\varepsilon'$ and $\|\tilde g_j'\|_\infty\le1$,
\begin{equation}\label{eq:eq-lawshift}
 \Bigl|\int\tilde g_j\,dL'_Q-\int\tilde g_j\,d\LPq{Q}\Bigr|\le L_M\varepsilon'
 \qquad\text{for every } Q\in[0,1] .
\end{equation}

Now let $m$ satisfy $q_m\in I_i$, $\FTAP(m)\ge N(\feq-t)$ and
$\Wone(\mu_m,\LPq{q_m})\ge\kappa$. By \eqref{eq:eq-net} there is $j$ with
$\int\tilde g_j\,d\mu_m\ge\int\tilde g_j\,d\LPq{q_m}+\kappa/2
\ge\int\tilde g_j\,d\LPq{Q_i}-\omega_{\mathrm P}(\rho)+\kappa/2$, and since
$\sum_k\tilde g_j(m_k)=N\int\tilde g_j\,d\mu_m$,
\[
 \FTAP(m)+s\sum_{k=1}^N\tilde g_j(m_k)
 \ge N\Bigl[\feq-t+s\int\tilde g_j\,d\LPq{Q_i}-s\,\omega_{\mathrm P}(\rho)+\frac{s\kappa}2\Bigr].
\]
Compare with \eqref{eq:eq-perturbed} on $I_i$, bounded through \eqref{eq:eq-basecost} and
\eqref{eq:eq-lawshift}: the corresponding random variable $\mathcal T_{ij}$ satisfies
\[
 \mathcal T_{ij}\ge N\Bigl[\frac{s\kappa}2-t-s\,\omega_{\mathrm P}(\rho)-K\varepsilon'^2-sL_M\varepsilon'
 -C_2s^2-\varepsilon_A\Bigr]
 \ge N s\kappa\Bigl(\frac{16}{32}-\frac2{32}-\frac2{32}-\frac2{32}-\frac4{32}-\frac1{32}\Bigr)
 =\frac5{32}Ns\kappa ,
\]
using $C_2s^2\le s\kappa/8$. Markov's inequality, with $C_R\rho\le\frac2{32}s\kappa$, bounds the
probability of this event by $(N+1)e^{-\frac3{32}N\theta s\kappa}$. A union bound over the
finitely many pairs $(i,j)$ gives, for all large $N$,
\begin{equation}\label{eq:eq-step2}
 \Pout\bigl(\exists m\in\cube:\ \FTAP(m)\ge N(\feq-t),\ \Wone(\mu_m,\LPq{q_m})\ge\kappa\bigr)
 \le e^{-c_2N},
\end{equation}
where $c_2:=\theta s\kappa/16=\varepsilon's\kappa/32\le t$.

\emph{Step 3: conclusion.} Let $c=\min\{c_1,c_2\}/2$. In each of Steps~1 and~2 the constant
$c_k$ does not exceed the threshold parameter $t$ of that step, so the event in
Corollary~\ref{cor:eq} at the threshold $N(\feq-c)$ is contained in the union of the events of
\eqref{eq:eq-step1} and \eqref{eq:eq-step2}. Its outer
probability is therefore at most $e^{-c_1N}+e^{-c_2N}\le2e^{-2cN}\le e^{-cN}$ for all large $N$.

\emph{Step 4: a fixed radius.} Let $Q\in[0,1]\setminus\CP$ and apply the bound of Step~3 with
$\kappa=\dist(Q,\CP)$, which is positive because $\CP$ is closed; call the resulting constant
$c(Q)$. If $\sup\{\FTAP(m):q_m=Q\}>N(\feq-c(Q))$, some $m$ with $q_m=Q$ has
$\FTAP(m)\ge N(\feq-c(Q))$ and $\dist(q_m,\CP)=\kappa$, so this event has outer probability at
most $e^{-c(Q)N}$. It is measurable, since $\FTAP-\HN$ is deterministic and the supremum is a
$1$-Lipschitz function of $\HN$ in the uniform norm. So for all large $N$,
$\sup\{\FTAP(m):q_m=Q\}\le N(\feq-c(Q))$ with probability at least $1-e^{-c(Q)N}$.
\end{proof}

\subsection{Relation to the TAP representation}\label{sec:eq-remarks}

\begin{remark}[Two consistency checks]\label{rem:eq-checks}
\emph{Replica symmetry at high temperature.} Let $\xi''(1)<1$. For $\zeta=\delta_0$ the
coefficient is $1$, $\Phi_{\delta_0}(t,x)=\log\cosh x+\frac12(\xi'(1)-\xi'(t))$, and $\tanh X_t$
is a martingale with $d\tanh X_t=\sqrt{\xi''}\sech^2X_t\,dB_t$. Hence
$\Gz{\delta_0}(t)=\int_0^t\xi''\E\sech^4X_r\,dr\le\xi'(t)\le\xi''(1)t<t$ for $t\in(0,1]$, so
$\Obs_{\delta_0}$ is strictly increasing on $[0,1]$. By \eqref{eq:loc-first-order} and
\eqref{eq:loc-moving-mass}, $P_\mu-P_{\delta_0}\ge\int\Obs_{\delta_0}\,d(\mu-\delta_0)\ge0$ for
every $\mu\in\Pspace([0,1])$, so $\delta_0$ is a Parisi measure, with $\CP=\{0\}$ and
$\LPq{0}=\delta_0$. Corollary~\ref{cor:eq} then says that, outside an event of outer probability at
most $e^{-cN}$, every $m$ with
$\FTAP(m)\ge N(\feq-c)$ has $q_m<\kappa$, so that most of its coordinates are small, in agreement
with \eqref{eq:eq-tapzero}.

\emph{One step of replica symmetry breaking.} If $\zP=a\delta_0+(1-a)\delta_q$ with
$a\in(0,1)$, Lemma~\ref{lem:eq-obstacle} gives $\CP\supset\{0,q\}$ and Lemma~\ref{lem:eq-zero}
adds nothing. Corollary~\ref{cor:eq} localizes near $\CP$, consistently with \eqref{eq:eq-tapzero}
at $0$ and, for $\xi\in\classC$, with \cite[Theorem~2]{CPS} at $q$.
\end{remark}

\begin{remark}[Relation to the TAP representation]\label{rem:eq-cps}
\emph{The representation theorems.} For $\xi\in\classC$, the TAP representation at equilibrium,
\cite[Theorem~2]{CPS} and \cite[Theorems~2, 3 and~5]{CPS2}, concerns radii $q\in\supp\zP$; at
every radius the maximum of the TAP functional on the sphere is at most the free energy up to
$o(N)$, with probability tending to one \cite[Section~5]{CPS}. Corollary~\ref{cor:eq} makes this
comparison strict off $\CP\supset\supp\zP$: there the maximum is smaller by a linear amount,
with an exponential bound uniform over the cube. At massless contacts, radii in
$\CP\setminus\supp\zP$, it does not decide whether the representation holds, although its field-law
clause applies there (Section~\ref{sec:open}). Section~\ref{sec:intro-results}
compares it with Subag's characterization for spherical models \cite[Theorems~5 and~10]{Sub24}.

\emph{An exact formula at every radius.} For $\xi\in\classC$ and every radius $Q$,
\cite[Theorem~36]{CPS}, which Chen, Panchenko and Subag import from \cite[Theorem~5]{CP18},
gives an exact variational formula, over a
Lagrange multiplier $\lambda\in\R$ and a coefficient on $[0,Q]$, for the limit of the expected
maximum over $\{q_m=Q\}$ of the TAP functional with its order parameter frozen at $\zP$; this
maximum bounds that of $\FTAP$ on the sphere from above. They use the formula at $\lambda=0$
with the coefficient $\cdf{\zP}$ on $[0,Q]$. Moving $\lambda$ away from $0$ changes the formula
at first order by $\Gz{\zP}(Q)-Q$ per unit of $\lambda$; moving the coefficient toward that
of $\delta_0$ changes it by $(1-\zP([0,Q)))(\min\OP-\OP(Q))$, which is negative for
$Q\notin\CP$ below $\max\supp\zP$ by Lemmas~\ref{lem:eq-obstacle}\ref{lem:eq-obstacle-a}
and~\ref{lem:eq-zero}. With the Borell--TIS inequality, as in \cite[Section~5]{CPS}, and the
normalization shift of Section~\ref{sec:intro-setting}, this gives, for $\xi\in\classC$, the
fixed-radius clause of
Corollary~\ref{cor:eq} in a few lines at every radius with $\Gz{\zP}(Q)\ne Q$ and at every
$Q\notin\CP$ below $\max\supp\zP$. Neither these first-order moves nor the field-law clause
alone reach the remaining radii off $\CP$, those above $\max\supp\zP$ with
$\Gz{\zP}(Q)=Q$, if there are any: since $\int y^2\,\mu_m(dy)=q_m$,
$\int y^2\,\LPq{q}(dy)=\Gz{\zP}(q)$ and $y\mapsto y^2$ is $2$-Lipschitz on $[-1,1]$, that
clause gives only $|\Gz{\zP}(q_m)-q_m|\le2\Wone(\mu_m,\LPq{q_m})<2\kappa$, which places
$q_m$ near the set where $\Gz{\zP}(q)=q$, that is, where $\OP'$ vanishes. Excluding
stationary points of $\OP$ that are not minima is the work of the glued comparison, through the
downward jump of Lemma~\ref{lem:eq-gap}, which at those radii goes from $1+\varepsilon$ to $1$.

\emph{The field law.} In our normalization $\LPq{q}$ is the law that \cite[(1.29)]{CPS2},
recalling \cite[Chapter~4]{Pan13} and \cite{AJ19a,AJ19b}, attaches informally to ancestor states
at $q\in\supp\zP$. The field-law clause of Corollary~\ref{cor:eq} is a precise statement of this
kind for every $m$ with $\FTAP(m)\ge N(\feq-c)$, at any radius; it concerns empirical laws of TAP
near-maximizers, not the Gibbs spin distributions of \cite{AJ19b}.
\end{remark}

%% file: sections/sec7-open.tex
\section{Open problems}\label{sec:open}\label{sec:intro-open}

Several questions remain open.
\begin{itemize}
\item \emph{Outliers.} Theorem~\ref{thm:struct}\ref{thm:struct-iii} controls the bulk of the
 Hessian, not its largest eigenvalues; \cite[Remark~B.1]{Bou26b} also notes that the outliers
 require a separate argument.
\item \emph{Massless contacts.} Corollary~\ref{cor:eq} does not decide whether the massless contacts
 $Q\in\CP\setminus\supp\zP$ carry the TAP representation, although its field-law clause applies
 there. Does $N^{-1}\sup\{\FTAP(m):q_m=Q\}$ converge in probability to $\feq$ there? Boursier
 indicates a construction when the contact set is finite, under a strict stability assumption
 \cite[abstract and Section~1.3]{Bou26b}; we know of no written proof. For spherical models the
 representation holds: every point
 of Talagrand's set $S_\beta$, which contains $\supp\mu_\beta$ \cite[Proposition~9]{Sub24}, is
 multi-samplable \cite[Theorem~10]{Sub24}, hence carries the representation \cite[Theorem~5]{Sub24}, and
 massless contacts occur, for instance in the $2$-spin model at low temperature
 \cite[Remark~13]{Sub24}. For SK at low temperature ($\xi''(0)>1$), $\supp\zP=[0,q_{\mathrm P}]$
 with $q_{\mathrm P}>0$ by the recent preprint \cite[Theorem~1.1]{Lop26} (near criticality
 \cite{Zho26a}), so massless contacts could only lie in $(q_{\mathrm P},1)$. The analogue at an
 upper level $f$: is the rate $\Sigma(f)$ attained on spheres of every radius in $\Cf$, beyond the
 first contact $\qf$, near which the Kac--Rice construction attains it
 (Proposition~\ref{prop:kr-main})?
\item \emph{Marginality beyond the $2+4$ family.} For SK at low temperature the first contact is
 marginal at all levels close to $\feq$ (the companion paper \cite{PeiB}, not used here); whether
 every contact is marginal there, so that Theorem~\ref{thm:stab}\ref{thm:stab-b} applies, is not
 addressed here.
\item \emph{Other models.} An external field. Non-convex mixtures: the Kac--Rice lower bound
 extends under a sign condition on $\xi$, provided the variational inputs hold
 (Remark~\ref{rem:kr-window}), while the matching upper bound would need a step that is only
 sketched in \cite{CGKLM}. Whether $\xi\in\classC$ can be replaced by \condR\ throughout; by
 Remark~\ref{rem:rpde} it suffices that the facts on the Parisi equation and the TAP correction
 that we cite hold under \condR.
\item \emph{Counting.} The ordinary annealed and quenched counts of TAP critical points at level $f$,
 for which \cite[Conjectures~1 and~2]{Bou26a} proposed Legendre formulas that
 \cite[Section~1.2 and Appendix~C]{Bou26b} corrects (the annealed count of SUSY states has rate
 $\Sigma(f)$, Proposition~\ref{prop:stab-susy}; at upper levels Boursier reads the supersymmetric
 quenched complexity, heuristically, as the probability of existence \cite[Section~1.2]{Bou26b}, which
 Theorem~\ref{thm:exist}\ref{thm:exist-d} computes),
 and their organization \cite[Conjecture~3]{Bou26a}.
\item \emph{Hessians at the equilibrium level.} We prove no Hessian statement at $\feq$:
 Theorem~\ref{thm:struct}\ref{thm:struct-iii}--\ref{thm:struct-iv} concerns levels
 $f\in(\feq,\fone)$ only.
\end{itemize}

%% file: sections/appA-parisi.tex
\section{The Parisi functional with non-monotone coefficients}\label{app:pde}

The glued coefficient of Proposition~\ref{prop:glued-comparison} is a cascade coefficient below
the radius, which may exceed one, and a distribution function above it, with no ordering at the
radius; the Ruelle cascade in its proof has these values as parameters. Parisi functionals with such
coefficients also appear in the extended variational principle of \cite{EAMS21} and in the bounds
of \cite[Section~6]{HS25}. Appendix~\ref{app:pde-toolbox} develops their calculus, used in
Section~\ref{sec:loc} and Appendix~\ref{app:glued}; Appendix~\ref{app:var-constrained} proves the
regularity of every level (Theorem~\ref{thm:exist}\ref{thm:exist-b}) and
Lemma~\ref{lem:loc-obstacle}; Appendix~\ref{app:var-gap} completes Proposition~\ref{prop:loc-gap};
and Appendix~\ref{app:pde-standing} proves the regularity at the contact set used in
Section~\ref{sec:hess} and Appendices~\ref{app:matrix}--\ref{app:nc}.

\paragraph{Standing assumption.}
Throughout, $\xi$ satisfies \condR. Appendix~\ref{app:pde-toolbox} uses $\xi$ only through the
following facts on $[0,1]$, which hold because $\xi\not\equiv0$ has nonnegative coefficients:
$\xi'(0)=0$; $\xi''$ and $\xi'''$ are continuous and nonnegative; $\xi''>0$ on $(0,1]$; so $\xi'$ is
strictly increasing and $\xi''\le\xi''(1)$. Among the facts on the Parisi equation, only
Corollary~\ref{cor:pde-majorant} (Appendix~\ref{app:glued-approx}) assumes \condH, since convexity
of $\xi$ on $[-1,1]$ enters only through Proposition~\ref{prop:glued-comparison}, by way of
Slepian's lemma. In Appendices~\ref{app:var-constrained}--\ref{app:var-gap},
Lemma~\ref{lem:loc-obstacle}, Remark~\ref{rem:var-level} and Proposition~\ref{prop:loc-gap} assume
\condH, as Section~\ref{sec:loc} does. Appendix~\ref{app:pde-standing} assumes in addition that
$\xi\in\classC$.

\subsection{Cole--Hopf solutions, sandwich inequalities and stability}\label{app:pde-toolbox}

Everything rests on the verification identity of Lemma~\ref{lem:pde-ac}, which compares the
Cole--Hopf solution with any controlled diffusion; it yields comparisons in the coefficient and in
the terminal value, and from these the calculus of the Parisi functional
(Corollary~\ref{cor:pde-step}). Fix $G\ge1$ (in the applications $G=3$) and an interval
$[t_0,t_1]\subset[0,1]$.

\paragraph{Coefficients, terminal values and Cole--Hopf solutions.}
A \emph{step coefficient} on $[t_0,t_1)$ is a function
$c=\sum_{j=1}^Jc_j\ind_{[\tau_{j-1},\tau_j)}$ with $t_0=\tau_0<\tau_1<\dots<\tau_J=t_1$ and
$c_j\in[0,G]$; no monotonicity is assumed. The \emph{terminal class} is
\[
 \mathcal T_G=\bigl\{\psi\in C^2(\R)\ \text{convex}:\ |\psi'|<1,\ \psi''\le G(1-\psi'^2)\bigr\}.
\]
Thus $\mathcal T_G\subset\mathcal T_{G'}$ for $G\le G'$, and $\log\cosh\in\mathcal T_1$. For a step
coefficient $c$ and $\psi\in\mathcal T_G$, the \emph{Cole--Hopf solution} $W=W[c,\psi]$ on
$[t_0,t_1]\times\R$ is defined backward: $W(t_1,\cdot)=\psi$, and for $t\in[\tau_{j-1},\tau_j)$,
with $\sigma_j(t)^2=\xi'(\tau_j)-\xi'(t)>0$ and $Z$ a standard Gaussian variable,
\begin{equation}\label{eq:pde-colehopf}
 W(t,y)=\begin{cases}
 \dfrac1{c_j}\log\E\exp\bigl(c_jW(\tau_j,y+\sigma_j(t)Z)\bigr) & \text{if } c_j>0,\\[1ex]
 \E\,W(\tau_j,y+\sigma_j(t)Z) & \text{if } c_j=0 .
 \end{cases}
\end{equation}
The expectations are finite because $W(\tau_j,\cdot)$ is $1$-Lipschitz, by
Lemma~\ref{lem:pde-class}\ref{lem:pde-class-d} below applied on the later steps. We write
$v=\partial_yW$. When $[t_0,t_1]=[0,1]$, $\psi=\log\cosh$ and $c=\cdf\zeta$ for a finitely
supported $\zeta\in\Pspace([0,1])$, \eqref{eq:pde-colehopf} is the classical Cole--Hopf solution
$\Phi_\zeta$ of the Parisi equation \cite[(11)--(12)]{AC15a}, and $v=\Phi_{\zeta,x}$.

\begin{lemma}[Class preservation and regularity]\label{lem:pde-class}
Let $c$ be a step coefficient on $[t_0,t_1)$ with values in $[0,G]$, let $\psi\in\mathcal T_G$ and
$W=W[c,\psi]$.
\begin{thmlist}
\item\label{lem:pde-class-a} $W(t,\cdot)\in\mathcal T_G$ for every $t\in[t_0,t_1]$; in particular
 $|v|<1$ and $0\le\partial_yv\le G(1-v^2)\le G$. If $\psi''>0$ on $\R$, then $\partial_yv>0$.
\item\label{lem:pde-class-b} On each closed step $[\tau_{j-1},\tau_j]$, $W$ is of class $C^{1,2}$,
 with one-sided $t$-derivatives at the ends, and
 $\partial_tW+\frac{\xi''}2\bigl(\partial_{yy}W+c_j(\partial_yW)^2\bigr)=0$.
\item\label{lem:pde-class-c} On each open step $v$ is of class $C^{1,2}$ and
 $\partial_tv+\frac{\xi''}2\partial_{yy}v+\xi''c_j\,v\,\partial_yv=0$. Moreover, $v$ is continuous
 on $[t_0,t_1]\times\R$.
\item\label{lem:pde-class-d} $|W(t,y)-W(t,0)|\le|y|$.
\end{thmlist}
\end{lemma}

\begin{proof}
Work on one step $[\tau_{j-1},\tau_j)$, write $c=c_j$, $\sigma=\sigma_j$ and
$\psi_j=W(\tau_j,\cdot)$, and assume $\psi_j\in\mathcal T_G$. This holds for $j=J$, and once (a)
is proved on the step it holds for $j-1$; so backward induction covers every step. Since
$|\psi_j'|<1$, we have $|\psi_j(z)-\psi_j(0)|\le|z|$.

\emph{The case $c>0$.} Let $p_s(z)=(2\pi s)^{-1/2}e^{-z^2/(2s)}$ and
\[
 \Theta(t,y)=\E e^{c\psi_j(y+\sigma(t)Z)}=\int e^{c\psi_j(z)}p_{\sigma(t)^2}(z-y)\,dz .
\]
The integrand is at most $e^{c(\psi_j(0)+|z|)}p_{\sigma(t)^2}(z-y)$, and $\sigma(t)>0$ for
$t<\tau_j$. So differentiation under the integral is legitimate on
$[\tau_{j-1},\tau_j)\times\R$: there $\Theta$ is positive, $C^1$ in $t$ and $C^\infty$ in $y$, and,
since $\partial_sp_s=\frac12\partial_{zz}p_s$ and $\frac{d}{dt}\sigma(t)^2=-\xi''(t)$, it solves
$\partial_t\Theta=-\frac{\xi''}2\partial_{yy}\Theta$. Hence $W=c^{-1}\log\Theta$ satisfies
$\partial_tW=-\frac{\xi''}2(\partial_{yy}W+c(\partial_yW)^2)$ on the half-open step.
Differentiating $\Theta(t,y)=\E e^{c\psi_j(y+\sigma Z)}$ in $y$ under the expectation, which is
legitimate because $\psi_j'$ and $\psi_j''$ are bounded and continuous, gives
\begin{equation}\label{eq:pde-tilt}
 \partial_yW(t,y)=\E_{\mu}\psi_j'(y+\sigma\,\cdot),\qquad
 \partial_{yy}W(t,y)=\E_{\mu}\psi_j''(y+\sigma\,\cdot)+c\Var_{\mu}\psi_j'(y+\sigma\,\cdot),
\end{equation}
where $\mu=\mu_{t,y}$ is the probability measure on $\R$ whose density with respect to the
standard Gaussian measure is proportional to $e^{c\psi_j(y+\sigma(t)z)}$.

\emph{Continuity at $t=\tau_j$.} This density equals
$e^{c(\psi_j(y+\sigma z)-\psi_j(y))}/\E e^{c(\psi_j(y+\sigma Z)-\psi_j(y))}$, whose numerator is at
most $e^{c\sigma|z|}$ and whose denominator is at least $\E e^{-c\sigma|Z|}\ge e^{-G}$ for
$\sigma\le1$; so these densities are uniformly integrable for $\sigma\le1$. Since
$\psi_j,\psi_j',\psi_j''$ are continuous and $\psi_j',\psi_j''$ are bounded, dominated convergence
gives $W\to\psi_j$, $\partial_yW\to\psi_j'$ and $\partial_{yy}W\to\psi_j''$ as
$(t,y')\to(\tau_j,y)$, locally uniformly in $y$. By the equation, $\partial_tW$ extends
continuously to $t=\tau_j$, and since $W(\tau_j,y)-W(t,y)=\int_t^{\tau_j}\partial_tW(s,y)\,ds$,
the left $t$-derivative at $\tau_j$ exists and equals this extension. This proves (b) when $c>0$.

\emph{The case $c=0$.} Then $W(t,y)=\E\psi_j(y+\sigma(t)Z)$ is a heat flow, and the same argument
without the tilt gives (b), with $\partial_tW=-\frac{\xi''}2\partial_{yy}W$; the formulas
\eqref{eq:pde-tilt} hold with $c=0$ and $\mu$ the standard Gaussian measure.

\emph{(c).} On the open step $W$ is $C^\infty$ in $y$, and the $t$-derivatives of its
$y$-derivatives are continuous (differentiate $\Theta$ under the integral). Differentiating the
equation in $y$ gives the equation for $v$. At $\tau_j$ the left limit of $v$ is
$\psi_j'=v(\tau_j,\cdot)$, which is the value used on the next step; so $v$ is continuous on
$[t_0,t_1]\times\R$.

\emph{(a).} By \eqref{eq:pde-tilt}, $\psi_j''\le G(1-\psi_j'^2)$, $c\le G$ and Jensen's inequality,
\[
 \partial_{yy}W\le G-(G-c)\E_\mu\psi_j'^2-c\,(\E_\mu\psi_j')^2
 \le G\bigl(1-(\E_\mu\psi_j')^2\bigr)=G\bigl(1-(\partial_yW)^2\bigr).
\]
Both terms of $\partial_{yy}W$ in \eqref{eq:pde-tilt} are nonnegative, and the first is positive
if $\psi_j''>0$. Finally $|\partial_yW|<1$, since it is a $\mu$-average of values of $\psi_j'$ in
$(-1,1)$. Hence $W(t,\cdot)\in\mathcal T_G$, and $\partial_{yy}W>0$ if $\psi_j''>0$.

\emph{(d)} follows from $|\partial_yW|<1$.
\end{proof}

Part~(b) of the next lemma is the stochastic-control representation of the Parisi equation
\cite[Theorem~3]{AC15b}, \cite[Lemma~18]{JT16}, written for coefficients that need not be monotone or
bounded by one; its proof uses only $c\ge0$.

\begin{lemma}[Auffinger--Chen process and verification]\label{lem:pde-ac}
Let $c$, $\psi$, $W$ and $v$ be as in Lemma~\ref{lem:pde-class}, let
$(t,y)\in[t_0,t_1]\times\R$, and let $B$ be a standard Brownian motion with filtration
$(\mathcal F_r)$. The equation
\begin{equation}\label{eq:pde-ac}
 dX_r=\xi''(r)c(r)v(r,X_r)\,dr+\sqrt{\xi''(r)}\,dB_r,\qquad X_t=y,
\end{equation}
has a unique strong solution on $[t,t_1]$.
\begin{thmlist}
\item\label{lem:pde-ac-a} For $r\in[t,t_1]$,
 \[
  v(r,X_r)=v(t,y)+\int_t^r\sqrt{\xi''(s)}\,\partial_yv(s,X_s)\,dB_s .
 \]
 In particular $v(r,X_r)=\E[\psi'(X_{t_1})\mid\mathcal F_r]$ and $v(t,y)=\E\psi'(X_{t_1})$.
\item\label{lem:pde-ac-b} \textup{(Verification.)} Let $u$ be a bounded progressively
 measurable process, and let $X^u$ solve $dX^u_r=\xi''(r)c(r)u_r\,dr+\sqrt{\xi''(r)}\,dB_r$,
 $X^u_t=y$. Then
 \begin{equation}\label{eq:pde-verification}
  \E\psi(X^u_{t_1})=W(t,y)+\frac12\E\int_t^{t_1}\xi''c\,u_r^2\,dr
  -\frac12\E\int_t^{t_1}\xi''c\,\bigl(u_r-v(r,X^u_r)\bigr)^2\,dr .
 \end{equation}
 Hence $W(t,y)\ge\E\psi(X^u_{t_1})-\frac12\E\int_t^{t_1}\xi''c\,u_r^2\,dr$, with equality for
 $u_r=v(r,X_r)$.
\end{thmlist}
\end{lemma}

\begin{proof}
By Lemma~\ref{lem:pde-class}\ref{lem:pde-class-a}, $c\le G$ and $\xi''\le\xi''(1)$, the drift in
\eqref{eq:pde-ac} is bounded and $G^2\xi''(1)$-Lipschitz in the space variable, so strong existence
and uniqueness are standard.

(a) On $[\tau_{j-1}+\epsilon',\tau_j-\epsilon]$ the function $v$ is $C^{1,2}$
(Lemma~\ref{lem:pde-class}\ref{lem:pde-class-c}), and It\^o's formula together with the equation
for $v$ gives $dv(r,X_r)=\sqrt{\xi''(r)}\,\partial_yv(r,X_r)\,dB_r$: the drift terms cancel. Let
$\epsilon,\epsilon'\downarrow0$, using the continuity of $v$ and of the paths of $X$ and the bound
$|\partial_yv|\le G$, and chain over the steps. The stochastic integral is a bounded martingale,
and $v(t_1,\cdot)=\psi'$; this gives the conditional expectations.

(b) On each closed step, Lemma~\ref{lem:pde-class}\ref{lem:pde-class-b}, It\^o's formula and the
equation for $W$ give
\[
 dW(r,X^u_r)=\tfrac12\xi''(r)c(r)\bigl[u_r^2-\bigl(u_r-v(r,X^u_r)\bigr)^2\bigr]dr
 +\sqrt{\xi''(r)}\,v(r,X^u_r)\,dB_r .
\]
The stochastic integral is a true martingale because $|v|<1$, and $\psi$ grows at most linearly
while $X^u_{t_1}$ has Gaussian tails. Taking expectations and summing over the steps gives
\eqref{eq:pde-verification}; the inequality follows from $c\ge0$.
\end{proof}

The verification identity turns comparisons of coefficients and of terminal values into
comparisons of solutions: one runs the optimal control of one system in the other.

\begin{lemma}[Coefficient sandwich]\label{lem:pde-coef}
Let $c_1,c_2$ be step coefficients on $[t_0,t_1)$ with values in $[0,G]$, let $\psi\in\mathcal T_G$,
$W_i=W[c_i,\psi]$ and $v_i=\partial_yW_i$. For $y\in\R$ let $X^i$ solve \eqref{eq:pde-ac} for
$c_i$ from $(t_0,y)$, and put $\Gamma_i(r)=\E\,v_i(r,X^i_r)^2$. Then
\begin{equation}\label{eq:pde-coef}
 \frac12\int_{t_0}^{t_1}\xi''(c_2-c_1)\Gamma_1\,dr
 \ \le\ W_2(t_0,y)-W_1(t_0,y)\ \le\
 \frac12\int_{t_0}^{t_1}\xi''(c_2-c_1)\Gamma_2\,dr .
\end{equation}
In particular $|W_2-W_1|\le\frac12\xi''(1)\|c_1-c_2\|_{L^1}$ on $[t_0,t_1]\times\R$.
\end{lemma}

\begin{proof}
Apply Lemma~\ref{lem:pde-ac}\ref{lem:pde-ac-b} to the system with coefficient $c_2$, with the
control $u_r=v_1(r,X^1_r)$ and the Brownian motion that drives $X^1$. Then
$X^u_{t_1}=X^1_{t_1}+Y$ with $Y=\int_{t_0}^{t_1}\xi''(c_2-c_1)u_r\,dr$, and convexity of
$\psi$ gives
\[
 W_2(t_0,y)\ge\E\psi(X^1_{t_1})+\E\bigl[\psi'(X^1_{t_1})Y\bigr]
 -\frac12\E\int_{t_0}^{t_1}\xi''c_2u_r^2\,dr .
\]
The equality case of Lemma~\ref{lem:pde-ac}\ref{lem:pde-ac-b} for the first system gives
$W_1(t_0,y)=\E\psi(X^1_{t_1})-\frac12\E\int_{t_0}^{t_1}\xi''c_1u_r^2\,dr$. By
Lemma~\ref{lem:pde-ac}\ref{lem:pde-ac-a} and Fubini's theorem,
$\E[\psi'(X^1_{t_1})u_r]=\E u_r^2=\Gamma_1(r)$, so
$\E[\psi'(X^1_{t_1})Y]=\int_{t_0}^{t_1}\xi''(c_2-c_1)\Gamma_1\,dr$. Subtracting gives the left
inequality in \eqref{eq:pde-coef}; exchanging the two systems gives the right one. The last claim
follows from $0\le\Gamma_i\le1$ and $\xi''\le\xi''(1)$, applied on $[t,t_1]$ for every $t$.
\end{proof}

\begin{lemma}[Terminal sandwich]\label{lem:pde-terminal}
Let $c$ be a step coefficient on $[t_0,t_1)$ with values in $[0,G]$, let $\psi_1,\psi_2\in\mathcal T_G$,
$W_i=W[c,\psi_i]$, and let $X^i$ solve \eqref{eq:pde-ac} for $(c,\psi_i)$ from $(t_0,y)$. Then
\[
 \E[\psi_2-\psi_1](X^1_{t_1})\ \le\ W_2(t_0,y)-W_1(t_0,y)\ \le\ \E[\psi_2-\psi_1](X^2_{t_1}).
\]
In particular $\sup|W_2-W_1|\le\sup|\psi_2-\psi_1|$.
\end{lemma}

\begin{proof}
Apply Lemma~\ref{lem:pde-ac}\ref{lem:pde-ac-b} to the second system with the control
$u_r=v_1(r,X^1_r)$. The coefficient is the same, so $X^u=X^1$, and the equality case for the first
system gives the lower bound. Exchanging the systems gives the upper bound.
\end{proof}

\begin{lemma}[Stability]\label{lem:pde-stability}
For $i=1,2$ let $c_i$ be a step coefficient on $[t_0,t_1)$ with values in $[0,G]$, let
$\psi_i\in\mathcal T_G$, $W_i=W[c_i,\psi_i]$ and $v_i=\partial_yW_i$, and put $E_0=e^{G^2\xi'(1)}$.
Then
\[
 D:=\sup_{[t_0,t_1]\times\R}|v_1-v_2|
 \le e^{G^2E_0\xi'(1)}\bigl(\|\psi_1'-\psi_2'\|_\infty+GE_0\xi''(1)\|c_1-c_2\|_{L^1}\bigr),
\]
and the solutions $X^1,X^2$ of \eqref{eq:pde-ac} for $(c_i,\psi_i)$, started at a common point and
driven by the same Brownian motion, satisfy almost surely
\[
 \sup_r|X^1_r-X^2_r|\le E_0\bigl(\xi''(1)\|c_1-c_2\|_{L^1}+G\xi'(1)D\bigr).
\]
\end{lemma}

\begin{proof}
Let $D(t)=\sup_y|v_1-v_2|(t,y)$, which is lower semicontinuous, hence measurable, and let
$e(r)=|X^1_r-X^2_r|$ for solutions started at $(t,y)$. Since $|v_1|\le1$, $c_2\le G$ and
$|\partial_yv_2|\le G$,
\[
 |c_1v_1(r,X^1_r)-c_2v_2(r,X^2_r)|\le|c_1-c_2|(r)+G\,D(r)+G^2e(r).
\]
Hence $e(r)\le a+G^2\int_t^r\xi''e$ with $a=\int_t^{t_1}\xi''(|c_1-c_2|+GD)$, and Gronwall's lemma
with $\int_0^1\xi''=\xi'(1)$ gives $e\le E_0a$. By Lemma~\ref{lem:pde-ac}\ref{lem:pde-ac-a},
$v_i(t,y)=\E\psi_i'(X^i_{t_1})$, and $\psi_2'$ is $G$-Lipschitz, so
$|v_1-v_2|(t,y)\le\|\psi_1'-\psi_2'\|_\infty+G\,\E e(t_1)$. Therefore
$D(t)\le\|\psi_1'-\psi_2'\|_\infty+GE_0\int_t^{t_1}\xi''(|c_1-c_2|+GD)$, and the backward Gronwall
lemma gives the bound on $D$, since $\int\xi''|c_1-c_2|\le\xi''(1)\|c_1-c_2\|_{L^1}$. Inserting it
into $e\le E_0a$ gives the bound on the processes.
\end{proof}

\paragraph{Constants.} We record, for later use,
\begin{equation}\label{eq:pde-constants}
\begin{gathered}
 L_v=GE_0\xi''(1)e^{G^2E_0\xi'(1)},\qquad L_X=E_0\bigl(\xi''(1)+G\xi'(1)L_v\bigr),\qquad
 L_M=L_v+GL_X,\\
 L_\Gamma=2L_M,\qquad L_\psi=E_0G\xi'(1)e^{G^2E_0\xi'(1)} .
\end{gathered}
\end{equation}
If $\psi_1=\psi_2$, Lemma~\ref{lem:pde-stability} gives $D\le L_v\|c_1-c_2\|_{L^1}$,
$\sup_r|X^1_r-X^2_r|\le L_X\|c_1-c_2\|_{L^1}$, hence
$\E|v_1(r,X^1_r)-v_2(r,X^2_r)|\le L_M\|c_1-c_2\|_{L^1}$ and, with $\Gamma_i$ as in
Lemma~\ref{lem:pde-coef}, $|\Gamma_1-\Gamma_2|\le L_\Gamma\|c_1-c_2\|_{L^1}$. If $c_1=c_2$, it gives
$\sup_r|X^1_r-X^2_r|\le L_\psi\|\psi_1'-\psi_2'\|_\infty$. These constants depend only on $G$,
$\xi'(1)$ and $\xi''(1)$.

\paragraph{The Parisi functional.}
For a step coefficient $c$ on $[0,1)$ with values in $[0,G]$, let $X^c$ solve \eqref{eq:pde-ac} for
$(c,\log\cosh)$ from $(0,0)$, let $v_c=\partial_yW[c,\log\cosh]$, and define
\begin{equation}\label{eq:pde-P}
\begin{gathered}
 P(c)=W[c,\log\cosh](0,0)-\frac12\int_0^1t\xi''(t)c(t)\,dt,\qquad
 \Gamma_c(t)=\E\,v_c(t,X^c_t)^2,\\
 \Obs_c(t)=\frac12\int_t^1\xi''(s)\bigl(\Gamma_c(s)-s\bigr)\,ds,\\
 DP(c)[h]=\frac12\int_0^1\xi''(t)h(t)\bigl(\Gamma_c(t)-t\bigr)\,dt=-\int_0^1h\,\Obs_c'\,dt,
\end{gathered}
\end{equation}
for $h\in L^1[0,1]$.
By Lemma~\ref{lem:pde-measurable}\ref{lem:pde-measurable-iv} below, $P(\cdf\zeta)=P_\zeta$,
$\Gamma_{\cdf\zeta}=\Gz\zeta$ and $\Obs_{\cdf\zeta}=\Obs_\zeta$ for $\zeta\in\Pspace([0,1])$.

\begin{corollary}[The Parisi functional on step coefficients]\label{cor:pde-step}
Let $c,c_1,c_2$ be step coefficients on $[0,1)$ with values in $[0,G]$.
\begin{thmlist}
\item\label{cor:pde-step-a} \textup{(Sandwich.)}
 $DP(c_1)[c_2-c_1]\le P(c_2)-P(c_1)\le DP(c_2)[c_2-c_1]$.
\item\label{cor:pde-step-b} \textup{(Convexity.)} $P$ is convex.
\item\label{cor:pde-step-c} \textup{(Lipschitz bound.)}
 $|P(c_1)-P(c_2)|\le\frac12\xi''(1)\|c_1-c_2\|_{L^1}$.
\item\label{cor:pde-step-d} \textup{(Second moments and field laws.)}
 $\sup_t|\Gamma_{c_1}-\Gamma_{c_2}|(t)\le L_\Gamma\|c_1-c_2\|_{L^1}$, and, if $X^{c_1}$ and $X^{c_2}$
 are driven by the same Brownian motion,
 $\E|v_{c_1}(t,X^{c_1}_t)-v_{c_2}(t,X^{c_2}_t)|\le L_M\|c_1-c_2\|_{L^1}$ for every $t$.
\item\label{cor:pde-step-e} \textup{(Second order.)} For every step function $h$ and $\varepsilon\in\R$
 such that $c+\varepsilon h$ takes values in $[0,G]$,
 \[
  P(c+\varepsilon h)\le P(c)+\varepsilon\,DP(c)[h]+K\varepsilon^2\|h\|_{L^1}^2,\qquad
  K:=\tfrac12\xi''(1)L_\Gamma .
 \]
\end{thmlist}
The constant $K$ depends only on $G$, $\xi'(1)$ and $\xi''(1)$.
\end{corollary}

\begin{proof}
(a) is Lemma~\ref{lem:pde-coef} with $[t_0,t_1]=[0,1]$, $\psi=\log\cosh$ and $y=0$, minus the
linear terms of $P$. (b) follows from the left half of (a), which gives an affine minorant of $P$
touching it at $c_1$. (c) follows from (a) and $|\Gamma-t|\le1$. (d) is the paragraph on constants.
For (e), the right half of (a) gives
$P(c+\varepsilon h)-P(c)\le\varepsilon\,DP(c)[h]
+\frac\varepsilon2\int_0^1\xi''h\,(\Gamma_{c+\varepsilon h}-\Gamma_c)\,dt$, and by (d) the last term is
at most $\frac12\xi''(1)L_\Gamma\varepsilon^2\|h\|_{L^1}^2$.
\end{proof}

The Parisi functional of a probability measure has a coefficient that is a distribution function,
which in general is not a step function, and the variational arguments of Section~\ref{sec:loc}
and Appendices~\ref{app:var-constrained}--\ref{app:var-gap} move through coefficients that are
neither. The comparison lemmas extend by continuity. Let $\mathcal B_G=\mathcal B_G([t_0,t_1))$ be the set of measurable
$c:[t_0,t_1)\to[0,G]$, modulo equality almost everywhere, with the $L^1$ distance; step
coefficients are dense in $\mathcal B_G$.

\begin{lemma}[Measurable coefficients]\label{lem:pde-measurable}
Let $c\in\mathcal B_G$ and $\psi\in\mathcal T_G$, and let $c_k$ be step coefficients with values in
$[0,G]$ such that $\|c_k-c\|_{L^1}\to0$. Put $W_k=W[c_k,\psi]$ and $v_k=\partial_yW_k$.
\begin{romanlist}
\item\label{lem:pde-measurable-i} \textup{(Extension.)} $W_k$ and $v_k$ converge uniformly on
 $[t_0,t_1]\times\R$ to limits $W[c,\psi]$ and $v_c$ that do not depend on the sequence. The
 function $v_c$ is continuous, $v_c=\partial_yW[c,\psi]$, $|v_c|\le1$, and $v_c(t,\cdot)$ is
 nondecreasing and $G$-Lipschitz; it is $1$-Lipschitz if $c\le1$ and $\psi\in\mathcal T_1$. If
 $\psi$ is even, then $W[c,\psi](t,\cdot)$ is even, $v_c(t,\cdot)$ is odd and $v_c(t,0)=0$.
\item\label{lem:pde-measurable-ii} \textup{(Auffinger--Chen process.)} For
 $(t,y)\in[t_0,t_1]\times\R$ the equation
 $dX_r=\xi''(r)c(r)v_c(r,X_r)\,dr+\sqrt{\xi''(r)}\,dB_r$, $X_t=y$, has a unique strong solution
 $X^c$, and the solutions $X^k$ for $c_k$, driven by the same Brownian motion, satisfy
 \[
  \sup_r|X^k_r-X^c_r|\le E_0\bigl(\xi''(1)\|c_k-c\|_{L^1}+G\xi'(1)\|v_k-v_c\|_\infty\bigr)\to0,
 \]
 uniformly in the sample point.
\item\label{lem:pde-measurable-iii} \textup{(Consequences.)} $v_k(r,X^k_r)\to v_c(r,X^c_r)$
 uniformly in $r$ and in the sample point. Hence $\Gamma_c(r):=\E\,v_c(r,X^c_r)^2$ is the uniform
 limit of $\Gamma_{c_k}$ and is continuous, and $\Law(v_c(r,X^c_r))$ is the $\Wone$-limit of
 $\Law(v_k(r,X^k_r))$. When $[t_0,t_1)=[0,1)$, $\psi=\log\cosh$ and $(t,y)=(0,0)$, define $P(c)$,
 $\Obs_c$ and $DP(c)$ by \eqref{eq:pde-P}; then $P(c)=\lim_kP(c_k)$,
 $\Obs_c\in C^1([0,1])$ with $\Obs_c'(t)=-\frac12\xi''(t)(\Gamma_c(t)-t)$, and
 Corollary~\ref{cor:pde-step} holds for coefficients in $\mathcal B_G$ (with $h$ measurable in
 \ref{cor:pde-step-e}).
\item\label{lem:pde-measurable-iv} \textup{(Identification.)} Let $\zeta\in\Pspace([0,1])$. Then
 $W[\cdf\zeta,\log\cosh]$ on $[0,1]\times\R$ is the solution $\Phi_\zeta$ of the Parisi equation in
 the sense of \cite[p.~7]{AC15a}, $v_{\cdf\zeta}=\Phi_{\zeta,x}$, the process $X^{\cdf\zeta}$ started
 at $(0,0)$ is the Auffinger--Chen process of $\zeta$, and
 \[
  P(\cdf\zeta)=P_\zeta,\qquad\Gamma_{\cdf\zeta}=\Gz\zeta,\qquad\Obs_{\cdf\zeta}=\Obs_\zeta .
 \]
 For $Q\in[0,1)$ and $\eta\in\Pspace([Q,1])$, the same construction on $[Q,1]$ with the
 coefficient $\cdf\eta(t)=\eta([Q,t])$ defines $\Phi_\eta=W[\cdf\eta,\log\cosh]$.
\end{romanlist}
\end{lemma}

\begin{proof}
(i) By Lemma~\ref{lem:pde-coef}, $\sup|W_k-W_l|\le\frac12\xi''(1)\|c_k-c_l\|_{L^1}$, and by
Lemma~\ref{lem:pde-stability} with $\psi_1=\psi_2$, $\sup|v_k-v_l|\le L_v\|c_k-c_l\|_{L^1}$. Both
sequences are therefore uniformly Cauchy, and interleaving two approximating sequences shows that
the limits do not depend on the sequence. The function $v_c$ is continuous as a uniform limit of
continuous functions (Lemma~\ref{lem:pde-class}\ref{lem:pde-class-c}), and letting $k\to\infty$
in $W_k(t,y)-W_k(t,0)=\int_0^yv_k(t,z)\,dz$ gives $v_c=\partial_yW[c,\psi]$. Monotonicity, the
bound $|v_c|\le1$ and the Lipschitz constants pass to the limit from
Lemma~\ref{lem:pde-class}\ref{lem:pde-class-a}; for the $1$-Lipschitz clause choose the $c_k$ with
values in $[0,1]$, which does not change the limit. If $\psi$ is even, every $W_k(t,\cdot)$ is even,
because the Gaussian convolutions in \eqref{eq:pde-colehopf} preserve evenness, and so is the
limit.

(ii) The drift is bounded, measurable in $r$ and $G^2\xi''(1)$-Lipschitz in the space variable by
(i), so the strong solution exists and is unique. With $e_k(r)=|X^k_r-X^c_r|$,
$|c_kv_k(r,X^k_r)-cv_c(r,X^c_r)|\le|c_k-c|(r)+G\|v_k-v_c\|_\infty+G^2e_k(r)$, and Gronwall's lemma
gives the bound.

(iii) $|v_k(r,X^k_r)-v_c(r,X^c_r)|\le\|v_k-v_c\|_\infty+G\,e_k(r)$, which tends to zero uniformly.
The statements on $\Gamma_c$ and on the laws follow, the latter because the $\Wone$ distance of
two laws is at most the expected distance of a coupling. $P(c_k)\to P(c)$ by (i). Since $\Gamma_c$
is continuous, $\Obs_c\in C^1$. The quantities in Corollary~\ref{cor:pde-step}, namely $P(c)$,
$DP(c)[h]$, $\Gamma_c$ and the coupled variables $v_c(t,X^c_t)$, are continuous in $c$ (and $DP(c)[h]$
in $h$) for the $L^1$ distance, by (i) and (iii). Approximating $c$, $c_1$, $c_2$ and
$c+\varepsilon h$ by step coefficients with values in $[0,G]$, the statements of the corollary pass
from step coefficients to $\mathcal B_G$.

(iv) For finitely supported $\zeta$ the two definitions agree, by \eqref{eq:pde-colehopf}. For
general $\zeta$, Auffinger and Chen \cite[p.~7]{AC15a} define $\Phi_\zeta$ as the continuous
extension, in the supremum norm, of the map $\zeta\mapsto\Phi_\zeta$ from finitely supported
measures, with the distance $\int_0^1|\cdf\zeta-\cdf{\zeta'}|\,dt$; the Lipschitz property on which
this definition rests is Lemma~\ref{lem:pde-coef}. Hence $\Phi_\zeta=W[\cdf\zeta,\log\cosh]$ by (i),
and $\Phi_{\zeta,x}=v_{\cdf\zeta}$. By strong uniqueness, $X^{\cdf\zeta}$ is the Auffinger--Chen
process of $\zeta$, and the three identities are the definitions of $P_\zeta$, $\Gz\zeta$ and
$\Obs_\zeta$ in Sections~\ref{sec:intro-setting}--\ref{sec:intro-structure}.
\end{proof}

The comparison of Proposition~\ref{prop:glued-comparison} uses, as terminal value at the radius
$Q$, the Fenchel transform of the TAP correction with a separable perturbation. The next lemma
places it in the terminal class.

\begin{lemma}[Perturbed Fenchel terminal values]\label{lem:pde-fenchel}
Let $Q\in[0,1)$, let $\eta\in\Pspace([Q,1])$ have a step distribution function, and put
$\Phi=\Phi_\eta(Q,\cdot)$. As in Section~\ref{sec:intro-setting}, let
$\Phistar_\eta(Q,m)=\inf_{y\in\R}\{\Phi(y)-my\}$ for $m\in(-1,1)$, extended continuously to
$m=\pm1$.
\begin{thmlist}
\item\label{lem:pde-fenchel-a} $\Phi\in\mathcal T_1$, and $\Phi'$ is a strictly increasing $C^1$
 bijection of $\R$ onto $(-1,1)$ with $0<\Phi''\le1-\Phi'^2$.
\item\label{lem:pde-fenchel-b} Let $b$ be the inverse of $\Phi'$ and $h=-\Phistar_\eta(Q,\cdot)$.
 Then $h$ is $C^2$ on $(-1,1)$, with $h'(m)=b(m)\to\pm\infty$ as $m\to\pm1$ and
 $h''(m)=1/\Phi''(b(m))\ge1/(1-m^2)\ge1$.
\item\label{lem:pde-fenchel-c} Let $g_0\in C^2([-1,1])$ and $s\in\R$ with
 $|s|\,\|g_0''\|_\infty\le\frac12$, and put
 \[
  \Psi_s(z)=\sup_{|m|\le1}\bigl[zm+\Phistar_\eta(Q,m)+s\,g_0(m)\bigr],\qquad z\in\R .
 \]
 The supremum is attained at a unique point $m^*_s(z)\in(-1,1)$, $\Psi_s'=m^*_s$,
 $\Psi_s''>0$, and $\Psi_s\in\mathcal T_2$. Moreover $\Psi_0=\Phi$ and $m^*_0=\Phi'$.
\end{thmlist}
\end{lemma}

\begin{proof}
(a) By Lemma~\ref{lem:pde-class}\ref{lem:pde-class-a} with $G=1$ (the coefficient $\cdf\eta$ takes
values in $[0,1]$ and $\log\cosh''>0$), $\Phi\in\mathcal T_1$ and $\Phi''>0$. By
Lemma~\ref{lem:pde-ac}\ref{lem:pde-ac-a}, $\Phi'(y)=\E\tanh(X_1)$ for the solution $X$ of
\eqref{eq:pde-ac} from $(Q,y)$, and $X_1\ge y-\xi'(1)+\int_Q^1\sqrt{\xi''}\,dB$; so $\Phi'(y)\to1$ as
$y\to\infty$, and $\Phi'(y)\to-1$ as $y\to-\infty$ in the same way.

(b) For $|m|<1$ the infimum defining $\Phistar_\eta(Q,m)$ is attained at $y=b(m)$, by convexity, so
$h(m)=mb(m)-\Phi(b(m))$. Since $b$ is $C^1$ with $b'=1/\Phi''(b)$, differentiation gives
$h'=b$ and $h''=b'=1/\Phi''(b)$, and $\Phi''(b(m))\le1-\Phi'(b(m))^2=1-m^2$.

(c) Put $k=h-sg_0$. By (b), $k''(m)\ge(1-m^2)^{-1}-\frac12\ge\frac12$, and $k'(m)\to\pm\infty$ as
$m\to\pm1$ because $g_0'$ is bounded. So $k'$ is an increasing $C^1$ bijection of $(-1,1)$ onto
$\R$, and the strictly concave function $m\mapsto zm-k(m)$ has a unique maximizer
$m^*_s(z)=(k')^{-1}(z)$ in $(-1,1)$; it is strictly monotone on either side of $m^*_s(z)$, so the
continuous extension of $\Phistar_\eta(Q,\cdot)$ to $m=\pm1$ does not change the supremum. Hence
$\Psi_s(z)=zm^*_s(z)-k(m^*_s(z))$, and since $k'(m^*_s(z))=z$ and $m^*_s$ is $C^1$ with
$(m^*_s)'=1/k''(m^*_s)$, we get $\Psi_s'=m^*_s\in(-1,1)$ and $\Psi_s''=1/k''(m^*_s)>0$. The bound
$k''(m)\ge(1+m^2)/(2(1-m^2))$ gives
$\Psi_s''\le2(1-m_s^{*2})/(1+m_s^{*2})\le2(1-\Psi_s'^2)$, so $\Psi_s\in\mathcal T_2$. For $s=0$,
$m^*_0=(h')^{-1}=\Phi'$, and $b(m^*_0(z))=z$ gives
$\Psi_0(z)=zm^*_0(z)-h(m^*_0(z))=\Phi(z)$.
\end{proof}

\subsection{The constrained problem and the regularity of every level}
\label{app:var-constrained}\label{app:var}

This subsection proves Theorem~\ref{thm:exist}\ref{thm:exist-b}, that every level of
$(\feq,\fone)$ is regular, and Lemma~\ref{lem:loc-obstacle} on the obstacle and its contact set.
The dual function $\Lambda(\theta)=\theta p(\theta)$ is the value of a family of convex problems:
for $\theta\in[0,1]$ one minimizes $\zeta\mapsto P_\zeta$ over
\[
 \mathcal Z_\theta=\bigl\{\zeta\in\Pspace([0,1]):\ \zeta(\{0\})\ge\theta\bigr\}.
\]
The set $\mathcal Z_\theta$ is convex, and it is weakly closed, since
$\limsup_n\zeta_n(\{0\})\le\zeta(\{0\})$ whenever $\zeta_n\to\zeta$ weakly; so it is compact.
Note that $\mathcal Z_0=\Pspace([0,1])$ and $\mathcal Z_1=\{\delta_0\}$. We use the distance
\[
 d(\mu,\mu')=\int_0^1\bigl|\cdf\mu(s)-\cdf{\mu'}(s)\bigr|\,ds,\qquad \mu,\mu'\in\Pspace([0,1]),
\]
which is the Wasserstein distance on $[0,1]$ and induces the weak topology. The logic is short.
Strict convexity of $P$ makes each constrained minimizer unique and forces the constraint to be
active above $\thstar$. A quadratic remainder in the first variation then makes $\Lambda$
continuously differentiable, and uniqueness makes $\Lambar$ strictly convex on $[\thstar,1]$, so
that $\Lambda'$ is a homeomorphism onto $(\feq,\fone)$ and every level has exactly one dual
minimizer. Three properties of the Parisi functional carry these arguments.
\begin{itemize}[leftmargin=3.2em,topsep=3pt,itemsep=2pt]
\item[(V1)] \emph{Strict convexity.} The map $\zeta\mapsto P_\zeta$ is strictly convex on
 $\Pspace([0,1])$.
\item[(V2)] \emph{Lipschitz continuity.} There is $L_\xi<\infty$ such that
 $|P_\mu-P_{\mu'}|\le L_\xi\,d(\mu,\mu')$ for all $\mu,\mu'$.
\item[(V3)] \emph{First variation with a quadratic remainder.} There is $C_\xi<\infty$ such that
 \begin{equation}\label{eq:var-quadratic}
  \Bigl|P_{\tilde\mu}-P_\mu-\int\Obs_\mu\,d(\tilde\mu-\mu)\Bigr|\le C_\xi\,d(\tilde\mu,\mu)^2
  \qquad\text{for all } \mu,\tilde\mu\in\Pspace([0,1]).
 \end{equation}
\end{itemize}
Under \condR\ all three hold. Property (V1) is due to Auffinger and Chen \cite[Theorem~2]{AC15b},
who state it under a summability condition such as~\eqref{hyp:C}. Jagannath and Tobasco prove it
for every $\xi$ whose radius of convergence exceeds one \cite[Theorem~2]{JT16}; their strictness step \cite[Lemma~21]{JT16} uses
that $\xi''>0$ on $(0,1]$, which holds here because $\xi\not\equiv0$ and every degree is at least
two. Property (V2) holds with $L_\xi=\frac32\xi''(1)$: the solutions of the Parisi equation
satisfy $\sup|\Phi_\mu-\Phi_{\mu'}|\le\xi''(1)\,d(\mu,\mu')$ \cite[Lemma~14]{JT16}, and the integral
terms of $P_\mu$ and $P_{\mu'}$ differ by at most $\frac12\xi''(1)\,d(\mu,\mu')$. Property (V3) is the
inequality on which the proof of \cite[Lemma~3.2]{JT17} rests: the first variation $G_\mu$ of
\cite[(1.1.1)]{JT17}, at inverse temperature one and zero external field, is the obstacle
$\Obs_\mu$, and~\eqref{eq:var-quadratic} is a quantitative form of the first-variation formula
\cite[Theorem~2]{Che17}. The papers \cite{JT16,JT17} assume only that the radius of convergence of
$\xi$ exceeds one. Under \condR, the convexity of $P$, (V2) and (V3) also follow from
Appendix~\ref{app:pde-toolbox}, with the smaller constant $\frac12\xi''(1)$ in place of $L_\xi$
and $C_\xi=K$, the constant
of~\eqref{eq:loc-second-order}: apply Corollary~\ref{cor:pde-step}, extended to measurable
coefficients by Lemma~\ref{lem:pde-measurable}, to $c=\cdf\mu$ and $h=\cdf{\tilde\mu}-\cdf\mu$, and
use~\eqref{eq:loc-moving-mass}. Appendix~\ref{app:pde-toolbox} does not give strict convexity.
The solution $u_\mu$ used in \cite{JT16,JT17} agrees with $\Phi_\mu$ under \condR: for $\mu$ with
finitely many atoms the Cole--Hopf solution is a weak solution with $|\partial_x\Phi_\mu|\le1$,
hence equals $u_\mu$ by the uniqueness in \cite[Proposition~8.1]{JT17}, and both are continuous in $\mu$ for the sup
norm (\cite[Lemma~14]{JT16} and Lemma~\ref{lem:pde-coef}).

\paragraph{The dual function.} Except for the value of $P_{\delta_0}$ and the lower bound on
$\fone$, the facts in the next lemma are in \cite[Lemma~3.2]{Bou26b} and its proof. We include
the short argument because the proof of Lemma~\ref{lem:loc-obstacle} below needs it with the
calculus of Appendix~\ref{app:pde-toolbox} alone.

\begin{lemma}[The dual function]\label{lem:var-dual}
Assume \condR. For every $\theta\in[0,1]$, $P$ attains its infimum over $\mathcal Z_\theta$. The
function $\Lambar$ is convex and continuous on $[0,1]$, $\theta\feq\le\Lambar(\theta)\le\theta
P_{\delta_0}$ there, and
\[
 \Lambar(1)=\lim_{\theta\uparrow1}\Lambda(\theta)=P_{\delta_0}=\frac{\xi(1)}2,\qquad
 \Lambar'_+(0)=\feq,\qquad P_{\delta_0}\le\fone\le P_{\delta_0}+L_\xi .
\]
\end{lemma}
The proof uses only the convexity of $P$ and (V2).

\begin{proof}
Minimizers exist because $\mathcal Z_\theta$ is compact and $P$ is continuous by (V2). Since
$\delta_0\in\mathcal Z_\theta$ and $P\ge\feq$, we have $\theta\feq\le\Lambda(\theta)\le\theta
P_{\delta_0}$ for $\theta\in[0,1)$.

\emph{The endpoint $1$.} If $0<h\le1$ and $\mu\in\mathcal Z_{1-h}$, then $\cdf\mu\ge1-h$, so
$d(\mu,\delta_0)=\int_0^1(1-\cdf\mu)\le h$ and $P_\mu\ge P_{\delta_0}-L_\xi h$ by (V2). Hence
\begin{equation}\label{eq:var-endpoint}
 (1-h)\bigl(P_{\delta_0}-L_\xi h\bigr)\le\Lambda(1-h)\le(1-h)P_{\delta_0},\qquad 0<h\le1,
\end{equation}
so $\Lambda(\theta)\to P_{\delta_0}$ as $\theta\uparrow1$, and $\Lambar(1)=P_{\delta_0}$.

\emph{Convexity.} Let $\theta_1,\theta_2\in[0,1]$, $t\in(0,1)$ and $\theta=t\theta_1+(1-t)\theta_2$;
the cases $\theta=0$ and $\theta=1$ are trivial. Let $\zeta^{(i)}$ minimize $P$ over
$\mathcal Z_{\theta_i}$, so that $\theta_iP_{\zeta^{(i)}}=\Lambar(\theta_i)$ (for $\theta_i=1$ this
is $\zeta^{(i)}=\delta_0$ and $\Lambar(1)=P_{\delta_0}$). The probability measure
\begin{equation}\label{eq:var-mixture}
\begin{split}
 &\tilde\mu=\frac{t\theta_1\zeta^{(1)}+(1-t)\theta_2\zeta^{(2)}}{\theta}\quad\text{satisfies}\\
 &\tilde\mu(\{0\})\ge\frac{t\theta_1^2+(1-t)\theta_2^2}{\theta}
 =\theta+\frac{t(1-t)(\theta_1-\theta_2)^2}{\theta}\ge\theta,
\end{split}
\end{equation}
so $\tilde\mu\in\mathcal Z_\theta$, and convexity of $P$ with the weights $t\theta_1/\theta$ and
$(1-t)\theta_2/\theta$ gives
$\Lambda(\theta)\le\theta P_{\tilde\mu}\le t\Lambar(\theta_1)+(1-t)\Lambar(\theta_2)$.

\emph{Continuity and slopes.} The convex function $\Lambar$ is continuous on $(0,1)$, at $1$
by~\eqref{eq:var-endpoint}, and at $0$ because $|\Lambda(\theta)|\le\theta\max|P|$. Let $\zeta'$
minimize $P$ over $\Pspace([0,1])$. For $s\in(0,1)$ the measure $(1-s)\zeta'+s\delta_0$ lies in
$\mathcal Z_s$ at distance at most $s$ from $\zeta'$, so $\feq\le p(s)\le\feq+L_\xi s$ and
$\Lambar'_+(0)=\lim_{s\downarrow0}p(s)=\feq$. The left derivative $\fone$ at $1$ is at least the
chord slope $\Lambar(1)-\Lambar(0)=P_{\delta_0}$, and it is the limit of the chord slopes
$(\Lambar(1)-\Lambda(1-h))/h$, which are at most $P_{\delta_0}+L_\xi$
by~\eqref{eq:var-endpoint}.

\emph{The value $P_{\delta_0}$.} For $\zeta=\delta_0$ the coefficient $\zeta([0,t])$ equals one,
and $\Phi(t,x)=\log\cosh x+\frac12(\xi'(1)-\xi'(t))$ solves the Parisi equation, because
$\Phi_{xx}+\Phi_x^2=1$. With $\int_0^1t\xi''(t)\,dt=\xi'(1)-\xi(1)$ this gives
$P_{\delta_0}=\frac12\xi'(1)-\frac12(\xi'(1)-\xi(1))=\frac12\xi(1)$.
\end{proof}

\paragraph{Constrained minimizers.} Strict convexity turns the dual function into a smooth
object. The derivative of $\Lambda$ comes from an envelope argument (compare \cite[Lemma~5.3]{CGKLM} at
zero temperature): $\Lambda$ is convex, and the
quadratic remainder in (V3) makes $p=\Lambda/\theta$ semiconcave, so the two one-sided
derivatives are squeezed together.

\begin{lemma}[Constrained minimizers]\label{lem:var-basic}
Assume \condR. Recall that $\thstar=\zP(\{0\})$, where $\zP$ is the minimizer of $P$ over
$\Pspace([0,1])$, unique by~\ref{lem:var-basic-a}.
\begin{thmlist}
\item\label{lem:var-basic-a} For every $\theta\in[0,1]$, $P$ has exactly one minimizer
 $\zeta_\theta$ on $\mathcal Z_\theta$. In particular the Parisi measure $\zP=\zeta_0$ is unique, and
 $\zeta_1=\delta_0$.
\item\label{lem:var-basic-b} If $\theta\le\thstar$ and $\theta<1$, then $\zeta_\theta=\zP$ and
 $\Lambda(\theta)=\theta\feq$. If $\thstar<\theta<1$, then $\zeta_\theta(\{0\})=\theta$.
\item\label{lem:var-basic-c} The functions $p$ and $\Lambda$ are continuously differentiable on
 $(0,1)$, with
 \begin{equation}\label{eq:var-pprime}
  p'(\theta)=\frac{\Obs_{\zeta_\theta}(0)-\int\Obs_{\zeta_\theta}\,d\zeta_\theta}{1-\theta},
  \qquad \Lambda'(\theta)=p(\theta)+\theta p'(\theta).
 \end{equation}
 The derivative $\Lambda'$ is nondecreasing, and it equals $\feq$ on $(0,\thstar]\cap(0,1)$.
\item\label{lem:var-basic-d} $\Lambar$ is strictly convex on $[\thstar,1]$.
\end{thmlist}
\end{lemma}

\begin{proof}
\ref{lem:var-basic-a} Existence is Lemma~\ref{lem:var-dual}, and uniqueness follows from (V1)
because $\mathcal Z_\theta$ is convex.

\ref{lem:var-basic-b} If $\theta\le\thstar$, then $\zP\in\mathcal Z_\theta$ minimizes $P$ over all
of $\Pspace([0,1])$, so $\zeta_\theta=\zP$ by uniqueness, and $p(\theta)=\feq$. Let
$\thstar<\theta<1$ and suppose that $a=\zeta_\theta(\{0\})>\theta$. For $\mu\in\Pspace([0,1])$
and $0<s\le1-\theta/a$, the measure $(1-s)\zeta_\theta+s\mu$ has mass at least $(1-s)a\ge\theta$
at $0$, so it lies in $\mathcal Z_\theta$, and minimality and convexity give
$P_{\zeta_\theta}\le(1-s)P_{\zeta_\theta}+sP_\mu$, that is, $P_{\zeta_\theta}\le P_\mu$. So
$\zeta_\theta$ minimizes $P$ over $\Pspace([0,1])$, hence $\zeta_\theta=\zP$; but
$\zP(\{0\})=\thstar<\theta<\zeta_\theta(\{0\})$, a contradiction.

\ref{lem:var-basic-c} Fix $\theta\in[0,1)$ and let $D_\theta$ denote the right side of the first
identity in~\eqref{eq:var-pprime}. For $\theta'\in[0,1)$ put $\lambda=(\theta'-\theta)/(1-\theta)$
and consider
\[
 \zeta_\theta+\lambda(\delta_0-\zeta_\theta)
 =\theta'\delta_0+\frac{1-\theta'}{1-\theta}\bigl(\zeta_\theta-\theta\delta_0\bigr).
\]
Since $\zeta_\theta(\{0\})\ge\theta$, the right side is a probability measure with mass at least
$\theta'$ at the origin, so it lies in $\mathcal Z_{\theta'}$. Its distribution function differs from
$\cdf{\zeta_\theta}$ by $\lambda(1-\cdf{\zeta_\theta})$, and
$\int_0^1(1-\cdf{\zeta_\theta})\le1-\zeta_\theta(\{0\})\le1-\theta$, so its distance from
$\zeta_\theta$ is at most $|\theta'-\theta|$; and $\lambda\int\Obs_{\zeta_\theta}\,d(\delta_0-\zeta_\theta)
=(\theta'-\theta)D_\theta$. Comparing $p(\theta')$ with the value of $P$ at this measure and
applying~\eqref{eq:var-quadratic} at $\mu=\zeta_\theta$ gives
\begin{equation}\label{eq:var-semiconcave}
 p(\theta')\le p(\theta)+(\theta'-\theta)D_\theta+C_\xi(\theta'-\theta)^2,
 \qquad \theta,\theta'\in[0,1).
\end{equation}
Now let $\theta\in(0,1)$. The convex function $\Lambda$ has one-sided derivatives
$\Lambda'_-(\theta)\le\Lambda'_+(\theta)$, so $p=\Lambda/\theta$ has one-sided derivatives
$p'_\pm(\theta)=(\Lambda'_\pm(\theta)-p(\theta))/\theta$ with $p'_-(\theta)\le p'_+(\theta)$.
Taking $\theta'=\theta\pm h$ in~\eqref{eq:var-semiconcave} and letting $h\downarrow0$ gives
$p'_+(\theta)\le D_\theta\le p'_-(\theta)$. Hence $p$ is differentiable at $\theta$ with
$p'(\theta)=D_\theta$, and so is $\Lambda=\theta p$, with $\Lambda'=p+\theta p'$. The derivative
of the convex function $\Lambda$ is nondecreasing, and as a derivative it has the intermediate
value property, so it has no jumps: $\Lambda\in C^1(0,1)$, and then $p=\Lambda/\theta$ is
$C^1$ as well. By~\ref{lem:var-basic-b}, $\Lambda(\theta)=\theta\feq$ on $[0,\thstar]\cap[0,1)$,
so $\Lambda'=\feq$ on $(0,\thstar)$ and, by continuity, at $\thstar$ if $\thstar<1$.

\ref{lem:var-basic-d} Let $\thstar\le\theta_1<\theta_2\le1$, $t\in(0,1)$ and
$\theta=t\theta_1+(1-t)\theta_2\in(\thstar,1)$. Form $\tilde\mu$ from $\zeta_{\theta_1}$ and
$\zeta_{\theta_2}$ as in~\eqref{eq:var-mixture}. Since $\theta_1\ne\theta_2$,
\eqref{eq:var-mixture} and~\ref{lem:var-basic-b} give
$\tilde\mu(\{0\})>\theta=\zeta_\theta(\{0\})$, so $\tilde\mu\in\mathcal Z_\theta$ and
$\tilde\mu\ne\zeta_\theta$. Uniqueness of the minimizer gives $p(\theta)<P_{\tilde\mu}$, and convexity of $P$ gives
$\Lambda(\theta)<\theta P_{\tilde\mu}\le t\Lambar(\theta_1)+(1-t)\Lambar(\theta_2)$.
\end{proof}

The next proposition is the regularity of all levels. Recall that a level $f\in(\feq,\fone)$ is
regular if $\theta\mapsto\Lambar(\theta)-\theta f$ has a unique minimizer $\thf$ on $[0,1]$.

\begin{proposition}[Every level is regular]\label{prop:var-regular}
Assume \condR.
\begin{thmlist}
\item\label{prop:var-regular-a} If $\thstar=1$, then $\Lambar(\theta)=\theta\feq$ on $[0,1]$ and
 $\fone=\feq$. If $\thstar<1$, then $\Lambar'_+(\thstar)=\feq<\fone$, and $\Lambda'$ is a
 strictly increasing homeomorphism of $(\thstar,1)$ onto $(\feq,\fone)$.
\item\label{prop:var-regular-b} Every $f\in(\feq,\fone)$ is regular, and its dual minimizer is
 $\thf=(\Lambda')^{-1}(f)\in(\thstar,1)$. In particular $\Lambda'(\thf)=f$ and
 $\Sigma(f)=\Lambda(\thf)-\thf f$.
\item\label{prop:var-regular-c} The function $\Sigma$ is finite, concave, nonincreasing and
 $1$-Lipschitz on $\R$. It vanishes on $(-\infty,\feq]$, is strictly decreasing on
 $[\feq,\infty)$, and equals $\xi(1)/2-f$ for $f\ge\fone$. It is continuously differentiable on
 $(\feq,\fone)$ with $\Sigma'(f)=-\thf$, and strictly concave on $[\feq,\fone]$.
\item\label{prop:var-regular-d} The right derivative of $\Sigma$ at $\feq$ equals $-\thstar$.
\end{thmlist}
\end{proposition}

\begin{proof}
\ref{prop:var-regular-a} If $\thstar=1$, then $\zP=\delta_0$, so $\feq=P_{\delta_0}=\Lambar(1)$ by
Lemma~\ref{lem:var-dual}, and with Lemma~\ref{lem:var-basic}\ref{lem:var-basic-b},
$\Lambar(\theta)=\theta\feq$ on $[0,1]$; hence $\fone=\feq$. Let $\thstar<1$. If $\thstar=0$,
then $\Lambar'_+(0)=\feq$ by Lemma~\ref{lem:var-dual}; if $0<\thstar<1$, then
$\Lambar'_+(\thstar)=\Lambda'(\thstar)=\feq$ by Lemma~\ref{lem:var-basic}\ref{lem:var-basic-c}. Let
$\theta_m=(1+\thstar)/2$ be the midpoint of $[\thstar,1]$. Convexity and the strict convexity of
Lemma~\ref{lem:var-basic}\ref{lem:var-basic-d} give
\[
 \feq=\Lambar'_+(\thstar)\le\frac{\Lambar(\theta_m)-\Lambar(\thstar)}{\theta_m-\thstar}
 <\frac{\Lambar(1)-\Lambar(\theta_m)}{1-\theta_m}\le\Lambar'_-(1)=\fone .
\]
On $(\thstar,1)$ the derivative $\Lambda'$ is continuous and strictly increasing, by
Lemma~\ref{lem:var-basic}\ref{lem:var-basic-c}--\ref{lem:var-basic-d}. Its limit $L_0$ at
$\thstar+$ is at least $\Lambar'_+(\thstar)$; and for small $h>0$ the mean value theorem, together
with the continuity of $\Lambar$ at $\thstar$, writes the chord slope of $\Lambar$ on
$[\thstar,\thstar+h]$ as a value of $\Lambda'$ in $(\thstar,\thstar+h)$, hence at least $L_0$. Letting
$h\downarrow0$ gives $L_0\le\Lambar'_+(\thstar)$, so $L_0=\feq$. In the same way
$\Lambda'(\theta)\le(\Lambar(1)-\Lambar(\theta))/(1-\theta)\le\fone$ for $\theta<1$, and the mean
value theorem on $[1-h,1]$, with the continuity of $\Lambar$ at $1$, shows that the limit of
$\Lambda'$ at $1-$ is $\fone$. So $\Lambda'$ maps $(\thstar,1)$ onto $(\feq,\fone)$.

\ref{prop:var-regular-b} Let $f\in(\feq,\fone)$, so $\thstar<1$, and put
$G_f(\theta)=\Lambar(\theta)-\theta f$, a convex continuous function on $[0,1]$. By~\ref{prop:var-regular-a},
$\thf=(\Lambda')^{-1}(f)\in(\thstar,1)$ is an interior critical point of $G_f$, hence a minimizer.
For $0\le\theta<\thstar$, Lemma~\ref{lem:var-basic}\ref{lem:var-basic-b} gives
$G_f(\theta)=\theta(\feq-f)>\thstar(\feq-f)=G_f(\thstar)\ge G_f(\thf)$. On $[\thstar,1]$ the
function $G_f$ is strictly convex, so it has at most one minimizer there. Hence $\thf$ is the
unique minimizer.

\ref{prop:var-regular-c} As an infimum over $\theta\in[0,1]$ of affine functions of $f$ with slopes
$-\theta\in[-1,0]$, $\Sigma$ is concave, nonincreasing and $1$-Lipschitz once it is finite; and
$0=\Lambar(0)\ge\Sigma(f)\ge\min_{[0,1]}\Lambar-\max\{0,f\}$. If $f\le\feq$, then
$\Lambar(\theta)-\theta f\ge\theta(\feq-f)\ge0$ by Lemma~\ref{lem:var-dual}, with equality at
$\theta=0$, so $\Sigma(f)=0$. If $f\ge\fone$, then the supporting line
$\Lambar(\theta)\ge\Lambar(1)+\fone(\theta-1)$ gives
$\Lambar(\theta)-\theta f\ge\Lambar(1)-f+(1-\theta)(f-\fone)\ge\Lambar(1)-f$, with equality at
$\theta=1$, so $\Sigma(f)=\Lambar(1)-f=\xi(1)/2-f$. For $f,f+h\in(\feq,\fone)$,
\[
 \Sigma(f+h)\le G_{f+h}(\thf)=\Sigma(f)-h\thf,\qquad
 \Sigma(f)\le G_f(\theta_{f+h})=\Sigma(f+h)+h\theta_{f+h},
\]
and $f\mapsto\thf$ is continuous by~\ref{prop:var-regular-a}; hence $\Sigma'(f)=-\thf$. This
derivative is continuous, negative and strictly decreasing on $(\feq,\fone)$, so $\Sigma$ is
strictly decreasing and strictly concave there. Since $\Sigma$ is continuous, it is strictly
decreasing on $[\feq,\fone]$, and strictly concave there as well: if it were affine on a chord of
$[\feq,\fone]$, $\Sigma'$ would be constant on a nondegenerate subinterval of $(\feq,\fone)$. With
the slope $-1$ on $[\fone,\infty)$, $\Sigma$ is strictly decreasing on $[\feq,\infty)$.

\ref{prop:var-regular-d} If $\thstar=1$, then $\Sigma(f)=\feq-f$ for $f\ge\feq=\fone$ by~\ref{prop:var-regular-c}, and the
right derivative is $-1$. If $\thstar<1$, then $\thf\downarrow\thstar$ as $f\downarrow\feq$
by~\ref{prop:var-regular-a}, and the mean value theorem writes
$(\Sigma(f)-\Sigma(\feq))/(f-\feq)$ as $-\theta_{f'}$ for some $f'\in(\feq,f)$; this tends to
$-\thstar$.
\end{proof}

\begin{proof}[Proof of Theorem~\ref{thm:exist}\ref{thm:exist-b}]
Part~\ref{thm:exist-b} is
Lemma~\ref{lem:var-basic}\ref{lem:var-basic-c} together with
Proposition~\ref{prop:var-regular}\ref{prop:var-regular-a}--\ref{prop:var-regular-b}. Both use
only (V1)--(V3), which hold under \condR; so this part of Theorem~\ref{thm:exist} does not need
$\xi\in\classC$. In particular, under \condH\ the constrained minimizer $\zf=\zeta_{\thf}$ is
unique. Section~\ref{sec:loc} does not use this uniqueness: Lemma~\ref{lem:loc-obstacle} and
Theorem~\ref{thm:struct}\ref{thm:struct-i}--\ref{thm:struct-ii} hold for every constrained
minimizer, with proofs that rest on Appendix~\ref{app:pde-toolbox} alone.
\end{proof}

\begin{remark}[The rate at the equilibrium level]\label{rem:var-slope}
Since $\Sigma=0$ on $(-\infty,\feq]$,
Proposition~\ref{prop:var-regular}\ref{prop:var-regular-d} shows that $\feq$ is a corner of
$\Sigma$ exactly when the Parisi measure has an atom at the origin. The identity
$\Sigma'(\feq+)=-\zP(\{0\})$ is recorded in \cite[Remark~1.4]{Bou26a}, where it is matched with an
extreme-value heuristic for the largest TAP free energies, near $N\feq$.
\end{remark}

\begin{proof}[Proof of Lemma~\ref{lem:loc-obstacle}]
The arguments follow \cite[Lemma~3.2]{Bou26b}
and its proof, except for the uniqueness of $\thf$. We give them in full because under
\condH\ they rest only on the calculus of Appendix~\ref{app:pde-toolbox}, which provides the
convexity of $P$, the inequalities \eqref{eq:loc-first-order} and~\eqref{eq:loc-second-order},
and $\Obs_\zeta\in C^1([0,1])$ with $\Obs_\zeta'=-\frac12\xi''(\Gz\zeta-t)$ and $\Obs_\zeta(1)=0$
(Corollary~\ref{cor:pde-step} and Lemma~\ref{lem:pde-measurable}). We use neither the strict
convexity of $P$ nor the uniqueness of the constrained minimizer. Assume \condH\ and fix $f\in(\feq,\fone)$.

\ref{lem:loc-obstacle-a} The proof of Lemma~\ref{lem:var-dual} uses only the convexity of $P$
and (V2), so it applies: $\Lambar$ is convex and continuous on $[0,1]$, in particular lower
semicontinuous, $\Lambar'_+(0)=\feq$, and $\fone=\Lambar'_-(1)$ is finite. The convex continuous
function $G_f(\theta)=\Lambar(\theta)-\theta f$ has right derivative $\feq-f<0$ at $0$ and left
derivative $\fone-f>0$ at $1$, so its minimizers form a nonempty compact subset of $(0,1)$, on
which $\Sigma(f)=\Lambda(\theta)-\theta f$. Uniqueness is proved after~\ref{lem:loc-obstacle-b}.

\ref{lem:loc-obstacle-b} Let $\theta$ be any minimizer of $G_f$. By compactness and continuity,
$P$ attains its infimum over $\mathcal Z_\theta$; let $\zeta$ be any minimizer. Suppose that
$\theta'=\zeta(\{0\})>\theta$. For $\mu\in\Pspace([0,1])$ and $0<\lambda\le1-\theta/\theta'$,
the measure $(1-\lambda)\zeta+\lambda\mu$ lies in $\mathcal Z_\theta$, so minimality and convexity
give $P_\zeta\le P_\mu$. Thus $P_\zeta=\feq$, and for $s\in[\theta,\theta']\cap[0,1)$ the measure
$\zeta$ is admissible at $s$, so $\Lambda(s)=s\feq$ and $G_f(s)=s(\feq-f)$, which is strictly
decreasing in $s$ because $f>\feq$. This contradicts the minimality of $\theta$. Hence
$\zeta(\{0\})=\theta$, and $\Sigma(f)=\Lambda(\theta)-\theta f=\theta(P_\zeta-f)$.

\emph{Uniqueness of $\thf$.} Suppose that $\theta_1<\theta_2$ both minimize $G_f$. By convexity
$G_f$ is constant on $[\theta_1,\theta_2]$, so $\bar\theta=(\theta_1+\theta_2)/2$ is a minimizer
and $\Lambda(\bar\theta)=\frac12(\Lambda(\theta_1)+\Lambda(\theta_2))$. Let $\zeta_i$ minimize $P$
over $\mathcal Z_{\theta_i}$; by the saturation just proved, $\zeta_i(\{0\})=\theta_i$. The
probability measure $\tilde\mu=(\theta_1\zeta_1+\theta_2\zeta_2)/(2\bar\theta)$ has
$\tilde\mu(\{0\})=(\theta_1^2+\theta_2^2)/(2\bar\theta)>\bar\theta$, and by convexity
$\bar\theta P_{\tilde\mu}\le\frac12(\theta_1P_{\zeta_1}+\theta_2P_{\zeta_2})=\Lambda(\bar\theta)$.
So $\tilde\mu$ minimizes $P$ over $\mathcal Z_{\bar\theta}$ at the minimizer $\bar\theta$ of $G_f$, while
$\tilde\mu(\{0\})>\bar\theta$, which contradicts saturation. Hence the minimizer $\thf$ is unique.

\ref{lem:loc-obstacle-c} Let $\theta=\thf$, let $\zeta=\zf$ be a minimizer of $P$ over
$\mathcal Z_\theta$, and write $\Obs=\Obs_\zeta$ and $o=\inf_{(0,1]}\Obs$. The regularity of $\Obs$
was recalled above.

\emph{First-order condition.} For $\zeta'\in\mathcal Z_\theta$ and $\lambda\in(0,1]$ the measure
$(1-\lambda)\zeta+\lambda\zeta'$ lies in $\mathcal Z_\theta$, and its distribution function is
$\cdf\zeta+\lambda(\cdf{\zeta'}-\cdf\zeta)$, with $\|\cdf{\zeta'}-\cdf\zeta\|_{L^1}\le1$. By
minimality, \eqref{eq:loc-second-order} and~\eqref{eq:loc-moving-mass},
$0\le\lambda\int\Obs\,d(\zeta'-\zeta)+K\lambda^2$, and letting $\lambda\downarrow0$,
\begin{equation}\label{eq:var-first-order}
 \int\Obs\,d\zeta'\ge\int\Obs\,d\zeta\qquad\text{for every }\zeta'\in\mathcal Z_\theta .
\end{equation}
With $\zeta'=\theta\delta_0+(1-\theta)\delta_t$, $t\in(0,1]$, and $\zeta(\{0\})=\theta$, this reads
$(1-\theta)\Obs(t)\ge\int_{(0,1]}\Obs\,d\zeta$. Taking the infimum over $t$,
$o\ge(1-\theta)^{-1}\int_{(0,1]}\Obs\,d\zeta\ge o$, because $\zeta((0,1])=1-\theta$. Hence
$\Obs=o$ $\zeta$-almost everywhere on $(0,1]$, and by continuity on $\supp\zeta\cap(0,1]$; also
$\Obs(0)\ge o$ by continuity.

\emph{The strict inequality at the origin.} Suppose that $\Obs(0)=o$. Then $\Obs\ge o$ on
$[0,1]$ and $\int\Obs\,d\zeta=\theta\Obs(0)+(1-\theta)o=o$, so $\int\Obs\,d(\mu-\zeta)\ge0$ for
every $\mu\in\Pspace([0,1])$, and~\eqref{eq:loc-first-order} with~\eqref{eq:loc-moving-mass}
gives $P_\mu\ge P_\zeta$. Thus $P_\zeta=\feq$ and $\Lambda(\theta)=\theta\feq$. For $s\in(\theta,1)$,
the measure
\[
 \zeta_s=\frac{s-\theta}{1-\theta}\,\delta_0+\frac{1-s}{1-\theta}\,\zeta
\]
has mass exactly $s$ at the origin, and
$\cdf{\zeta_s}=\cdf\zeta+\frac{s-\theta}{1-\theta}(1-\cdf\zeta)$. Since
$\int\Obs\,d(\delta_0-\zeta)=\Obs(0)-o=0$, \eqref{eq:loc-second-order}
and~\eqref{eq:loc-moving-mass} give $P_{\zeta_s}\le\feq+K(s-\theta)^2/(1-\theta)^2$, and therefore
\[
 G_f(s)\le s\feq-sf+sK\frac{(s-\theta)^2}{(1-\theta)^2}
 =G_f(\theta)+(s-\theta)(\feq-f)+O\bigl((s-\theta)^2\bigr)<G_f(\theta)
\]
for $s-\theta>0$ small, because $f>\feq$. This contradicts the minimality of $\theta$. Hence
$\Obs(0)>o$.

\emph{The contact set.} The continuous function $\Obs$ attains its minimum over $[0,1]$, and not at
$0$, because $\Obs(0)>o$; so $o$ is attained in $(0,1]$, that is, $o=\Omin$. Next,
$\Gz\zeta(1)=\E\tanh^2(X_1)<1$, so $\Obs'(1)=\frac12\xi''(1)(1-\Gz\zeta(1))>0$; since $\Obs(1)=0$,
$\Obs<0$ just below $1$, and $\Omin<0=\Obs(1)$. Therefore
$\Cf=\argmin_{(0,1]}\Obs=\{t\in[0,1]:\Obs(t)=\Omin\}$ (as $\Obs(0)>\Omin$) is a nonempty compact
subset of $(0,1)$, and it contains $\supp\zeta\cap(0,1]$.
\end{proof}

\begin{remark}[The level and the height of the obstacle at the origin]\label{rem:var-level}
Assume \condH\ and let $f\in(\feq,\fone)$. By Lemma~\ref{lem:loc-obstacle}, $\zf(\{0\})=\thf$ and
$\Obs=\Omin$ on $\supp\zf\cap(0,1]$, so $\int\Obs\,d\zf=\thf\Obs(0)+(1-\thf)\Omin$ and the right
side of~\eqref{eq:var-pprime} at $\theta=\thf$ equals $\Obs(0)-\Omin$. Since $\Lambda'(\thf)=f$,
Lemma~\ref{lem:var-basic}\ref{lem:var-basic-c} gives
\[
 f=P_{\zf}+\thf\bigl(\Obs(0)-\Omin\bigr),\qquad
 \Sigma(f)=\thf\bigl(P_{\zf}-f\bigr)=-\thf^2\bigl(\Obs(0)-\Omin\bigr).
\]
So the strict inequality $\Obs(0)>\Omin$ of Lemma~\ref{lem:loc-obstacle}\ref{lem:loc-obstacle-c}
is equivalent to $\Sigma(f)<0$, and the height of the obstacle at the origin above its minimum,
together with $\thf$, determines the rate.
\end{remark}

\subsection{The relaxed variational gap: the remaining direction}\label{app:var-gap}

We complete the proof of Proposition~\ref{prop:loc-gap}. Section~\ref{sec:loc-gap} proved
parts~\ref{prop:loc-gap-a} and~\ref{prop:loc-gap-b}, and part~\ref{prop:loc-gap-c} when $\Cf$
meets $[0,Q)$: then mass is moved across the radius down to a contact point $r<Q$, by raising the
cascade coefficient on $[r,Q)$. When $\Cf$ lies above $Q$ the mass must move the other way, from
the radius up to a contact point $r>Q$. The cascade part is then left alone, and the band
distribution function is lowered on $[Q,r)$. This is possible only because the constrained
minimizer charges $[0,Q]$ with mass at least $\thf$. As before, $\zeta=\zf$ is a constrained
minimizer, $\theta=\thf$ and $\Obs=\Obs_\zeta$.

\begin{proof}[Proof of Proposition~\ref{prop:loc-gap}\ref{prop:loc-gap-c} when $\Cf\cap{[0,Q)}=\emptyset$]
Let $Q\in[0,1)\setminus\Cf$ and assume that $\Cf$ has no point in $[0,Q)$. By
Lemma~\ref{lem:loc-obstacle}\ref{lem:loc-obstacle-c}, $\Cf$ is a nonempty subset of $(0,1)$, so it
contains a point $r\in(Q,1)$. Put $d=\Obs(Q)-\Omin$, which is positive because $Q\notin\Cf$ when
$Q>0$, and because $\Obs(0)>\Omin$ when $Q=0$.

\emph{The data.} For $0<\varepsilon\le\theta$ let $\xc=\cdf\zeta|_{[0,Q)}$ (there is no cascade
coefficient if $Q=0$) and
\[
 \eta_\varepsilon=\Shift_Q\zeta-\varepsilon\delta_Q+\varepsilon\delta_r .
\]
The measure $\Shift_Q\zeta$ has an atom of mass $\zeta([0,Q])\ge\zeta(\{0\})\ge\theta\ge\varepsilon$
at $Q$, so $\eta_\varepsilon$ is a probability measure on $[Q,1]$; only the constraint
$\zeta(\{0\})\ge\theta$ is used here, not the saturation of Lemma~\ref{lem:loc-obstacle}. The cascade
coefficient $\xc$ is nondecreasing with values in $[\theta,1]$, so $(\xc,\eta_\varepsilon)\in
\mathcal A_Q$ (for $Q=0$, $\eta_\varepsilon=\zeta-\varepsilon\delta_0+\varepsilon\delta_r\in
\Pspace([0,1])=\mathcal A_0$). The distribution function of $\eta_\varepsilon$ on $[Q,1]$ is
$\cdf\zeta-\varepsilon\ind_{[Q,r)}$, so
\[
 \glue(\xc,\eta_\varepsilon)=\cdf\zeta+\varepsilon h,\qquad h=-\ind_{[Q,r)},\qquad
 \|h\|_{L^1}=r-Q\le1,
\]
and $\glue(\xc,\eta_\varepsilon)\ge\theta-\varepsilon\ge0$ on $[Q,r)$, so it is a coefficient. If
$Q>0$, the measure $\zeta$ has no atom at $Q\notin\Cf\supset\supp\zeta\cap(0,1]$, so the glued
coefficient jumps down by $\varepsilon$ at $Q$, as Lemma~\ref{lem:loc-monotone} requires of any gain.

\emph{The gain.} The first variation is
\[
 DP(\cdf\zeta)[h]=-\int_0^1h\,\Obs'\,dt=\int_Q^r\Obs'(t)\,dt=\Obs(r)-\Obs(Q)=-d,
\]
and~\eqref{eq:loc-second-order} gives $P(\glue(\xc,\eta_\varepsilon))\le
P_\zeta-\varepsilon d+K\varepsilon^2$. If $d\le2K\theta$, the choice $\varepsilon=d/(2K)\le\theta$
gives the gain $d^2/(4K)$. Otherwise take $\varepsilon=\theta$; then $K\theta<d/2$ and the gain
is at least $\theta d-K\theta^2\ge\theta d/2$. In both cases
$P(\glue(\xc,\eta_\varepsilon))\le P_\zeta-\varphi(d)$, with $\varphi$ as
in~\eqref{eq:loc-gap-function}.
\end{proof}

In this direction the cap $\varepsilon\le\thf$ cannot be relaxed: since
$\supp\zf\cap(0,1]\subset\Cf$ and $\Cf\cap[0,Q]=\emptyset$, only the atom of $\Shift_Q\zf$ at $Q$,
of mass exactly $\zf([0,Q])=\thf$, can move up to $r$. The cap is the source of the branch
$\thf d/2$ of $\varphi$.

\paragraph{The uniform constant.} The gap in Proposition~\ref{prop:loc-gap}\ref{prop:loc-gap-c}
depends on the radius $Q$ only through $\Obs(Q)-\Omin$, because the constant $K$
of~\eqref{eq:loc-second-order} depends neither on the coefficient nor on the direction. This is
Corollary~\ref{cor:pde-step}\ref{cor:pde-step-e}, extended to measurable coefficients by
Lemma~\ref{lem:pde-measurable}: the upper half of the sandwich in
Corollary~\ref{cor:pde-step}\ref{cor:pde-step-a} bounds $P(c+\varepsilon h)-P(c)-\varepsilon DP(c)[h]$ by
$\frac\varepsilon2\int_0^1\xi''h\,(\Gamma_{c+\varepsilon h}-\Gamma_c)\,dt$, and the Lipschitz bound
$\sup|\Gamma_{c+\varepsilon h}-\Gamma_c|\le L_\Gamma\varepsilon\|h\|_{L^1}$ of
Corollary~\ref{cor:pde-step}\ref{cor:pde-step-d} turns this into $K\varepsilon^2\|h\|_{L^1}^2$ with
$K=\frac12\xi''(1)L_\Gamma$. The constant $L_\Gamma$ of~\eqref{eq:pde-constants} depends only on
$\xi'(1)$, $\xi''(1)$ and the bound $3$ on the coefficients. This uniformity is what lets the
finitely many shells in the proof of Theorem~\ref{thm:struct}\ref{thm:struct-i} share a single
gain.

\subsection{Regularity at the contact set}\label{app:pde-standing}

Section~\ref{sec:hess} and Appendices~\ref{app:matrix}--\ref{app:nc} need, near contacts, bounds
on $\Phi_{\alpha,xx}$ and its next two derivatives uniform in $\alpha$, the identity
$\Obs''=\frac12\xi''\Plef$, and continuity of the field laws in the radius.
Lemma~\ref{lem:pde-standing} proves them for every $\xi\in\classC$, single-degree covariances
included: neither convexity of $\xi$ nor a second active degree is used, and its contact clause,
in~\ref{lem:pde-standing-c}, also uses $\xi''>0$ on $(0,1)$. The class $\classC$ enters through
\cite[Proposition~1(i)]{AC15a}, which is stated for it (compare Remark~\ref{rem:rpde}).

For $\alpha\in\Pspace([0,1])$ let $\Phi=\Phi_\alpha$ be as in
Lemma~\ref{lem:pde-measurable}\ref{lem:pde-measurable-iv}, with derivatives $\Phi_x$,
$\Phi_{xx}$, $\Phi_{xxx}$, $\Phi_{xxxx}$ in the space variable. Let $X$ be the Auffinger--Chen
process of $\alpha$ started at $(0,0)$, and put $M_t=\Phi_x(t,X_t)$ and $U_t=\Phi_{xx}(t,X_t)$; for
$\alpha=\zf$ these give the field law $\Lq{t}=\Law(M_t)$, the variable $\UQ$ at $t=Q$ and the
Plefka factor $\Plef(Q)=1-\xi''(Q)\E U_Q^2$ of Section~\ref{sec:intro-structure}. Put
$V_*=\xi'(1)$ and
\[
 \underline c=e^{-2V_*-2\sqrt{V_*}},\qquad\overline C=8e^{5V_*},\qquad
 C_{34}=2e^{2V_*}+6e^{3V_*}+12V_*e^{8V_*}.
\]

\begin{lemma}[Regularity at the contact set]\label{lem:pde-standing}
Assume $\xi\in\classC$ and let $\alpha\in\Pspace([0,1])$.
\begin{thmlist}
\item\label{lem:pde-standing-a} \textup{(Derivative bounds.)} For every $j\ge0$, $\partial_x^j\Phi$
 exists and is continuous on $[0,1]\times\R$, and for all $(t,b)\in[0,1]\times\R$,
 \[
  \underline c\,e^{-2|b|}\le\Phi_{xx}(t,b)\le\overline C\,e^{-2|b|},\qquad
  |\Phi_{xxx}(t,b)|+|\Phi_{xxxx}(t,b)|\le C_{34}\,\Phi_{xx}(t,b).
 \]
\item\label{lem:pde-standing-b} \textup{(Inverse field.)} For every $t\in[0,1]$, $\Phi_x(t,\cdot)$ is
 an odd, strictly increasing bijection of $\R$ onto $(-1,1)$. Its inverse $b_\alpha(t,m)$ is jointly
 continuous on $[0,1]\times(-1,1)$, and
 \[
  |b_\alpha(t,m)|\le\tfrac12\log(\overline C/2)+\tfrac12\log\frac1{1-|m|}.
 \]
 The function $u_\alpha(t,m)=\Phi_{xx}(t,b_\alpha(t,m))$, extended by $u_\alpha(t,\pm1)=0$, is
 jointly continuous on $[0,1]\times[-1,1]$, and $0\le u_\alpha(t,m)\le(2\overline C/\underline c)(1-|m|)$.
\item\label{lem:pde-standing-c} \textup{(Obstacle.)} $M_t=\int_0^t\sqrt{\xi''(s)}\,U_s\,dB_s$, the
 function $\Gz\alpha$ is in $C^1([0,1])$ with $\Gz\alpha'=\xi''\,\E U_t^2$, and
 $\Obs_\alpha\in C^2([0,1])$ with
 \[
  \Obs_\alpha'(t)=-\tfrac12\xi''(t)\bigl(\Gz\alpha(t)-t\bigr),\qquad
  \Obs_\alpha''(t)=-\tfrac12\xi'''(t)\bigl(\Gz\alpha(t)-t\bigr)-\tfrac12\xi''(t)\bigl(\Gz\alpha'(t)-1\bigr).
 \]
 Every $Q\in\argmin_{(0,1]}\Obs_\alpha$ lies in $(0,1)$ and satisfies $\Gz\alpha(Q)=Q$ and
 $\Obs_\alpha''(Q)=\frac12\xi''(Q)\bigl(1-\xi''(Q)\E U_Q^2\bigr)\ge0$. In particular, for
 $f\in(\feq,\fone)$ and $\alpha=\zf$, one has $\Plef=2\Obs_{\zf}''/\xi''\ge0$ on $\Cf$.
\item\label{lem:pde-standing-d} \textup{(Field laws.)} The maps $t\mapsto\Law(M_t)$ and
 $t\mapsto\Law(U_t)$ are $\Wone$-continuous on $[0,1]$, and $u_\alpha(t,\cdot)$ pushes $\Law(M_t)$
 forward to $\Law(U_t)$ for every $t\in[0,1]$.
\end{thmlist}
The constants $\underline c$, $\overline C$ and $C_{34}$ depend only on $\xi'(1)$; the bounds in
(a)--(b) are uniform in $\alpha$ and $t$.
\end{lemma}

\begin{proof}
Write $a_t=\xi''(t)$.

\emph{(a), finitely supported $\alpha$.} Then $\cdf\alpha$ is a step coefficient with values in
$[0,1]$, and $\Phi$ is given by \eqref{eq:pde-colehopf}. On each step of $\cdf\alpha$,
differentiating \eqref{eq:pde-colehopf} under the Gaussian expectation shows that $\Phi$ is
$C^\infty$ in $x$ with $\partial_x^j\Phi$ bounded for $j\ge1$ (all derivatives of $\log\cosh$ of
order at least one are bounded, and the derivatives produced by \eqref{eq:pde-colehopf} are
polynomials in tilted moments of the derivatives of the terminal value), and that
$\partial_t\partial_x^j\Phi$ is obtained by differentiating the Parisi equation; by dominated
convergence, as in the proof of Lemma~\ref{lem:pde-class}\ref{lem:pde-class-b}, these functions
extend continuously to the closed step. Put $L=\frac12a\,\partial_{xx}+a\cdf\alpha\,\Phi_x\partial_x$.
Differentiating $\partial_t\Phi+\frac12a(\Phi_{xx}+\cdf\alpha\,\Phi_x^2)=0$ in $x$ once, twice, three
and four times gives, on each step,
\begin{gather*}
 (\partial_t+L)\Phi_x=0,\qquad(\partial_t+L)\Phi_{xx}+a\cdf\alpha\,\Phi_{xx}^2=0,\qquad
 (\partial_t+L)\Phi_{xxx}+3a\cdf\alpha\,\Phi_{xx}\Phi_{xxx}=0,\\
 (\partial_t+L)\Phi_{xxxx}+4a\cdf\alpha\,\Phi_{xx}\Phi_{xxxx}+3a\cdf\alpha\,\Phi_{xxx}^2=0 .
\end{gather*}
Let $\tilde X$ be the Auffinger--Chen process of $\alpha$ started at $(t,b)$, whose generator is $L$, and
let $I_{ts}=\int_t^sa_r\cdf\alpha(r)\Phi_{xx}(r,\tilde X_r)\,dr$. By
Lemma~\ref{lem:pde-class}\ref{lem:pde-class-a} with $G=1$, $0<\Phi_{xx}\le1$, so
$0\le I_{ts}\le\int_0^1a=V_*$. It\^o's formula on each step, chained over the steps as in the proof of
Lemma~\ref{lem:pde-ac}\ref{lem:pde-ac-a}, gives the Feynman--Kac formulas
\begin{gather*}
 \Phi_x(t,b)=\E\tanh(\tilde X_1),\qquad
 \Phi_{xx}(t,b)=\E\bigl[\Phi_{xx}(s,\tilde X_s)e^{I_{ts}}\bigr]\quad(t\le s\le1),\\
 \Phi_{xxx}(t,b)=-2\,\E\bigl[\sech^2(\tilde X_1)\tanh(\tilde X_1)e^{3I_{t1}}\bigr],\\
 \Phi_{xxxx}(t,b)=\E\bigl[z(\tilde X_1)e^{4I_{t1}}\bigr]
 +3\,\E\int_t^1a_s\cdf\alpha(s)\Phi_{xxx}(s,\tilde X_s)^2e^{4I_{ts}}\,ds,
\end{gather*}
where $z(y)=4\sech^2y\tanh^2y-2\sech^4y$ is the third derivative of $\tanh$.
Write $\tilde X_1=b+D+J$ with $D=\int_t^1a_s\cdf\alpha(s)\Phi_x(s,\tilde X_s)\,ds$ and $J=\int_t^1\sqrt{a_s}\,dB_s$;
then $|D|\le V_*$, and $J$ is centred Gaussian with variance $\xi'(1)-\xi'(t)\le V_*$. We use
$e^{-2|y|}\le\sech^2y\le4e^{-2|y|}$ and $|b|-|D|-|J|\le|\tilde X_1|\le|b|+|D|+|J|$. Since $I\ge0$, Jensen's
inequality and $\E|J|\le\sqrt{V_*}$ give
\[
 \Phi_{xx}(t,b)\ge\E\sech^2(\tilde X_1)\ge e^{-2|b|-2V_*}\,\E e^{-2|J|}\ge\underline c\,e^{-2|b|},
\]
and since $\E e^{2|J|}\le2\E e^{2J}\le2e^{2V_*}$,
\[
 \Phi_{xx}(t,b)\le e^{V_*}\E\sech^2(\tilde X_1)\le4e^{3V_*}e^{-2|b|}\,\E e^{2|J|}\le\overline C\,e^{-2|b|}.
\]
Next, $e^{3I_{t1}}\le e^{2V_*}e^{I_{t1}}$ gives
$|\Phi_{xxx}(t,b)|\le2e^{2V_*}\E[\sech^2(\tilde X_1)e^{I_{t1}}]=2e^{2V_*}\Phi_{xx}(t,b)$. For the fourth
derivative, $|z|\le6\sech^2$ bounds the first term by $6e^{3V_*}\Phi_{xx}(t,b)$. In the second,
$\Phi_{xxx}^2\le4e^{4V_*}\Phi_{xx}^2\le4e^{4V_*}\Phi_{xx}$ and
$\E\Phi_{xx}(s,\tilde X_s)\le\E[\Phi_{xx}(s,\tilde X_s)e^{I_{ts}}]=\Phi_{xx}(t,b)$, so it is at most
$12e^{8V_*}\int_t^1a_s\,ds\,\Phi_{xx}(t,b)\le12V_*e^{8V_*}\Phi_{xx}(t,b)$. Together these give (a)
for finitely supported $\alpha$.

\emph{(a), general $\alpha$.} Take finitely supported $\alpha_n$ converging weakly to $\alpha$, so
that $\|\cdf{\alpha_n}-\cdf\alpha\|_{L^1}\to0$. By
\cite[Proposition~1(i)]{AC15a}, $\partial_x^j\Phi_\alpha$ exists and is continuous on
$[0,1]\times\R$ for every $j$, and $\partial_x^j\Phi_{\alpha_n}\to\partial_x^j\Phi_\alpha$ uniformly;
here $\Phi_\alpha$ of \cite{AC15a} is ours, by
Lemma~\ref{lem:pde-measurable}\ref{lem:pde-measurable-iv}. The inequalities proved above have
constants independent of $n$, so they pass to the limit.

\emph{(b).} For finitely supported $\alpha$ and $b\ge0$, $\Phi_x(t,b)=\E\tanh(\tilde X_1)\ge
\E\tanh(b-V_*+J)\to1$ as $b\to\infty$, so $1-\Phi_x(t,b)=\int_b^\infty\Phi_{xx}(t,y)\,dy$, and (a)
gives
\[
 \tfrac12\underline c\,e^{-2b}\le1-\Phi_x(t,b)\le\tfrac12\overline C\,e^{-2b}\qquad(b\ge0);
\]
by uniform convergence the same holds for every $\alpha$. The function $\Phi_x(t,\cdot)$ is odd
(Lemma~\ref{lem:pde-measurable}\ref{lem:pde-measurable-i}) and strictly increasing
($\Phi_{xx}>0$), and by the display its limits at $\pm\infty$ are $\pm1$; so it is a bijection onto
$(-1,1)$. For $m\in[0,1)$, $b=b_\alpha(t,m)\ge0$ satisfies $1-m=1-\Phi_x(t,b)$, and the display
gives $\tfrac12\log\frac{\underline c}{2(1-m)}\le b\le\tfrac12\log\frac{\overline C}2+\tfrac12\log\frac1{1-m}$,
with $\overline C/2>1$; oddness gives the bound for $m<0$. If $(t_n,m_n)\to(t,m)$ in
$[0,1]\times(-1,1)$, the numbers $b_\alpha(t_n,m_n)$ are bounded, and every limit point $b'$
satisfies $\Phi_x(t,b')=m$ by joint continuity of $\Phi_x$, so $b'=b_\alpha(t,m)$. Hence $b_\alpha$,
and with it $u_\alpha$, is continuous on $[0,1]\times(-1,1)$. Finally, with $b=b_\alpha(t,m)$,
$e^{-2|b|}\le2(1-|m|)/\underline c$ by the display and the evenness of $\Phi_{xx}(t,\cdot)$, so
$0<u_\alpha(t,m)\le\overline C e^{-2|b|}\le(2\overline C/\underline c)(1-|m|)$, and the extension by
zero at $m=\pm1$ is continuous on $[0,1]\times[-1,1]$.

\emph{(c).} Take $\alpha_n$ as in the proof of (a), with Auffinger--Chen processes $X^n$
started at $(0,0)$ and driven by the same Brownian motion as $X$. By
Lemma~\ref{lem:pde-ac}\ref{lem:pde-ac-a} and $\Phi_{\alpha_n,x}(0,0)=0$ (oddness),
$\Phi_{\alpha_n,x}(t,X^n_t)=\int_0^t\sqrt{a_s}\,\Phi_{\alpha_n,xx}(s,X^n_s)\,dB_s$. As $n\to\infty$,
the left side tends to $M_t$ uniformly, by
Lemma~\ref{lem:pde-measurable}\ref{lem:pde-measurable-iii}. The integrands tend to $U_s$
pointwise, by the uniform convergence in (a), the continuity of $\Phi_{xx}$ and
Lemma~\ref{lem:pde-measurable}\ref{lem:pde-measurable-ii}, and they are bounded by one; by the It\^o
isometry and bounded convergence the stochastic integrals converge in $L^2$. This gives the
representation of $M_t$, which is also \cite[Lemma~2]{AC15b}. By the It\^o isometry,
$\Gz\alpha(t)=\E M_t^2=\int_0^ta_s\E U_s^2\,ds$. The map $s\mapsto\E U_s^2$ is continuous, by
dominated convergence, since $\Phi_{xx}$ is continuous and bounded and $X$ has continuous paths. So
$\Gz\alpha\in C^1([0,1])$ with $\Gz\alpha'=a\,\E U^2$, and since $\xi''\in C^1([0,1])$, the formula
$\Obs_\alpha'=-\frac12\xi''(\Gz\alpha-t)$ of Lemma~\ref{lem:pde-measurable}\ref{lem:pde-measurable-iii}
shows that $\Obs_\alpha\in C^2([0,1])$ with the stated $\Obs_\alpha''$.

Let $Q\in\argmin_{(0,1]}\Obs_\alpha$. Since $\Gz\alpha(1)=\E\tanh^2(X_1)<1$,
$\Obs_\alpha'(1)=\frac12\xi''(1)(1-\Gz\alpha(1))>0$, so $\Obs_\alpha<\Obs_\alpha(1)$ just below $1$,
and $Q<1$. Thus $Q$ is an interior minimum of the $C^2$ function $\Obs_\alpha$, and
$\Obs_\alpha'(Q)=0\le\Obs_\alpha''(Q)$. Since $\xi''(Q)>0$, this gives $\Gz\alpha(Q)=Q$ and
$\Obs_\alpha''(Q)=-\frac12\xi''(Q)(\Gz\alpha'(Q)-1)=\frac12\xi''(Q)(1-\xi''(Q)\E U_Q^2)$. For
$\alpha=\zf$ the set $\argmin_{(0,1]}\Obs_{\zf}$ is $\Cf$, and $1-\xi''(Q)\E U_Q^2=\Plef(Q)$.

\emph{(d).} The functions $\Phi_x$ and $\Phi_{xx}$ are continuous and bounded by one, and $X$ has
continuous paths. If $t_n\to t$, then $M_{t_n}\to M_t$ and $U_{t_n}\to U_t$ almost surely, and by
bounded convergence $\Wone(\Law(M_{t_n}),\Law(M_t))\le\E|M_{t_n}-M_t|\to0$, and likewise for $U$.
Since $\Phi_x(t,\cdot)$ is injective by (b), $b_\alpha(t,M_t)=X_t$, so $u_\alpha(t,M_t)=U_t$, which
is the last claim.
\end{proof}

%% file: sections/appB-glued.tex
\section{Ruelle cascades and the glued comparison}\label{app:glued}

Proposition~\ref{prop:glued-comparison} is Boursier's comparison for the maximum of $\FTAP$ over
the cube \cite[Proposition~2.7]{Bou26b}, run on a sphere or thin shell, and two features of her
argument make this possible. The band measure above the radius $Q$ enters only through the
deterministic correction $\Phistar_\eta(Q,\cdot)$, whose Fenchel transform ($\Phi_\eta(Q,\cdot)$
when $s=0$) is the terminal value at $Q$ of a Cole--Hopf recursion for the cascade part, and no step
ties the cascade parameters to $\eta$ or bounds them by
one. On a sphere the self-overlap does not vary, so the cascade increments below $Q$ have
exponential moments given exactly by that recursion, and Doob's inequality is needed only across
the shell, not on every level. Appendices~\ref{app:pde-cascade}--\ref{app:glued-approx} prove the
cascade facts, the proposition and its extension to general data, \ref{app:glued-map} maps each
step onto \cite[Section~2]{Bou26b}, and \ref{app:var-field} proves the field law,
Theorem~\ref{thm:struct}\ref{thm:struct-ii}. Throughout, $\xi$ satisfies \condR;
Appendix~\ref{app:pde-cascade} uses $\xi$ only on $[0,1]$, through $\xi(0)=0$ and the facts in
the standing assumption of Appendix~\ref{app:pde}.

\subsection{Cascades and the Onsager shift}\label{app:pde-cascade}

\paragraph{Cascades and fields.}
Fix radial levels $0=r_0<r_1<\dots<r_L\le1$, $L\ge1$, and parameters $\xc_1,\dots,\xc_L>0$, and
write $\xc(t)=\xc_\ell$ for $t\in[r_{\ell-1},r_\ell)$. No ordering of the parameters is assumed
until Lemma~\ref{lem:cascade-domination}. Let $\N^0=\{\varnothing\}$ and
$\mathcal V=\bigcup_{\ell=1}^{L}\N^\ell$. A vertex $v\in\N^\ell$ has depth $|v|=\ell$, its
children are $vn$, $n\in\N$, and $v|k$ is its ancestor at depth $k\le\ell$. The leaves are the
$\alpha\in\N^L$, and $\alpha\wedge\beta=\max\{k\le L:\alpha|k=\beta|k\}$. To every vertex $v$ of
depth $\ell-1<L$ attach a Poisson point process on $\R$ with intensity
$\xc_\ell e^{-\xc_\ell u}\,du$, all of them independent; list its points in decreasing order and
give the $n$-th one to the edge from $v$ to $vn$ as its label $u_{vn}$. Set
\[
 u_\alpha=\sum_{\ell=1}^{L}u_{\alpha|\ell},\qquad
 U=\sup_{\alpha\in\N^L}u_\alpha\in(-\infty,\infty].
\]
Let $\vth(t)=t\xi'(t)-\xi(t)$, so that $\vth(0)=0$ and $\vth'(t)=t\xi''(t)\ge0$ on $[0,1]$.
Independently of the labels, let $W^i_v$, $v\in\mathcal V$, $0\le i\le N$, be independent standard
Brownian motions, and for $s\in[0,r_L]$ and leaves $\alpha$ put
\begin{align*}
 B^i(s,\alpha)&=\sum_{\ell=1}^{L}\bigl[W^i_{\alpha|\ell}(\xi'(s\wedge r_\ell))
 -W^i_{\alpha|\ell}(\xi'(s\wedge r_{\ell-1}))\bigr],\qquad1\le i\le N,\\
 Y(s,\alpha)&=\sqrt N\sum_{\ell=1}^{L}\bigl[W^0_{\alpha|\ell}(\vth(s\wedge r_\ell))
 -W^0_{\alpha|\ell}(\vth(s\wedge r_{\ell-1}))\bigr].
\end{align*}
These are centred Gaussian processes with continuous paths, independent of the labels and of each
other, with
\[
 \E\,B^i(s,\alpha)B^j(t,\beta)=\delta_{ij}\,\xi'\bigl(s\wedge t\wedge r_{\alpha\wedge\beta}\bigr),
 \qquad
 \E\,Y(s,\alpha)Y(t,\beta)=N\vth\bigl(s\wedge t\wedge r_{\alpha\wedge\beta}\bigr);
\]
they are the fields of \cite[Definition~2.3]{Bou26b} on our levels, and we call $Y$ the
\emph{Onsager field}. The restriction of the fields of a leaf $\alpha$ to $[0,r_\ell]$ depends only
on $\alpha|\ell$; we call it the path of $\alpha|\ell$. Let $\mathcal G_\ell$ be the
$\sigma$-field generated by the labels of depth at most $\ell$ and the paths of the vertices of
depth $\ell$. By construction, for a vertex $v$ of depth $\ell-1$, the labels $(u_{vn})_{n\ge1}$
and the families indexed by $n$ formed by the increments after $r_{\ell-1}$ of the path of $vn$
together with all labels and fields of the subtree below $vn$ are independent of
$\mathcal G_{\ell-1}$; these families are i.i.d.\ over $n$ and independent of $(u_{vn})_{n\ge1}$.
Finally, the \emph{Onsager shift} is
\[
 A(q)=\frac N2\int_0^qt\xi''(t)\xc(t)\,dt,\qquad q\in[0,r_L].
\]
A \emph{tree automorphism} is a depth-preserving bijection $\pi$ of
$\{\varnothing\}\cup\mathcal V$ with $\pi(v|k)=\pi(v)|k$ for every $v$ and every $k\le|v|$. It
permutes the leaves, and $\pi(\alpha)\wedge\pi(\beta)=\alpha\wedge\beta$.

The identity of Lemma~\ref{lem:onsager-shift} below is the first identity of
\cite[Remark~2.2]{Bou26b}, where it is stated for parameters $0<\xc_1<\dots<\xc_L\le1$ given by a
finite-step distribution function and is proved by applying the marking identity
\cite[(2.2)]{Bou26b} level by level. We use it with parameters that may exceed one
(Proposition~\ref{prop:glued-comparison}), so we give a complete proof. It needs one property of
Poisson processes, which we state first.

\begin{lemma}[Shifted Poisson points with independent marks]\label{lem:poisson-shift}
Let $\xc>0$, let $u_1>u_2>\cdots$ be the points of a Poisson point process on $\R$ with intensity
$\xc e^{-\xc u}\,du$, and let $\mathcal E$ be a Polish space. Let $(\gamma_n,S_n)_{n\ge1}$ be
i.i.d.\ random elements of $\R\times\mathcal E$, independent of $(u_n)_{n\ge1}$, such that
$\gamma_1$ is independent of $S_1$ and $c:=\xc^{-1}\log\E e^{\xc\gamma_1}<\infty$. Then almost
surely the numbers $u_n+\gamma_n$, $n\ge1$, are pairwise distinct and only finitely many of them
exceed any given level. If $k\mapsto n_k$ is the bijection of $\N$ that lists them in decreasing
order, then
\[
 \bigl(u_{n_k}+\gamma_{n_k}-c,\ S_{n_k}\bigr)_{k\ge1}\overset{d}{=}\bigl(u_k,\ S_k\bigr)_{k\ge1}
 \qquad\text{in }(\R\times\mathcal E)^{\N}.
\]
\end{lemma}

\begin{proof}
\emph{Marking formula.} Let $\mu(du)=\xc e^{-\xc u}\,du$, and let $P_\gamma$, $P_S$ be the laws of
$\gamma_1$, $S_1$. For every measurable $h:\R\times\R\times\mathcal E\to[0,\infty]$,
\begin{equation}\label{eq:marking-laplace}
 \E\exp\Bigl(-\sum_{n\ge1}h(u_n,\gamma_n,S_n)\Bigr)
 =\exp\Bigl(-\int\bigl(1-e^{-h}\bigr)\,d(\mu\otimes P_\gamma\otimes P_S)\Bigr).
\end{equation}
Indeed, given that the Poisson number $K_t$ of points in $(t,\infty)$, of mean $e^{-\xc t}$, is
$k$, these points with their marks are, up to order, $k$ i.i.d.\ triples of law
$\mu_t\otimes P_\gamma\otimes P_S$ with $\mu_t=e^{\xc t}\mu|_{(t,\infty)}$: the points are the
ordered values of $k$ i.i.d.\ draws from $\mu_t$, and the marks are i.i.d.\ and independent of
them. Averaging over $k$ gives \eqref{eq:marking-laplace} for the points in $(t,\infty)$; let
$t\to-\infty$ (dominated convergence on the left, monotone on the right).

\emph{Change of variables.} Let $\Xi=\sum_n\delta_{(u_n,S_n)}$ and
$\Xi'=\sum_n\delta_{(u_n+\gamma_n-c,\,S_n)}$. For measurable $\varphi\ge0$ on $\R\times\mathcal E$
vanishing on some $(-\infty,t]\times\mathcal E$, \eqref{eq:marking-laplace} with
$h(u,g,\varsigma)=\varphi(u+g-c,\varsigma)$, the substitution $w=u+g-c$ and
$\int e^{\xc(g-c)}P_\gamma(dg)=1$, then \eqref{eq:marking-laplace} with
$h(u,g,\varsigma)=\varphi(u,\varsigma)$, give
\[
 \E e^{-\Xi'(\varphi)}
 =\exp\Bigl(-\iint\bigl(1-e^{-\varphi(w,\varsigma)}\bigr)\,\xc e^{-\xc w}\,dw\,P_S(d\varsigma)\Bigr)
 =\E e^{-\Xi(\varphi)}.
\]
For $\varphi=\epsilon\ind_{(t,\infty)\times\mathcal E}$, $\epsilon\downarrow0$, this gives
$\Xi'((t,\infty)\times\mathcal E)<\infty$ almost surely for every $t$; for
$\varphi=\sum_{i\le m}\epsilon_i\ind_{D_i}$ it gives equal laws for $(\Xi(D_i))_{i\le m}$ and
$(\Xi'(D_i))_{i\le m}$ whenever the $D_i$ are Borel subsets of some $(t,\infty)\times\mathcal E$.
Such sets generate the $\sigma$-field of the maps $\omega\mapsto\omega(D)$ on the space $\mathbf M$
of counting measures finite on every $(t,\infty)\times\mathcal E$, since
$\omega(D)=\lim_{t\to-\infty}\omega(D\cap((t,\infty)\times\mathcal E))$; hence
$\Xi'\overset{d}{=}\Xi$ in $\mathbf M$ by a $\pi$-system argument.

\emph{Measurable enumeration.} Almost surely $\Xi$ lies in the set $\mathbf M^*\subset\mathbf M$ of
measures with infinitely many atoms, all of mass one and with pairwise distinct first
coordinates; $\mathbf M^*$ is measurable, and so is the map
$\operatorname{List}:\mathbf M^*\to(\R\times\mathcal E)^{\N}$ listing the atoms by decreasing first
coordinate (both facts are standard). Hence $\Xi'\in\mathbf M^*$ almost surely, which is the first
claim, and $\operatorname{List}(\Xi')=(u_{n_k}+\gamma_{n_k}-c,S_{n_k})_{k\ge1}$ has the law of
$\operatorname{List}(\Xi)=(u_k,S_k)_{k\ge1}$.
\end{proof}

\begin{lemma}[Onsager shift]\label{lem:onsager-shift}
Let $q\in[0,r_L]$ be deterministic.
\begin{thmlist}
\item\label{lem:onsager-shift-a} There is a tree automorphism $\pi$, measurable with respect to
 $\sigma(u,Y)$, such that
 \[
  \bigl(u_{\pi(\alpha)}+Y(q,\pi(\alpha))-A(q)\bigr)_{\alpha\in\N^L}
  \overset{d}{=}\bigl(u_\alpha\bigr)_{\alpha\in\N^L}\qquad\text{in }\R^{\N^L}.
 \]
 In particular $\sup_\alpha\bigl(u_\alpha+Y(q,\alpha)-A(q)\bigr)\overset{d}{=}U$, and the supremum
 on the left is attained with the same probability as the supremum defining $U$.
\item\label{lem:onsager-shift-b} Let $\mathcal H$ be a $\sigma$-field independent of
 $\sigma(u,Y)$, and let $\mathsf q$ be an $\mathcal H$-measurable random variable with values in
 $[0,r_L]$. Then $V=\sup_\alpha\bigl(u_\alpha+Y(\mathsf q,\alpha)-A(\mathsf q)\bigr)$ is independent
 of $\mathcal H$, and $V\overset{d}{=}U$.
\item\label{lem:onsager-shift-c} Let $(u^b,Y^b)_{b\ge1}$ be i.i.d.\ copies of $(u,Y)$, let
 $\mathcal G$ be a $\sigma$-field independent of $\sigma\bigl((u^b,Y^b)_{b\ge1}\bigr)$, and for each
 $b$ let $\mathsf q_b$ be a $\mathcal G$-measurable random variable with values in $[0,r_L]$. Then
 the variables
 $V_b=\sup_\alpha\bigl(u^b_\alpha+Y^b(\mathsf q_b,\alpha)-A(\mathsf q_b)\bigr)$, $b\ge1$, are
 i.i.d.\ with the law of $U$, and $(V_b)_{b\ge1}$ is independent of $\mathcal G$.
\end{thmlist}
\end{lemma}

\begin{proof}
\emph{Step 1: edge marks.} Put $t_\ell=r_\ell\wedge q$ for $0\le\ell\le L$, so that
$0=t_0\le t_1\le\dots\le t_L=q$. For a Borel set $I\subset[0,1]$ write
$|I|_\vth=\int_It\xi''(t)\,dt$. Since $\vth(0)=0$, we have
$\vth(s\wedge t\wedge r)=|[0,s]\cap[0,t]\cap[0,r]|_\vth$, and expanding the products gives, for
$0\le s\le s'\le r_L$, $0\le t\le t'\le r_L$ and leaves $\alpha,\beta$,
\begin{equation}\label{eq:increment-covariance}
 \E\bigl[(Y(s',\alpha)-Y(s,\alpha))(Y(t',\beta)-Y(t,\beta))\bigr]
 =N\,\bigl|(s,s']\cap(t,t']\cap[0,r_{\alpha\wedge\beta}]\bigr|_\vth .
\end{equation}
For $v\in\mathcal V$ with $|v|=\ell$ let $\alpha^v\in\N^L$ be the leaf $v1\cdots1$, and set
\[
 \gamma_v=Y(t_\ell,\alpha^v)-Y(t_{\ell-1},\alpha^v),\qquad
 \sigma_\ell^2=N\bigl(\vth(t_\ell)-\vth(t_{\ell-1})\bigr),\qquad
 A_\ell=\tfrac12\xc_\ell\sigma_\ell^2 .
\]
\begin{romanlist}
\item\emph{Almost surely, $Y(q,\alpha)=\sum_{\ell=1}^L\gamma_{\alpha|\ell}$ for every leaf
 $\alpha$.} First, $Y(0,\alpha)=0$ almost surely, since its variance is $N\vth(0)=0$. If
 $v=\alpha|\ell$, then $\alpha\wedge\alpha^v\ge\ell$, so $r_{\alpha\wedge\alpha^v}\ge r_\ell\ge t_\ell$,
 and for $t\le t_\ell$ the covariance formula gives
 $\E(Y(t,\alpha)-Y(t,\alpha^v))^2=N(2\vth(t)-2\vth(t))=0$. Hence
 $Y(t_\ell,\alpha)-Y(t_{\ell-1},\alpha)=\gamma_{\alpha|\ell}$ almost surely, and summing over $\ell$
 gives the claim for each $\alpha$, hence for all of the countably many leaves at once.
\item\emph{The $\gamma_v$, $v\in\mathcal V$, are independent centred Gaussian variables with
 $\Var\gamma_v=\sigma_{|v|}^2$, and they are independent of the labels.} They are linear
 functionals of $Y$, hence jointly Gaussian, centred and independent of the labels. Let $|v|=\ell$,
 $|w|=k$ and $j=\alpha^v\wedge\alpha^w$. By \eqref{eq:increment-covariance},
 $\E\gamma_v\gamma_w=N\bigl|(t_{\ell-1},t_\ell]\cap(t_{k-1},t_k]\cap[0,r_j]\bigr|_\vth$. If $\ell\ne k$
 the two intervals are disjoint. If $\ell=k$ and $v\ne w$, then $j\le\ell-1$; either
 $q\le r_{\ell-1}$ and $(t_{\ell-1},t_\ell]$ is empty, or $t_{\ell-1}=r_{\ell-1}\ge r_j$ and the
 intersection with $[0,r_j]$ is empty. If $v=w$, then $j=L$ and $r_L\ge t_\ell$, which gives
 $\sigma_\ell^2$. Uncorrelated jointly Gaussian variables are independent.
\item\emph{$\xc_\ell^{-1}\log\E e^{\xc_\ell\gamma_v}=A_\ell$ for $|v|=\ell$, and
 $\sum_{\ell=1}^LA_\ell=A(q)$.} The first identity is the Gaussian moment generating function. For
 the second, $[0,q)$ is the disjoint union of the intervals $[t_{\ell-1},t_\ell)$, and whenever such
 an interval is nonempty, $t_{\ell-1}=r_{\ell-1}$ and $t_\ell\le r_\ell$, so $\xc=\xc_\ell$ on it.
 Hence $\sum_\ell A_\ell=\frac N2\sum_\ell\xc_\ell\int_{t_{\ell-1}}^{t_\ell}t\xi''(t)\,dt=A(q)$.
\end{romanlist}

\emph{Step 2: sorting.} For $L'\ge1$ let $\mathcal V^{L'}=\bigcup_{\ell=1}^{L'}\N^\ell$. Call
$\lambda\in\R^{\mathcal V^{L'}}$ \emph{sortable} if for every vertex $v$ of depth less than $L'$,
including $\varnothing$, the numbers $\lambda_{vn}$, $n\in\N$, are pairwise distinct and only
finitely many of them exceed any given level. For sortable $\lambda$ define $\pi_\lambda$
recursively by $\pi_\lambda(\varnothing)=\varnothing$ and $\pi_\lambda(vn)=\pi_\lambda(v)m$, where
$\lambda_{\pi_\lambda(v)m}$ is the $n$-th largest of the numbers $\lambda_{\pi_\lambda(v)m'}$,
$m'\in\N$. Then $\pi_\lambda$ is a tree automorphism, and we set
$\operatorname{Sort}(\lambda)=\lambda\circ\pi_\lambda$. The sortable vectors form a Borel set, and
$\lambda\mapsto\pi_\lambda$ and $\lambda\mapsto\operatorname{Sort}(\lambda)$ are Borel on it, since
ranks are countable sums of indicators. The labels of a cascade of depth $L'$ are almost surely
sortable, with $\pi=\id$. For $n\in\N$ write $\lambda^{[n]}\in\R^{\mathcal V^{L'-1}}$ for
$w\mapsto\lambda_{nw}$. Then $\lambda$ is sortable if and only if the root values
$(\lambda_n)_{n\ge1}$ satisfy the condition above and every $\lambda^{[n]}$ is sortable, and in that
case
\begin{equation}\label{eq:sort-recursion}
 \operatorname{Sort}(\lambda)_k=\lambda_{n_k},\qquad
 \operatorname{Sort}(\lambda)^{[k]}=\operatorname{Sort}\bigl(\lambda^{[n_k]}\bigr),\qquad k\ge1,
\end{equation}
where $k\mapsto n_k$ lists the root values in decreasing order.

\emph{Claim.} Let $L'\ge1$, let $u'$ be the labels of a cascade of depth $L'$ with arbitrary
parameters $\xc'_1,\dots,\xc'_{L'}>0$, and let $(\gamma'_v)_{v\in\mathcal V^{L'}}$ be independent
of $u'$ and independent centred Gaussian variables with $\Var\gamma'_v=\sigma'^2_{|v|}$ for some
$\sigma'_1,\dots,\sigma'_{L'}\ge0$. Put
$\lambda_v=u'_v+\gamma'_v-\frac12\xc'_{|v|}\sigma'^2_{|v|}$. Then $\lambda$ is almost surely sortable
and $\operatorname{Sort}(\lambda)\overset{d}{=}u'$ in $\R^{\mathcal V^{L'}}$.

We argue by induction on $L'$. For $L'=1$ this is Lemma~\ref{lem:poisson-shift} with
$\gamma_n=\gamma'_n$, $c=\frac12\xc'_1\sigma'^2_1$, and $\mathcal E$ a one-point space. Let $L'\ge2$.
The triples $\bigl(\gamma'_n,u'^{[n]},\gamma'^{[n]}\bigr)$, $n\ge1$, are i.i.d., independent of the
root labels $(u'_n)_{n\ge1}$, and $\gamma'_n$ is independent of $(u'^{[n]},\gamma'^{[n]})$. Each
$u'^{[n]}$ is the label vector of a cascade of depth $L'-1$ with parameters
$\xc'_2,\dots,\xc'_{L'}$, and $\lambda^{[n]}$ is built from $(u'^{[n]},\gamma'^{[n]})$ as in the
claim. By the induction hypothesis, $S_n=\operatorname{Sort}(\lambda^{[n]})$ is almost surely
defined and $S_n\overset{d}{=}u'^{[1]}$. The pairs $(\gamma'_n,S_n)$ are i.i.d., independent of
$(u'_n)$, and $\gamma'_n$ is independent of $S_n$. Lemma~\ref{lem:poisson-shift} with
$\mathcal E=\R^{\mathcal V^{L'-1}}$, $\gamma_n=\gamma'_n$, marks $S_n$ and
$c=\frac12\xc'_1\sigma'^2_1$ shows that the root values $\lambda_n=u'_n+\gamma'_n-c$ are almost
surely pairwise distinct with finitely many above each level, so that $\lambda$ is almost surely
sortable (every $\lambda^{[n]}$ is almost surely sortable by the induction hypothesis), and that
\[
 \bigl(\lambda_{n_k},S_{n_k}\bigr)_{k\ge1}\overset{d}{=}\bigl(u'_k,S_k\bigr)_{k\ge1}
 \overset{d}{=}\bigl(u'_k,u'^{[k]}\bigr)_{k\ge1},
\]
the second equality because both sequences of marks are i.i.d.\ with the same law and independent
of $(u'_k)$. By \eqref{eq:sort-recursion}, $\operatorname{Sort}(\lambda)$ is the image of the left
side under the Borel map that sends a sequence $(y_k,\varsigma_k)_{k\ge1}$ in
$(\R\times\R^{\mathcal V^{L'-1}})^{\N}$ to the vector with root values $y_k$ and subtrees
$\varsigma_k$, and the same map sends the right side to $u'$. This proves the claim.

\emph{Step 3: proof of (a).} Apply the claim with $L'=L$, $\xc'=\xc$, $u'=u$, $\gamma'=\gamma$ and
$\sigma'_\ell=\sigma_\ell$, which Step~1(ii) allows; then $\lambda_v=u_v+\gamma_v-A_{|v|}$. Let
$\pi=\pi_\lambda$ on the event that $\lambda$ is sortable and $\pi=\id$ otherwise; it is
$\sigma(u,Y)$-measurable, because $\gamma$ is a function of $Y$. Let
$\operatorname{Leaf}:\R^{\mathcal V}\to\R^{\N^L}$ be the leaf-sum map,
$\operatorname{Leaf}(\lambda)_\alpha=\sum_{\ell=1}^L\lambda_{\alpha|\ell}$. Since $\pi$ is a tree
automorphism, $\operatorname{Leaf}(\lambda\circ\pi)_\alpha=\operatorname{Leaf}(\lambda)_{\pi(\alpha)}$,
and by Step~1(i) and (iii), $\operatorname{Leaf}(\lambda)_\alpha=u_\alpha+Y(q,\alpha)-A(q)$ for all
$\alpha$, almost surely. Therefore
\[
 \bigl(u_{\pi(\alpha)}+Y(q,\pi(\alpha))-A(q)\bigr)_\alpha
 =\operatorname{Leaf}\bigl(\operatorname{Sort}(\lambda)\bigr)
 \overset{d}{=}\operatorname{Leaf}(u)=(u_\alpha)_\alpha .
\]
The supremum of a family indexed by $\N^L$ is a Borel function of the family, and the set of
families whose supremum is attained is Borel; neither changes when the index set is permuted by
$\pi$. This gives the last assertion of (a).

\emph{Step 4: proof of (b).} Let $\mathsf Y=(Y(\cdot,\alpha))_\alpha$, a random element of
$C([0,r_L])^{\N^L}$. The map
$\mathsf V(s,\mathsf u,\mathsf y)
=\sup_\alpha\bigl(\sum_{\ell}\mathsf u_{\alpha|\ell}+\mathsf y_\alpha(s)-A(s)\bigr)$ is Borel on
$[0,r_L]\times\R^{\mathcal V}\times C([0,r_L])^{\N^L}$, because evaluation $(s,y)\mapsto y(s)$ is
continuous, $A$ is continuous and the supremum is countable. Let $h:(-\infty,\infty]\to\R$ be
bounded and Borel and let $H\in\mathcal H$. Since $(\mathsf q,\ind_H)$ is independent of
$(u,\mathsf Y)$, $V=\mathsf V(\mathsf q,u,\mathsf Y)$ and Fubini's theorem give
$\E[h(V)\ind_H]=\E[\varphi(\mathsf q)\ind_H]$ with $\varphi(s)=\E h(\mathsf V(s,u,\mathsf Y))$. By
(a), $\varphi(s)=\E h(U)$ for every $s$, so $\E[h(V)\ind_H]=\E h(U)\,\Prob(H)$.

\emph{Step 5: proof of (c).} The $\sigma$-fields $\mathcal G$ and $\sigma(u^b,Y^b)$, $b\ge1$, are
mutually independent. Fix $b$ and apply (b) to $(u^b,Y^b)$ with
$\mathcal H_b=\mathcal G\vee\sigma\bigl((u^{b'},Y^{b'}):b'\ne b\bigr)$: $V_b$ has the law of $U$ and
is independent of $\mathcal H_b$, which contains $\mathcal G$ and every $V_{b'}$ with $b'\ne b$.
For distinct $b_1,\dots,b_k$, Borel sets $D_1,\dots,D_k$ and $\Omega_0\in\mathcal G$, applying
this successively to $b_1,b_2,\dots$ gives
$\Prob\bigl(\Omega_0\cap\bigcap_{i\le k}\{V_{b_i}\in D_i\}\bigr)
=\Prob(\Omega_0)\prod_{i\le k}\Prob(U\in D_i)$.
\end{proof}

The second cascade fact bounds exponential moments of a cascade maximum to which marks are added
level by level. In the proof of Proposition~\ref{prop:glued-comparison} it plays the role of
\cite[Lemma~2.5]{Bou26b}; it needs strictly increasing parameters, but no bound by one.

\begin{lemma}[Cascade domination]\label{lem:cascade-domination}
Assume $0<\xc_1<\dots<\xc_L$. For $1\le\ell\le L$ and every vertex $w$ of depth $\ell$ let
$\delta_\ell(w)$ be a measurable function of the path of $w$, and let $c_\ell\in\R$ be
deterministic with
\[
 \E\bigl[e^{\xc_\ell\delta_\ell(vn)}\bigm|\mathcal G_{\ell-1}\bigr]\le e^{\xc_\ell c_\ell}
 \quad\text{almost surely, for every vertex $v$ of depth $\ell-1$ and every $n$.}
\]
Then, for every $\lambda\in(0,\xc_1)$,
\[
 \E\exp\Bigl(\lambda\sup_{\alpha\in\N^L}\Bigl(u_\alpha+\sum_{\ell=1}^L\delta_\ell(\alpha|\ell)\Bigr)\Bigr)
 \le e^{\lambda\sum_\ell c_\ell}\,\E e^{\lambda U}<\infty .
\]
\end{lemma}

\begin{proof}
For a vertex $w$ of depth $\ell<L$ let
$R_w=\sup_{\alpha:\,\alpha|\ell=w}\sum_{k>\ell}\bigl(u_{\alpha|k}+\delta_k(\alpha|k)\bigr)$, and
$R_w=0$ for leaves; $R_\varnothing$ is the supremum in the statement. For $0\le\ell\le L$ let
$U^{(\ell)}$ be the supremum of the leaf labels of a cascade with the levels $\ell+1,\dots,L$ and
parameters $\xc_{\ell+1},\dots,\xc_L$, so that $U^{(L)}=0$ and $U^{(0)}=U$, and put
$c_{>\ell}=\sum_{k>\ell}c_k$ and $\xc_{L+1}=\infty$. We show by downward induction on $\ell$ that,
for every $\lambda\in(0,\xc_{\ell+1})$, $\E e^{\lambda U^{(\ell)}}<\infty$ and
\[
 \E\bigl[e^{\lambda R_w}\bigm|\mathcal G_\ell\bigr]\le e^{\lambda c_{>\ell}}\,\E e^{\lambda U^{(\ell)}}
 \qquad\text{for every vertex $w$ of depth $\ell$.}
\]
For $\ell=L$ there is nothing to prove. Let $1\le\ell\le L$, assume the statement at $\ell$, and fix
a vertex $v$ of depth $\ell-1$. Then $R_v=\sup_n(u_{vn}+Y_n)$ with $Y_n=\delta_\ell(vn)+R_{vn}$. By
the independence structure of the construction, conditionally on $\mathcal G_{\ell-1}$ the $Y_n$
are i.i.d.\ and independent of the labels $(u_{vn})_n$, which form a Poisson point process with
intensity $\xc_\ell e^{-\xc_\ell u}\,du$. By the induction hypothesis at $\lambda=\xc_\ell<\xc_{\ell+1}$
and the assumption on $\delta_\ell$,
\[
 \E\bigl[e^{\xc_\ell Y_n}\bigm|\mathcal G_{\ell-1}\bigr]
 =\E\Bigl[e^{\xc_\ell\delta_\ell(vn)}\,\E\bigl[e^{\xc_\ell R_{vn}}\bigm|\mathcal G_\ell\bigr]\Bigm|\mathcal G_{\ell-1}\Bigr]
 \le e^{\xc_\ell c_{\ge\ell}}\,\E e^{\xc_\ell U^{(\ell)}}<\infty ,
\]
where $c_{\ge\ell}=c_\ell+c_{>\ell}$. Lemma~\ref{lem:poisson-shift} (with a one-point mark space),
applied conditionally on $\mathcal G_{\ell-1}$, shows that conditionally on $\mathcal G_{\ell-1}$
the supremum $R_v$ is attained and $R_v-\varkappa'_v$ has the law of the largest point $u_1$ of
that Poisson process, where
$\varkappa'_v=\xc_\ell^{-1}\log\E[e^{\xc_\ell Y_1}\mid\mathcal G_{\ell-1}]\le c_{\ge\ell}+\varkappa_\ell$
and $\varkappa_\ell=\xc_\ell^{-1}\log\E e^{\xc_\ell U^{(\ell)}}$. The same lemma, applied with i.i.d.\
copies of $U^{(\ell)}$ in place of the $Y_n$, gives $U^{(\ell-1)}\overset{d}{=}\varkappa_\ell+u_1$.
Since $\Prob(u_1>t)\le e^{-\xc_\ell t}$, we have $\E e^{\lambda u_1}<\infty$, hence
$\E e^{\lambda U^{(\ell-1)}}<\infty$, for $\lambda<\xc_\ell$, and then
\[
 \E\bigl[e^{\lambda R_v}\bigm|\mathcal G_{\ell-1}\bigr]=e^{\lambda\varkappa'_v}\,\E e^{\lambda u_1}
 \le e^{\lambda c_{\ge\ell}}e^{\lambda\varkappa_\ell}\,\E e^{\lambda u_1}
 =e^{\lambda c_{\ge\ell}}\,\E e^{\lambda U^{(\ell-1)}} .
\]
This is the statement at $\ell-1$. The case $\ell=0$, after taking expectations, is the lemma.
\end{proof}

\begin{remark}\label{rem:onsager-shift-scope}
The proof of Lemma~\ref{lem:onsager-shift} uses only $\xc_\ell>0$ and $\vth'\ge0$ on $[0,r_L]$; it
needs no ordering of the parameters, no bound $\xc_\ell\le1$ and no relation to a Parisi measure.
If $\xc_1<\dots<\xc_L$, then $U$ is almost surely finite and attained, so the supremum in
Lemma~\ref{lem:onsager-shift}\ref{lem:onsager-shift-a} is a maximum (for
\ref{lem:onsager-shift-b}--\ref{lem:onsager-shift-c} apply the argument of its Step~4 to the
indicator of attainment). By induction on $L$: for $L=1$, $U$ is the largest Poisson
point; for $L\ge2$, the subtree suprema $U^{[n]}$ are i.i.d., independent of the root labels, almost
surely finite and attained, and $\E e^{\xc_1U^{[1]}}<\infty$ by Lemma~\ref{lem:cascade-domination}
for the subtree (parameters $\xc_2<\dots<\xc_L$, exponent $\xc_1<\xc_2$); so by
Lemma~\ref{lem:poisson-shift} with $\gamma_n=U^{[n]}$ only finitely many of the $u_n+U^{[n]}$ exceed
any level, and $U=\sup_n(u_n+U^{[n]})$ is attained.
\end{remark}

\subsection{Proof of the glued comparison}\label{app:glued-proof}

\paragraph{Data and statement.}
The hypotheses and notation of Proposition~\ref{prop:glued-comparison} are in force (the levels
$q_\ell$ and values $\xc_\ell$ of $\xc$, the top value $y$, the band measure $\eta$, the
perturbation $(s,g_0)$), with $\Phi_\eta$ and $\Phistar_\eta(Q,\cdot)$ as in
Section~\ref{sec:loc-glued}, and $\mathcal S=\{m\in\cube:\ Q\le q_m\le Q+\rho\}$ is the shell, a
sphere if $\rho=0$. For a step coefficient $c$ on $[t_0,t_1)$ with values in $[0,3]$ and a terminal
function $\psi$ at time $t_1$, $W[c,\psi]$ is the Cole--Hopf solution~\eqref{eq:pde-colehopf}: it
equals $\psi$ at $t_1$, and on a step $[\tau',\tau)$ of $c$, where $c$ takes the value $c'>0$,
\begin{equation}\label{eq:glued-colehopf}
 W(t,z)=\frac1{c'}\log\E\exp\Bigl(c'\,W\bigl(\tau,z+\sqrt{\xi'(\tau)-\xi'(t)}\,Z\bigr)\Bigr),
 \qquad t\in[\tau',\tau),\ z\in\R,
\end{equation}
with $Z$ a standard Gaussian variable (the heat flow $\E W(\tau,z+\sqrt{\xi'(\tau)-\xi'(t)}\,Z)$
if $c'=0$). The objects of the proposition are
\begin{align}
 \Psi_s(z)&=\sup_{|m|\le1}\bigl[zm+\Phistar_\eta(Q,m)+s\,g_0(m)\bigr],\qquad z\in\R,
 \label{eq:glued-terminal}\\
 P_s(\xc,\eta)&=W_s(0,0)-\frac12\int_0^Qt\xi''(t)\xc(t)\,dt
 -\frac12\int_Q^1t\xi''(t)\cdf\eta(t)\,dt,\label{eq:glued-Ps}
\end{align}
where $W_s=W[\xc,\Psi_s]$ on $[0,Q]$, and $W_s(0,0)=\Psi_s(0)$ if $Q=0$. For $m\in\mathcal S$ let
\begin{equation}\label{eq:glued-T}
 T(m)=\HN(m)+\sum_{i=1}^N\bigl[\Phistar_\eta(Q,m_i)+s\,g_0(m_i)\bigr]-NW_s(0,0)+A(q_m),
\end{equation}
where
\[
 A(q)=\frac N2\int_0^qt\xi''(t)\bar\xc(t)\,dt,\qquad
 \bar\xc=\xc\ \text{on }[0,Q),\qquad \bar\xc=y\ \text{on }[Q,Q+\rho],
\]
and let
\begin{equation}\label{eq:glued-majorant}
 \mathcal T=\sup_{m\in\mathcal S}T(m)+\frac N2\xi''(1)\rho .
\end{equation}
Step~0 below shows that $T$ is almost surely continuous on the compact set $\mathcal S$, so
$\mathcal T$ is a supremum over a countable dense subset and is a random variable. We prove that,
for every $N\ge1$,
\begin{gather}
 \FTAP(m)+s\sum_{i=1}^Ng_0(m_i)-NP_s(\xc,\eta)\le\mathcal T
 \qquad\text{for every }m\in\mathcal S,\label{eq:glued-pointwise}\\
 \E e^{\theta\mathcal T}\le\exp\Bigl\{\theta\Bigl(NC_R\rho
 +\ind_{\{\rho>0\}}\frac{\log(N+1)}{y}\Bigr)\Bigr\},\qquad C_R=5\xi''(1),
 \label{eq:glued-moment}
\end{gather}
which is Proposition~\ref{prop:glued-comparison}.

\paragraph{Plan.}
Step~0 is deterministic. It proves \eqref{eq:glued-pointwise} and reduces
\eqref{eq:glued-moment} to the bound
\begin{equation}\label{eq:glued-goal}
 \log\E\exp\Bigl(\theta\sup_{m\in\mathcal S}T(m)\Bigr)\le\theta d_N,\qquad
 d_N=\ind_{\{\rho>0\}}\Bigl[\frac{\log(N+1)}y+N\kappa_\rho+\frac{y(N+1)\kappa_\rho}2\Bigr],
\end{equation}
where $\kappa_\rho=\xi'(Q+\rho)-\xi'(Q)$. To prove it we compare two families of Gaussian
processes, indexed by a copy $b$ of the system, a point $m$ of a finite set $D\subset\mathcal S$
and a leaf $\alpha$ of a Ruelle cascade whose levels are those of $\xc$, extended to the top of
the shell. The first family, $X_b$, adds to $\HN(m)$ an Onsager field at self-overlap $q_m$; the
second, $\widetilde X_b$, replaces both by a cavity field that is linear in $m$. An outer
Poisson level $(w_b)$ of parameter $\theta$ turns expectations of maxima into exponential
moments of order $\theta$. With $U$ the supremum of the cascade labels, the proof establishes
\begin{equation}\label{eq:glued-chain}
 \begin{split}
 \E\max_bw_b+\frac{\log\E e^{\theta\max_DT}+\log\E e^{\theta U}}{\theta}
 &\le\E\sup_{b,m,\alpha}X_b\le\E\sup_{b,m,\alpha}\widetilde X_b\\
 &\le\E\max_bw_b+d_N+\frac{\log\E e^{\theta U}}{\theta}.
 \end{split}
\end{equation}
The middle inequality is Slepian's lemma (Step~3). The right one combines Fenchel's inequality
with an exact computation below $Q$ and Doob's inequality on the shell (Steps~4--5). The left
one is the Onsager shift of Lemma~\ref{lem:onsager-shift} (Step~6). The common terms cancel, and
$D$ increases to a dense subset of $\mathcal S$ (Step~7). We first record the properties of the
perturbed terminal $\Psi_s$ that the recursion needs.

\begin{lemma}[The perturbed terminal]\label{lem:glued-terminal}
In the setting above, $\Psi_0=\Phi_\eta(Q,\cdot)$. For every $(s,g_0)$ as above, the function
$\Psi_s$ is convex and twice continuously differentiable, and
\[
 |\Psi_s'|<1,\qquad 0<\Psi_s''\le2\bigl(1-(\Psi_s')^2\bigr).
\]
Consequently Lemma~\ref{lem:pde-class} applies to $W_s=W[\xc,\Psi_s]$ with $G=3$: for every
$t\in[0,Q]$, $W_s(t,\cdot)$ is convex, $|\partial_zW_s(t,\cdot)|<1$ and
$|W_s(t,z)-W_s(t,0)|\le|z|$.
\end{lemma}

\begin{proof}
The first two assertions are Lemma~\ref{lem:pde-fenchel}\ref{lem:pde-fenchel-c}, since
$\mathcal T_2$ is the class of convex $C^2$ functions $\psi$ with $|\psi'|<1$ and
$\psi''\le2(1-\psi'^2)$. As $\mathcal T_2\subset\mathcal T_3$ and $\xc$ takes values in $(0,3]$,
the last assertion is Lemma~\ref{lem:pde-class}\ref{lem:pde-class-a}
and~\ref{lem:pde-class-d}.
\end{proof}

\begin{proof}[Proof of Proposition~\ref{prop:glued-comparison}]
\emph{Step 0: deterministic reduction.}
Let $m\in\mathcal S$ and $Q'=q_m\in[Q,Q+\rho]$. The measure
$\Shift_{Q'}\eta=\eta|_{(Q',1]}+\eta([Q,Q'])\delta_{Q'}\in\Pspace([Q',1])$ has distribution
function $\cdf\eta$ on $[Q',1]$, so its solution of the Parisi equation coincides with
$\Phi_\eta$ on $[Q',1]$. Since $\FTAP$ is an infimum over $\Pspace([q_m,1])$,
\[
 \FTAP(m)\le F_{\Shift_{Q'}\eta}(m)
 =\HN(m)+\sum_{i=1}^N\Phistar_\eta(Q',m_i)-\frac N2\int_{Q'}^1t\xi''(t)\cdf\eta(t)\,dt .
\]
By Lemma~\ref{lem:pde-class}, on each step
$\partial_t\Phi_\eta=-\frac{\xi''}2(\partial_{xx}\Phi_\eta+\cdf\eta\,(\partial_x\Phi_\eta)^2)\le0$,
so $\Phi_\eta(Q',\cdot)\le\Phi_\eta(Q,\cdot)$ and therefore
$\Phistar_\eta(Q',\cdot)\le\Phistar_\eta(Q,\cdot)$. Moreover
$\int_Q^{Q'}t\xi''\cdf\eta\le\xi''(1)\rho$, since $\xi''$ is nondecreasing on $[0,1]$, and
$A(q_m)\ge A(Q)=\frac N2\int_0^Qt\xi''\xc$, since $\bar\xc\ge0$. Comparing with
\eqref{eq:glued-Ps} and \eqref{eq:glued-T},
\begin{equation}\label{eq:glued-step0}
 \FTAP(m)+s\sum_{i=1}^Ng_0(m_i)-NP_s(\xc,\eta)\le T(m)+\frac N2\xi''(1)\rho\le\mathcal T,
\end{equation}
which is \eqref{eq:glued-pointwise}. The function $\Phistar_\eta(Q,\cdot)$ is concave, hence
lower semicontinuous on the closed interval $[-1,1]$, and upper semicontinuous as an infimum of
affine functions, so it is continuous. Since $\HN$ has almost surely continuous paths under
\condH, $T$ is almost surely continuous on $\mathcal S$. Given \eqref{eq:glued-goal},
\eqref{eq:glued-moment} follows from \eqref{eq:glued-majorant}, $\kappa_\rho\le\xi''(1)\rho$
and $y\le3$:
\[
 \frac1N\log\E e^{\theta\mathcal T}\le\theta\Bigl[\frac12\xi''(1)\rho+\frac{d_N}N\Bigr]
 \le\theta\Bigl[\Bigl(\frac12+1+\frac y2\Bigl(1+\frac1N\Bigr)\Bigr)\xi''(1)\rho
 +\ind_{\{\rho>0\}}\frac{\log(N+1)}{yN}\Bigr],
\]
and $\frac12+1+\frac y2(1+\frac1N)\le\frac92$. It remains to prove \eqref{eq:glued-goal}.

\medskip\noindent\emph{Step 1: the cascade.}
If $Q>0$, set $r_\ell=q_\ell$ for $0\le\ell<L$ and $r_L=Q+\rho$. Level $\ell$ is
$[r_{\ell-1},r_\ell]$ with parameter $\xc_\ell$, so the top level $[q_{L-1},Q+\rho]$ has
parameter $y$. If $Q=0$, use a single level: $L=1$, $r_1=\rho$, with parameter $\xc_1=y=1$. In
both cases $\theta<\xc_1<\dots<\xc_L\le3$, $\bar\xc=\xc_\ell$ on $[r_{\ell-1},r_\ell)$, and
$r_L\le1$. Let $(u_\alpha)_{\alpha\in\N^L}$ be the labels of the additive Ruelle probability
cascade with these levels and parameters, let $Y$ be the Onsager field, both as in
Appendix~\ref{app:pde-cascade} (they are the objects of \cite[Definitions~2.2
and~2.3]{Bou26b}), and let $U=\sup_\alpha u_\alpha$. Let $B^i_s(\alpha):=B^i(s,\alpha)$, $1\le i\le N$,
$s\in[0,r_L]$, be the cavity fields of Appendix~\ref{app:pde-cascade}
(\cite[Definition~2.3]{Bou26b}), and take the cascade and all fields independent of $\HN$. The
$B^i$ are centred Gaussian processes with continuous paths and
\[
 \E B^i_s(\alpha)B^j_t(\beta)=\delta_{ij}\,\xi'(s\wedge t\wedge r_{\alpha\wedge\beta}).
\]
They are built like $Y$, from independent Brownian motions attached to the vertices of the
tree, with $\vth$ replaced by $\xi'$ and one family for each $i$. Along a leaf, each
$B^i(\alpha)$ is a Gaussian martingale with $\langle B^i(\alpha)\rangle_t=\xi'(t)$, and two
leaves with $\alpha\wedge\beta\ge\ell$ have the same paths on $[0,r_\ell]$. Let $\mathcal G_\ell$
be the $\sigma$-field generated by the labels at levels at most $\ell$ and all paths on
$[0,r_\ell]$, as in Lemma~\ref{lem:cascade-domination}. Given $\mathcal G_{\ell-1}$, the
increments of the paths on $[r_{\ell-1},r_\ell]$ along the children of a vertex at level
$\ell-1$ are independent of each other and of $\mathcal G_{\ell-1}$, and the cavity-field
components of each are $N$ independent continuous Gaussian processes with independent increments and quadratic variation
$\xi'(t)-\xi'(r_{\ell-1})$.

\medskip\noindent\emph{Step 2: the outer level.}
Fix a finite set $D\subset\mathcal S$. Let $(w_b)_{b\ge1}$ be the points of a Poisson point
process on $\R$ with intensity $\theta e^{-\theta w}\,dw$, and let $(\HN^b,u^b,Y^b,B^b)$,
$b\ge1$, be independent copies of $(\HN,u,Y,B)$, independent of $(w_b)$. This outer level is
Boursier's device for exponential moments \cite[Proposition~2.7, Steps~1--2]{Bou26b}, the
mechanism of \cite{CGKLM}, as she notes. For $b\ge1$, $m\in D$ and $\alpha\in\N^L$ set
\begin{align*}
 X_b(m,\alpha)&=w_b+\HN^b(m)+Y^b(q_m,\alpha)+u^b_\alpha+\Upsilon(m),\\
 \widetilde X_b(m,\alpha)&=w_b+Z^b(m,\alpha)+u^b_\alpha+\Upsilon(m),
\end{align*}
where $\Upsilon(m)=\sum_{i}[\Phistar_\eta(Q,m_i)+s\,g_0(m_i)]-NW_s(0,0)$ and
$Z^b(m,\alpha)=\sum_im_iB^{b,i}_{q_m}(\alpha)$. These are defined because $q_m\le Q+\rho=r_L$.

\medskip\noindent\emph{Step 3: Slepian's lemma.}
Condition on all labels $(w_b)$ and $(u^b)$. The two families are then Gaussian with the same
mean $w_b+u^b_\alpha+\Upsilon(m)$. Their variances agree:
$\E(\HN(m)+Y(q_m,\alpha))^2=N\xi(q_m)+N\vth(q_m)=Nq_m\xi'(q_m)=\E Z(m,\alpha)^2$. For two index
points in the same copy, with $R_{12}=\ov{m^1}{m^2}\in[-1,1]$ and
$Q_{12}=q_{m^1}\wedge q_{m^2}\wedge r_{\alpha^1\wedge\alpha^2}\in[0,1]$, the covariance of the
first family exceeds that of the second by
\[
 N\bigl[\xi(R_{12})+\vth(Q_{12})-R_{12}\xi'(Q_{12})\bigr]
 =N\bigl[\xi(R_{12})-\xi(Q_{12})-\xi'(Q_{12})(R_{12}-Q_{12})\bigr]\ge0 .
\]
This is where convexity of $\xi$ enters; since $R_{12}$ may be negative, it is needed on all of
$[-1,1]$, and odd degrees are allowed. For index points in different copies both covariances
vanish. By Slepian's lemma, on every finite set of indices the maximum of the first family is
stochastically dominated by the maximum of the second. Take finite index sets increasing to
$\{(b,m,\alpha):b\ge1,\ m\in D,\ \alpha\in\N^L\}$ and use monotone convergence (the maxima are
bounded below by the value at one index point, which is integrable):
\begin{equation}\label{eq:glued-slepian}
 \E\sup_{b,\,m\in D,\,\alpha}X_b(m,\alpha)\le\E\sup_{b,\,m\in D,\,\alpha}\widetilde X_b(m,\alpha).
\end{equation}

\medskip\noindent\emph{Step 4: Fenchel's inequality.}
By \eqref{eq:glued-terminal}, $m_iz+\Phistar_\eta(Q,m_i)+s\,g_0(m_i)\le\Psi_s(z)$ for every
$z\in\R$. Let
\[
 S_t(\alpha)=\sum_{i=1}^NW_s\bigl(t,B^i_t(\alpha)\bigr)\quad(0\le t\le Q),\qquad
 S_t(\alpha)=\sum_{i=1}^N\Psi_s\bigl(B^i_t(\alpha)\bigr)\quad(Q\le t\le Q+\rho).
\]
The two definitions agree at $t=Q$ because $W_s(Q,\cdot)=\Psi_s$, and $S_0=NW_s(0,0)$ because
$B_0=0$. Then $Z(m,\alpha)+\Upsilon(m)\le S_{q_m}(\alpha)-S_0$, and since
$q_m\in[Q,Q+\rho]$,
\begin{equation}\label{eq:glued-fenchel}
 \widetilde X_b(m,\alpha)\le w_b+u^b_\alpha+\sup_{Q\le t\le Q+\rho}\bigl(S^b_t(\alpha)-S^b_0\bigr)
 =w_b+u^b_\alpha+\sum_{\ell=1}^L\delta^b_\ell(\alpha|\ell),
\end{equation}
where $\delta_\ell=S_{r_\ell}-S_{r_{\ell-1}}$ for $\ell<L$ and
$\delta_L=\sup_{Q\le t\le Q+\rho}S_t-S_{r_{L-1}}$. Leaves through the vertex $\alpha|\ell$ share
their paths on $[0,r_\ell]$, so $\delta_\ell$ is a measurable functional of the path of
$\alpha|\ell$ on $[0,r_\ell]$.

\medskip\noindent\emph{Step 5: the cascade side.}
We check the hypothesis of Lemma~\ref{lem:cascade-domination}: for every vertex $v$ at level
$\ell-1$ and every child $vn$,
\begin{equation}\label{eq:glued-levels}
 \E\bigl[e^{\xc_\ell\,\delta_\ell(vn)}\bigm|\mathcal G_{\ell-1}\bigr]\le e^{\xc_\ell c_\ell},
 \qquad c_\ell=0\ \ (\ell<L),\qquad c_L=d_N .
\end{equation}

\emph{Exact levels.} Let $\ell<L$; or let $\ell=L$ and $Q>0$, with $r_L$ replaced by $Q$ in what
follows. Given $\mathcal G_{\ell-1}$, the starting points $B^i_{r_{\ell-1}}(vn)$ are known, and
the increments $B^i_{r_\ell}(vn)-B^i_{r_{\ell-1}}(vn)$, $1\le i\le N$, are independent
$N(0,\xi'(r_\ell)-\xi'(r_{\ell-1}))$ variables, independent of $\mathcal G_{\ell-1}$. On
$[r_{\ell-1},r_\ell)$ the coefficient $\xc$ equals $\xc_\ell$, so the definition
\eqref{eq:glued-colehopf} of $W_s$, applied coordinate by coordinate, gives
\begin{equation}\label{eq:glued-exact}
 \E\bigl[e^{\xc_\ell(S_{r_\ell}-S_{r_{\ell-1}})(vn)}\bigm|\mathcal G_{\ell-1}\bigr]=1 ;
\end{equation}
the exponential moments are finite because $|W_s(t,z)-W_s(t,0)|\le|z|$
(Lemma~\ref{lem:glued-terminal}). This proves \eqref{eq:glued-levels} for $\ell<L$, and for
$\ell=L$ when $\rho=0$. Identity \eqref{eq:glued-exact} is the level-wise identity behind the
second identity of \cite[Remark~2.2]{Bou26b} (her martingale step, p.~9), here read off from the Cole--Hopf recursion with an arbitrary
terminal value at $Q$.

\emph{Doob's inequality on the shell.} Let $\rho>0$ and $\ell=L$; at this level $vn$ is a leaf
$\alpha$. Let $\mathcal G_Q$ be generated by $\mathcal G_{L-1}$ and all paths on $[0,Q]$, and
$\Delta=\sup_{Q\le t\le Q+\rho}(S_t-S_Q)$, so that $\delta_L=(S_Q-S_{r_{L-1}})+\Delta$. (If
$Q=0$, then $S_Q-S_{r_{L-1}}=0$ and $\mathcal G_Q$ is trivial.) By \eqref{eq:glued-exact} on
$[q_{L-1},Q]$, it suffices to show $\E[e^{y\Delta}\mid\mathcal G_Q]\le e^{yd_N}$. Given
$\mathcal G_Q$, the processes $(B^i_t-B^i_Q)_{t\in[Q,Q+\rho]}$ along $\alpha$ are independent
continuous Gaussian martingales with $d\langle B^i\rangle_t=\xi''(t)\,dt$. By It\^o's formula,
\[
 S_t-S_Q=M_t+\frac12\int_Q^t\xi''(r)\sum_{i=1}^N\Psi_s''(B^i_r)\,dr,\qquad
 M_t=\sum_{i=1}^N\int_Q^t\Psi_s'(B^i_r)\,dB^i_r .
\]
By Lemma~\ref{lem:glued-terminal}, $|\Psi_s'|<1$ and $0<\Psi_s''\le2$, so
$\langle M\rangle_t\le N\kappa_\rho$ and the drift term is at most $N\kappa_\rho$, where
$\kappa_\rho=\xi'(Q+\rho)-\xi'(Q)\le\xi''(1)\rho$. Let
$\mathcal E_t=\exp(yM_t-\frac{y^2}2\langle M\rangle_t)$, a martingale since $\langle M\rangle$
is bounded, and $p=1+\frac1N$. Then
$e^{y\Delta}\le e^{yN\kappa_\rho+y^2N\kappa_\rho/2}\sup_t\mathcal E_t$. Writing
$\mathcal E^p=\exp(pyM-\frac{p^2y^2}2\langle M\rangle)\exp(\frac{(p^2-p)y^2}2\langle M\rangle)$,
whose first factor has conditional expectation at most one, Doob's $L^p$ inequality gives
\[
 \E\Bigl[\sup_t\mathcal E_t\Bigm|\mathcal G_Q\Bigr]
 \le\frac p{p-1}\Bigl(\E\bigl[\mathcal E_{Q+\rho}^p\bigm|\mathcal G_Q\bigr]\Bigr)^{1/p}
 \le(N+1)\,e^{(p-1)y^2N\kappa_\rho/2}=(N+1)\,e^{y^2\kappa_\rho/2}.
\]
Hence $\E[e^{y\Delta}\mid\mathcal G_Q]\le(N+1)\exp\{yN\kappa_\rho+\frac{y^2(N+1)\kappa_\rho}2\}
=e^{yd_N}$, which is \eqref{eq:glued-levels} for $\ell=L$. This is the computation in the proof
of \cite[Lemma~2.5]{Bou26b}, performed on the shell only.

With \eqref{eq:glued-levels} and $\theta<\xc_1$, Lemma~\ref{lem:cascade-domination} gives, for
$U''=\sup_\alpha\bigl(u_\alpha+\sum_{\ell=1}^L\delta_\ell(\alpha|\ell)\bigr)$,
\begin{equation}\label{eq:glued-cascade}
 \E e^{\theta U''}\le e^{\theta d_N}\,\E e^{\theta U}<\infty .
\end{equation}

\medskip\noindent\emph{Step 6: the lower side.}
By \eqref{eq:glued-T},
$X_b(m,\alpha)=w_b+T^b(m)+\bigl[u^b_\alpha+Y^b(q_m,\alpha)-A(q_m)\bigr]$, where $T^b$ is $T$
with $\HN$ replaced by $\HN^b$. Let $m^b_*$ be the first maximizer of $T^b$ over $D$ in a fixed
enumeration of $D$, and $q^b_*=q_{m^b_*}\in[Q,Q+\rho]\subset[0,r_L]$. Then
\[
 \sup_{m\in D,\,\alpha}X_b(m,\alpha)\ge w_b+\max_DT^b+V_b,\qquad
 V_b=\sup_\alpha\bigl(u^b_\alpha+Y^b(q^b_*,\alpha)-A(q^b_*)\bigr).
\]
The function $A$ is the shift of Lemma~\ref{lem:onsager-shift} for our cascade, because
$\bar\xc=\xc_\ell$ on $[r_{\ell-1},r_\ell)$. Apply part~(c) of that lemma with
$\mathcal G=\sigma\bigl((\HN^b)_{b\ge1},(w_b)_{b\ge1}\bigr)$, which is independent of the copies
$(u^b,Y^b)_{b\ge1}$, and with the $\mathcal G$-measurable radii $\mathsf q_b=q^b_*$: the
variables $V_b$ are i.i.d.\ with the law of $U$ and independent of $\mathcal G$. Since the
$\HN^b$ are i.i.d.\ and independent of $(w_b)$, the variables $\max_DT^b+V_b$, $b\ge1$, are
i.i.d., independent of $(w_b)$, and $\max_DT^b$ is independent of $V_b$.

\medskip\noindent\emph{Step 7: conclusion.}
We use the outer marking identity: if $(\chi_b)_{b\ge1}$ are i.i.d., independent of $(w_b)$,
and $\E e^{\theta\chi_1}<\infty$, then
\begin{equation}\label{eq:glued-outer}
 \E\sup_b(w_b+\chi_b)=\E\max_bw_b+\frac1\theta\log\E e^{\theta\chi_1},
\end{equation}
by Lemma~\ref{lem:poisson-shift} with $\xc=\theta$ and a one-point mark space; this is
\cite[(2.11)]{Bou26b}, and $\E\max_bw_b$ is finite. On the upper side take $\chi_b=U''_b$, the
variable $U''$ of Step~5 for the $b$-th copy. The $U''_b$ are i.i.d.\ and independent of
$(w_b)$, so \eqref{eq:glued-fenchel}, \eqref{eq:glued-outer} and \eqref{eq:glued-cascade} give
\[
 \E\sup_{b,\,m\in D,\,\alpha}\widetilde X_b
 \le\E\max_bw_b+d_N+\frac1\theta\log\E e^{\theta U}.
\]
On the lower side take $\chi_b=\max_DT^b+V_b$. Here $\E e^{\theta\max_DT}<\infty$ because $D$ is
finite, and $\E e^{\theta U}<\infty$ by \eqref{eq:glued-cascade}; Step~6 and
\eqref{eq:glued-outer} give
\[
 \E\sup_{b,\,m\in D,\,\alpha}X_b
 \ge\E\max_bw_b+\frac1\theta\log\E e^{\theta\max_DT}+\frac1\theta\log\E e^{\theta U}.
\]
These two bounds and \eqref{eq:glued-slepian} form the chain \eqref{eq:glued-chain}. All its
terms are finite, so $\log\E e^{\theta\max_DT}\le\theta d_N$. Let $D$ increase along finite sets
to a countable dense subset of $\mathcal S$. By the continuity of $T$ (Step~0) and monotone
convergence, $\log\E e^{\theta\sup_{\mathcal S}T}\le\theta d_N$, which is
\eqref{eq:glued-goal}.
\end{proof}

\subsection{Step approximation and the measurable majorant}\label{app:glued-approx}

Proposition~\ref{prop:glued-comparison} is proved for step data ($\xc$ and $\cdf\eta$ step
functions), with $(s,g_0)$ as there; $P_s$, given by~\eqref{eq:glued-Ps}, makes sense for every step
$\xc$ with values in $[0,3]$. Since $\Psi_0=\Phi_\eta(Q,\cdot)$ (Lemma~\ref{lem:glued-terminal}) and
the recursion~\eqref{eq:pde-colehopf} on $[0,Q]$ with this terminal value continues the one on
$[Q,1]$, $P_0(\xc,\eta)=P(\glue(\xc,\eta))$ for step data. The next lemma approximates the data of
Section~\ref{sec:loc}, which need not be step data.

\begin{lemma}[Approximation by step data]\label{lem:pde-approx}
Let $\theta\in(0,1)$, $Q\in[0,1)$, $\xc:[0,Q)\to[\theta,2]$ nondecreasing, $\eta\in\Pspace([Q,1])$,
and $(s,g_0)$ as above. For $k\ge1$ let $\xc_k(t)=\xc(jQ/k)+1/k$ for $t\in[jQ/k,(j+1)Q/k)$,
$0\le j<k$, and let $\eta_k$ be the image of $\eta$ under
$t\mapsto Q+(1-Q)\lceil k(t-Q)/(1-Q)\rceil/k$.
\begin{thmlist}
\item\label{lem:pde-approx-a} $\xc_k$ is a nondecreasing step function with values in
 $(\theta,3]$, and $\|\xc_k-\xc\|_{L^1}\le\frac Qk\bigl(\xc(Q-)-\xc(0)+1\bigr)$.
\item\label{lem:pde-approx-b} $\eta_k\in\Pspace([Q,1])$ has a step distribution function, and
 $\|\cdf{\eta_k}-\cdf\eta\|_{L^1}=\Wone(\eta_k,\eta)\le(1-Q)/k$.
\item\label{lem:pde-approx-c} For step data $(\xc_1',\eta_1')$ and $(\xc_2',\eta_2')$ with cascade
 coefficients taking values in $[0,3]$ (not necessarily monotone),
 \[
  |P_s(\xc_1',\eta_1')-P_s(\xc_2',\eta_2')|
  \le\xi''(1)\bigl(\|\xc_1'-\xc_2'\|_{L^1}+\|\cdf{\eta_1'}-\cdf{\eta_2'}\|_{L^1}\bigr),
 \]
 uniformly in $(s,g_0)$.
\end{thmlist}
Consequently $P_s$ extends uniquely to general data $(\xc,\eta)$ as a continuous function for the
$L^1$ distances, $|P_s(\xc_k,\eta_k)-P_s(\xc,\eta)|\le\xi''(1)\bigl(\|\xc_k-\xc\|_{L^1}
+\|\cdf{\eta_k}-\cdf\eta\|_{L^1}\bigr)\le3\xi''(1)/k$, and
$P_0(\xc,\eta)=P(\glue(\xc,\eta))$ in the sense of Lemma~\ref{lem:pde-measurable}.
\end{lemma}

\begin{proof}
(a) The values are $\xc(jQ/k)+1/k\in(\theta,2+1/k]$, nondecreasing in $j$. On the $j$-th cell,
$|\xc-\xc_k|\le\xc((j+1)Q/k-)-\xc(jQ/k)+1/k$ by monotonicity, and summing over the cells
telescopes. (b) Each point moves by at most $(1-Q)/k$, and on the line the $\Wone$ distance is the
$L^1$ distance of the distribution functions.

(c) For the cascade coefficient: $\Psi_s\in\mathcal T_2\subset\mathcal T_3$ by
Lemma~\ref{lem:pde-fenchel}, so Lemma~\ref{lem:pde-coef} with $G=3$ bounds the change of
$W[\cdot,\Psi_s](0,0)$ by $\frac12\xi''(1)\|\xc_1'-\xc_2'\|_{L^1}$, and the linear term changes by
at most as much. For the band measure: by Lemma~\ref{lem:pde-terminal},
$|W[\xc,\Psi^{(1)}_s]-W[\xc,\Psi^{(2)}_s]|\le\sup|\Psi^{(1)}_s-\Psi^{(2)}_s|$, where $\Psi^{(i)}_s$
is built from $\eta_i'$. The supremum defining $\Psi_s$ changes by at most
$\sup_m|\Phistar_{\eta_1'}(Q,m)-\Phistar_{\eta_2'}(Q,m)|\le\sup_y|\Phi_{\eta_1'}(Q,y)-\Phi_{\eta_2'}(Q,y)|$,
which is at most $\frac12\xi''(1)\|\cdf{\eta_1'}-\cdf{\eta_2'}\|_{L^1}$ by Lemma~\ref{lem:pde-coef}
on $[Q,1]$; the linear term changes by at most as much. The consequences follow from density of
step data, from (a)--(b) with $\xc(Q-)-\xc(0)\le2$, and, for $s=0$, from
Corollary~\ref{cor:pde-step}\ref{cor:pde-step-c} and Lemma~\ref{lem:pde-measurable}, since
$\|\glue(\xc_1',\eta_1')-\glue(\xc_2',\eta_2')\|_{L^1}
=\|\xc_1'-\xc_2'\|_{L^1}+\|\cdf{\eta_1'}-\cdf{\eta_2'}\|_{L^1}$.
\end{proof}

Combined with Proposition~\ref{prop:glued-comparison}, the lemma gives the comparison in the form
in which it is used, \eqref{eq:loc-glued-general}.

\begin{corollary}[Measurable majorant]\label{cor:pde-majorant}
Assume \condH. Let $\rho\in[0,1-Q]$, with $\rho>0$ if $Q=0$, and let the other data be as in
Lemma~\ref{lem:pde-approx}. For every $\varepsilon_A>0$ and $N\ge1$ there is a
random variable $\mathcal T'$ such that
\[
 \FTAP(m)+s\sum_{i=1}^Ng_0(m_i)\le NP_s(\xc,\eta)+N\varepsilon_A+\mathcal T'
 \qquad\text{for every } m\in\cube \text{ with } q_m\in[Q,Q+\rho],
\]
and $\E e^{\theta\mathcal T'}\le(N+1)\,e^{N\theta C_R\rho}$.
The bound on $\E e^{\theta\mathcal T'}$ depends neither on $Q$ nor on the data.
\end{corollary}

\begin{proof}
Choose $k\ge3\xi''(1)/\varepsilon_A$, so that $P_s(\xc_k,\eta_k)\le P_s(\xc,\eta)+\varepsilon_A$ by
Lemma~\ref{lem:pde-approx}. After merging equal neighbouring values, $\xc_k$ is a step function with
strictly increasing values in $(\theta,3]$, so $(\xc_k,\eta_k)$ are step data for
Proposition~\ref{prop:glued-comparison}, with top parameter $y_k=\xc_k(Q-)>\theta$ (and $y=1>\theta$
if $Q=0$). Let $\mathcal T'$ be the random variable that the proposition provides for these data.
Then the pointwise bound holds, and
$\E e^{\theta\mathcal T'}\le(N+1)^{\theta/y_k}e^{N\theta C_R\rho}\le(N+1)e^{N\theta C_R\rho}$.
\end{proof}

\subsection{Comparison with Boursier's cube bound}\label{app:glued-map}

Table~\ref{tab:glued-map} maps each step of the proof above onto \cite[Section~2]{Bou26b}, where
the cascade parameters $0<\theta_1<\dots<\theta_L\le1$ are the values of the distribution function
of a finite-step $\zeta$ with $\zeta(\{0\})>\theta$.

\begin{table}[ht]
\centering
\small
\caption{The proof of Proposition~\ref{prop:glued-comparison} step by step, against the proof of
the cube bound \cite[Proposition~2.7]{Bou26b}. No step uses monotonicity of the coefficient
across the radius $Q$ or a bound by one on the cascade parameters. Rows follow the order of
\cite[Section~2]{Bou26b}.}
\label{tab:glued-map}
\begin{tabular}{@{}p{0.24\textwidth}p{0.25\textwidth}p{0.45\textwidth}@{}}
\toprule
\raggedright Here & \raggedright In \cite{Bou26b} & \raggedright What changes\tabularnewline
\midrule
\raggedright Step 1: cascade on the levels of $\xc$, top level $[q_{L-1},Q+\rho]$
 & \raggedright Definitions 2.2 and 2.3, (2.3); pp.~7--8
 & \raggedright The parameters are the values of $\xc$, any increasing positive numbers. The
 bound $\theta_L\le1$ (p.~7) is a property of her $\zeta$, not of the construction.
 \tabularnewline\addlinespace[3pt]
\raggedright Step 6: Onsager shift
 & \raggedright Marking identity (2.2); Remark 2.2, first identity; pp.~8--9
 & \raggedright Lemma~\ref{lem:onsager-shift}, proved for any positive parameters, used in
 its form~(c) for a countable family of copies.\tabularnewline\addlinespace[3pt]
\raggedright Step 5: exact levels below $Q$
 & \raggedright Remark 2.2, second identity; pp.~8--9
 & \raggedright Read off from the Cole--Hopf recursion, with an arbitrary terminal value
 $\Psi_s$ at $Q$.\tabularnewline\addlinespace[3pt]
\raggedright Step 0: reduction
 & \raggedright (2.4); p.~9
 & \raggedright The band correction is frozen at $Q$; the shell costs
 $\frac N2\xi''(1)\rho$.\tabularnewline\addlinespace[3pt]
\raggedright Step 3: Slepian's lemma
 & \raggedright (2.5), (2.6), Lemma 2.4; pp.~9--10; Proposition~2.7, Step~3; p.~14
 & \raggedright Unchanged: common mean, equal variances, convexity of $\xi$ on $[-1,1]$;
 the points $m$ range over a finite set $D$.\tabularnewline\addlinespace[3pt]
\raggedright Step 4: Fenchel's inequality
 & \raggedright (2.7); p.~10
 & \raggedright $\Psi_s$ in place of $\Phi_\zeta(q_m,\cdot)$.\tabularnewline\addlinespace[3pt]
\raggedright Step 5: the cascade side
 & \raggedright Lemma 2.5; pp.~10--12; Proposition~2.7, Step~5; p.~14
 & \raggedright Doob's inequality only on $[Q,Q+\rho]$. The induction over levels is
 Lemma~\ref{lem:cascade-domination}, which needs only $\xc_1<\dots<\xc_L$.\tabularnewline\addlinespace[3pt]
\raggedright Step 6: lower side
 & \raggedright Lemma 2.6; p.~12; Proposition~2.7, Step~6; p.~14
 & \raggedright The first maximizer over a finite set replaces a measurable maximizer over the
 cube.\tabularnewline\addlinespace[3pt]
\raggedright Steps 2 and 7: outer level
 & \raggedright Proposition~2.7, Steps~1, 2 and~7; pp.~13--14
 & \raggedright Unchanged; needs only $\theta<\xc_1$.\tabularnewline\addlinespace[3pt]
\raggedright Non-step data
 & \raggedright Proposition~2.7, end of the proof; pp.~14--15
 & \raggedright Lemma~\ref{lem:pde-approx} and Corollary~\ref{cor:pde-majorant}: continuity in
 the data in $L^1$, with no monotonicity of the glued coefficient.\tabularnewline
\bottomrule
\end{tabular}
\end{table}

\paragraph{What is used.}
The proof uses that $\xc$ is nondecreasing on $[0,Q)$, strictly after merging equal values, for
Lemma~\ref{lem:cascade-domination}; that $\xc_1>\theta$, for the outer level; that $\eta$ is a
probability measure on $[Q,1]$, for $\FTAP\le F_{\Shift_{q_m}\eta}$ in Step~0; that $|s|\,\|g_0''\|_\infty\le\frac12$, for
Lemma~\ref{lem:glued-terminal}; that $Q+\rho\le1$, so that the top level of the cascade ends at
$r_L\le1$ (Step~1); and that $\xi$ is
convex on $[-1,1]$, in Slepian's lemma. The cap $\xc\le3$ only fixes $C_R$ and the class of
Lemma~\ref{lem:pde-class} containing $W_s$; any finite cap would do.

\paragraph{What is not used.}
Each of Boursier's steps goes through for parameters above one and a downward jump at $Q$, once the
bound $|\partial_x\Phi|\le1$ she uses is available (Lemma~\ref{lem:pde-class}); by
Lemma~\ref{lem:loc-monotone}, for $Q>0$ a nondecreasing glued coefficient gains nothing over
$P_{\zf}$, so the jump is the only source of gain.

\paragraph{Consistency with the cube bound.}
If $\zeta(\{0\})\ge\theta$, the data $(\cdf\zeta|_{[0,Q)},\Shift_Q\zeta)$ have glued coefficient
$\cdf\zeta$, so $P_0=P_\zeta$. Summing Corollary~\ref{cor:pde-majorant} over a finite cover of
$[0,1]$ by shells, and letting the approximation error and $\rho$ tend to zero, gives back
\cite[Proposition~2.7]{Bou26b}.

\subsection{Proof of the field law}\label{app:var-field}

This subsection completes the proof of Theorem~\ref{thm:struct}\ref{thm:struct-ii} outlined in
Section~\ref{sec:loc-field}: test functions detect $\Wone$-deviations (Lemma~\ref{lem:var-net}); by
the tangent bound (Lemma~\ref{lem:var-tangent}, stated for a general base measure for
Section~\ref{sec:eq}), the perturbation $s\sum_ig(m_i)$ raises the bound of
Proposition~\ref{prop:glued-comparison} by $s$ times the field-law mean of $g$, up to $O(s^2)$; and
an exponent budget lets the linear gain beat all errors.

\begin{lemma}[Test functions]\label{lem:var-net}
For every $\kappa>0$ there are $\tilde g_1,\dots,\tilde g_J\in C^\infty(\R)$ with
$\|\tilde g_j'\|_\infty\le1$ and $\|\tilde g_j''\|_\infty<\infty$ such that, for all probability
measures $\mu,\mu'$ on $[-1,1]$,
\[
 \Wone(\mu,\mu')\ge\kappa\quad\Longrightarrow\quad
 \int\tilde g_j\,d(\mu-\mu')\ge\frac\kappa2\ \text{ for some } j\le J .
\]
\end{lemma}

\begin{proof}
The set $\mathcal L$ of $1$-Lipschitz functions on $[-1,1]$ that vanish at $0$ is compact for the
supremum norm, by the Arzel\`a--Ascoli theorem. Let $g_1,\dots,g_J$ be a $\kappa/16$-net of
$\mathcal L$. Extend each $g_j$ to a $1$-Lipschitz function on $\R$ and convolve it with a smooth
probability density supported in $[-\kappa/16,\kappa/16]$; the result $\tilde g_j$ is smooth, with
$\|\tilde g_j'\|_\infty\le1$, bounded second derivative, and $\sup_{[-1,1]}|\tilde g_j-g_j|\le\kappa/16$.
So every $g\in\mathcal L$ is within $\kappa/8$ of some $\tilde g_j$ on $[-1,1]$. By the
Kantorovich--Rubinstein duality, $\Wone(\mu,\mu')=\sup_{g\in\mathcal L}\int g\,d(\mu-\mu')$. If
$\Wone(\mu,\mu')\ge\kappa$, choose $g\in\mathcal L$ with $\int g\,d(\mu-\mu')\ge3\kappa/4$ and $j$
with $\sup_{[-1,1]}|g-\tilde g_j|\le\kappa/8$; then
$\int\tilde g_j\,d(\mu-\mu')\ge3\kappa/4-2\cdot\kappa/8=\kappa/2$.
\end{proof}

The tangent bound is the heart of the argument. Recall from Lemma~\ref{lem:pde-approx} the
functional $P_s(\xc,\eta)$ of~\eqref{eq:glued-Ps}, extended to general data, and the constant
$L_\psi$ of~\eqref{eq:pde-constants} (with $G=3$), which controls the Auffinger--Chen process
under a change of the terminal value.

\begin{lemma}[Tangent bound]\label{lem:var-tangent}
Assume \condH. Let $\theta\in(0,1)$ and $\zeta'\in\Pspace([0,1])$ with $\zeta'(\{0\})\ge\theta$,
let $X'$ be the Auffinger--Chen process of $\zeta'$, and put
$L'_q=\Law(\Phi_{\zeta',x}(q,X'_q))$ for $q\in[0,1]$. Let $g\in C^2([-1,1])$ with
$\|g'\|_\infty\le1$, and let $|s|\le(2\|g''\|_\infty)^{-1}$. Put $C_2=1+L_\psi$, which depends
only on $\xi'(1)$ and $\xi''(1)$.
\begin{thmlist}
\item\label{lem:var-tangent-a} For every $Q\in[0,1)$,
 \begin{equation}\label{eq:var-tangent}
  P_s\bigl(\cdf{\zeta'}|_{[0,Q)},\Shift_Q\zeta'\bigr)\le P_{\zeta'}+s\int g\,dL'_Q+C_2s^2 .
 \end{equation}
\item\label{lem:var-tangent-b} For every $Q\in[0,1)$, $\rho'\in(0,1-Q]$, $\varepsilon_A>0$ and
 $N\ge1$ there is a random variable $\mathcal T$ with
 $\E e^{\theta\mathcal T}\le(N+1)e^{N\theta C_R\rho'}$ such that
 \begin{multline*}
  \sup\Bigl\{\FTAP(m)+s\sum_{i=1}^Ng(m_i):\ m\in\cube,\ q_m\in[Q,Q+\rho']\Bigr\}\\
  \le N\Bigl[P_{\zeta'}+s\int g\,dL'_Q+C_2s^2+\varepsilon_A\Bigr]+\mathcal T .
 \end{multline*}
\end{thmlist}
\end{lemma}

The proof uses only $\zeta'(\{0\})\ge\theta$, which makes the base data admissible, and not any
minimality of $\zeta'$. In words: the perturbation shifts the terminal value at the radius by
$s\,g$ evaluated at the optimal magnetization $\Phi_{\zeta',x}(Q,\cdot)$ of the TAP correction, to
first order; the comparison carries this shift down to the origin along the Auffinger--Chen
process, where it averages to $s\int g\,dL'_Q$.

\begin{proof}
\ref{lem:var-tangent-a} \emph{Step data.} Let $\xc$ be a step function on $[0,Q)$ with values in
$[0,3]$ (none if $Q=0$), and let $\eta\in\Pspace([Q,1])$ have a step distribution function. Let $c$ be their
glued coefficient~\eqref{eq:loc-glue}, $v_c=\partial_yW[c,\log\cosh]$, and $X^c$ the solution
of~\eqref{eq:pde-ac} for $(c,\log\cosh)$ from $(0,0)$. We claim that
\begin{equation}\label{eq:var-tangent-step}
 P_s(\xc,\eta)\le P_0(\xc,\eta)+s\,\E g\bigl(v_c(Q,X^c_Q)\bigr)+C_2s^2 .
\end{equation}
Let $\Psi_s$ and its maximizer $m^*_s$ be as in Lemma~\ref{lem:pde-fenchel}, and
$h=-\Phistar_\eta(Q,\cdot)$, so that $h''\ge1$ and $\Psi_0'=m^*_0=\Phi_\eta'(Q,\cdot)$.

First, the perturbation of the terminal value. By the definition of $\Psi_s$ as a supremum,
$\Psi_s(z)\ge zm^*_0(z)-h(m^*_0(z))+s\,g(m^*_0(z))=\Psi_0(z)+s\,g(m^*_0(z))$. In the other
direction, the optimality conditions $h'(m^*_s(z))-s\,g'(m^*_s(z))=z=h'(m^*_0(z))$ and $h''\ge1$
give $|m^*_s(z)-m^*_0(z)|\le|s|\,\|g'\|_\infty\le|s|$, and therefore
\[
 \Psi_s(z)=zm^*_s-h(m^*_s)+s\,g(m^*_s)\le\Psi_0(z)+s\,g(m^*_s(z))
 \le\Psi_0(z)+s\,g(m^*_0(z))+s^2 .
\]
So $0\le\Psi_s-\Psi_0-s\,g\circ m^*_0\le s^2$, and $\|\Psi_s'-\Psi_0'\|_\infty=\sup|m^*_s-m^*_0|\le|s|$.

Next, the transport to the origin. If $Q=0$, then $P_s-P_0=\Psi_s(0)-\Psi_0(0)$, $m^*_0(0)=0$ and
$v_c(0,X^c_0)=0$, so~\eqref{eq:var-tangent-step} follows from the last display. Let $Q>0$. The
integral terms of $P_s$ and $P_0$ coincide, so $P_s-P_0=W[\xc,\Psi_s](0,0)-W[\xc,\Psi_0](0,0)$. By
Lemma~\ref{lem:pde-fenchel}, $\Psi_0,\Psi_s\in\mathcal T_2\subset\mathcal T_3$. Let $X^0$ and $X^s$
solve~\eqref{eq:pde-ac} on $[0,Q]$ for $(\xc,\Psi_0)$ and $(\xc,\Psi_s)$ from $(0,0)$, driven by
the same Brownian motion. The upper half of the terminal sandwich,
Lemma~\ref{lem:pde-terminal}, gives
\[
 P_s-P_0\le\E\bigl[\Psi_s-\Psi_0\bigr](X^s_Q)\le s\,\E g\bigl(m^*_0(X^s_Q)\bigr)+s^2 .
\]
Lemma~\ref{lem:pde-stability} with equal coefficients gives
$\sup_r|X^s_r-X^0_r|\le L_\psi\|\Psi_s'-\Psi_0'\|_\infty\le L_\psi|s|$. The function $g\circ m^*_0$
is $1$-Lipschitz, because $g$ is and $(m^*_0)'=\Phi_\eta''(Q,\cdot)\le1$
(Lemma~\ref{lem:pde-fenchel}\ref{lem:pde-fenchel-a}); hence
$s\,\E g(m^*_0(X^s_Q))\le s\,\E g(m^*_0(X^0_Q))+L_\psi s^2$. Finally, the recursion
\eqref{eq:pde-colehopf} on $[0,Q]$ with terminal value $\Psi_0=\Phi_\eta(Q,\cdot)$ continues the one
on $[Q,1]$, so $W[\xc,\Psi_0]$ and $\Phi_\eta$ together form $W[c,\log\cosh]$; thus $X^0=X^c$ on
$[0,Q]$ and $m^*_0(X^0_Q)=v_c(Q,X^c_Q)$. This proves~\eqref{eq:var-tangent-step}.

\emph{General data.} Let $\xc=\cdf{\zeta'}|_{[0,Q)}$ and $\eta=\Shift_Q\zeta'$. These data are
admissible in Lemma~\ref{lem:pde-approx}, since $\xc$ is nondecreasing with values in
$[\theta,1]$, and their glued coefficient is $\cdf{\zeta'}$. Let $(\xc_k,\eta_k)$ be the step data
of Lemma~\ref{lem:pde-approx}, with glued coefficients $c_k$. Then
$\|c_k-\cdf{\zeta'}\|_{L^1}\le3/k$, and as $k\to\infty$:
$P_s(\xc_k,\eta_k)\to P_s(\xc,\eta)$ by Lemma~\ref{lem:pde-approx};
$P_0(\xc_k,\eta_k)=P(c_k)\to P(\cdf{\zeta'})=P_{\zeta'}$ by Corollary~\ref{cor:pde-step}\ref{cor:pde-step-c}
and Lemma~\ref{lem:pde-measurable}\ref{lem:pde-measurable-iv}; and, since $g$ is $1$-Lipschitz,
\[
 \bigl|\E g\bigl(v_{c_k}(Q,X^{c_k}_Q)\bigr)-\textstyle\int g\,dL'_Q\bigr|
 \le\E\bigl|v_{c_k}(Q,X^{c_k}_Q)-v_{\cdf{\zeta'}}(Q,X^{\cdf{\zeta'}}_Q)\bigr|
 \le L_M\|c_k-\cdf{\zeta'}\|_{L^1}\to0
\]
by Corollary~\ref{cor:pde-step}\ref{cor:pde-step-d}, extended to measurable coefficients by
Lemma~\ref{lem:pde-measurable}, which also identifies $v_{\cdf{\zeta'}}(Q,X^{\cdf{\zeta'}}_Q)$ with
$\Phi_{\zeta',x}(Q,X'_Q)$. Applying~\eqref{eq:var-tangent-step} to $(\xc_k,\eta_k)$ and letting
$k\to\infty$ gives~\eqref{eq:var-tangent}.

\ref{lem:var-tangent-b} The data $(\cdf{\zeta'}|_{[0,Q)},\Shift_Q\zeta')$ and the perturbation
$(s,g)$ are admissible in Corollary~\ref{cor:pde-majorant}, which is~\eqref{eq:loc-glued-general};
it gives $\mathcal T$ with the stated moment bound and
$\sup[\FTAP(m)+s\sum_ig(m_i)]\le NP_s(\cdf{\zeta'}|_{[0,Q)},\Shift_Q\zeta')+N\varepsilon_A+\mathcal T$
over the shell. Insert~\eqref{eq:var-tangent}.
\end{proof}

\begin{proof}[Proof of Theorem~\ref{thm:struct}\ref{thm:struct-ii}]
Assume \condH, fix $f\in(\feq,\fone)$ and $\kappa>0$, and write $\theta=\thf$, $\zeta=\zf$ (any
constrained minimizer) and $X$ for the Auffinger--Chen process of $\zeta$, so that
$\Lq{q}=\Law(\Phi_{\zeta,x}(q,X_q))$. By Lemma~\ref{lem:loc-obstacle}, $\zeta(\{0\})=\theta$ and
$\Sigma(f)=\theta(P_\zeta-f)$.

\emph{Test functions and continuity in the radius.} Let $\tilde g_1,\dots,\tilde g_J$ be given by
Lemma~\ref{lem:var-net} for $\kappa$; they are the functions $g_j$ of
Section~\ref{sec:loc-field}. The map $q\mapsto\Phi_{\zeta,x}(q,X_q)$ is almost surely continuous on
$[0,1]$ and bounded by one (Lemma~\ref{lem:pde-measurable}), so by dominated convergence each
$q\mapsto\int\tilde g_j\,d\Lq{q}$ is continuous, hence uniformly continuous, on $[0,1]$. Let
$\omega$ be a common modulus of continuity for $j\le J$.

\emph{Parameters.} Choose, in this order,
\[
 s=\min\Bigl\{\frac{\kappa}{8C_2},\ \bigl(2\max_j\|\tilde g_j''\|_\infty\bigr)^{-1}\Bigr\},
 \qquad \rho\in(0,1]\ \text{with}\ C_R\rho\le\frac{s\kappa}{16}\ \text{and}\
 \omega(\rho)\le\frac\kappa{16},
\]
and put $\delta=s\kappa/16$ and $\varepsilon_A=s\kappa/32$. Then $C_2s^2\le s\kappa/8$, and
$(s,\tilde g_j)$ is admissible in Lemma~\ref{lem:var-tangent}.

\emph{One shell and one test function.} Let $Q_i=i\rho$ and $I_i=[Q_i,\min(Q_i+\rho,1)]$ for the
integers $i\ge0$ with $Q_i<1$. Fix $i$ and $j$. Suppose that $m\in\cube$ has $q_m\in I_i$,
$\FTAP(m)\ge N(f-\delta)$ and $N^{-1}\sum_k\tilde g_j(m_k)\ge\int\tilde g_j\,d\Lq{q_m}+\kappa/2$.
Since $|\int\tilde g_j\,d\Lq{q_m}-\int\tilde g_j\,d\Lq{Q_i}|\le\omega(\rho)$,
\[
 \FTAP(m)+s\sum_{k=1}^N\tilde g_j(m_k)\ge
 N\Bigl[f-\delta+s\Bigl(\int\tilde g_j\,d\Lq{Q_i}-\omega(\rho)+\frac\kappa2\Bigr)\Bigr].
\]
Lemma~\ref{lem:var-tangent}\ref{lem:var-tangent-b}, with $\zeta'=\zeta$, $Q=Q_i$,
$\rho'=\min(\rho,1-Q_i)$ and $g=\tilde g_j$, provides a random variable $\mathcal T_{ij}$ with
$\E e^{\theta\mathcal T_{ij}}\le(N+1)e^{N\theta C_R\rho}$ that bounds the left side by
$N[P_\zeta+s\int\tilde g_j\,d\Lq{Q_i}+C_2s^2+\varepsilon_A]+\mathcal T_{ij}$. The terms
$s\int\tilde g_j\,d\Lq{Q_i}$ cancel, so such an $m$ can exist only on the event
\[
 \Bigl\{\mathcal T_{ij}\ge N\Bigl[f-P_\zeta-\delta-s\,\omega(\rho)+\frac{s\kappa}2-C_2s^2
 -\varepsilon_A\Bigr]\Bigr\}.
\]
By Markov's inequality applied to $e^{\theta\mathcal T_{ij}}$, and $\Sigma(f)=\theta(P_\zeta-f)$,
this event has probability at most
\[
 (N+1)\exp N\Bigl\{\Sigma(f)+\theta\Bigl[C_2s^2+\varepsilon_A+C_R\rho+\delta+s\,\omega(\rho)
 -\frac{s\kappa}2\Bigr]\Bigr\}.
\]

\emph{The exponent budget.} By the choice of parameters, the bracket is at most
\begin{equation}\label{eq:var-budget}
 s\kappa\Bigl(\frac4{32}+\frac1{32}+\frac2{32}+\frac2{32}+\frac2{32}-\frac{16}{32}\Bigr)
 =-\frac5{32}\,s\kappa .
\end{equation}
The gain $s\kappa/2$ is linear in $s$, the cost $C_2s^2$ of the perturbation is quadratic, and the
remaining errors were chosen small compared with $s\kappa$.

\emph{Union bound.} If $\Wone(\mu_m,\Lq{q_m})\ge\kappa$, then $q_m\in I_i$ for some $i$, and
Lemma~\ref{lem:var-net} gives a $j$ with $N^{-1}\sum_k\tilde g_j(m_k)\ge\int\tilde g_j\,d\Lq{q_m}
+\kappa/2$. So the event of~\ref{thm:struct-ii} is contained in the union of the events above over
the at most $(\lceil1/\rho\rceil+1)J$ pairs $(i,j)$, a finite union of measurable events. Its
outer probability is at most
$(\lceil1/\rho\rceil+1)J(N+1)e^{N(\Sigma(f)-5\theta s\kappa/32)}\le e^{N(\Sigma(f)-\theta s\kappa/8)}$
for all large $N$. This proves~\ref{thm:struct-ii} with $\delta=s\kappa/16$ and
$c=\thf s\kappa/8$.
\end{proof}

The constants $\delta=s\kappa/16$ and $c=\thf s\kappa/8$ depend on $\kappa$ through $s$, hence
through $C_2$ and the second derivatives of the test functions; the modulus $\omega$ and the number
$J$ of test functions enter only through $\rho$ and the threshold beyond which the bound holds.
Taking the smaller of these values of $\delta$ and of those in the proof of
Theorem~\ref{thm:struct}\ref{thm:struct-i}, and half the smaller of the two values of $c$, gives the
joint form stated after the outline in Section~\ref{sec:loc-field}.

%% file: sections/appC-gaussian.tex
\section{Gaussian bounds and the TAP maximum}\label{app:gauss}

This appendix supplies the inputs of Section~\ref{sec:max} that are not part of the Ramsey
argument. Appendix~\ref{app:gauss-bounds} bounds the field and the moments of $\ZN$ and $\MN$ and
compares $\xi$ with its truncations in the degree; Appendix~\ref{app:gauss-truncation} deduces the
band theorem at every exponential rate under \condR, which is where Theorem~\ref{thm:max} leaves
$\classC$; Appendix~\ref{app:gauss-critical} proves Theorem~\ref{thm:exist}\ref{thm:exist-a}.
Appendices~\ref{app:var-moments}--\ref{app:var-overlap} give the upper tail of $\MN$ and a lower
bound on the self-overlap of near-maximizers; Appendix~\ref{app:gauss-rpde} explains where
$\classC$ enters. Here $K_\xi=\xi'(1)+\xi''(1)$.

\subsection{Gaussian bounds}\label{app:gauss-bounds}

\paragraph{The field.} We realize the field as the series
\[
 \HN=\sum_{p\ge2}\beta_p\,H_{N,p},\qquad
 H_{N,p}(x)=N^{-(p-1)/2}\sum_{i_1,\dots,i_p=1}^Ng^{(p)}_{i_1\dots i_p}\,x_{i_1}\cdots x_{i_p},
\]
with independent standard Gaussian coefficients $g^{(p)}_{i_1\dots i_p}$. Then
$\E H_{N,p}(x)H_{N,p}(y)=N\ov{x}{y}^p$, so that $\E\HN(x)\HN(y)=N\xi(\ov{x}{y})$ wherever the
series converges, and $\HN(0)=0$ because $p\ge2$. Throughout this appendix
\[
 K_\xi=\xi'(1)+\xi''(1)=\sum_{p\ge2}p^2\beta_p^2 ,
\]
which is finite under \condR. For $D\ge2$ we split $\xi=\xi_{\le D}+\xi_{>D}$ and
$\HN=\HN^{\le D}+\HN^{>D}$ according to $p\le D$ and $p>D$; the two parts of the field are
independent, with covariances $N\xi_{\le D}(\ov{x}{y})$ and $N\xi_{>D}(\ov{x}{y})$.

\paragraph{Gaussian tools.} We use two classical facts about centered Gaussian processes
$(X_s)_{s\in S}$ indexed by a countable set $S$; they apply to processes with continuous paths on
subsets of $\R^N$ through countable dense subsets. The \emph{Sudakov--Fernique inequality}: if
$(Y_s)_{s\in S}$ is another such process with $\E\sup_SY<\infty$ and
$\E(X_s-X_{s'})^2\le\E(Y_s-Y_{s'})^2$ for all $s,s'$, then $\E\sup_SX\le\E\sup_SY$. \emph{Gaussian
concentration of the supremum} (see \cite{Ledoux}; the tail bound, the Borell--TIS inequality,
is \cite[(1.24)]{Ledoux}): if $\sigma^2=\sup_s\E X_s^2<\infty$ and $\sup_SX<\infty$ almost surely,
then $M=\sup_SX$ is integrable and
\begin{equation}\label{eq:gauss-conc}
 \E e^{\lambda(M-\E M)}\le e^{\lambda^2\sigma^2/2}\quad(\lambda\in\R),\qquad
 \Prob\bigl(M-\E M\ge u\bigr)\le e^{-u^2/(2\sigma^2)}\quad(u>0).
\end{equation}

\needspace{5\baselineskip}
\begin{lemma}[The field]\label{lem:gauss-field}
Assume \condR\ and fix $N$.
\begin{thmlist}
\item\label{lem:gauss-field-a} There is $r>1$, depending only on $\xi$, such that almost surely the
 series defining $\HN$ and the series of its first partial derivatives converge uniformly on
 $[-r,r]^N$. In particular $\HN$ is almost surely continuously differentiable on a neighbourhood
 of $\cube$, and $\sup_{\cube}|\partial_i\HN|<\infty$ for every $i$.
\item\label{lem:gauss-field-b} For $x,y\in\cube$, $\E(\HN(x)-\HN(y))^2\le K_\xi|x-y|^2$, and
 \[
  \E\sup_{\cube}\HN\le\sqrt{2K_\xi/\pi}\;N .
 \]
\end{thmlist}
\end{lemma}

\begin{proof}
\emph{One degree.} Let $\rho\ge1$ and $x,y\in[-\rho,\rho]^N$, so that $|x|,|y|\le\rho\sqrt N$.
Writing $x^{\otimes p}-y^{\otimes p}=\sum_{k=0}^{p-1}x^{\otimes k}\otimes(x-y)\otimes
y^{\otimes(p-1-k)}$ gives $|x^{\otimes p}-y^{\otimes p}|\le p(\rho\sqrt N)^{p-1}|x-y|$, hence
\begin{equation}\label{eq:gauss-increment}
 \E\bigl(H_{N,p}(x)-H_{N,p}(y)\bigr)^2=N^{1-p}\,|x^{\otimes p}-y^{\otimes p}|^2
 \le p^2\rho^{2(p-1)}|x-y|^2 .
\end{equation}
Compare $H_{N,p}$ on $[-\rho,\rho]^N$ with the process $p\rho^{p-1}g\cdot x$, where $g$ is a
standard Gaussian vector of $\R^N$, and note that
\[
 \E\sup_{x\in[-\rho,\rho]^N}g\cdot x=\rho\,\E\sum_{i=1}^N|g_i|=\rho N\sqrt{2/\pi}.
\]
The Sudakov--Fernique inequality gives $\E\sup_{[-\rho,\rho]^N}H_{N,p}\le\sqrt{2/\pi}\,p\rho^pN$,
and the same bound for $-H_{N,p}$. Both suprema are nonnegative because $H_{N,p}(0)=0$, so
\begin{equation}\label{eq:gauss-sup-p}
 \E\sup_{[-\rho,\rho]^N}|H_{N,p}|\le2\sqrt{2/\pi}\;p\rho^pN .
\end{equation}

\ref{lem:gauss-field-a} Let $R_\xi>1$ be the radius of convergence of $\xi$ and fix
$r\in(1,R_\xi^{1/2})$. For fixed $i$ and fixed coordinates $x_j$, $j\ne i$, in $[-r,r]$, the map
$s\mapsto H_{N,p}(x_1,\dots,x_{i-1},s,x_{i+1},\dots,x_N)$ is a polynomial of degree at most $p$.
A.~A.~Markov's inequality for the derivative of a polynomial, $\max_{[-r,r]}|P'|\le
p^2r^{-1}\max_{[-r,r]}|P|$ for $\deg P\le p$, therefore gives
$\sup_{[-r,r]^N}|\partial_iH_{N,p}|\le p^2r^{-1}\sup_{[-r,r]^N}|H_{N,p}|$. With
\eqref{eq:gauss-sup-p},
\[
 \E\sum_{p\ge2}|\beta_p|\sup_{[-r,r]^N}\Bigl(|H_{N,p}|+\sum_{i=1}^N|\partial_iH_{N,p}|\Bigr)
 \le2\sqrt{2/\pi}\,N(1+N)\sum_{p\ge2}p^3|\beta_p|r^p<\infty,
\]
because $\limsup_p(|\beta_p|r^p)^{1/p}=rR_\xi^{-1/2}<1$. So almost surely both series converge
absolutely and uniformly on $[-r,r]^N$, and $\HN$ is continuously differentiable on $(-r,r)^N$.

\ref{lem:gauss-field-b} By independence of the $H_{N,p}$ and \eqref{eq:gauss-increment} with
$\rho=1$, $\E(\HN(x)-\HN(y))^2\le\sum_p\beta_p^2p^2|x-y|^2=K_\xi|x-y|^2$. By~\ref{lem:gauss-field-a}
the field is continuous on $\cube$, and the Sudakov--Fernique inequality against $\sqrt{K_\xi}\,g\cdot x$
gives $\E\sup_{\cube}\HN\le\sqrt{K_\xi}\,\E\sup_{\cube}g\cdot x=\sqrt{2K_\xi/\pi}\,N$.
\end{proof}

The next lemma collects the a priori bounds of Section~\ref{sec:max}, including the comparison with
truncated covariances used in Appendix~\ref{app:gauss-truncation}. For $D$ at least the smallest
active degree, $\xi_{\le D}$ is a nonzero polynomial with nonnegative coefficients, so
Appendix~\ref{app:pde} applies to it. Quantities built from $(\HN^{\le D},\xi_{\le D})$ in place of
$(\HN,\xi)$ carry a superscript $(D)$: $\ZN^{(D)}$, $\Phi^{(D)}_\zeta$, $P^{(D)}_\zeta$,
$\feq^{(D)}$, $\Lambda^{(D)}$, $\FTAP^{(D)}$, $\MN^{(D)}$, and $\TAPcps^{(D)}$ for the functional
$\TAPcps$ of Section~\ref{sec:intro-setting}. Put
\[
 T_D=\sup_{\cube}|\HN^{>D}|,\qquad K_D=\xi_{>D}'(1)+\xi_{>D}''(1),\qquad
 \varepsilon_D=2\sqrt{2K_D/\pi} .
\]
Under \condR, $\xi_{>D}(1)$, $\xi'_{>D}(1)$, $\xi''_{>D}(1)$ and $\varepsilon_D$ tend to zero as
$D\to\infty$.

\begin{lemma}[A priori bounds]\label{lem:gauss-bounds}
Assume \condR. Put $C^{\max}_\xi=\sqrt{2K_\xi/\pi}+\xi'(1)/2$ and, for $\theta>0$,
$C_\theta=\theta\sqrt{2K_\xi/\pi}+\theta^2\xi(1)$.
\begin{thmlist}
\item\label{lem:gauss-bounds-a} \textup{(Partition function.)} For $\theta>0$, $\E\ZN^\theta\ge1$
 and $\E\ZN^{2\theta}\le e^{2C_\theta N}$; if moreover $\theta\le1$, then
 $\E\ZN^{2\theta}\le e^{2\theta N\xi(1)}$. For $\theta<0$,
 $\E\ZN^\theta\ge e^{-|\theta|\sqrt{2K_\xi/\pi}\,N}$ and $\E\ZN^{2\theta}\le e^{2\theta^2N\xi(1)}$.
\item\label{lem:gauss-bounds-b} \textup{(Deterministic bounds.)} $\FTAP(0)=N\feq$, and
 $\FTAP(x)\le\HN(x)+N\xi'(1)/2$ for every $x\in\cube$. Hence
 $N\feq\le\MN\le\sup_{\cube}\HN+N\xi'(1)/2$.
\item\label{lem:gauss-bounds-c} \textup{(The maximum.)} For $\lambda>0$,
 $\E e^{\lambda\MN}\le\exp\{N[C^{\max}_\xi\lambda+\lambda^2\xi(1)/2]\}$; for $\lambda<0$,
 $\E e^{\lambda\MN}\le e^{\lambda N\feq}$.
\item\label{lem:gauss-bounds-d} \textup{(Truncation.)} Let $D$ be at least the smallest active
 degree.
 \begin{enumerate}[label=\textup{(d\arabic*)}, ref=(d\arabic*), leftmargin=*, topsep=2pt,
  itemsep=1pt, beginpenalty=10000]
 \item\label{lem:gauss-bounds-d1} For every $\zeta\in\Pspace([0,1])$,
  $0\le\Phi_\zeta-\Phi^{(D)}_\zeta\le\xi'_{>D}(1)$ on $[0,1]\times\R$ and
  $|P_\zeta-P^{(D)}_\zeta|\le\frac32\xi'_{>D}(1)$. Consequently
  $|\feq-\feq^{(D)}|\le\frac32\xi'_{>D}(1)$,
  $\sup_{\theta\in[0,1)}|\Lambda(\theta)-\Lambda^{(D)}(\theta)|\le\frac32\xi'_{>D}(1)$,
  $\sup_{\mu\in\Pspace([-1,1])}|\TAPcps(\mu)-\TAPcps^{(D)}(\mu)|\le\frac32\xi'_{>D}(1)$, and
  \[
   |\log\ZN-\log\ZN^{(D)}|\le T_D,\qquad
   \sup_{\cube}|\FTAP-\FTAP^{(D)}|\le T_D+\tfrac32N\xi'_{>D}(1),
  \]
  and the last bound also holds for $|\MN-\MN^{(D)}|$.
 \item\label{lem:gauss-bounds-d2} $\E T_D\le\varepsilon_DN$, and
  $\Prob\bigl(T_D>(\varepsilon_D+u)N\bigr)\le\exp\bigl(-u^2N/(2\xi_{>D}(1))\bigr)$ for $u>0$ if
  $\xi_{>D}\not\equiv0$ (otherwise $T_D=0$).
 \item\label{lem:gauss-bounds-d3} For $\theta>0$,
  $\bigl|N^{-1}\log\E\ZN^\theta-N^{-1}\log\E(\ZN^{(D)})^\theta\bigr|
  \le\theta\varepsilon_D+\theta^2\xi_{>D}(1)/2$.
 \end{enumerate}
\item\label{lem:gauss-bounds-e} \textup{(Tilts.)} For $\theta>0$ and every event $A$,
 $\Qth(A)\le e^{C_\theta N}\Prob(A)^{1/2}$; for $\theta\in(0,1]$ also
 $\Qth(A)\le e^{\theta N\xi(1)}\Prob(A)^{1/2}$.
\end{thmlist}
In particular, \eqref{eq:max-crude} holds for every $\theta\ne0$ and every $N$ with
$C(\xi,\theta)=2|\theta|(C^{\max}_\xi+|\feq|)+2\theta^2\xi(1)$.
\end{lemma}

Here $\MN$ is a random variable by Lemma~\ref{lem:gauss-interior}\ref{lem:gauss-interior-a} below,
whose proof uses only Lemma~\ref{lem:gauss-field} and part~\ref{lem:gauss-bounds-d1}.

\begin{proof}
Write $S=\sup_{\cube}\HN$. By Lemma~\ref{lem:gauss-field}, $S$ is finite and
$\E S\le\sqrt{2K_\xi/\pi}\,N$. Since $\Var\HN(x)=N\xi(q_x)\le N\xi(1)$, the bounds
\eqref{eq:gauss-conc} apply to $S$ with $\sigma^2=N\xi(1)$.

\ref{lem:gauss-bounds-a} By Jensen's inequality for the average over $\sigma$,
$\log\ZN\ge L:=2^{-N}\sum_{\sigma\in\hcube}\HN(\sigma)$, a centered Gaussian variable with
$\Var L=2^{-2N}\sum_{\sigma,\sigma'}N\xi(\ov{\sigma}{\sigma'})\le N\xi(1)$, since
$|\xi(t)|\le\xi(1)$ on $[-1,1]$. Also $\log\ZN\le\max_{\hcube}\HN\le S$, so $\log\ZN$ is
integrable. For $\theta>0$, Jensen's inequality gives
$\E\ZN^\theta\ge e^{\theta\E\log\ZN}\ge e^{\theta\E L}=1$, and \eqref{eq:gauss-conc} with
$\lambda=2\theta$ gives $\E\ZN^{2\theta}\le\E e^{2\theta S}\le e^{2\theta\E S+2\theta^2N\xi(1)}\le
e^{2C_\theta N}$. If $\theta\le1$, then $\E\ZN^{2\theta}\le(\E\ZN^2)^\theta$ and
$\E\ZN^2=2^{-2N}\sum_{\sigma,\sigma'}e^{N\xi(1)+N\xi(\ov{\sigma}{\sigma'})}\le e^{2N\xi(1)}$. For
$\theta<0$, $\E\ZN^\theta\ge e^{\theta\E\log\ZN}\ge e^{\theta\E S}\ge e^{-|\theta|\sqrt{2K_\xi/\pi}N}$,
and $\ZN^{2\theta}\le e^{2\theta L}$ gives $\E\ZN^{2\theta}\le e^{2\theta^2\Var L}\le
e^{2\theta^2N\xi(1)}$.

\ref{lem:gauss-bounds-b} Since $\HN(0)=0$ and $q_0=0$, $F_\zeta(0)=N\Phistar_\zeta(0,0)-\frac
N2\int_0^1t\xi''(t)\cdf\zeta(t)\,dt$. The function $\Phi_\zeta(0,\cdot)$ is even and convex
(Lemma~\ref{lem:pde-measurable}\ref{lem:pde-measurable-i}), so
$\Phistar_\zeta(0,0)=\inf_b\Phi_\zeta(0,b)=\Phi_\zeta(0,0)$. Hence $F_\zeta(0)=NP_\zeta$ for every
$\zeta$, and $\FTAP(0)=N\feq$ by \eqref{eq:intro-ftap}. For the upper bound, $\delta_1$ is
admissible in \eqref{eq:intro-ftap} at every $x$. Its distribution function vanishes on $[0,1)$, so
the integral term of $F_{\delta_1}$ vanishes and $\Phi_{\delta_1}$ solves the heat equation:
$\Phi_{\delta_1}(q,b)=\E\log\cosh\bigl(b+(\xi'(1)-\xi'(q))^{1/2}Z\bigr)$ with $Z$ standard
Gaussian, by \eqref{eq:pde-colehopf}. Taking $b=0$ in the infimum and using
$\log\cosh z\le z^2/2$ gives $\Phistar_{\delta_1}(q,m)\le\Phi_{\delta_1}(q,0)\le\xi'(1)/2$ for
$|m|<1$, and for $m=\pm1$ by continuity. So $\FTAP(x)\le F_{\delta_1}(x)\le\HN(x)+N\xi'(1)/2$.

\ref{lem:gauss-bounds-c} For $\lambda>0$, part~\ref{lem:gauss-bounds-b} and \eqref{eq:gauss-conc}
give $\E e^{\lambda\MN}\le e^{\lambda N\xi'(1)/2}\E e^{\lambda S}\le
\exp\{\lambda N\xi'(1)/2+\lambda\sqrt{2K_\xi/\pi}N+\lambda^2N\xi(1)/2\}$. For $\lambda<0$ use
$\MN\ge N\feq$.

\ref{lem:gauss-bounds-d1} First let $\zeta$ have finitely many atoms, and put $u=\Phi_\zeta$,
$w=\Phi^{(D)}_\zeta$ and $\delta=u-w$. By Lemma~\ref{lem:pde-class}, applied with $\xi$ and with
$\xi_{\le D}$ and $G=1$, both functions are of class $C^{1,2}$ on each closed step of $\cdf\zeta$,
with $|u_x|,|w_x|<1$ and $0\le u_{xx},w_{xx}\le1$. Subtracting the two Parisi equations and using
$u_x^2-w_x^2=(u_x+w_x)\delta_x$,
\[
 \partial_t\delta+\frac{\xi_{\le D}''}2\Bigl(\delta_{xx}+\cdf\zeta\,(u_x+w_x)\,\delta_x\Bigr)
 =-\frac{\xi_{>D}''}2\bigl(u_{xx}+\cdf\zeta\,u_x^2\bigr).
\]
Fix $(q,b)$ and let $X$ solve $dX_r=\frac12\xi_{\le D}''(r)\cdf\zeta(r)(u_x+w_x)(r,X_r)\,dr+
\sqrt{\xi_{\le D}''(r)}\,dB_r$ on $[q,1]$ with $X_q=b$; the drift is bounded and Lipschitz in the
space variable. It\^o's formula on each closed step, chained over the steps as in the proof of
Lemma~\ref{lem:pde-ac}, and $\delta(1,\cdot)=0$ give
\begin{equation}\label{eq:gauss-delta}
 \delta(q,b)=\frac12\,\E\int_q^1\xi_{>D}''(r)\bigl(u_{xx}+\cdf\zeta\,u_x^2\bigr)(r,X_r)\,dr ;
\end{equation}
the stochastic integral is a martingale because $|\delta_x|\le2$. Since
$0\le u_{xx}+\cdf\zeta u_x^2\le2$, we get $0\le\delta(q,b)\le\int_q^1\xi''_{>D}\le\xi'_{>D}(1)$.
For general $\zeta$, approximate $\cdf\zeta$ in $L^1$ by distribution functions of measures with
finitely many atoms; $\Phi_\zeta$ and $\Phi^{(D)}_\zeta$ are the uniform limits
(Lemma~\ref{lem:pde-measurable}\ref{lem:pde-measurable-i}), and the bound passes to the limit.

The conjugates then satisfy $|\Phistar_\zeta(q,m)-\Phi^{*,(D)}_\zeta(q,m)|\le\xi'_{>D}(1)$ for
$|m|<1$, and for $m=\pm1$ by continuity. The integral terms differ by
$\frac12\int_q^1t\xi''_{>D}(t)\cdf\zeta(t)\,dt\le\frac12\vth_{>D}(1)\le\frac12\xi'_{>D}(1)$, where
$\vth_{>D}(1)=\xi'_{>D}(1)-\xi_{>D}(1)$. Hence, for every $\zeta$, $|P_\zeta-P^{(D)}_\zeta|$ and,
for every $x\in\cube$, $N^{-1}\bigl|(F_\zeta(x)-\HN(x))-(F^{(D)}_\zeta(x)-\HN^{\le D}(x))\bigr|$
are at most $\frac32\xi'_{>D}(1)$. Taking infima over $\zeta$, over all of $\Pspace([0,1])$ or over
$\{\zeta(\{0\})\ge\theta\}$, gives the bounds on $\feq$, on $\Lambda=\theta p(\theta)$ and, since
$\FTAP(x)-\HN(x)=N(\TAPcps(\mu_x)-\log2)$, on $\TAPcps$ at every empirical law. The last
bound extends to every $\mu\in\Pspace([-1,1])$, because the same computation applies to the
infimum over $\zeta$ in \cite[(2.24)--(2.25)]{CPS} at any $\mu$. Finally
$|\HN-\HN^{\le D}|\le T_D$ on $\cube$ gives $e^{-T_D}\ZN^{(D)}\le\ZN\le e^{T_D}\ZN^{(D)}$ and the
bounds on $\FTAP$ and $\MN$.

\ref{lem:gauss-bounds-d2} If $\xi_{>D}\not\equiv0$, Lemma~\ref{lem:gauss-field} applies to the
covariance $\xi_{>D}$, whose radius of convergence is at least that of $\xi$, and to the fields
$\pm\HN^{>D}$. Both suprema are nonnegative since $\HN^{>D}(0)=0$, so
$\E T_D\le\E\sup\HN^{>D}+\E\sup(-\HN^{>D})\le\varepsilon_DN$. The variable $T_D$ is the supremum
of the process $(x,s)\mapsto s\HN^{>D}(x)$ on $\cube\times\{-1,1\}$, whose variance is at most
$N\xi_{>D}(1)$, and \eqref{eq:gauss-conc} gives the tail bound.

\ref{lem:gauss-bounds-d3} By~\ref{lem:gauss-bounds-d1}, $\ZN^\theta\le(\ZN^{(D)})^\theta e^{\theta
T_D}$, and $T_D$ is independent of $\ZN^{(D)}$. So \eqref{eq:gauss-conc}
and~\ref{lem:gauss-bounds-d2} give $\E\ZN^\theta\le\E(\ZN^{(D)})^\theta\,e^{\theta\varepsilon_DN
+\theta^2N\xi_{>D}(1)/2}$. Similarly $\E\ZN^\theta\ge\E(\ZN^{(D)})^\theta\,\E e^{-\theta T_D}
\ge\E(\ZN^{(D)})^\theta\,e^{-\theta\varepsilon_DN}$ by Jensen's inequality.

\ref{lem:gauss-bounds-e} By the Cauchy--Schwarz inequality and~\ref{lem:gauss-bounds-a},
$\Qth(A)=\E[\ZN^\theta\ind_A]/\E\ZN^\theta\le(\E\ZN^{2\theta})^{1/2}\Prob(A)^{1/2}$.

For the last assertion, compare the exponents in~\ref{lem:gauss-bounds-a}
and~\ref{lem:gauss-bounds-c} (with $\lambda=2\theta$) with $C(\xi,\theta)$, using
$C^{\max}_\xi\ge\sqrt{2K_\xi/\pi}$.
\end{proof}

\subsection{The band theorem at every exponential rate}\label{app:gauss-truncation}

We prove Lemma~\ref{lem:max-truncation}. Recall the band free energy
$F^{\mathrm{band}}_n(m;\varepsilon,\delta)$ of \eqref{eq:max-band}, and that
$F^{\mathrm{band}}_n(m)-\FTAP(m)=N\bigl(\mathrm{TAP}_{N,n}(m,\varepsilon,\delta)-\TAPcps(\mu_m)\bigr)$
in the notation of \cite[(2.14), (2.24)--(2.25)]{CPS} (Section~\ref{sec:max-ramsey}). As in
Section~\ref{sec:intro-structure}, the exceptional events are bounded in outer probability.

\paragraph{Quantifiers in the band theorem.} For $\xi\in\classC$, \cite[Theorem~1]{CPS} states:
for all $c,t>0$, if $\varepsilon,\delta>0$ are small enough and $n$ is large enough, then for all
large $N$, outside an event of probability at most $e^{-cN}$, every $m\in\cube$ satisfies
$|F^{\mathrm{band}}_n(m;\varepsilon,\delta)-\FTAP(m)|<tN$. Lemma~\ref{lem:max-truncation} asks in
addition for $n\ge L/\min(\varepsilon,\delta)$, so the order of the choices matters. The proof in
\cite[Section~4]{CPS} allows it. The bound $F^{\mathrm{band}}_n<\FTAP+tN$ is proved in
\cite[Section~4.1]{CPS} by contradiction along sequences $\varepsilon_N,\delta_N\to0$ and
$n_N\to\infty$, so it holds for all sufficiently small $\varepsilon,\delta$ and all sufficiently
large $n$. The bound $F^{\mathrm{band}}_n>\FTAP-tN$ is proved in \cite[Section~4.2]{CPS} for every
triple $(\varepsilon,\delta,n)$ that admits an auxiliary parameter $\eta<\varepsilon/6$, small
enough for the uniform continuity of $\TAPcps$, with
$(t/4)^2c_\xi/(1/n+\delta+\varepsilon/2)>2c+\log(2/\eta)$; here $c_\xi$ is the constant of the
concentration bound \cite[Lemma~18]{CPS}. With $\delta=\varepsilon$ and $\eta$ a fixed multiple of
$\varepsilon$, the left side grows like $1/\varepsilon$ and the right side like
$\log(1/\varepsilon)$, so the condition holds once $\varepsilon$ is small enough and $n$ is large
enough depending on $\varepsilon$. Hence, for every $L\ge1$, one may take $\varepsilon=\delta$
small and then $n$ large, with $n\ge L$ and $n\ge L/\varepsilon$. For a polynomial $\xi$ this is
Lemma~\ref{lem:max-truncation}.

\paragraph{Truncation.} Under \condR\ alone we apply this to the polynomial $\xi_{\le D}$, and the
remainder $\HN^{>D}$ is controlled by Lemma~\ref{lem:gauss-bounds}\ref{lem:gauss-bounds-d}. The
degree $D$ has to grow with the rate $c$: by
Lemma~\ref{lem:gauss-bounds}\ref{lem:gauss-bounds-d2} the remainder $T_D$ exceeds its mean by $uN$
with probability at most $e^{-u^2N/(2\xi_{>D}(1))}$, and this is smaller than $e^{-cN}$ only when
$\xi_{>D}(1)$ is small compared with $u^2/c$. This is why the band theorem at every exponential
rate under \condR\ needs a truncation level that depends on the rate, chosen before the band
parameters.

\begin{proof}[Proof of Lemma~\ref{lem:max-truncation}]
Fix $c,t>0$ and $L\ge1$.

\emph{Choice of $D$.} Choose $D$, at least the smallest active degree, with
\[
 2\varepsilon_D+\tfrac32\xi'_{>D}(1)\le t/4\qquad\text{and}\qquad
 (t/8)^2\ge2(c+1)\,\xi_{>D}(1).
\]
This is possible because $\varepsilon_D$, $\xi'_{>D}(1)$ and $\xi_{>D}(1)$ tend to zero. The
choice depends only on $\xi$, $c$ and $t$.

\emph{The truncated band theorem.} Apply the polynomial case to $\xi_{\le D}$, with rate $c+1$,
accuracy $t/2$ and the given $L$. This gives $\varepsilon,\delta>0$ and $n$ with $n\ge L$ and
$n\ge L/\min(\varepsilon,\delta)$ such that, for all large $N$, the event
\[
 \Omega^{(D)}=\Bigl\{\bigl|F^{\mathrm{band},(D)}_n(m;\varepsilon,\delta)-\FTAP^{(D)}(m)\bigr|
 <tN/2\ \text{for every } m\in\cube\Bigr\}
\]
has a complement of outer probability at most $e^{-(c+1)N}$. Here $F^{\mathrm{band},(D)}_n$ is
\eqref{eq:max-band} with $\HN^{\le D}$ in place of $\HN$.

\emph{The remainder.} Let $\Omega'=\{T_D\le(\varepsilon_D+t/8)N\}$. By
Lemma~\ref{lem:gauss-bounds}\ref{lem:gauss-bounds-d2} with $u=t/8$ and the choice of $D$,
$\Prob((\Omega')^c)\le e^{-(c+1)N}$.

\emph{Comparison.} On $\Omega^{(D)}\cap\Omega'$, fix $m\in\cube$. Since
$F^{\mathrm{band},(D)}_n(m)>\FTAP^{(D)}(m)-tN/2>-\infty$, the band $\Band$ is nonempty. The band
does not depend on the field, and for every $(\sigma^1,\dots,\sigma^n)$ in it the exponents
$\sum_{i\le n}\HN(\sigma^i)$ and $\sum_{i\le n}\HN^{\le D}(\sigma^i)$ differ by at most $nT_D$; by
\eqref{eq:max-band}, $|F^{\mathrm{band}}_n(m)-F^{\mathrm{band},(D)}_n(m)|\le T_D$. By
Lemma~\ref{lem:gauss-bounds}\ref{lem:gauss-bounds-d1},
$|\FTAP(m)-\FTAP^{(D)}(m)|\le T_D+\frac32N\xi'_{>D}(1)$. Therefore
\[
 \bigl|F^{\mathrm{band}}_n(m)-\FTAP(m)\bigr|<\frac{tN}2+2T_D+\frac32N\xi'_{>D}(1)
 \le\frac{tN}2+\frac{tN}4+\Bigl(2\varepsilon_D+\frac32\xi'_{>D}(1)\Bigr)N\le tN .
\]
So the event in \eqref{eq:max-band-thm} is contained in $(\Omega^{(D)})^c\cup(\Omega')^c$, whose
outer probability is at most $2e^{-(c+1)N}\le e^{-cN}$ for $N\ge\log2$.
\end{proof}

Propositions~\ref{prop:max-lower} and~\ref{prop:max-upper} are deduced from
Lemma~\ref{lem:max-truncation} in Section~\ref{sec:max-ramsey}: the upper band uses the half
$\FTAP<F^{\mathrm{band}}_n+tN$ together with \eqref{eq:max-band-le}, and the lower band uses the
other half at a barycenter supplied by Ramsey's theorem.

\subsection{Maximizers are interior critical points}\label{app:gauss-critical}

Theorem~\ref{thm:exist}\ref{thm:exist-a} has two parts. Entropy pushes the maximizers of $\FTAP$
off the faces of the cube (Lemma~\ref{lem:gauss-interior}). At an interior point the generalized
TAP free energy and the TAP function with the minimizing order parameter frozen have the same
directional derivatives, and the frozen function is continuously differentiable
(Lemma~\ref{lem:gauss-critical}). The second part is the envelope argument by which Chen,
Panchenko and Subag compute the gradient of the TAP correction \cite[Theorem~5 and
Section~11]{CPS}; we apply it at the maximizers and observe that it also applies to the frozen
function.

For $\xi\in\classC$ we use the inverse field of
Lemma~\ref{lem:pde-standing}\ref{lem:pde-standing-b}: for $\zeta\in\Pspace([0,1])$ and
$q\in[0,1]$, $b_\zeta(q,\cdot)$ is the inverse of the bijection
$\partial_x\Phi_\zeta(q,\cdot)\colon\R\to(-1,1)$, and
$u_\zeta(q,m)=\partial_x^2\Phi_\zeta(q,b_\zeta(q,m))$. In the notation of \cite[(2.22)]{CPS},
$b_\zeta(q,m)=\Psi(q,m,\zeta)$. For $\zeta\in\Pspace([0,1])$, $q\in[0,1]$ and $z\in\cube$ put
\[
 G_\zeta(q,z)=\sum_{j=1}^N\Phistar_\zeta(q,z_j)-\frac N2\int_q^1t\,\xi''(t)\,\cdf\zeta(t)\,dt,
 \qquad\text{so that}\qquad F_\zeta(z)=\HN(z)+G_\zeta(q_z,z).
\]

\begin{lemma}[Continuity, attainment and interiority]\label{lem:gauss-interior}
\leavevmode
\begin{thmlist}
\item\label{lem:gauss-interior-a} Assume \condR. Almost surely, for every $N$, $\HN$ is continuously
 differentiable on a neighbourhood of $\cube$, $\FTAP$ is continuous on $\cube$, and $\MN$ is
 attained. In particular $\MN$ is a random variable.
\item\label{lem:gauss-interior-b} Assume $\xi\in\classC$. Fix $N$ and a realization in which $\HN$ is
 continuously differentiable on a neighbourhood of $\cube$. Let
 $\gamma(b)=\E\tanh(b+\sqrt{\xi'(1)}\,|Z|)$, with $Z$ standard Gaussian, which is an increasing
 bijection of $\R$ onto $(-1,1)$, and put $\psi_\xi(m)=\gamma^{-1}(m)-\xi''(1)$ for $m\in(-1,1)$,
 so that $\psi_\xi$ is increasing and $\psi_\xi(m)\to\infty$ as $m\uparrow1$. If $x\in\cube$ has
 $x_i=1$ and $y=x-se_i$ with $s\in(0,1)$, then
 \begin{equation}\label{eq:gauss-push}
  \FTAP(y)-\FTAP(x)\ \ge\ s\Bigl[\psi_\xi(1-s)-3\xi''(1)-\sup_{\cube}|\partial_i\HN|\Bigr],
 \end{equation}
 and the same holds if $x_i=-1$ and $y=x+se_i$. Consequently every maximizer of $\FTAP$ lies in
 $\ocube$.
\end{thmlist}
\end{lemma}

\begin{proof}
\ref{lem:gauss-interior-a} The first assertion is Lemma~\ref{lem:gauss-field}\ref{lem:gauss-field-a},
for each $N$ outside a null set. By Section~\ref{sec:intro-setting},
$\FTAP(x)=\HN(x)+N\TAPcps(\mu_x)-N\log2$. The functional $\TAPcps(\mu)$ depends on $\mu$ only
through the law of $|m|$ under $\mu$, because $\Phistar_\zeta(q,\cdot)$ is even, and $m\mapsto|m|$
does not increase $\Wone$. For a polynomial covariance, $\TAPcps$ is therefore $\Wone$-continuous on
$\Pspace([-1,1])$ by \cite[Theorem~12(ii)]{CPS}. Under \condR, $\TAPcps$ is the uniform limit of
$\TAPcps^{(D)}$ as $D\to\infty$ by Lemma~\ref{lem:gauss-bounds}\ref{lem:gauss-bounds-d1}, hence
continuous. Since $\Wone(\mu_x,\mu_y)\le N^{-1}\sum_i|x_i-y_i|$, the map $x\mapsto\FTAP(x)$ is
continuous on the compact set $\cube$, so its maximum is attained and equals its supremum over a
countable dense subset.

\ref{lem:gauss-interior-b} \emph{Step 1: a lower bound on the inverse field.} Let
$\zeta\in\Pspace([0,1])$, $q\in[0,1]$ and $b\in\R$. By
Lemma~\ref{lem:pde-ac}\ref{lem:pde-ac-a} for measures with finitely many atoms, and by
Lemma~\ref{lem:pde-measurable}\ref{lem:pde-measurable-i}--\ref{lem:pde-measurable-ii} in general,
$\partial_x\Phi_\zeta(q,b)=\E\tanh(X_1)$, where $X$ solves
$dX_r=\xi''(r)\cdf\zeta(r)\,\partial_x\Phi_\zeta(r,X_r)\,dr+\sqrt{\xi''(r)}\,dB_r$ on $[q,1]$ with
$X_q=b$. Since $0\le\cdf\zeta\le1$ and $|\partial_x\Phi_\zeta|\le1$, the drift is at most
$\xi''(r)\le\xi''(1)$ on an interval of length at most one, so $X_1\le b+\xi''(1)+W$ with
$W=\int_q^1\sqrt{\xi''(r)}\,dB_r$ centered Gaussian of variance $\xi'(1)-\xi'(q)\le\xi'(1)$. Since
$\tanh$ is increasing, $\partial_x\Phi_\zeta(q,b)\le\E\tanh(b+\xi''(1)+|W|)\le\gamma(b+\xi''(1))$.
Applied at $b=b_\zeta(q,m)$ this gives $m\le\gamma(b_\zeta(q,m)+\xi''(1))$, that is,
\begin{equation}\label{eq:gauss-inverse-lower}
 b_\zeta(q,m)\ge\psi_\xi(m)\qquad\text{for all } \zeta\in\Pspace([0,1]),\ q\in[0,1],\ m\in(-1,1).
\end{equation}

\emph{Step 2: Lipschitz dependence on the self-overlap.} For $\zeta$ with finitely many atoms,
$\Phi_\zeta$ is continuous in $t$ and solves the Parisi equation on each closed step
(Lemma~\ref{lem:pde-class}\ref{lem:pde-class-b}), with
$|\partial_t\Phi_\zeta|=\frac{\xi''}2|\partial_x^2\Phi_\zeta+\cdf\zeta(\partial_x\Phi_\zeta)^2|\le
\xi''(1)$. By uniform approximation, $|\Phi_\zeta(q,b)-\Phi_\zeta(q',b)|\le\xi''(1)|q-q'|$ for every
$\zeta$. Taking infima, $|\Phistar_\zeta(q,m)-\Phistar_\zeta(q',m)|\le\xi''(1)|q-q'|$ for
$|m|<1$, and for $m=\pm1$ by continuity.

\emph{Step 3: moving a coordinate inward.} Let $x_i=1$, $y=x-se_i$ and $\zeta\in\Pspace([0,1])$.
The function $\Phistar_\zeta(q,\cdot)$ is concave and continuous on $[-1,1]$, and for
$m,m'\in(-1,1)$,
$\Phistar_\zeta(q,m')\le\Phi_\zeta(q,b_\zeta(q,m))-m'b_\zeta(q,m)=\Phistar_\zeta(q,m)-(m'-m)\,
b_\zeta(q,m)$. So $-b_\zeta(q,m)$ is a supergradient at $m$, and since $b_\zeta(q,\cdot)$ is
continuous, $\Phistar_\zeta(q,\cdot)$ is continuously differentiable on $(-1,1)$ with derivative
$-b_\zeta(q,\cdot)$ (compare \cite[(2.21)--(2.23)]{CPS}). By \eqref{eq:gauss-inverse-lower} and the
monotonicity of $\psi_\xi$, with an improper integral at $m=1$,
\[
 G_\zeta(q_y,y)-G_\zeta(q_y,x)=\Phistar_\zeta(q_y,1-s)-\Phistar_\zeta(q_y,1)
 =\int_{1-s}^1b_\zeta(q_y,m)\,dm\ge s\,\psi_\xi(1-s).
\]
Next, $|q_y-q_x|=(2s-s^2)/N\le2s/N$, and Step~2 together with $t\,\xi''(t)\cdf\zeta(t)\le\xi''(1)$
gives $|G_\zeta(q_y,x)-G_\zeta(q_x,x)|\le\frac32N\xi''(1)|q_y-q_x|\le3\xi''(1)s$. Finally
$\HN(y)-\HN(x)\ge-s\sup_{\cube}|\partial_i\HN|$. Adding the three bounds,
\[
 F_\zeta(y)-F_\zeta(x)\ge s\Bigl[\psi_\xi(1-s)-3\xi''(1)-\sup_{\cube}|\partial_i\HN|\Bigr]
 \qquad\text{for every } \zeta .
\]
Since $\FTAP(x)\le F_\zeta(x)$ for every $\zeta$, $\FTAP(y)-\FTAP(x)\ge\inf_\zeta(F_\zeta(y)-
F_\zeta(x))$, which proves \eqref{eq:gauss-push}. The case $x_i=-1$ is the same, because
$\Phistar_\zeta(q,\cdot)$ is even. The supremum of $|\partial_i\HN|$ is finite and
$\psi_\xi(1-s)\to\infty$ as $s\downarrow0$, so the right side of \eqref{eq:gauss-push} is positive
for small $s$, and $x$ is not a maximizer.
\end{proof}

For $y\in\ocube$, $v\in\R^N$ and $|t|$ small, so that $y+tv\in\ocube$, put $q_t=q_{y+tv}$ and,
for $\zeta\in\Pspace([0,1])$ and $h\in\R^N$,
\begin{equation}\label{eq:gauss-fy}
 f_y(t;\zeta,h)=\sum_{i=1}^N\bigl[\Phi_\zeta(q_t,h_i)-(y_i+tv_i)h_i\bigr]
 -\frac N2\int_{q_t}^1s\,\xi''(s)\,\cdf\zeta(s)\,ds .
\end{equation}
This is $N$ times the functional \cite[(11.1)]{CPS} at the point $y+tv$, minus $N\log2$. By the
definition of $\Phistar_\zeta$,
\[
 F_\zeta(y+tv)=\HN(y+tv)+\inf_{h\in\R^N}f_y(t;\zeta,h),\qquad
 \FTAP(y+tv)=\HN(y+tv)+\inf_{\zeta,h}f_y(t;\zeta,h).
\]

\begin{lemma}[Exact criticality]\label{lem:gauss-critical}
Assume $\xi\in\classC$. Fix $N$ and a realization in which $\HN$ is continuously differentiable on a
neighbourhood of $\cube$.
\begin{thmlist}
\item\label{lem:gauss-critical-a} For every $\zeta\in\Pspace([0,1])$ the function $F_\zeta$ is
 continuously differentiable on $\ocube$, with
 \begin{equation}\label{eq:gauss-gradient}
  \partial_iF_\zeta(y)=\partial_i\HN(y)-b_\zeta(q_y,y_i)
  -y_i\,\xi''(q_y)\,\frac1N\sum_{j=1}^Nu_\zeta(q_y,y_j),\qquad y\in\ocube .
 \end{equation}
\item\label{lem:gauss-critical-b} Let $x\in\ocube$, and let $\zeta_x$ be the unique minimizer of
 $\mu\mapsto F_\mu(x)$ over $\Pspace([q_x,1])$ \cite[Theorem~12(iii) and Section~11]{CPS}. For every
 $v\in\R^N$, the function $t\mapsto\FTAP(x+tv)$ is differentiable at $t=0$, with derivative
 $\nabla F_{\zeta_x}(x)\cdot v$.
\item\label{lem:gauss-critical-c} If $x\in\ocube$ is a local maximizer of $\FTAP$, then
 $\nabla F_{\zeta_x}(x)=0$ and $F_{\zeta_x}(x)=\FTAP(x)$.
\end{thmlist}
\end{lemma}

\begin{proof}
\emph{Minimizers in $h$.} By Lemma~\ref{lem:pde-standing}\ref{lem:pde-standing-b}, for
$m\in(-1,1)$ the function $b\mapsto\Phi_\zeta(q,b)-mb$ is strictly convex with critical point
$b_\zeta(q,m)$. Hence, for fixed $\zeta$, the unique minimizer of $f_y(t;\zeta,\cdot)$ is
$h_t=(b_\zeta(q_t,y_i+tv_i))_{i\le N}$, which is continuous in $t$ because $b_\zeta$ is jointly
continuous on $[0,1]\times(-1,1)$.

\emph{One-sided derivatives.} For fixed $(\zeta,h)$, the right and left derivatives in $q$ of
$q\mapsto\sum_i\Phi_\zeta(q,h_i)-\frac N2\int_q^1s\xi''(s)\cdf\zeta(s)\,ds$ are computed in
\cite[Section~11]{CPS} from the Parisi equation:
\begin{equation}\label{eq:gauss-qderiv}
 -\frac{\xi''(q)}2\sum_{i=1}^N\partial_x^2\Phi_\zeta(q,h_i)
 +\frac{\xi''(q)}2\,\cdf\zeta^{\pm}(q)\Bigl(Nq-\sum_{i=1}^N\bigl(\partial_x\Phi_\zeta(q,h_i)\bigr)^2\Bigr),
\end{equation}
where $\cdf\zeta^+(q)=\zeta([0,q])=\cdf\zeta(q)$ and $\cdf\zeta^-(q)=\zeta([0,q))$. Since
$t\mapsto q_t$ is smooth, with $\frac{d}{dt}q_t=2(y+tv)\cdot v/N$, and $q\mapsto\Phi_\zeta(q,b)$ is
Lipschitz (Step~2 of the previous proof), $f_y(\cdot;\zeta,h)$ has right and left derivatives in $t$: the chain rule applies
with the one-sided $q$-derivative \eqref{eq:gauss-qderiv} taken on the side toward which $q_t$
moves, and the term $-\sum_iv_ih_i$ is added. The only factor that may be discontinuous is
$\cdf\zeta^{\pm}(q)$, and it multiplies $Nq-\sum_i(\partial_x\Phi_\zeta(q,h_i))^2$. At $(t,h)=(0,h_0)$,
where $h_{0,i}=b_\zeta(q_y,y_i)$ and hence $\partial_x\Phi_\zeta(q_y,h_{0,i})=y_i$, this factor
equals $Nq_y-|y|^2=0$ \cite[(11.9)]{CPS}. As $\partial_x\Phi_\zeta$ and $\partial_x^2\Phi_\zeta$ are
continuous on $[0,1]\times\R$ (Lemma~\ref{lem:pde-standing}\ref{lem:pde-standing-a}), both
one-sided $t$-derivatives of $f_y(\cdot;\zeta,\cdot)$ are therefore continuous at $(0,h_0)$, as
functions of $(t,h)$, with the common value
\begin{equation}\label{eq:gauss-common}
 \sum_{i=1}^Nv_i\Bigl[-b_\zeta(q_y,y_i)-y_i\,\xi''(q_y)\,\frac1N\sum_{j=1}^Nu_\zeta(q_y,y_j)\Bigr].
\end{equation}

\ref{lem:gauss-critical-a} The envelope lemma \cite[Lemma~37]{CPS}, applied on the metric space
$\R^N$ to $f_y(\cdot;\zeta,\cdot)$ and its minimizers $h_t$, for right and for left derivatives,
shows that $t\mapsto\inf_hf_y(t;\zeta,h)$ is differentiable at $t=0$ with derivative
\eqref{eq:gauss-common}. Adding $\nabla\HN(y)\cdot v$, every directional derivative of $F_\zeta$
at $y$ exists and is the linear form in $v$ given by \eqref{eq:gauss-gradient}. The right side of
\eqref{eq:gauss-gradient} is continuous on $\ocube$, because $b_\zeta$ and $u_\zeta$ are jointly
continuous on $[0,1]\times(-1,1)$ (Lemma~\ref{lem:pde-standing}\ref{lem:pde-standing-b}). A function
whose partial derivatives exist and are continuous on an open set is continuously differentiable
there.

\ref{lem:gauss-critical-b} Let $y=x$. In \cite[Section~11]{CPS} it is shown that the minimizer
$(\zeta_{x+tv},h_{x+tv})$ of $f_x(t;\cdot,\cdot)$ over $\Pspace([0,1])\times\R^N$, normalized by
$\zeta([0,q_t))=0$, is unique and continuous at $t=0$, and that the one-sided $t$-derivatives of
$f_x$ are continuous at $(0,(\zeta_x,h_x))$ jointly in $(t,\zeta,h)$, by the cancellation just
described; here $h_{x,i}=b_{\zeta_x}(q_x,x_i)$. The envelope lemma \cite[Lemma~37]{CPS} on
$\Pspace([0,1])\times\R^N$ then shows that $t\mapsto\FTAP(x+tv)-\HN(x+tv)$ has right and left
derivatives at $t=0$ equal to those of $f_x(\cdot;\zeta_x,\cdot)$ at $(0,h_x)$. By the proof
of~\ref{lem:gauss-critical-a} with $\zeta=\zeta_x$, both equal \eqref{eq:gauss-common}, which is the
derivative of $F_{\zeta_x}(x+tv)-\HN(x+tv)$ at $t=0$.

\ref{lem:gauss-critical-c} At a local maximizer in the open cube, the derivative in~\ref{lem:gauss-critical-b}
vanishes for every $v$. So $\nabla F_{\zeta_x}(x)=0$, and $F_{\zeta_x}(x)=\FTAP(x)$ by the choice of
$\zeta_x$.
\end{proof}

\begin{proof}[Proof of Theorem~\ref{thm:exist}\ref{thm:exist-a}]
Let $\xi\in\classC$. By Lemma~\ref{lem:gauss-interior}\ref{lem:gauss-interior-a}, there is an event
of probability one on which, for every $N$, $\HN$ is continuously differentiable on a neighbourhood
of $\cube$, $\FTAP$ is continuous and $\MN$ is attained. Fix a realization in this event and a
maximizer $x^*$. By Lemma~\ref{lem:gauss-interior}\ref{lem:gauss-interior-b}, $x^*\in\ocube$. With
$\zeta^*=\zeta_{x^*}$, Lemma~\ref{lem:gauss-critical} shows that $F_{\zeta^*}$ is continuously
differentiable on $\ocube$, that $\nabla F_{\zeta^*}(x^*)=0$, and that
$F_{\zeta^*}(x^*)=\FTAP(x^*)=\MN$.
\end{proof}

\begin{remark}[The generalized TAP equations]\label{rem:gauss-tap}
By \cite[(11.12)]{CPS}, at the minimizing order parameter the factor
$\xi''(q_x)N^{-1}\sum_ju_{\zeta_x}(q_x,x_j)$ in the last term of \eqref{eq:gauss-gradient}
equals $\xi''(q_x)\int_{q_x}^1\zeta_x([0,s])\,ds$, so
\eqref{eq:gauss-gradient} at $(x,\zeta_x)$ is the gradient formula \cite[Theorem~5]{CPS} and
$\nabla F_{\zeta^*}(x^*)=0$ are the generalized TAP equations \cite[Remark~6]{CPS}.
Lemma~\ref{lem:gauss-critical} does not use this: away from $x$, the gradient of $F_{\zeta_x}$
keeps the empirical average.
\end{remark}

\subsection{Exponential moments and the upper tail of the maximum}\label{app:var-moments}

Section~\ref{sec:max-existence} obtained the local rate of $\MN$ at a level $f\in(\feq,\fone)$,
Theorem~\ref{thm:exist}\ref{thm:exist-d}, by tilting the law of the disorder with $e^{\thf\MN}$;
the inputs were the exponential moments of $\MN$ from Theorem~\ref{thm:max}\ref{thm:max-b} and
\condFM, and the differentiability of $\Lambda$ at $\thf$. The same inputs give the whole upper
tail of $\MN$ below $\fone$ and, through Theorem~\ref{thm:max}, that of the free energy.

\begin{proposition}[The upper tail of the TAP maximum]\label{prop:var-tail}
Assume \condH\ and $\xi\in\classC$. Let $I(f)=-\Sigma(f)$ for $f\ge\feq$ and $I(f)=+\infty$ for
$f<\feq$.
\begin{thmlist}
\item\label{prop:var-tail-a} For every $\theta\in[0,1]$,
 $\lim_N N^{-1}\log\E e^{\theta\MN}=\Lambar(\theta)$.
\item\label{prop:var-tail-b} For every $f\in\R$,
 $\limsup_N N^{-1}\log\Prob(\MN\ge Nf)\le\Sigma(f)$, and for every $f<\fone$,
 \[
  \lim_{N\to\infty}\frac1N\log\Prob\bigl(\MN\ge Nf\bigr)=\Sigma(f).
 \]
\item\label{prop:var-tail-c} For every closed set $A\subset\R$ and every open set $U\subset\R$,
 \[
  \limsup_{N\to\infty}\frac1N\log\Prob\bigl(N^{-1}\MN\in A\bigr)\le-\inf_AI,\qquad
  \liminf_{N\to\infty}\frac1N\log\Prob\bigl(N^{-1}\MN\in U\bigr)\ge-\inf_{U\cap(-\infty,\fone)}I .
 \]
\item\label{prop:var-tail-d} $\MN\ge N\feq$, and for every $\delta>0$ there is $c_\delta>0$ such
 that $\Prob(\MN\ge N(\feq+\delta))\le e^{-c_\delta N}$ for all large $N$.
\item\label{prop:var-tail-e} For every $N\ge1$ and $f\in\R$,
 $\Prob(\log\ZN\ge Nf)\le e^{N\Sigma(f)}$, and for every $f<\fone$,
 \[
  \lim_{N\to\infty}\frac1N\log\Prob\bigl(N^{-1}\log\ZN\ge f\bigr)=\Sigma(f).
 \]
\end{thmlist}
\end{proposition}

Parts~\ref{prop:var-tail-a}--\ref{prop:var-tail-d} use only \condR, $\xi\in\classC$ and
\condFM; convexity enters them only through Theorem~\ref{thm:exist}\ref{thm:exist-e}, which
supplies \condFM. The bound for every $N$ in~\ref{prop:var-tail-e} uses Guerra's bound, which needs
convexity of $\xi$ on $[-1,1]$. The open-set bound in~\ref{prop:var-tail-c} stops at $\fone$:
the argument has no lower bound for levels above $\fone$.

\begin{proof}
\ref{prop:var-tail-a} For $\theta\in(0,1)$ this is Theorem~\ref{thm:exist}\ref{thm:exist-c}
and~\ref{thm:exist-e}, and for $\theta=1$ it is the last assertion of
Theorem~\ref{thm:exist}\ref{thm:exist-c}, since $\Lambar(1)=\xi(1)/2$. For $\theta=0$ both sides
vanish.

\ref{prop:var-tail-b} Markov's inequality gives
$\Prob(\MN\ge Nf)\le e^{-\theta Nf}\E e^{\theta\MN}$ for $\theta\in[0,1]$, so
by~\ref{prop:var-tail-a} the upper limit is at most $\Lambar(\theta)-\theta f$; the infimum over
$\theta\in[0,1]$ is $\Sigma(f)$. For the lower bound, let $f<\fone$. If $f\le\feq$, then
$\Prob(\MN\ge Nf)=1$ by~\eqref{eq:max-crude-M} and $\Sigma(f)=0$. Let $\feq<f<\fone$ and
$0<\varepsilon<(\fone-f)/2$, and put $f'=f+\varepsilon\in(\feq,\fone)$. Since
$\{|N^{-1}\MN-f'|\le\varepsilon\}\subset\{\MN\ge Nf\}$, and since the lower limit in
Theorem~\ref{thm:exist}\ref{thm:exist-d} is nondecreasing in the window and tends to $\Sigma(f')$
as the window shrinks,
\[
 \liminf_{N\to\infty}\frac1N\log\Prob\bigl(\MN\ge Nf\bigr)
 \ge\liminf_{N\to\infty}\frac1N\log\Prob\bigl(|N^{-1}\MN-f'|\le\varepsilon\bigr)
 \ge\Sigma(f+\varepsilon).
\]
Let $\varepsilon\downarrow0$ and use the continuity of $\Sigma$.

\ref{prop:var-tail-c} By Proposition~\ref{prop:var-regular}\ref{prop:var-regular-c}, $I$ is
nondecreasing on $[\feq,\infty)$ and infinite below $\feq$. Let $A$ be closed. If $A$ does not
meet $[\feq,\infty)$, then $\Prob(N^{-1}\MN\in A)=0$ by~\eqref{eq:max-crude-M}. Otherwise let $a$
be the least point of $A\cap[\feq,\infty)$. Then $\{N^{-1}\MN\in A\}\subset\{\MN\ge Na\}$ almost
surely, and~\ref{prop:var-tail-b} gives the upper limit $\Sigma(a)=-I(a)=-\inf_AI$. Let $U$ be
open and $f\in U\cap[\feq,\fone)$, and choose $\delta>0$ with $(f-\delta,f+\delta)\subset U$. Then
\[
 \Prob\bigl(N^{-1}\MN\in U\bigr)\ge\Prob\bigl(\MN\ge Nf\bigr)-\Prob\bigl(\MN\ge N(f+\delta)\bigr).
\]
By Proposition~\ref{prop:var-regular}\ref{prop:var-regular-c}, $\Sigma(f+\delta)<\Sigma(f)$; fix
$\gamma>0$ with $3\gamma<\Sigma(f)-\Sigma(f+\delta)$. By~\ref{prop:var-tail-b}, for all large $N$
the first probability is at least $e^{N(\Sigma(f)-\gamma)}$, and the second is at most
$e^{N(\Sigma(f+\delta)+\gamma)}\le\frac12e^{N(\Sigma(f)-\gamma)}$. So the lower limit is at least
$\Sigma(f)-\gamma=-I(f)-\gamma$. Let $\gamma\downarrow0$ and take the supremum over
$f\in U\cap[\feq,\fone)$; the points of $U$ below $\feq$ have $I=+\infty$.

\ref{prop:var-tail-d} The first claim is~\eqref{eq:max-crude-M}. By
Proposition~\ref{prop:var-regular}\ref{prop:var-regular-c}, $\Sigma(\feq+\delta)<0$, and
\ref{prop:var-tail-b} gives the claim with $c_\delta=-\Sigma(\feq+\delta)/2$.

\ref{prop:var-tail-e} Let $\theta\in(0,1)$ and $\alpha\in\Pspace([0,1])$ with
$\alpha(\{0\})\ge\theta$. Guerra's bound, $(N\theta)^{-1}\log\E\ZN^\theta\le P_\alpha$ for every
$N$ (Lemma~\ref{lem:moments-guerra}; \cite[Theorem~6.1]{Tal07}), and Markov's inequality give
$\Prob(\log\ZN\ge Nf)\le e^{-\theta Nf}\E\ZN^\theta\le e^{N\theta(P_\alpha-f)}$. Taking the
infimum over $\alpha$ gives $e^{N(\Lambda(\theta)-\theta f)}$, and $\theta=0$ gives the trivial
bound $1$. The infimum of $\Lambda(\theta)-\theta f$ over $\theta\in[0,1)$ equals $\Sigma(f)$,
because $\Lambda(\theta)\to\Lambar(1)$ as $\theta\uparrow1$ (Lemma~\ref{lem:var-dual}). This
proves the bound for every $N$, and hence the upper half of the limit. For the lower half, let
$f<\fone$, $0<t<\fone-f$ and $\gamma>0$. Then
\[
 \Prob\bigl(\log\ZN\ge Nf\bigr)\ge\Prob\bigl(\MN\ge N(f+t)\bigr)
 -\Prob\bigl(\MN-\log\ZN>tN\bigr).
\]
By~\ref{prop:var-tail-b} the first term is at least $e^{N(\Sigma(f+t)-\gamma)}$ for all large $N$,
and by Theorem~\ref{thm:max}\ref{thm:max-a} the second is at most $e^{N(\Sigma(f+t)-\gamma-1)}$
for all large $N$. Hence the lower limit is at least $\Sigma(f+t)-\gamma$. Let $\gamma\downarrow0$
and $t\downarrow0$.
\end{proof}

\begin{remark}[Consequences for the free energy]\label{rem:var-tail}\label{rem:max-subsumed}
The limit in~\ref{prop:var-tail-e} follows from \condFM\ and is not new where \condFM\ is known;
the bound for every $N$ is Guerra's. For even $\xi$, \condFM\ is Talagrand's
\cite[Section~7]{Tal07}, and for $\xi$ finite on $\R$ it is stated in \cite[Theorem~2.1]{CGKLM};
Boursier notes that it makes $-\Sigma$ the upper-tail rate function of $N^{-1}\log\ZN$
\cite[Remark~1.3]{Bou26b} and, ``with a little more work'', gives a large deviation principle on
$[\feq,\infty)$ with rate function $-\Sigma$, the Legendre transform of $\Lambar$
\cite[Section~1.4]{Bou26a}. The upper bound in~\ref{prop:var-tail-b} is hers
\cite[Section~2.2]{Bou26b}; through TAP states she proves the lower half of~\ref{prop:var-tail-e}
at regular levels whose first positive contact is strictly Plefka stable \cite[Remark~1.4 and
Proposition~3.12]{Bou26b}. Here it holds at every level below $\fone$, whatever the stability of
the contacts; Theorem~\ref{thm:max} transfers every speed-$N$ statement about $N^{-1}\log\ZN$
to $N^{-1}\MN$ and, through Theorem~\ref{thm:exist}\ref{thm:exist-a}, its local rates at levels
of $(\feq,\fone)$ to the existence event. We do not claim a full large deviation principle for the
free energy: the open-set lower bound in~\ref{prop:var-tail-c} stops at $\fone$, and \condFM\
enters only through Theorem~\ref{thm:exist}\ref{thm:exist-e}.
\end{remark}

\subsection{Positive overlap of near-maximizers}\label{app:var-overlap}

Theorem~\ref{thm:max} says nothing about where the maximizers of $\FTAP$ sit. Some information is
cheap and holds for every mixture: a point of small self-overlap has TAP free energy
close to $\feq$ or below, so near-maximizers at a level above $\feq$ have self-overlap bounded
below. With $K_\xi=\xi'(1)+\xi''(1)$ as above, write $\vth(q)=q\xi'(q)-\xi(q)$, let
$p(\theta)=\inf\{P_\zeta:\zeta(\{0\})\ge\theta\}$, so that $\Lambda(\theta)=\theta p(\theta)$, and
for $\theta>0$ let $\ell_\theta=\liminf_N(N\theta)^{-1}\log\E\ZN^\theta$. The events below are
measurable because $\FTAP$ is continuous on the cube.

\begin{proposition}[Positive overlap of near-maximizers]\label{rem:max-overlap}
\leavevmode
\begin{thmlist}
\item\label{rem:max-overlap-a} Assume \condR. For every $f>\feq$ and $c>0$ there is
 $\kappa\in(0,1]$ such that, for every $N$,
 \[
  \Prob\bigl(\exists x\in\cube:\ q_x\le\kappa,\ \FTAP(x)\ge Nf\bigr)\le e^{-cN};
 \]
 any $\kappa$ with $\vth(\kappa)/2+\sqrt{K_\xi\kappa}<(f-\feq)/2$ and $(f-\feq)^2\ge8c\,\xi(\kappa)$
 works.
\item\label{rem:max-overlap-b} Assume \condR, and let $\theta,t>0$ with $\feq+t<\ell_\theta$.
 There are $\kappa,c'>0$ such that, for all large $N$,
 \[
  \Qth\bigl(\exists x\in\cube:\ q_x\le\kappa,\ \FTAP(x)\ge\log\ZN-tN\bigr)\le e^{-c'N}.
 \]
 Together with Corollary~\ref{cor:max-tilt}: with $\Qth$-probability at least $1-2e^{-c'N}$, the
 set $\{x:\FTAP(x)\ge\log\ZN-tN\}$ contains every maximizer of $\FTAP$ and lies in
 $\{q_x>\kappa\}$. If moreover $\theta\in(\thstar,1)$ and $\FMge{\theta}$ holds,
 then $\ell_\theta\ge p(\theta)>\feq$.
\item\label{rem:max-overlap-c} Assume $\xi\in\classC$ and \condFM, and let $f\in(\feq,\fone)$,
 $\varepsilon\in(0,(f-\feq)/4)$ and $c>0$. There is $\kappa>0$ such that, for all large $N$,
 \[
  \Prob\Bigl(\exists x\in\cube:\ q_x\le\kappa,\ \FTAP(x)\ge\MN-\varepsilon N\ \Bigm|\
  |N^{-1}\MN-f|\le\varepsilon\Bigr)\le e^{-cN}.
 \]
\end{thmlist}
\end{proposition}

Part~\ref{rem:max-overlap-a} is a comparison with the Parisi measure (Lemma~\ref{lem:var-small});
\ref{rem:max-overlap-b} follows by the change of measure of Corollary~\ref{cor:max-tilt}
and~\ref{rem:max-overlap-c} by dividing by the lower bound of Theorem~\ref{thm:exist}\ref{thm:exist-d}.

\begin{lemma}[TAP values at small self-overlap]\label{lem:var-small}
Assume \condR.
\begin{thmlist}
\item\label{lem:var-small-a} For every $N$ and every $x\in\cube$,
 \[
  \FTAP(x)\le\HN(x)+N\feq+\frac N2\,\vth(q_x).
 \]
\item\label{lem:var-small-b} For every $\kappa\in(0,1]$ and $u>0$,
 \[
  \Prob\Bigl(\sup\bigl\{\HN(x):\ x\in\cube,\ q_x\le\kappa\bigr\}>N\bigl(\sqrt{K_\xi\kappa}+u\bigr)\Bigr)
  \le\exp\Bigl(-\frac{u^2N}{2\xi(\kappa)}\Bigr).
 \]
\end{thmlist}
\end{lemma}

\begin{proof}
\ref{lem:var-small-a} Let $q=q_x$ and $\zeta\in\Pspace([0,1])$. By~\eqref{eq:intro-ftap},
$\FTAP(x)\le F_\zeta(x)$. Three bounds compare $F_\zeta(x)$ with $NP_\zeta$.
\begin{steplist}
\item $\Phistar_\zeta(q,a)\le\Phi_\zeta(q,0)$ for every $a\in[-1,1]$: take $b=0$ in the infimum
 that defines $\Phistar_\zeta(q,a)$ for $|a|<1$, and pass to the continuous extension at $a=\pm1$.
\item $\Phi_\zeta(q,0)\le\Phi_\zeta(0,0)$. When $\zeta$ has finitely many atoms,
 $\Phi_\zeta(t,\cdot)$ is convex and
 $\partial_t\Phi_\zeta=-\frac{\xi''}2(\partial_x^2\Phi_\zeta+\zeta([0,t])(\partial_x\Phi_\zeta)^2)\le0$
 on each step of the Cole--Hopf recursion, so $t\mapsto\Phi_\zeta(t,0)$ is nonincreasing; this
 passes to general $\zeta$ by the uniform continuity of $\zeta\mapsto\Phi_\zeta$
 \cite[Lemma~14]{JT16}.
\item Since $\zeta([0,t])\le1$ and $\int_0^qt\xi''(t)\,dt=\vth(q)$,
 $\int_q^1t\xi''(t)\zeta([0,t])\,dt\ge\int_0^1t\xi''(t)\zeta([0,t])\,dt-\vth(q)$.
\end{steplist}
Together, $F_\zeta(x)\le\HN(x)+NP_\zeta+\frac N2\vth(q)$. Take the infimum over $\zeta$.

\ref{lem:var-small-b} The set $S_\kappa=\{x\in\cube:q_x\le\kappa\}$ is compact and $\HN$ is
almost surely continuous, so the supremum is a supremum over a countable dense subset, and its law
depends only on the covariance of $\HN$. For $x,y\in\cube$,
$\E(\HN(x)-\HN(y))^2\le K_\xi|x-y|^2$ by Lemma~\ref{lem:gauss-field}\ref{lem:gauss-field-b}. The
Sudakov--Fernique inequality, comparing with $\sqrt{K_\xi}\,g\cdot x$ for a standard Gaussian vector
$g$, gives
\[
 \E\sup_{S_\kappa}\HN\le\sqrt{K_\xi}\,\E\sup_{|x|\le\sqrt{\kappa N}}g\cdot x
 =\sqrt{K_\xi\kappa N}\,\E|g|\le\sqrt{K_\xi\kappa}\,N .
\]
On $S_\kappa$ the variance of $\HN(x)$ is $N\xi(q_x)\le N\xi(\kappa)$, since $\xi$ is
nondecreasing on $[0,1]$. The Borell--TIS inequality~\eqref{eq:gauss-conc} with deviation $uN$
gives the claim.
\end{proof}

\begin{proof}[Proof of Proposition~\ref{rem:max-overlap}]
\ref{rem:max-overlap-a} Let $u=(f-\feq)/2$ and let $\kappa\in(0,1]$ satisfy
$\vth(\kappa)/2+\sqrt{K_\xi\kappa}<u$ and $u^2\ge2c\,\xi(\kappa)$; such $\kappa$ exist because
$\vth(\kappa)\le\kappa\xi'(1)$ and $\xi(\kappa)\le\kappa^2\xi(1)$. Outside the event of
Lemma~\ref{lem:var-small}\ref{lem:var-small-b}, which has probability at most
$e^{-u^2N/(2\xi(\kappa))}\le e^{-cN}$, every $x\in\cube$ with $q_x\le\kappa$ satisfies, by
Lemma~\ref{lem:var-small}\ref{lem:var-small-a} and the monotonicity of $\vth$ on $[0,1]$,
\[
 \FTAP(x)\le N\bigl(\sqrt{K_\xi\kappa}+u\bigr)+N\feq+\frac N2\vth(\kappa)<N(\feq+2u)=Nf .
\]
Both bounds hold for every $N$.

\ref{rem:max-overlap-b} Let $C=C(\xi,\theta)$ be the constant of~\eqref{eq:max-crude}, so that
$\Qth(B)\le e^{2CN}\Prob(B)^{1/2}$ for every event $B$ (proof of Corollary~\ref{cor:max-tilt}). The
level $\ell_\theta$ is finite because $\E\ZN^\theta\le(\E\ZN^{2\theta})^{1/2}\le e^{CN/2}$. Fix
$a\in(\feq+t,\ell_\theta)$, put $\gamma=(a-t-\feq)/2>0$, and choose $\kappa\in(0,1]$ with
$\vth(\kappa)/2+\sqrt{K_\xi\kappa}<\gamma$ and $\gamma^2\ge4(2C+1)\xi(\kappa)$. Consider
\[
 B_1=\{\log\ZN\ge Na\},\qquad
 B_2=\bigl\{\sup\{\HN(x):x\in\cube,\ q_x\le\kappa\}\le N(\sqrt{K_\xi\kappa}+\gamma)\bigr\}.
\]
For all large $N$, $(N\theta)^{-1}\log\E\ZN^\theta\ge(\ell_\theta+a)/2$, so
$\Qth(B_1^c)\le e^{\theta Na}/\E\ZN^\theta\le e^{-\theta N(\ell_\theta-a)/2}$. By
Lemma~\ref{lem:var-small}\ref{lem:var-small-b},
$\Qth(B_2^c)\le e^{2CN}e^{-\gamma^2N/(4\xi(\kappa))}\le e^{-N}$. On $B_1\cap B_2$, every $x$ with
$q_x\le\kappa$ satisfies, by Lemma~\ref{lem:var-small}\ref{lem:var-small-a},
\[
 \FTAP(x)\le N\bigl(\sqrt{K_\xi\kappa}+\gamma\bigr)+N\feq+\frac N2\vth(\kappa)<N(\feq+2\gamma)
 =N(a-t)\le\log\ZN-tN .
\]
So the event of~\ref{rem:max-overlap-b} lies in $B_1^c\cup B_2^c$, whose $\Qth$-probability is at
most $e^{-c'N}$ with $c'=\frac12\min\{\theta(\ell_\theta-a)/2,1\}$, for all large $N$. By
Corollary~\ref{cor:max-tilt}, $\Qth(\MN<\log\ZN-tN)\le e^{-c'N}$ for all large $N$. Off both events,
the maximum $\MN$ is attained, every maximizer lies in $\{x:\FTAP(x)\ge\log\ZN-tN\}$, and this set
lies in $\{q_x>\kappa\}$. Finally, if $\theta\in(\thstar,1)$ and $\FMge{\theta}$
holds, then $\ell_\theta\ge\Lambda(\theta)/\theta=p(\theta)$; and $p(\theta)>\feq$, since
otherwise $\zeta_\theta$ would minimize $P$ and equal $\zP$ by
Lemma~\ref{lem:var-basic}\ref{lem:var-basic-a}, while
$\zeta_\theta(\{0\})\ge\theta>\thstar=\zP(\{0\})$.

\ref{rem:max-overlap-c} The lower limit in Theorem~\ref{thm:exist}\ref{thm:exist-d} is
nondecreasing in $\varepsilon$ and tends to $\Sigma(f)$ as $\varepsilon\downarrow0$, so it is at
least $\Sigma(f)$ for each fixed $\varepsilon$. Hence, for all large $N$, the conditioning event
$\{|N^{-1}\MN-f|\le\varepsilon\}$ has probability at least $e^{N(\Sigma(f)-1)}$. On it, $\FTAP(x)\ge\MN-\varepsilon N$ implies
$\FTAP(x)\ge N(f-2\varepsilon)$, and $f-2\varepsilon>\feq$ because $\varepsilon<(f-\feq)/4$.
Apply~\ref{rem:max-overlap-a} at the level $f-2\varepsilon$ with the rate $c+1-\Sigma(f)$, which
is positive because $\Sigma\le0$; it gives $\kappa$ such that the unconditional probability of the
event is at most $e^{-(c+1-\Sigma(f))N}$ for every $N$. Dividing gives the bound $e^{-cN}$.
\end{proof}

\subsection{The class of covariances}\label{app:gauss-rpde}

\begin{remark}[Scope of the cited facts on the Parisi equation]\label{rem:rpde}
The class $\classC$ enters our proofs only through facts on the Parisi equation and on the
functionals $\TAPcps(\mu,\zeta)$ that we cite from \cite{CPS} (Theorems~5 and~12, Lemma~37, and
Section~11 there), from \cite{AC15a,AC15b} and from \cite{Che17}. These papers state their models
under a summability condition such as $\sum_p2^p\beta_p^2<\infty$, which there serves to define the
Hamiltonian, while the cited statements involve $\xi$ only through $\xi'$ and $\xi''$ on $[0,1]$.
If they hold for every $\xi$ satisfying \condR, then every statement of this paper made for
$\xi\in\classC$ holds under \condR\ (together with \condH\ where it is assumed). We have not
checked their proofs in that generality, and we do not claim this extension. It would add the
covariances with $\sum_p2^p\beta_p^2=\infty$, all of which have radius of convergence at most two;
covariances of radius two with $\sum_p2^p\beta_p^2<\infty$ already lie in $\classC$. The rest of
this subsection gives details.
\end{remark}

\paragraph{Where the class enters.} The papers \cite{CPS}, \cite{AC15a}, \cite{AC15b} and
\cite{Che17} state their models under a summability condition such as $\sum_p2^p\beta_p^2<\infty$
(for \cite{CPS}, see (1.2) there). Two groups of statements from them are used only for
$\xi\in\classC$. The first consists of the continuity of the generalized TAP correction and the
uniqueness of the minimizing order parameter \cite[Theorem~12(ii)--(iii)]{CPS}, and the continuity
of the minimizers and the one-sided derivatives of \cite[Section~11]{CPS}, with the envelope lemma
\cite[Lemma~37]{CPS}. They enter Lemma~\ref{lem:gauss-critical}, hence
Theorem~\ref{thm:exist}\ref{thm:exist-a}; the uniqueness, with \cite[(2.26) and
Theorem~12(i)]{CPS}, also enters the order parameters of near-optimal states in Section~\ref{sec:hess}
(Proposition~\ref{prop:hess-own}, Appendix~\ref{app:matrix-own}). The second is the regularity
theory of the Parisi equation \cite[Proposition~1(i)]{AC15a}, which enters through
Lemma~\ref{lem:pde-standing}: in Lemmas~\ref{lem:gauss-interior}\ref{lem:gauss-interior-b}
and~\ref{lem:gauss-critical}, in Sections~\ref{sec:hess} and~\ref{sec:stab} and in
Appendices~\ref{app:matrix}
to~\ref{app:nc}. The strict convexity of the Parisi functional \cite[Theorem~2]{AC15b} and its first
variation \cite[Theorem~2]{Che17} are also stated in $\classC$, but they impose no restriction:
Theorem~\ref{thm:exist}\ref{thm:exist-b}, the uniqueness of the constrained minimizers (in
particular of $\zf$) and the regularity of $\Lambda$ rest only on properties (V1)--(V3) of
Appendix~\ref{app:var-constrained}, which hold under \condR\ by \cite{JT16,JT17}. Within
Theorem~\ref{thm:exist}, the class $\classC$ is therefore used only outside part~\ref{thm:exist-b}.

Each statement of the two groups concerns the Parisi equation, the Parisi functional, or the
functionals $\TAPcps(\mu,\zeta)$ of \cite[(2.24)]{CPS} and their minimizers; none involves $\HN$.
They depend on $\xi$ only through $\xi'$ and $\xi''$ on $[0,1]$, namely the coefficient $\xi''$ of
the Parisi equation and the integral $\int_q^1s\xi''(s)\cdf\zeta(s)\,ds$. Under \condR\ these
functions are continuous on $[0,1]$, and $\xi''>0$ on $(0,1]$. Appendices~\ref{app:pde-toolbox}
and~\ref{app:pde-cascade}, which construct the Parisi solutions for general order parameters and
prove the comparison inequalities that we use, use $\xi$ only on $[0,1]$ and assume no summability
condition; and Appendices~\ref{app:gauss-bounds}--\ref{app:gauss-truncation} and
Lemma~\ref{lem:gauss-interior}\ref{lem:gauss-interior-a} use \cite[Theorems~1 and~12(ii)]{CPS}
only for the polynomial covariances $\xi_{\le D}$.

\paragraph{Why truncation does not help.} Theorem~\ref{thm:max} meets the same restriction, since
\cite[Theorem~1]{CPS} is stated in $\classC$, and Lemma~\ref{lem:max-truncation} removes it by
truncating $\xi$ in the degree: its conclusion allows an error $tN$, and the remainder $\HN^{>D}$
changes $F^{\mathrm{band}}_n-\FTAP$ by at most $2T_D+\frac32N\xi'_{>D}(1)$.
Theorem~\ref{thm:exist}\ref{thm:exist-a} and the uniqueness statements of the first group are exact
for the given $\xi$: truncation moves $\FTAP$ by at most $T_D+\frac32N\xi'_{>D}(1)$, a small
multiple of $N$ outside an exponentially rare event
(Lemma~\ref{lem:gauss-bounds}\ref{lem:gauss-bounds-d1}--\ref{lem:gauss-bounds-d2}), but it does
not carry critical points of the truncated TAP functions to critical points of $F_\zeta$, and an
approximation in $\xi$ does not give uniqueness.

\paragraph{Other hypotheses.} Remarks~\ref{rem:kr-window} and~\ref{rem:kr-summable} discuss the
convexity and summability hypotheses of the Kac--Rice construction.

%% file: sections/appD-matrix.tex
\section{Random-matrix inputs}\label{app:matrix}

This appendix proves the random-matrix statements of Section~\ref{sec:hess} and the Gaussian
matrix lemmas used in Appendix~\ref{app:kr}. Appendices~\ref{app:matrix-gaussian}--\ref{app:matrix-pastur}
treat Gaussian matrices, the one-sided Hessians of $F_\alpha$, and reflected Pastur laws with their
edge. The decisive input is Lemma~\ref{lem:matrix-uniform} in Appendix~\ref{app:matrix-uniform}: a
Lipschitz spectral statistic of $\Hess\HN(x)-D$ concentrates with probability $1-e^{-cN^2}$, and
this beats the $e^{CN}$ elements of a net of points and diagonal matrices. Only the last two
subsections concern the rare event: Appendix~\ref{app:matrix-own} proves
Proposition~\ref{prop:hess-own} through Theorem~\ref{thm:struct}\ref{thm:struct-i}--\ref{thm:struct-ii},
and Appendix~\ref{app:matrix-index} proves the index form of the bulk through
Theorem~\ref{thm:struct}\ref{thm:struct-iii} and shows, with the Kac--Rice construction, that its
window cannot be widened.

\paragraph{Standing assumptions.}
Appendix~\ref{app:matrix-gaussian} and Lemma~\ref{lem:matrix-uniform} use only \condR. Everywhere
else we assume $\xi\in\classC$, which implies \condR: the regularity of the Parisi equation in
Lemma~\ref{lem:pde-standing} is proved for this class. Wherever $\zf$ and the objects defined from
it appear ($\Cf$, $\Lq{Q}$, $\UQ$, $\chiQ$, $\Plef$, $\rhoQ$, $\nuQ$, $\lamp$), the standing
assumptions of Section~\ref{sec:hess} are in force: \condH, $\xi\in\classC$, $f\in(\feq,\fone)$,
and $r_0\in(0,1/4)$ is fixed with $\Cf\subset[2r_0,1-2r_0]$. We write $\underline c$,
$\overline C$ and $C_{34}$ for the constants of Lemma~\ref{lem:pde-standing}; they depend only on
$\xi'(1)$.

\paragraph{Notation and elementary facts.}
The conventions of Section~\ref{sec:hess-pastur} are used throughout: $\GOE$ is a GOE matrix,
eigenvalues are ordered decreasingly, the empirical spectral distribution of a compression to
$x^\perp$ is normalized by $N-1$, and $\|A\|$ is the operator norm of a matrix $A$. For
$\mu,\mu'\in\Pspace(\R)$ let $d_{\mathrm K}(\mu,\mu')=\sup_t|\mu((-\infty,t])-\mu'((-\infty,t])|$
be the Kolmogorov distance; $\dL\le d_{\mathrm K}$. For $K>0$ let $\mathcal D_K$ be the set of
diagonal $N\times N$ matrices with entries in $[0,K]$. We use five elementary facts. For
symmetric $N\times N$ matrices $A,B,E$,
\begin{gather}
 \dL\bigl(\ESD(A),\ESD(A+E)\bigr)\le\|E\|,\label{eq:matrix-weyl}\\
 \dL\bigl(\ESD(A),\ESD(B)\bigr)\le d_{\mathrm K}\bigl(\ESD(A),\ESD(B)\bigr)\le\frac{\rank(A-B)}N .
 \label{eq:matrix-rank}
\end{gather}
The first holds because $|\lambda_i(A+E)-\lambda_i(A)|\le\|E\|$ for every $i$ (Weyl), the second
because $\lambda_{i+k}(A)\le\lambda_i(B)\le\lambda_{i-k}(A)$ with $k=\rank(A-B)$ (interlacing). For
real random variables $Y,Y'$ on one probability space, Markov's inequality at level
$s=(\E|Y-Y'|)^{1/2}$ gives $\Prob(Y\le t)\le\Prob(Y'\le t+s)+s$, and symmetrically; hence
\begin{equation}\label{eq:matrix-coupling}
 \dL\bigl(\Law(Y),\Law(Y')\bigr)\le\bigl(\E|Y-Y'|\bigr)^{1/2}.
\end{equation}
If $\varphi:\R\to[0,1]$ is monotone, then, writing $\varphi$ as an average of indicators of
half-lines,
\begin{equation}\label{eq:matrix-monotone}
 \Bigl|\int\varphi\,d\mu-\int\varphi\,d\mu'\Bigr|\le d_{\mathrm K}(\mu,\mu').
\end{equation}
Finally, $\dL$ metrizes weak convergence on $\Pspace(\R)$, so $\dL(\mu,\mu')\le\liminf_n\dL(\mu_n,\mu_n')$
whenever $\mu_n\to\mu$ and $\mu_n'\to\mu'$ weakly.

\subsection{Gaussian matrix estimates}\label{app:matrix-gaussian}

Three Gaussian facts enter. For $\sqrt a\,\GOE-D$ with $D$ deterministic and diagonal, the expected
top eigenvalue is bounded by an explicit one-dimensional infimum, and the resolvent satisfies a
self-consistent equation up to $O(N^{-1})$. At a fixed point, $\Hess\HN$ is a multiple of a GOE
matrix up to rank two, and its oscillation over the cube is controlled by the third derivative.

\begin{lemma}[Gaussian comparison and the resolvent equation]\label{lem:matrix-comparison}
Let $a\ge0$ and let $D_N$ be a deterministic real diagonal $N\times N$ matrix.
\begin{thmlist}
\item\label{lem:matrix-comparison-a} \textup{(Comparison.)}
 \begin{equation}\label{eq:matrix-comparison}
  \E\lambda_{\max}\bigl(\sqrt a\,\GOE-D_N\bigr)
  \le\inf_{t<\min_iD_{N,ii}}\Bigl\{-t+\frac aN\sum_{i\le N}\frac1{D_{N,ii}-t}\Bigr\}.
 \end{equation}
 In particular $\E\lambda_{\max}(\GOE)\le2$, and $\Prob(\|\GOE\|>2+\delta)\le2e^{-\delta^2N/4}$ for
 every $\delta>0$.
\item\label{lem:matrix-comparison-b} \textup{(Resolvent.)} For $z=E+\mathrm i\varepsilon$ with
 $\varepsilon>0$ let $G_N(z)=(D_N-\sqrt a\,\GOE-z)^{-1}$, $m_N(z)=N^{-1}\tr G_N(z)$ and
 $\bar m_N(z)=\E m_N(z)$. For deterministic $v,w\in\R^N$,
 \begin{equation}\label{eq:matrix-resolvent-var}
  \E\bigl|v^{\mathsf T}(G_N-\E G_N)w\bigr|^2\le\frac{4a|v|^2|w|^2}{N\varepsilon^4},\qquad
  \E|m_N-\bar m_N|^2\le\frac{4a}{N^2\varepsilon^4}.
 \end{equation}
 Moreover $(\E G_N)_{ij}=0$ for $i\ne j$, and
 \begin{equation}\label{eq:matrix-approx-dyson}
  \E G_{N,ii}=\frac1{D_{N,ii}-z-a\bar m_N}+O\Bigl(\frac1N\Bigr),\qquad
  \bar m_N=\frac1N\sum_{i\le N}\frac1{D_{N,ii}-z-a\bar m_N}+O\Bigl(\frac1N\Bigr),
 \end{equation}
 where the errors are bounded by $C_{a,\varepsilon}/N$, uniformly in $i$ and in $D_N$.
\item\label{lem:matrix-comparison-c} \textup{(Resolvent equation.)} For $\rho\in\Pspace(\R)$ and
 $z\in\Cplx^+$ the equation
 \begin{equation}\label{eq:matrix-dyson}
  m=\int\frac{\rho(dD)}{D-z-am},\qquad\operatorname{Im}m>0,
 \end{equation}
 has at most one solution. If $\ESD(D_N)\to\rho$ weakly, it has a solution $m(z)$, and
 $\bar m_N(z)\to m(z)$ and $m_N(z)\to m(z)$ in probability.
\end{thmlist}
\end{lemma}

\begin{proof}
\ref{lem:matrix-comparison-a} For unit vectors $v,w$ with $v\cdot w=r$, put
$X_v=\sqrt a\,v^{\mathsf T}\GOE v$ and $Y_v=2\sqrt{a/N}\,\gamma\cdot v$, with $\gamma$ a standard
Gaussian vector. From the covariance of $\GOE$,
\[
 \E(X_v-X_w)^2=\frac{4a}N(1-r^2)\le\frac{8a}N(1-r)=\E(Y_v-Y_w)^2 .
\]
The Sudakov--Fernique inequality compares expected suprema of two Gaussian processes with this
ordering of increment variances. Its interpolation proof applies unchanged when the same
deterministic function, here $-v^{\mathsf T}D_Nv$, is added to both processes, because that
function does not contribute to the derivative along the interpolation. On finite subsets of the
sphere, and then by continuity on the whole sphere,
\[
 \E\lambda_{\max}\bigl(\sqrt a\,\GOE-D_N\bigr)=\E\sup_{|v|=1}\{X_v-v^{\mathsf T}D_Nv\}
 \le\E\sup_{|v|=1}\{Y_v-v^{\mathsf T}D_Nv\}.
\]
For $t<\min_iD_{N,ii}$ and $|v|=1$, completing the square in each coordinate gives
$2\sqrt{a/N}\,\gamma\cdot v-v^{\mathsf T}D_Nv\le-t+\frac aN\sum_i\gamma_i^2/(D_{N,ii}-t)$, and
taking expectations proves~\eqref{eq:matrix-comparison}. For $D_N=0$ and $a=1$ the infimum is $2$,
attained at $t=-1$. As a function of the independent standard Gaussian variables that define
$\GOE$, the largest eigenvalue is $\sqrt{2/N}$-Lipschitz, so Gaussian concentration
\cite{Ledoux} gives $\Prob(\lambda_{\max}(\GOE)>2+\delta)\le e^{-\delta^2N/4}$, and the same holds
for $\lambda_{\max}(-\GOE)$ since $-\GOE$ has the law of $\GOE$.

\ref{lem:matrix-comparison-b} Write $\GOE=(W_{kl})$. The entries $W_{kl}$, $k\le l$, are independent, with variance
$1/N$ for $k<l$ and $2/N$ for $k=l$. Differentiating the inverse,
$\partial G/\partial W_{kl}=\sqrt a\,GE^{kl}G$ with $E^{kl}=e_ke_l^{\mathsf T}+e_le_k^{\mathsf T}$
for $k<l$ and $E^{kk}=e_ke_k^{\mathsf T}$. Since $G$ is symmetric,
$\partial(v^{\mathsf T}Gw)/\partial W_{kl}=\sqrt a\,[(Gv)_k(Gw)_l+(Gv)_l(Gw)_k]$ for $k<l$, and the
Gaussian Poincar\'e inequality, applied to real and imaginary parts, gives
$\Var(v^{\mathsf T}Gw)\le\frac{4a}N|Gv|^2|Gw|^2\le4a|v|^2|w|^2/(N\varepsilon^4)$, because
$\|G\|\le\varepsilon^{-1}$. For the trace, $\partial m_N/\partial W_{kl}=2\sqrt a\,(G^2)_{kl}/N$
for $k<l$ and $\sqrt a\,(G^2)_{kk}/N$ for $k=l$, so
$\Var m_N\le\frac{4a}{N^3}\|G^2\|_{\mathrm F}^2\le4a/(N^2\varepsilon^4)$. Conjugating by a diagonal
sign matrix preserves the law of $\GOE$ and fixes $D_N$, so $(\E G)_{ij}=0$ for $i\ne j$. From
$(D_{N,ii}-z)G_{ii}=1+\sqrt a\sum_kW_{ik}G_{ki}$ and Gaussian integration by parts, in which the
diagonal variance $2/N$ supplies the second copy of $G_{ii}^2$,
\begin{equation}\label{eq:matrix-ibp}
 (D_{N,ii}-z)\E G_{ii}=1+a\E[m_NG_{ii}]+\frac aN\E(G^2)_{ii}.
\end{equation}
By~\eqref{eq:matrix-resolvent-var} with $v=w=e_i$ and Cauchy--Schwarz,
$|\E[m_NG_{ii}]-\bar m_N\E G_{ii}|\le4a/(N^{3/2}\varepsilon^4)$, and
$|\frac aN(G^2)_{ii}|\le a/(N\varepsilon^2)$. Since $\operatorname{Im}\bar m_N\ge0$,
$|D_{N,ii}-z-a\bar m_N|\ge\varepsilon$, and~\eqref{eq:matrix-approx-dyson} follows by dividing and
averaging over $i$.

\ref{lem:matrix-comparison-c} Let $m_1,m_2$ be two solutions. Taking imaginary parts
in~\eqref{eq:matrix-dyson},
$a\int|D-z-am_j|^{-2}\rho(dD)=a\operatorname{Im}m_j/(\operatorname{Im}z+a\operatorname{Im}m_j)<1$.
Subtracting the two equations, if $m_1\ne m_2$ then
$1=a\int\rho(dD)/[(D-z-am_1)(D-z-am_2)]$, which contradicts Cauchy--Schwarz. Now let
$\ESD(D_N)\to\rho$. Since $|\bar m_N|\le\varepsilon^{-1}$, every subsequence has a further
subsequence along which $\bar m_N\to m^*$ with $\operatorname{Im}m^*\ge0$. In the average
in~\eqref{eq:matrix-approx-dyson} the integrand $D\mapsto(D-z-a\bar m_N)^{-1}$ is bounded by
$\varepsilon^{-1}$ and converges uniformly in $D$ to $(D-z-am^*)^{-1}$, which is bounded and
continuous; so $m^*=\int\rho(dD)/(D-z-am^*)$. The imaginary part of the integrand is positive, so
$\operatorname{Im}m^*>0$, and $m^*$ solves~\eqref{eq:matrix-dyson}. By uniqueness the whole
sequence $\bar m_N$ converges to this solution, and~\eqref{eq:matrix-resolvent-var} gives
convergence of $m_N$ in probability.
\end{proof}

The next lemma is used only in Appendix~\ref{app:kr}, at profiles that are strictly Plefka stable
in the sense $1-aN^{-1}\sum_iU_i^2>0$. It shows that the matrix $D_N-\sqrt a\,\GOE$ then has a
spectral gap, and that its inverse acts on fixed vectors like the diagonal matrix of the $U_i$.

\begin{lemma}[Spectral gap and inverse at a strictly stable profile]\label{lem:matrix-gap}
Let $a>0$, let $0<U_i\le1$ for $i\le N$, put $\chi_N=N^{-1}\sum_iU_i$ and
$D_N=\diag(U_i^{-1}+a\chi_N)$.
\begin{thmlist}
\item\label{lem:matrix-gap-a} If $1-aN^{-1}\sum_iU_i^2\ge d>0$, there are $g_0,c_M>0$ depending
 only on $a$ and $d$ such that
 \begin{equation}\label{eq:matrix-gap}
  \Prob\bigl(D_N-\sqrt a\,\GOE\not\succeq g_0I\bigr)\le e^{-c_MN}.
 \end{equation}
 The same holds, with $(g_0,c_M)$ replaced by $(g_0/2,c_M/4)$, after $U_i^{-1}$ is replaced by
 $\min(U_i^{-1},K)$, for every fixed $K$ larger than a constant depending only on $a$ and $d$.
\item\label{lem:matrix-gap-b} Suppose that the empirical laws of the $U_i$ converge weakly to the
 law of a $(0,1]$-valued variable $U$ with $\chi=\E U$ and $1-a\E U^2>0$. Then, for deterministic
 vectors $v_N,w_N$ of bounded Euclidean norm,
 \begin{equation}\label{eq:matrix-inverse}
  v_N^{\mathsf T}\bigl(D_N-\sqrt a\,\GOE\bigr)^{-1}w_N-\sum_{i\le N}U_iv_{N,i}w_{N,i}\longrightarrow0,
  \qquad
  \frac1N\tr\bigl(D_N-\sqrt a\,\GOE\bigr)^{-1}\longrightarrow\chi
 \end{equation}
 in probability, where the left sides are set to zero outside the event of~\ref{lem:matrix-gap-a}.
 The same holds for random vectors of bounded norm independent of $\GOE$.
\item\label{lem:matrix-gap-c} \textup{(General profiles.)} Let $D_N$ be deterministic real diagonal
 matrices with $\ESD(D_N)\to\rho'$ weakly, and let $m_0\in\R$ satisfy
 \[
  m_0=\int\frac{\rho'(dD)}{D-am_0},\qquad a\varsigma<1,\qquad
  \varsigma:=\int\frac{\rho'(dD)}{(D-am_0)^2},
 \]
 where $D_{N,ii}-am_0\ge c_\circ>0$ for all $i$ and all large $N$. If
 $\Prob(D_N-\sqrt a\,\GOE\succeq g'I)\to1$ for some $g'>0$, then, for $v_N,w_N$ as
 in~\ref{lem:matrix-gap-b},
 \begin{equation}\label{eq:matrix-variant}
 \begin{gathered}
  v_N^{\mathsf T}\bigl(D_N-\sqrt a\,\GOE\bigr)^{-1}w_N-\sum_{i\le N}\frac{v_{N,i}w_{N,i}}{D_{N,ii}-am_0}
  \longrightarrow0,\\
  \frac1N\tr\bigl(D_N-\sqrt a\,\GOE\bigr)^{-1}\longrightarrow m_0
 \end{gathered}
 \end{equation}
 in probability, with the left sides set to zero outside the event
 $\{D_N-\sqrt a\,\GOE\succeq g'I\}$.
\item\label{lem:matrix-gap-d} \textup{(Capped profiles.)} Let $U$ be as in~\ref{lem:matrix-gap-b},
 let $\rho_K$ be the law of $\min(U^{-1},K)+a\chi$ and $V=\max(U,K^{-1})$. For a small fixed
 $\delta>0$ and all large $K$, the map $m\mapsto\int\rho_K(dD)/(D-am)=\E[V/(1-a(m-\chi)V)]$ is a
 contraction of $[\chi-\delta,\chi+\delta]$; its fixed point $m_0=\chi+O(K^{-1})$ satisfies the
 hypotheses of~\ref{lem:matrix-gap-c} with $\rho'=\rho_K$ and $c_\circ=\frac12$, for every profile
 $D_{N,ii}=\min(u_i^{-1},K)+a\chi_N'$ with $0<u_i\le1$, $\chi_N'\to\chi$ and $\ESD(D_N)\to\rho_K$,
 except for the gap condition, which must be checked separately.
\end{thmlist}
\end{lemma}

\begin{proof}
\ref{lem:matrix-gap-a} Put $t=a\chi_N+s$ with $s=\min\{\frac12,d/(4a)\}$. Then
$D_{N,ii}-t=U_i^{-1}-s>0$, and the right side of~\eqref{eq:matrix-comparison} equals
\[
 -s+\frac{as}N\sum_i\frac{U_i^2}{1-sU_i}\le-ds+2as^2\le-\frac{ds}2 ,
\]
since $U_i^2/(1-sU_i)\le U_i^2+2s$ for $s\le\frac12$. The largest eigenvalue of
$\sqrt a\,\GOE-D_N$ is $\sqrt{2a/N}$-Lipschitz in the standard Gaussian variables, so Gaussian
concentration at distance $ds/4$ gives~\eqref{eq:matrix-gap} with $g_0=ds/4$ and
$c_M=d^2s^2/(64a)$. For the cap put $V_i=\max(U_i,K^{-1})$; at the same $t$ the right side
of~\eqref{eq:matrix-comparison} becomes
$\frac aN\sum_i(V_i-U_i)-s+\frac{as}N\sum_iV_i^2/(1-sV_i)$, which differs from the uncapped value
by at most $(a+10as)/K$, uniformly in the profile, because $y\mapsto y^2/(1-sy)$ is
$10$-Lipschitz on $[0,1]$. For large $K$ the bound $-ds/4$ survives, and concentration at distance
$ds/8$ gives the capped statement.

\ref{lem:matrix-gap-b} Since $U$ has no atom at $0$, the empirical laws of $U_i^{-1}$ converge to
that of $U^{-1}$, and $\chi_N\to\chi$; so $\ESD(D_N)\to\rho:=\Law(U^{-1}+a\chi)$. Fix
$\varepsilon>0$ and put $z=\mathrm i\varepsilon$. By
Lemma~\ref{lem:matrix-comparison}\ref{lem:matrix-comparison-b}, $\E v^{\mathsf T}Gw$ equals
$\sum_iv_iw_i(D_{N,ii}-z-a\bar m_N)^{-1}+O(|v||w|/N)$, the denominators have modulus at least
$\varepsilon$, and $\bar m_N\to m(z)$ by part~\ref{lem:matrix-comparison-c}. Together with the
variance bound,
\begin{equation}\label{eq:matrix-resolvent-limit}
 v_N^{\mathsf T}G_N(z)w_N-\sum_i\frac{v_{N,i}w_{N,i}}{D_{N,ii}-z-am(z)}\longrightarrow0
\end{equation}
in probability; no delocalization is needed, since $\sum_i|v_iw_i|\le|v||w|$. We identify $m(z)$
for small $\varepsilon$. At $z=0$ the equation $m=\int\rho(dD)/(D-z-am)$, without the sign
condition, has the solution $m=\chi$, because $D-a\chi=U^{-1}$; the derivative of its right side in
$m$ there is $a\E U^2<1$, and all denominators are at least $1$ near $(0,\chi)$. The analytic
implicit function theorem gives a local solution $m_{\mathrm{loc}}(z)$ with
$m_{\mathrm{loc}}(0)=\chi$ and $m_{\mathrm{loc}}'(0)=\E U^2/(1-a\E U^2)>0$. Hence
$\operatorname{Im}m_{\mathrm{loc}}(\mathrm i\varepsilon)>0$ for small $\varepsilon$, and uniqueness
in~\eqref{eq:matrix-dyson} gives $m(\mathrm i\varepsilon)=m_{\mathrm{loc}}(\mathrm i\varepsilon)\to\chi$.
On the event of~\ref{lem:matrix-gap-a}, which has probability tending to one because
$N^{-1}\sum_iU_i^2\to\E U^2$,
$\|G_N(\mathrm i\varepsilon)-(D_N-\sqrt a\,\GOE)^{-1}\|\le\varepsilon/g_0^2$. Finally
$D_{N,ii}-\mathrm i\varepsilon-am(\mathrm i\varepsilon)=U_i^{-1}+\delta_N(\varepsilon)$ with
$\delta_N(\varepsilon)=a(\chi_N-\chi)+a\chi-\mathrm i\varepsilon-am(\mathrm i\varepsilon)$
independent of $i$, and $|(U_i^{-1}+\delta)^{-1}-U_i|\le2|\delta|$ when $|\delta|\le\frac12$,
because $U_i^{-1}\ge1$. Letting first $N\to\infty$ and then $\varepsilon\downarrow0$
in~\eqref{eq:matrix-resolvent-limit} proves the first claim. For the trace,
$m_N(\mathrm i\varepsilon)\to m(\mathrm i\varepsilon)$ in probability, the gap event gives
$|m_N(\mathrm i\varepsilon)-N^{-1}\tr(D_N-\sqrt a\,\GOE)^{-1}|\le\varepsilon/g_0^2$, and
$m(\mathrm i\varepsilon)\to\chi$. Random vectors independent of $\GOE$ are handled by
conditioning on them.

\ref{lem:matrix-gap-c} The proof of~\ref{lem:matrix-gap-b} applies with $(\rho',m_0)$ in place of
$(\rho,\chi)$. Lemma~\ref{lem:matrix-comparison} uses only that the diagonal is real and that its
empirical laws converge weakly, with errors uniform in the profile. The implicit function theorem
at $(0,m_0)$ uses the displayed conditions, since $\supp\rho'\subset[am_0+c_\circ,\infty)$, and gives
$m'(0)=\varsigma/(1-a\varsigma)>0$, which identifies the branch. The lower bound $c_\circ$ makes
the reciprocals $(D_{N,ii}-\mathrm i\varepsilon-am(\mathrm i\varepsilon))^{-1}$ converge to
$(D_{N,ii}-am_0)^{-1}$ uniformly in $i$ as $\varepsilon\downarrow0$. Only empirical laws and
weighted averages of reciprocals enter, not operator-norm convergence of the diagonals, so sparse
changes of coordinates are allowed.

\ref{lem:matrix-gap-d} Since $D-am=V^{-1}(1-a(m-\chi)V)$ for $D=V^{-1}+a\chi$, the formula for the
map holds. Its value at $m=\chi$ is $\E V=\chi+\E(K^{-1}-U)_+=\chi+O(K^{-1})$, and on
$|m-\chi|\le\delta$ its derivative is at most $a\E V^2/(1-a\delta)^2$, where
$a\E V^2\le a\E U^2+aK^{-2}<1$ for large $K$. For small $\delta$ and large $K$ the map is therefore
a contraction of $[\chi-\delta,\chi+\delta]$ into itself, and its fixed point $m_0$ satisfies
$|m_0-\chi|=O(K^{-1})$ and $a\varsigma<1$. The denominators satisfy
$D_{N,ii}-am_0\ge1-a|\chi_N'-\chi|-a|m_0-\chi|\ge\frac12$ for large $N$ and $K$.
\end{proof}

We turn to the Hessian of $\HN$. Recall that $\E\HN(x)\HN(y)=N\xi(\ov xy)$ and that, under
\condR, $\xi$ and all its derivatives are finite on $[-1,1]$, with $|\xi^{(k)}(t)|\le\xi^{(k)}(1)$
there. In particular $\HN$ is almost surely smooth on a neighbourhood of the cube. For
$y\in\cube$ we write $\|D^3\HN(y)\|$ for the supremum of $|D^3\HN(y)[\phi_1,\phi_2,\phi_3]|$ over
unit vectors $\phi_1,\phi_2,\phi_3$.

\begin{lemma}[The Hessian of the Gaussian field]\label{lem:matrix-third}
Assume \condR\ and put $C_H=2\xi''(1)+4\xi'''(1)+\xi''''(1)$.
\begin{thmlist}
\item\label{lem:matrix-third-a} For $y\in\cube$ with $r=q_y$ and every symmetric $N\times N$
 matrix $S$,
 \begin{equation}\label{eq:matrix-hess-cov}
  \Var\tr\bigl(S\,\Hess\HN(y)\bigr)=\frac{2\xi''(r)}N\tr S^2+\frac{4\xi'''(r)}{N^2}|Sy|^2
  +\frac{\xi''''(r)}{N^3}\bigl(y^{\mathsf T}Sy\bigr)^2\le\frac{C_H}N\tr S^2 .
 \end{equation}
 Consequently, at a fixed $y$,
 \[
  \Hess\HN(y)\overset{d}{=}\sqrt{\xi''(r)}\,\GOE+\frac1N\bigl(y\gamma^{\mathsf T}+\gamma y^{\mathsf T}\bigr)
  +\frac{Z}{N^{3/2}}\,yy^{\mathsf T},
 \]
 with $\GOE$, $\gamma\sim N(0,\xi'''(r)I_N)$ and $Z\sim N(0,\xi''''(r))$ independent. The last two
 terms have rank at most two.
\item\label{lem:matrix-third-b} For every $M>0$ there are $C^H_2,C^H_3$, depending only on $\xi$ and
 $M$, such that for every $N$
 \begin{equation}\label{eq:matrix-sup}
 \begin{gathered}
  \Prob\Bigl(\sup_{y\in\cube}\|\Hess\HN(y)\|>C^H_2\Bigr)\le e^{-MN},\\
  \Prob\Bigl(\sup_{y\in\cube}\|D^3\HN(y)\|>\frac{C^H_3}{\sqrt N}\Bigr)\le e^{-MN}.
 \end{gathered}
 \end{equation}
\end{thmlist}
\end{lemma}

\begin{proof}
\ref{lem:matrix-third-a} Differentiating $N\xi(\ov xy)$ twice in $x$ and twice in $y$ and setting
$x=y$ gives
\[
 \E\,\partial_{ij}\HN(y)\,\partial_{kl}\HN(y)
 =\frac{\xi''}N(\delta_{ik}\delta_{jl}+\delta_{il}\delta_{jk})
 +\frac{\xi'''}{N^2}(y_iy_k\delta_{jl}+y_iy_l\delta_{jk}+y_jy_k\delta_{il}+y_jy_l\delta_{ik})
 +\frac{\xi''''}{N^3}y_iy_jy_ky_l,
\]
with the derivatives of $\xi$ at $r$. The three terms are the covariances of the three independent
summands in the displayed decomposition, which is therefore an identity in law between centred
Gaussian matrices. Contracting with $S$ gives~\eqref{eq:matrix-hess-cov}; the bound uses
$|Sy|^2\le\|S\|^2|y|^2\le N\tr S^2$, $(y^{\mathsf T}Sy)^2\le N^2\tr S^2$ and $r\le1$.

\ref{lem:matrix-third-b} For $k\in\{2,3\}$ let $T_k(y,\phi)=D^k\HN(y)[\phi_1,\dots,\phi_k]$ on the
index set $\mathcal I_k=\cube\times(S^{N-1})^k$. By the series
$N\xi(\ov xy)=N\sum_p\beta_p^2(x\cdot y/N)^p$, the $j$-th derivative in $x$ along unit vectors
$\phi_1,\dots,\phi_j$ is $N^{1-j}\xi^{(j)}(\ov xy)\prod_i(y\cdot\phi_i)$, and each further
derivative in $y$ along a unit vector either differentiates $\xi^{(j)}(\ov xy)$, producing a factor
$x\cdot\psi/N$, or replaces a factor $y\cdot\phi_i$ by $\phi_i\cdot\psi$. Since
$|y\cdot\phi|\le\sqrt N$ and $|x\cdot\psi|/N\le N^{-1/2}$, every mixed derivative of order $j$ in
$x$ and $l$ in $y$, contracted with unit vectors, is bounded by $CN^{1-(j+l)/2}$, with $C$
depending on $\xi^{(i)}(1)$ for $i\le j+l$. Hence $\Var T_k(y,\phi)\le CN^{1-k}$. For the canonical
distance $d$ of $T_k$, the mean value theorem in $y$ (using the bound for $D^{k+1}\HN$) and
multilinearity in $\phi$ give
\[
 d\bigl((y,\phi),(y',\phi')\bigr)\le CN^{(1-k)/2}\Bigl(\frac{|y-y'|}{\sqrt N}+\sum_{i\le k}|\phi_i-\phi_i'|\Bigr).
\]
The set $N^{-1/2}\cube$ lies in the unit ball, which, like the unit sphere, is covered by
$(3/s)^N$ balls of radius $s\le1$. For the distance in brackets, the product index set is
therefore covered by $(3(k+1)/s)^{(k+1)N}$ sets of radius $s\le k+1$, and Dudley's entropy bound
gives
\[
 \E\sup_{\mathcal I_k}T_k\le CN^{(1-k)/2}\int_0^{2(k+1)}\sqrt{(k+1)N\log(3(k+1)/s)}\,ds\le C'N^{1-k/2},
\]
and the Borell--TIS inequality, with variance $\sigma^2\le CN^{1-k}$, gives
$\Prob(\sup T_k>\E\sup T_k+u)\le e^{-u^2/(2\sigma^2)}$. Taking $u$ a large constant for $k=2$ and a
large multiple of $N^{-1/2}$ for $k=3$ proves~\eqref{eq:matrix-sup}; the process is continuous, so
the suprema are measurable. Derivatives of $\xi$ up to order eight enter, through the canonical
distance of $T_3$.
\end{proof}

\begin{remark}[Conditioned fields]\label{rem:matrix-conditional}
Appendix~\ref{app:kr} uses the third-derivative bound for $\HN$ conditioned on its value and its
gradient at a point. Gaussian conditioning decreases variances and increment variances, so the
centred part of the conditioned field satisfies Lemma~\ref{lem:matrix-third}\ref{lem:matrix-third-b}
with the same constants. Its conditional mean is a separate matter, which Appendix~\ref{app:kr}
treats under its hypothesis of two active degrees. Theorem~\ref{thm:struct} uses only the
unconditional bounds of Lemma~\ref{lem:matrix-third}, which hold for every $\xi$ satisfying \condR,
single-degree covariances included.
\end{remark}

\subsection{The TAP function: correction, susceptibilities and Hessian}\label{app:matrix-hessian}

From now on $\xi\in\classC$. The TAP function splits into the Gaussian field and a deterministic
function of the self-overlap and the magnetization law. We study how this function and the
susceptibilities $u_\alpha$ of~\eqref{eq:hess-u} depend on their arguments, and compute the
derivatives of $F_\alpha$.

With $V_\alpha$ as in~\eqref{eq:hess-V}, put
\[
 V_{\mathrm{TAP}}(r,\mu)=\inf_{\beta\in\Pspace([0,1])}V_\beta(r,\mu),\qquad r\in[0,1],\ \mu\in\Pspace([-1,1]).
\]
Then for every $x\in\cube$
\begin{equation}\label{eq:matrix-split}
 F_\alpha(x)=\HN(x)+NV_\alpha(q_x,\mu_x),\qquad \FTAP(x)=\HN(x)+NV_{\mathrm{TAP}}(q_x,\mu_x).
\end{equation}
The first identity is the definition, since $\sum_i\Phistar_\alpha(q_x,x_i)=N\int\Phistar_\alpha(q_x,m)\,\mu_x(dm)$.
For the second, the function $\Phi_\beta(r,\cdot)$ is obtained by solving the Parisi equation
backward from time one, so $V_\beta(r,\mu)$ depends on $\beta$ only through $\cdf\beta$ on $[r,1]$;
hence $V_\beta(r,\mu)=V_{\Shift_r\beta}(r,\mu)$, and the infimum over $\Pspace([r,1])$ that defines
$\FTAP$ equals the infimum over $\Pspace([0,1])$.

\paragraph{Time derivatives.}
For $\alpha\in\Pspace([0,1])$ and $t\in[0,1]$ put $\cdf\alpha^+(t)=\alpha([0,t])=\cdf\alpha(t)$
and $\cdf\alpha^{-}(t)=\alpha([0,t))$. For $j\le2$, $b\in\R$ and $t\in(0,1)$, the function
$t\mapsto\partial_x^j\Phi_\alpha(t,b)$ has derivatives from the right and from the left,
\begin{equation}\label{eq:matrix-onesided}
 \partial_t^\pm\partial_x^j\Phi_\alpha(t,b)
 =-\frac{\xi''(t)}2\,\partial_x^j\bigl(\Phi_{\alpha,xx}+\cdf\alpha^{\pm}(t)\Phi_{\alpha,x}^2\bigr)(t,b).
\end{equation}
Indeed, for finitely supported $\alpha$ the integrated form
$\partial_x^j\Phi_\alpha(t,b)-\partial_x^j\Phi_\alpha(s,b)
=-\int_s^t\frac{\xi''}2\partial_x^j(\Phi_{\alpha,xx}+\cdf\alpha\,\Phi_{\alpha,x}^2)(\tau,b)\,d\tau$ holds
on each step of $\cdf\alpha$, as in the proof of Lemma~\ref{lem:pde-standing}\ref{lem:pde-standing-a},
and across steps by continuity. For general $\alpha$ it follows by approximation, since
$\partial_x^i\Phi_{\alpha_n}\to\partial_x^i\Phi_\alpha$ uniformly for every $i$ when
$\cdf{\alpha_n}\to\cdf\alpha$ in $L^1$ \cite[Proposition~1(i)]{AC15a}. The integrand is jointly
continuous apart from the factor $\cdf\alpha$, which is right-continuous with left limits; this
gives~\eqref{eq:matrix-onesided}. Since $0<\Phi_{\alpha,xx}\le1$ and $|\Phi_{\alpha,x}|\le1$,
in particular $|\partial_t^\pm\Phi_\alpha|\le\xi''(1)$.

For $|m|<1$ write $b=b_\alpha(t,m)$, and put
\[
 v_\alpha(t,m)=\frac{\Phi_{\alpha,xxx}}{\Phi_{\alpha,xx}}(t,b),\qquad
 T_\alpha(t,m)=\Bigl(\Phi_{\alpha,xxxx}-\frac{\Phi_{\alpha,xxx}^2}{\Phi_{\alpha,xx}}\Bigr)(t,b).
\]
By Lemma~\ref{lem:pde-standing}\ref{lem:pde-standing-a}, the bounds $|v_\alpha|\le C_{34}$ and
$|T_\alpha|\le C_{34}(1+C_{34})$ hold. We differentiate the identity
$\Phi_{\alpha,x}(t,b_\alpha(t,m))=m$ in $m$ and in $t$, using~\eqref{eq:matrix-onesided}; this gives
the inverse-field identities
\begin{equation}\label{eq:matrix-inverse-field}
\begin{gathered}
 \partial_mb_\alpha=\frac1{u_\alpha},\qquad
 \partial_t^\pm b_\alpha=\frac{\xi''(t)}2\bigl(v_\alpha+2\cdf\alpha^{\pm}(t)\,m\bigr),\\
 \partial_mu_\alpha=v_\alpha,\qquad
 \partial_t^\pm u_\alpha=-\frac{\xi''(t)}2T_\alpha-\xi''(t)\cdf\alpha^{\pm}(t)\,u_\alpha^2 .
\end{gathered}
\end{equation}
For the last one, $\partial_t^\pm u_\alpha=\partial_t^\pm\Phi_{\alpha,xx}+\Phi_{\alpha,xxx}\partial_t^\pm b_\alpha$,
and the terms $\xi''\cdf\alpha^{\pm}m\Phi_{\alpha,xxx}$ cancel.

\paragraph{The TAP correction.}
Put $\varrho(t)=t\log(2e/t)$ for $t\in(0,2]$ and $\varrho(0)=0$; it is concave and nondecreasing
on $[0,2]$, and $\varrho(t)\ge t$ there. Put
\[
 C_*=\log2+\xi'(1),\qquad L_1=\tfrac32\xi''(1),\qquad L_2=\tfrac12\bigl(1+\log(\overline C/2)\bigr).
\]

\begin{lemma}[The TAP correction]\label{lem:matrix-correction}
Assume $\xi\in\classC$. For all $\alpha,\beta\in\Pspace([0,1])$, $r,r'\in[0,1]$ and
$\mu,\mu'\in\Pspace([-1,1])$:
\begin{thmlist}
\item\label{lem:matrix-correction-a} $|\Phistar_\alpha(r,m)|\le C_*$ for every $m\in[-1,1]$;
\item\label{lem:matrix-correction-b}
 $|V_\alpha(r,\mu)-V_\alpha(r',\mu')|\le L_1|r-r'|+L_2\,\varrho\bigl(\Wone(\mu,\mu')\bigr)$, and the same
 holds for $V_{\mathrm{TAP}}$;
\item\label{lem:matrix-correction-c}
 $|V_\alpha(r,\mu)-V_\beta(r,\mu)|\le\xi''(1)\int_r^1|\cdf\alpha(t)-\cdf\beta(t)|\,dt$.
\end{thmlist}
\end{lemma}

\begin{proof}
\ref{lem:matrix-correction-a} Let $\Phi_0$ be the solution on $[r,1]$ with coefficient zero and
terminal value $\log\cosh$. By the coefficient sandwich of Lemma~\ref{lem:pde-coef}, which extends
to measurable coefficients by Lemma~\ref{lem:pde-measurable}, and since the second moments
appearing there lie in $[0,1]$, $\Phi_\alpha(r,b)$ lies between $\Phi_0(r,b)$ and
$\Phi_0(r,b)+\frac12(\xi'(1)-\xi'(r))$. Here $\Phi_0(r,b)=\E\log\cosh(b+\sigma Z)$ with $Z$ standard Gaussian and
$\sigma^2=\xi'(1)-\xi'(r)$, so Jensen's inequality gives $\Phi_\alpha(r,b)\ge\log\cosh b\ge|b|-\log2$,
and $\log\cosh y\le y^2/2$ gives $\Phi_\alpha(r,0)\le\xi'(1)$. Hence
$-\log2\le\inf_b\{\Phi_\alpha(r,b)-mb\}\le\Phi_\alpha(r,0)\le\xi'(1)$.

\ref{lem:matrix-correction-b} By~\eqref{eq:matrix-onesided},
$\sup_b|\Phi_\alpha(r,b)-\Phi_\alpha(r',b)|\le\xi''(1)|r-r'|$, and an infimum preserves this bound;
the integral term of $V_\alpha$ is $\frac12\xi''(1)$-Lipschitz in $r$. For $|m|<1$ the infimum
defining $\Phistar_\alpha(r,m)$ is attained at $b_\alpha(r,m)$, and $\partial_m\Phistar_\alpha=-b_\alpha$.
By Lemma~\ref{lem:pde-standing}\ref{lem:pde-standing-b},
$|b_\alpha(r,m)|\le\frac12\log(\overline C/2)+\frac12\log\frac1{1-|m|}$. The function
$s\mapsto\log\frac1{1-|s|}$ on $[-1,1]$ has decreasing rearrangement $y\mapsto\log(2/y)$ on $[0,2]$,
so its integral over any interval of length $h$ is at most $\int_0^h\log(2/y)\,dy=\varrho(h)$.
Therefore, using $h\le\varrho(h)$ and continuity at $m=\pm1$,
$|\Phistar_\alpha(r,m)-\Phistar_\alpha(r,m')|\le L_2\varrho(|m-m'|)$. Integrating over an optimal
coupling of $\mu$ and $\mu'$ and using the concavity of $\varrho$ gives the bound in $\mu$. An
infimum over $\alpha$ of functions with a common modulus has the same modulus.

\ref{lem:matrix-correction-c} By Lemma~\ref{lem:pde-coef} on $[r,1]$, extended to measurable
coefficients by Lemma~\ref{lem:pde-measurable},
$\sup_b|\Phi_\alpha(r,b)-\Phi_\beta(r,b)|\le\frac12\xi''(1)\int_r^1|\cdf\alpha-\cdf\beta|$. The
conjugates inherit this bound, and the integral terms of $V_\alpha$ and $V_\beta$ differ by at most
the same amount.
\end{proof}

\paragraph{Susceptibilities.}
Up to a common shift, the diagonal entries of the Hessian are the reciprocals
$u_\alpha(q_x,x_i)^{-1}$ of the susceptibilities.
Proposition~\ref{prop:hess-bulk} replaces the order parameter $\alpha$ of a point by $\zf$, and
for this it needs $u_\alpha$ to depend continuously on $\alpha$, uniformly in the magnetization.
We use the constants $E_0$ and $L_v$ of Lemma~\ref{lem:pde-stability} and~\eqref{eq:pde-constants}
with $G=1$, that is, $E_0=e^{\xi'(1)}$ and $L_v=E_0\,\xi''(1)\,e^{E_0\xi'(1)}$.

\begin{lemma}[Stability in the order parameter]\label{lem:matrix-orderstab}
Assume $\xi\in\classC$. For $\alpha,\beta\in\Pspace([0,1])$ and $t\in[0,1]$ put
$\epsilon=\int_t^1|\cdf\alpha(s)-\cdf\beta(s)|\,ds$. Then
\[
 \sup_{b\in\R}\bigl|\Phi_{\alpha,x}(t,b)-\Phi_{\beta,x}(t,b)\bigr|\le L_v\epsilon,\qquad
 \sup_{b\in\R}\bigl|\Phi_{\alpha,xx}(t,b)-\Phi_{\beta,xx}(t,b)\bigr|\le2\bigl(C_{34}L_v\epsilon\bigr)^{1/2},
\]
and there is a nondecreasing function $\omega_S$ on $[0,1]$ with $\omega_S(\epsilon)\to0$ as
$\epsilon\downarrow0$, depending only on $\xi'(1)$ and $\xi''(1)$, such that
\[
 \sup_{m\in[-1,1]}\bigl|u_\alpha(t,m)-u_\beta(t,m)\bigr|\le\omega_S(\epsilon).
\]
\end{lemma}

\begin{proof}
On $[t,1]$ the functions $\Phi_{\alpha,x}$ and $\Phi_{\beta,x}$ are the space derivatives of the
solutions with terminal value $\log\cosh$ and coefficients $\cdf\alpha$ and $\cdf\beta$ restricted
to $[t,1)$, with values in $[0,1]$. Lemma~\ref{lem:pde-stability} with $G=1$ and equal terminal
values, extended to measurable coefficients by Lemma~\ref{lem:pde-measurable}\ref{lem:pde-measurable-i},
gives the first bound. For the second, $w=\Phi_{\alpha,x}(t,\cdot)-\Phi_{\beta,x}(t,\cdot)$ satisfies
$\|w\|_\infty\le L_v\epsilon$ and $\|w''\|_\infty\le2C_{34}$, because
$|\Phi_{\alpha,xxx}|\le C_{34}\Phi_{\alpha,xx}\le C_{34}$ by
Lemma~\ref{lem:pde-standing}\ref{lem:pde-standing-a}. Taylor's formula at $b\pm h$ gives
$|w'(b)|\le\|w\|_\infty/h+h\|w''\|_\infty/2$, and the choice $h=(2\|w\|_\infty/\|w''\|_\infty)^{1/2}$
gives $\|w'\|_\infty\le(2\|w\|_\infty\|w''\|_\infty)^{1/2}$.

For the susceptibilities fix $\tau\in(0,1)$. If $|m|\le1-\tau$, then by
Lemma~\ref{lem:pde-standing}\ref{lem:pde-standing-b} both $b_\alpha=b_\alpha(t,m)$ and
$b_\beta=b_\beta(t,m)$ lie in $[-B_\tau,B_\tau]$, $B_\tau=\frac12\log(\overline C/2)+\frac12\log\frac1\tau$,
where $\Phi_{\alpha,xx}\ge c_\tau:=\underline c\,e^{-2B_\tau}$. Since
$\Phi_{\alpha,x}(t,b_\alpha)=m=\Phi_{\beta,x}(t,b_\beta)$, the mean value theorem gives
$c_\tau|b_\alpha-b_\beta|\le|\Phi_{\alpha,x}(t,b_\beta)-\Phi_{\beta,x}(t,b_\beta)|\le L_v\epsilon$.
Hence
\[
 |u_\alpha(t,m)-u_\beta(t,m)|\le C_{34}|b_\alpha-b_\beta|+\bigl|\Phi_{\alpha,xx}(t,b_\beta)-\Phi_{\beta,xx}(t,b_\beta)\bigr|
 \le\frac{C_{34}L_v\epsilon}{c_\tau}+2\bigl(C_{34}L_v\epsilon\bigr)^{1/2}.
\]
If $|m|>1-\tau$, Lemma~\ref{lem:pde-standing}\ref{lem:pde-standing-b} gives
$0\le u_\alpha,u_\beta\le(2\overline C/\underline c)\tau$. The infimum over $\tau\in(0,1)$ of the
sum of the two bounds is the required $\omega_S(\epsilon)$.
\end{proof}

\paragraph{Gradient and Hessian.}
We now compute the derivatives of $F_\alpha$. The self-overlap enters $F_\alpha$ in two places, and
in the gradient the two resulting terms cancel; in the Hessian the factor $\cdf\alpha^{\pm}$
survives, which is why an atom of $\alpha$ at $q_x$ produces two one-sided Hessians.

\begin{lemma}[Gradient and one-sided Hessians of the TAP function]\label{lem:matrix-hessian}
Assume $\xi\in\classC$. Let $\alpha\in\Pspace([0,1])$ and $x\in\ocube\setminus\{0\}$, and put
$r=q_x\in(0,1)$, $a_r=\xi''(r)$, $u_i=u_\alpha(r,x_i)$, $b_i=b_\alpha(r,x_i)$, $v_i=v_\alpha(r,x_i)$,
$T_i=T_\alpha(r,x_i)$, $\chi_\alpha(x)=N^{-1}\sum_iu_i$ and
$\Plef_N^\alpha(x)=1-a_rN^{-1}\sum_iu_i^2$.
\begin{thmlist}
\item\label{lem:matrix-hessian-a} $F_\alpha$ is continuously differentiable on
 $\ocube\setminus\{0\}$, and
 \begin{equation}\label{eq:matrix-gradient}
  \nabla F_\alpha(x)=\nabla\HN(x)-(b_i)_{i\le N}-a_r\chi_\alpha(x)\,x .
 \end{equation}
\item\label{lem:matrix-hessian-b} For $\sigma\in\{+,-\}$ let
 \begin{equation}\label{eq:matrix-ambient}
  \nabla^2_\sigma F_\alpha(x)=\Hess\HN(x)-\diag\bigl(u_1^{-1},\dots,u_N^{-1}\bigr)-a_r\chi_\alpha(x)\,I
  +R_\sigma(x),
 \end{equation}
 \[
  R_\sigma(x)=-\frac{a_r}N\bigl(xv^{\mathsf T}+vx^{\mathsf T}\bigr)
  -\frac2N\Bigl\{\xi'''(r)\chi_\alpha(x)-\frac{a_r^2}{2N}\sum_{i\le N}T_i
  +a_r\cdf\alpha^{\sigma}(r)\,\Plef_N^\alpha(x)\Bigr\}xx^{\mathsf T}.
 \]
 Then, for every $e\in\R^N$,
 $\lim_{s\downarrow0}s^{-1}(\nabla F_\alpha(x+se)-\nabla F_\alpha(x))$ equals
 $\nabla^2_+F_\alpha(x)e$ if $x\cdot e\ge0$ and $\nabla^2_-F_\alpha(x)e$ if $x\cdot e\le0$. The
 matrix $R_\sigma(x)$ is symmetric, of rank at most two, with range in
 $\operatorname{span}\{x,v\}$, and
 $R_+(x)-R_-(x)=-\frac{2a_r}N\,\alpha(\{r\})\,\Plef_N^\alpha(x)\,xx^{\mathsf T}$. In particular
 $\nabla^2_+F_\alpha(x)=\nabla^2_-F_\alpha(x)$ when $\alpha(\{r\})=0$, and the compressions of the
 two one-sided Hessians to $x^\perp$ always coincide.
\item\label{lem:matrix-hessian-c} If $\nabla F_\alpha(x)=0$, this common compression is the
 Riemannian Hessian at $x$ of the restriction of $F_\alpha$ to the sphere $\{y:q_y=r\}$.
\end{thmlist}
\end{lemma}

\begin{proof}
\ref{lem:matrix-hessian-a} Write $r_y=q_y$. In $F_\alpha(y)=\HN(y)+\sum_i\Phistar_\alpha(r_y,y_i)
-\frac N2\int_{r_y}^1t\xi''(t)\cdf\alpha(t)\,dt$, the derivative in $y_j$ of the middle term at fixed
$r_y$ is $-b_\alpha(r_y,y_j)$. The infimum defining $\Phistar_\alpha$ is attained at a unique point
that depends continuously on $r$, so~\eqref{eq:matrix-onesided} gives the one-sided derivatives
$\partial_r^\pm\Phistar_\alpha(r,m)=-\frac{a_r}2(u_\alpha+\cdf\alpha^{\pm}(r)m^2)$. With
$\partial_{y_j}r_y=2y_j/N$, the terms produced by the dependence on $r_y$ are
\[
 \frac{2y_j}N\Bigl[-\frac{a_r}2\sum_i\bigl(u_i+\cdf\alpha^{\pm}(r)y_i^2\bigr)
 +\frac N2\,r\,a_r\,\cdf\alpha^{\pm}(r)\Bigr]=-a_r\chi_\alpha(y)\,y_j ,
\]
because $\sum_iy_i^2=Nr$: the terms containing $\cdf\alpha^{\pm}$ cancel. So both one-sided
partial derivatives exist and are given by~\eqref{eq:matrix-gradient}, which is continuous in $y$;
hence $F_\alpha$ is continuously differentiable.

\ref{lem:matrix-hessian-b} Let $e\in\R^N$ and $y=x+se$. Then
$r_y-r=2s(x\cdot e)/N+s^2|e|^2/N$, so for small $s>0$ the self-overlap moves to the right of $r$ if
$x\cdot e>0$, to the left if $x\cdot e<0$, and by $O(s^2)$ if $x\cdot e=0$. By the
identities~\eqref{eq:matrix-inverse-field}, the joint continuity of $u_\alpha$, $v_\alpha$ and
$b_\alpha$ (Lemma~\ref{lem:pde-standing}\ref{lem:pde-standing-b}), and the mean value theorem in
the variable $m$,
\[
 b_\alpha(r_y,y_j)-b_\alpha(r,x_j)=\frac{se_j}{u_j}+\frac{a_r}2\bigl(v_j+2\cdf\alpha^{\sigma}(r)x_j\bigr)(r_y-r)+o(s),
\]
with $\sigma$ the sign of $x\cdot e$, and similarly for $u_\alpha(r_y,y_i)$. Differentiating the $j$-th component
$\partial_j\HN(y)-b_\alpha(r_y,y_j)-\xi''(r_y)\chi_\alpha(y)y_j$ of~\eqref{eq:matrix-gradient}
along $e$ with these expansions gives $\sum_k H_{jk}e_k$, where the contributions to $H_{jk}$ of the
last two terms are
\begin{align*}
 \partial_k\bigl[b_\alpha(r_y,y_j)\bigr]&=\frac{\delta_{jk}}{u_j}+\frac{a_r}N\bigl(v_jx_k+2\cdf\alpha^{\sigma} x_jx_k\bigr),\\
 \partial_k\bigl[\xi''(r_y)\chi_\alpha(y)y_j\bigr]&=a_r\chi_\alpha\delta_{jk}+\frac{a_r}Nx_jv_k
 +\frac{2\xi'''(r)\chi_\alpha}Nx_jx_k\\
 &\qquad-\frac{a_r^2}{N^2}\Bigl(\sum_iT_i\Bigr)x_jx_k
 -\frac{2a_r^2\cdf\alpha^{\sigma}}{N^2}\Bigl(\sum_iu_i^2\Bigr)x_jx_k ,
\end{align*}
where $\partial_k\chi_\alpha=N^{-1}[v_k+(2x_k/N)\sum_i\partial_r^\sigma u_\alpha(r,x_i)]$ was used.
Collecting the multiples of $x_jx_k$ gives~\eqref{eq:matrix-ambient}. When $x\cdot e=0$ the terms
containing $r_y-r$ are $o(s)$, and $R_+(x)e=R_-(x)e=-\frac{a_r}Nx(v\cdot e)$. The remaining claims
are read off from the formula for $R_\sigma$, since $\cdf\alpha^{+}(r)-\cdf\alpha^{-}(r)=\alpha(\{r\})$
and $Px=0$ for the orthogonal projection $P$ onto $x^\perp$.

\ref{lem:matrix-hessian-c} On the sphere $\{q_y=r\}$ the function $F_\alpha$ coincides with the
smooth function $\tilde F(y)=\HN(y)+\sum_i\Phistar_\alpha(r,y_i)-c$, for a constant $c$. The
Riemannian Hessian of a function on a sphere centred at the origin is
$P(\Hess\tilde F(x)-(x\cdot\nabla\tilde F(x)/|x|^2)I)P$. Here
$\Hess\tilde F=\Hess\HN-\diag(u_i^{-1})$ and $\nabla\tilde F=\nabla\HN-(b_i)_i$, and
$\nabla F_\alpha(x)=0$ gives $x\cdot\nabla\tilde F(x)=a_r\chi_\alpha(x)|x|^2$ by~\eqref{eq:matrix-gradient}.
Since $PR_\sigma(x)P=0$ on $x^\perp$, the Riemannian Hessian equals the compression of
$\nabla^2_\sigma F_\alpha(x)$.
\end{proof}

For the uniform estimates below we cap the diagonal. Coordinates with small susceptibility have
large diagonal entries; capping them changes only those coordinates.

\begin{corollary}[Capped diagonal]\label{cor:matrix-truncated}
In the setting of Lemma~\ref{lem:matrix-hessian}, for $K\ge1$ put
$D^\alpha_K(x)=\diag\bigl(\min(u_i^{-1},K)+a_r\chi_\alpha(x)\bigr)_{i\le N}$ and
$p^\alpha_K(x)=N^{-1}\#\{i:u_i<1/K\}$. Then $D^\alpha_K(x)\in\mathcal D_{K+\xi''(1)}$, and for
$\sigma\in\{+,-\}$
\[
 \dL\Bigl(\ESD\bigl(\nabla^2_\sigma F_\alpha(x)\bigr),\ESD\bigl(\Hess\HN(x)-D^\alpha_K(x)\bigr)\Bigr)
 \le p^\alpha_K(x)+\frac2N .
\]
\end{corollary}

\begin{proof}
The entries of $D^\alpha_K(x)$ lie in $[1,K+\xi''(1)]$ because $0<u_i\le1$ and $0\le\chi_\alpha\le1$.
The difference of the two matrices is $-\diag((u_i^{-1}-K)_+)+R_\sigma(x)$, of rank at most
$Np^\alpha_K(x)+2$; apply~\eqref{eq:matrix-rank}.
\end{proof}

\subsection{Reflected Pastur laws and the edge}\label{app:matrix-pastur}

We prove Lemma~\ref{lem:hess-pastur}, then the edge hypothesis~\eqref{eq:hess-fhyp} for the
variables $\UQ$, then parts~\ref{lem:hess-edge-c}--\ref{lem:hess-edge-e} of
Lemma~\ref{lem:hess-edge}, and finally a version of the square-root edge that is uniform under
bounded changes of measure, which Appendix~\ref{app:matrix-index} needs.

\begin{proof}[Proof of Lemma~\ref{lem:hess-pastur}]
\ref{lem:hess-pastur-a} \emph{Uniqueness.} Let $w\in\Cplx^+$ and let $g_1,g_2$ solve
\eqref{eq:hess-pastur} with $\operatorname{Im}g_j<0$. Taking imaginary parts,
$-\operatorname{Im}g_j=(\operatorname{Im}w+a|\operatorname{Im}g_j|)I_j$ with
$I_j=\int\rho(dD)|w+D-ag_j|^{-2}$, so $aI_j<1$. If $g_1\ne g_2$, subtracting the two equations
gives $1=a\int\rho(dD)/[(w+D-ag_1)(w+D-ag_2)]$, which contradicts Cauchy--Schwarz.

\emph{Existence and convergence.} Let $D_N$ be real diagonal matrices with $\ESD(D_N)\to\rho$, put
$X_N=D_N-\sqrt a\,\GOE$, and let $m_N(z)=N^{-1}\tr(X_N-z)^{-1}$ be the Stieltjes transform of
$\ESD(X_N)$. By Lemma~\ref{lem:matrix-comparison}\ref{lem:matrix-comparison-c}, $m_N(z)\to m(z)$ in
probability for every $z\in\Cplx^+$, where $m(z)$ is the solution of~\eqref{eq:matrix-dyson}. Realize
the matrices $\GOE$, $N\ge1$, on one probability space. By
Lemma~\ref{lem:matrix-comparison}\ref{lem:matrix-comparison-a} and the Borel--Cantelli lemma,
$\|\GOE\|\le3$ for all large $N$ almost surely, and then $\dL(\ESD(X_N),\ESD(D_N))\le3\sqrt a$
by~\eqref{eq:matrix-weyl}; so $(\ESD(X_N))_N$ is almost surely tight. Given a subsequence, a
diagonal argument yields a further subsequence along which, almost surely, $m_N(z_k)\to m(z_k)$ for
every $z_k$ in a countable dense subset of $\Cplx^+$. Along it, almost surely, every weak limit
point of $\ESD(X_N)$ is a probability measure whose Stieltjes transform agrees with $m$ on that set,
hence on $\Cplx^+$; since a Stieltjes transform determines its measure, all limit points coincide
with one probability measure $\mu$, which depends only on $(a,\rho)$. So $\ESD(X_N)\to\mu$ weakly
in probability, and the expected empirical spectral distributions converge to $\mu$ as well. Let
$\nu$ be the image of $\mu$ under $\lambda\mapsto-\lambda$, so that
$\ESD(\sqrt a\,\GOE-D_N)\to\nu$. For $w\in\Cplx^+$, the Cauchy transform of $\nu$ is
$g(w)=\int\mu(d\lambda)/(w+\lambda)=\overline{m(-\bar w)}$. Conjugating~\eqref{eq:matrix-dyson} at
the point $-\bar w\in\Cplx^+$ gives~\eqref{eq:hess-pastur}, and $\operatorname{Im}g(w)<0$. Applied
to the diagonal matrices of the quantiles of $\rho$, this shows that~\eqref{eq:hess-pastur} has a
solution for every $\rho\in\Pspace([0,\infty))$; by uniqueness, and since a Cauchy transform
determines its measure, $\nu=\Pastur{a}{\rho}$ is well defined.

\ref{lem:hess-pastur-b} Take $D_N$ with nonnegative entries and $\ESD(D_N)\to\rho$, and let
$\delta>0$; by~\ref{lem:hess-pastur-a} the expected empirical spectral distributions of
$\sqrt a\,\GOE-D_N$ converge weakly to $\nu$. On $\{\|\GOE\|\le2+\delta\}$, $\sqrt a\,\GOE-D_N\preceq\sqrt a(2+\delta)I$, so
$\ESD(\sqrt a\,\GOE-D_N)$ gives no mass to $(\sqrt a(2+\delta),\infty)$. Write
$D_N=D_N\ind_{\{D_N\le s\}}+D_N\ind_{\{D_N>s\}}$. The matrix $\sqrt a\,\GOE-D_N\ind_{\{D_N\le s\}}$
has spectrum in $[-\sqrt a(2+\delta)-s,\infty)$, and subtracting the positive semidefinite matrix
$D_N\ind_{\{D_N>s\}}$, of rank $\#\{i:D_{N,ii}>s\}$, moves at most that many eigenvalues below this
point. Hence
$\ESD(\sqrt a\,\GOE-D_N)((-\infty,-\sqrt a(2+\delta)-s))\le\ESD(D_N)((s,\infty))$ on the same event.
Taking expectations, using $\Prob(\|\GOE\|>2+\delta)\to0$, and applying the portmanteau theorem
(open sets on the left, the closed set $[s,\infty)$ on the right), we get
$\nu((\sqrt a(2+\delta),\infty))=0$ and $\nu((-\infty,-\sqrt a(2+\delta)-s))\le\rho([s,\infty))$;
let $\delta\downarrow0$.

\ref{lem:hess-pastur-c} Let $\nu_n=\Pastur{a_n}{\rho_n}$ with Cauchy transforms $g_n$. By
\ref{lem:hess-pastur-b} and the tightness of $(\rho_n)$, the family $(\nu_n)$ is tight. Let $\nu^*$
be a subsequential weak limit; its Cauchy transform $g^*$ is the pointwise limit of $g_n$. For
$w\in\Cplx^+$ we have $\operatorname{Im}(w-a_ng_n(w))\ge\operatorname{Im}w$, so the integrands
$D\mapsto(w+D-a_ng_n(w))^{-1}$ are bounded by $1/\operatorname{Im}w$ and converge uniformly on
$[0,\infty)$ to $(w+D-ag^*(w))^{-1}$. Since $\rho_n\to\rho$ weakly, $g^*$
solves~\eqref{eq:hess-pastur} for $(a,\rho)$, and $\operatorname{Im}g^*<0$; by uniqueness
$\nu^*=\Pastur{a}{\rho}$. The uniform continuity follows, because
$[0,A]\times\Pspace([0,B])$ is a compact metric space.

\ref{lem:hess-pastur-d} Take independent copies $(D_i,D_i')$, $i\ge1$, of the pair $(D,D')$,
independent of the matrices $\GOE$, and put $D_N=\diag(D_i)_{i\le N}$ and
$D_N'=\diag(D_i')_{i\le N}$. Almost surely $\ESD(D_N)\to\rho$ and $\ESD(D_N')\to\rho'$, so by
\ref{lem:hess-pastur-a}, applied conditionally on the diagonals, $\ESD(\sqrt a\,\GOE-D_N)$ and
$\ESD(\sqrt a\,\GOE-D_N')$ converge in probability to $\Pastur{a}{\rho}$ and $\Pastur{a}{\rho'}$. If
$|D-D'|\le\varepsilon$ almost surely, their L\'evy distance is at most
$\max_i|D_i-D_i'|\le\varepsilon$, by~\eqref{eq:matrix-weyl}. If $\Prob(D\ne D')\le p$, it is at most
$N^{-1}\#\{i\le N:D_i\ne D_i'\}\to\Prob(D\ne D')$, by~\eqref{eq:matrix-rank} and the law of large
numbers. Similarly
$\dL(\ESD(\sqrt a\,\GOE-D_N),\ESD(\sqrt{a'}\,\GOE-D_N))\le|\sqrt a-\sqrt{a'}|(2+\delta)$ on
$\{\|\GOE\|\le2+\delta\}$. In each case the bound passes to the limits by the lower
semicontinuity of $\dL$ under weak convergence, along a subsequence with almost sure convergence.
\end{proof}

The edge function $z_Q$ of Section~\ref{sec:hess-edge} reaches its minimum inside its domain
because $\UQ$ puts enough mass near its essential supremum. This is where the diffusion $X$ enters:
its law at time $Q$ is equivalent to a Gaussian law, so $X_Q$ visits every neighbourhood of the
maximizer of $\Phi_{\zf,xx}(Q,\cdot)$.

\begin{lemma}[The edge hypothesis for the Parisi variables]\label{lem:matrix-hyp}
Let $Q\in(0,1)$, and let $u_{\max}$ be the maximum of $b\mapsto\Phi_{\zf,xx}(Q,b)$ over $\R$.
Then $u_{\max}$ is the essential supremum of $\UQ$, and there are $c,\varepsilon_0>0$ such that
\[
 \Prob\bigl(\UQ^{-1}-u_{\max}^{-1}<\varepsilon\bigr)\ge c\,\varepsilon\qquad(0<\varepsilon\le\varepsilon_0).
\]
In particular $U=\UQ$ and $a=\aQ$ satisfy~\eqref{eq:hess-fhyp}.
\end{lemma}

\begin{proof}
Write $\Phi=\Phi_{\zf}$. By Lemma~\ref{lem:pde-standing}\ref{lem:pde-standing-a} the function
$b\mapsto\Phi_{xx}(Q,b)$ is positive and continuous and tends to $0$ at $\pm\infty$, so it attains
its maximum $u_{\max}$ at some $b_0$; it is twice differentiable with
$|\partial_b^2\Phi_{xx}|=|\Phi_{xxxx}|\le C_{34}$. Hence
$u_{\max}-\Phi_{xx}(Q,b)\le\frac12C_{34}(b-b_0)^2$, and if $|b-b_0|\le r_1$, where $r_1>0$ is such
that $\Phi_{xx}(Q,\cdot)\ge u_{\max}/2$ on $[b_0-r_1,b_0+r_1]$,
\[
 \Phi_{xx}(Q,b)^{-1}-u_{\max}^{-1}=\frac{u_{\max}-\Phi_{xx}(Q,b)}{\Phi_{xx}(Q,b)\,u_{\max}}
 \le\frac{C_{34}}{u_{\max}^2}(b-b_0)^2 .
\]
Let $X$ be the Auffinger--Chen process of $\zf$. Its drift is $\sqrt{\xi''(t)}\,\theta_t$ with
$\theta_t=\sqrt{\xi''(t)}\cdf{\zf}(t)\Phi_x(t,X_t)$ and $|\theta_t|\le\sqrt{\xi''(1)}$. Let
$Z=\exp(-\int_0^Q\theta\,dB-\frac12\int_0^Q\theta^2\,dt)$. By Girsanov's theorem, under the measure
$Z\,d\Prob$ the variable $X_Q$ is centred Gaussian with variance $\xi'(Q)>0$, and
$\E Z^2\le e^{\xi''(1)}$, because $\exp(-2\int\theta\,dB-2\int\theta^2)$ is a martingale. By
Cauchy--Schwarz, for every Borel set $A$,
$\Prob(X_Q\in A)\ge e^{-\xi''(1)}\Prob(\sqrt{\xi'(Q)}\,Z'\in A)^2$ with $Z'$ standard Gaussian. So
$X_Q$ has full support, which gives the first claim, and
$\Prob(|X_Q-b_0|<s)\ge c_1s^2$ for $s\le1$. Since $\UQ=\Phi_{xx}(Q,X_Q)$, the event
$\{|X_Q-b_0|<\min(r_1,(\varepsilon/C_{34})^{1/2}u_{\max})\}$ is contained in
$\{\UQ^{-1}-u_{\max}^{-1}<\varepsilon\}$, whose probability is therefore at least $c\varepsilon$ for
small $\varepsilon$. With $D=\UQ^{-1}+\aQ\chiQ$ we have $D-\essinf D=\UQ^{-1}-u_{\max}^{-1}$, and
\[
 \E\bigl[(D-\essinf D)^{-2}\bigr]=\int_0^\infty\Prob\bigl(D-\essinf D<y^{-1/2}\bigr)\,dy
 \ge\int_{\varepsilon_0^{-2}}^\infty c\,y^{-1/2}\,dy=\infty .
 \qedhere
\]
\end{proof}

\begin{proof}[Proof of Lemma~\ref{lem:hess-edge}\ref{lem:hess-edge-c}--\ref{lem:hess-edge-e}]
We use the notation of Section~\ref{sec:hess-edge} and parts~\ref{lem:hess-edge-a}--\ref{lem:hess-edge-b},
proved there. Recall that $U\le u_{\max}\le1$, so $U^{-1}\ge1$.

\ref{lem:hess-edge-c} Put $z_2=z''(-a\chi)=2a\E U^3$ and $I=[-a\chi-\frac12,-a\chi+\frac12]$, which
lies in $(s_-,\infty)$ because $-a\chi-s_-=u_{\max}^{-1}\ge1$. For $s\in I$,
$s+D=U^{-1}+(s+a\chi)\ge\frac12U^{-1}$, so
$z''(s)=2a\E[U^3(1+(s+a\chi)U)^{-3}]\ge2a(2/3)^3\E U^3=:c_I$ and
$|z'''(s)|=6a\E(s+D)^{-4}\le96a$. Let $c_0=c_I/2$ and $|\Plef|\le c_0$. Since $z'(-a\chi)=\Plef$ and
$z''\ge c_I$ on $I$, the critical point lies in $I$, and $h:=s^*+a\chi$ satisfies
$|h|\le|\Plef|/c_I\le\frac12$. Taylor's formula at $-a\chi$ gives
$0=z'(s^*)=\Plef+z_2h+R_1$ with $|R_1|\le48ah^2$, and
$z(s^*)=\Plef h+\frac12z_2h^2+R_2$ with $|R_2|\le16a|h|^3$. Since
$\Plef h+\frac12z_2h^2=-\Plef^2/(2z_2)+\frac{z_2}2(h+\Plef/z_2)^2$ and $|h+\Plef/z_2|=|R_1|/z_2$,
\[
 \Bigl|\lamp+\frac{\Plef^2}{4a\E U^3}\Bigr|=\Bigl|z(s^*)+\frac{\Plef^2}{2z_2}\Bigr|
 \le\frac{(48a)^2h^4}{2z_2}+16a|h|^3\le C|\Plef|^3 ,
\]
with $C$ depending only on $a$ and $\E U^3$.

\ref{lem:hess-edge-d} Let $g$ be the Cauchy transform of $\nu$, and $G(\sigma)=\E(\sigma+D)^{-1}$
for $\sigma>s_-$. The proof of part~\ref{lem:hess-edge-a} in Section~\ref{sec:hess-edge} shows that
every $w>\lamp$ satisfies $g(w)=G(s(w))$ and $z(s(w))=w$, where $s(w)=w-ag(w)>s^*$; thus $s(w)$ is
the inverse of the increasing branch of $z$, and $s(w)\downarrow s^*$ as $w\downarrow\lamp$. Since
$G'=(z'-1)/a$ and $s'(w)=1/z'(s(w))$,
\begin{equation}\label{eq:matrix-gprime}
 -g'(w)=\int\frac{\nu(d\lambda)}{(w-\lambda)^2}=\frac{1-z'(s(w))}{a\,z'(s(w))}\qquad(w>\lamp).
\end{equation}
Because $z'''<0$, $z''$ is decreasing, so $c_2\le z''\le C_2$ on $[s^*,s^*+1]$. Let $t=w-\lamp$ and
$s=s(w)$, and assume $t\le t_0$, so that $s\le s^*+1$. Then $t=\int_{s^*}^sz'\le(s-s^*)z'(s)$ and
$z'(s)\ge c_2(s-s^*)$, so $t\le z'(s)^2/c_2$; and $t\ge c_2(s-s^*)^2/2$, $z'(s)\le C_2(s-s^*)$.
Hence
\begin{equation}\label{eq:matrix-zprime}
 \sqrt{c_2t}\le z'(s(w))\le C_2\sqrt{2t/c_2}\qquad(0<t=w-\lamp\le t_0).
\end{equation}
\emph{Upper bound.} For $0<\kappa\le t_0$ take $w=\lamp+\kappa$; every
$\lambda\in[\lamp-\kappa,\lamp]$ has $w-\lambda\le2\kappa$, so
by~\eqref{eq:matrix-gprime}--\eqref{eq:matrix-zprime}
\[
 \nu([\lamp-\kappa,\lamp])\le4\kappa^2(-g'(w))\le\frac{4\kappa^2}{a\sqrt{c_2\kappa}}=C_4\kappa^{3/2},
 \qquad C_4:=\frac4{a\sqrt{c_2}} .
\]
\emph{Lower bound.} Put $F(y)=\nu([\lamp-y,\lamp])$. Then $F(0)=0$, since an atom at $\lamp$ would
make $g(w)$ unbounded as $w\downarrow\lamp$, while $g(w)=G(s(w))\to G(s^*)<\infty$. As $\nu$ is
carried by $(-\infty,\lamp]$, integration by parts gives, for $w=\lamp+t$,
$-g'(w)=\int_0^\infty2F(y)(t+y)^{-3}\,dy$. If moreover $t\le c_2/(8C_2^2)$, then $z'(s)\le\frac12$
by~\eqref{eq:matrix-zprime}, and~\eqref{eq:matrix-gprime} gives $-g'(w)\ge1/(2az'(s))\ge c_3t^{-1/2}$
with $c_3=\sqrt{c_2/2}/(2aC_2)$. Split the integral at $\Theta t$ and at $y_0=t_0$, where
$\Theta=(16C_4/c_3)^2$; using $F\le1$, $F(y)\le C_4y^{3/2}$ for $y\le t_0$, and monotonicity of $F$,
\[
 c_3t^{-1/2}\le\frac{F(\Theta t)}{t^2}+\frac{4C_4}{\sqrt{\Theta t}}+\frac1{t_0^2}
 =\frac{F(\Theta t)}{t^2}+\frac{c_3}{4\sqrt t}+\frac1{t_0^2}
\]
(if $\Theta t>t_0$ the middle range is empty and the bound still holds). For
$t\le(c_3t_0^2/4)^2$ the last term is at most $c_3/(4\sqrt t)$, so $F(\Theta t)\ge\frac{c_3}2t^{3/2}$.
With
\[
 t_1=\min\Bigl\{t_0,\ \frac{c_2}{8C_2^2},\ \Bigl(\frac{c_3t_0^2}4\Bigr)^2\Bigr\},\qquad
 \kappa_0=\min\{t_0,\Theta t_1\},\qquad c=\frac{c_3}2\Theta^{-3/2},\qquad C=C_4,
\]
this gives $c\kappa^{3/2}\le\nu([\lamp-\kappa,\lamp])\le C\kappa^{3/2}$ for $0<\kappa\le\kappa_0$.
The constants $c$, $C$, $\kappa_0$ are continuous, strictly positive functions of
$(a,C_2,c_2,t_0)\in(0,\infty)^4$.

\ref{lem:hess-edge-e} Lemma~\ref{lem:matrix-hyp} gives~\eqref{eq:hess-fhyp} for $U=\UQ$, and with
$a=\aQ$ the objects of Section~\ref{sec:hess-edge} are $\nu=\nuQ$, $\lamp=\lamp(Q)$ and
$\Plef=\Plef(Q)$, by definition. For continuity, write $z_Q$, $s_-(Q)$, $s^*_Q$ for the objects
attached to $Q$. By Lemma~\ref{lem:pde-standing}\ref{lem:pde-standing-a},
$\Phi_{\zf,xx}(Q,b)\le\overline Ce^{-2|b|}$ and $\Phi_{\zf,xx}(Q,0)\ge\underline c$, so
$u_{\max}(Q)$ is a maximum over $|b|\le\frac12\log(\overline C/\underline c)$ of a jointly
continuous function, hence continuous in $Q$. The map $Q\mapsto\Law(\UQ)$ is weakly continuous
(Lemma~\ref{lem:pde-standing}\ref{lem:pde-standing-d}), so $\chiQ$ and $s_-(Q)$ are continuous.
Fix $(Q_0,s_0)$ with $s_0>s_-(Q_0)$, and write $\sigma=s+\aQ\chiQ$, so that
$\E(s+D_Q)^{-j}=\E[\UQ^j(1+\sigma\UQ)^{-j}]$. For $(Q,s)$ near $(Q_0,s_0)$ we have
$\sigma\ge-1/(u_{\max}(Q_0)+2\delta_1)$ and $u_{\max}(Q)\le u_{\max}(Q_0)+\delta_1$ for some
$\delta_1>0$; so the laws of $\UQ$ are carried by $[0,u_{\max}(Q_0)+\delta_1]$, where
$1+\sigma U\ge\delta_1/(u_{\max}(Q_0)+2\delta_1)$. The integrand is therefore bounded and jointly
continuous there, and weak continuity of $\Law(\UQ)$ shows that $z_Q(s)$, $z_Q'(s)$ and $z_Q''(s)$
are jointly continuous near $(Q_0,s_0)$. Since $z_Q'$ is strictly increasing with the unique zero
$s^*_Q$, the strict signs of $z_{Q_0}'(s^*_{Q_0}\pm\delta)$ persist for $Q$ near $Q_0$, so
$Q\mapsto s^*_Q$ is continuous. Hence $\lamp(Q)=z_Q(s^*_Q)$ and the quantities $C_2=z_Q''(s^*_Q)$,
$c_2=z_Q''(s^*_Q+1)$ and $t_0=z_Q(s^*_Q+1)-\lamp(Q)$ are continuous and positive, and so are the
constants of~\ref{lem:hess-edge-d}. For $Q\in\Cf$,
Lemma~\ref{lem:pde-standing}\ref{lem:pde-standing-c} gives $\Plef(Q)=2\Obs_{\zf}''(Q)/\xi''(Q)\ge0$.
Substituting into~\ref{lem:hess-edge-b} gives $\lamp(Q)=0$ if and only if $\Obs_{\zf}''(Q)=0$, and
into~\ref{lem:hess-edge-c}, for $|\Plef(Q)|\le c_0$,
$\lamp(Q)=-\Obs_{\zf}''(Q)^2/(\xi''(Q)^3\E\UQ^3)+O(\Obs_{\zf}''(Q)^3)$, with a constant
depending only on $\xi''(Q)$ and $\E\UQ^3$.
\end{proof}

Appendix~\ref{app:matrix-index} applies the square-root edge to laws obtained from $\Law(X_Q)$ by a
bounded change of density, as in the tilt of Section~\ref{sec:stab-kr}, and needs its constants to
be uniform in the size of the change. The substitution $\sigma=s+a\chi$ fixes the domain of the edge
function, and continuity does the rest.

\begin{lemma}[Uniformity under bounded tilts]\label{lem:matrix-tilt}
Let $a>0$, and let $U$ be a $(0,1]$-valued random variable on a probability space
$(\Omega,\mathcal F,P_0)$ with $P_0$-essential supremum $u_{\max}$ and
$\E_0[(U^{-1}-u_{\max}^{-1})^{-2}]=\infty$. Let $\mathcal R$ be a bounded random variable with
$\E_0\mathcal R=0$, let $\eta_0=1/(2\|\mathcal R\|_\infty)$ (any $\eta_0>0$ if $\mathcal R=0$), and
for $\eta\in[0,\eta_0]$ let $dP_\eta=(1-\eta\mathcal R)\,dP_0$, $\chi_\eta=\E_\eta U$,
$\rho_\eta=\Law_{P_\eta}(U^{-1}+a\chi_\eta)$, $\nu^{(\eta)}=\Pastur{a}{\rho_\eta}$ and
$\lamp^{(\eta)}=\sup\supp\nu^{(\eta)}$.
\begin{thmlist}
\item\label{lem:matrix-tilt-a} Each $P_\eta$ satisfies~\eqref{eq:hess-fhyp}, so
 Lemma~\ref{lem:hess-edge} applies to $\nu^{(\eta)}$.
\item\label{lem:matrix-tilt-b} The maps $\eta\mapsto\lamp^{(\eta)}$ and
 $\eta\mapsto(C_2,c_2,t_0)(\eta)$, the latter defined as in
 Lemma~\ref{lem:hess-edge}\ref{lem:hess-edge-d}, are continuous on $[0,\eta_0]$, the latter with
 values in $(0,\infty)^3$.
\item\label{lem:matrix-tilt-c} There are $0<c_\sharp\le C_\sharp$ and $\kappa_\sharp>0$ such that,
 for all $\eta\in[0,\eta_0]$ and $\kappa\in(0,\kappa_\sharp]$,
 \[
  c_\sharp\kappa^{3/2}\le\nu^{(\eta)}\bigl([\lamp^{(\eta)}-\kappa,\lamp^{(\eta)}]\bigr)\le C_\sharp\kappa^{3/2}.
 \]
\end{thmlist}
\end{lemma}

\begin{proof}
\ref{lem:matrix-tilt-a} The density $1-\eta\mathcal R$ lies in $[\frac12,\frac32]$, so $P_\eta$
and $P_0$ have the same null sets, $U$ has essential supremum $u_{\max}$ under $P_\eta$, and with
$D=U^{-1}+a\chi_\eta$,
$\E_\eta[(D-\essinf D)^{-2}]=\E_\eta[(U^{-1}-u_{\max}^{-1})^{-2}]\ge\frac12\E_0[(U^{-1}-u_{\max}^{-1})^{-2}]=\infty$.

\ref{lem:matrix-tilt-b} Let $z_\eta$ be the edge function of $\nu^{(\eta)}$. The substitution
$\sigma=s+a\chi_\eta$ maps its domain $s>-(u_{\max}^{-1}+a\chi_\eta)$ onto $\sigma>-u_{\max}^{-1}$,
for every $\eta$, and
\[
 z_\eta(\sigma-a\chi_\eta)=\hat z_\eta(\sigma)-a\chi_\eta,\qquad
 \hat z_\eta(\sigma):=\sigma+a\,\E_0\Bigl[\frac{1-\eta\mathcal R}{\sigma+U^{-1}}\Bigr].
\]
For $k\le2$, $\partial_\sigma^k\hat z_\eta(\sigma)$ is affine in $\eta$, with coefficients that are
bounded continuous functions of $\sigma\ge\sigma_1$ for every $\sigma_1>-u_{\max}^{-1}$, by dominated
convergence with the bound $\|\mathcal R\|_\infty^j(\sigma_1+u_{\max}^{-1})^{-1-k}$. So
$\hat z_\eta$, $\hat z_\eta'$ and $\hat z_\eta''$ are jointly continuous on
$[0,\eta_0]\times(-u_{\max}^{-1},\infty)$, and $\chi_\eta=\chi_0-\eta\E_0[\mathcal RU]$ is
continuous. By Lemma~\ref{lem:hess-edge}\ref{lem:hess-edge-a}, $\hat z_\eta'$ is strictly
increasing with a unique zero $\sigma^*_\eta$, and the strict signs of
$\hat z_{\eta_1}'(\sigma^*_{\eta_1}\pm\delta)$ persist for $\eta$ near $\eta_1$; so
$\eta\mapsto\sigma^*_\eta$ is continuous. Derivatives in $s$ and in $\sigma$ coincide, so
$\lamp^{(\eta)}=\hat z_\eta(\sigma^*_\eta)-a\chi_\eta$, $C_2(\eta)=\hat z_\eta''(\sigma^*_\eta)$,
$c_2(\eta)=\hat z_\eta''(\sigma^*_\eta+1)$ and
$t_0(\eta)=\hat z_\eta(\sigma^*_\eta+1)-\hat z_\eta(\sigma^*_\eta)$ are continuous; the last three
are positive because $\hat z_\eta''>0$.

\ref{lem:matrix-tilt-c} The constants $c$, $C$, $\kappa_0$ of
Lemma~\ref{lem:hess-edge}\ref{lem:hess-edge-d} are continuous positive functions of
$(a,C_2,c_2,t_0)$. On the compact interval $[0,\eta_0]$ put $c_\sharp=\min_\eta c(\eta)$,
$C_\sharp=\max_\eta C(\eta)$ and $\kappa_\sharp=\min_\eta\kappa_0(\eta)$.
\end{proof}

\subsection{Uniform spectral estimates over the cube}\label{app:matrix-uniform}

At a fixed point $x$ and for a fixed diagonal $D$, Lemmas~\ref{lem:matrix-third}\ref{lem:matrix-third-a}
and~\ref{lem:hess-pastur}\ref{lem:hess-pastur-a} show that $\ESD(\Hess\HN(x)-D)$ is close to
$\Pastur{\xi''(q_x)}{\ESD(D)}$. Section~\ref{sec:hess} needs this at random points, so we prove it
for all points of the cube and all bounded diagonal matrices at once. For the
Sherrington--Kirkpatrick model, estimates of this kind are proved in \cite[Section~4]{JSS}.

\begin{lemma}[Uniform spectral estimate]\label{lem:matrix-uniform}
Assume \condR. For every $K,\eta,M>0$ there is $N_0$ such that, for $N\ge N_0$, outside an event of
outer probability at most $e^{-MN}$,
\[
 \dL\Bigl(\ESD\bigl(\Hess\HN(x)-D\bigr),\Pastur{\xi''(q_x)}{\ESD(D)}\Bigr)\le\eta
 \qquad\text{for every } x\in\cube \text{ and } D\in\mathcal D_K .
\]
\end{lemma}

\begin{proof}
The proof has four steps: two uniform events, a concentration bound at a fixed pair $(x,D)$ for a
finite family of test functions, a union bound over a net, and interpolation between net points.

\emph{Step 1: two uniform events.} Let $C^H_2,C^H_3$ be given by
Lemma~\ref{lem:matrix-third}\ref{lem:matrix-third-b} with $M+1$ in place of $M$, and let
$\Omega_2=\{\sup_{\cube}\|\Hess\HN\|\le C^H_2\}$ and $\Omega_3=\{\sup_{\cube}\|D^3\HN\|\le C^H_3N^{-1/2}\}$.
Each fails with probability at most $e^{-(M+1)N}$.

\emph{Step 2: a fixed pair.} Let $h\in(0,1)$, $J=[-\ell,\ell]$ with $\ell=C^H_2+K+2\sqrt{\xi''(1)}+1$, and
$\varphi_y(\lambda)=\min\{1,\max\{0,(y+h-\lambda)/h\}\}$ for $y$ in a grid $\mathcal G\subset[-\ell,\ell+h]$ of
mesh $h$. If two probability measures carried by $J$ satisfy
$|\int\varphi_y\,d\mu-\int\varphi_y\,d\mu'|\le h$ for every $y\in\mathcal G$, then
$\mu((-\infty,t])\le\int\varphi_y\,d\mu\le\int\varphi_y\,d\mu'+h\le\mu'((-\infty,t+2h])+h$ for the grid
point $y\in[t,t+h]$, and symmetrically; so $\dL(\mu,\mu')\le2h$.

Fix $x\in\cube$ and $D\in\mathcal D_K$, and put $a=\xi''(q_x)$. The Hessian $\Hess\HN(x)$ is a
centred Gaussian vector in the space of symmetric matrices with the Frobenius inner product, and
by~\eqref{eq:matrix-hess-cov} its covariance operator has norm at most $C_H/N$; so it is the image of a
standard Gaussian vector under a linear map of operator norm at most $(C_H/N)^{1/2}$. By the
Hoffman--Wielandt inequality \cite{HW} and Cauchy--Schwarz, $A\mapsto N^{-1}\tr\varphi_y(A-D)$ is
$h^{-1}N^{-1/2}$-Lipschitz in the Frobenius norm. Gaussian concentration \cite{Ledoux} therefore
gives
\[
 \Prob\Bigl(\Bigl|\frac1N\tr\varphi_y\bigl(\Hess\HN(x)-D\bigr)-\E[\,\cdot\,]\Bigr|>\frac h2\Bigr)
 \le2\exp\Bigl(-\frac{h^4N^2}{8C_H}\Bigr).
\]
By the decomposition of Lemma~\ref{lem:matrix-third}\ref{lem:matrix-third-a},
\eqref{eq:matrix-rank} and~\eqref{eq:matrix-monotone}, the expectation differs by at most $2/N$
from $\E N^{-1}\tr\varphi_y(\sqrt a\,\GOE-D)$. We claim that for $N\ge N_1(h,K)$, uniformly in
$a\in[0,\xi''(1)]$, $D\in\mathcal D_K$ and $y\in\mathcal G$,
\[
 \Bigl|\E\frac1N\tr\varphi_y\bigl(\sqrt a\,\GOE-D\bigr)-\int\varphi_y\,d\Pastur{a}{\ESD(D)}\Bigr|\le\frac h4 .
\]
Otherwise there are $N_k\to\infty$, $a_k$, $D_k\in\mathcal D_K$ and $y_k$ violating it; passing to a
subsequence, $a_k\to a$, $\ESD(D_k)\to\rho\in\Pspace([0,K])$ and $y_k=y$. By
Lemma~\ref{lem:hess-pastur}\ref{lem:hess-pastur-a} the expected empirical spectral distribution of
$\sqrt a\,\GOE-D_k$ converges to $\Pastur{a}{\rho}$, and by~\eqref{eq:matrix-weyl} and
$\Prob(\|\GOE\|>3)\le2e^{-N/4}$ the same holds for $\sqrt{a_k}\,\GOE-D_k$. By
Lemma~\ref{lem:hess-pastur}\ref{lem:hess-pastur-c}, $\Pastur{a_k}{\ESD(D_k)}\to\Pastur{a}{\rho}$.
Since $\varphi_y$ is bounded and continuous, both integrals converge to the same limit, a
contradiction.

\emph{Step 3: a net.} Let $\mathcal N_x\subset\cube$ be an $\epsilon_0\sqrt N$-net of the cube in the
Euclidean norm, with $|\mathcal N_x|\le(6/\epsilon_0)^N$, and let $\mathcal N_D\subset\mathcal D_K$ be
an $\epsilon_0$-net in the supremum norm, with $|\mathcal N_D|\le(K/\epsilon_0+1)^N$. By Step 2, for
$N\ge\max(N_1,8/h)$, the probability that some $(y,D')\in\mathcal N_x\times\mathcal N_D$ and some
$y_0\in\mathcal G$ have
$|N^{-1}\tr\varphi_{y_0}(\Hess\HN(y)-D')-\int\varphi_{y_0}\,d\Pastur{\xi''(q_y)}{\ESD(D')}|>h$ is at most
$e^{CN}e^{-h^4N^2/(8C_H)}\le e^{-(M+1)N}$ for large $N$. Outside this event and on $\Omega_2$, both
$\ESD(\Hess\HN(y)-D')$ and $\Pastur{\xi''(q_y)}{\ESD(D')}$ are carried by $J$, the latter by
Lemma~\ref{lem:hess-pastur}\ref{lem:hess-pastur-b} (its support lies in
$[-K-2\sqrt{\xi''(1)},2\sqrt{\xi''(1)}]$); so their L\'evy distance is at most $2h$ for every net
pair.

\emph{Step 4: interpolation.} Let $x\in\cube$ and $D\in\mathcal D_K$, and pick a net pair $(y,D')$
with $|x-y|\le\epsilon_0\sqrt N$ and $\|D-D'\|\le\epsilon_0$. On $\Omega_3$,
$\|(\Hess\HN(x)-D)-(\Hess\HN(y)-D')\|\le(C^H_3+1)\epsilon_0$, which by~\eqref{eq:matrix-weyl} moves
the empirical spectral distribution by at most $(C^H_3+1)\epsilon_0$ in L\'evy distance. For the
targets, Lemma~\ref{lem:hess-pastur}\ref{lem:hess-pastur-d}, with the diagonal entries of $D$ and
$D'$ coupled through a uniform index, gives
$\dL(\Pastur{\xi''(q_y)}{\ESD(D')},\Pastur{\xi''(q_x)}{\ESD(D)})\le\epsilon_0+2|\sqrt{\xi''(q_y)}-\sqrt{\xi''(q_x)}|$,
and $|q_x-q_y|\le2\epsilon_0$. Choose $h$ and $\epsilon_0$ such that
$(C^H_3+2)\epsilon_0+2h+2\omega(2\epsilon_0)\le\eta$, where $\omega$ is the modulus of continuity of
$\sqrt{\xi''}$ on $[0,1]$. The bad event is contained in the union of $\Omega_2^c$, $\Omega_3^c$ and
the event of Step 3, whose probability is at most $3e^{-(M+1)N}\le e^{-MN}$.
\end{proof}

\begin{proof}[Proof of Proposition~\ref{prop:hess-bulk}]
Write $u=u_{\zf}$ and $a_t=\xi''(t)$. By Lemma~\ref{lem:pde-standing}\ref{lem:pde-standing-b},
$u$ is uniformly continuous on $[0,1]\times[-1,1]$; let $\omega_u$ be its modulus for the distance
$|t-t'|+|m-m'|$. For $K\ge1$ the function $\psi_K=1/\max(u,K^{-1})=\min(u^{-1},K)$ has modulus at
most $K^2\omega_u$, because $y\mapsto1/\max(y,K^{-1})$ is $K^2$-Lipschitz on $[0,1]$. For
$(Q,\mu)\in\mathcal K$ let $\rho_K(Q,\mu)$ be the law of $\psi_K(Q,m)+a_Q\int u(Q,\cdot)\,d\mu$ under
$m\sim\mu$. Fix $\kappa>0$ and choose, in order:
\begin{steplist}
\item\label{step:matrix-bulk-K} $K\ge1$ with
 $\sup_{(Q,\mu)\in\mathcal K}\mu(\{m:u(Q,m)\le4/K\})\le\kappa/8$. This is possible because
 $(Q,\mu)\mapsto\mu(\{u(Q,\cdot)\le\epsilon\})$ is upper semicontinuous (the set
 $\{(t,m):u(t,m)\le\epsilon\}$ is closed) and decreases, as $\epsilon\downarrow0$, to
 $\mu(\{u(Q,\cdot)=0\})=\mu(\{-1,1\})=0$; on the compact set $\mathcal K$ the convergence is uniform
 (Dini).
\item\label{step:matrix-bulk-E} The event of Lemma~\ref{lem:matrix-uniform} with
 $(K+\xi''(1),\kappa/8,M+1)$ in place of $(K,\eta,M)$; let $\Omega_E$ be its complement.
\item\label{step:matrix-bulk-w} $w\in(0,r_0/2]$ so small that $(K^2+\xi''(1))\,\omega_S(w)\le\kappa/8$
 and $\omega_S(w)\le1/K$, with $\omega_S$ from Lemma~\ref{lem:matrix-orderstab}; that
 $\tau_2(w):=K(\omega_u(2w^{1/2})+w^{1/2})\le\kappa/8$; and that
 \[
  \tau_1(w):=\Bigl(K^2\omega_u(2w^{1/2})+Kw^{1/2}+\xi''(1)\bigl(\omega_u(2w^{1/2})+w^{1/2}\bigr)+\xi'''(1)w\Bigr)^{1/2}
 \]
 and $\xi'''(1)w$ have sum so small that the modulus of uniform continuity in
 Lemma~\ref{lem:hess-pastur}\ref{lem:hess-pastur-c} on $[0,\xi''(1)]\times\Pspace([0,K+\xi''(1)])$
 is at most $\kappa/8$.
\end{steplist}
Now let $N\ge8/\kappa$ be large, and on $\Omega_E$ take $(Q,\mu)\in\mathcal K$, $x\in\ocube$ and
$\alpha$ satisfying the three window conditions; put $\sigma\in\{+,-\}$. Note
$|q_x-Q|\le w\le r_0/2$, so $q_x\in(0,1)$.
\begin{enumerate}[label=\textup{(\alph*)},leftmargin=*,widest=f,topsep=3pt,itemsep=2pt]
\item By Corollary~\ref{cor:matrix-truncated},
 $\dL(\ESD(\nabla^2_\sigma F_\alpha(x)),\ESD(\Hess\HN(x)-D^\alpha_K(x)))\le p^\alpha_K(x)+2/N$.
\item By Lemma~\ref{lem:matrix-uniform}, since $D^\alpha_K(x)\in\mathcal D_{K+\xi''(1)}$,
 $\dL(\ESD(\Hess\HN(x)-D^\alpha_K(x)),\Pastur{a_{q_x}}{\ESD(D^\alpha_K(x))})\le\kappa/8$.
\item By Lemma~\ref{lem:matrix-orderstab} at $t=q_x$,
 $\sup_m|u_\alpha(q_x,m)-u(q_x,m)|\le\omega_S(w)$. Since $y\mapsto\min(y^{-1},K)$ is
 $K^2$-Lipschitz on $(0,1]$ and $|\chi_\alpha(x)-\chi_{\zf}(x)|\le\omega_S(w)$, the diagonal
 matrices $D^\alpha_K(x)$ and $D^{\zf}_K(x)$ differ by at most $(K^2+\xi''(1))\omega_S(w)\le\kappa/8$
 entrywise, and Lemma~\ref{lem:hess-pastur}\ref{lem:hess-pastur-d} moves the target by at most
 $\kappa/8$.
\item Couple $m\sim\mu_x$ and $m'\sim\mu$ with $\E|m-m'|\le w$. On $\{|m-m'|\le w^{1/2}\}$ the
 points $(q_x,m)$ and $(Q,m')$ are within distance $2w^{1/2}$, and the complement has probability
 at most $w^{1/2}$. Hence $\E|\psi_K(q_x,m)-\psi_K(Q,m')|\le K^2\omega_u(2w^{1/2})+Kw^{1/2}$ and
 $\E|u(q_x,m)-u(Q,m')|\le\omega_u(2w^{1/2})+w^{1/2}$; also $|a_{q_x}-a_Q|\le\xi'''(1)w$. Since
 $\ESD(D^{\zf}_K(x))$ is the law of $\psi_K(q_x,m)+a_{q_x}\int u(q_x,\cdot)\,d\mu_x$,
 \eqref{eq:matrix-coupling} gives $\dL(\ESD(D^{\zf}_K(x)),\rho_K(Q,\mu))\le\tau_1(w)$, and the
 choice of $w$ and Lemma~\ref{lem:hess-pastur}\ref{lem:hess-pastur-c} give
 $\dL(\Pastur{a_{q_x}}{\ESD(D^{\zf}_K(x))},\Pastur{a_Q}{\rho_K(Q,\mu)})\le\kappa/8$.
\item The variables $\psi_K(Q,m')+c$ and $u(Q,m')^{-1}+c$, with $c=a_Q\int u(Q,\cdot)\,d\mu$, differ
 only when $u(Q,m')<1/K$. Lemma~\ref{lem:hess-pastur}\ref{lem:hess-pastur-d} and the choice of $K$
 give $\dL(\Pastur{a_Q}{\rho_K(Q,\mu)},\Pastur{a_Q}{\rho_{Q,\mu}})\le\kappa/8$.
\item Since $\sup_m|u_\alpha(q_x,m)-u(q_x,m)|\le1/K$, every $i$ with $u_\alpha(q_x,x_i)<1/K$ has
 $u(q_x,x_i)<2/K$. Let $\phi(t,m)=\min\{1,\max\{0,3-Ku(t,m)\}\}$, which equals one where
 $u\le2/K$, vanishes where $u\ge3/K$, and is $K$-Lipschitz in $u$. With the coupling of (d),
 \[
  p^\alpha_K(x)\le\int\phi(q_x,\cdot)\,d\mu_x\le\E\phi(Q,m')+K\E\bigl|u(q_x,m)-u(Q,m')\bigr|
 \]
 and the last sum is at most $\mu(\{u(Q,\cdot)\le3/K\})+\tau_2(w)\le\kappa/4$.
\end{enumerate}
By the triangle inequality,
$\dL(\ESD(\nabla^2_\sigma F_\alpha(x)),\Pastur{\xi''(Q)}{\rho_{Q,\mu}})\le\frac\kappa4+\frac2N+\frac{4\kappa}8<\kappa$.
The exceptional event is $\Omega_E^c$, of outer probability at most $e^{-(M+1)N}\le e^{-MN}$.
\end{proof}

\subsection{Order parameters of near-optimal states}\label{app:matrix-own}

Proposition~\ref{prop:hess-own} rests on one deterministic fact: at a contact $Q$, the band problem
$\beta\mapsto V_\beta(Q,\Lq{Q})$ has the unique minimizer $\Shift_Q\zf$. Uniqueness is a theorem of
Chen, Panchenko and Subag, and a first variation identifies the minimizer. At the first
contact this minimality is stated in the proof of \cite[Lemma~3.7]{Bou26b}, through
\cite[Lemma~A.4]{Bou26b}; we prove it at every contact. Compactness makes this
quantitative (Lemma~\ref{lem:matrix-own}), and Theorem~\ref{thm:struct}\ref{thm:struct-i}--\ref{thm:struct-ii}
brings every near-maximizer into its range.

\begin{lemma}[The band problem at a contact]\label{lem:matrix-band}
Let $Q\in\Cf$ and $\hat\zeta=\Shift_Q\zf$. Then $\hat\zeta$ is the unique minimizer of
$\beta\mapsto V_\beta(Q,\Lq{Q})$ over $\Pspace([Q,1])$. In particular
$V_{\mathrm{TAP}}(Q,\Lq{Q})=V_{\zf}(Q,\Lq{Q})$.
\end{lemma}

\begin{proof}
\emph{Uniqueness.} At a contact $\Gz{\zf}(Q)=Q$
(Lemma~\ref{lem:pde-standing}\ref{lem:pde-standing-c}), that is, $\int m^2\,\Lq{Q}(dm)=Q$. For such
$\mu$, Chen, Panchenko and Subag define the TAP correction $\TAPcps(\mu,\beta)$ \cite[(2.24)]{CPS}
with the terminal value $\log(2\cosh)$; hence $\TAPcps(\mu,\beta)=V_\beta(Q,\mu)+\log2$. They prove
that its infimum over $\Pspace([Q,1])$ is attained at a unique point
\cite[(2.26) and Theorem~12(i), (iii)]{CPS}. Since the conjugate is even in $m$, only the law of
$|m|$ under $\mu$ enters there; under $\Lq{Q}$ this law is not carried by $\{0,1\}$, because $X_Q$
has full support (proof of Lemma~\ref{lem:matrix-hyp}) and $\Phi_{\zf,x}(Q,\cdot)$ maps $\R$ onto
$(-1,1)$. Their model class \cite[(1.2)]{CPS} contains every mixture with
$\sum_p2^p\beta_p^2<\infty$, that is, every $\xi\in\classC$. So $\beta\mapsto V_\beta(Q,\Lq{Q})$ has
at most one minimizer over $\Pspace([Q,1])$.

\emph{Convexity.} For a step coefficient $c$ on $[Q,1)$ with values in $[0,1]$, the verification
identity of Lemma~\ref{lem:pde-ac}\ref{lem:pde-ac-b} writes $W[c,\log\cosh](Q,y)$ as the maximum
over bounded progressively measurable controls $u$ of
$\E\log\cosh\bigl(y+\int_Q^1\xi''cu\,dr+\int_Q^1\sqrt{\xi''}\,dB\bigr)-\frac12\E\int_Q^1\xi''cu^2\,dr$.
For fixed $u$ this is jointly convex in $(c,y)$, so the maximum is jointly convex, and by
Lemma~\ref{lem:pde-measurable}\ref{lem:pde-measurable-i} this persists for measurable coefficients.
A partial infimum of a jointly convex function is convex, so $c\mapsto\inf_y\{W[c](Q,y)-my\}$ is
convex. Since $\beta\mapsto\cdf\beta$ is affine and
$\Phi_\beta(Q,\cdot)=W[\cdf\beta|_{[Q,1)},\log\cosh](Q,\cdot)$, the map $\beta\mapsto V_\beta(Q,\mu)$ is
convex on $\Pspace([Q,1])$ for every $\mu$.

\emph{First variation.} Let $\mu=\Lq{Q}$, $\alpha\in\Pspace([Q,1])$, and
$\beta_\epsilon=\hat\zeta+\epsilon(\alpha-\hat\zeta)$ for $\epsilon\in(0,1]$. By convexity,
$V_\alpha(Q,\mu)-V_{\hat\zeta}(Q,\mu)\ge\epsilon^{-1}(V_{\beta_\epsilon}(Q,\mu)-V_{\hat\zeta}(Q,\mu))$. For
$y\in\R$ let $Y^y$ be the Auffinger--Chen process of the coefficient $\cdf{\hat\zeta}=\cdf{\zf}$ on
$[Q,1]$ started at $(Q,y)$, and $\Gamma^y(r)=\E\,\Phi_{\zf,x}(r,Y^y_r)^2$. The left inequality of the
coefficient sandwich (Lemma~\ref{lem:pde-coef}, extended to measurable coefficients by
Lemma~\ref{lem:pde-measurable}) gives
$\Phi_{\beta_\epsilon}(Q,y)-\Phi_{\hat\zeta}(Q,y)\ge\frac\epsilon2\int_Q^1\xi''(r)(\cdf\alpha-\cdf{\hat\zeta})(r)\Gamma^y(r)\,dr$.
Evaluating at $y=b_\epsilon(m):=b_{\beta_\epsilon}(Q,m)$, where the infimum defining
$\Phistar_{\beta_\epsilon}(Q,m)$ is attained, and bounding $\Phistar_{\hat\zeta}(Q,m)$ from above by
the value at the same point,
\[
 \frac{V_{\beta_\epsilon}(Q,\mu)-V_{\hat\zeta}(Q,\mu)}\epsilon
 \ge\frac12\int\mu(dm)\int_Q^1\xi''(r)(\cdf\alpha-\cdf{\hat\zeta})(r)\,\Gamma^{b_\epsilon(m)}(r)\,dr
 -\frac12\int_Q^1r\xi''(r)(\cdf\alpha-\cdf{\hat\zeta})(r)\,dr .
\]
As $\epsilon\downarrow0$, $b_\epsilon(m)\to b_{\hat\zeta}(Q,m)$ for $|m|<1$, by
Lemma~\ref{lem:matrix-orderstab} and the strict monotonicity of $\Phi_{\hat\zeta,x}(Q,\cdot)$; and
$y\mapsto\Gamma^y(r)$ is continuous, since the drift of $Y^y$ is Lipschitz in space. By bounded
convergence the right side tends to
$\frac12\int_Q^1\xi''(r)(\cdf\alpha-\cdf{\hat\zeta})(r)\bigl(\int\Gamma^{b_{\hat\zeta}(Q,m)}(r)\,\mu(dm)-r\bigr)dr$.
Now $\Phi_{\hat\zeta}(Q,\cdot)=\Phi_{\zf}(Q,\cdot)$, so $b_{\hat\zeta}(Q,M_Q)=X_Q$, where $X$ is the
Auffinger--Chen process of $\zf$ and $M_Q=\Phi_{\zf,x}(Q,X_Q)\sim\mu$. By strong uniqueness,
$(X_r)_{r\ge Q}$ is the process $Y^{X_Q}$, so $\int\Gamma^{b_{\hat\zeta}(Q,m)}(r)\,\mu(dm)=\Gz{\zf}(r)$.
With $\Obs=\Obs_{\zf}$, $\Obs'=-\frac12\xi''(\Gz{\zf}-t)$, $\Obs(1)=0$, and
$(\cdf\alpha-\cdf{\hat\zeta})(r)=(\alpha-\hat\zeta)([Q,r])$ for $r\ge Q$, Fubini's theorem turns the
limit into
\[
 -\int_Q^1\Obs'(r)\,(\alpha-\hat\zeta)([Q,r])\,dr=\int_{[Q,1]}\Obs\,d(\alpha-\hat\zeta)
 =\int_{[Q,1]}\bigl(\Obs-\Obs(Q)\bigr)\,d\alpha\ \ge0 .
\]
Here we used Lemma~\ref{lem:loc-obstacle}\ref{lem:loc-obstacle-c}: $\Obs\ge\Omin=\Obs(Q)$ on
$[Q,1]\subset(0,1]$, and $\hat\zeta$ is carried by $\{Q\}\cup(\supp\zf\cap(Q,1])\subset\Cf$, where
$\Obs=\Omin$. Hence $V_\alpha(Q,\mu)\ge V_{\hat\zeta}(Q,\mu)$ for every $\alpha\in\Pspace([Q,1])$, and
$\hat\zeta$ is the minimizer. Finally $V_{\zf}=V_{\hat\zeta}$ at $(Q,\mu)$, by the remark
after~\eqref{eq:matrix-split}.
\end{proof}

For a configuration whose magnetization law is a perturbation of $\Lq{Q}$, the lemma bounds the
excess of $F_{\zf}$ over $\FTAP$. Appendix~\ref{app:kr} uses this for the tilted laws of the
Kac--Rice construction; the total variation distance below is $d_{\mathrm{TV}}(\pi,\pi')=\sup_A|\pi(A)-\pi'(A)|$.

\begin{corollary}[The fixed order parameter is nearly optimal near a contact]\label{cor:matrix-tilted}
Let $Q\in\Cf$, $M=\Phi_{\zf,x}(Q,\cdot)$ and $\pi_0=\Law(X_Q)$. For every $x\in\cube$ and every
probability measure $\pi$ on $\R$,
\[
 0\le\frac{F_{\zf}(x)-\FTAP(x)}N\le4C_*\,d_{\mathrm{TV}}(\pi,\pi_0)+2L_1|q_x-Q|+2L_2\,\varrho\bigl(\Wone(\mu_x,M_\#\pi)\bigr).
\]
If moreover $x_i=M(b_i)$ for all $i$, $\int M^2\,d\pi=Q$, and the empirical law of the $b_i$ is
within quadratic Wasserstein distance $w\le2$ of $\pi$, then $\Wone(\mu_x,M_\#\pi)\le w$ and
$|q_x-Q|\le2w$, so the right side is at most $4C_*d_{\mathrm{TV}}(\pi,\pi_0)+(4L_1+2L_2)\varrho(w)$.
\end{corollary}

\begin{proof}
By~\eqref{eq:matrix-split}, $F_{\zf}(x)-\FTAP(x)=N[V_{\zf}(q_x,\mu_x)-V_{\mathrm{TAP}}(q_x,\mu_x)]\ge0$.
Lemma~\ref{lem:matrix-correction}\ref{lem:matrix-correction-b}, applied to $V_{\zf}$ and to
$V_{\mathrm{TAP}}$, moves $(q_x,\mu_x)$ to $(Q,M_\#\pi)$ at cost
$2L_1|q_x-Q|+2L_2\varrho(\Wone(\mu_x,M_\#\pi))$. By
Lemma~\ref{lem:matrix-correction}\ref{lem:matrix-correction-a},
$|V_\beta(Q,M_\#\pi)-V_\beta(Q,\Lq{Q})|\le2C_*d_{\mathrm{TV}}(M_\#\pi,M_\#\pi_0)\le2C_*d_{\mathrm{TV}}(\pi,\pi_0)$
for every $\beta$, and the same holds for $V_{\mathrm{TAP}}$. Since
$V_{\zf}(Q,\Lq{Q})=V_{\mathrm{TAP}}(Q,\Lq{Q})$ by Lemma~\ref{lem:matrix-band}, the first claim
follows. For the second, $M$ is $1$-Lipschitz because $0<\Phi_{\zf,xx}\le1$, so
$\Wone(\mu_x,M_\#\pi)\le\Wone(N^{-1}\sum_i\delta_{b_i},\pi)\le w$; and $|m^2-m'^2|\le2|m-m'|$ on
$[-1,1]$ gives $|q_x-Q|\le2w$. Finally $w\le\varrho(w)$ on $[0,2]$.
\end{proof}

\begin{lemma}[Near-minimizers of the band problem]\label{lem:matrix-own}
For every $\tau>0$ there is $\gamma>0$ such that, for all $Q\in\Cf$, $r\in[0,1]$,
$\mu\in\Pspace([-1,1])$ and $\alpha\in\Pspace([0,1])$ with
\[
 |r-Q|\le\gamma,\qquad\Wone(\mu,\Lq{Q})\le\gamma,\qquad V_\alpha(r,\mu)\le V_{\mathrm{TAP}}(r,\mu)+\gamma ,
\]
one has $\int_r^1|\cdf\alpha(t)-\cdf{\zf}(t)|\,dt\le\tau$.
\end{lemma}

\begin{proof}
Suppose not. Then there are $\tau>0$ and, for every $n$, $Q_n\in\Cf$ and $(r_n,\mu_n,\alpha_n)$
satisfying the three conditions with $\gamma=1/n$ but
$\int_{r_n}^1|\cdf{\alpha_n}-\cdf{\zf}|>\tau$. Passing to a subsequence, $Q_n\to Q\in\Cf$, because
$\Cf$ is compact (Lemma~\ref{lem:loc-obstacle}\ref{lem:loc-obstacle-c}), and $\alpha_n\to\alpha$
weakly, so that $\cdf{\alpha_n}\to\cdf\alpha$ in $L^1$. Then $r_n\to Q$ and, since
$t\mapsto\Lq{t}$ is $\Wone$-continuous (Lemma~\ref{lem:pde-standing}\ref{lem:pde-standing-d}),
$\Wone(\mu_n,\Lq{Q})\to0$. By Lemma~\ref{lem:matrix-correction}\ref{lem:matrix-correction-b}--\ref{lem:matrix-correction-c},
$V_{\alpha_n}(r_n,\mu_n)\to V_\alpha(Q,\Lq{Q})$ and
$V_{\mathrm{TAP}}(r_n,\mu_n)\to V_{\mathrm{TAP}}(Q,\Lq{Q})$. Hence $V_\alpha(Q,\Lq{Q})\le V_{\mathrm{TAP}}(Q,\Lq{Q})$.
Since $V_\alpha(Q,\cdot)=V_{\Shift_Q\alpha}(Q,\cdot)$, the measure $\Shift_Q\alpha$ minimizes the band
problem, and Lemma~\ref{lem:matrix-band} gives $\Shift_Q\alpha=\Shift_Q\zf$, that is,
$\cdf\alpha=\cdf{\zf}$ on $[Q,1]$. Therefore
$\int_{r_n}^1|\cdf{\alpha_n}-\cdf{\zf}|\le\int_0^1|\cdf{\alpha_n}-\cdf\alpha|+|r_n-Q|\to0$, a
contradiction.
\end{proof}

\begin{proof}[Proof of Proposition~\ref{prop:hess-own}]
Let $\kappa>0$, let $\gamma>0$ be given by Lemma~\ref{lem:matrix-own} with $\tau=\kappa/2$, and
choose $\kappa_1\in(0,\min\{\kappa,\gamma/2,r_0\}]$ such that $\Wone(\Lq{Q},\Lq{Q'})\le\gamma/2$
whenever $|Q-Q'|\le\kappa_1$; this is possible because $t\mapsto\Lq{t}$ is uniformly continuous on
$[0,1]$ (Lemma~\ref{lem:pde-standing}\ref{lem:pde-standing-d}). Let $\delta_1,c>0$ be such that
both parts \ref{thm:struct-i} and~\ref{thm:struct-ii} of Theorem~\ref{thm:struct} hold at accuracy
$\kappa_1$ with these constants, and put $\varepsilon=\min\{\delta_1,\gamma/2\}$. For large $N$,
outside the union of the two exceptional events, of outer probability at most
$2e^{N(\Sigma(f)-c)}\le e^{N(\Sigma(f)-c/2)}$, let $x\in\cube$ and $\alpha$ satisfy
$\FTAP(x)\ge N(f-\varepsilon)$ and $F_\alpha(x)\le\FTAP(x)+2\varepsilon N$. Then
$\dist(q_x,\Cf)<\kappa_1\le\kappa$ and $\Wone(\mu_x,\Lq{q_x})<\kappa_1\le\kappa$. Pick $Q\in\Cf$
with $|q_x-Q|<\kappa_1\le\gamma$; since $\kappa_1\le r_0$, also $q_x\in(r_0,1-r_0)$. Then
$\Wone(\mu_x,\Lq{Q})<\kappa_1+\gamma/2\le\gamma$, and by~\eqref{eq:matrix-split}
$V_\alpha(q_x,\mu_x)-V_{\mathrm{TAP}}(q_x,\mu_x)=(F_\alpha(x)-\FTAP(x))/N\le2\varepsilon\le\gamma$.
Lemma~\ref{lem:matrix-own} gives $\int_{q_x}^1|\cdf\alpha-\cdf{\zf}|\le\kappa/2<\kappa$. The
proposition holds with $\varepsilon$ and $c/2$.
\end{proof}

\subsection{The index form and the index window}\label{app:matrix-index}

We prove the index form of the bulk, Proposition~\ref{prop:hess-index} and
Corollary~\ref{cor:hess-marg}, with $\kappa_*$, $\Cf^\kappa$, $\kappa'(\kappa)$
of~\eqref{eq:hess-kappaprime} and $\mathcal H_\varepsilon(x)$ as in Section~\ref{sec:hess-proof}.
Both turn the L\'evy bound of Theorem~\ref{thm:struct}\ref{thm:struct-iii} into eigenvalue counts,
and the square-root edge of Lemma~\ref{lem:hess-edge}\ref{lem:hess-edge-d} gives the order
$\kappa^{3/2}$ of the index. We then show that one parameter cannot serve for both index and window;
that proof uses the construction of Appendix~\ref{app:kr} through
properties~\ref{step:stab-k1}--\ref{step:stab-k3} of Section~\ref{sec:stab-kr}.

\begin{proof}[Proof of Proposition~\ref{prop:hess-index}]
\ref{prop:hess-index-a} The map $Q\mapsto\nuQ$ is weakly continuous on $(0,1)$: $\aQ$ and
$\chiQ=\E\UQ$ are continuous, $\Law(\UQ)$ is weakly continuous
(Lemma~\ref{lem:pde-standing}\ref{lem:pde-standing-d}) and
$\UQ>0$, so $\rhoQ$ is weakly continuous, and Lemma~\ref{lem:hess-pastur}\ref{lem:hess-pastur-c}
applies. The edge $\lamp$ is continuous by Lemma~\ref{lem:hess-edge}\ref{lem:hess-edge-e}. Let
$Q_n\to Q$ in $\Cf^\kappa$. For large $n$,
$(\lamp(Q_n)-\kappa/2,\infty)\supset(\lamp(Q)-\kappa/4,\infty)$, and the portmanteau theorem gives
$\liminf_n\nu_{Q_n}((\lamp(Q_n)-\kappa/2,\infty))\ge\nuQ((\lamp(Q)-\kappa/4,\infty))>0$, the last
because $\lamp(Q)$ lies in the support. By compactness of $\Cf^\kappa$ the infimum
in~\eqref{eq:hess-kappaprime} is positive. For the order of magnitude, the constants $c,C,\kappa_0$
of Lemma~\ref{lem:hess-edge}\ref{lem:hess-edge-d} for $\nuQ$ depend continuously on $Q$ and are
positive, so there are $c_\flat,C_\flat,\kappa_\flat>0$ with $c\ge c_\flat$, $C\le C_\flat$ and
$\kappa_0\ge\kappa_\flat$ for every $Q$ in the compact set $\Cf^{\kappa_*}$. Let
$\kappa_{\mathrm{ix}}=\min(\kappa_*,2\kappa_\flat)$. Since $\nuQ$ is
carried by $(-\infty,\lamp(Q)]$, for $\kappa\le\kappa_{\mathrm{ix}}$ and $Q\in\Cf^\kappa$,
\[
 c_\flat(\kappa/4)^{3/2}\le\nuQ\bigl([\lamp(Q)-\kappa/4,\lamp(Q)]\bigr)
 \le\nuQ\bigl((\lamp(Q)-\kappa/2,\infty)\bigr)\le C_\flat(\kappa/2)^{3/2},
\]
which gives~\ref{prop:hess-index-a} with $c_{\mathrm{ix}}=c_\flat4^{-3/2}/2$ and
$C_{\mathrm{ix}}=C_\flat2^{-3/2}/2$.

\ref{prop:hess-index-b} Let $\kappa'=\kappa'(\kappa)$ and $\eta=\min(\kappa/2,\kappa'/2)$, and
let $\varepsilon,c$ be given by Theorem~\ref{thm:struct}\ref{thm:struct-iii} at accuracy $\eta$.
Outside its exceptional event, let $x\in\ocube$ with $\FTAP(x)\ge N(f-\varepsilon)$. Then
$\dist(q_x,\Cf)<\eta\le\kappa$, so $q_x\in\Cf^\kappa$. Write $\lamp=\lamp(q_x)$ and
$\nu=\nu_{q_x}$, and let $H\in\mathcal H_\varepsilon(x)$ be one of the one-sided Hessians, so that
$\dL(\ESD(H),\nu)\le\eta$. Since $\nu$ is carried by $(-\infty,\lamp]$ and $\eta<\kappa$, the
right inequality of~\eqref{eq:hess-levy} gives
\[
 \ESD(H)\bigl((\lamp+\kappa,\infty)\bigr)\le\nu\bigl([\lamp+\kappa-\eta,\infty)\bigr)+\eta
 =\eta<\kappa' .
\]
Since $\eta\le\kappa/2$ and $q_x\in\Cf^\kappa$, the left inequality gives
\[
 \ESD(H)\bigl((\lamp-\kappa,\infty)\bigr)\ge\nu\bigl((\lamp-\kappa/2,\infty)\bigr)-\eta
 \ge2\kappa'-\kappa'/2=\tfrac32\kappa' .
\]
The eigenvalue counts are integers. So $\#\{i:\lambda_i(H)>\lamp-\kappa\}>\kappa'N$ gives
$\lambda_{\lceil\kappa'N\rceil}(H)>\lamp-\kappa$, and
$\#\{i:\lambda_i(H)>\lamp+\kappa\}<\kappa'N\le\lceil\kappa'N\rceil$ gives
$\lambda_{\lceil\kappa'N\rceil}(H)\le\lamp+\kappa$. For a compression $H'$ of such an $H$, the
interlacing inequalities before~\eqref{eq:hess-interlace} change each count by at most one, and
cannot increase the count above $\lamp+\kappa$. The count above $\lamp-\kappa$ remains at least
$\frac32\kappa'N-1\ge\lceil\kappa'N\rceil$ once $\kappa'N\ge4$, and the same conclusions follow.
\end{proof}

\begin{proof}[Proof of Corollary~\ref{cor:hess-marg}]
By Lemma~\ref{lem:hess-edge}\ref{lem:hess-edge-e}, $\lamp$ is uniformly continuous on the compact
set $\Cf^{\kappa_*}$. Take $\kappa_{\rm loc}\in(0,\kappa/4]$ with
$|\lamp(Q)-\lamp(Q')|\le\kappa/2$ whenever $Q,Q'\in\Cf^{\kappa_*}$ and
$|Q-Q'|\le\kappa_{\rm loc}$. Let $\kappa'=\kappa'(\kappa/2)$, and apply
Theorem~\ref{thm:struct}\ref{thm:struct-iii} at accuracy
$\eta=\min(\kappa/4,\kappa'/2,\kappa_{\rm loc})$; this gives $\varepsilon$ and $c$. Outside its
exceptional event, let $x$ be as in the statement, and pick $Q\in\Cf$ with $|q_x-Q|<\eta$. The
proof of Proposition~\ref{prop:hess-index}\ref{prop:hess-index-b}, run with $\kappa/2$ in place of
$\kappa$, uses exactly $\eta\le\kappa/4$, $\eta\le\kappa'(\kappa/2)/2$ and
$q_x\in\Cf^{\kappa/2}$. It shows that every $H\in\mathcal H_\varepsilon(x)$ has
$\lambda_{\lceil\kappa'N\rceil}(H)\in[\lamp(q_x)-\kappa/2,\lamp(q_x)+\kappa/2]$ and fewer than
$\kappa'N$ eigenvalues above $\lamp(q_x)+\kappa/2$, for large $N$. Since
$|\lamp(q_x)-\lamp(Q)|\le\kappa/2$, the interval lies in $[\lamp(Q)-\kappa,\lamp(Q)+\kappa]$, and
$\lamp(Q)+\kappa\ge\lamp(q_x)+\kappa/2$. The description of $\lamp(Q)$ is
Lemma~\ref{lem:hess-edge}\ref{lem:hess-edge-e}.
\end{proof}

\paragraph{The index window.}
Proposition~\ref{prop:hess-index} controls the $\lceil\kappa'N\rceil$-th eigenvalue at accuracy
$\kappa$ only with $\kappa'$ of order $\kappa^{3/2}$. Remark~\ref{rem:stab-window} in
Section~\ref{sec:hess-proof} states that a single parameter for both the index and the window is
impossible; we prove it here, using the points retained by the Kac--Rice construction.

We use the notation of Section~\ref{sec:stab-kr}: $q=\qf$, $a=\xi''(q)$, $M=\Phi_{\zf,x}(q,\cdot)$,
$U=\Phi_{\zf,xx}(q,\cdot)$, $\pi_0=\Law(X_q)$, the tilts $\pi_\eta$ of~\eqref{eq:stab-tilt}, the
admissible tilts, the marked counts $\widehat C_N(\eta,\delta_c,w)$ with overlap cutoff $\delta_c$,
and the properties~\ref{step:stab-k1}--\ref{step:stab-k3} of the retained points. For $x\in\ocube$
let $\mathcal H^{\mathrm{fix}}(x)$ consist of $\Hpm F_{\zf}(x)$ and their (common) compression to
$x^\perp$.

\begin{proposition}[The index window cannot be widened]\label{prop:matrix-window}
Assume the hypotheses of Theorem~\ref{thm:stab}\ref{thm:stab-a} and let $f\in(\feq,\fone)$. There is
$\kappa_{\mathrm w}\in(0,\frac12]$ with the following property. For $\kappa\in(0,\kappa_{\mathrm w}]$ and $\delta>0$ let
$\mathcal B_N(\kappa,\delta)$ be the event that some $x\in\ocube$ satisfies
\[
 \nabla F_{\zf}(x)=0,\qquad\FTAP(x)\ge N(f-\delta),\qquad
 \lambda_{\lfloor\kappa N\rfloor}(H)<\lamp(q_x)-\kappa\ \text{ for every } H\in\mathcal H^{\mathrm{fix}}(x).
\]
Then for every $c'>0$ and all large $N$, $\mathcal B_N(\kappa,\delta)$ contains an event of
probability at least $e^{N(\Sigma(f)-c')}$. Consequently, for no $\kappa\in(0,\kappa_{\mathrm w}]$ are there
$\delta,c>0$ such that, for all large $N$, outside an event of outer probability at most
$e^{N(\Sigma(f)-c)}$, every $m\in\ocube$ with $\FTAP(m)\ge N(f-\delta)$ has
$\lambda_{\lfloor\kappa N\rfloor}(H)\ge\lamp(q_m)-\kappa$ for some $H\in\mathcal H^{\mathrm{fix}}(m)$.
\end{proposition}

\begin{proof}
The parameters are chosen in the order: $\kappa_{\mathrm w}$; then $\kappa$, $\delta$ and $c'$; then the tilt
$\eta$; then the cutoff $\delta_c$; then the window $w$; then $N$. The point is that $\kappa_{\mathrm w}$ must
not depend on $\eta$, and Lemma~\ref{lem:matrix-tilt} provides this.

\emph{Step 0: uniform edge constants.} If $\Plef(q)>0$ only $\eta=0$ is admissible; set
$\mathcal R=0$. If $\Plef(q)=0$, let $\mathcal R$ and $\eta_0=(2\|\mathcal R\|_\infty)^{-1}$ be as
in Section~\ref{sec:stab-kr}, so that every admissible $\eta$ lies in $(0,\eta_0]$. Under $\pi_0$
the variable $U$ has the law of $U_q$, so Lemma~\ref{lem:matrix-hyp} gives the hypothesis of
Lemma~\ref{lem:matrix-tilt} with $P_0=\pi_0$. In the notation of Section~\ref{sec:stab-kr},
$\nu^{(\eta)}=\Pastur{a}{\Law_{\pi_\eta}(U^{-1}+a\chi_\eta)}$, and
Lemma~\ref{lem:matrix-tilt}\ref{lem:matrix-tilt-c} gives $C_\sharp,\kappa_\sharp>0$ with
\[
 \nu^{(\eta)}\bigl([\lamp^{(\eta)}-y,\lamp^{(\eta)}]\bigr)\le C_\sharp y^{3/2}
 \qquad(\eta\in[0,\eta_0],\ 0<y\le\kappa_\sharp).
\]
Choose $\kappa_{\mathrm w}\in(0,\frac12]$ with $\frac32\kappa_{\mathrm w}\le\kappa_\sharp$ and
$C_\sharp(\frac32)^{3/2}\kappa_{\mathrm w}^{1/2}\le\frac14$.

\emph{Step 1: the bulk at the retained points.} The law $M_\#\pi_\eta$ puts no mass at $\pm1$, and
$u_{\zf}(q,M(b))=U(b)$; so, in the notation of Proposition~\ref{prop:hess-bulk},
$\Pastur{a}{\rho_{q,M_\#\pi_\eta}}=\nu^{(\eta)}$. Moreover $\lamp^{(\eta)}\le\lamp(q)$: for
$\eta=0$ they are equal, and if $\Plef(q)=0$ then $\lamp(q)=0$ while $\lamp^{(\eta)}\le0$, both by
Lemma~\ref{lem:hess-edge}\ref{lem:hess-edge-b}, applied under $\pi_\eta$ with
Lemma~\ref{lem:matrix-tilt}\ref{lem:matrix-tilt-a}.

\emph{Step 2: parameters.} Fix $\kappa\in(0,\kappa_{\mathrm w}]$, $\delta>0$ and $c'\in(0,1)$.
\begin{itemize}
\item Choose an admissible $\eta$ with $\varepsilon_3(\eta)\le c'/4$ and
 $\varepsilon_1(\eta)\le\delta/4$, in the notation of
 \ref{step:stab-k2}--\ref{step:stab-k3}.
\item Choose $\delta_c\in(0,q/2)$ with $\frac12\log(1-\delta_c/q)\ge-c'/4$. By \ref{step:stab-k3}
 there is $w_2>0$ such that
 $\liminf_NN^{-1}\log\Prob(\widehat C_N(\eta,\delta_c,w)\ge1)\ge\Sigma(f)-3c'/4$ for every
 $w\in(0,w_2]$.
\item Apply Proposition~\ref{prop:hess-bulk} to the one-point set
 $\mathcal K=\{(q,M_\#\pi_\eta)\}$, which lies in $[r_0,1-r_0]\times\Pspace([-1,1])$ because
 $q\in\Cf$, with accuracy $\kappa/4$ and exponent $2-\Sigma(f)>0$. This gives $w_H>0$ and, for
 large $N$, a measurable set $E_N$ of probability at most $2e^{N(\Sigma(f)-2)}$, containing the
 exceptional set of the proposition (the union of $\Omega_2^c$, $\Omega_3^c$ and the event of
 Step~3 of Lemma~\ref{lem:matrix-uniform}), outside of which the conclusion of the proposition
 holds.
\item By Lemma~\ref{lem:hess-edge}\ref{lem:hess-edge-e} choose $w_\lambda>0$ with
 $|\lamp(q')-\lamp(q)|\le\kappa/4$ whenever $|q'-q|\le2w_\lambda$.
\item Fix $w\in(0,\min\{w_2,w_H/2,w_\lambda\}]$ with $\varepsilon_2(w)\le\delta/4$.
\end{itemize}

\emph{Step 3: the retained points.} By \ref{step:stab-k1}--\ref{step:stab-k2}, every point $x$
retained by $\widehat C_N(\eta,\delta_c,w)$ lies in $\ocube$, satisfies $\nabla F_{\zf}(x)=0$,
$\Wone(\mu_x,M_\#\pi_\eta)\le w$ and $|q_x-q|\le2w$, and
$|N^{-1}\FTAP(x)-f|\le\delta/2+2N^{-2}\le\delta$ for large $N$.

\emph{Step 4: the spectrum at the retained points.} On $\{\widehat C_N\ge1\}\setminus E_N$ let $x$
be a retained point. Since $|q_x-q|\le w_H$ and $\Wone(\mu_x,M_\#\pi_\eta)\le w_H$,
Proposition~\ref{prop:hess-bulk} with $\alpha=\zf$ and Step 1 give
$\dL(\ESD(\Hpm F_{\zf}(x)),\nu^{(\eta)})\le\kappa/4$. By the choice of $w_\lambda$ and Step 1,
$\lamp(q_x)-\kappa\ge\lamp(q)-\frac54\kappa\ge\lamp^{(\eta)}-\frac54\kappa$. If
$\dL(\mu_1,\mu_2)\le\epsilon$, then $\mu_1([t,\infty))\le\mu_2([t-\epsilon,\infty))+\epsilon$; and
$\nu^{(\eta)}$ is carried by $(-\infty,\lamp^{(\eta)}]$. Hence, by Step 0 and
$\frac32\kappa\le\kappa_\sharp$,
\[
 \ESD\bigl(\Hpm F_{\zf}(x)\bigr)\bigl([\lamp(q_x)-\kappa,\infty)\bigr)
 \le\nu^{(\eta)}\bigl([\lamp^{(\eta)}-\tfrac32\kappa,\lamp^{(\eta)}]\bigr)+\frac\kappa4
 \le C_\sharp\bigl(\tfrac32\kappa\bigr)^{3/2}+\frac\kappa4\le\frac\kappa2 .
\]
So at most $\kappa N/2$ eigenvalues of either one-sided Hessian are at least $\lamp(q_x)-\kappa$.
For $\kappa N>2$ we have $\kappa N/2<\lfloor\kappa N\rfloor$, hence
$\lambda_{\lfloor\kappa N\rfloor}(\Hpm F_{\zf}(x))<\lamp(q_x)-\kappa$. Cauchy interlacing gives
$\lambda_i(H')\le\lambda_i(H)$ for the compression $H'$ and $i\le N-1$, and
$\lfloor\kappa N\rfloor\le N-1$, so the same holds for the compression. With Step 3, $x$ witnesses
$\mathcal B_N(\kappa,\delta)$.

\emph{Step 5: probability.} The event $\{\widehat C_N\ge1\}\setminus E_N$ is contained in
$\mathcal B_N(\kappa,\delta)$, and by Step~2 its probability is at least
$e^{N(\Sigma(f)-7c'/8)}-2e^{N(\Sigma(f)-2)}\ge e^{N(\Sigma(f)-c')}$ for large $N$, since $c'<1$.
This proves the first claim for $c'\in(0,1)$, and hence for every $c'>0$.

For the consequence, suppose that for some $\kappa\le\kappa_{\mathrm w}$ there were such $\delta$ and $c$.
Every witness of $\mathcal B_N(\kappa,\delta)$ violates the property in question, so the event of
the first claim with $c'=c/2$, of probability at least $e^{N(\Sigma(f)-c/2)}$, would be contained in
an event of outer probability at most $e^{N(\Sigma(f)-c)}$, which is impossible for large $N$.
\end{proof}

The obstruction is the square-root edge: uniformly in the tilt, $\nu^{(\eta)}$ puts mass of order
$y^{3/2}$ within distance $y$ of its edge, so a window of width $\kappa$ below the edge holds only
about $\kappa^{3/2}N$ eigenvalues, fewer than $\kappa N$. The argument uses only the construction
of Appendix~\ref{app:kr}, Proposition~\ref{prop:hess-bulk} and Lemmas~\ref{lem:hess-edge},
\ref{lem:matrix-hyp} and~\ref{lem:matrix-tilt}, at the fixed order parameter $\zf$.

%% file: sections/appE-kacrice.tex
\section{The Kac--Rice construction}\label{app:kr}

This appendix proves Proposition~\ref{prop:kr-main}: at every level $f\in(\feq,\fone)$ there are
critical points of the \emph{fixed} function $F_{\zf}$ at the rate $\Sigma(f)$, around each of
which $F_{\zf}$ is strongly concave, with a margin that depends on the accuracy
(Section~\ref{sec:stab-kr}), on a ball of radius of order $\sqrt N$.
Theorem~\ref{thm:stab}\ref{thm:stab-a} is read off from its proof (Section~\ref{sec:stab-kr}), and
for convex mixtures that are not even it gives the lower half of the fractional-moment formula
(Proposition~\ref{prop:max-fm}\ref{prop:max-fm-d}). The proof uses Appendix~\ref{app:pde},
Lemmas~\ref{lem:loc-obstacle}, \ref{lem:poisson-shift} and~\ref{lem:hess-pastur} and
Appendix~\ref{app:matrix}, and none of the main results; in particular it does not use \condFM.

We say that $\xi$ has \emph{two active degrees} if $\beta_p\neq0$ for at least two values of
$p$. A single monomial that is convex on $[-1,1]$ is even, so every convex $\xi$ that is not
even has two active degrees.

\begin{proposition}[Critical points of the fixed TAP function]\label{prop:kr-main}
Assume \condH, $\xi\in\classC$ and that $\xi$ has two active degrees, and let
$f\in(\feq,\fone)$. Then, for every $\varepsilon>0$,
\[
 \liminf_{N\to\infty}\frac1N\log\Prob\bigl(\Eex{\varepsilon}{f}\bigr)\ge\Sigma(f),
\]
and the witnesses of $\Eex{\varepsilon}{f}$ can be taken of the form $(\zf,x)$. More precisely,
the retained points of Appendix~\ref{app:kr-count} have the properties
\ref{step:stab-k1}--\ref{step:stab-k3} listed in Section~\ref{sec:stab-kr}.
\end{proposition}

\paragraph{The obstruction and the two adjustments.}
The proof counts critical points $x$ of $F_{\zf}$ whose self-overlap is near the first positive
contact $\qf=\min\Cf$ and whose coordinates, written as $x_i=\Phi_{\zf,x}(\qf,b_i)$, have
empirical field law $N^{-1}\sum_i\delta_{b_i}$ near the law $\pi_0$ of the Auffinger--Chen process
at time $\qf$. At such a point, conditionally on criticality, the Hessian of $F_{\zf}$
transverse to $x$ is $\sqrt{\xi''(\qf)}$ times a GOE matrix minus a diagonal matrix with entries
$\Phi_{\zf,xx}(\qf,b_i)^{-1}+\xi''(\qf)\chi$, where $\chi$ is the average of the
$\Phi_{\zf,xx}(\qf,b_i)$. Its spectral gap is governed by the Plefka factor $\Plef(\qf)$. If $\Plef(\qf)>0$, the gap is strict, and the isolation and spherical-code
argument of Boursier \cite[Section~3]{Bou26b}, which adapts the isolated-critical-point argument
of Huang and Sellke \cite{HS}, applies. If $\Plef(\qf)=0$, the bulk of that matrix reaches zero
(Theorem~\ref{thm:struct}\ref{thm:struct-iv}), and neither the local strong concavity nor the
limit of the normalized log-determinant can be read off.

We change the data selected inside the Kac--Rice integral, not the function. First, a bounded
tilt $\pi_\eta=(1-\eta\mathcal R)\pi_0$ of the field law keeps the self-overlap
$\E_{\pi_\eta}\Phi_{\zf,x}(\qf,b)^2=\qf$ and makes the Plefka factor of the selected data
positive, of order $\eta$ (Lemma~\ref{lem:kr-tilt}). The tilted law no longer satisfies the
integration-by-parts identity of $\pi_0$, which is what makes the radial part of the Hessian and
the bound on competitors at smaller overlap vanish to the right order. Both acquire errors of
order $\eta$. Second, with two active degrees the value $\HN(x)$ keeps a conditional variance of
order $N$ after its gradient is fixed (Lemma~\ref{lem:kr-energy}), so the selected energy can be
raised by $L\eta N$ at a cost of order $\eta$ in the rate. A large fixed $L$ makes both the radial
Schur complement and the competitor bound strictly negative, of order $\eta$
(Appendices~\ref{app:kr-hessian} and~\ref{app:kr-gap}). Then retained critical points are isolated:
two of them have overlap less than a cutoff $\delta$, so their number is bounded by a spherical
code, and dividing the expected number by this deterministic bound bounds the existence
probability from below (Appendix~\ref{app:kr-count}). The rate obtained is
$\Sigma(f)-O(\eta)+\frac12\log(1-\delta/\qf)$, and we let $\delta\downarrow0$ and then
$\eta\downarrow0$. If $\Plef(\qf)>0$, the construction runs at $\eta=0$ and is Boursier's.

Appendix~\ref{app:kr-tilt} constructs the two adjustments, \ref{app:kr-hessian}
and~\ref{app:kr-gap} prove the Hessian bound and the competitor gap, \ref{app:kr-count} counts
and proves Proposition~\ref{prop:kr-main}, \ref{app:kr-susy} proves the margin bound and the SUSY
count of Section~\ref{sec:stab-kr}, and \ref{app:kr-remarks} bounds the retained states given
existence and locates the hypotheses on $\xi$.

\paragraph{Standing assumptions and notation.}
Throughout this appendix $\xi$ satisfies \condH, $\xi\in\classC$, $\xi$ has two active degrees,
and $f\in(\feq,\fone)$ is fixed. We write
\begin{gather*}
 \theta=\thf,\qquad \zeta=\zf,\qquad q=\qf,\qquad \Phi=\Phi_\zeta,\\
 a=\xi''(q),\qquad c=\xi'(q),\qquad a_1=\xi'''(q),\qquad a_2=\xi''''(q),
\end{gather*}
and
\[
 M(b)=\Phi_x(q,b),\qquad U(b)=\Phi_{xx}(q,b),\qquad \pi_0=\Law(X_q),
\]
where $X$ is the Auffinger--Chen process of $\zeta$. For a law $\pi$ on $\R$, $\E_\pi$ denotes
expectation in $b\sim\pi$, and $\E_0=\E_{\pi_0}$. The Plefka factor at the first contact is
$\Plef_0=\Plef(q)=1-a\E_0U^2$. For $x\in\ocube$ we write $\pi_x=N^{-1}\sum_i\delta_{M^{-1}(x_i)}$ for
the empirical law of its \emph{fields} at the fixed radius $q$ (Lemma~\ref{lem:kr-contact} below
shows that $M$ is a bijection of $\R$ onto $(-1,1)$), and
$\mu_x=N^{-1}\sum_i\delta_{x_i}=M_\#\pi_x$.

\subsection{The field law at the first contact, its tilt and the energy shift}\label{app:kr-tilt}

We first collect what the construction uses about the constrained minimizer and its first
contact. All of it is proved in Appendix~\ref{app:pde} or in Section~\ref{sec:loc}.

\begin{lemma}[The first contact]\label{lem:kr-contact}
\leavevmode
\begin{thmlist}
\item\label{lem:kr-contact-a} $\theta\in(0,1)$, $\zeta(\{0\})=\theta$ and $\zeta((0,q))=0$, so that
 $\cdf\zeta(t)=\theta$ for $t\in[0,q)$. Moreover $\Sigma(f)=\theta(P_\zeta-f)$.
\item\label{lem:kr-contact-b} $q\in(0,1)$, $\Obs_\zeta(r)>\Obs_\zeta(q)$ for $r\in(0,q)$,
 $\Gz\zeta(q)=q$, and $\Plef_0=1-a\E_0U^2=2\Obs_\zeta''(q)/a\ge0$.
\item\label{lem:kr-contact-c} \textup{(Level identity.)}
 $f=P_\zeta+\theta\bigl(\Obs_\zeta(0)-\Obs_\zeta(q)\bigr)$.
\item\label{lem:kr-contact-d} \textup{(Regularity.)} For every $t\in[0,1]$, $\Phi_x(t,\cdot)$ is an
 odd, strictly increasing bijection of $\R$ onto $(-1,1)$, and $0<\Phi_{xx}\le1$. There are
 $0<c_1\le C_1$ and $C_{34}$, depending only on $\xi'(1)$, with
 \[
  c_1e^{-2|b|}\le\Phi_{xx}(t,b)\le C_1e^{-2|b|},\qquad
  |\Phi_{xxx}(t,b)|+|\Phi_{xxxx}(t,b)|\le C_{34}\Phi_{xx}(t,b),
 \]
 for all $(t,b)\in[0,1]\times\R$.
\item\label{lem:kr-contact-e} \textup{(Plateau law.)} The law $\pi_0$ has the density
 \begin{equation}\label{eq:kr-law}
  \frac{d\pi_0}{d\mathcal N(0,c)}(b)=\exp\bigl\{\theta\bigl[\Phi(q,b)-\Phi(0,0)\bigr]\bigr\}
 \end{equation}
 with respect to the centred Gaussian law of variance $c=\xi'(q)$. It is positive and has Gaussian
 tails, and for every $C^1$ function $F$ with $|F|+|F'|\le C(1+e^{C|b|})$,
 \begin{equation}\label{eq:kr-ward}
  \E_0\bigl[bF(b)/c\bigr]=\E_0\bigl[F'(b)+\theta M(b)F(b)\bigr].
 \end{equation}
\end{thmlist}
\end{lemma}

\begin{proof}
\ref{lem:kr-contact-a}--\ref{lem:kr-contact-b} By
Lemma~\ref{lem:loc-obstacle}\ref{lem:loc-obstacle-a}--\ref{lem:loc-obstacle-b}, $\theta\in(0,1)$,
$\zeta(\{0\})=\theta$ and $\Sigma(f)=\theta(P_\zeta-f)$, and $\zeta$ is unique
(Lemma~\ref{lem:var-basic}\ref{lem:var-basic-a}). By Lemma~\ref{lem:loc-obstacle}\ref{lem:loc-obstacle-c}, $\Cf$ is a
compact subset of $(0,1)$ containing $\supp\zeta\cap(0,1]$, and $\Obs_\zeta>\Omin$ on
$(0,1]\setminus\Cf$. Since $q=\min\Cf$, this gives $\zeta((0,q))=0$ and
$\Obs_\zeta(r)>\Omin=\Obs_\zeta(q)$ for $0<r<q$. The identities $\Gz\zeta(q)=q$ and
$\Plef_0=2\Obs_\zeta''(q)/a\ge0$ are Lemma~\ref{lem:pde-standing}\ref{lem:pde-standing-c}.

\ref{lem:kr-contact-c} This is Remark~\ref{rem:var-level}: by
Lemma~\ref{lem:var-basic}\ref{lem:var-basic-c}, $\Lambda'=p+\theta p'$ with $p(\theta)=P_\zeta$ and
$p'(\theta)=\bigl(\Obs_\zeta(0)-\int\Obs_\zeta\,d\zeta\bigr)/(1-\theta)=\Obs_\zeta(0)-\Obs_\zeta(q)$, since
$\zeta(\{0\})=\theta$ and $\Obs_\zeta=\Obs_\zeta(q)$ on $\supp\zeta\cap(0,1]$; and $\Lambda'(\theta)=f$
by Proposition~\ref{prop:var-regular}\ref{prop:var-regular-b}.

\ref{lem:kr-contact-d} This is Lemma~\ref{lem:pde-standing}\ref{lem:pde-standing-a}--\ref{lem:pde-standing-b},
together with $\Phi_{xx}\le1$, which holds for finitely supported order parameters by
Lemma~\ref{lem:pde-class}\ref{lem:pde-class-a} with $G=1$ and passes to the limit by the uniform
convergence in the proof of Lemma~\ref{lem:pde-standing}.

\ref{lem:kr-contact-e} On $[0,q)$ the coefficient $\cdf\zeta$ equals $\theta$, so
$w(t,b)=e^{\theta\Phi(t,b)}$ solves the backward heat equation $\partial_tw+\frac{\xi''}2w_{xx}=0$
there (the Cole--Hopf form \eqref{eq:pde-colehopf} of the Parisi equation on a step). On $[0,q]$
the Auffinger--Chen process, whose drift is $\theta\xi''\Phi_x=\xi''\partial_x\log w$, is the
Doob transform of the Gaussian martingale $\int_0^t\sqrt{\xi''}\,dB$ by $w$. Its law at time $q$
therefore has density $w(q,b)/w(0,0)$ with respect to the law of that martingale at time $q$,
which is $\mathcal N(0,\xi'(q))$ because $\xi'(0)=0$. This is \eqref{eq:kr-law}, the plateau
Cole--Hopf identity of \cite[Lemma~A.3, (A.11)]{Bou26b} with lower endpoint zero. Since
$|\Phi_x|<1$, $|\Phi(q,b)-\Phi(q,0)|\le|b|$, so the density is positive and $\pi_0$ has Gaussian
tails. Differentiating the logarithm of the density gives $-b/c+\theta M(b)$, and
\eqref{eq:kr-ward} follows by integration by parts.
\end{proof}

\paragraph{Identities at the contact.}
Put $\chi_0=\E_0U$ and $k_0=\E_0\,bM(b)$. By \eqref{eq:kr-ward} with $F=M$ and
$\E_0M^2=\Gz\zeta(q)=q$,
\begin{equation}\label{eq:kr-k0}
 k_0=c\,(\chi_0+\theta q).
\end{equation}
Along the Auffinger--Chen process, It\^o's formula and the Parisi equation on the plateau give
$\frac d{dt}\E\Phi(t,X_t)=\frac\theta2\xi''(t)\Gz\zeta(t)$ for $0<t<q$. Recall
$\vth(q)=q\xi'(q)-\xi(q)=\int_0^qt\xi''(t)\,dt$. Since
$\Obs_\zeta(0)-\Obs_\zeta(q)=\frac12\int_0^q\xi''(t)(\Gz\zeta(t)-t)\,dt$, the level identity gives
\begin{equation}\label{eq:kr-ito}
 \E_0\Phi(q,b)-\Phi(0,0)=\frac\theta2\int_0^q\xi''(t)\Gz\zeta(t)\,dt
 =f-P_\zeta+\frac\theta2\vth(q).
\end{equation}
Define the \emph{energy density at the contact}
\begin{equation}\label{eq:kr-h0}
 h_0=c\chi_0+\theta\xi(q).
\end{equation}
Then, since $\cdf\zeta=\theta$ on $[0,q)$, the definition of $P_\zeta$, \eqref{eq:kr-k0} and
\eqref{eq:kr-ito} give
\begin{equation}\label{eq:kr-level}
 h_0+\E_0\bigl[\Phi(q,b)-bM(b)\bigr]-\frac12\int_q^1t\xi''(t)\cdf\zeta(t)\,dt=f .
\end{equation}
By the Fenchel identity $\Phistar(q,M(b))=\Phi(q,b)-bM(b)$, the middle term is
$\int\Phistar(q,m)\,\mu_0(dm)$ with $\mu_0=M_\#\pi_0$, which is the field law $\Lq{q}$. So
\eqref{eq:kr-level} says that a point
$x$ with $\HN(x)=Nh_0$ and $\mu_x=\mu_0$ has $F_\zeta(x)=Nf$; this is how $h_0$ is chosen.

\paragraph{The field-law tilt.}
The following lemma opens the Plefka inequality without changing the self-overlap.

\begin{lemma}[Field-law tilt]\label{lem:kr-tilt}
Assume $\Plef_0=0$. Let $\mathcal R$ be the residual of $U^2$ after orthogonal projection in
$L^2(\pi_0)$ onto $\operatorname{span}\{1,M^2\}$. Then $\mathcal R$ is bounded and not identically
zero, and $\E_0\mathcal R=\E_0\mathcal RM^2=0$ and $\E_0\mathcal RU^2=\E_0\mathcal R^2$. For
$0<\eta\le\eta_0:=(2\|\mathcal R\|_\infty)^{-1}$, the measure $\pi_\eta=(1-\eta\mathcal R)\pi_0$ is a
probability measure whose density with respect to $\pi_0$ lies in $[\frac12,\frac32]$, and
\[
 \E_{\pi_\eta}M^2=q,\qquad
 \Plef_\eta:=1-a\E_{\pi_\eta}U^2=a\eta\,\E_0\mathcal R^2>0 .
\]
\end{lemma}

\begin{proof}
The functions $1$ and $M^2$ are linearly independent in $L^2(\pi_0)$, because $M$ is strictly
increasing and $\pi_0$ has a positive density. All three functions $1,M^2,U^2$ are bounded, so the
projection has bounded coefficients and $\mathcal R$ is bounded. Orthogonality gives
$\E_0\mathcal R=\E_0\mathcal RM^2=0$ and $\E_0\mathcal RU^2=\E_0\mathcal R^2$, and these give the
normalization, $\E_{\pi_\eta}M^2=\E_0M^2=q$ and $\E_{\pi_\eta}U^2=\E_0U^2-\eta\E_0\mathcal R^2$.
Suppose $\mathcal R=0$ in $L^2(\pi_0)$. By continuity and the full support of $\pi_0$,
$U^2=\bar A+\bar BM^2$ at every field. As $b\to\pm\infty$, $M(b)\to\pm1$ and $U(b)\to0$
(Lemma~\ref{lem:kr-contact}\ref{lem:kr-contact-d}), so $\bar B=-\bar A$, and $\bar A>0$ since
$U>0$. Then $M'=U=\sqrt{\bar A(1-M^2)}$, whose solution reaches $1$ at a finite field, which
contradicts $|M|<1$.
\end{proof}

At $\Plef_0=0$ the tilted law $\pi_\eta$ has a bounded density with respect to $\pi_0$, hence
finite exponential moments of every order. We also put $\pi_\eta=\pi_0$ for $\eta=0$, and, for
every law $\pi$ with $\E_\pi M^2=q$,
\[
 \chi_\pi=\E_\pi U,\qquad k_\pi=\E_\pi\,bM(b),\qquad d_\pi=k_\pi+aq\chi_\pi,\qquad
 \upsilon_\pi=\E_\pi\bigl[UM^2\bigr]>0 .
\]
We abbreviate $\chi_\eta=\chi_{\pi_\eta}$, and similarly for $k$, $d$, $\upsilon$.

\paragraph{The energy shift.}
The second adjustment raises the selected value of $\HN$. It is possible because, with two
active degrees, the value of $\HN$ at a point is not determined by its gradient. Put
$s_j=\sum_pp^j\beta_p^2q^p$ for $j\ge0$, so that $s_0=\xi(q)$, $s_1=qc$ and $s_2=q(c+aq)$, and
\[
 S_0=\begin{pmatrix}s_0&s_1\\ s_1&s_2\end{pmatrix},\qquad
 A(r)=\frac{s_2\xi(r)-s_1r\xi'(r)}{\det S_0},\qquad
 B(r)=\frac{s_0r\xi'(r)-s_1\xi(r)}{\det S_0}.
\]
At a point $x$ with $q_x=q$, $(\HN(x),x\cdot\nabla\HN(x))$ has covariance $NS_0$, the projection of
$\nabla\HN(x)$ onto $x^\perp$ is independent of it, and for every $y$,
$N^{-1}\Cov(\HN(y),(\HN(x),x\cdot\nabla\HN(x)))=(\xi(r),r\xi'(r))$ with $r=\ov xy$. So the
regression of $\HN(y)$ on this pair is $A(r)\HN(x)+B(r)\,x\cdot\nabla\HN(x)$, and $A(q)=1$,
$A'(q)=0$, $B(q)=0$, $B'(q)=q^{-1}$.

\begin{lemma}[Energy shift]\label{lem:kr-energy}
$\det S_0>0$, and
\begin{equation}\label{eq:kr-collar}
 1-A(r)\ge(1-r/q)^2\quad(0\le r\le q),\qquad qA''(q)<0 .
\end{equation}
At a point $x$ with $q_x=q$, the conditional variance of $\HN(x)$ given $\nabla\HN(x)$ is
$N\sigma_H^2$, where
\begin{equation}\label{eq:kr-sigmaH}
 \sigma_H^2=\xi(q)-\frac{qc^2}{c+aq}=\frac{\det S_0}{s_2}>0 .
\end{equation}
\end{lemma}

\begin{proof}
Let $w_p=\beta_p^2q^p/s_0$, a probability vector on the degrees. Then
$\det S_0=s_0^2\Var_w(p)>0$, since two degrees are active. Grouping the numerator of $A$ into
unordered pairs of degrees gives, for $z\in[0,1]$,
\[
 A(qz)=\sum_{p<s}\lambda_{ps}\,\frac{sz^p-pz^s}{s-p},\qquad
 \lambda_{ps}=\frac{w_pw_s(s-p)^2}{\Var_w(p)},\qquad\sum_{p<s}\lambda_{ps}=1 .
\]
For $2\le p<s$, polynomial multiplication verifies
\[
 1-\frac{sz^p-pz^s}{s-p}=(1-z)^2Q_{ps}(z),\qquad
 Q_{ps}(z)=\sum_{j=0}^{p-1}(j+1)z^j+\sum_{j=p}^{s-2}\frac{p(s-1-j)}{s-p}\,z^j .
\]
All coefficients of $Q_{ps}$ are nonnegative and its constant coefficient is one, so
$Q_{ps}\ge1$ on $[0,1]$; summing with the weights $\lambda_{ps}$ proves the first inequality in
\eqref{eq:kr-collar}. Each pair has second derivative $-ps$ at $z=1$, so
$qA''(q)=-q^{-1}\sum_{p<s}\lambda_{ps}ps<0$. The series may be differentiated termwise because
$\xi$ has radius of convergence greater than one. Finally, $\HN(x)$ is independent of the part of
$\nabla\HN(x)$ orthogonal to $x$, and its regression on $x\cdot\nabla\HN(x)$ leaves the variance
$N(s_0-s_1^2/s_2)=N\det S_0/s_2$, which is \eqref{eq:kr-sigmaH}.
\end{proof}

For $\eta\ge0$ we select the energy density
\begin{equation}\label{eq:kr-heta}
 h_\eta=h_0+L\eta,
\end{equation}
with a constant $L\ge0$ fixed once and for all in \eqref{eq:kr-L} below. Both the covariance $\xi$
and the order parameter $\zeta$ defining $F_\zeta$ stay fixed throughout the construction; only
the data $(\pi_\eta,h_\eta)$ selected inside the Kac--Rice integral change.

\subsection{A strict conditional Hessian bound}\label{app:kr-hessian}

This subsection proves Proposition~\ref{prop:kr-local}: given criticality at a point whose
selected data are close to $(\pi,h)$, $F_\zeta$ is strongly concave on a ball of radius of order
$\sqrt N$ and its normalized log-determinant has a deterministic limit, if $\Plef_\pi>0$ and the
radial Schur complement has a negative limit $\mathcal S(\pi,h)$. Both fail at a marginal contact
\eqref{eq:kr-schur}; the tilt restores the first, the energy shift the second
\eqref{eq:kr-schur-shift}.

\paragraph{Inverse fields, the gradient and the Hessian.}
For $t\in[0,1]$ and $m\in(-1,1)$ let $b_t(m)=b_\zeta(t,m)$ be the inverse of $\Phi_x(t,\cdot)$, and
let $u(t,m)=u_\zeta(t,m)$, $v(t,m)=v_\zeta(t,m)$ and $T(t,m)=T_\zeta(t,m)$ be the susceptibility and
the two ratios of Appendix~\ref{app:matrix-hessian}:
\[
 u(t,m)=\Phi_{xx}(t,b_t(m)),\qquad
 v(t,m)=\frac{\Phi_{xxx}}{\Phi_{xx}}(t,b_t(m)),\qquad
 T(t,m)=\Bigl(\Phi_{xxxx}-\frac{\Phi_{xxx}^2}{\Phi_{xx}}\Bigr)(t,b_t(m)).
\]
At the radius $q$ these are functions of the field: $b_q(M(b))=b$, $u(q,M(b))=U(b)$, and we write
$v(b)=U'(b)/U(b)$ and $T(b)=U''(b)-U'(b)^2/U(b)$ for $v(q,M(b))$ and $T(q,M(b))$. We write
$\cdf\zeta^+(t)=\zeta([0,t])$ and $\cdf\zeta^-(t)=\zeta([0,t))$, as in
Appendix~\ref{app:matrix-hessian}; they differ only at the atoms of $\zeta$, and
$\cdf\zeta^\pm\ge\theta$ everywhere because $\zeta(\{0\})=\theta$. By \eqref{eq:matrix-inverse-field},
for $m\in(-1,1)$ and $t\in(0,1)$,
\begin{equation}\label{eq:kr-inverse-field}
 \partial_t^\pm b_t(m)=\frac{\xi''(t)}2\bigl(v(t,m)+2\cdf\zeta^\pm(t)\,m\bigr),\qquad
 \partial_t^\pm u(t,m)=-\frac{\xi''(t)}2T(t,m)-\xi''(t)\cdf\zeta^\pm(t)\,u(t,m)^2 ,
\end{equation}
together with $\partial_mb_t=1/u$ and $\partial_mu=v$. The Hessian bounds below are evaluated at
points whose coordinates may approach $\pm1$, so we also need $v$ and $T$ up to the boundary.

\begin{lemma}[Inverse fields up to the boundary]\label{lem:kr-endpoints}
The functions $u,v,T$ are jointly continuous on $[0,1]\times(-1,1)$, with $|v|\le C_{34}$ and
$|T|\le(C_{34}+C_{34}^2)\,u$. As $m\to\pm1$, uniformly in $t\in[0,1]$,
\[
 u(t,m)\to0,\qquad v(t,m)\to\mp2,\qquad T(t,m)\to0,
\]
so $u,v,T$ extend to jointly continuous functions on $[0,1]\times[-1,1]$. Moreover $b_t(m)$,
$u(t,m)$ and $\log u(t,m)$ are Lipschitz in $t\in[0,1]$, uniformly in $m\in(-1,1)$.
\end{lemma}

\begin{proof}
The bounds on $v$ and $T$ follow from Lemma~\ref{lem:kr-contact}\ref{lem:kr-contact-d}, since
$|T|\le|\Phi_{xxxx}|+\Phi_{xxx}^2/\Phi_{xx}$. Joint continuity on $[0,1]\times(-1,1)$ follows from
the continuity of the spatial derivatives of $\Phi$ and of $b_t(m)$
(Lemma~\ref{lem:pde-standing}\ref{lem:pde-standing-a}--\ref{lem:pde-standing-b}). The limit of $u$
is Lemma~\ref{lem:pde-standing}\ref{lem:pde-standing-b}, and then $T\to0$ by the bound on $T$.

For $v$, by oddness of $\Phi_x(t,\cdot)$ it suffices to show $v(t,m)\to-2$ as $m\uparrow1$.
Since $b_t(m)\ge\frac12\log\frac{c_1}{2(1-m)}$ for $m\in[0,1)$ (proof of
Lemma~\ref{lem:pde-standing}\ref{lem:pde-standing-b}), it suffices to show
$\Phi_{xxx}(t,b)/\Phi_{xx}(t,b)\to-2$ as $b\to\infty$, uniformly in $t$. First let $\zeta'$ be
finitely supported, let $\tilde X$ be its Auffinger--Chen process started at $(t,b)$ and put
$I=\int_t^1\xi''(r)\cdf{\zeta'}(r)\Phi_{\zeta',xx}(r,\tilde X_r)\,dr\in[0,V_*]$, $V_*=\xi'(1)$. The
Feynman--Kac formulas in the proof of Lemma~\ref{lem:pde-standing} give
\[
 \frac{\Phi_{\zeta',xxx}(t,b)}{\Phi_{\zeta',xx}(t,b)}+2
 =\frac{2\,\E\bigl[\sech^2(\tilde X_1)\,e^{I}\,\bigl(1-\tanh(\tilde X_1)e^{2I}\bigr)\bigr]}
  {\E\bigl[\sech^2(\tilde X_1)\,e^{I}\bigr]} .
\]
Write $\tilde X_r=b+D_r+J_r$ with $|D_r|\le V_*$ and $J_r=\int_t^r\sqrt{\xi''(s)}\,dB_s$, and let
$R>0$ and $b>V_*+R$. On $\{\sup_r|J_r|\le R\}$ the path stays above $b-V_*-R$, so
$I\le V_*C_1e^{-2(b-V_*-R)}$ and $1-\tanh\tilde X_1\le2e^{-2(b-V_*-R)}$; there the last factor in
the numerator is at most $\epsilon(b,R)$, with $\epsilon(b,R)\to0$ as $b\to\infty$ for fixed $R$.
Elsewhere it is at most $1+e^{2V_*}$. The denominator is at least $\E\sech^2(\tilde X_1)\ge
e^{-2b-2V_*}\E e^{-2|J_1|}$, and the numerator restricted to $\{\sup_r|J_r|>R\}$ is at most
$Ce^{-2b}\E\bigl[e^{2|J_1|}\ind\{\sup_r|J_r|>R\}\bigr]$, with $C$ depending only on $V_*$. By
Cauchy--Schwarz and the Gaussian maximal inequality for $J$, whose variance is at most $V_*$, the
ratio of these two bounds tends to zero as $R\to\infty$, uniformly in $b$, $t$ and $\zeta'$. Taking
$b\to\infty$ and then $R\to\infty$ gives the claim for finitely supported $\zeta'$, uniformly in
$\zeta'$. For $\zeta$ it follows by the uniform convergence of the spatial derivatives along
finitely supported approximations (proof of Lemma~\ref{lem:pde-standing}\ref{lem:pde-standing-a}).

Finally, the one-sided derivatives in \eqref{eq:kr-inverse-field} are bounded, and
$|\partial_t^\pm\log u|\le\xi''(1)\bigl(\frac12(C_{34}+C_{34}^2)+1\bigr)$; a continuous function
with bounded one-sided derivatives is Lipschitz.
\end{proof}

For $x\in\ocube$ with $r=q_x>0$, write $b^{(r)}=(b_r(x_i))_i$, $u_i=u(r,x_i)$ and
$\chi_x=N^{-1}\sum_iu_i$. By Lemma~\ref{lem:matrix-hessian}\ref{lem:matrix-hessian-a},
\begin{equation}\label{eq:kr-gradient}
 \nabla F_\zeta(x)=\nabla\HN(x)-b^{(r)}-\xi''(r)\chi_x\,x ;
\end{equation}
the terms produced by $\cdf\zeta^\pm$ cancel, so $F_\zeta$ is continuously differentiable on
$\ocube\setminus\{0\}$, and by \eqref{eq:kr-inverse-field} and Lemma~\ref{lem:kr-endpoints} its
gradient is locally Lipschitz there, also across the spheres $\{q_x=t\}$ at atoms $t$ of $\zeta$.
By Lemma~\ref{lem:matrix-hessian}\ref{lem:matrix-hessian-b}, on either side $\sigma\in\{+,-\}$ of
such a sphere,
\begin{equation}\label{eq:kr-hessian}
\begin{split}
 \nabla^2_\sigma F_\zeta(y)={}&\nabla^2\HN(y)-\diag(u_i^{-1})-\xi''(r)\chi_yI
 -\frac{\xi''(r)}N\bigl(y\mathbf v_y^{\mathsf T}+\mathbf v_yy^{\mathsf T}\bigr)\\
 &-\frac2N\Bigl\{\xi'''(r)\chi_y-\frac{\xi''(r)^2}{2N}\sum_iT(r,y_i)
 +\xi''(r)\cdf\zeta^\sigma(r)\Plef_N(y)\Bigr\}yy^{\mathsf T},
\end{split}
\end{equation}
where $r=q_y$, $u_i=u(r,y_i)$, $\mathbf v_y=(v(r,y_i))_i$ and
$\Plef_N(y)=1-\xi''(r)N^{-1}\sum_iu_i^2$ is the empirical Plefka factor. The reciprocal $u_i^{-1}$
appears only on the diagonal. The two one-sided Hessians differ by
$-\frac2N\xi''(r)\zeta(\{r\})\Plef_N(y)\,yy^{\mathsf T}$, and the right derivative of $\nabla F_\zeta$ at $y$ in a
direction $e$ is $\nabla^2_+F_\zeta(y)e$ if $y\cdot e\ge0$ and $\nabla^2_-F_\zeta(y)e$ if $y\cdot e\le0$
(Lemma~\ref{lem:matrix-hessian}\ref{lem:matrix-hessian-b}).

\paragraph{Exact conditional blocks at an anchor.}
Consider an \emph{anchor} $x\in\ocube$ with $q_x=q$, so that $x_i=M(b_i)$ with $b=(b_i)$ the
fields of $x$. Write $\langle\varphi\rangle_N=N^{-1}\sum_i\varphi(b_i)$ and
\[
 \chi_N=\langle U\rangle_N,\qquad k_N=\langle bM\rangle_N,\qquad d_N=k_N+aq\chi_N,\qquad
 e=\frac{x}{\sqrt{Nq}},\qquad P=I-ee^{\mathsf T}.
\]
By \eqref{eq:kr-gradient}, the critical equation $\nabla F_\zeta(x)=0$ reads
$\nabla\HN(x)=b+a\chi_Nx$. It holds for every empirical field law of self-overlap $q$, and it
uses no integration-by-parts identity for that law. We condition on
\begin{equation}\label{eq:kr-conditioning}
 \HN(x)=Nh,\qquad \nabla\HN(x)=b+a\chi_Nx .
\end{equation}
With $S_0$ and $A,B$ from Lemma~\ref{lem:kr-energy}, define
$\mathbf t=(qa,\ q(2a+qa_1))^{\mathsf T}$ and
\begin{equation}\label{eq:kr-block-coefficients}
 A_2=qA''(q),\qquad B_2=qB''(q),\qquad
 \tau=a+qa_1-\frac{qa^2}c,\qquad
 \varsigma^2=2a+4qa_1+q^2a_2-\mathbf t^{\mathsf T}S_0^{-1}\mathbf t .
\end{equation}
Then $(A_2,B_2)=\mathbf t^{\mathsf T}S_0^{-1}$, both $\tau/N$ and $\varsigma^2/N$ are conditional
Gaussian variances, and, with $s_3=\sum_pp^3\beta_p^2q^p$,
\[
 \tau=\frac{s_1s_3-s_2^2}{q^2s_1}>0,\qquad A_2=-\frac{qs_1}{\det S_0}\,\tau<0 ;
\]
the strict signs use two active degrees and $q>0$. Finally, with $v$ and $T$ as above, put
\begin{gather*}
 \dot u(b)=-\frac a2T(b)-a\theta U(b)^2,\qquad
 \ell(b)=\frac bc-v(b),\\
 \psi(b)=-\frac{M(b)}{\sqrt q\,U(b)}+a\sqrt q\,\ell(b),\qquad
 \psi_N=\frac1{\sqrt N}\bigl(\psi(b_i)\bigr)_i .
\end{gather*}

\begin{lemma}[Conditional block representation]\label{lem:kr-block}
Under \eqref{eq:kr-conditioning}, relative to the decomposition $\R e\oplus e^\perp$, the left
Hessian has the law
\[
 \nabla^2_-F_\zeta(x)\overset{d}{=}
 \begin{pmatrix} r_N+\varsigma Z/\sqrt N & z_N^{\mathsf T}\\ z_N & C_N\end{pmatrix},
\]
where
\begin{align}
 C_N&=P\bigl[\sqrt a\,\GOE-\diag(U(b_i)^{-1}+a\chi_N)\bigr]P\big|_{e^\perp},\label{eq:kr-transverse}\\
 z_N&=P\psi_N+\sqrt{\tau/N}\,g,\label{eq:kr-cross}\\
 r_N&=A_2h+B_2d_N-\frac{\langle M^2/U\rangle_N}q-a\chi_N-2a\langle Mv\rangle_N
 -2q\bigl[a\theta+a_1\chi_N+a\langle\dot u\rangle_N\bigr].\label{eq:kr-radial}
\end{align}
Here $\GOE$ is a GOE matrix, $g$ a standard Gaussian vector in $e^\perp$ and $Z$ a standard
Gaussian variable, all independent. The right Hessian has the same representation with $\theta$
replaced by $\zeta([0,q])$ in $r_N$ and in $\dot u$.
\end{lemma}

\begin{proof}
Write $\mathcal G=\nabla\HN(x)$ and $\mathcal H=\nabla^2\HN(x)$. Differentiating the covariance
$N\xi(\ov yz)$ gives $\Cov(\HN(x),x\cdot\mathcal G)=NS_0$ and $\Cov(P\mathcal G)=cP$, the two
being independent, and, for unit vectors $u',v'\perp e$,
\begin{gather*}
 \Cov(u'^{\mathsf T}\mathcal He,v'\cdot\mathcal G)=\frac{a\sqrt q}{\sqrt N}\,u'\cdot v',\qquad
 \Cov(u'^{\mathsf T}\mathcal He,v'^{\mathsf T}\mathcal He)=\frac{a+qa_1}N\,u'\cdot v',\\
 \Cov\bigl(e^{\mathsf T}\mathcal He,(\HN(x),x\cdot\mathcal G)\bigr)=\mathbf t^{\mathsf T},
\end{gather*}
and $\Var(e^{\mathsf T}\mathcal He)=(2a+4qa_1+q^2a_2)/N$. The transverse block
$P\mathcal HP|_{e^\perp}$ has GOE covariance with off-diagonal variance $a/N$ and is uncorrelated
with all conditioned variables, with the cross entries and with the radial entry; the cross
entries are uncorrelated with $\HN(x)$, $x\cdot\mathcal G$ and the radial entry, and the radial
entry with $P\mathcal G$. Gaussian regression, with $P\mathcal G=Pb$ under
\eqref{eq:kr-conditioning}, gives
\[
 P\mathcal He=\frac{a\sqrt q}{c\sqrt N}\,Pb+\sqrt{\tau/N}\,g,\qquad
 e^{\mathsf T}\mathcal He=A_2h+B_2d_N+\varsigma Z/\sqrt N .
\]
Now take $y=x$ and $\sigma=-$ in \eqref{eq:kr-hessian}, where $\cdf\zeta^-(q)=\theta$ and
$\xi'''(q)\chi_N-\frac{a^2}{2N}\sum_iT(b_i)+a\theta\Plef_N(x)=a\theta+a_1\chi_N+a\langle\dot u\rangle_N$.
Its transverse, cross and radial blocks are \eqref{eq:kr-transverse}--\eqref{eq:kr-radial}; for
the cross block note $P\diag(U^{-1})e=P(M/U)(b_i)_i/\sqrt{Nq}$ and
$\frac aNP\mathbf v_x(x\cdot e)=a\sqrt{q/N}\,P\mathbf v_x$. Complete the transverse GOE block to a
full GOE matrix by independent entries; they do not appear in $P\GOE P$. The right Hessian
differs only through $\cdf\zeta^+(q)=\zeta([0,q])$.
\end{proof}

\paragraph{The radial Schur complement.}
Selected fields must have controlled exponential moments. Put $t_0=8$. For a law $\pi$ on $\R$
with finite exponential moments of every order fix $L_0>\E_\pi e^{t_0|b|}$; the \emph{cap} on the
fields of $x$ is
\begin{equation}\label{eq:kr-cap}
 \frac1N\sum_ie^{t_0|b_i|}\le L_0 .
\end{equation}
By Lemma~\ref{lem:kr-contact}\ref{lem:kr-contact-d}, $M^2/U\le Ce^{2|b|}$, $|v|\le C$, $|\dot u|\le C$
and $\psi(b)^2\le C(e^{4|b|}+b^2+1)$. So under the cap $|\psi_N|\le B_\psi$, with $B_\psi$ depending
only on $L_0$, and every continuous function of growth at most $Ce^{4|b|}$ is uniformly integrable
along empirical field laws that satisfy the cap. For a law $\pi$ with $\E_\pi M^2=q$ and $h\in\R$
define
\begin{equation}\label{eq:kr-schur-generic}
\begin{split}
 \mathcal S(\pi,h)={}&A_2h+B_2d_\pi+(\tau-a)\chi_\pi-\frac{2ak_\pi}c
 -2q\bigl\{a\theta+a_1\chi_\pi+a\E_\pi\dot u\bigr\}\\
 &+qa^2\E_\pi\bigl[U\ell^2\bigr]-\frac q{\upsilon_\pi}\bigl[-1+a\E_\pi(UM\ell)\bigr]^2 .
\end{split}
\end{equation}
It is a continuous function of finitely many field moments, and $\partial_h\mathcal S=A_2<0$.

The formulas above were written at an anchor with $q_x=q$. For an anchor of actual self-overlap
$Q$ near $q$ we use the actual inverse fields $b_Q(x_i)$, take $M,U$ at time $Q$ (so
$U_i=u(Q,x_i)$), replace every covariance coefficient $q,a,c,a_1,a_2,s_j,A,B$ by its value at
$Q$, and in the rank-one term of \eqref{eq:kr-hessian} keep the coefficient $\theta$ in place of
$\cdf\zeta^\pm(Q)$; we call the result the \emph{frozen} left Hessian. The block algebra of
Lemma~\ref{lem:kr-block} holds verbatim for it, with $\dot u=-\frac{\xi''(Q)}2T-\xi''(Q)\theta U^2$
used as an algebraic abbreviation.

\begin{lemma}[Radial Schur complement]\label{lem:kr-schur}
Let $\pi$ have finite exponential moments of every order, with $\E_\pi M^2=q$ and
$\Plef_\pi=1-a\E_\pi U^2>0$, and let $h\to h_\infty$ as $N\to\infty$.
\begin{thmlist}
\item\label{lem:kr-schur-a} Consider anchors whose fields satisfy the cap \eqref{eq:kr-cap}, with
 $W_2(\pi_x,\pi)\to0$ as $N\to\infty$. Under \eqref{eq:kr-conditioning}, with probability tending
 to one, $C_N\preceq-g_0I$ for a constant $g_0>0$ depending only on $a$ and $\Plef_\pi$, and
 \[
  r_N+\varsigma Z/\sqrt N-z_N^{\mathsf T}C_N^{-1}z_N\longrightarrow\mathcal S(\pi,h_\infty)
  \qquad\text{in probability.}
 \]
\item\label{lem:kr-schur-b} If $\mathcal S(\pi,h_\infty)<0$, there is $g>0$ such that, for the
 anchors of~\ref{lem:kr-schur-a}, both one-sided Hessians satisfy $\nabla^2_\pm F_\zeta(x)\preceq-gI$
 with probability tending to one.
\item\label{lem:kr-schur-c} \textup{(Anchors off the sphere.)} Consider anchors $x$ with
 $q_x\to q$ whose fields satisfy \eqref{eq:kr-cap}, with $W_2(\pi_x,\pi)\to0$ as $N\to\infty$, and
 condition on $\HN(x)=Nh$ and $\nabla F_\zeta(x)=0$. Then the conclusions of~\ref{lem:kr-schur-a}
 hold for the blocks of the frozen left Hessian, and, if $\mathcal S(\pi,h_\infty)<0$, the conclusion
 of~\ref{lem:kr-schur-b} holds for the frozen left Hessian.
\end{thmlist}
\end{lemma}

\begin{proof}
\ref{lem:kr-schur-a} Let $\widetilde C_N=\sqrt a\,\GOE-\diag(U(b_i)^{-1}+a\chi_N)$ on $\R^N$. The empirical laws of the
$U(b_i)$ converge to $\Law_\pi(U)$, a law on $(0,1]$ with mean $\chi_\pi$ and
$1-a\E_\pi U^2=\Plef_\pi>0$. Lemma~\ref{lem:matrix-gap}\ref{lem:matrix-gap-a}--\ref{lem:matrix-gap-b} gives a gap
$\widetilde C_N\preceq-g_0I$ with probability at least $1-e^{-c_MN}$, the inverse equivalent
$v_N^{\mathsf T}\widetilde C_N^{-1}w_N+\sum_iU(b_i)v_{N,i}w_{N,i}\to0$ in probability for vectors
of bounded norm that are deterministic or independent of $\GOE$, and
$N^{-1}\tr\widetilde C_N^{-1}\to-\chi_\pi$ (all set to zero off the gap event). The inverse of the
compression $C_N$, extended by zero on $\R e$, is
\[
 K_N=\widetilde C_N^{-1}-\frac{\widetilde C_N^{-1}ee^{\mathsf T}\widetilde C_N^{-1}}
 {e^{\mathsf T}\widetilde C_N^{-1}e},
\]
since $K_Ne=0$ and $P\widetilde C_NPK_N=P$. The denominator converges to $-\upsilon_\pi/q<0$. By
the inverse equivalent applied to $\psi_N$ and $e$,
$\psi_N^{\mathsf T}K_N\psi_N\to-\E_\pi(U\psi^2)+[\E_\pi(UM\psi)]^2/\upsilon_\pi$. The rank-one
term of $K_N$ has trace at most $g_0^{-1}$ in absolute value on the gap event, so
$N^{-1}\tr K_N\to-\chi_\pi$. Conditionally on $\GOE$, the Gaussian cross term and the centred
quadratic form in $g$ have variance $O(N^{-1})$. Hence
\[
 z_N^{\mathsf T}C_N^{-1}z_N\longrightarrow-\E_\pi(U\psi^2)+\frac{[\E_\pi(UM\psi)]^2}{\upsilon_\pi}
 -\tau\chi_\pi ,
\]
while $\varsigma Z/\sqrt N\to0$ and $r_N$ converges to the corresponding $\pi$-average, by uniform
integrability under the cap. The apparently singular moments of $1/U$ cancel: since
$U\psi^2=M^2/(qU)-2aM\ell+qa^2U\ell^2$, $UM\psi=\sqrt q(-M^2/q+aUM\ell)$ and
$\E_\pi(M\ell)=k_\pi/c-\E_\pi(Mv)$, the limit of the Schur complement is \eqref{eq:kr-schur-generic}.

\ref{lem:kr-schur-b} Take an event of probability tending to one on which
$C_N\preceq-g_0I$, the Schur complement is at most $-s<0$ and $|g|^2\le2N$. There
$\|C_N^{-1}z_N\|\le K_*:=g_0^{-1}(B_\psi+\sqrt{2\tau})$. For $(t,v)\in\R\oplus e^\perp$, with
$\hat r_N=r_N+\varsigma Z/\sqrt N$ and $w=v+tC_N^{-1}z_N$,
\[
 \hat r_Nt^2+2tz_N^{\mathsf T}v+v^{\mathsf T}C_Nv
 =(\hat r_N-z_N^{\mathsf T}C_N^{-1}z_N)t^2+w^{\mathsf T}C_Nw
 \le-\min(s,g_0)(t^2+|w|^2),
\]
and $t^2+|v|^2\le\max(2,1+2K_*^2)(t^2+|w|^2)$. This is a gap for $\nabla^2_-F_\zeta(x)$. The right
Hessian differs by $-\frac2N a\,\zeta(\{q\})\Plef_N(x)\,xx^{\mathsf T}$, which is negative
semidefinite once $\Plef_N(x)\ge0$, and $\Plef_N(x)\to\Plef_\pi>0$.

\ref{lem:kr-schur-c} By Lemma~\ref{lem:kr-endpoints}, $|b_Q(m)-b_q(m)|\le C|Q-q|$ uniformly in
$m$, so the empirical laws of the fields $b_Q(x_i)$ and $b_q(x_i)=M^{-1}(x_i)$ are within
$W_2$-distance $C|Q-q|$, and \eqref{eq:kr-cap} persists with $L_0$ replaced by
$e^{t_0C|Q-q|}L_0$. The coefficients are smooth in $Q$, and $u,v,T$ are uniformly continuous on
$[q/2,(1+q)/2]\times[-1,1]$. Hence empirical windows formulated with the fixed-radius fields
$\pi_x$ give limiting moments, and values of \eqref{eq:kr-schur-generic}, within $o(1)$ as
$Q\to q$, and the proofs of~\ref{lem:kr-schur-a} and of the left-Hessian part
of~\ref{lem:kr-schur-b} apply to the frozen left Hessian.
\end{proof}

\paragraph{The Schur complement at the contact.}
At the data of the contact the Schur complement is explicit. Apply \eqref{eq:kr-ward} to $F=UM$,
$F=bU/c$ and $F=U'$ (all bounded with bounded derivatives, by
Lemma~\ref{lem:kr-contact}\ref{lem:kr-contact-d}):
\[
 \E_0(UM\ell)=\E_0U^2+\theta\upsilon_0,\qquad
 \E_0(U\ell^2)+\E_0T=\frac{\chi_0}c+\theta\E_0U^2+\theta^2\upsilon_0 .
\]
The identity $(A_2,B_2)S_0=\mathbf t^{\mathsf T}$ gives $A_2c+B_2(c+aq)=2a+qa_1$ and
$A_2\xi(q)+B_2cq=qa$. With $h_0=c\chi_0+\theta\xi(q)$ and $d_0=(c+aq)\chi_0+c\theta q$ (by
\eqref{eq:kr-k0}), the terms of \eqref{eq:kr-schur-generic} before the last square add up to
$-3qa\theta\Plef_0+qa^2\theta^2\upsilon_0$, and the last square is
$-q(-\Plef_0+a\theta\upsilon_0)^2/\upsilon_0$. Therefore
\begin{equation}\label{eq:kr-schur}
 \mathcal S(\pi_0,h_0)=-qa\theta\Plef_0-\frac{q\Plef_0^2}{\upsilon_0} .
\end{equation}
At $\Plef_0=0$ this is the value of the continuous expression \eqref{eq:kr-schur-generic}; it says
nothing about the inversion of a marginal finite-dimensional matrix. If $\Plef_0>0$, then
$\mathcal S(\pi_0,h_0)<0$, and the untilted data satisfy both hypotheses of
Lemma~\ref{lem:kr-schur}.

Now let $\Plef_0=0$ and $0<\eta\le\eta_0$. Every field moment appearing in
\eqref{eq:kr-schur-generic} is the $\pi_\eta$-average of a function of growth at most
$C(1+b^2)$, hence affine in $\eta$ with a bounded slope, and $\upsilon_\eta\ge\upsilon_0/2$ because
the density of $\pi_\eta$ is at least $\frac12$. So there is $C_S$ with
$|\mathcal S(\pi_\eta,h_0)-\mathcal S(\pi_0,h_0)|\le C_S\eta$ for $0\le\eta\le\eta_0$. Since
$\mathcal S(\pi_0,h_0)=0$ and $\partial_h\mathcal S=A_2<0$, every
$L\ge L_S:=(C_S+1)/|A_2|$ gives
\begin{equation}\label{eq:kr-schur-shift}
 \mathcal S(\pi_\eta,h_\eta)\le-\eta\qquad(0<\eta\le\eta_0),\qquad\text{while}\quad
 \Plef_\eta=a\eta\,\E_0\mathcal R^2>0 .
\end{equation}

\paragraph{Third derivatives under the conditioning.}
Two active degrees enter once more, through the conditional mean of the third derivative of
$\HN$.

\begin{lemma}[Conditional third derivatives]\label{lem:kr-c3}
For every $H_0<\infty$ there are $\varepsilon_0,c_3>0$ and $C=C(H_0)$ such that the following holds
for every $N$, every $x\in\ocube$ with
$|q_x-q|\le\varepsilon_0$ and every value $(Nh',\mathbf y)$ of $(\HN(x),\nabla\HN(x))$ with
$|h'|+N^{-1/2}|\mathbf y|\le H_0$: conditionally on $(\HN(x),\nabla\HN(x))=(Nh',\mathbf y)$,
\[
 \Prob_{\rm cond}\Bigl(\sup_{z\in\cube}\|D^3\HN(z)\|_{\rm op}>\frac C{\sqrt N}\Bigr)\le e^{-c_3N}.
\]
In particular, on the complement, $\|\nabla^2\HN(y)-\nabla^2\HN(x)\|\le
C|y-x|/\sqrt N$ for all $y\in\cube$.
\end{lemma}

\begin{proof}
Write $\HN=m_{\rm cond}+G$, with $m_{\rm cond}$ the conditional mean and $G$ the centred conditional
field. Gaussian conditioning decreases pointwise variances and variances of increments, so the
bounds on the pointwise variance and the canonical metric of the trilinear process
$(z,\varphi_1,\varphi_2,\varphi_3)\mapsto D^3\HN(z)[\varphi_1,\varphi_2,\varphi_3]$ used in the proof of
Lemma~\ref{lem:matrix-third}\ref{lem:matrix-third-b} hold for $G$, and that proof gives the bound for
$G$ (Remark~\ref{rem:matrix-conditional}). For the mean,
let $\mathbf Z=(\HN(x)/\sqrt N,\nabla\HN(x))$. Its covariance is uniformly nondegenerate for
$|q_x-q|\le\varepsilon_0$: the part of $\nabla\HN(x)$ orthogonal to $x$ has covariance
$\xi'(q_x)$ times the identity, and the pair $(\HN(x)/\sqrt N,\,x\cdot\nabla\HN(x)/|x|)$ has a
covariance matrix whose determinant is $\det S_0(q_x)/q_x>0$, where $S_0(r)$ is $S_0$ with $q$
replaced by $r$ (Lemma~\ref{lem:kr-energy}). On the conditioning event
$|\Cov(\mathbf Z)^{-1/2}\mathbf Z|\le C\sqrt N$. By Gaussian regression and Cauchy--Schwarz, each
$m_{\rm cond}$-derivative $D^3m_{\rm cond}(z)[\varphi_1,\varphi_2,\varphi_3]$, $|\varphi_j|=1$, is at most the
unconditional standard deviation of $D^3\HN(z)[\varphi_1,\varphi_2,\varphi_3]$, which is at most $C/N$
by the variance bound in the proof of Lemma~\ref{lem:matrix-third}\ref{lem:matrix-third-b}, times
$C\sqrt N$. This is the
only place where the nondegeneracy of $\Cov(\mathbf Z)$, and hence the second active degree, is
used in this lemma.
\end{proof}

\paragraph{The strict Hessian bound.}
We now state the result of this section. For a law $\pi$ and $h\in\R$, data $(x,h')$, with
$x\in\ocube$ and $h'\in\R$, are \emph{admissible of width $w$} if the fields of $x$ satisfy
\eqref{eq:kr-cap} and $W_2(\pi_x,\pi)+|h'-h|\le w$. Since $M$ is $1$-Lipschitz and
$|m^2-m'^2|\le2|m-m'|$ on $[-1,1]$, admissible data have $|q_x-\E_\pi M^2|\le2w$. We write
$\Prob_{x,h'}$ for the conditional law given $\HN(x)=Nh'$ and $\nabla F_\zeta(x)=0$, and
\[
 \mathfrak d_\pi=-\E_\pi\log U+\frac a2\chi_\pi^2 .
\]

\begin{proposition}[Strict Hessian bound around a conditioned critical point]\label{prop:kr-local}
Let $\pi$ be a law on $\R$ with finite exponential moments of every order and $h\in\R$, such that
\[
 \E_\pi M^2=q,\qquad \Plef_\pi=1-a\E_\pi U^2>0,\qquad \mathcal S(\pi,h)<0,
\]
and fix $L_0>\E_\pi e^{t_0|b|}$. There are $g_*,\rho_*>0$ such that the following hold, where the
suprema are over admissible data of width at most $w$.
\begin{thmlist}
\item\label{prop:kr-local-a} \textup{(Strong concavity on a ball.)}
 \[
  \lim_{w\downarrow0}\limsup_{N\to\infty}\sup_{(x,h')}\Prob_{x,h'}\Bigl(\exists\,y\in\ocube,\
  |y-x|\le\rho_*\sqrt N,\ \sigma\in\{\pm\}:\ \nabla^2_\sigma F_\zeta(y)\not\preceq-g_*I\Bigr)=0 .
 \]
\item\label{prop:kr-local-b} \textup{(Typical determinant.)} For every $t>0$,
 \[
  \lim_{w\downarrow0}\limsup_{N\to\infty}\sup_{(x,h')}\Prob_{x,h'}\Bigl(\exists\,\sigma\in\{\pm\}:\
  \Bigl|\frac1N\log\bigl|\det\nabla^2_\sigma F_\zeta(x)\bigr|-\mathfrak d_\pi\Bigr|>t\Bigr)=0 .
 \]
\end{thmlist}
The constants $g_*,\rho_*$ depend on $\pi$, $h$ and $L_0$ only.
\end{proposition}

For the data $(\pi_0,h_0)$ at a strictly Plefka-stable first contact, part~\ref{prop:kr-local-a}
is, for the left Hessian, Boursier's conditional Hessian estimate
\cite[Proposition~B.1(ii)]{Bou26b}. The proof below works
for general data $(\pi,h)$, with constants depending on them, and evaluates no matrix inverse at
the marginal data $(\pi_0,h_0)$. Part~\ref{prop:kr-local-b} is a law of large numbers in
conditional probability; the count of Appendix~\ref{app:kr-count} needs neither an annealed
determinant formula nor a lower-tail estimate for the determinant (compare
\cite[Lemma~B.2]{Bou26b}).

\begin{proof}
By the contradiction argument at the end of the proof it suffices to treat sequences of data
whose widths tend to zero; along such a sequence $q_x\to q$, the empirical field laws converge to
$\pi$ in $W_2$, \eqref{eq:kr-cap} holds and $h'\to h$. We work under $\Prob_{x,h'}$, and all
statements hold with conditional probability tending to one.

\emph{Step 1: a continuous upper envelope.} For $K\ge1$ and $y\in\ocube$ with $r=q_y>0$, let
$\mathcal M_{N,K}(y)$ be the matrix obtained from \eqref{eq:kr-hessian} by replacing $u_i^{-1}$ with
$\min(u_i^{-1},K)$ and $\cdf\zeta^\sigma(r)$ with $\theta$. By Lemma~\ref{lem:kr-endpoints}, all
coefficient functions of $\mathcal M_{N,K}$ are uniformly continuous functions of $(r,y_i)$ on
$[q/2,(1+q)/2]\times[-1,1]$, across every atom of $\zeta$. If $|y-x|\le\rho\sqrt N$, at most a
fraction $\rho$ of the coordinates have $|y_i-x_i|>\sqrt\rho$, so averages $N^{-1}\sum_i\phi(r,y_i)$
of such functions change by at most $\omega_\phi(\sqrt\rho+|q_y-q_x|)+2\|\phi\|_\infty\rho$. In
particular $\Plef_N(y)\ge\Plef_\pi/2$ for $|y-x|\le\rho_1\sqrt N$ and narrow windows. Since
$\cdf\zeta^\sigma\ge\theta$, this gives the exact matrix inequality
\begin{equation}\label{eq:kr-envelope}
 \nabla^2_\sigma F_\zeta(y)\preceq\mathcal M_{N,K}(y)\qquad(\sigma=\pm,\ |y-x|\le\rho_1\sqrt N).
\end{equation}

\emph{Step 2: the cap keeps the anchor gap.} The matrix $\mathcal M_{N,\infty}(x)$ is the frozen left
Hessian at $x$, so Lemma~\ref{lem:kr-schur}\ref{lem:kr-schur-c} gives $-\mathcal M_{N,\infty}(x)\succeq g_1I$. Write $-\mathcal M_{N,\infty}(x)=A_N-F_N$ with
$A_N=\diag(u_i^{-1}+a\chi_N)-\sqrt a\,\GOE$ (the transverse GOE completed as in
Lemma~\ref{lem:kr-block}) and $F_N$ of bounded rank. By Lemma~\ref{lem:kr-block}, the vectors
spanning the range of $F_N$ are, after normalization by $\sqrt N$, combinations of $x$, $b$,
$(M/U)(b_i)$ and bounded functions of the fields, so \eqref{eq:kr-cap} gives $\|F_N\|\le C$. By
Lemma~\ref{lem:matrix-gap}\ref{lem:matrix-gap-a}, $A_N\succeq g_0I$. A convex combination of
$A_N-F_N\succeq g_1I$ and $A_N-F_N\succeq A_N-CI$ gives
\[
 A_N-F_N\succeq t_*A_N,\qquad t_*=\frac{g_1}{C+g_1}.
\]
Let $A_{N,K}$ be $A_N$ with $u_i^{-1}$ capped at $K$. The capped statement of
Lemma~\ref{lem:matrix-gap}\ref{lem:matrix-gap-a} gives $A_{N,K}\succeq g_0'I$ for all large fixed $K$. On $\{\|\sqrt a\,\GOE\|\le B\}$ we have
\[
 \|A_N^{-1}-A_{N,K}^{-1}\|\le C'/K :
\]
split the indices into those with $u_i^{-1}\le K$ and the rest; on the second block both inverses
have norm at most $(K-B)^{-1}$, the Schur complements of the first block differ by at most
$2B^2/(K-B)$ and have inverses of norm at most $1/g_0'$, and the block inverse formula bounds every
block of the difference by $O(K^{-1})$. Hence $\|A_N^{-1/2}-A_{N,K}^{-1/2}\|=O(K^{-1/2})$.
Conjugating the last display by $A_N^{-1/2}$ and using $\|F_N\|\le C$ gives
$I-A_{N,K}^{-1/2}F_NA_{N,K}^{-1/2}\succeq\frac{t_*}2I$ for large fixed $K$. Since capping changes
only the diagonal, $\mathcal M_{N,K}(x)=F_N-A_{N,K}\preceq-\frac{t_*}2g_0'I$. No independence of
$F_N$ and $\GOE$ is used.

\emph{Step 3: a uniform negative median near the anchor.} This step cannot follow from
operator-norm continuity of the diagonal, since a single coordinate of $y$ may move by a fixed
amount. Under $\Prob_{x,h'}$ with $w\le1$ we have $|h'|\le|h|+1$, and
$\nabla\HN(x)=b^{(q_x)}+\xi''(q_x)\chi_xx$ by~\eqref{eq:kr-gradient}, so $N^{-1/2}|\nabla\HN(x)|$ is
bounded by~\eqref{eq:kr-cap} and Lemma~\ref{lem:kr-endpoints}; this fixes $H_0$, depending only on
$h$, $L_0$ and $\xi$, and with it the constants $C=C(H_0)$ and $c_3$ of Lemma~\ref{lem:kr-c3},
used here and in Step~4. By that lemma, outside an event of conditional probability $e^{-c_3N}$ the
conditional Hessian of $\HN$ moves by at most $C|y-x|/\sqrt N\le C\rho$ in operator norm over the
ball. Freeze it at $x$. Apart from the capped diagonal, every deterministic term of
$\mathcal M_{N,K}(y)$ then changes by $o_\rho(1)$ in operator norm, uniformly for
$|y-x|\le\rho\sqrt N$, by Step~1; normalized rank-one terms are controlled by the normalized
Euclidean norms of their vectors. The remaining matrix, expressed through the conditional
transverse GOE and the independent cross vector of Lemma~\ref{lem:kr-block}, has transverse block
$P[\sqrt a\,\GOE-\diag(\min(u(r,y_i)^{-1},K)+a\chi_N)]P$, cross vector
\[
 P\Bigl\{-\diag\bigl(\min(u(r,y_i)^{-1},K)\bigr)e+\frac{a\sqrt q}{\sqrt N}\bigl(b/c-\mathbf v_x\bigr)\Bigr\}
 +\sqrt{\tau/N}\,g,
\]
and radial entry obtained from \eqref{eq:kr-radial} by replacing $q^{-1}\langle M^2/U\rangle_N$
with $e^{\mathsf T}\diag(\min(u(r,y_i)^{-1},K))e$, plus the independent residual
$\varsigma Z/\sqrt N$. Its Schur vectors are deterministic vectors plus the same Gaussian $g$,
independent of the auxiliary full GOE. The finite-$N$ comparison
\eqref{eq:matrix-comparison} still gives a transverse gap: at its fixed trial point the
reciprocal denominators stay bounded below and their average changes by $o_\rho(1)$. The
compressed-inverse formula uses finitely many bilinear forms of the transverse inverse with
deterministic profile vectors and the independent cross vector, and the proof of
Lemma~\ref{lem:matrix-gap}\ref{lem:matrix-gap-b}--\ref{lem:matrix-gap-c}, which works at a fixed
spectral parameter through Lemma~\ref{lem:matrix-comparison}\ref{lem:matrix-comparison-b}--\ref{lem:matrix-comparison-c},
applies to the changed bounded diagonal; the changed cross weight satisfies $|\min(u(r,y_i)^{-1},K)x_i|\le K|x_i|$, so no
integrability of the fields of $y$ is needed. For weights involving the fields of $x$ first restrict
to $|b_i|\le L'$ and use uniform continuity; \eqref{eq:kr-cap} makes the remaining squared norms
uniformly small as $L'\to\infty$. Gaussian quadratic forms of the cross vector have conditional
variance $O(N^{-1})$. Hence every compressed Schur entry converges to a deterministic weighted
average, which is continuous under the profile changes above and tends to its value at the anchor
as the radius and the windows shrink. If no uniform negative median existed in any small ball and
window, we could select radii and widths tending to zero, $N\to\infty$, and failing points $y$; along
this sequence the transverse gap and every Schur entry converge to their capped anchor values, the
capped anchor gap of Step~2 makes the limiting Schur complement strictly negative, and the block
inverse identity gives a full gap, a contradiction. So for some $g_2,\rho_0>0$, depending on $K$,
\[
 \operatorname{median}_{\Prob_{x,h'}}\lambda_{\max}\bigl(\mathcal M_{N,K}(y)\bigr)\le-4g_2
 \qquad(|y-x|\le\rho_0\sqrt N)
\]
for all sufficiently narrow windows and large $N$.

\emph{Step 4: concentration and a net.} For a deterministic symmetric $S$, covariance
differentiation gives the bound \eqref{eq:matrix-hess-cov} of
Lemma~\ref{lem:matrix-third}\ref{lem:matrix-third-a},
$\Var\tr(S\nabla^2\HN(y))\le\frac CN\tr S^2$, and conditioning decreases this variance. So
$\lambda_{\max}(\mathcal M_{N,K}(y))$ is a $C/\sqrt N$-Lipschitz function of the remaining standard
Gaussian variables, and, uniformly over deterministic $y$ in the ball,
\begin{equation}\label{eq:kr-pointwise}
 \Prob_{x,h'}\bigl(\lambda_{\max}(\mathcal M_{N,K}(y))>-3g_2\bigr)\le2e^{-c_KN}.
\end{equation}
Fix a mesh $s>0$ and round each coordinate of $y-x$ toward zero to a multiple of $s$. The rounded
point stays in the cube and in the ball, is within $s$ of $y$ in every coordinate, and has at most
$\rho^2N/s^2$ nonzero rounded coordinates. For $\rho^2/s^2<\frac12$ the set of rounded points has
cardinality at most
\[
 \sum_{j\le\rho^2N/s^2}\binom Nj\lceil4/s\rceil^j
 \le\exp\Bigl\{N\Bigl[H_2(\rho^2/s^2)+\frac{\rho^2}{s^2}\log\lceil4/s\rceil+o(1)\Bigr]\Bigr\},
\]
with $H_2(p)=-p\log p-(1-p)\log(1-p)$. Between a point and its rounding the Hessian of $\HN$
moves by at most $Cs$ (Lemma~\ref{lem:kr-c3}) and the other terms by at most $\omega_K(s)$, a
modulus from Step~1. Choose, in this order, the cap $K$ of Step~2 (which fixes $g_2$, $\rho_0$,
$c_K$), a mesh $s$ with $Cs+\omega_K(s)<g_2$, a radius $\rho\le\min(\rho_0,\rho_1)$ with
$\rho^2/s^2<\frac12$ and entropy below $c_K/2$, and then narrow windows and large $N$. The union
bound over the rounded points and \eqref{eq:kr-pointwise} give $\mathcal M_{N,K}(y)\preceq-2g_2I$ for
all $|y-x|\le\rho\sqrt N$, with exponentially small failure probability, and
\eqref{eq:kr-envelope} proves~\ref{prop:kr-local-a} with $\rho_*=\rho$ and $g_*=2g_2$.

\emph{Step 5: the typical determinant.} Let $\rho_K$ be the limit law of the capped diagonal
entries $\min(U^{-1},K)+a\chi_\pi$, and let $\mu_K$ be the reflection of
$\Pastur{a}{\rho_K}$ (Lemma~\ref{lem:hess-pastur}), the limit of the empirical spectral
distribution of $A_{N,K}$. Its Stieltjes transform $m_K(z)=\int\mu_K(d\lambda)/(\lambda-z)$ solves
$m_K(z)=\int\rho_K(dD)/(D-z-am_K(z))$
(Lemma~\ref{lem:matrix-comparison}\ref{lem:matrix-comparison-c}), and $\supp\mu_K\subset[g_0',\infty)$ by the gap. On the real
half-line below $\supp\mu_K$ put
\[
 L_K(z)=\int\log\bigl(D-z-am_K(z)\bigr)\,\rho_K(dD)+\frac a2m_K(z)^2 .
\]
The logarithms are defined: with $\varpi(z)=z+am_K(z)$, which is real and strictly increasing
there and tends to $-\infty$ as $z\to-\infty$, the imaginary part of the equation at
$z+\mathrm i\epsilon$ and Fatou's lemma give $a\int(D-\varpi(z))^{-2}\rho_K(dD)\le1$; integrating
over $z$ in a compact interval $J$ and substituting $\varpi(z)$, whose derivative is bounded on
$J$, shows that $\rho_K$ charges no point of the interior of $\varpi(J)$, since
$\int_I(D-x)^{-2}dx=\infty$ for $D$ in the interior of an interval $I$. So $\rho_K$-almost every
$D$ exceeds $\varpi(z)$, and the equation extends to the real half-line by dominated convergence.
Differentiating and using the equation gives $L_K'=-m_K$, and $L_K(z)=\log(-z)+o(1)$ as
$z\to-\infty$; these are also the derivative and the normalization of
$\int\log(\lambda-z)\,\mu_K(d\lambda)$, so the two functions coincide. Spectral convergence and the
gap give $N^{-1}\log\det A_{N,K}\to L_K(0)$. To remove the cap, increase the capped entries one at a
time: by the matrix determinant lemma and the gap,
\[
 0\le\log\det A_N-\log\det A_{N,K}\le\sum_{i:\,u_i^{-1}>K}\log\Bigl(1+\frac{u_i^{-1}-K}{g_0'}\Bigr).
\]
Each summand is at most $C(1+|b_i|)$ and vanishes unless $|b_i|\ge c'\log K-C$, so \eqref{eq:kr-cap}
makes the normalized sum uniformly small for large $K$. Real roots of the capped equation need not
be unique, so we identify $m_K(0)$ through traces: $N^{-1}\tr A_{N,K}^{-1}\to m_K(0)$, and by the
capped form of the inverse equivalent
(Lemma~\ref{lem:matrix-gap}\ref{lem:matrix-gap-c}--\ref{lem:matrix-gap-d}) the same traces converge
to a fixed point $\chi_\pi+O(K^{-1})$. Hence $m_K(0)\to\chi_\pi$, the logarithm arguments
$\min(U^{-1},K)+a(\chi_\pi-m_K(0))$ are at least $\frac12$ for large $K$, and dominated convergence,
with $|\log U|\le C(1+|b|)$, gives $L_K(0)\to\mathfrak d_\pi$ as $K\to\infty$. Finally, the negative of
the conditional Hessian differs from $A_N$ by a matrix of bounded rank $r$, and both have a gap
$g_*'>0$. By \eqref{eq:kr-cap}, $\max_i|b_i|\le t_0^{-1}\log(NL_0)$, so both have operator norm at
most $CN^C$, and eigenvalue interlacing bounds the difference of their log-determinants by
$2r\log(CN^C/g_*')=O(\log N)$. The other one-sided Hessian differs by a rank-one term and keeps the
gap. This proves~\ref{prop:kr-local-b} along the sequence.

\emph{Step 6: uniformity.} If either double limit in the proposition were positive, there would
be $\delta_*>0$, widths $w_k\downarrow0$, dimensions $N_k\uparrow\infty$ and admissible data of
width at most $w_k$ at dimension $N_k$ whose conditional failure probabilities exceed $\delta_*$.
These data obey \eqref{eq:kr-cap} with the same $L_0$, and their radii, empirical field laws and
energies converge to $q$, $\pi$ and $h$. Steps~1--5 apply along this sequence, so the failure
probabilities tend to zero, a contradiction.
\end{proof}

\subsection{A conditional competitor gap}\label{app:kr-gap}

Proposition~\ref{prop:kr-local} isolates a conditioned critical point from critical points at
overlap close to $q$. Here we bound $F_\zeta$ at competitors at a fixed overlap $r\in(0,q)$, for an
arbitrary empirical field law, since $\pi_\eta$ violates \eqref{eq:kr-ward}; energy and field law
enter only through an endpoint term. This is the one-level version of the comparison behind
\cite[Proposition~3.9]{Bou26b}, with windows of positive width.

For $0<r<q$ put
\[
 \lambda_r=\frac{\xi'(r)}c,\qquad
 \mathsf v_r=\xi'(r)-\frac{\xi'(r)^2}c,\qquad
 \bar{\mathsf v}_r=c-\xi'(r),\qquad
 \mathsf w_r=\int_r^qt\xi''(t)\,dt=\vth(q)-\vth(r),
\]
and $\mathcal B(r,b)=\E\,\Phi(r,\lambda_rb+Z)$ with $Z\sim\mathcal N(0,\mathsf v_r)$. For a law $\pi$
with finite second moment and $\E_\pi M^2=q$, and $h\in\R$, define
\begin{align}
 \mathcal C_{\pi,h}(r)&=(A(r)-1)h+B(r)d_\pi+\Bigl(1-\frac{\lambda_rr}q\Bigr)k_\pi,
 \label{eq:kr-endpoint}\\
 J_{\pi,h}(r)&=\E_\pi\bigl[\mathcal B(r,b)-\Phi(q,b)\bigr]+\mathcal C_{\pi,h}(r)
 -\frac\theta2\,\mathsf w_r .\label{eq:kr-J}
\end{align}
The same formulas define $J$ for an empirical law. Conditioning means Gaussian conditioning on the
displayed values of $\HN(x)$ and $\nabla\HN(x)$.

\begin{lemma}[Conditional comparison and a quadratic gap]\label{lem:kr-comparison}
\leavevmode
\begin{thmlist}[beginpenalty=10000]
\item\label{lem:kr-comparison-a} \textup{(Exact slices.)} Let $x\in\ocube$ with $q_x=q$ and
 condition on \eqref{eq:kr-conditioning}. Then, for $0<r<q$,
 \begin{equation}\label{eq:kr-exact-comparison}
  \frac1N\E_{\rm cond}\sup\bigl\{F_\zeta(y)-F_\zeta(x):\ y\in\cube,\ q_y=q,\ \ov xy=r\bigr\}
  \le J_{\pi_x,h}(r),
 \end{equation}
 where the supremum over the empty set is $-\infty$.
\item\label{lem:kr-comparison-b} \textup{(Windows.)} Fix a compact interval $I\subset(0,q)$, a law
 $\pi$ as above and $h\in\R$. For $\omega>0$, consider anchors with $|q_x-q|\le\omega$,
 $W_2(\pi_x,\pi)\le\omega$, conditioned on $\HN(x)=Nh'$ with $|h'-h|\le\omega$ and on the gradient
 value given by $\nabla F_\zeta(x)=0$, and replace the competitor constraints by $|q_y-q|\le\omega$
 and $|\ov xy-r|\le\omega$. Then the left side of \eqref{eq:kr-exact-comparison} is at most
 $J_{\pi,h}(r)+o_\omega(1)$, where $o_\omega(1)$ depends only on $I$, $\xi$, $\zeta$, $\pi$ and $h$,
 tends to zero as $\omega\downarrow0$, and is uniform in $N$, $r\in I$ and the conditional data.
\item\label{lem:kr-comparison-c} \textup{(Gap.)} For $0<r<q$,
 \begin{equation}\label{eq:kr-baseline-J}
  J_{\pi_0,h_0}(r)=-\theta\bigl(\Obs_\zeta(r)-\Obs_\zeta(q)\bigr).
 \end{equation}
 If $\Plef_0=0$, there is $L_J$ such that, for every $L\ge L_J$ in \eqref{eq:kr-heta} and every
 $\eta\in(0,\eta_0]$,
 \begin{equation}\label{eq:kr-tilted-gap}
  J_{\pi_\eta,h_\eta}(r)\le-\theta\bigl(\Obs_\zeta(r)-\Obs_\zeta(q)\bigr)-\eta\,(q-r)^2,
  \qquad 0<r<q .
 \end{equation}
\end{thmlist}
\end{lemma}

\begin{proof}[Proof of \ref{lem:kr-comparison-a}]
\emph{Regression.} The covariance of $(\HN(x),x\cdot\nabla\HN(x))$ is $NS_0$, the projection of
$\nabla\HN(x)$ onto $x^\perp$ is independent of this pair with covariance $c$ times the identity
on $x^\perp$, and for $\ov xy=r$
\[
 N^{-1}\Cov\bigl(\HN(y),(\HN(x),x\cdot\nabla\HN(x))\bigr)=\bigl(\xi(r),r\xi'(r)\bigr).
\]
Inverting $S_0$ and then regressing on the transverse gradient gives
\begin{equation}\label{eq:kr-condmean}
 \frac1N\E_{\rm cond}\HN(y)=A(r)h+B(r)d_N+\lambda_r\Bigl(\frac{b\cdot y}N-\frac rqk_N\Bigr),
\end{equation}
with the actual empirical moments $k_N=\E_{\pi_x}bM$ and $d_N=k_N+aq\chi_N$. Let $G_x$ be the centred
conditional field. On the slice $\{q_y=q,\ \ov xy=r\}$ its covariance is
$N[\xi(\ov yz)-\xi'(r)^2\ov yz/c]$ up to a constant independent of $y,z$. With
\[
 K_r(t)=\xi(t)-\frac{\xi'(r)^2}ct+r\xi'(r)-\xi(r),
\]
this gives
\begin{equation}\label{eq:kr-increments}
 \frac1{2N}\E\bigl(G_x(y)-G_x(z)\bigr)^2=K_r(q)-K_r(\ov yz).
\end{equation}
Since $0<\xi'(r)<c$, the numbers $\mathsf v_r=K_r'(r)$, $\bar{\mathsf v}_r=K_r'(q)-K_r'(r)$ and
$\mathsf w_r=qK_r'(q)-K_r(q)$ are positive, and $K_r(r)=rK_r'(r)$.

\emph{A Poisson comparison field.} Let $(\Pi_\alpha)$ be the points of a Poisson process on $\R$ with
intensity $\theta e^{-\theta u}\,du$. Independently, let $Z_i\sim\mathcal N(0,\mathsf v_r)$,
$g_{i\alpha}\sim\mathcal N(0,\bar{\mathsf v}_r)$ and $Y_\alpha\sim\mathcal N(0,N\mathsf w_r)$ be
independent; the $Z_i$ are common to all leaves $\alpha$. With
$\Psi(y)=\sum_i\{\lambda_rb_iy_i+\Phistar(q,y_i)\}$, compare
\[
 X(y,\alpha)=G_x(y)+Y_\alpha+\Pi_\alpha+\Psi(y),\qquad
 \widetilde X(y,\alpha)=\sum_iy_i(Z_i+g_{i\alpha})+\Pi_\alpha+\Psi(y).
\]
For two indices put $q_{\alpha\beta}=q$ if $\alpha=\beta$ and $q_{\alpha\beta}=r$ otherwise. By
\eqref{eq:kr-increments}, the squared increment of $\widetilde X$ minus that of $X$ is
\begin{equation}\label{eq:kr-bregman}
 2N\bigl\{K_r(R)-K_r(q_{\alpha\beta})-K_r'(q_{\alpha\beta})(R-q_{\alpha\beta})\bigr\}\ge0,
 \qquad R=\ov yz .
\end{equation}
For equal leaves this is $2NK_r'(q)(q-R)-2N[K_r(q)-K_r(R)]$; for distinct leaves it uses
$\mathsf w_r=qK_r'(q)-K_r(q)$ and $K_r(r)=rK_r'(r)$. The inequality holds because
$K_r''=\xi''\ge0$ on $[-1,1]$, including negative values of $R$; this is the only place where the
proof of Proposition~\ref{prop:kr-main} uses the convexity of $\xi$ on $[-1,0]$
(Proposition~\ref{prop:kr-retained}\ref{prop:kr-retained-b} also uses it, through
Proposition~\ref{prop:max-fm}\ref{prop:max-fm-a}; see Remark~\ref{rem:kr-window}). Conditionally on the scores $\Pi_\alpha$, apply the
Sudakov--Fernique inequality to finite subsets of the slice and finitely many leaves; the means of
$X$ and $\widetilde X$ agree. Passing to all leaves is legitimate because, after subtracting one
fixed variable, the increasing maxima are nonnegative and their limits have finite expectation by
\eqref{eq:kr-marking} below; passing to the whole slice uses continuity and the integrability of
Gaussian suprema. So $\E\sup_{y,\alpha}X\le\E\sup_{y,\alpha}\widetilde X$.

\emph{Marking.} If $V_\alpha$ are i.i.d., independent of $(\Pi_\alpha)$, with
$\E e^{\theta V}<\infty$, then Lemma~\ref{lem:poisson-shift} (with a trivial mark space) shows that
$\sup_\alpha(\Pi_\alpha+V_\alpha)-\theta^{-1}\log\E e^{\theta V}$ has the law of $\sup_\alpha \Pi_\alpha$,
a Gumbel variable with finite mean. Hence
\begin{equation}\label{eq:kr-marking}
 \E\sup_\alpha(\Pi_\alpha+V_\alpha)=\E\sup_\alpha \Pi_\alpha+\theta^{-1}\log\E e^{\theta V} .
\end{equation}
On the left, the supremum separates into $\sup_y\{G_x(y)+\Psi(y)\}$ plus $\sup_\alpha(Y_\alpha+\Pi_\alpha)$,
whose marking cost is $\theta N\mathsf w_r/2$. On the right, enlarge the slice to the cube and use
$\sup_{m\in[-1,1]}\{tm+\Phistar(q,m)\}=\Phi(q,t)$ coordinatewise. Conditionally on the root vector
$(Z_i)$, the leaf marks $\sum_i\Phi(q,Z_i+\lambda_rb_i+g_{i\alpha})$ are i.i.d., and since
$\cdf\zeta=\theta$ on $(r,q)$, the Cole--Hopf form of the Parisi equation gives
\[
 \frac1\theta\log\E_g\exp\Bigl\{\theta\sum_i\Phi(q,Z_i+\lambda_rb_i+g_i)\Bigr\}
 =\sum_i\Phi(r,Z_i+\lambda_rb_i);
\]
the exponential moments are finite since $|\Phi(q,t)|\le C+|t|$. Applying \eqref{eq:kr-marking}
on both sides and cancelling $\E\sup_\alpha \Pi_\alpha$,
\[
 \E\sup_y\bigl\{G_x(y)+\Psi(y)\bigr\}\le\sum_i\mathcal B(r,b_i)-\frac{N\theta}2\mathsf w_r .
\]
Finally, \eqref{eq:kr-condmean}, the identity $\Phistar(q,x_i)=\Phi(q,b_i)-b_ix_i$ and the
cancellation of the equal Onsager integrals at $q_y=q_x$ give \eqref{eq:kr-exact-comparison}.
\end{proof}

\begin{proof}[Proof of \ref{lem:kr-comparison-b}]
Let $Q=q_x$ and $t_y=\ov xy$. The same regression gives, for all $y,z$,
\[
 \frac1N\E\,G_x(y)G_x(z)=\xi(\ov yz)-\mathbf s(t_y)^{\mathsf T}S_0(Q)^{-1}\mathbf s(t_z)
 -\frac{\xi'(t_y)\xi'(t_z)}{\xi'(Q)}\Bigl(\ov yz-\frac{t_yt_z}Q\Bigr),
\]
with $\mathbf s(t)=(\xi(t),t\xi'(t))^{\mathsf T}$ and $S_0(Q)$ the covariance matrix of
$N^{-1/2}(\HN(x),x\cdot\nabla\HN(x))$ (as in the proof of Lemma~\ref{lem:kr-c3}). Both $S_0(Q)$ and $\xi'(Q)$ are uniformly
invertible near $q$, so on the competitor window
\[
 \Bigl|\frac1{2N}\E\bigl(G_x(y)-G_x(z)\bigr)^2-K_r(q)+K_r(\ov yz)\Bigr|\le C\omega ;
\]
here and below $C$ depends on $I$, $\xi$, $\zeta$, $\pi$ and $h$ only. The conditional mean is controlled
without any identity specific to the law: the selected gradient is
$\nabla\HN(x)=b^{(Q)}+\xi''(Q)\chi_xx$, and by Lemma~\ref{lem:kr-endpoints}
$\|b^{(Q)}-b\|/\sqrt N\le C\omega$ and $|\chi_x-\chi_N|\le C\omega$, while the Wasserstein window
bounds $N^{-1}|b|^2$. Replacing $Q,t_y$ by $q,r$ and the actual-radius quantities by their
fixed-radius values in the conditional mean costs at most $CN\omega$, and the conjugate
$\Phistar(\cdot,m)$ and the Onsager integral change by at most $CN\omega$ when the radius is frozen,
because $|\partial_t\Phi|$ is bounded and the infimum defining $\Phistar$ preserves this bound.

Fix $\rho\in(0,1)$ and a Euclidean $\rho\sqrt N$-net $y_1,\dots,y_n$ of the competitor window, with
$n\le(1+2/\rho)^N$. On the net the increment difference \eqref{eq:kr-bregman} now has an error of
at most $CN\omega$, which includes the discrepancies in the self-overlaps of the competitors. Add
independent Gaussians $\epsilon_j+\epsilon'_\alpha$, each of variance $C_0N\omega$, to
$\widetilde X(y_j,\alpha)$. Distinct indices differ in the physical or the leaf index, so every
squared increment gains at least $2C_0N\omega$, and the comparison holds without error for $C_0$
large. The physical noise costs at most $\E\max_j\epsilon_j\le\sqrt{2C_0N\omega\log n}$. The leaf
noise has an exact cost: replacing $\Pi_\alpha$ by $\Pi_\alpha+\epsilon'_\alpha$ translates the Poisson
process by $\theta C_0N\omega/2$, by \eqref{eq:kr-marking}, and the other leaf marks keep their
law, also conditionally on the root field. After division by $N$ the two noises cost at most
$C\sqrt{\omega\log(1+2/\rho)}+C\omega$, whatever the number of leaves.

It remains to remove the net. Conditioning reduces increment variances, and the canonical metric
of $\HN$ is at most $C|y-z|$ on the cube, so the Gaussian entropy integral gives
$\E\sup_y|G_x(y)-G_x(y_{j(y)})|\le CN\rho\sqrt{\log(C/\rho)}$. The linear part of the endpoint
changes by at most $CN\rho$. For the conjugate part, let $\Omega$ be a modulus of continuity of
$\Phistar(q,\cdot)$ on $[-1,1]$ and $C_\Phi$ a bound on its absolute value; splitting coordinates
at $|y_i-z_i|=\sqrt\rho$ gives
$N^{-1}\sum_i|\Phistar(q,y_i)-\Phistar(q,z_i)|\le\Omega(\sqrt\rho)+2C_\Phi\rho$ when
$|y-z|\le\rho\sqrt N$. On the comparison side we again enlarge the net to the cube. Finally the
endpoint averages are close to their $\pi$-averages: $b\mapsto\mathcal B(r,b)-\Phi(q,b)$ is
uniformly Lipschitz, $U$ is bounded and Lipschitz, $b\mapsto bM(b)$ has derivative at most
$C(1+|b|)$, and Cauchy--Schwarz under a $W_2$-coupling bounds the differences by $C\omega$.
Altogether the normalized error is at most
\[
 C\omega+C\sqrt{\omega\log(1+2/\rho)}+C\rho\sqrt{\log(C/\rho)}+C\rho+\Omega(\sqrt\rho)+2C_\Phi\rho ,
\]
which is small if $\rho$ is small and then $\omega$ is small. All estimates precede the limit in
$N$.
\end{proof}

\begin{proof}[Proof of \ref{lem:kr-comparison-c}]
\emph{The baseline.} The regression coefficients satisfy
$cA(r)+(c+aq)B(r)=r\xi'(r)/q$ and $\xi(q)A(r)+qcB(r)=\xi(r)$. With \eqref{eq:kr-k0} and
\eqref{eq:kr-h0}, every $\chi_0$-term of \eqref{eq:kr-endpoint} cancels and
\[
 \mathcal C_{\pi_0,h_0}(r)=\theta\bigl\{qc-\xi(q)-r\xi'(r)+\xi(r)\bigr\}=\theta\,\mathsf w_r .
\]
On the plateau, $e^{\theta\Phi}$ is space-time harmonic for the heat flow with clock $\xi'$, and
the Auffinger--Chen process on $[0,q]$ is its Doob transform. The density change depends only on
the endpoint, so given $X_q=b$ the law of $X_r$ is the Gaussian bridge with mean $\lambda_rb$ and
variance $\mathsf v_r$. Hence $\E_0\mathcal B(r,b)=\E\Phi(r,X_r)$, and It\^o's formula gives
$\frac d{dt}\E\Phi(t,X_t)=\frac\theta2\xi''(t)\Gz\zeta(t)$ on $(0,q)$. Therefore
\[
 J_{\pi_0,h_0}(r)=-\frac\theta2\int_r^q\xi''(t)\bigl(\Gz\zeta(t)-t\bigr)\,dt
 =-\theta\bigl(\Obs_\zeta(r)-\Obs_\zeta(q)\bigr),
\]
which is \eqref{eq:kr-baseline-J}.

\emph{The double contact at $q$.} Recall from before Lemma~\ref{lem:kr-energy} that $A(q)=1$,
$A'(q)=0$, $B(q)=0$ and $B'(q)=q^{-1}$. Although $\mathsf v_q=0$, the heat-semigroup rule, which differentiates in the
variance and not in its square root, gives
$\partial_r\mathcal B(r,b)=\E\bigl[\Phi_t+\lambda_r'b\,\Phi_x+\frac12\mathsf v_r'\Phi_{xx}\bigr](r,\lambda_rb+Z)$.
As $r\uparrow q$, $\mathsf v_r'\to-a$ and $\lambda_r'\to a/c$, and the Parisi equation on the plateau
gives $\Phi_t(q-,b)=-\frac a2(U+\theta M^2)$. So
\[
 \partial_r\mathcal B(q-,b)=-aU(b)-\frac{a\theta}2M(b)^2+\frac acbM(b),\qquad
 \mathcal C_{\pi,h}'(q)=a\chi_\pi-\frac ack_\pi ,
\]
and for every $\pi$ with $\E_\pi M^2=q$ and every $h$,
\begin{equation}\label{eq:kr-double-contact}
 J_{\pi,h}(q)=0,\qquad J_{\pi,h}'(q-)=0 .
\end{equation}
A second heat differentiation expresses $\partial_r^2\mathcal B(r,b)$ as an average at
$(r,\lambda_rb+Z)$ of $\Phi_{tt},\Phi_{tx},\Phi_{txx}$ and of $\Phi_x,\dots,\Phi_{xxxx}$, with
coefficients built from $\lambda_r',\lambda_r'',\mathsf v_r',\mathsf v_r''$ and $b$, at most quadratic in
$b$; on the plateau the Parisi equation expresses the time derivatives through spatial derivatives
of order at most four. By Lemma~\ref{lem:kr-contact}\ref{lem:kr-contact-d} these are bounded on
$[0,q]\times\R$, so $|\partial_r^2\mathcal B(r,b)|\le C(1+b^2)$ for $0<r<q$, which also justifies the
left derivatives at $q$ by dominated convergence.

\emph{The tilt and the shift.} The map $\pi\mapsto J_{\pi,h_0}$ is affine, and
$\pi_\eta-\pi_0=-\eta\mathcal R\pi_0$ with $\mathcal R$ bounded, while $\pi_0$ has moments of all
orders. So $\sup_{0<r<q}|(J_{\pi_\eta,h_0}-J_{\pi_0,h_0})''(r)|\le C_J\eta$ for $0<\eta\le\eta_0$.
Integrating twice from $q$ and using \eqref{eq:kr-double-contact},
$|J_{\pi_\eta,h_0}(r)-J_{\pi_0,h_0}(r)|\le C_J\eta(q-r)^2$. The energy shift contributes exactly
$L\eta(A(r)-1)\le-L\eta(1-r/q)^2$, by Lemma~\ref{lem:kr-energy}. With $L_J=q^2(C_J+1)$,
\eqref{eq:kr-baseline-J} gives \eqref{eq:kr-tilted-gap}.
\end{proof}

We now fix the energy shift:
\begin{equation}\label{eq:kr-L}
 L=\max\{L_S,L_J\},
\end{equation}
with $L_S$ from \eqref{eq:kr-schur-shift}. A tilt $\eta$ is \emph{admissible} if $\eta=0$ when
$\Plef_0>0$, and $\eta\in(0,\eta_0]$ when $\Plef_0=0$. For every admissible $\eta$ the data
$(\pi_\eta,h_\eta)$ satisfy the hypotheses of Proposition~\ref{prop:kr-local}, by
\eqref{eq:kr-schur} and \eqref{eq:kr-schur-shift}, and
$J_{\pi_\eta,h_\eta}(r)\le-\theta(\Obs_\zeta(r)-\Obs_\zeta(q))<0$ for $0<r<q$, by
Lemma~\ref{lem:kr-comparison}\ref{lem:kr-comparison-c} and
Lemma~\ref{lem:kr-contact}\ref{lem:kr-contact-b}.

\paragraph{The order of the parameters.}
At a fixed admissible tilt, Proposition~\ref{prop:kr-local} first selects a radius $\rho_\eta$.
Then a cutoff $\delta$ is chosen so small that competitors with overlap at least $q-\delta$ lie in
the ball of that radius. Overlaps below $\delta$ are not excluded; the spherical code of
Appendix~\ref{app:kr-count} handles them. On the remaining compact interval $[\delta,q-\delta]$,
$J_{\pi_\eta,h_\eta}$ has a strictly negative maximum. Only then are the competitor width $\omega$
and the empirical windows chosen. Conditional variances are at most $N\xi(1)$, so Gaussian
concentration and a finite cover of $[\delta,q-\delta]$ turn the gap in conditional expectation
into an exponentially small probability for competitors (Lemma~\ref{lem:kr-marks}).

\subsection{The marked count and existence}\label{app:kr-count}

We count critical points of $F_\zeta$ with data close to an admissible pair $(\pi_\eta,h_\eta)$
that pass two marks: strong concavity on a ball, and a strict gap against competitors at overlaps
in $[\delta,q-\delta]$. The expected count has exponent $\Sigma(f)-O(\eta)$, the marks give distinct
counted points overlaps below the cutoff $\delta'$ of~\ref{step:kr-par4}, so a one-sided spherical
code bounds their number, and the ratio bounds the existence probability.

\paragraph{The value of $F_\zeta$ at the selected data.}
By \eqref{eq:matrix-split}, $F_\alpha(x)=\HN(x)+NV_\alpha(q_x,\mu_x)$ and
$\FTAP(x)=\HN(x)+NV_{\mathrm{TAP}}(q_x,\mu_x)$, with $V_\alpha$ as in \eqref{eq:hess-V} and
$V_{\mathrm{TAP}}=\inf_\beta V_\beta$. We use the modulus $\varrho(t)=t\log(2e/t)$ and the constants
$C_*,L_1,L_2$ of Lemma~\ref{lem:matrix-correction}; $\varrho$ is nondecreasing on $[0,2]$ and
$\varrho(t)\ge t$ there. Put $\|\mathcal R\|_\infty=0$ if $\Plef_0>0$, and
\[
 C_V=4L_1+2L_2+4C_*\|\mathcal R\|_\infty .
\]

\begin{lemma}[Values at the selected data]\label{lem:kr-correction}
Let $0\le\eta\le\eta_0$, with $\eta=0$ if $\Plef_0>0$, let $x\in\ocube$, and let
$w=W_2(\pi_x,\pi_\eta)\le2$. Then
$\Wone(\mu_x,M_\#\pi_\eta)\le w$, $|q_x-q|\le2w$, and
\begin{thmlist}
\item\label{lem:kr-correction-a}
 $|V_\zeta(q_x,\mu_x)-V_\zeta(q,M_\#\pi_\eta)|\le C_V\varrho(w)$;
\item\label{lem:kr-correction-b} $0\le\bigl(F_\zeta(x)-\FTAP(x)\bigr)/N\le C_V\bigl(\eta+\varrho(w)\bigr)$.
\end{thmlist}
Moreover the value at the selected data,
\begin{equation}\label{eq:kr-feta}
 f_\eta=h_\eta+V_\zeta(q,M_\#\pi_\eta),\qquad\text{satisfies}\qquad
 f_0=f,\qquad |f_\eta-f|\le(L+C_V)\,\eta .
\end{equation}
\end{lemma}

\begin{proof}
The first two assertions are in Corollary~\ref{cor:matrix-tilted} (with $Q=q$ and $\pi=\pi_\eta$),
and~\ref{lem:kr-correction-a} follows from them and
Lemma~\ref{lem:matrix-correction}\ref{lem:matrix-correction-b}, since $2L_1w+L_2\varrho(w)\le
C_V\varrho(w)$. Corollary~\ref{cor:matrix-tilted} also gives~\ref{lem:kr-correction-b}, because the
total variation distance between $\pi_\eta=(1-\eta\mathcal R)\pi_0$ and $\pi_0$ is at most
$\eta\|\mathcal R\|_\infty$; the corollary rests on the identity
$V_{\mathrm{TAP}}(q,M_\#\pi_0)=V_\zeta(q,M_\#\pi_0)$ at the contact $q$
(Lemma~\ref{lem:matrix-band}). Finally, by the Fenchel identity
$\Phistar(q,M(b))=\Phi(q,b)-bM(b)$, \eqref{eq:kr-level} says $h_0+V_\zeta(q,M_\#\pi_0)=f$; and
$|V_\zeta(q,M_\#\pi_\eta)-V_\zeta(q,M_\#\pi_0)|\le2C_*\eta\|\mathcal R\|_\infty$ because the conjugates
are bounded by $C_*$ (Lemma~\ref{lem:matrix-correction}\ref{lem:matrix-correction-a}), while
$h_\eta-h_0=L\eta$.
\end{proof}

\paragraph{A marked Kac--Rice formula.}
The atoms of $\zeta$ may be dense, so $F_\zeta$ need not be twice differentiable everywhere. Its
gradient is locally Lipschitz on $\{x\in\ocube:q_x>0\}$ (after \eqref{eq:kr-gradient}), and that
suffices.

\begin{lemma}[Marked Kac--Rice formula]\label{lem:kr-marked}
Let $\mathcal A$ be a Borel subset of $\{x\in\ocube:q_x>0\}$ and let $\mathsf Q(x,\HN)\ge0$ be a
measurable mark. Then
\[
 \E\sum_{x\in\mathcal A:\,\nabla F_\zeta(x)=0}\mathsf Q(x,\HN)
 =\int_{\mathcal A}p_{\nabla F_\zeta(x)}(0)\,
 \E\bigl[|\det\nabla^2F_\zeta(x)|\,\mathsf Q(x,\HN)\bigm|\nabla F_\zeta(x)=0\bigr]\,dx,
\]
also when both sides are infinite. Almost surely no critical point in $\{q_x>0\}$ lies on a sphere
$\{q_x=t\}$ with $t$ an atom of $\zeta$, so the Hessian on the right is the ordinary one for almost
every $x$. The conditional expectation is the canonical one: given $\nabla F_\zeta(x)=0$, the field
$\HN$ has the law of the Gaussian regression of $\HN$ on $\nabla\HN(x)$ evaluated at the gradient
value \eqref{eq:kr-gradient}.
\end{lemma}

\begin{proof}
Fix $x_*\neq0$. The covariance $\Sigma_*$ of $\nabla\HN(x_*)$ is positive definite. Put
$\mathbf Z=\Sigma_*^{-1/2}\nabla\HN(x_*)$, $\psi_j(x)=\E[\HN(x)\mathbf Z_j]$, so that
$\HN(x)=\psi(x)\cdot\mathbf Z+\widetilde H_N(x)$ with $\widetilde H_N$ independent of $\mathbf Z$.
The matrix $\mathcal J(x)_{ij}=\partial_i\psi_j(x)$ equals $\Sigma_*^{1/2}$ at $x_*$, hence is invertible
on a deterministic neighbourhood $U_*$. Conditionally on $\widetilde H_N$, let
$b_*(x)=\nabla\widetilde H_N(x)+\nabla(F_\zeta-\HN)(x)$ and $T_*(x)=-\mathcal J(x)^{-1}b_*(x)$, a
locally Lipschitz map; the critical equation is $\mathbf Z=T_*(x)$, and at almost every $x$,
$DT_*(x)=-\mathcal J(x)^{-1}\nabla^2F_\zeta(x)$ evaluated at $\mathbf Z=T_*(x)$. The area formula
for Lipschitz maps \cite[Section~3]{Federer}, extended to the nonnegative weight equal to the
Gaussian density at $T_*(x)$ times the mark by monotone convergence, equates the Gaussian average
of the marked number of solutions with the integral of that weight times $|\det DT_*|$. Averaging
over $\widetilde H_N$ and changing variables from $\mathbf Z$ to the gradient gives the formula on
$U_*$; countably many such charts and Tonelli's theorem give it on $\mathcal A$. The union of the
atom spheres is Lebesgue-null, a locally Lipschitz map between spaces of equal dimension maps null
sets to null sets, and $\mathbf Z$ is independent of $\widetilde H_N$; this gives the second
assertion. Finally, the conditional law produced by the argument is Gaussian with mean affine in
the gradient value and covariance independent of it, as is the regression of $\HN$ on
$\nabla\HN(x)$; two such versions agree for almost every gradient value, since the law of
$\nabla F_\zeta(x)$ is nondegenerate, hence for every value. Conditioning in addition on $\HN(x)$,
as below, is the same canonical regression on $(\HN(x),\nabla\HN(x))$.
\end{proof}

\paragraph{Parameters and marks.}
Fix the following, in this order.
\begin{enumerate}[label=\textup{(Par\arabic*)},ref=(Par\arabic*),leftmargin=*,widest=6,topsep=3pt,itemsep=2pt]
\item\label{step:kr-par1} An admissible tilt $\eta$, and $(\pi,h)=(\pi_\eta,h_\eta)$.
\item\label{step:kr-par2} A cap $L_0>\E_\pi e^{t_0|b|}$, $t_0=8$.
\item\label{step:kr-par3} The margin and the radius $g_\eta=g_*$, $\rho_\eta=\min\{\rho_*,\sqrt q/2\}$, with
 $g_*,\rho_*$ given by Proposition~\ref{prop:kr-local} for $(\pi,h,L_0)$.
\item\label{step:kr-par4} For $\delta\in(0,q/2)$, the effective cutoff
 $\delta'=\min\{\delta,\,q/4,\,\rho_\eta^2/4\}$ and the gap
 \[
  \gamma=\theta\min\bigl\{\Obs_\zeta(r)-\Obs_\zeta(q):\ \delta'\le r\le q-\delta'\bigr\}>0
 \]
 (Lemma~\ref{lem:kr-contact}\ref{lem:kr-contact-b}).
\item\label{step:kr-par5} A competitor width $\omega_c\in(0,\min\{\delta'/4,q/4\}]$ such that the error
 of Lemma~\ref{lem:kr-comparison}\ref{lem:kr-comparison-b}, for $I=[\delta',q-\delta']$ and $(\pi,h)$,
 is at most $\gamma/4$ for all widths in $(0,\omega_c]$.
\item\label{step:kr-par6} A window width $w\in(0,w_1]$, where $w_1=\min\{\omega_c/2,w_*\}$ and
 $C_V\varrho(w_*)\le\gamma/8$, and a measurable set $\mathcal W_N(w)\subset\ocube$.
\end{enumerate}
The \emph{retained points} are the $x\in\ocube$ satisfying:
\begin{enumerate}[label=\textup{(Ret\arabic*)},ref=(Ret\arabic*),leftmargin=*,widest=5,topsep=3pt,itemsep=2pt]
\item\label{step:kr-ret1} $\nabla F_\zeta(x)=0$;
\item\label{step:kr-ret2} $W_2(\pi_x,\pi)\le w$, the fields of $x$ satisfy \eqref{eq:kr-cap}, and
 $x\in\mathcal W_N(w)$;
\item\label{step:kr-ret3} $|\HN(x)/N-h|\le N^{-2}$;
\item\label{step:kr-ret4} $\nabla^2_\sigma F_\zeta(y)\preceq-g_\eta I$ for $\sigma\in\{\pm\}$ and every
 $y\in\ocube$ with $|y-x|\le\rho_\eta\sqrt N$;
\item\label{step:kr-ret5} $F_\zeta(y)\le F_\zeta(x)-\gamma N/2$ for every $y\in\mathcal K(x)$, where
 $\mathcal K(x)=\{y\in\cube:\ |q_y-q|\le\omega_c,\ \delta'\le\ov xy\le q-\delta'\}$.
\end{enumerate}
Let $\widehat C_N=\widehat C_N(\eta,\delta,w)$ be the number of retained points ($\widehat C_N=0$ on
the null event where $\HN$ is not of class $C^3$ near the cube). Mark~\ref{step:kr-ret4} is
measurable where~\ref{step:kr-ret2} holds; there $q_x\ge3q/4$ (Lemma~\ref{lem:kr-correction}
and~\ref{step:kr-par5}), so $\rho_\eta\sqrt N<|x|$ by~\ref{step:kr-par3} and the ball does not
contain $0$. By \eqref{eq:kr-hessian} and Lemma~\ref{lem:kr-endpoints}, at every point $y$ of the ball
other than its two points of extreme radius $y_\pm=x(1\pm\rho_\eta\sqrt N/|x|)$, each one-sided
Hessian is a limit of Hessians at points of the ball in a fixed countable dense subset of $\ocube$ whose radii are
not atoms of $\zeta$, approaching from the corresponding side; at $y_\pm$ (when they lie in $\ocube$)
both one-sided Hessians are given by \eqref{eq:kr-hessian} and are measurable functions of
$(x,\HN)$. So it suffices to impose~\ref{step:kr-ret4} on that set and at $y_\pm$. Mark~\ref{step:kr-ret5} is a supremum of a continuous function over a compact set. By
Lemma~\ref{lem:kr-correction}, \ref{step:kr-ret2} gives $|q_x-q|\le2w\le\omega_c$,
so no point belongs to its own competitor set $\mathcal K(x)$. For a fixed admissible $\eta$, the
margin and the radius in~\ref{step:kr-ret4} do not depend on $\delta$ or $w$.

\begin{lemma}[The marks]\label{lem:kr-marks}
Let $x\in\ocube$ satisfy~\ref{step:kr-ret2}.
\begin{thmlist}
\item\label{lem:kr-marks-a} If $x$ satisfies~\ref{step:kr-ret4}, then
 $\langle\nabla F_\zeta(y)-\nabla F_\zeta(y'),y-y'\rangle\le-g_\eta|y-y'|^2$ for all $y,y'$ in
 $B(x)=\{z\in\ocube:|z-x|\le\rho_\eta\sqrt N\}$.
\item\label{lem:kr-marks-b} If $|h'-h|\le N^{-2}$ and $N\ge\omega_c^{-1/2}$, then
 \[
  \Prob_{x,h'}\bigl(\text{\ref{step:kr-ret5} fails}\bigr)\le\frac q{\omega_c}
  \exp\Bigl\{-\frac{N\gamma^2}{32\,\xi(1)}\Bigr\}.
 \]
\item\label{lem:kr-marks-c} $\displaystyle\lim_{w\downarrow0}\limsup_{N\to\infty}\sup_{x,h'}
 \Prob_{x,h'}\bigl(\text{\ref{step:kr-ret4} fails}\bigr)=0$, the supremum being over $x$
 satisfying~\ref{step:kr-ret2} and $|h'-h|\le N^{-2}$.
\end{thmlist}
\end{lemma}

\begin{proof}
\ref{lem:kr-marks-a} The set $B(x)$ is convex, and for $z\in B(x)$,
$\sqrt{q_z}\ge\sqrt{q_x}-\rho_\eta\ge\sqrt{3q/4}-\sqrt q/2>0$. Let $y\neq y'$ in $B(x)$,
$z_t=y'+t(y-y')$ and $\phi(t)=\langle\nabla F_\zeta(z_t),y-y'\rangle$ for $t\in[0,1]$. The segment is a
compact subset of $\ocube\setminus\{0\}$, where $\nabla F_\zeta$ is locally Lipschitz, so $\phi$ is
Lipschitz and $\phi(1)-\phi(0)=\int_0^1\phi'(t)\,dt$. By
Lemma~\ref{lem:matrix-hessian}\ref{lem:matrix-hessian-b}, at every $t$ the right derivative of
$\phi$ is $\langle\nabla^2_\sigma F_\zeta(z_t)(y-y'),y-y'\rangle$, with $\sigma$ the sign of
$z_t\cdot(y-y')$, and by~\ref{step:kr-ret4} it is at most $-g_\eta|y-y'|^2$. At almost every $t$ the
derivative $\phi'(t)$ exists and equals the right derivative, and integration gives the claim.

\ref{lem:kr-marks-b} Put $r_j=\min\{\delta'+2j\omega_c,\,q-\delta'\}$ for
$0\le j\le J=\lceil(q-2\delta')/(2\omega_c)\rceil$, so that $\mathcal K(x)\subset\bigcup_jS_j$ with
$S_j=\{y\in\cube:|q_y-q|\le\omega_c,\,|\ov xy-r_j|\le\omega_c\}$, and $J+1\le q/\omega_c$. The anchor
data have width at most $\omega_c$, because $|q_x-q|\le2w\le\omega_c$, $W_2(\pi_x,\pi)\le w$ and
$|h'-h|\le N^{-2}\le\omega_c$. By Lemma~\ref{lem:kr-comparison}\ref{lem:kr-comparison-b}
and~\ref{lem:kr-comparison-c}, and the definition of $\gamma$,
$Y_j=\sup_{y\in S_j}\{F_\zeta(y)-F_\zeta(x)\}$ satisfies $\E_{x,h'}Y_j\le-\frac34\gamma N$ whenever
$S_j\neq\varnothing$. Under $\Prob_{x,h'}$, $F_\zeta(x)$ is a constant, and $Y_j$ is the supremum over
the compact set $S_j$ of a continuous function plus a centred Gaussian field $G_x$ with continuous
paths and $\Var G_x(y)\le\Var\HN(y)\le N\xi(1)$. Gaussian concentration for Lipschitz functions
\cite[(1.22)--(1.23)]{Ledoux}, applied to maxima over finite subsets increasing to a dense set and
passed to the limit by monotone convergence and Fatou's lemma, gives
$\Prob_{x,h'}(Y_j>-\frac12\gamma N)\le\exp\{-(\gamma N/4)^2/(2N\xi(1))\}$. If~\ref{step:kr-ret5}
fails, some $Y_j$ exceeds $-\frac12\gamma N$; the union bound over $J+1\le q/\omega_c$ indices
proves the claim.

\ref{lem:kr-marks-c} The data $(x,h')$ are admissible of width at most $w+N^{-2}$ in the sense of
Proposition~\ref{prop:kr-local}, and $\rho_\eta\le\rho_*$; part~\ref{prop:kr-local-a} of that
proposition gives the claim.
\end{proof}

\paragraph{The count exponent.}
For a law $\pi$ with $\E_\pi M^2=q$ and finite second moment, and $h\in\R$, put
\begin{equation}\label{eq:kr-SigmaKR}
 \Sigma_{\rm KR}(\pi,h)=-D\bigl(\pi\,\Vert\,\mathcal N(0,c)\bigr)
 +\frac{a(k_\pi-c\chi_\pi)^2}{2c(c+aq)}
 -\frac1{2\sigma_H^2}\Bigl(h-\frac{c\,d_\pi}{c+aq}\Bigr)^2 ,
\end{equation}
where $D$ is the relative entropy and $\sigma_H^2$ is \eqref{eq:kr-sigmaH}. The three terms come
from the field law, from the Gaussian density of the gradient together with the typical
determinant, and from the conditional density of the energy.

\begin{lemma}[Rate of the selected data]\label{lem:kr-normalization}
$\Sigma_{\rm KR}(\pi_0,h_0)=\Sigma(f)$, and there is $C_\Sigma$ with
$|\Sigma_{\rm KR}(\pi_\eta,h_\eta)-\Sigma(f)|\le C_\Sigma\eta$ for every admissible $\eta$.
\end{lemma}

\begin{proof}
By \eqref{eq:kr-law} and \eqref{eq:kr-ito},
$D(\pi_0\Vert\mathcal N(0,c))=\theta\,\E_0[\Phi(q,b)-\Phi(0,0)]=\theta(f-P_\zeta)+\frac{\theta^2}2\vth(q)$.
By \eqref{eq:kr-k0}, $k_0-c\chi_0=c\theta q$, and by \eqref{eq:kr-h0},
$h_0-cd_0/(c+aq)=\theta(\xi(q)-qc^2/(c+aq))=\theta\sigma_H^2$. So the two quadratic terms of
\eqref{eq:kr-SigmaKR} add up to
$\frac{\theta^2}2\bigl(\frac{acq^2}{c+aq}-\xi(q)+\frac{qc^2}{c+aq}\bigr)=\frac{\theta^2}2\vth(q)$, and
$\Sigma_{\rm KR}(\pi_0,h_0)=\theta(P_\zeta-f)=\Sigma(f)$ by
Lemma~\ref{lem:kr-contact}\ref{lem:kr-contact-a}. For $\eta>0$, the log-density ratio
$\log(1-\eta\mathcal R)$ is $O(\eta)$ uniformly, the original log-density grows at most linearly
in $|b|$, every other moment in \eqref{eq:kr-SigmaKR} is affine in $\eta$ with bounded slope, and
$h_\eta-h_0=L\eta$.
\end{proof}

\begin{lemma}[Expected marked count]\label{lem:kr-count}
Let $\mathcal W_N(w)$ be the set of $x\in\ocube$ whose fields satisfy
$|N^{-1}\sum_i\varphi(b_i)-\E_\pi\varphi|\le w$ for the five functions
$\varphi(b)=b^2,\ \log U(b),\ \log\frac{d\pi}{db}(b),\ U(b),\ bM(b)$. Then
\[
 \liminf_{w\downarrow0}\liminf_{N\to\infty}\frac1N\log\E\widehat C_N(\eta,\delta,w)
 \ge\Sigma_{\rm KR}(\pi_\eta,h_\eta).
\]
\end{lemma}

\begin{proof}
Let $\mathcal S_N$ be the set of $x$ satisfying~\ref{step:kr-ret2}; it is a Borel subset of
$\{q_x\ge q-2w\}$. For $x\in\mathcal S_N$ let $r=q_x$, $c_r=\xi'(r)$, $a_r=\xi''(r)$,
$k_x=N^{-1}b^{(r)}\cdot x$, $d_x=k_x+a_rr\chi_x$ and $\mathbf g(x)=b^{(r)}+a_r\chi_xx$, so that
$\nabla F_\zeta(x)=\nabla\HN(x)-\mathbf g(x)$. Covariance differentiation gives
$\Cov(\nabla\HN(x))=c_rI+\frac{a_r}Nxx^{\mathsf T}$, with determinant $c_r^N(1+a_rr/c_r)$, and
\[
 \log p_{\nabla\HN(x)}(\mathbf g(x))=-\frac N2\log(2\pi c_r)-\frac12\log\Bigl(1+\frac{a_rr}{c_r}\Bigr)
 -\frac{|b^{(r)}|^2}{2c_r}-N\Bigl[\frac{a_r\chi_xk_x}{c_r}+\frac{a_r^2r\chi_x^2}{2c_r}
 -\frac{a_rd_x^2}{2c_r(c_r+a_rr)}\Bigr],
\]
using $|\mathbf g(x)|^2=|b^{(r)}|^2+2Na_r\chi_xk_x+Na_r^2r\chi_x^2$ and $x\cdot\mathbf g(x)=Nd_x$.
Given $\nabla\HN(x)=\mathbf g(x)$, $\HN(x)$ is Gaussian with mean $Nc_rd_x/(c_r+a_rr)$ and variance
$N\sigma_H^2(r)$, where $\sigma_H^2(r)=\xi(r)-rc_r^2/(c_r+a_rr)>0$ is continuous near $q$
(Lemma~\ref{lem:kr-energy} at radius $r$). By Lemma~\ref{lem:kr-marked}, disintegrating over
$\HN(x)=Nh'$,
\[
 \E\widehat C_N=\int_{\mathcal S_N}p_{\nabla\HN(x)}(\mathbf g(x))\int_{|h'-h|\le N^{-2}}
 Np_x(Nh')\,\E_{x,h'}\bigl[|\det\nabla^2F_\zeta(x)|\,\ind\{\text{\ref{step:kr-ret4},\ \ref{step:kr-ret5}}\}\bigr]\,dh'\,dx ,
\]
with $p_x$ the conditional density of $\HN(x)$. Fix $t>0$. By Proposition~\ref{prop:kr-local}
and Lemma~\ref{lem:kr-marks}\ref{lem:kr-marks-b}--\ref{lem:kr-marks-c}, there is $w(t)>0$ such
that for $w\le w(t)$ and $N\ge N(t,w)$, uniformly over $x\in\mathcal S_N$ and $|h'-h|\le N^{-2}$,
the event that~\ref{step:kr-ret4} and~\ref{step:kr-ret5} hold and
$N^{-1}\log|\det\nabla^2F_\zeta(x)|\ge\mathfrak d_\pi-t$ has conditional probability at least
$\frac12$, so that
\[
 \E_{x,h'}\bigl[|\det\nabla^2F_\zeta(x)|\,\ind\{\text{\ref{step:kr-ret4},\ \ref{step:kr-ret5}}\}\bigr]
 \ge\tfrac12e^{N(\mathfrak d_\pi-t)} .
\]
Up to a constant factor, the energy window contributes at least
\[
 N^{-3/2}\exp\Bigl\{-\frac N{2\sigma_H^2(r)}\Bigl[\Bigl(h-\frac{c_rd_x}{c_r+a_rr}\Bigr)^2+O(N^{-2})\Bigr]\Bigr\}.
\]
Change
variables $x_i=M(b_i)$, a diffeomorphism of $\R^N$ onto $\ocube$ with Jacobian $\prod_iU(b_i)$, and
write $dx=\prod_i[U(b_i)/\pi'(b_i)]\,\pi^{\otimes N}(db)$ with $\pi'=d\pi/db$. On $\mathcal S_N$,
Lemma~\ref{lem:kr-endpoints} gives $|b^{(r)}_i-b_i|\le C|r-q|\le2Cw$ and $|\chi_x-\chi_N|+|k_x-k_N|\le Cw$,
hence $N^{-1}\bigl||b^{(r)}|^2-|b|^2\bigr|\le Cw$; the windows of $\mathcal W_N(w)$ control the
averages of $b^2,\log U,\log\pi',U$ and $bM$. Collecting the exponents,
$N^{-1}\sum_i\log U(b_i)$ cancels $-\E_\pi\log U$ in $\mathfrak d_\pi$ up to $w$; the terms
$-\E_\pi\log\pi'-\frac12\log(2\pi c)-\E_\pi b^2/(2c)$ make $-D(\pi\Vert\mathcal N(0,c))$; and the
quadratic identity
\[
 -\frac{a\chi k}c-\frac{a^2q\chi^2}{2c}+\frac{a(k+aq\chi)^2}{2c(c+aq)}+\frac a2\chi^2
 =\frac{a(k-c\chi)^2}{2c(c+aq)}
\]
turns the remaining terms of the gradient density and of $\mathfrak d_\pi$ into the second term of
\eqref{eq:kr-SigmaKR}. Hence, for $w\le w(t)$ and $N\ge N(t,w)$,
\[
 \E\widehat C_N\ge\tfrac12\,N^{-C}\,\pi^{\otimes N}(\mathcal S_N)\,
 \exp\bigl\{N\bigl(\Sigma_{\rm KR}(\pi,h)-t-\epsilon(w)\bigr)\bigr\},
\]
with $\epsilon(w)\to0$ as $w\downarrow0$. Under $\pi^{\otimes N}$ the empirical law of the fields
converges to $\pi$ in $W_2$ almost surely, $N^{-1}\sum_ie^{t_0|b_i|}\to\E_\pi e^{t_0|b|}<L_0$, and the
five averages converge, all by the law of large numbers; so $\pi^{\otimes N}(\mathcal S_N)\to1$. Let
$N\to\infty$, then $w\downarrow0$ and $t\downarrow0$.
\end{proof}

\paragraph{Separation and the spherical code.}
For $s\in(0,1)$ let $n_N(s)$ be the largest cardinality of a set of unit vectors in $\R^N$ with
pairwise inner products at most $s$.

\begin{lemma}[Spherical code]\label{lem:kr-code}
For $N\ge3$ and $\sqrt s<t<t'<1$,
\[
 n_N(s)\le\frac{2N}{t'-t}\bigl(1-t'^2\bigr)^{-(N-3)/2};\qquad\text{hence}\qquad
 \limsup_{N\to\infty}\frac1N\log n_N(s)\le-\frac12\log(1-s).
\]
\end{lemma}

\begin{proof}
Let $\mathcal A\subset S^{N-1}$ have pairwise inner products at most $s$, and let $K(u)=\{v\in S^{N-1}:\langle v,u\rangle\ge t\}$. We claim that
$|\mathcal A\cap K(u)|\le N$ for every unit $u$. If $u\in\mathcal A$, every other $v\in\mathcal A$ has
$\langle v,u\rangle\le s<t$, so the intersection is $\{u\}$. Otherwise each $v$ in it differs from
$\pm u$, so $v-\langle v,u\rangle u$ is a nonzero vector of $u^\perp$, and for distinct $v,v'$ in it
$\langle v-\langle v,u\rangle u,\,v'-\langle v',u\rangle u\rangle\le s-t^2<0$. In $\R^d$ at most $d+1$
vectors have pairwise negative inner products: if $v_1,\dots,v_{d+2}$ did, an affine dependence
$\sum_j\lambda_jv_j=0$, $\sum_j\lambda_j=0$, $\lambda\neq0$, would give
$z=\sum_{\lambda_j>0}\lambda_jv_j=\sum_{\lambda_k<0}|\lambda_k|v_k$ with
$|z|^2=\sum_{\lambda_j>0>\lambda_k}\lambda_j|\lambda_k|\langle v_j,v_k\rangle<0$. With $d=N-1$ this
proves the claim. For $u$ uniform on $S^{N-1}$, rotation invariance gives
$\E|\mathcal A\cap K(u)|=|\mathcal A|\Prob(u_1\ge t)$, so $|\mathcal A|\le N/\Prob(u_1\ge t)$. The
density of $u_1$ is proportional to $(1-z^2)^{(N-3)/2}\le1$ on $(-1,1)$, with normalizing integral
at most $2$, so $\Prob(u_1\ge t)\ge\frac12(t'-t)(1-t'^2)^{(N-3)/2}$. Let $N\to\infty$, $t'\downarrow t$
and $t\downarrow\sqrt s$.
\end{proof}

This is the bound of \cite[Lemma~3.11]{Bou26b}, where the counting step uses Rankin's theorem; we
include the short argument above.

\begin{lemma}[Separation of retained points]\label{lem:kr-separation}
Let $w\le w_1$ and $N\ge3\gamma^{-1/2}$. Then distinct retained points $x,y$ satisfy
$\ov xy<\delta'$, and $\widehat C_N(\eta,\delta,w)\le n_N(s_w)$ with $s_w=\delta'/(q-2w)\le\frac13$.
This holds for every realization of $\HN$ and every choice of the sets $\mathcal W_N(w)$.
\end{lemma}

\begin{proof}
Let $x\neq y$ be retained. By~\ref{step:kr-ret2} and~\ref{step:kr-ret3} and
Lemma~\ref{lem:kr-correction}\ref{lem:kr-correction-a},
$|F_\zeta(z)/N-f_\eta|\le N^{-2}+C_V\varrho(w)$ for $z\in\{x,y\}$. If $\ov xy\ge q-\delta'$, then
$|x-y|^2/N=q_x+q_y-2\ov xy\le4w+2\delta'\le3\delta'<\rho_\eta^2$, so $y\in B(x)$, and
Lemma~\ref{lem:kr-marks}\ref{lem:kr-marks-a} with~\ref{step:kr-ret1} at both points gives
$0=\langle\nabla F_\zeta(y)-\nabla F_\zeta(x),y-x\rangle\le-g_\eta|y-x|^2<0$, a contradiction. If
$\delta'\le\ov xy<q-\delta'$, then $y\in\mathcal K(x)$, so~\ref{step:kr-ret5} gives
$F_\zeta(y)-F_\zeta(x)\le-\gamma N/2$; but
$F_\zeta(y)-F_\zeta(x)\ge-2N(N^{-2}+C_V\varrho(w))>-\gamma N/2$, since $C_V\varrho(w)\le\gamma/8$ and
$N^{-2}<\gamma/8$, again a contradiction. So $\ov xy<\delta'$. The unit vectors $x/|x|$, for
retained $x$, are then pairwise distinct with inner products $\ov xy/\sqrt{q_xq_y}<\delta'/(q-2w)$,
since $q_x,q_y\ge q-2w$; and $s_w\le(q/4)/(3q/4)$.
\end{proof}

\begin{proof}[Proof of Proposition~\ref{prop:kr-main}]
Let $\eta$ be admissible and choose the parameters~\ref{step:kr-par1}--\ref{step:kr-par6}, with
$\mathcal W_N(w)$ as in Lemma~\ref{lem:kr-count}, and call the points counted by $\widehat C_N$
retained. We verify~\ref{step:stab-k1}--\ref{step:stab-k3} of Section~\ref{sec:stab-kr}, in which
$\pi_x$, $M$ and $\widehat C_N(\eta,\delta,w)$ are as here.

\emph{\ref{step:stab-k1}.} By~\ref{step:kr-ret1} and~\ref{step:kr-ret4}, every retained point is a
critical point of $F_\zeta$ in $\ocube$, and $\Hpm F_\zeta(y)\preceq-g_\eta I$ whenever
$|y-x|\le\rho_\eta\sqrt N$; $g_\eta$ and $\rho_\eta$ depend only on $\eta$ (through $L_0$, which we fix
as a function of $\eta$).

\emph{\ref{step:stab-k2}.} By Lemma~\ref{lem:kr-correction},
$\Wone(\mu_x,M_\#\pi_\eta)\le W_2(\pi_x,\pi_\eta)\le w$ and $|q_x-q|\le2w$. As in the proof of
Lemma~\ref{lem:kr-separation}, $|F_\zeta(x)/N-f_\eta|\le N^{-2}+C_V\varrho(w)$, and by
Lemma~\ref{lem:kr-correction}\ref{lem:kr-correction-b} the value $\FTAP(x)/N$ lies within
$C_V(\eta+\varrho(w))$ below $F_\zeta(x)/N$. With \eqref{eq:kr-feta},
\[
 \Bigl|\frac{F_\zeta(x)}N-f\Bigr|+\Bigl|\frac{\FTAP(x)}N-f\Bigr|
 \le\varepsilon_1(\eta)+\varepsilon_2(w)+2N^{-2},\qquad
 \varepsilon_1(\eta)=(2L+3C_V)\eta,\quad \varepsilon_2(w)=3C_V\varrho(w),
\]
and $\varepsilon_1(0)=0$, as required.

\emph{\ref{step:stab-k3}.} Since $\widehat C_N$ is a nonnegative integer,
$\E\widehat C_N\le\Prob(\widehat C_N\ge1)\max\widehat C_N\le\Prob(\widehat C_N\ge1)\,n_N(s_w)$ for
$N\ge3\gamma^{-1/2}$, by Lemma~\ref{lem:kr-separation}. Lemmas~\ref{lem:kr-count}
and~\ref{lem:kr-code}, the continuity of $w\mapsto s_w$ and $s_w\to\delta'/q$ give
\[
 \liminf_{w\downarrow0}\liminf_{N\to\infty}\frac1N\log\Prob\bigl(\widehat C_N(\eta,\delta,w)\ge1\bigr)
 \ge\Sigma_{\rm KR}(\pi_\eta,h_\eta)+\frac12\log\Bigl(1-\frac{\delta'}q\Bigr)
 \ge\Sigma(f)-\varepsilon_3(\eta)+\frac12\log\Bigl(1-\frac\delta q\Bigr),
\]
with $\varepsilon_3(\eta)=|\Sigma(f)-\Sigma_{\rm KR}(\pi_\eta,h_\eta)|\le C_\Sigma\eta$ by
Lemma~\ref{lem:kr-normalization}, and $\delta'\le\delta$.

\emph{Existence.} Fix $\varepsilon>0$. Choose an admissible $\eta$ with
$\varepsilon_1(\eta)\le\varepsilon/3$ (this is $\eta=0$ if $\Plef_0>0$), then $\delta\in(0,q/2)$, then
$w$ so small that $\varepsilon_2(w)\le\varepsilon/3$, and let $N^{-2}\le\varepsilon/6$. By~\ref{step:stab-k2}
every retained point $x$ satisfies \eqref{eq:event-exist} with $\zeta=\zf$, so
$\Prob(\Eex{\varepsilon}{f})\ge\Prob(\widehat C_N\ge1)$ for all large $N$, for every such $w$. By
\ref{step:stab-k3},
\[
 \liminf_{N\to\infty}\frac1N\log\Prob\bigl(\Eex{\varepsilon}{f}\bigr)
 \ge\Sigma(f)-C_\Sigma\eta+\frac12\log\Bigl(1-\frac\delta q\Bigr).
\]
The left side depends neither on $\eta$ nor on $\delta$. Let $\delta\downarrow0$ and then
$\eta\downarrow0$ through admissible values.
\end{proof}

If $\Plef_0>0$, the tilt is $\eta=0$ throughout: the Schur complement \eqref{eq:kr-schur} is strictly
negative, the obstacle is strictly above its contact value on $(0,q)$, and the argument is
that of \cite[Section~3]{Bou26b}, with her isolation conditions replaced by the marks~\ref{step:kr-ret4}--\ref{step:kr-ret5}; no limit in $\eta$ is needed.

\subsection{Margins and the annealed count of SUSY states}\label{app:kr-susy}

This subsection proves two statements of Section~\ref{sec:stab-kr} about the construction above:
Remark~\ref{rem:stab-margin}, which bounds its margins by the edge of the tilted bulk, and
Proposition~\ref{prop:stab-susy}, which compares the existence probability with the expected
number of SUSY states (Definition~\ref{def:stab-susy}). The upper bound of the latter
(Proposition~\ref{prop:kr-susy}) needs no spectral gap.

\paragraph{Margins.}
The tilted bulk law $\nu^{(\eta)}$ and its edge $\lamp^{(\eta)}$ are defined before
Remark~\ref{rem:stab-margin}.

\begin{proof}[Proof of Remark~\ref{rem:stab-margin}]
\emph{The bound $g_\eta\le-\lamp^{(\eta)}$.} Fix $\kappa>0$ and put
$\tau':=\frac12\min\{\kappa,\nu^{(\eta)}([\lamp^{(\eta)}-\kappa,\infty))\}$, which is positive because
$\lamp^{(\eta)}\in\supp\nu^{(\eta)}$. Apply Proposition~\ref{prop:hess-bulk} to the single pair
$(q,M_\#\pi_\eta)$, which puts no mass at $\pm1$, with the order parameter $\zf$ and accuracy $\tau'$.
Since $u(q,M(b))=\Phi_{\zf,xx}(q,b)=U(b)$, the bulk law attached to this pair is $\nu^{(\eta)}$. This gives
$w_H>0$ such that, for every $L'>0$ and all large $N$, outside an event of outer probability at most
$e^{-L'N}$, every $x\in\ocube$ with $|q_x-q|\le w_H$ and $\Wone(\mu_x,M_\#\pi_\eta)\le w_H$ satisfies
$\dL(\ESD(\Hpm F_{\zf}(x)),\nu^{(\eta)})\le\tau'$. Fix $\delta\in(0,q/2)$. By \ref{step:stab-k3} we
may choose $w\le w_H/2$ and $L'>0$ such that $\Prob(\widehat C_N(\eta,\delta,w)\ge1)\ge e^{-L'N/2}$ for
all large $N$. Then, for large $N$, with positive probability there is a retained point $x$ outside
the exceptional event; by \ref{step:stab-k2} it satisfies the two window conditions. For
$H=\Hpm F_{\zf}(x)$, since $\tau'<\kappa$,
\[
 \ESD(H)\bigl((\lamp^{(\eta)}-2\kappa,\infty)\bigr)
 \ge\nu^{(\eta)}\bigl((\lamp^{(\eta)}-2\kappa+\tau',\infty)\bigr)-\tau'
 \ge\nu^{(\eta)}\bigl([\lamp^{(\eta)}-\kappa,\infty)\bigr)-\tau'>0 .
\]
So $\lambda_{\max}(H)>\lamp^{(\eta)}-2\kappa$, while \ref{step:stab-k1} gives
$\lambda_{\max}(H)\le-g_\eta$. Hence $g_\eta<-\lamp^{(\eta)}+2\kappa$ for every $\kappa>0$.

\emph{The expansion.} Lemma~\ref{lem:hess-edge} applies to the law of $U$ under $\pi_\eta$. Its
hypothesis holds under $\pi_0$ (Lemma~\ref{lem:matrix-hyp}), and it is preserved by a change of
measure whose density $1-\eta\mathcal R$ lies in $[\frac12,\frac32]$, since such a change preserves
null sets and essential suprema and changes expectations of nonnegative functions by a factor at
least $\frac12$. The quadratic expansion of the edge in that lemma states that
$|\lamp+\Plef^2/(4a\E U^3)|\le C|\Plef|^3$ when $|\Plef|\le c_0$, with $c_0,C$ depending only on $a$
and $\E U^3$. Under $\pi_\eta$ one has $\E_{\pi_\eta}U^3\in[\frac12\E_{\pi_0}U^3,\frac32\E_{\pi_0}U^3]$,
so $c_0$ and $C$ can be taken uniform in $\eta\in(0,\eta_0]$ (they are explicit and monotone in
$\E U^3$; see the proof in Appendix~\ref{app:matrix-pastur}). Insert $\Plef_\eta$ from
\eqref{eq:stab-tilt} and $\E_{\pi_\eta}U^3=\E_{\pi_0}U^3-\eta\E_{\pi_0}[\mathcal RU^3]$.
\end{proof}

\paragraph{The annealed count.}
\begin{proof}[Proof of Proposition~\ref{prop:stab-susy}]
\emph{Lower bounds.} Take the construction of Theorem~\ref{thm:stab}\ref{thm:stab-a} at an
admissible tilt $\eta$.
The coupling that moves the mass of $\eta\mathcal R_+\pi_0$ onto that of $\eta\mathcal R_-\pi_0$
gives $W_2(\pi_\eta,\pi_0)^2\le2\eta\|\mathcal R\|_\infty\E_{\pi_0}b^2$. By \ref{step:stab-k1} and
\ref{step:stab-k2}, once $\eta$, $w$ and $N^{-1}$ are small enough in terms of $\varepsilon$, every
retained point is an $\varepsilon$-SUSY state. Hence
$\E\mathsf N_{N,\varepsilon}\ge\Prob(\mathsf N_{N,\varepsilon}\ge1)\ge\Prob(\widehat C_N\ge1)$, and
\ref{step:stab-k3} gives both lower bounds after $\delta\downarrow0$ and $\eta\downarrow0$.

\emph{Upper bound.} The upper bound of part~\ref{prop:stab-susy-i} is
Proposition~\ref{prop:kr-susy} below, whose hypotheses are those of
Proposition~\ref{prop:stab-susy}. The $\liminf_N$ and the $\limsup_N$ of
$N^{-1}\log\E\mathsf N_{N,\varepsilon}$ are nondecreasing in $\varepsilon$, since
$\mathsf N_{N,\varepsilon}$ is; with the lower bound this gives the last assertion of
part~\ref{prop:stab-susy-i}.

The first assertion of part~\ref{prop:stab-susy-ii} was obtained with the lower bounds. The second
follows from part~\ref{prop:stab-susy-i}, the first assertion, and
$1\le\E[\mathsf N_{N,\varepsilon}\mid\mathsf N_{N,\varepsilon}\ge1]
=\E\mathsf N_{N,\varepsilon}/\Prob(\mathsf N_{N,\varepsilon}\ge1)$.
\end{proof}

\paragraph{The upper bound.}
\begin{proposition}[Annealed upper bound for SUSY states]\label{prop:kr-susy}
Assume the hypotheses of Proposition~\ref{prop:kr-main}, and let $\mathsf N_{N,\varepsilon}$ be the
number of $\varepsilon$-SUSY states at level $f$. Then
\[
 \lim_{\varepsilon\downarrow0}\limsup_{N\to\infty}\frac1N\log\E\mathsf N_{N,\varepsilon}\le\Sigma(f).
\]
\end{proposition}

The proof uses neither the tilt, nor the energy shift, nor the cap \eqref{eq:kr-cap}, and its
only stability input is the non-strict Plefka inequality $\Plef_0\ge0$, which holds at every
contact. The Kac--Rice formula bounds $\E\mathsf N_{N,\varepsilon}$ by an integral over the window,
in which the condition on $\FTAP$ is dropped. The change of variables $x_i=M(b_i)$ has Jacobian
$\prod_iU(b_i)$ and turns this integral into an expectation under $\pi_0^{\otimes N}$. The Gaussian
density factors are explicit; the only factor that is not is the conditional expectation of
$|\det\nabla^2F_\zeta(x)|$ given the gradient and the value of $\HN$ at $x$. At a marginal contact
the limiting spectral measure of the transverse Hessian touches zero, while the empirical field
laws in the window may violate the Plefka inequality, so the estimate below inverts no matrix and
uses no spectral gap. It replaces $\log|\lambda|$ by the $\iota^{-1}$-Lipschitz function
$\log\max(|\lambda|,\iota)$, for which Gaussian concentration and a law of large numbers for
deformed GOE matrices hold uniformly over the window, and lets $\iota\downarrow0$ only in the
limiting logarithmic potential. That potential is finite and explicit because $\Plef_0\ge0$: the
law $\nu_{\qf}$ lives on $(-\infty,0]$ (Theorem~\ref{thm:struct}\ref{thm:struct-iv}), and
\[
 \int\log|\lambda|\,\nu_{\qf}(d\lambda)=-\E\log U_{\qf}+\frac{\xi''(\qf)}2\chi_{\qf}^2 ,
\]
also when the edge touches zero. This identity is \eqref{eq:kr-edge} at $\rho=\rho_{\qf}$ and
$w_0=a\chi_0$. It is an instance of the Hopf--Lax formula for the logarithmic potential of a free
convolution with a semicircle law \cite[Lemma~C.3(2), p.~67]{Bou26a}, stated there, for measures
with support bounded below, at every point at or below the left edge. Lemma~\ref{lem:kr-edge}
covers the case used here, under a logarithmic moment, together with the fact that the measure
puts no mass on $(-\infty,0]$; its only stability input is $a\int\rho(dy)(y-w_0)^{-2}\le1$, which
at $\rho_{\qf}$ is the non-strict Plefka inequality. With it, the Jacobian cancels the term
$-\E\log U_{\qf}$, and the remaining exponent is $\Sigma(f)$ by the identities that normalize the
count (Lemma~\ref{lem:kr-normalization}). Two active degrees enter the upper bound as well, through
the positive conditional variance of the energy.

\paragraph{Truncated log-determinants.}
For $\iota\in(0,1]$ and a real symmetric $n\times n$ matrix $A$ put
$l_\iota(\lambda)=\log\max(|\lambda|,\iota)$ and $\operatorname{ld}_\iota(A)=\sum_il_\iota(\lambda_i(A))$.
Then $|\det A|\le e^{\operatorname{ld}_\iota(A)}$, and $l_\iota$ is even, nondecreasing in $|\lambda|$,
at least $\log\iota$, and $\iota^{-1}$-Lipschitz.

\begin{lemma}[Spectral facts]\label{lem:kr-ld}
Let $A,B$ be real symmetric $n\times n$ matrices and $\iota\in(0,1]$.
\begin{thmlist}
\item\label{lem:kr-ld-a} $|\operatorname{ld}_\iota(A)-\operatorname{ld}_\iota(B)|\le\iota^{-1}\sqrt n\,\|A-B\|_{\rm HS}$.
\item\label{lem:kr-ld-b} If $n\ge2$ and $\check A$ is the compression of $A$ to a hyperplane, then
 $\operatorname{ld}_\iota(\check A)\le\operatorname{ld}_\iota(A)+\log\frac1\iota$ and
 $\operatorname{ld}_\iota(A)\le\operatorname{ld}_\iota(\check A)+\log\frac1\iota+2\log\max(\|A\|,1)$.
\item\label{lem:kr-ld-c} If $B=\mathsf D-E$ with $\mathsf D\succeq0$ diagonal and $E$ symmetric, then
 $\sum_i\varphi(|\lambda_i(B)|)\le\sum_i\varphi(\mathsf D_{ii}+\|E\|)$ for every nondecreasing
 $\varphi:[0,\infty)\to[0,\infty)$.
\end{thmlist}
\end{lemma}

\begin{proof}
\ref{lem:kr-ld-a} By the Hoffman--Wielandt inequality \cite{HW},
$\sum_i(\lambda_i(A)-\lambda_i(B))^2\le\|A-B\|_{\rm HS}^2$, with eigenvalues in decreasing order;
combine it with the Lipschitz bound for $l_\iota$ and Cauchy--Schwarz.

\ref{lem:kr-ld-b} Let $\mu_1\ge\dots\ge\mu_{n-1}$ be the eigenvalues of $\check A$. By
Courant--Fischer, $\lambda_{j+1}\le\mu_j\le\lambda_j$. If $\mu_j\ge0$ then
$l_\iota(\mu_j)\le l_\iota(\lambda_j)$; if $\mu_j<0$ then $l_\iota(\mu_j)\le l_\iota(\lambda_{j+1})$.
With $t=\#\{j:\mu_j\ge0\}$ the indices used are $\{1,\dots,t\}\cup\{t+2,\dots,n\}$; only
$\lambda_{t+1}$ is unused, and $l_\iota(\lambda_{t+1})\ge\log\iota$. This is the first inequality. For
the second, let $s=\#\{i:\lambda_i\ge0\}$. If $2\le i\le s$ then $0\le\lambda_i\le\mu_{i-1}$, and if
$s+1\le i\le n-1$ then $\mu_i\le\lambda_i<0$; so $l_\iota(\lambda_i)$ is at most $l_\iota$ of the
indicated $\mu$, and these indices are distinct. At most two eigenvalues of $A$ are unassigned, each
with $l_\iota\le\log\max(\|A\|,1)$, and at most one eigenvalue of $\check A$ is unused, with
$l_\iota\ge\log\iota$.

\ref{lem:kr-ld-c} By Courant--Fischer, $|\lambda_i(B)-\lambda_i(\mathsf D)|\le\|E\|$, and
$\lambda_i(\mathsf D)\ge0$, so $|\lambda_i(B)|\le\lambda_i(\mathsf D)+\|E\|$; apply $\varphi$ and sum.
\end{proof}

For a probability measure $\nu$ on $\R$ we write $\check\nu$ for its reflection, the law of
$-\lambda$ when $\lambda\sim\nu$. So $\check\nu[a,\rho]$, the reflection of the reflected Pastur law
of Lemma~\ref{lem:hess-pastur}, is the limit of the empirical spectral distributions of
$D_N-\sqrt a\,\GOE$. By Lemma~\ref{lem:hess-pastur}\ref{lem:hess-pastur-a}, after reflection and
complex conjugation, its Stieltjes transform $\mathsf s(z)=\int\check\nu[a,\rho](d\lambda)/(\lambda-z)$,
$z\in\Cplx^+$, is the unique solution with $\operatorname{Im}\mathsf s>0$ of
\begin{equation}\label{eq:kr-pastur}
 \mathsf s=\int\frac{\rho(dy)}{y-z-a\mathsf s} .
\end{equation}

\begin{lemma}[Logarithmic law of large numbers]\label{lem:kr-deformed}
Let $a>0$ and $\rho\in\Pspace([0,\infty))$ with $\int\log(1+y)\,\rho(dy)<\infty$. Let $N_k\to\infty$,
$a_k\to a$, and let $\mathsf D^{(k)}=\diag(d_1,\dots,d_{N_k})$ be deterministic, with $d_i\ge0$,
empirical laws converging weakly to $\rho$, and
\begin{equation}\label{eq:kr-tails}
 \lim_{K\to\infty}\limsup_{k\to\infty}\frac1{N_k}\sum_i\log_+\frac{d_i}K=0,\qquad
 \frac1{N_k}\log\bigl(1+\max_id_i\bigr)\to0 .
\end{equation}
Then $\int\log(1+|\lambda|)\,\check\nu[a,\rho](d\lambda)<\infty$, and for every continuous $\varphi$ with
$|\varphi(\lambda)|\le C_\varphi(1+\log(1+|\lambda|))$,
\[
 \frac1{N_k}\E\sum_{i=1}^{N_k}\varphi\bigl(\lambda_i(\mathsf D^{(k)}-\sqrt{a_k}\,W_{N_k})\bigr)
 \longrightarrow\int\varphi\,d\check\nu[a,\rho] .
\]
\end{lemma}

\begin{proof}
Write $\lambda_i=\lambda_i(\mathsf D^{(k)}-\sqrt{a_k}W)$. By
Lemma~\ref{lem:matrix-comparison}\ref{lem:matrix-comparison-a}, $\Prob(\|W\|>6)\le2e^{-4N}$. By
Lemma~\ref{lem:hess-pastur}\ref{lem:hess-pastur-a}, the expected empirical spectral distribution of
$\mathsf D^{(k)}-\sqrt aW$ converges weakly to $\check\nu[a,\rho]$; by Weyl's inequality its
eigenvalues move by at most $|\sqrt{a_k}-\sqrt a|\,\|W\|$ when $a$ is replaced by $a_k$, so the
expected empirical spectral distribution $\bar\mu_k$ of $\mathsf D^{(k)}-\sqrt{a_k}W$ has the same
limit.

For $R\ge1$ let $\mathfrak t_R(y)=\log(1+y)\ind\{y>R\}$. By Lemma~\ref{lem:kr-ld}\ref{lem:kr-ld-c},
$\frac1N\sum_i\mathfrak t_R(|\lambda_i|)\le\frac1N\sum_i\mathfrak t_R(d_i+\sqrt{a_k}\|W\|)$. Take
$R\ge\max(12\sqrt{\bar a},18)$, $\bar a=\sup_ka_k$, and work on $\{\|W\|\le6\}$. If
$d_i+6\sqrt{a_k}>R$, then $d_i>R/2$ and $d_i+6\sqrt{a_k}\le2d_i$, and with $K=\sqrt{R/2}\ge3$ and
$d_i>K^2$, $\log(1+2d_i)\le3\log(d_i/K)$. So the right side is at most $3N^{-1}\sum_i\log_+(d_i/K)$. On
$\{\|W\|>6\}$ it is at most $\log(1+\max_id_i)+\log(1+\sqrt{a_k}\|W\|)$, and
$\E\log^2(1+\sqrt{a_k}\|W\|)\le a_k\E\|W\|_{\rm HS}^2=a_k(N+1)$. By Cauchy--Schwarz and
\eqref{eq:kr-tails},
\[
 \lim_{R\to\infty}\limsup_{k\to\infty}\frac1{N_k}\E\sum_i\mathfrak t_R(|\lambda_i|)=0,\qquad
 \sup_k\frac1{N_k}\E\sum_i\log(1+|\lambda_i|)<\infty .
\]
Weak lower semicontinuity of integrals of nonnegative continuous functions gives the logarithmic
moment of $\check\nu[a,\rho]$. Finally, let $\mathsf k_R$ be continuous, $0\le\mathsf k_R\le1$, equal to
one on $[-R,R]$ and to zero off $[-R-1,R+1]$. Then
$\int\varphi\mathsf k_R\,d\bar\mu_k\to\int\varphi\mathsf k_R\,d\check\nu[a,\rho]$, while
$|\varphi|(1-\mathsf k_R)\le2C_\varphi\mathfrak t_R(|\cdot|)$ for $R\ge e$, which is small uniformly in large
$k$ for both measures.
\end{proof}

\begin{lemma}[Logarithmic potential at the edge]\label{lem:kr-edge}
Let $a>0$ and $\rho\in\Pspace([0,\infty))$ with $\int\log(1+y)\,\rho(dy)<\infty$, and let $w_0\in\R$
satisfy
\[
 \inf\supp\rho\ge w_0+1,\qquad a\int\frac{\rho(dy)}{y-w_0}=w_0,\qquad
 a\int\frac{\rho(dy)}{(y-w_0)^2}\le1 .
\]
Then $\check\nu=\check\nu[a,\rho]$ satisfies $\check\nu((-\infty,0])=0$ and
\begin{equation}\label{eq:kr-edge}
 \int\log\lambda\,\check\nu(d\lambda)=\int\log(y-w_0)\,\rho(dy)+\frac{w_0^2}{2a}\in\R .
\end{equation}
\end{lemma}

\begin{proof}
Put $y_\rho=\inf\supp\rho>w_0$, $\mathfrak g(w)=\int\rho(dy)/(y-w)$ and $\mathcal Z(w)=w-a\mathfrak g(w)$,
holomorphic off $[y_\rho,\infty)$. On $(-\infty,y_\rho)$, $\mathfrak g,\mathfrak g',\mathfrak g''>0$, so
$\mathcal Z'=1-a\mathfrak g'$ is strictly decreasing and less than one. By hypothesis $\mathcal Z(w_0)=0$
and $\mathcal Z'(w_0)\ge0$, so $0<\mathcal Z'<1$ on $(-\infty,w_0)$, and $\mathcal Z$ is an increasing
homeomorphism of $(-\infty,w_0]$ onto $(-\infty,0]$. Let $\Omega$ be its inverse, which is $C^1$ on
$(-\infty,0)$ with $\Omega'=1/\mathcal Z'(\Omega)>1$, and put $\mathsf m(z)=(\Omega(z)-z)/a=\mathfrak g(\Omega(z))$
for $z\le0$.

\emph{No mass below zero.} Fix $E<0$. Since $\mathcal Z'(\Omega(E))>0$, the holomorphic inverse
function theorem gives a disc $U_E$ about $E$ and a holomorphic $\Omega_E$ on $U_E$, with values in
a disc about $\Omega(E)$ that avoids $[y_\rho,\infty)$, such that $\mathcal Z\circ\Omega_E=\id$ and
$\Omega_E=\Omega$ on $U_E\cap\R$. Then $\widetilde{\mathsf m}=(\Omega_E-z)/a=\mathfrak g(\Omega_E)$ solves
\eqref{eq:kr-pastur} on $U_E$, is real on $U_E\cap\R$ with $\widetilde{\mathsf m}'=(\Omega'-1)/a>0$
there, and after shrinking $U_E$, $\operatorname{Re}\widetilde{\mathsf m}'\ge\kappa_E>0$ on $U_E$. By the
Cauchy--Riemann equations $\operatorname{Im}\widetilde{\mathsf m}(t+\mathrm iv)\ge\kappa_Ev>0$ for
$t+\mathrm iv\in U_E$, $v>0$, so by uniqueness $\mathsf s=\widetilde{\mathsf m}$ on $U_E\cap\Cplx^+$. Hence
$\operatorname{Im}\mathsf s(\lambda+\mathrm iv)\to0$ as $v\downarrow0$, uniformly on an interval $I_E$
about $E$. Since $\frac1\pi\operatorname{Im}\mathsf s(\lambda+\mathrm iv)$ is the density of $Y+v\mathsf C$,
with $Y\sim\check\nu$ and $\mathsf C$ an independent standard Cauchy variable, and $Y+v\mathsf C\to Y$ in law,
$\check\nu(I_E)\le\liminf_{v\downarrow0}\Prob(Y+v\mathsf C\in I_E)=0$. Countably many $I_E$ cover
$(-\infty,0)$, so $\check\nu((-\infty,0))=0$, and for real $z<0$, by dominated convergence,
$\int\check\nu(d\lambda)/(\lambda-z)=\mathsf m(z)$.

\emph{Two potentials.} For $z<0$ put $\mathcal L(z)=\int\log(\lambda-z)\,\check\nu(d\lambda)$ and
$\widetilde{\mathcal L}(z)=\int\log(y-\Omega(z))\,\rho(dy)+(\Omega(z)-z)^2/(2a)$. Both are finite
($y-\Omega(z)\ge1$ on $\supp\rho$, and $\check\nu$ has a logarithmic moment by
Lemma~\ref{lem:hess-pastur}\ref{lem:hess-pastur-b} applied to the reflection,
$\check\nu((2\sqrt a+s,\infty))\le\rho([s,\infty))$, and the hypothesis
$\int\log(1+y)\,\rho(dy)<\infty$), and
differentiation under the integrals gives $\mathcal L'=-\mathsf m$ and
$\widetilde{\mathcal L}'=-\Omega'\mathfrak g(\Omega)+\mathsf m(\Omega'-1)=-\mathsf m$. As $z\to-\infty$,
$\mathcal L(z)-\log(-z)\to0$ by dominated convergence; also $0<\Omega(z)-z=a\mathfrak g(\Omega(z))\le a$
and $\mathfrak g(\Omega(z))\to0$, so $(y-\Omega(z))/(-z)\to1$ with values in $[\frac12,1+y]$ for
$|z|\ge\max(2a,1)$, and $\widetilde{\mathcal L}(z)-\log(-z)\to0$. So $\mathcal L=\widetilde{\mathcal L}$ on
$(-\infty,0)$.

\emph{The edge.} As $z\uparrow0$, $\log(\lambda-z)\downarrow\log\lambda$ on $[0,\infty)$, and monotone
convergence gives $\mathcal L(z)\to\int\log\lambda\,d\check\nu\in[-\infty,\infty)$. Also
$\Omega(z)\uparrow w_0$, $\log(y-\Omega(z))\downarrow\log(y-w_0)\ge0$ on $\supp\rho$, and
$(\Omega(z)-z)^2\to w_0^2$. So $\widetilde{\mathcal L}(z)$ converges to the right side of \eqref{eq:kr-edge},
which is finite. In particular $\check\nu(\{0\})=0$.
\end{proof}

\paragraph{The window.}
Fix $\varepsilon_0\in(0,\min(q,1-q)/2)$ and $\mathcal J_0=[q-\varepsilon_0,q+\varepsilon_0]$. For
$x\in\ocube$ with $r=q_x\in\mathcal J_0$ let $a_r,c_r,k_x,d_x$ and $\mathbf g(x)$ be as in the proof
of Lemma~\ref{lem:kr-count}, with $\chi_x$ as before \eqref{eq:kr-gradient}, and put
\[
 \mathsf D_x=\diag\bigl(u(r,x_i)^{-1}+a_r\chi_x\bigr),
\]
so that $\nabla F_\zeta(x)=\nabla\HN(x)-\mathbf g(x)$ and $F_\zeta(x)/N=\HN(x)/N+V_\zeta(q_x,\mu_x)$. For
$0<\varepsilon\le\varepsilon_0$ let
$\mathcal W_\varepsilon=\{x\in\ocube:\ |q_x-q|\le\varepsilon,\ W_2(\pi_x,\pi_0)\le\varepsilon\}$ and let
$\mathcal W_\varepsilon^\circ$ be the set of $x\in\mathcal W_\varepsilon$ whose radius is not an atom of
$\zeta$; the difference is contained in countably many spheres. Write $b_i=M^{-1}(x_i)$ for the
fields of $x$. Constants below depend only on $\xi$, $\zeta$, $q$ and $\varepsilon_0$.

\begin{lemma}[Window facts]\label{lem:kr-susywindow}
For $0<\varepsilon\le\varepsilon_0$ and $x\in\mathcal W_\varepsilon$:
\begin{thmlist}
\item\label{lem:kr-susywindow-a} $N^{-1}|b|^2\le C$, $\max_i|b_i|\le C\sqrt N$ and
 $\max_iu(q_x,x_i)^{-1}\le Ce^{C\sqrt N}$;
\item\label{lem:kr-susywindow-b} $|N^{-1}\sum_i\log U(b_i)-\E_0\log U|\le C_{34}\varepsilon$ and
 $|N^{-1}\sum_i\Phi(q,b_i)-\E_0\Phi(q,b)|\le\varepsilon$;
\item\label{lem:kr-susywindow-c} $|\chi_x-\chi_0|+|k_x-k_0|\le C\varepsilon$ and
 $N^{-1}\bigl||b^{(q_x)}|^2-|b|^2\bigr|\le C\varepsilon$;
\item\label{lem:kr-susywindow-d} $|V_\zeta(q_x,\mu_x)-(f-h_0)|\le C_V\varrho(\varepsilon)$.
\end{thmlist}
\end{lemma}

\begin{proof}
Throughout, $\Wone(\pi_x,\pi_0)\le W_2(\pi_x,\pi_0)\le\varepsilon$, and $r=q_x$.
\ref{lem:kr-susywindow-a} Under an optimal coupling,
$(N^{-1}|b|^2)^{1/2}\le\varepsilon+(\E_0b^2)^{1/2}$, and $|b_i|\le|b|$. By
Lemma~\ref{lem:kr-contact}\ref{lem:kr-contact-d} and Lemma~\ref{lem:kr-endpoints},
$u(r,x_i)^{-1}\le c_1^{-1}e^{2|b_i^{(r)}|}\le c_1^{-1}e^{2|b_i|+C\varepsilon}$.
\ref{lem:kr-susywindow-b} $\log U$ is $C_{34}$-Lipschitz and $\Phi(q,\cdot)$ is $1$-Lipschitz; integrate
over an optimal coupling. \ref{lem:kr-susywindow-c} By Lemma~\ref{lem:kr-endpoints},
$|u(r,x_i)-U(b_i)|\le C|r-q|$ and $|b_i^{(r)}-b_i|\le C|r-q|$; $U$ is $C_{34}$-Lipschitz and $b\mapsto bM(b)$
is $C$-Lipschitz, since $|M+bU|\le1+C|b|e^{-2|b|}$; and
$\bigl||b^{(r)}|^2-|b|^2\bigr|\le\sum_i|b_i^{(r)}-b_i|(2|b_i|+C\varepsilon)\le C\varepsilon N$ by
\ref{lem:kr-susywindow-a} and Cauchy--Schwarz. \ref{lem:kr-susywindow-d} This is
Lemma~\ref{lem:kr-correction}\ref{lem:kr-correction-a} at $\eta=0$, together with
$h_0+V_\zeta(q,M_\#\pi_0)=f_0=f$ from \eqref{eq:kr-feta}.
\end{proof}

\begin{lemma}[Structure of the conditional Hessian]\label{lem:kr-susystructure}
Let $x\in\ocube$ with $r=q_x\in\mathcal J_0$ not an atom of $\zeta$, and $h\in\R$. Conditionally on
$\nabla\HN(x)=\mathbf g(x)$ and $\HN(x)=Nh$, the matrix $\nabla^2F_\zeta(x)$ has the law of
$\mathbf A=\mathbf A_0+\mathsf K\mathbf Y$, where $\mathbf Y$ is a standard Gaussian vector,
$\mathsf K$ is linear and $\mathbf A_0$ is a deterministic symmetric matrix, and:
\begin{thmlist}
\item\label{lem:kr-susystructure-a} there are an orthogonal matrix $O$ with first column $x/|x|$ and a
 GOE matrix $W'$ whose independent entries are coordinates of $\mathbf Y$ such that the compression
 of $\mathbf A$ to $x^\perp$ equals that of $\sqrt{a_r}\,W-\mathsf D_x$, where $W=OW'O^{\mathsf T}$ is a GOE
 matrix;
\item\label{lem:kr-susystructure-b} $\|\mathsf K\mathbf y\|_{\rm HS}^2\le C_{\mathsf K}N^{-1}|\mathbf y|^2$ for every
 $\mathbf y$;
\item\label{lem:kr-susystructure-c} $\|\mathbf A_0\|\le\max_iu(r,x_i)^{-1}+C(1+|h|+|b^{(r)}|/\sqrt N)$.
\end{thmlist}
\end{lemma}

\begin{proof}
Let $e=x/|x|$ and complete it to an orthonormal basis $e,f_2,\dots,f_N$, the columns of $O$. The
covariances in the proof of Lemma~\ref{lem:kr-block}, at radius $r$, show that the transverse
family $(f_j^{\mathsf T}\nabla^2\HN(x)f_k)$ has GOE covariance with off-diagonal variance $a_r/N$ and
is uncorrelated with $\HN(x)$, $\nabla\HN(x)$, the cross entries and the radial entry; so in these
coordinates it is $\sqrt{a_r}$ times the lower right block of a GOE matrix $W'$, which we complete
by independent entries that $\mathsf K$ ignores. The radial row is conditionally Gaussian with
covariance at most the unconditional one, of operator norm at most $C/N$; this gives
\ref{lem:kr-susystructure-a} for the Hamiltonian part and \ref{lem:kr-susystructure-b}. Its
conditional mean is $\frac{a_r\sqrt r}{c_r\sqrt N}Pb^{(r)}$ for the cross part and a linear
combination of $h$ and $d_x$ for the radial entry, with coefficients bounded on $\mathcal J_0$ because
$\det S_0(r)>0$ there (Lemma~\ref{lem:kr-energy}); and $|d_x|\le|b^{(r)}|/\sqrt N+C$. The
deterministic part of $\nabla^2F_\zeta(x)-\nabla^2\HN(x)$ is read off from \eqref{eq:kr-hessian}: its
compression to $x^\perp$ is that of $-\mathsf D_x$, and its norm is at most $\max_iu_i^{-1}+C$, since
$|\mathbf v_x|\le C_{34}\sqrt N$ and the scalar coefficient of $xx^{\mathsf T}$ is bounded.
\end{proof}

Recall that $\rho_{\qf}$ is the law of $U^{-1}+a\chi_0$ under $\pi_0$. Then
$\supp\rho_{\qf}\subset[1+a\chi_0,\infty)$ and $\int\log(1+y)\,\rho_{\qf}(dy)\le C+2\E_0|b|<\infty$, and
$\Pastur{a}{\rho_{\qf}}=\nu_{\qf}$.

\begin{lemma}[Diagonals in shrinking windows]\label{lem:kr-susydiag}
Let $\varepsilon_k\downarrow0$, $N_k\to\infty$ and $x^{(k)}\in\mathcal W_{\varepsilon_k}\subset(-1,1)^{N_k}$, with
radii $r_k$. Then $a_{r_k}\to a$, and the matrices $\mathsf D_{x^{(k)}}$ satisfy the hypotheses of
Lemma~\ref{lem:kr-deformed} with $\rho=\rho_{\qf}$.
\end{lemma}

\begin{proof}
We have $r_k\to q$, the diagonal entries are at least one, and $\chi_{x^{(k)}}\to\chi_0$ by
Lemma~\ref{lem:kr-susywindow}\ref{lem:kr-susywindow-c}. Since
$|\log(A+B)-\log(A'+B')|\le|\log A-\log A'|+|B-B'|$ for $A,A'\ge1$, $B,B'\ge0$, and $\log u$ is
Lipschitz in the radius (Lemma~\ref{lem:kr-endpoints}),
$\max_i|\log d_i^{(k)}-\log(U(b_i)^{-1}+a\chi_0)|\le C|r_k-q|+|a_{r_k}\chi_{x^{(k)}}-a\chi_0|\to0$. The
empirical field laws converge weakly to $\pi_0$, so the empirical laws of the $d_i^{(k)}$ converge to
$\rho_{\qf}$. Next, $\log d_i^{(k)}\le C_2+2|b_i|$, so $\log_+(d_i^{(k)}/K)\le(C_2+2|b_i|-\log K)_+$, a
$2$-Lipschitz function of $b_i$ whose $\pi_0$-mean tends to zero as $K\to\infty$; and
$\log(1+\max_id_i^{(k)})=O(\sqrt{N_k})$ by Lemma~\ref{lem:kr-susywindow}\ref{lem:kr-susywindow-a}.
\end{proof}

\begin{lemma}[An annealed determinant bound without a spectral gap]\label{lem:kr-susydet}
For $H_0<\infty$ and $0<\varepsilon\le\varepsilon_0$ let $\mathfrak D_N(\varepsilon,H_0)$ be the supremum of
$N^{-1}\log\E\bigl[|\det\nabla^2F_\zeta(x)|\bigm|\nabla\HN(x)=\mathbf g(x),\ \HN(x)=Nh\bigr]$ over
$x\in\mathcal W_\varepsilon^\circ$ and $|h|\le H_0$. Then
\[
 \lim_{\varepsilon\downarrow0}\limsup_{N\to\infty}\mathfrak D_N(\varepsilon,H_0)\le-\E_0\log U+\frac a2\chi_0^2 .
\]
\end{lemma}

\begin{proof}
Fix $\iota\in(0,1]$ and let $\mathbf A$ be as in Lemma~\ref{lem:kr-susystructure}. By
Lemmas~\ref{lem:kr-ld}\ref{lem:kr-ld-a} and~\ref{lem:kr-susystructure}\ref{lem:kr-susystructure-b},
$\Psi(\mathbf y)=\operatorname{ld}_\iota(\mathbf A_0+\mathsf K\mathbf y)$ is $L_\iota$-Lipschitz with
$L_\iota=\sqrt{C_{\mathsf K}}/\iota$. For a smooth $L_\iota$-Lipschitz function, Gaussian concentration
\cite[pp.~135--136]{Ledoux} gives $\E e^{\Psi(\mathbf Y)}\le e^{\E\Psi(\mathbf Y)+L_\iota^2/2}$; a
mollification $\Psi*\phi_\tau$, which is $L_\iota$-Lipschitz and within $L_\iota\tau$ of $\Psi$, extends
this to $\Psi$. Hence $\E|\det\mathbf A|\le\exp\{\E\operatorname{ld}_\iota(\mathbf A)+C_{\mathsf K}/(2\iota^2)\}$. By
Lemma~\ref{lem:kr-susystructure}\ref{lem:kr-susystructure-a}, the compression of $\mathbf A$ to
$x^\perp$ is that of $-\mathbf B_x$, $\mathbf B_x=\mathsf D_x-\sqrt{a_r}W$, so by
Lemma~\ref{lem:kr-ld}\ref{lem:kr-ld-b}, applied twice, and the evenness of $l_\iota$,
\[
 \operatorname{ld}_\iota(\mathbf A)\le\operatorname{ld}_\iota(\mathbf B_x)+2\log\tfrac1\iota
 +2\log\max(\|\mathbf A\|,1).
\]
By Jensen's inequality and Lemmas~\ref{lem:kr-susystructure}\ref{lem:kr-susystructure-b}--\ref{lem:kr-susystructure-c}
and~\ref{lem:kr-susywindow}\ref{lem:kr-susywindow-a}, $\E\log\max(\|\mathbf A\|,1)\le C(1+H_0+\sqrt N)$.
Therefore
\[
 \mathfrak D_N(\varepsilon,H_0)\le\sup_{x\in\mathcal W_\varepsilon}\frac1N\E\operatorname{ld}_\iota\bigl(\mathsf D_x-\sqrt{a_{q_x}}W\bigr)
 +\frac{C_{\mathsf K}/(2\iota^2)+2\log(1/\iota)+C(1+H_0+\sqrt N)}N .
\]
We claim that
$\lim_{\varepsilon\downarrow0}\limsup_N\sup_{x\in\mathcal W_\varepsilon}N^{-1}\E\operatorname{ld}_\iota(\mathsf D_x-\sqrt{a_{q_x}}W)\le\int l_\iota\,d\check\nu_{\qf}$,
where $\check\nu_{\qf}=\check\nu[a,\rho_{\qf}]$. Otherwise there are $\epsilon_*>0$, $\varepsilon_k\downarrow0$,
$N_k\uparrow\infty$ and $x^{(k)}\in\mathcal W_{\varepsilon_k}$ with the normalized expectation at least
$\int l_\iota\,d\check\nu_{\qf}+\epsilon_*$; but by Lemmas~\ref{lem:kr-susydiag} and~\ref{lem:kr-deformed},
with $\varphi=l_\iota$, which is continuous with $|l_\iota|\le\log(1/\iota)+\log(1+|\lambda|)$, the normalized
expectations converge to $\int l_\iota\,d\check\nu_{\qf}$. So the left side of the lemma is at most
$\int l_\iota\,d\check\nu_{\qf}$ for every $\iota\in(0,1]$. Apply Lemma~\ref{lem:kr-edge} with $\rho=\rho_{\qf}$ and
$w_0=a\chi_0$: $\inf\supp\rho_{\qf}\ge w_0+1$, $a\int\rho_{\qf}(dy)/(y-w_0)=a\E_0U=w_0$, and
$a\int\rho_{\qf}(dy)/(y-w_0)^2=a\E_0U^2=1-\Plef_0\le1$. So $\check\nu_{\qf}$ lives on $(0,\infty)$ and
$\int\log\lambda\,d\check\nu_{\qf}=-\E_0\log U+\frac a2\chi_0^2$. As $\iota\downarrow0$,
$l_\iota\downarrow\log|\lambda|$ and $l_1$ is integrable, so monotone convergence finishes the proof.
\end{proof}

\begin{proof}[Proof of Proposition~\ref{prop:kr-susy}]
Let $0<\varepsilon\le\varepsilon_0$. Apply Lemma~\ref{lem:kr-marked} on the compact shells
$\{x\in[-1+1/n,1-1/n]^N:|q_x-q|\le\varepsilon\}$ with the mark
$\ind_{\mathcal W_\varepsilon}(x)\ind\{|F_\zeta(x)/N-f|\le\varepsilon\}$, which drops the condition on $\FTAP$, and
let $n\to\infty$. Disintegrating over $\HN(x)=Nh$,
\[
 \E\mathsf N_{N,\varepsilon}\le\int_{\mathcal W_\varepsilon^\circ}\Upsilon(x)\,dx,\qquad
 \Upsilon(x)=p_{\nabla\HN(x)}(\mathbf g(x))\int_{I_x}Np_x(Nh)\,
 \E_{x,h}\bigl|\det\nabla^2F_\zeta(x)\bigr|\,dh,
\]
where $p_x$ is the conditional density of $\HN(x)$ given $\nabla\HN(x)=\mathbf g(x)$, $\E_{x,h}$ the
conditional expectation given in addition $\HN(x)=Nh$, and
$I_x=[f-V_\zeta(q_x,\mu_x)-\varepsilon,\,f-V_\zeta(q_x,\mu_x)+\varepsilon]\subset[-H_0,H_0]$, since the
conjugates and the Onsager integral are bounded.

\emph{Gaussian factors.} As in the proof of Lemma~\ref{lem:kr-count},
\begin{gather*}
 \log p_{\nabla\HN(x)}(\mathbf g(x))\le-\frac N2\log(2\pi c_r)-\frac{|b^{(r)}|^2}{2c_r}-N\mathcal Q_r(\chi_x,k_x),\\
 \mathcal Q_r(\chi,k)=\frac{a_r\chi k}{c_r}+\frac{a_r^2r\chi^2}{2c_r}-\frac{a_r(k+a_rr\chi)^2}{2c_r(c_r+a_rr)},
\end{gather*}
and $p_x(Nh)=(2\pi N\sigma_H^2(r))^{-1/2}\exp\{-N(h-h^\star_x)^2/(2\sigma_H^2(r))\}$ with
$h^\star_x=c_rd_x/(c_r+a_rr)$. Hence, for $x\in\mathcal W_\varepsilon^\circ$,
\[
 \int_{I_x}Np_x(Nh)\,\E_{x,h}|\det\nabla^2F_\zeta(x)|\,dh
 \le\frac{2\varepsilon N}{\sqrt{2\pi N\sigma_H^2(r)}}\exp\Bigl\{N\mathfrak D_N(\varepsilon,H_0)
 -N\inf_{h\in I_x}\frac{(h-h^\star_x)^2}{2\sigma_H^2(r)}\Bigr\}.
\]

\emph{Change of variables.} With $x_i=M(b_i)$ and the density
$\pi_0'(b)=(2\pi c)^{-1/2}\exp\{-b^2/(2c)+\theta[\Phi(q,b)-\Phi(0,0)]\}$ of \eqref{eq:kr-law},
\[
 \int_{\mathcal W_\varepsilon^\circ}\Upsilon(x)\,dx=\int\ind_{\mathcal W_\varepsilon^\circ}(x)\,\Upsilon(x)
 \prod_i\frac{U(b_i)}{\pi_0'(b_i)}\,\pi_0^{\otimes N}(db)
 \le\exp\Bigl\{N\sup_{x\in\mathcal W_\varepsilon^\circ}\Xi_N(x)\Bigr\},
\]
where, for $x\in\mathcal W_\varepsilon^\circ$,
\begin{align*}
 \Xi_N(x)\le{}&\mathfrak D_N(\varepsilon,H_0)+\frac1N\sum_i\log U(b_i)-\mathcal Q_r(\chi_x,k_x)
 -\theta\Bigl[\frac1N\sum_i\Phi(q,b_i)-\Phi(0,0)\Bigr]\\
 &+\frac12\log\frac c{c_r}+\frac{|b|^2}{2Nc}-\frac{|b^{(r)}|^2}{2Nc_r}
 -\inf_{h\in I_x}\frac{(h-h^\star_x)^2}{2\sigma_H^2(r)}+\frac{C\log N}N .
\end{align*}

\emph{Limits in the window.} Put $h^\star_0=cd_0/(c+aq)$. By \eqref{eq:kr-k0} and \eqref{eq:kr-h0},
$h_0-h^\star_0=\theta\sigma_H^2>0$, as in the proof of Lemma~\ref{lem:kr-normalization}. For $h\in I_x$,
Lemma~\ref{lem:kr-susywindow}\ref{lem:kr-susywindow-c}--\ref{lem:kr-susywindow-d} and the continuity of
$c_r,a_r$ give $|h-h^\star_x|\ge\theta\sigma_H^2-o_\varepsilon(1)$ uniformly in $x\in\mathcal W_\varepsilon$, and
$\sigma_H^2(r)\to\sigma_H^2$. With Lemma~\ref{lem:kr-susywindow}\ref{lem:kr-susywindow-b}--\ref{lem:kr-susywindow-c}
for the other terms and Lemma~\ref{lem:kr-susydet} for the first, in which $-\E_0\log U$ cancels the
average of $\log U(b_i)$,
\[
 \lim_{\varepsilon\downarrow0}\limsup_{N\to\infty}\sup_{x\in\mathcal W_\varepsilon^\circ}\Xi_N(x)
 \le\frac a2\chi_0^2-\mathcal Q_q(\chi_0,k_0)-\theta\bigl[\E_0\Phi(q,b)-\Phi(0,0)\bigr]-\frac{\theta^2\sigma_H^2}2 .
\]
Finally, $\frac a2\chi^2-\mathcal Q_q(\chi,k)=a(k-c\chi)^2/(2c(c+aq))$, which is the quadratic identity
in the proof of Lemma~\ref{lem:kr-count}, and $\theta[\E_0\Phi(q,b)-\Phi(0,0)]=D(\pi_0\Vert\mathcal N(0,c))$
by \eqref{eq:kr-law}. So the right side is $\Sigma_{\rm KR}(\pi_0,h_0)=\Sigma(f)$, by
Lemma~\ref{lem:kr-normalization}.
\end{proof}

At a strictly Plefka-stable first contact, Proposition~\ref{prop:kr-susy} is part of
\cite[Proposition~3.6]{Bou26b}, whose determinant estimate \cite[Lemma~B.2]{Bou26b} is stated
under the strict Plefka condition. The bound above uses no separation of the spectrum from zero,
no cutoff of the coordinates away from $\pm1$ (the diagonal entries are only bounded by
$Ce^{C\sqrt N}$), and no Plefka inequality for the empirical laws in the window. Boursier's earlier
annealed count \cite[Proposition~2.2]{Bou26a} uses a window that keeps the coordinates in
$[-1+\varepsilon^{1/2},1-\varepsilon^{1/2}]$ \cite[Definition~2.1]{Bou26a}; see the discussion after
Proposition~\ref{prop:stab-susy}.

\subsection{Retained states and the hypotheses on the covariance}\label{app:kr-remarks}

We bound the number of retained points given that one exists. Boursier observes that the
annealed and the quenched complexities of a count $\mathcal N$ may differ, by the exponential rate
of $\E[\mathcal N\mid\mathcal N\ge1]=\E\mathcal N/\Prob(\mathcal N\ge1)$ \cite[p.~2]{Bou26b}. For the
retained points this rate is controlled by the spherical code, also at a marginal first contact.

\begin{proposition}[Retained states are subexponentially many given existence]\label{prop:kr-retained}
Assume the hypotheses of Proposition~\ref{prop:kr-main}. Let $\eta$ be admissible and
$\delta\in(0,q/2)$, choose the parameters~\ref{step:kr-par1}--\ref{step:kr-par6}, and write
$\widehat C_N=\widehat C_N(\eta,\delta,w)$ and $s_w=\delta'/(q-2w)$. We put
$\widehat C_N^{\,0}=\ind\{\widehat C_N\ge1\}$.
\begin{thmlist}
\item\label{prop:kr-retained-a} \textup{(Deterministic bounds.)} For every choice of
 $\mathcal W_N(w)$, $w\in(0,w_1]$, $N\ge3\gamma^{-1/2}$ and $\lambda\ge0$,
 \[
  \Prob(\widehat C_N\ge1)\le\E\widehat C_N^{\,\lambda}\le n_N(s_w)^\lambda\,\Prob(\widehat C_N\ge1),
 \]
 and, if $\Prob(\widehat C_N\ge1)>0$, then $1\le\E[\widehat C_N\mid\widehat C_N\ge1]\le n_N(s_w)$.
\item\label{prop:kr-retained-b} \textup{(Exponents.)} Let $L_\lambda^-$ and $L_\lambda^+$ be
 $\liminf_{w\downarrow0}\liminf_N$ and $\limsup_{w\downarrow0}\limsup_N$ of
 $N^{-1}\log\E\widehat C_N^{\,\lambda}$. For every choice of $\mathcal W_N(w)$,
 \[
  L_\lambda^+\le\Sigma(f)+\theta\bigl(|f_\eta-f|+C_V\eta\bigr)-\frac\lambda2\log\Bigl(1-\frac{\delta'}q\Bigr),
 \]
 and for the choice of Lemma~\ref{lem:kr-count},
 $L_\lambda^-\ge\Sigma_{\rm KR}(\pi_\eta,h_\eta)+\frac12\log(1-\delta'/q)$, where the last term may be
 omitted when $\lambda\ge1$.
\item\label{prop:kr-retained-c} \textup{(Limits.)} For the choice of Lemma~\ref{lem:kr-count} there is
 $C_{\mathrm{ret}}$, independent of $\eta$, $\delta$ and $\lambda$, with
 \[
  \bigl|L_\lambda^\pm-\Sigma(f)\bigr|\le C_{\mathrm{ret}}\eta+\frac{1\vee\lambda}2\log\frac1{1-\delta'/q}\qquad(\lambda\ge0).
 \]
 So every iterated limit of $N^{-1}\log\E\widehat C_N^{\,\lambda}$, with $N\to\infty$ innermost and then
 $w\downarrow0$, tends to $\Sigma(f)$ as $\delta\downarrow0$ and then $\eta\downarrow0$, uniformly for
 $\lambda$ in compact subsets of $[0,\infty)$. Moreover, at every admissible $\eta$,
 $\Prob(\widehat C_N\ge1)>0$ for all sufficiently small $w$ and all large $N$, and
 \[
  \lim_{\delta\downarrow0}\limsup_{w\downarrow0}\limsup_{N\to\infty}\frac1N
  \log\E\bigl[\widehat C_N\bigm|\widehat C_N\ge1\bigr]=0 .
 \]
\end{thmlist}
\end{proposition}

\begin{proof}
\ref{prop:kr-retained-a} By Lemma~\ref{lem:kr-separation}, $\widehat C_N\le n_N(s_w)$ for every
realization. Since $\widehat C_N$ is a nonnegative integer,
$\ind\{\widehat C_N\ge1\}\le\widehat C_N^{\,\lambda}\le n_N(s_w)^\lambda\ind\{\widehat C_N\ge1\}$; take
expectations, and divide by $\Prob(\widehat C_N\ge1)$ for $\lambda=1$.

\ref{prop:kr-retained-b} Let $\MN=\max_{\cube}\FTAP$. As in the proof of~\ref{step:stab-k2},
every retained point $x$ has
\[
 \FTAP(x)/N\ge f-|f_\eta-f|-C_V\eta-2C_V\varrho(w)-N^{-2},
\]
so
$\{\widehat C_N\ge1\}\subset\{\MN\ge N(f-\varepsilon_N)\}$ with
$\varepsilon_N=|f_\eta-f|+C_V\eta+2C_V\varrho(w)+N^{-2}$. By Markov's inequality at $\theta=\thf$,
$\Prob(\widehat C_N\ge1)\le e^{-\theta N(f-\varepsilon_N)}\E e^{\theta\MN}$. By
Theorem~\ref{thm:max}\ref{thm:max-b} and Proposition~\ref{prop:max-fm}\ref{prop:max-fm-a},
$\limsup_NN^{-1}\log\E e^{\theta\MN}\le\Lambda(\theta)$, and $\Lambda(\theta)-\theta f=\Sigma(f)$. With
\ref{prop:kr-retained-a}, Lemma~\ref{lem:kr-code} and $s_w\to\delta'/q$ as $w\downarrow0$, this gives the
bound on $L_\lambda^+$; only the deterministic facts of Lemma~\ref{lem:kr-separation} were used. For the
lower bound, \ref{prop:kr-retained-a} with $\lambda=1$ gives
$\E\widehat C_N^{\,\lambda}\ge\Prob(\widehat C_N\ge1)\ge\E\widehat C_N/n_N(s_w)$; apply
Lemmas~\ref{lem:kr-count} and~\ref{lem:kr-code}. For $\lambda\ge1$, $\widehat C_N^{\,\lambda}\ge\widehat C_N$.

\ref{prop:kr-retained-c} By Lemma~\ref{lem:kr-normalization} and \eqref{eq:kr-feta},
$|\Sigma_{\rm KR}(\pi_\eta,h_\eta)-\Sigma(f)|\le C_\Sigma\eta$ and $|f_\eta-f|\le(L+C_V)\eta$; put
$C_{\mathrm{ret}}=C_\Sigma+\theta(L+2C_V)$. Substituting in~\ref{prop:kr-retained-b}, using
$\frac\lambda2\le\frac{1\vee\lambda}2$ on the upper side, the fact that the logarithmic term appears on
the lower side only for $\lambda<1$, and $L_\lambda^-\le L_\lambda^+$, gives the two-sided bound. Its right
side is nondecreasing in $\lambda$, which gives the uniformity on compact sets, and $\delta'=\delta$ for
small $\delta$. The number $\Sigma_{\rm KR}(\pi_\eta,h_\eta)$ is finite, so by Lemma~\ref{lem:kr-count}
there is $w_2\in(0,w_1]$ such that, for $w\le w_2$, $\E\widehat C_N>0$, hence
$\Prob(\widehat C_N\ge1)>0$, for all large $N$. For such $w$ and $N$,
\ref{prop:kr-retained-a} and Lemma~\ref{lem:kr-code} give
$0\le\limsup_NN^{-1}\log\E[\widehat C_N\mid\widehat C_N\ge1]\le-\frac12\log(1-s_w)$, which tends to
$-\frac12\log(1-\delta'/q)$ as $w\downarrow0$, and to zero as $\delta\downarrow0$.
\end{proof}

\begin{remark}[Scope of Proposition~\ref*{prop:kr-retained}]\label{rem:kr-retained}
At a strictly Plefka-stable first contact, Boursier's isolated count $\widehat N$ satisfies
$\E\widehat N\le n_N\,\Prob(\widehat N\neq0)$ (Step~(2) of \cite[Section~3.2]{Bou26b}, with
\cite[Lemma~3.11 and Section~3.7]{Bou26b}), which with \cite[Proposition~3.6]{Bou26b} and the proof
of \cite[Proposition~3.3]{Bou26b} gives the analogue of Proposition~\ref{prop:kr-retained}. Here
the contact may be marginal, the data are tilted, and $\Sigma(f)$ is reached only as
$\eta\downarrow0$. Part~\ref{prop:kr-retained-c} covers the moments $0<\lambda<1$ behind the
heuristic $\log n=\lim_{\epsilon\downarrow0}(n^\epsilon-1)/\epsilon$ of \cite[p.~2]{Bou26b}. At a
fixed cutoff the bound of~\ref{prop:kr-retained-a} is exponential, with rate $-\frac12\log(1-\delta'/(q-2w))$.
Nothing is claimed about all critical points of $F_{\zf}$ in the windows (without the
marks~\ref{step:kr-ret4}--\ref{step:kr-ret5}) or all TAP states given that one exists.
\end{remark}

We finally record where the hypotheses on $\xi$ enter the construction.

\begin{remark}[A sign window instead of convexity]\label{rem:kr-window}
Convexity of $\xi$ on $[-1,0]$ enters the proof of Proposition~\ref{prop:kr-main} only in
\eqref{eq:kr-bregman}, whose left side
is $2N\mathcal D_\xi(R,q_{\alpha\beta})$, $q_{\alpha\beta}\in\{q,r\}$, with
$\mathcal D_\xi(R,Q)=\xi(R)-\xi(Q)-\xi'(Q)(R-Q)$, since $K_r-\xi$ is affine. On the slice
$\{q_y=q,\ \ov xy=r\}$, $\ov yz\ge2r^2/q-q>-q$ (decompose along $x$ and $x^\perp$), and for
$R\le0\le Q\le1$, $\partial_Q\mathcal D_\xi(R,Q)=\xi''(Q)(Q-R)\ge0$, so
$\mathcal D_\xi(R,Q)\ge\mathcal D_\xi(R,0)=\xi(R)$. Hence
Lemma~\ref{lem:kr-comparison}\ref{lem:kr-comparison-a} holds whenever $\xi\ge0$ on $[-q,0]$,
and~\ref{lem:kr-comparison-b} on $I=[\delta,q-\delta]$ once also $\omega\le\delta^2/(4q)$, which
keeps the overlaps of the net points above $-q$. So the lower bound of
Proposition~\ref{prop:kr-main} extends to a level $f$ of a mixture with two active degrees that is
not convex on $[-1,1]$ but has $\xi\ge0$ on $[-\qf,0]$, provided the variational inputs of
Lemmas~\ref{lem:kr-contact} and~\ref{lem:matrix-band}, proved or cited only under \condH, hold at
$f$; for $\xi(t)=\beta_2^2t^2+\beta_3^2t^3$ this reads $\qf\le\beta_2^2/\beta_3^2$, while convexity
on $[-\qf,0]$ would need $\qf\le\beta_2^2/(3\beta_3^2)$. The matching upper bound needs a positivity
principle in the spirit of \cite{PanPos}, and the corresponding step of \cite[Theorem~2.1]{CGKLM}
is only sketched there, so we state no two-sided result for non-convex mixtures.
\end{remark}

\begin{remark}[Where the covariance enters]\label{rem:kr-summable}
Besides the variational inputs and convexity on $[-1,0]$ (Remark~\ref{rem:kr-window}), $\xi$ enters
through two active degrees and bounded derivatives. Two active degrees give $\det S_0>0$ and
$\sigma_H^2>0$ (Lemma~\ref{lem:kr-energy}), hence the energy shift, the strict signs of $A_2$ and
$\tau$ in Lemma~\ref{lem:kr-block}, the nondegeneracy in Lemma~\ref{lem:kr-c3} and the conditional
density of the energy in Lemma~\ref{lem:kr-count}. Derivatives of $\xi$ on $[-1,1]$ of order at
most eight bound the pointwise variance and the canonical metric of $D^3\HN$
(Lemmas~\ref{lem:matrix-third} and~\ref{lem:kr-c3}). The class $\classC$ enters only through the
cited facts on the Parisi equation and on the TAP correction (Remark~\ref{rem:rpde}); a time change reducing them to a
polynomial covariance is a natural route to replacing $\xi\in\classC$ in
Proposition~\ref{prop:kr-main} by $\sum_pp^8\beta_p^2<\infty$, under strict convexity and two active
degrees. We do not claim this extension: a complete proof must also revisit the uses of the power
series of $\xi$ (as in Lemma~\ref{lem:kr-energy}) and the regularity of $\HN$ behind the Kac--Rice
formula. No statement of this paper depends on it.
\end{remark}

%% file: sections/appF-nearcritical.tex
\section{The near-critical \texorpdfstring{$2+4$}{2+4} family}\label{app:nc}

Theorem~\ref{thm:stab}\ref{thm:stab-b} needs every contact to be marginal, a degenerate condition.
We show that it holds at all levels close to $\feq$ for the mixtures $\xi(t)=bt^2+et^4$ slightly
above the critical Sherrington--Kirkpatrick covariance $\xi_{\mathrm c}(t)=t^2/2$:
Proposition~\ref{prop:nc-example} describes their constrained minimizers, and
Appendix~\ref{app:nc-marginal} proves Lemma~\ref{lem:stab-family}, which gives
Theorem~\ref{thm:stab}\ref{thm:stab-c} (Section~\ref{sec:stab-family}). No finite-$N$ estimate
enters.

\subsection{Statement and outline}\label{app:nc-example}

For $\theta\in[0,1]$, the \emph{constrained minimizer at $\theta$}, written $\zeta_\theta$, is the
minimizer of $\zeta\mapsto P_\zeta$ over $\{\zeta\in\Pspace([0,1]):\zeta(\{0\})\ge\theta\}$. For
the covariances of this appendix it exists and is unique
(Lemma~\ref{lem:var-basic}\ref{lem:var-basic-a}); in particular $\zeta_0=\zP$, and
$p(\theta)=P_{\zeta_\theta}$. When
$\argmin_{(0,1]}\Obs_{\zeta_\theta}$ is a nonempty compact subset of $(0,1)$, its least point is
the \emph{first positive contact} of $\zeta_\theta$; for $\theta=\thf$ it is the point $\qf$ of
Section~\ref{sec:intro-structure}.

\begin{proposition}[Marginal first contacts near criticality]\label{prop:nc-example}
There is $\varepsilon_*\in(0,\frac12]$ such that every covariance
\begin{equation}\label{eq:nc-family}
 \xi(t)=bt^2+et^4,\qquad b>\tfrac12,\qquad e>0,\qquad b-\tfrac12+e<\varepsilon_*,
\end{equation}
has the following properties, where $a_0=\xi''(0)=2b$.
\begin{thmlist}
\item\label{prop:nc-example-a} \textup{(Equilibrium.)} $\supp\zP=[0,v_{\mathrm P}]$ for some
 $v_{\mathrm P}\in(0,1)$, $\zP(\{0\})=0$, and $U_{\mathrm P}:=\Phi_{\zP,xx}(0,0)=a_0^{-1/2}$.
\item\label{prop:nc-example-b} \textup{(Marginal first contact.)} There is $\theta_0\in(0,1)$ such
 that for every $\theta\in(0,\theta_0)$
 \[
  \zeta_\theta(\{0\})=\theta,\qquad
  \supp\zeta_\theta=\{0\}\cup[q_\theta,v_\theta],\qquad 0<q_\theta<v_\theta<1,
 \]
 $q_\theta$ is the first positive contact of $\zeta_\theta$, and it is marginal:
 \[
  \xi''(q_\theta)\,\E\,\Phi_{\zeta_\theta,xx}(q_\theta,X_{q_\theta})^2=1,
 \]
 where $X$ is the Auffinger--Chen process of $\zeta_\theta$.
\item\label{prop:nc-example-c} \textup{(Levels.)} The number $\fzero:=p(\theta_0)$ satisfies
 $\feq<\fzero<\fone$, and $\thf\in(0,\theta_0)$ for every $f\in(\feq,\fzero)$. Hence, for every
 such $f$, $\zf=\zeta_{\thf}$ and its first positive contact $\qf=q_{\thf}$ has $\Plef(\qf)=0$.
\item\label{prop:nc-example-d} \textup{(Onset.)} Put
 $D_{\mathrm P}=\xi''''(0)\,U_{\mathrm P}^2+a_0^3\,\Phi_{\zP,xxxx}(0,0)^2$. Then
 $D_{\mathrm P}>0$, and
 \begin{equation}\label{eq:nc-onset}
  \lim_{\theta\downarrow0}\frac{q_\theta}{\theta}=\frac{3\sqrt{a_0}}{D_{\mathrm P}},\qquad
  \lim_{\theta\downarrow0}\frac{\zeta_\theta(\{q_\theta\})}{\theta}=\frac12 .
 \end{equation}
\end{thmlist}
\end{proposition}

\begin{remark}[The hypotheses]\label{rem:nc-hyp}
The covariances \eqref{eq:nc-family} are polynomials with $\xi''(t)=2b+12et^2>0$ on $[-1,1]$, so
they satisfy \condH, lie in $\classC$ and have two active degrees: the hypotheses of
Theorems~\ref{thm:exist}, \ref{thm:struct} and~\ref{thm:stab}\ref{thm:stab-a}. The condition
$b>\frac12$, that is $\xi''(0)>\xi_{\mathrm c}''(0)=1$, forces $\zP\neq\delta_0$; $e>0$ gives
$\xi'''>0$ on $(0,1]$, which puts the origin in $\supp\zP$, and $\xi''''(0)>0$, which makes
$D_{\mathrm P}$ positive.
\end{remark}

\paragraph{Outline of the proof.} We perturb the critical pair $(\xi_{\mathrm c},\delta_0)$; write
$\Gamma=\Gz\zeta$. At the critical pair $\Gamma<\id$ on $(0,1]$, so near criticality every contact
lies in a small-time region $(0,t_{\mathrm c})$. On a plateau of $\zeta$ (an open interval of zero
mass), $\Gamma'''$ is an explicit expectation of the spatial derivatives of $\Phi_\zeta$ of orders
two to four along the Auffinger--Chen process (Appendix~\ref{app:nc-plateau}); at the critical pair
and at the origin it is at least $4$, because the third and fourth derivatives of $\log\cosh$ at $0$
are $0$ and $-2$. So near criticality $\Gamma'$ is strictly convex on every plateau below
$t_{\mathrm c}$. At a contact $t$, $\Gamma(t)=t$ and $\Gamma'(t)\le1$; if two contacts $t_1<t_2$
enclosed a plateau, the mean of $\Gamma'$ over $(t_1,t_2)$ would equal one and would be strictly
smaller than $\frac12(\Gamma'(t_1)+\Gamma'(t_2))\le1$. This adapts the near-critical gap-exclusion
mechanism of Zhou \cite{Zho26a} to contacts that may carry no mass (Appendix~\ref{app:nc-gap}), and
it applies to every constrained minimizer, since all of them stay close to $\delta_0$
(Appendix~\ref{app:nc-confine}). So the positive support of $\zeta_\theta$ is an interval
$[q_\theta,v_\theta]$ of contacts, whence $\Gamma'(q_\theta)=1$; the onset limits
\eqref{eq:nc-onset} come from the cubic behaviour of $\Gamma-\id$ on $(0,q_\theta)$
(Appendix~\ref{app:nc-proof}).

\subsection{Plateau calculus}\label{app:nc-plateau}

\paragraph{Notation.} In Appendices~\ref{app:nc-plateau}--\ref{app:nc-confine} the covariance varies
in the set
\[
 \mathcal K=\bigl\{\xi_{b,e}(t)=bt^2+et^4:\ \tfrac12\le b\le1,\ 0\le e\le1\bigr\},\qquad
 d(\xi_{b,e})=\bigl(b-\tfrac12\bigr)+e .
\]
It contains $\xi_{\mathrm c}=\xi_{1/2,0}$ and, since $\varepsilon_*\le\frac12$, the family
\eqref{eq:nc-family}, on which $d(\xi)<\varepsilon_*$. Every $\xi\in\mathcal K$ lies in $\classC$, and
on $[0,1]$
\begin{equation}\label{eq:nc-K}
 1\le\xi''\le14,\qquad |\xi''-1|\le12\,d(\xi),\qquad 0\le\xi'''\le24\,d(\xi),\qquad
 \xi''''=24e\le24\,d(\xi).
\end{equation}
For $\xi\in\mathcal K$ and $\zeta\in\Pspace([0,1])$ let $\Phi=\Phi^\xi_\zeta$, $X$, $\Gamma=\Gz\zeta^\xi$,
$\Obs=\Obs^\xi_\zeta$ and $P^\xi_\zeta$ be the objects of
Sections~\ref{sec:intro-setting}--\ref{sec:intro-structure} computed with the covariance $\xi$, and
write
\[
 a=\xi'',\qquad \alpha=\cdf\zeta,\qquad u=\partial_x\Phi,\qquad
 U_t=\partial_x^2\Phi(t,X_t),\quad V_t=\partial_x^3\Phi(t,X_t),\quad W_t=\partial_x^4\Phi(t,X_t).
\]
A tilde marks the same objects for a second pair $(\tilde\xi,\tilde\zeta)$. Let
$m_1(\zeta)=\int t\,\zeta(dt)=\|1-\alpha\|_{L^1}$, the Wasserstein distance from $\zeta$ to
$\delta_0$. A \emph{plateau} of $\zeta$ is an open interval $I\subset(0,1)$ with $\zeta(I)=0$; on
it $\alpha$ is constant, and the constant is the \emph{height} of $I$.

For a fixed covariance the bounds below are due to Auffinger and Chen \cite{AC15a}.
Lemma~\ref{lem:pde-standing}\ref{lem:pde-standing-a} bounds $\partial_x^ku$ for $1\le k\le3$
uniformly on $\mathcal K$, where $\xi'(1)\le6$ (its upper bounds increase with $\xi'(1)$);
we need every order, and stability in the pair
$(\xi,\zeta)$.

\begin{lemma}[Jets and their stability]\label{lem:nc-jets}
For every $k\ge0$ there is $C_k<\infty$, and there is $C<\infty$, such that the following hold for
all $\xi,\tilde\xi\in\mathcal K$ and $\zeta,\tilde\zeta\in\Pspace([0,1])$.
\begin{thmlist}[beginpenalty=10000]
\item\label{lem:nc-jets-a} $\partial_x^j\Phi$ is continuous on $[0,1]\times\R$ for every
 $j\ge0$, $\Phi(t,\cdot)$ is even, $|u|\le1$, $0<\partial_xu\le1$, and
 $\sup_{[0,1]\times\R}|\partial_x^ku|\le C_k$.
\item\label{lem:nc-jets-b} \textup{(Jet stability.)}
 \begin{equation}\label{eq:nc-jet}
  \|\partial_x^k(u-\tilde u)\|_\infty
  \le C_k\bigl(\|a-\tilde a\|_{L^1}+\|a\alpha-\tilde a\tilde\alpha\|_{L^1}\bigr).
 \end{equation}
\item\label{lem:nc-jets-c} \textup{(Small times.)} $\E X_t^2\le2\|a\|_\infty t+2\|a\|_\infty^2t^2$
 for $t\in[0,1]$.
\item\label{lem:nc-jets-d} \textup{(Continuity in the pair.)}
 \begin{equation}\label{eq:nc-continuity}
  \sup_{[0,1]}\bigl(|\Gamma-\tilde\Gamma|+|\Gamma'-\tilde\Gamma'|\bigr)
  +\bigl|P^\xi_\zeta-P^{\tilde\xi}_{\tilde\zeta}\bigr|
  \le C\bigl(\|a-\tilde a\|_\infty+\|a\alpha-\tilde a\tilde\alpha\|_{L^1}\bigr).
 \end{equation}
\end{thmlist}
\end{lemma}

\begin{proof}
\ref{lem:nc-jets-a} The continuity of $\partial_x^j\Phi$ is
Lemma~\ref{lem:pde-standing}\ref{lem:pde-standing-a}. For a fixed covariance, its uniform
convergence along weakly convergent sequences of measures is \cite[Proposition~1(i)]{AC15a}; here
$\Phi$ is the solution of \cite{AC15a} by Lemma~\ref{lem:pde-measurable}\ref{lem:pde-measurable-iv}.
Evenness is Lemma~\ref{lem:pde-measurable}\ref{lem:pde-measurable-i}, and $|u|\le1$,
$0<\partial_xu\le1$ are Lemma~\ref{lem:pde-class}\ref{lem:pde-class-a} and, for general $\zeta$,
Lemma~\ref{lem:pde-measurable}\ref{lem:pde-measurable-i}, with
Lemma~\ref{lem:pde-standing}\ref{lem:pde-standing-a} for the strict positivity (see also \cite[Proposition~2]{AC15a}). By the
uniform convergence it suffices to prove the bounds of \ref{lem:nc-jets-a}
and~\ref{lem:nc-jets-b} for finitely supported $\zeta$ and $\tilde\zeta$, which we now assume. Then
$\Phi$ is given on each step of $\alpha$ by the Cole--Hopf formula \eqref{eq:pde-colehopf}, all
$\partial_x^ju$ are continuous across the steps, and $u$ solves
\begin{equation}\label{eq:nc-u}
 \partial_tu+\tfrac a2\,\partial_x^2u+a\alpha\,u\,\partial_xu=0,\qquad u(1,x)=\tanh x .
\end{equation}
For $k\ge2$, differentiating \eqref{eq:nc-u} $k$ times in $x$ gives a linear backward equation for
$\partial_x^ku$ with second-order coefficient $a/2$, drift $a\alpha u$, zeroth-order coefficient
$(k+1)a\alpha\,\partial_xu$, and a forcing term that is a polynomial in
$\partial_xu,\dots,\partial_x^{k-1}u$ whose coefficients are bounded in terms of $k$. Its terminal value
$\partial_x^k\tanh$ is bounded. The Feynman--Kac formula along the Auffinger--Chen process started
at $(t,x)$, chained over the steps (the maximum principle), bounds $\partial_x^ku$ in terms of $k$,
$\|a\|_\infty$ and $C_1,\dots,C_{k-1}$, and induction on $k$ gives \ref{lem:nc-jets-a}.

\ref{lem:nc-jets-b} The difference $w=u-\tilde u$ vanishes at $t=1$ and solves
\begin{equation}\label{eq:nc-w}
 \partial_tw+\tfrac a2\,\partial_x^2w+a\alpha\,u\,\partial_xw+a\alpha\,(\partial_x\tilde u)\,w
 =-\tfrac{a-\tilde a}2\,\partial_x^2\tilde u-(a\alpha-\tilde a\tilde\alpha)\,\tilde u\,\partial_x\tilde u .
\end{equation}
Its drift and zeroth-order coefficients are bounded, and by \ref{lem:nc-jets-a} the supremum over
$x$ of the right side at time $s$ is at most $C_2|a-\tilde a|(s)+|a\alpha-\tilde a\tilde\alpha|(s)$.
The maximum principle bounds $\|w\|_\infty$ by a constant times the time integral of this
supremum, which is \eqref{eq:nc-jet} for $k=0$. Differentiating \eqref{eq:nc-w} $k$ times gives an
equation of the same form for $\partial_x^kw$, with bounded coefficients, whose forcing is
$\partial_x^k$ of the right side of \eqref{eq:nc-w}, bounded as before by \ref{lem:nc-jets-a},
plus products of $\partial_x^jw$, $j<k$, with derivatives of $u$ and $\tilde u$. Induction on $k$
gives \eqref{eq:nc-jet}. All constants depend only on $k$ and on the bounds \eqref{eq:nc-K}.

\ref{lem:nc-jets-c} This follows from $X_t=\int_0^ta\alpha\,u(s,X_s)\,ds+\int_0^t\sqrt a\,dB_s$ and
$|\alpha u|\le1$.

\ref{lem:nc-jets-d} Drive $X$ and $\tilde X$ by the same Brownian motion and let
$D=\sup_{t\le1}|X_t-\tilde X_t|$. Since $|\partial_xu|\le1$, $\alpha\le1$ and $|\tilde u|\le1$,
\[
 \bigl|a\alpha\,u(s,X_s)-\tilde a\tilde\alpha\,\tilde u(s,\tilde X_s)\bigr|
 \le\|a\|_\infty\bigl(|X_s-\tilde X_s|+\|u-\tilde u\|_\infty\bigr)+|a\alpha-\tilde a\tilde\alpha|(s).
\]
Gronwall's inequality, applied pathwise, and Doob's inequality for the martingale
$\int_0^t(\sqrt a-\sqrt{\tilde a})\,dB_s$ give
\[
 (\E D^2)^{1/2}\le e^{\|a\|_\infty}\Bigl(\|a\|_\infty\|u-\tilde u\|_\infty
 +\|a\alpha-\tilde a\tilde\alpha\|_{L^1}+2\|\sqrt a-\sqrt{\tilde a}\|_{L^2}\Bigr),
\]
and $|\sqrt a-\sqrt{\tilde a}|\le\frac12|a-\tilde a|$ because $a,\tilde a\ge1$. Now
$|\Gamma(t)-\tilde\Gamma(t)|\le2\E|u(t,X_t)-\tilde u(t,\tilde X_t)|\le2(\|u-\tilde u\|_\infty+\E D)$.
By Lemma~\ref{lem:pde-standing}\ref{lem:pde-standing-c}, $\Gamma'=a\,\E U_t^2$, and since
$0<U\le1$ and $|\partial_x^2u|\le C_2$,
$|\Gamma'(t)-\tilde\Gamma'(t)|\le|a-\tilde a|(t)+2\|a\|_\infty(\|\partial_x(u-\tilde u)\|_\infty+C_2\E D)$.
For the Parisi functional, integrating the Parisi equation in time at $x=0$ (on each step for
finitely supported $\zeta$, and in general by approximation, Lemma~\ref{lem:pde-measurable}) gives
\[
 P^\xi_\zeta=\frac12\int_0^1a(t)\bigl(\partial_xu(t,0)+\alpha(t)\,u(t,0)^2-t\,\alpha(t)\bigr)\,dt ,
\]
whose integrand changes by at most
$|a-\tilde a|+\|a\|_\infty(\|\partial_x(u-\tilde u)\|_\infty+2\|u-\tilde u\|_\infty)+2|a\alpha-\tilde a\tilde\alpha|$.
Together with \eqref{eq:nc-jet} for $k=0,1$, these bounds give \eqref{eq:nc-continuity}.
\end{proof}

\begin{lemma}[Calculus of $\Gamma$]\label{lem:nc-plateau}
Let $\xi\in\mathcal K$ and $\zeta\in\Pspace([0,1])$.
\begin{thmlist}
\item\label{lem:nc-plateau-a} $\Gamma\in C^1([0,1])$ with $\Gamma'(t)=\xi''(t)\,\E U_t^2$;
 moreover $\Gamma(0)=0$, $U_0=\Phi_{xx}(0,0)>0$ and $V_0=0$.
\item\label{lem:nc-plateau-b} Let $I$ be a plateau of $\zeta$ with height $\omega$. Then
 $\Gamma\in C^3(I)$ and, on $I$, with $U=U_t$, $V=V_t$, $W=W_t$,
 \begin{align}
  \Gamma''&=\xi'''\,\E U^2+\xi''^2\,\E\bigl(V^2-2\omega U^3\bigr),\label{eq:nc-second}\\
  \Gamma'''&=\xi''''\,\E U^2+3\xi''\xi'''\,\E\bigl(V^2-2\omega U^3\bigr)
   +\xi''^3\,\E\bigl(W^2-12\omega UV^2+6\omega^2U^4\bigr).\label{eq:nc-third}
 \end{align}
 For fixed $\omega$ the right sides are continuous functions of $t\in[0,1]$. In particular, at an
 endpoint of $I$ the one-sided derivative of $\Gamma'$ from inside $I$ exists and equals the right
 side of \eqref{eq:nc-second} there.
\item\label{lem:nc-plateau-c} For $t\in[0,1)$ the right derivative $\Gamma''(t+)$ of $\Gamma'$ at
 $t$ exists and equals the right side of \eqref{eq:nc-second} at $t$ with $\omega=\zeta([0,t])$;
 for $t\in(0,1]$ the left derivative $\Gamma''(t-)$ exists and equals the same expression with
 $\omega=\zeta([0,t))$. In particular, for $t\in(0,1)$,
 \begin{equation}\label{eq:nc-jump}
  \Gamma''(t-)-\Gamma''(t+)=2\,\xi''(t)^2\,\zeta(\{t\})\,\E U_t^3 .
 \end{equation}
\end{thmlist}
\end{lemma}

\begin{proof}
\ref{lem:nc-plateau-a} The first claim is
Lemma~\ref{lem:pde-standing}\ref{lem:pde-standing-c}. Since $X_0=0$ and $\Phi(0,\cdot)$ is even
(Lemma~\ref{lem:nc-jets}\ref{lem:nc-jets-a}), $\Gamma(0)=u(0,0)^2=0$ and
$V_0=\partial_x^3\Phi(0,0)=0$, and $U_0>0$ by Lemma~\ref{lem:nc-jets}\ref{lem:nc-jets-a}.

\ref{lem:nc-plateau-b} On $I$ the Parisi equation has the constant coefficient $\omega$. Taking, in
Lemma~\ref{lem:pde-measurable}, step approximations of $\alpha$ that equal $\omega$ on a compact
subinterval $[s,t']\subset I$, one sees that on $[s,t']\times\R$ the function $\Phi$ is the
Cole--Hopf solution \eqref{eq:pde-colehopf} with coefficient $\omega$ and terminal value
$\Phi(t',\cdot)$. So $u$ and its spatial derivatives are $C^1$ in $t$ on $I$; $u$ solves
\eqref{eq:nc-u} with $\alpha=\omega$, and its spatial derivatives solve the differentiated
equations, classically. It\^o's formula
gives, on $I$,
\[
 dU_t=-\omega\xi''(t)\,U_t^2\,dt+\sqrt{\xi''(t)}\,V_t\,dB_t,\qquad
 dV_t=-3\omega\xi''(t)\,U_tV_t\,dt+\sqrt{\xi''(t)}\,W_t\,dB_t,
\]
with bounded integrands, so
\[
 \frac{d}{dt}\E U_t^2=\xi''\,\E\bigl(V^2-2\omega U^3\bigr),\qquad
 \frac{d}{dt}\E\bigl(V_t^2-2\omega U_t^3\bigr)=\xi''\,\E\bigl(W^2-12\omega UV^2+6\omega^2U^4\bigr).
\]
The right sides are continuous in $t$, by the continuity of $\partial_x^j\Phi$, the continuity of
the paths of $X$ and bounded convergence. Differentiating $\Gamma'=\xi''\,\E U^2$ twice gives
\eqref{eq:nc-second} and \eqref{eq:nc-third}. The statement at the endpoints follows from the mean
value theorem.

\ref{lem:nc-plateau-c} This is \cite[Proposition~3, (21)--(22)]{AC15a}. There the expectations are
taken along the time-changed Brownian motion $M(t)=\int_0^t\sqrt{\xi''(s)}\,dB_s$, with the
exponential weight of \cite[(13)]{AC15a}. By It\^o's formula and the Parisi equation, the logarithm
of that weight is the Girsanov exponent of the drift $\xi''\cdf\zeta\,\Phi_x$, so by Girsanov's
theorem the weight is the density of the law of $(X_s)_{s\le t}$ with respect to that of
$(M(s))_{s\le t}$. Hence the statements of \cite{AC15a} hold along $X$; in particular $\Gamma$ is
the function of \cite[(18)]{AC15a}. Subtracting the two one-sided derivatives gives
\eqref{eq:nc-jump}.
\end{proof}

\subsection{Gap exclusion near criticality}\label{app:nc-gap}

\paragraph{The critical pair.} For $\xi_{\mathrm c}(t)=t^2/2$ and $\zeta=\delta_0$ one has
$\xi_{\mathrm c}''\equiv1$ and $\cdf{\delta_0}\equiv1$, and
\begin{equation}\label{eq:nc-critical}
 \Phi_{\mathrm c}(t,x):=\Phi^{\xi_{\mathrm c}}_{\delta_0}(t,x)=\log\cosh x+\tfrac12(1-t),\qquad
 \bigl(\partial_x^2,\partial_x^3,\partial_x^4\bigr)\Phi_{\mathrm c}(t,0)=(1,0,-2).
\end{equation}
Its Auffinger--Chen process solves $dX_t=\tanh(X_t)\,dt+dB_t$, the Doob transform of Brownian
motion by the space-time harmonic function $e^{-t/2}\cosh x$, so
$\E F(X_t)=e^{-t/2}\,\E\bigl[F(\sqrt tZ)\cosh(\sqrt tZ)\bigr]$ for a standard Gaussian variable $Z$.
Hence $\Gamma_{\mathrm c}:=\Gz{\delta_0}^{\xi_{\mathrm c}}$ satisfies $\Gamma_{\mathrm c}(0)=0$ and
\begin{equation}\label{eq:nc-critical-gamma}
 \Gamma_{\mathrm c}'(t)=\E\sech^4(X_t)=e^{-t/2}\,\E\sech^3(\sqrt tZ)<1\qquad(0<t\le1),
\end{equation}
so $\Gamma_{\mathrm c}(t)<t$ for $0<t\le1$. The obstacle $\Obs_{\mathrm c}$ of $\delta_0$ at
$\xi_{\mathrm c}$ therefore has
$\Obs_{\mathrm c}'=\frac12(t-\Gamma_{\mathrm c})>0$ on $(0,1]$, so
$\Obs_{\mathrm c}\ge\Obs_{\mathrm c}(0)=\int\Obs_{\mathrm c}\,d\delta_0$, and
\eqref{eq:loc-first-order} with \eqref{eq:loc-moving-mass} give $P^{\xi_{\mathrm c}}_\mu\ge
P^{\xi_{\mathrm c}}_{\delta_0}$ for every $\mu\in\Pspace([0,1])$. By the strict convexity of the
Parisi functional \cite[Theorem~2]{AC15b}, $\delta_0$ is the unique Parisi measure of
$\xi_{\mathrm c}$.

Both features of the critical pair used in the outline survive a small perturbation.

\begin{lemma}[Two estimates near the critical pair]\label{lem:nc-near}
There are $t_{\mathrm c}\in(0,1)$ and $\varepsilon_2>0$ such that, for every $\xi\in\mathcal K$
with $d(\xi)<\varepsilon_2$ and every $\zeta\in\Pspace([0,1])$ with $m_1(\zeta)<\varepsilon_2$, the
following hold.
\begin{thmlist}
\item\label{lem:nc-near-a} \textup{(Third derivative.)} For every plateau $I$ of $\zeta$,
 \begin{equation}\label{eq:nc-third-positive}
  \Gamma'''(t)\ge2\qquad\text{for every } t\in I\cap(0,t_{\mathrm c}).
 \end{equation}
\item\label{lem:nc-near-b} \textup{(Below the diagonal.)} $\Gamma(t)<t$ and $\Obs'(t)>0$ for
 every $t\in[t_{\mathrm c},1]$.
\end{thmlist}
\end{lemma}

The bound \eqref{eq:nc-third-positive} is pointwise: it applies to $I\cap(0,t_{\mathrm c})$ also
when the plateau $I$ extends beyond $t_{\mathrm c}$.

\begin{proof}
We compare with the critical pair, for which $\tilde a\equiv1$, $\tilde\alpha\equiv1$ and
$\tilde u=\tanh$. By \eqref{eq:nc-K}, $\|a-1\|_\infty\le12\,d(\xi)$ and
$\|a\alpha-1\|_{L^1}\le\|a-1\|_\infty+\|\alpha-1\|_{L^1}\le12\,d(\xi)+m_1(\zeta)$, so
\eqref{eq:nc-jet} gives
\begin{equation}\label{eq:nc-near-jets}
 \|\partial_x^ku-\partial_x^k\tanh\|_\infty\le C\bigl(d(\xi)+m_1(\zeta)\bigr),\qquad k=1,2,3,
\end{equation}
where here and below $C$ depends only on $\mathcal K$.

\ref{lem:nc-near-a} Let $t\in I\cap(0,t_{\mathrm c})$, where $I$ has height $\omega\in[0,1]$, and
put $R_\omega(y_1,y_2,y_3)=y_3^2-12\omega y_1y_2^2+6\omega^2y_1^4$, so that the last expectation in
\eqref{eq:nc-third} is $\E R_\omega(U_t,V_t,W_t)$. At the critical values,
$R_\omega(1,0,-2)=4+6\omega^2\ge4$. On the cube $[-C_\sharp,C_\sharp]^3$ with
$C_\sharp=\max\{C_1,C_2,C_3,2\}$, which contains $(U_t,V_t,W_t)$ by
Lemma~\ref{lem:nc-jets}\ref{lem:nc-jets-a}, $R_\omega$ is Lipschitz with a constant $L$ that does
not depend on $\omega\in[0,1]$. By \eqref{eq:nc-critical}, $(1,0,-2)$ is the value at $x=0$ of
$(\partial_x\tanh,\partial_x^2\tanh,\partial_x^3\tanh)$, and these functions are Lipschitz. So
\eqref{eq:nc-near-jets} and Lemma~\ref{lem:nc-jets}\ref{lem:nc-jets-c} give
\[
 \E\bigl|(U_t,V_t,W_t)-(1,0,-2)\bigr|\le C\bigl(d(\xi)+m_1(\zeta)+\E|X_t|\bigr)
 \le C_\flat\bigl(d(\xi)+m_1(\zeta)+\sqrt t\bigr)
\]
for a constant $C_\flat$. Hence $\E R_\omega(U_t,V_t,W_t)\ge4-LC_\flat(d(\xi)+m_1(\zeta)+\sqrt t)$.
The other two terms of \eqref{eq:nc-third} are at most $Cd(\xi)$ in absolute value, by
\eqref{eq:nc-K}. Choose $t_{\mathrm c}\in(0,1)$ with $LC_\flat\sqrt{t_{\mathrm c}}\le1$, and then
$\varepsilon_2>0$ with $2LC_\flat\varepsilon_2+C\varepsilon_2\le1$. Then
$\E R_\omega(U_t,V_t,W_t)\ge2$, so $\xi''\ge1$ gives
$\xi''(t)^3\,\E R_\omega(U_t,V_t,W_t)\ge4-1-LC_\flat(d(\xi)+m_1(\zeta))$, and
$\Gamma'''(t)\ge3-LC_\flat(d(\xi)+m_1(\zeta))-Cd(\xi)\ge2$.

\ref{lem:nc-near-b} By \eqref{eq:nc-critical-gamma}, $\kappa_F:=\min_{[t_{\mathrm c},1]}(t-\Gamma_{\mathrm c}(t))>0$.
By \eqref{eq:nc-continuity} with the critical pair,
$\sup_{[0,1]}|\Gamma-\Gamma_{\mathrm c}|\le C(d(\xi)+m_1(\zeta))$. Decreasing $\varepsilon_2$ so that
$2C\varepsilon_2<\kappa_F$, we get $\Gamma(t)<t$ on $[t_{\mathrm c},1]$, and then
$\Obs'=\frac12\xi''(t-\Gamma)>0$ there.
\end{proof}

At every local minimum of the obstacle, whether or not it carries mass, $\Gamma$ meets the diagonal
from the correct side.

\begin{lemma}[Contact conditions]\label{lem:nc-contact}
Let $\xi\in\mathcal K$ and $\zeta\in\Pspace([0,1])$. If $\Obs$ has a local minimum at
$t\in(0,1)$, then
\begin{equation}\label{eq:nc-contact}
 \Gamma(t)=t,\qquad\Gamma'(t)\le1 .
\end{equation}
If $\Obs$, as a function on $[0,1]$, has a local minimum at $0$, then \eqref{eq:nc-contact} holds at
$t=0$.
\end{lemma}

\begin{proof}
By Lemma~\ref{lem:pde-standing}\ref{lem:pde-standing-c}, $\Obs\in C^2([0,1])$ with
$\Obs'=-\frac12\xi''(\Gamma-t)$ and $\Obs''=-\frac12\xi'''(\Gamma-t)-\frac12\xi''(\Gamma'-1)$. At a
local minimum $t\in(0,1)$, $\Obs'(t)=0\le\Obs''(t)$, and $\xi''(t)>0$ gives \eqref{eq:nc-contact}.
At $0$, $\Gamma(0)=0$ by Lemma~\ref{lem:nc-plateau}\ref{lem:nc-plateau-a}, so
$\Gamma(s)-s=(\Gamma'(0)-1)s+o(s)$ and
\[
 \Obs(t)-\Obs(0)=-\frac12\int_0^t\xi''(s)\bigl(\Gamma(s)-s\bigr)\,ds
 =-\frac{\xi''(0)}4\bigl(\Gamma'(0)-1\bigr)t^2+o(t^2).
\]
A local minimum at $0$ makes the coefficient of $t^2$ nonnegative, and $\xi''(0)>0$ gives
$\Gamma'(0)\le1$.
\end{proof}

Zhou \cite{Zho26a} proves that, slightly below the critical temperature, the Parisi measure of the
Sherrington--Kirkpatrick model is supported on an interval starting at the origin. To exclude a gap
$(q,q')$ near the origin between two support points, he compares the mean of $\Gamma'$ over the gap
with $\frac12(\Gamma'(q)+\Gamma'(q'))$, which is at most one by the conditions of Auffinger and Chen
at support points \cite[Theorem~5]{AC15a}, and obtains the strict inequality, for short gaps, from a perturbation of
the critical pair involving the fourth spatial derivative of the Parisi solution, and for long gaps
from a global inequality at $\delta_0$ \cite[Section~1.5, Theorems~5 and~7]{Zho26a}. The next lemma runs
this comparison for an arbitrary measure, through
the pointwise bound \eqref{eq:nc-third-positive}; its endpoints need only satisfy
\eqref{eq:nc-contact}, so they may be contacts that carry no mass (Lemma~\ref{lem:nc-contact}), and
the measure may be a constrained minimizer. Zhou's theorem itself is not used; support points away
from the origin are handled by confinement (Appendix~\ref{app:nc-confine}).

\begin{lemma}[Gap exclusion]\label{lem:nc-gap}
Let $\xi\in\mathcal K$, $\zeta\in\Pspace([0,1])$ and $t_{\mathrm c}\in(0,1)$, and assume that
$\Gamma'''\ge2$ at every point of $(0,t_{\mathrm c})$ that lies in a plateau of $\zeta$. Then there
are no $0\le t_1<t_2<t_{\mathrm c}$ with $\zeta((t_1,t_2))=0$ such that \eqref{eq:nc-contact} holds
at $t_1$ and at $t_2$.
\end{lemma}

\begin{proof}
Suppose that such $t_1,t_2$ exist. Then $(t_1,t_2)$ is a plateau contained in $(0,t_{\mathrm c})$,
so $\Gamma'''\ge2$ on it and $\Gamma'$ is strictly convex on $(t_1,t_2)$. Since $\Gamma'$ is
continuous on $[t_1,t_2]$ (Lemma~\ref{lem:nc-plateau}\ref{lem:nc-plateau-a}), it is strictly convex
on $[t_1,t_2]$ and lies strictly below its chord on $(t_1,t_2)$. By \eqref{eq:nc-contact} at both
ends,
\[
 1=\frac{\Gamma(t_2)-\Gamma(t_1)}{t_2-t_1}
 =\frac1{t_2-t_1}\int_{t_1}^{t_2}\Gamma'(t)\,dt
 <\frac{\Gamma'(t_1)+\Gamma'(t_2)}2\le1,
\]
a contradiction.
\end{proof}

\subsection{Constrained minimizers stay near criticality}\label{app:nc-confine}

For $\xi\in\mathcal K$ and $\theta\in[0,1]$ let $\zeta^\xi_\theta$ be the constrained minimizer at
$\theta$ for the covariance $\xi$. Since $\delta_0$ is admissible for every constraint and is the
unique minimizer at $\xi_{\mathrm c}$, every constrained minimizer near criticality is close to
$\delta_0$, uniformly in the constraint.

\begin{lemma}[Confinement]\label{lem:nc-confine}
For every $\varepsilon>0$ there is $\varepsilon_3>0$ such that $m_1(\zeta^\xi_\theta)<\varepsilon$
for every $\xi\in\mathcal K$ with $d(\xi)<\varepsilon_3$ and every $\theta\in[0,1]$.
\end{lemma}

\begin{proof}
For $\xi\in\mathcal K$ and $\theta\in[0,1]$ the measure $\delta_0$ is feasible at $\theta$, so
\begin{equation}\label{eq:nc-trial}
 \inf_{\zeta}P^\xi_\zeta\le P^\xi_{\zeta^\xi_\theta}\le P^\xi_{\delta_0}.
\end{equation}
Suppose that the lemma fails for some $\varepsilon>0$. Then there are $\xi_n\in\mathcal K$ with
$d(\xi_n)\to0$ and $\theta_n\in[0,1]$ such that $\zeta_n:=\zeta^{\xi_n}_{\theta_n}$ has
$m_1(\zeta_n)\ge\varepsilon$. By compactness a subsequence converges weakly to some $\zeta_\infty$,
and $m_1(\zeta_\infty)\ge\varepsilon$, so $\zeta_\infty\ne\delta_0$. Weak convergence on $[0,1]$
gives $\|\cdf{\zeta_n}-\cdf{\zeta_\infty}\|_{L^1}\to0$, and $\|\xi_n''-1\|_\infty\to0$ by
\eqref{eq:nc-K}. By \eqref{eq:nc-continuity},
$P^{\xi_n}_{\zeta_n}\to P^{\xi_{\mathrm c}}_{\zeta_\infty}$ and
$P^{\xi_n}_{\delta_0}\to P^{\xi_{\mathrm c}}_{\delta_0}$ along the subsequence. The right half of
\eqref{eq:nc-trial} gives in the limit $P^{\xi_{\mathrm c}}_{\zeta_\infty}\le
P^{\xi_{\mathrm c}}_{\delta_0}=\min_\zeta P^{\xi_{\mathrm c}}_\zeta$, so $\zeta_\infty$ is a Parisi
measure of $\xi_{\mathrm c}$, and $\zeta_\infty=\delta_0$ by uniqueness. This is a contradiction.
\end{proof}

Combining Lemmas~\ref{lem:nc-near}--\ref{lem:nc-confine} gives the estimates on which the proof of
Proposition~\ref{prop:nc-example} runs.

\begin{corollary}[Standing estimates]\label{cor:nc-standing}
Let $t_{\mathrm c}$ and $\varepsilon_2$ be as in Lemma~\ref{lem:nc-near}, let $\varepsilon_3$ be
given by Lemma~\ref{lem:nc-confine} for $\varepsilon=\varepsilon_2$, and put
$\varepsilon_*=\min\{\frac12,\varepsilon_2,\varepsilon_3\}$. Let $\xi\in\mathcal K$ with
$d(\xi)<\varepsilon_*$, let $\theta\in[0,1]$, and let $\zeta=\zeta^\xi_\theta$, with $\Gamma$ and
$\Obs$ its second-moment function and obstacle. Then the following hold.
\begin{thmlist}
\item\label{cor:nc-standing-a} $\Gamma'''\ge2$ at every point of $(0,t_{\mathrm c})$ that lies in a
 plateau of $\zeta$; and $\Gamma<\id$ and $\Obs'>0$ on $[t_{\mathrm c},1]$.
\item\label{cor:nc-standing-b} Every minimum point of $\Obs$ over $(0,1]$ lies in
 $(0,t_{\mathrm c})$, and every minimum point of $\Obs$ over $[0,1]$ lies in $[0,t_{\mathrm c})$.
\item\label{cor:nc-standing-c} No two points $t_1<t_2$ of $\argmin_{(0,1]}\Obs$ satisfy
 $\zeta((t_1,t_2))=0$. If $\theta=0$, the same holds for any two points of $\argmin_{[0,1]}\Obs$.
\end{thmlist}
\end{corollary}

\begin{proof}
\ref{cor:nc-standing-a} By Lemma~\ref{lem:nc-confine}, $m_1(\zeta)<\varepsilon_2$, and
Lemma~\ref{lem:nc-near} applies.
\ref{cor:nc-standing-b} The obstacle is strictly increasing on $[t_{\mathrm c},1]$ by
\ref{cor:nc-standing-a}. So at every $t\in[t_{\mathrm c},1]$ it takes smaller values just below $t$,
at points of $(0,1]$; hence $t$ is not a minimum point over $(0,1]$ or over $[0,1]$.
\ref{cor:nc-standing-c} By \ref{cor:nc-standing-b}, the points of $\argmin_{(0,1]}\Obs$ are local
minima in $(0,t_{\mathrm c})$, and a point of $\argmin_{[0,1]}\Obs$ is either such a local minimum
or the point $0$, at which $\Obs$ restricted to $[0,1]$ is minimal. By Lemma~\ref{lem:nc-contact},
\eqref{eq:nc-contact} holds at all these points, and Lemma~\ref{lem:nc-gap}, whose hypothesis is
\ref{cor:nc-standing-a}, excludes the configuration.
\end{proof}

\subsection{Proof of Proposition~\ref*{prop:nc-example}}\label{app:nc-proof}

Let $\varepsilon_*$ be as in Corollary~\ref{cor:nc-standing}, and fix $\xi$ as in
\eqref{eq:nc-family}. Then $\xi\in\mathcal K$ and $d(\xi)<\varepsilon_*$, so
Corollary~\ref{cor:nc-standing} applies to every constrained minimizer $\zeta_\theta$,
$\theta\in[0,1]$. We write $\Gamma_\theta$ and $\Obs_\theta$ for the second-moment function and
the obstacle of $\zeta_\theta$, and $\Gamma_{\mathrm P}$, $\OP$ when $\theta=0$; thus
$\CP=\argmin_{[0,1]}\OP$. Besides Corollary~\ref{cor:nc-standing}, the proof uses three properties
of the covariance: $\xi''(0)=2b>1$, $\xi'''(t)=24et>0$ for $t>0$, and $\xi''''(0)=24e>0$.

\emph{Step 1: the equilibrium measure.} We prove part~\ref{prop:nc-example-a}.
\begin{steplist}
\item \emph{The origin is in the support.} Suppose that $r:=\min\supp\zP>0$. Then $r\in\CP$ by
 Lemma~\ref{lem:eq-obstacle}\ref{lem:eq-obstacle-a}, and $r<t_{\mathrm c}$ by
 Corollary~\ref{cor:nc-standing}\ref{cor:nc-standing-b}, so \eqref{eq:nc-contact} holds at $r$ by
 Lemma~\ref{lem:nc-contact}. The interval $(0,r)$ is a plateau of height $0$, and by
 \eqref{eq:nc-second} the function $h=\Gamma_{\mathrm P}-\id$ satisfies
 $h''=\xi'''\,\E U^2+\xi''^2\,\E V^2>0$ there, because $\xi'''>0$ on $(0,r)$ and $U>0$. So $h$ is
 strictly convex on $[0,r]$ with $h(0)=h(r)=0$; it is negative on $(0,r)$, and therefore
 $h'(r)\ge(h(r)-h(s))/(r-s)>0$ for $s\in(0,r)$, that is $\Gamma_{\mathrm P}'(r)>1$. This
 contradicts \eqref{eq:nc-contact}.
\item \emph{The measure is not $\delta_0$.} For $\delta_0$ the Parisi equation is solved by
 $\log\cosh x+\frac12(\xi'(1)-\xi'(t))$ \cite[(11)]{AC15a}, so by
 Lemma~\ref{lem:nc-plateau}\ref{lem:nc-plateau-a} its function $\Gamma$ has
 $\Gamma'(0)=\xi''(0)\sech^4(0)=2b>1$. If $\zP=\delta_0$, then $0\in\CP$ by
 Lemma~\ref{lem:eq-obstacle}\ref{lem:eq-obstacle-a}, and Lemma~\ref{lem:nc-contact} would give
 $\Gamma'(0)\le1$.
\item \emph{The support is an interval.} By Lemma~\ref{lem:eq-obstacle}\ref{lem:eq-obstacle-a} and
 Corollary~\ref{cor:nc-standing}\ref{cor:nc-standing-b}, $\supp\zP\subset\CP\subset[0,t_{\mathrm c})$.
 If $\supp\zP$ were not an interval, a connected component $(s_1,s_2)$ of
 $(\min\supp\zP,\max\supp\zP)\setminus\supp\zP$ would have $s_1,s_2\in\supp\zP\subset\CP$ and
 $\zP((s_1,s_2))=0$, which Corollary~\ref{cor:nc-standing}\ref{cor:nc-standing-c} excludes. With
 (1) and (2), $\supp\zP=[0,v_{\mathrm P}]$ for some $v_{\mathrm P}\in(0,t_{\mathrm c})$.
\item \emph{No atom at the origin.} Every point of $[0,v_{\mathrm P}]$ lies in $\CP$, so
 $\Gamma_{\mathrm P}(t)=t$ on $[0,v_{\mathrm P}]$ by Lemma~\ref{lem:nc-contact}. Hence
 $\Gamma_{\mathrm P}'\equiv1$ on $[0,v_{\mathrm P}]$, and $\Gamma_{\mathrm P}''(0+)=0$. By
 Lemma~\ref{lem:nc-plateau}\ref{lem:nc-plateau-c} at $t=0$, where $\xi'''(0)=0$ and $V_0=0$,
 \[
  0=\Gamma_{\mathrm P}''(0+)=-2a_0^2\,\zP(\{0\})\,U_{\mathrm P}^3 .
 \]
 Since $U_{\mathrm P}>0$, $\zP(\{0\})=0$; in the notation of
 Theorem~\ref{thm:exist}\ref{thm:exist-b}, $\thstar=0$. Finally
 $1=\Gamma_{\mathrm P}'(0)=a_0U_{\mathrm P}^2$ by Lemma~\ref{lem:nc-plateau}\ref{lem:nc-plateau-a}, so
 $U_{\mathrm P}=a_0^{-1/2}$.
\end{steplist}

\emph{Step 2: the constraint is active and the first contact is separated from the origin.} Let
$\theta\in(0,1)$. The facts proved here are those of
Lemma~\ref{lem:loc-obstacle}\ref{lem:loc-obstacle-b}--\ref{lem:loc-obstacle-c}, at a general
$\theta$ instead of $\thf$.
\begin{steplist}
\item \emph{Active constraint.} Since $\thstar=\zP(\{0\})=0$ by Step~1,
 Lemma~\ref{lem:var-basic}\ref{lem:var-basic-b} gives $\zeta_\theta(\{0\})=\theta$: a larger atom
 would leave room to move mass in every direction, and by convexity $\zeta_\theta$ would be the
 unconstrained minimizer $\zP$, which has no atom at $0$.
\item \emph{The contacts.} Write $\zeta_\theta=\theta\delta_0+(1-\theta)\rho_\theta$ with
 $\rho_\theta\in\Pspace((0,1])$, and put $c_\theta=\int\Obs_\theta\,d\rho_\theta$. As in the proof of
 Lemma~\ref{lem:eq-obstacle}\ref{lem:eq-obstacle-a}, the minimality of $\zeta_\theta$ on the convex
 set $\{\zeta(\{0\})\ge\theta\}$, \eqref{eq:loc-second-order} and \eqref{eq:loc-moving-mass} give
 $\int\Obs_\theta\,d(\mu-\zeta_\theta)\ge0$ for every feasible $\mu$. With
 $\mu=\theta\delta_0+(1-\theta)\delta_t$, $t\in(0,1]$, this gives $\Obs_\theta(t)\ge c_\theta$, with
 equality $\rho_\theta$-almost everywhere and hence, by continuity, on
 $\supp\zeta_\theta\cap(0,1]$. With $\mu=\delta_0$ it gives $\Obs_\theta(0)\ge c_\theta$. If
 $\Obs_\theta(0)=c_\theta$, then $\min_{[0,1]}\Obs_\theta=c_\theta=\int\Obs_\theta\,d\zeta_\theta$,
 so $\int\Obs_\theta\,d(\mu-\zeta_\theta)\ge0$ for every $\mu\in\Pspace([0,1])$, and
 \eqref{eq:loc-first-order} with \eqref{eq:loc-moving-mass} gives $P_\mu\ge P_{\zeta_\theta}$; then
 $\zeta_\theta=\zP$, which (1) excludes. Hence $\Obs_\theta(0)>c_\theta=\min_{(0,1]}\Obs_\theta$.
\item \emph{The first contact.} By (2), $\argmin_{(0,1]}\Obs_\theta=\{t\in(0,1]:\Obs_\theta(t)=c_\theta\}$
 contains $\supp\zeta_\theta\cap(0,1]$. It is nonempty and closed, it is bounded away from $0$ by
 the continuity of $\Obs_\theta$, and it lies in $(0,t_{\mathrm c})$ by
 Corollary~\ref{cor:nc-standing}\ref{cor:nc-standing-b}. Its least point $q_\theta$ is the first
 positive contact of $\zeta_\theta$, $\Obs_\theta>c_\theta$ on $[0,q_\theta)$, and
 $\zeta_\theta((0,q_\theta))=0$.
\end{steplist}

\emph{Step 3: the support of $\zeta_\theta$ and the marginal contact.} We prove
part~\ref{prop:nc-example-b}.
\begin{steplist}
\item \emph{Limits as $\theta\downarrow0$.} The function $p$ is nondecreasing, since the feasible sets
 decrease. The measures $\theta\delta_0+(1-\theta)\zP$ are feasible, so by convexity
 $\feq\le p(\theta)\le\theta P_{\delta_0}+(1-\theta)\feq\to\feq$. Every weak limit point $\zeta'$ of
 $\zeta_\theta$ as $\theta\downarrow0$ therefore has $P_{\zeta'}=\feq$, by the continuity of $P$
 (Lemma~\ref{lem:nc-jets}\ref{lem:nc-jets-d}), so $\zeta'=\zP$. By compactness,
 $\zeta_\theta\Rightarrow\zP$.
\item \emph{Choice of $\theta_0$.} There is $\theta_0\in(0,1)$ such that
 $\supp\zeta_\theta\cap(0,1]$ has at least two points for every $\theta\in(0,\theta_0)$. This set is
 nonempty because $\zeta_\theta((0,1])=1-\theta>0$. If it were a single point along a sequence
 $\theta_n\downarrow0$, then $\zeta_{\theta_n}=\theta_n\delta_0+(1-\theta_n)\delta_{r_n}$, whose weak
 limit points are Dirac masses, whereas $\zeta_{\theta_n}\Rightarrow\zP$ and $\supp\zP$ has positive
 length by Step~1.
\item \emph{The support.} Let $\theta\in(0,\theta_0)$ and $S_\theta=\supp\zeta_\theta\cap(0,1]$. By
 Step~2, $S_\theta$ is a closed subset of $[q_\theta,t_{\mathrm c})$ contained in
 $\argmin_{(0,1]}\Obs_\theta$. If $S_\theta$ were not an interval, two of its points would enclose
 an interval of zero $\zeta_\theta$-mass (take a component of the complement, as in Step~1), which
 Corollary~\ref{cor:nc-standing}\ref{cor:nc-standing-c} excludes. So $S_\theta=[r_\theta,v_\theta]$,
 with $r_\theta<v_\theta$ by (2). If $q_\theta<r_\theta$, the points $q_\theta$ and $r_\theta$ of
 $\argmin_{(0,1]}\Obs_\theta$ would enclose $(q_\theta,r_\theta)$, which has zero mass; so
 $r_\theta=q_\theta$. Together with Step~2,
 $\supp\zeta_\theta=\{0\}\cup[q_\theta,v_\theta]$ with $0<q_\theta<v_\theta<t_{\mathrm c}<1$.
\item \emph{Marginality.} Every point of $[q_\theta,v_\theta]$ is a positive contact and lies in
 $(0,1)$, so $\Gamma_\theta(t)=t$ on $[q_\theta,v_\theta]$ by Lemma~\ref{lem:nc-contact}. Since
 $v_\theta>q_\theta$ and $\Gamma_\theta'$ is continuous, $\Gamma_\theta'(q_\theta)=1$, which by
 Lemma~\ref{lem:nc-plateau}\ref{lem:nc-plateau-a} is the marginality identity of
 part~\ref{prop:nc-example-b}.
\end{steplist}

\emph{Step 4: the interval of levels.} We prove part~\ref{prop:nc-example-c}.
\begin{steplist}
\item \emph{Strict monotonicity of $p$.} If $p(\theta_1)=p(\theta_2)$ for some
 $0\le\theta_1<\theta_2<1$, then $\zeta_{\theta_2}$ is feasible at $\theta_1$ and attains
 $p(\theta_1)$, so $\zeta_{\theta_2}=\zeta_{\theta_1}$ by uniqueness; but their masses at $0$ are
 $\theta_2$ and $\theta_1$ (Steps~1 and~2). So $p$ is strictly increasing on $[0,1)$, and
 $p(\theta)>\feq$ for $\theta>0$.
\item \emph{Strict convexity of $\Lambda$.} Since $\thstar=0$,
 Lemma~\ref{lem:var-basic}\ref{lem:var-basic-d} shows that $\Lambar$ is strictly convex on $[0,1]$.
 (The mechanism: for $\theta=\lambda\theta_1+(1-\lambda)\theta_2$ with $\theta_1\ne\theta_2$, the
 mixture $\theta^{-1}(\lambda\theta_1\zeta_{\theta_1}+(1-\lambda)\theta_2\zeta_{\theta_2})$ has an
 atom $\theta+\lambda(1-\lambda)(\theta_1-\theta_2)^2/\theta>\theta$ at $0$, so it is feasible at
 $\theta$ but is not $\zeta_\theta$, and uniqueness makes the convexity inequality for $\Lambda$
 strict.)
\item \emph{The levels.} Let $\fzero=p(\theta_0)$. The function $\Lambar$ is convex by (2) and
 finite on $[0,1]$, since $\feq\le p\le P_{\delta_0}$, so it has a subgradient $f$ at the interior
 point $\theta_0$. Comparing with $\Lambar(0)=0$ gives
 $f\ge\Lambda(\theta_0)/\theta_0=\fzero$, and by (2) and the monotonicity of chord slopes of convex
 functions, $f<(\Lambar(s)-\Lambar(\theta_0))/(s-\theta_0)\le\Lambar'_-(1)=\fone$ for
 $s\in(\theta_0,1)$. Together with (1), $\feq<\fzero<\fone$. Now let $f\in(\feq,\fzero)$. For
 $s\in[\theta_0,1)$, $\Lambar(s)-sf=s(p(s)-f)\ge s(\fzero-f)>0$, and
 $\Lambar(1)-f=\liminf_{s\uparrow1}(\Lambda(s)-sf)\ge\fzero-f>0$. On the other hand $p(s)<f$ for
 small $s>0$ by Step~3(1), so $\Lambar(s)-sf<0=\Lambar(0)$. Hence the minimizer $\thf$ of
 $\theta\mapsto\Lambar(\theta)-\theta f$, which is unique by
 Lemma~\ref{lem:loc-obstacle}\ref{lem:loc-obstacle-a}, lies in $(0,\theta_0)$. Then
 $\zf=\zeta_{\thf}$, and by Step~3 its first positive contact $\qf=q_{\thf}$ satisfies
 $\Plef(\qf)=1-\xi''(\qf)\,\E U_{\qf}^2=0$.
\end{steplist}

\emph{Step 5: the onset.} We prove part~\ref{prop:nc-example-d}. Let $\theta\in(0,\theta_0)$, and
write $q=q_\theta$, $\Gamma=\Gamma_\theta$, and $U,V,W$ for the derivatives of
$\Phi_{\zeta_\theta}$ along the Auffinger--Chen process of $\zeta_\theta$. Put
$W_{\mathrm P}=\Phi_{\zP,xxxx}(0,0)$, so that $D_{\mathrm P}=\xi''''(0)U_{\mathrm P}^2+a_0^3W_{\mathrm P}^2$;
it is positive because $\xi''''(0)=24e>0$ and $U_{\mathrm P}>0$.

Since $\zeta_\theta\Rightarrow\zP$ and every interval $(0,\epsilon)$ with
$0<\epsilon<v_{\mathrm P}$ has positive $\zP$-mass, the portmanteau theorem gives
$\zeta_\theta((0,\epsilon))>0$ for small $\theta$; as $\zeta_\theta((0,q))=0$, we get $q\to0$ as
$\theta\downarrow0$. The interval $(0,q)$ is a plateau of $\zeta_\theta$ of height $\theta$. As
$\theta\downarrow0$,
\[
 \sup_{0<t<q}\E\bigl(|U_t-U_{\mathrm P}|+|V_t|+|W_t-W_{\mathrm P}|\bigr)\longrightarrow0 .
\]
Indeed, \eqref{eq:nc-jet} with $\tilde\xi=\xi$ and $\tilde\zeta=\zP$ makes the derivatives of
$\Phi_{\zeta_\theta}$ uniformly close to those of $\Phi_{\zP}$, because
$\|\xi''(\cdf{\zeta_\theta}-\cdf{\zP})\|_{L^1}\to0$; the derivatives of $\Phi_{\zP}$ are continuous
at $(0,0)$, where $\partial_x^3\Phi_{\zP}$ vanishes; and $\E X_t^2\le Cq$ for $t\le q$ by
Lemma~\ref{lem:nc-jets}\ref{lem:nc-jets-c}. In \eqref{eq:nc-third} with $\omega=\theta$, moreover,
$\xi''''\equiv24e=\xi''''(0)$, $\xi'''(t)\to0$ and $\xi''(t)\to a_0$ uniformly on $(0,q)$, and the
terms that carry the height $\theta$ tend to zero. Therefore
\begin{equation}\label{eq:nc-onset-third}
 \sup_{0<t<q}\bigl|\Gamma'''(t)-D_{\mathrm P}\bigr|\longrightarrow0\qquad(\theta\downarrow0).
\end{equation}
At the left end of the plateau, Lemma~\ref{lem:nc-plateau}\ref{lem:nc-plateau-b} with
$\xi'''(0)=0$ and $V_0=0$ gives the exact value
\begin{equation}\label{eq:nc-onset-zero}
 \Gamma''(0+)=-2a_0^2\,\theta\,\Phi_{\zeta_\theta,xx}(0,0)^3 .
\end{equation}
Since $\Gamma(0)=0$, $\Gamma(q)=q$ and $\Gamma'(q)=1$ (Step~3), two integrations by parts on the
plateau give the exact identity
\[
 0=q\Gamma'(q)-\Gamma(q)=\int_0^qt\,\Gamma''(t)\,dt
 =\frac{q^2}2\,\Gamma''(0+)+\frac12\int_0^q(q^2-s^2)\,\Gamma'''(s)\,ds .
\]
By \eqref{eq:nc-onset-third} the last integral is $\frac23q^3(D_{\mathrm P}+o(1))$, so
$\Gamma''(0+)=-\frac23D_{\mathrm P}\,q+o(q)$. Comparing with \eqref{eq:nc-onset-zero}, and using
$\Phi_{\zeta_\theta,xx}(0,0)\to U_{\mathrm P}=a_0^{-1/2}$, we obtain
\[
 \frac{q_\theta}\theta\longrightarrow\frac{3a_0^2U_{\mathrm P}^3}{D_{\mathrm P}}
 =\frac{3\sqrt{a_0}}{D_{\mathrm P}} .
\]
One more integration gives the left derivative of $\Gamma'$ at $q$,
\[
 \Gamma''(q-)=\Gamma''(0+)+\int_0^q\Gamma'''(s)\,ds=\frac{D_{\mathrm P}}3\,q+o(q),
\]
while $\Gamma'\equiv1$ on $[q,v_\theta]$ gives $\Gamma''(q+)=0$. By \eqref{eq:nc-jump},
\[
 2\,\xi''(q)^2\,\zeta_\theta(\{q\})\,\E U_q^3=\Gamma''(q-)-\Gamma''(q+)=\frac{D_{\mathrm P}}3\,q+o(q),
\]
and $\xi''(q)^2\,\E U_q^3\to a_0^2U_{\mathrm P}^3$ by the convergences above. With the limit of
$q_\theta/\theta$ this gives
\[
 \frac{\zeta_\theta(\{q_\theta\})}\theta\longrightarrow
 \frac{D_{\mathrm P}}3\cdot\frac{3a_0^2U_{\mathrm P}^3}{D_{\mathrm P}}\cdot\frac1{2a_0^2U_{\mathrm P}^3}
 =\frac12 .
\]
This proves part~\ref{prop:nc-example-d} and completes the proof of
Proposition~\ref{prop:nc-example}.\qed

\begin{remark}[The onset and its prediction]\label{rem:nc-onset}
As $\zeta_\theta((0,q_\theta))=0$, Step~5 gives $\zeta_\theta([0,q_\theta])=\frac32\theta+o(\theta)$.
Parisi and Rizzo \cite[Section~IV and Appendix~B, (B2)]{PR} found this picture for the
Sherrington--Kirkpatrick model in a perturbative replica computation, with the replica number $n$
for $\theta$: their order parameter $q(x)$ has a first plateau of height $\frac32n\dot q(0)+o(n)$
ending at $x_{\mathrm c}=\frac32n+O(n^2)$, where $\dot q(0)$ is the equilibrium slope at $x=0$. For
that model at low temperature, the companion paper \cite{PeiB}, not used here, shows by the same
calculus that the first positive contact is marginal for small $\theta$, with position and atom
given to leading order by \eqref{eq:nc-onset} at $\xi''''(0)=0$.
\end{remark}

\subsection{Every contact is marginal}\label{app:nc-marginal}

Lemma~\ref{lem:stab-family} (Section~\ref{sec:stab-family}) adds to
Proposition~\ref{prop:nc-example} that $[q_\theta,v_\theta]$ is the whole set of positive contacts
and that each of its points is marginal. The proof is a Rolle argument on the two plateaus
$(0,q_\theta)$ and $(v_\theta,1)$.

\begin{proof}[Proof of Lemma~\ref{lem:stab-family}]
Let $\xi$ be as in \eqref{eq:nc-family}, $\theta\in(0,\theta_0)$ and $\zeta=\zeta_\theta$, and write
$\Gamma=\Gz\zeta$, $\Obs=\Obs_\zeta$, $q=q_\theta$, $v=v_\theta$ and $h=\Gamma-\id$. We use the
following facts.
\begin{itemize}
\item By Proposition~\ref{prop:nc-example}\ref{prop:nc-example-b} and Step~2 of its proof,
 $\supp\zeta=\{0\}\cup[q,v]$, every point of $[q,v]$ minimizes $\Obs$ over $(0,1]$, $\Gamma'(q)=1$,
 and $\Obs(r)>\Obs(q)$ for $0\le r<q$.
\item By Corollary~\ref{cor:nc-standing}, $q,v\in(0,t_{\mathrm c})$, $h'''\ge2$ at every point of
 $(0,t_{\mathrm c})$ that lies in a plateau of $\zeta$, and $h<0$ on $[t_{\mathrm c},1]$.
\item By Lemma~\ref{lem:nc-plateau}, $h\in C^1([0,1])$, $h(0)=0$, and $h$ is $C^3$ on every
 plateau. By Lemma~\ref{lem:pde-standing}\ref{lem:pde-standing-c}, $\Obs'=-\frac12\xi''h$ and
 $\Obs''=-\frac12\xi'''h-\frac12\xi''(\Gamma'-1)$, with $\xi''>0$.
\end{itemize}

\emph{The middle.} Every point of $[q,v]$ is a minimum point of $\Obs$ over $(0,1]$ and lies in
$(0,1)$, so $\Obs'=0$ there, that is $h=0$ on $[q,v]$. Hence $h'=0$ on $(q,v)$ and, by the
continuity of $\Gamma'$, at $q$ and $v$.

\emph{The left plateau.} The interval $(0,q)$ is a plateau inside $(0,t_{\mathrm c})$, so
$h'''\ge2$ on it, and $h(0)=h(q)=h'(q)=0$. Suppose that $h(r)=0$ for some $r\in(0,q)$. Rolle's
theorem gives zeros of $h'$ at some $s_1\in(0,r)$ and $s_2\in(r,q)$, besides the zero at $q$; then
zeros of $h''$ in $(s_1,s_2)$ and in $(s_2,q)$; then a zero of $h'''$ in $(0,q)$, which is
impossible. ($h'$ is continuous on $[0,q]$ and $h''$ is $C^1$ on the open plateau, which is all
Rolle's theorem needs.) So $h$ has no zero in $(0,q)$ and has constant sign there. Since
$\Obs(r)-\Obs(q)=\frac12\int_r^q\xi''h>0$ for $r\in(0,q)$, the sign is positive: $h>0$ on $(0,q)$.

\emph{The right plateau.} On $(v,1]$ the distribution function of $\zeta$ equals one, so
$(v,t_{\mathrm c})$ lies in a plateau inside $(0,t_{\mathrm c})$, where $h'''\ge2$; also
$h(v)=h'(v)=0$, and $h(t_{\mathrm c})<0$. Two distinct zeros $r<r'$ of $h$ in $(v,t_{\mathrm c})$
would give zeros of $h'$ at $v$, in $(v,r)$ and in $(r,r')$, hence two zeros of $h''$ and one of
$h'''$, which is impossible. If $h$ had exactly one zero $r\in(v,t_{\mathrm c})$ with $h'(r)=0$, then
$h'$ would vanish at $v$, in $(v,r)$ and at $r$, with the same contradiction. If $h$ had exactly one
zero $r\in(v,t_{\mathrm c})$ with $h'(r)\ne0$, then $h$ would change sign at $r$, and $h(t_{\mathrm c})<0$
would force $h>0$ on $(v,r)$; then $\Obs(t)-\Obs(v)=-\frac12\int_v^t\xi''h<0$ for $t\in(v,r)$,
contradicting the minimality of $\Obs(v)$. So $h$ has no zero in $(v,t_{\mathrm c})$, and
$h(t_{\mathrm c})<0$ gives $h<0$ on $(v,t_{\mathrm c})$, hence on $(v,1]$.

\emph{Conclusion.} Since $\Obs'=-\frac12\xi''h$ with $\xi''>0$, the obstacle is strictly decreasing
on $[0,q]$, constant on $[q,v]$ and strictly increasing on $[v,1]$. Hence
$\argmin_{(0,1]}\Obs=[q,v]=\supp\zeta\cap(0,1]$, and also $\{t\in(0,1]:\Gamma(t)=t\}=[q,v]$. On
$[q,v]$ one has $h=0$ and $\Gamma'=1$, so $\Obs''=0$ there.
\end{proof}

%% file: sections/appG-moments.tex
\section{Fractional moments}\label{app:moments}

Theorem~\ref{thm:exist} reduces the large deviations of $\MN$, and the existence of generalized
TAP critical points, to \condFM. Here we prove \condFM\ for every convex $\xi\in\classC$
(Proposition~\ref{prop:max-fm}, hence Theorem~\ref{thm:exist}\ref{thm:exist-e}). Its upper half is
Guerra's bound, a finite-$N$ inequality for every convex $\xi$. Its lower half is elementary for
$\theta\le\thstar$; for $\theta\in(\thstar,1)$ it rests on Talagrand's computation
\cite[Section~7]{Tal07} if $\xi$ is even, and otherwise on the Kac--Rice construction of
Appendix~\ref{app:kr} through Theorem~\ref{thm:max}\ref{thm:max-b} (Table~\ref{tab:moments-inputs}).
Appendix~\ref{app:moments-pure} gives a second route for pure even $\xi$ that avoids
\cite[Section~7]{Tal07}.

\begin{table}[ht]
\centering
\small
\begin{tabular}{@{}p{0.25\textwidth}p{0.17\textwidth}p{0.15\textwidth}p{0.33\textwidth}@{}}
\toprule
Covariance & Upper half, & Lower half, & Lower half, \\
 & $0<\theta<1$ & $\theta\le\thstar$ & $\thstar<\theta<1$\\
\midrule
$\xi$ even
 & Lemma~\ref{lem:moments-guerra}
 & \raggedright elementary (also covered by \cite[\S7]{Tal07})
 & \raggedright \cite[\S7]{Tal07} with truncation
 (Proposition~\ref{prop:max-fm}\ref{prop:max-fm-c})\tabularnewline
\addlinespace
$\xi$ convex, not even
 & Lemma~\ref{lem:moments-guerra}
 & elementary
 & Proposition~\ref{prop:kr-main} and Theorem~\ref{thm:max}\ref{thm:max-b}\\
\addlinespace
$\xi=\beta^2t^p$, $p$ even (second route)
 & Lemma~\ref{lem:moments-guerra}
 & elementary
 & Proposition~\ref{prop:moments-pure}: Proposition~\ref{prop:kr-main} for $\xi+\gamma^2t^l$,
 transfer of $\MN$, $\gamma\downarrow0$\\
\bottomrule
\end{tabular}
\caption{Inputs for \condFM\ when $\xi$ satisfies \condH\ and $\xi\in\classC$.
Lemma~\ref{lem:moments-guerra} is Guerra's bound \cite[Theorem~6.1]{Tal07}; ``elementary'' is
Proposition~\ref{prop:max-fm}\ref{prop:max-fm-b}. Rows one and two prove
Proposition~\ref{prop:max-fm}; row three is the alternative of Appendix~\ref{app:moments-pure}.}
\label{tab:moments-inputs}
\end{table}

\subsection{Guerra's bound and the inputs for each class}\label{app:moments-guerra}

Guerra's interpolation \cite{Gue03} bounds the free energy above by Parisi functionals. The same
scheme, with the disorder averaged at exponent $\theta$, bounds the fractional moments when $\xi$
is convex \cite[Theorem~6.1]{Tal07}. In our normalization it reads as follows.

\begin{lemma}[Guerra's bound for fractional moments]\label{lem:moments-guerra}
Assume \condH. For every $N\ge1$, every $\theta\in(0,1)$ and every $\alpha\in\Pspace([0,1])$
with $\alpha(\{0\})\ge\theta$,
\begin{equation}\label{eq:moments-guerra}
 \frac1{N\theta}\log\E\ZN^\theta\le P_\alpha .
\end{equation}
Consequently $N^{-1}\log\E\ZN^\theta\le\Lambda(\theta)$ for every $N\ge1$ and every
$\theta\in(0,1)$.
\end{lemma}

\begin{proof}
Talagrand works with a Gaussian Hamiltonian on $\hcube$ whose covariance is $N\xi(\ov{\sigma}{\tau})$,
under the summability condition $\sum_pp^2\beta_p^2<\infty$, and allows an external field
\cite[(3.1)--(3.4) and Theorem~3.1]{Tal07}; we take the field to be zero. The summability
condition holds under \condR, and $\E\ZN^\theta$ depends only on the law of the centered
Gaussian vector $(\HN(\sigma))_{\sigma\in\hcube}$, whose covariance is $N\xi(\ov{\sigma}{\tau})$. So
\cite[Theorem~6.1]{Tal07}, which assumes $\xi$ convex, applies. Its parameters are numbers
$\theta=m_0<m_1<\dots<m_k<1$ and $0=q_0<q_1<\dots<q_k<q_{k+1}=1$; let $\alpha$ be the probability
measure on $[0,1]$ with $\alpha([0,t])=m_l$ for $t\in[q_l,q_{l+1})$, $0\le l\le k$. Talagrand's
partition function is $2^N\ZN$, and the terminal value of his Parisi recursion is $\log\cosh$. In
our normalization his bound therefore reads
\[
 \frac1{N\theta}\log\E\ZN^\theta\le X_0-\frac12\sum_{l=0}^km_l\bigl(\vth(q_{l+1})-\vth(q_l)\bigr),
\]
where $X_{k+1}=\log\cosh(z_0+\dots+z_k)$ for independent centered Gaussian variables $z_l$ of
variance $\xi'(q_{l+1})-\xi'(q_l)$, and $X_l=m_l^{-1}\log\E_l\exp(m_lX_{l+1})$ for
$l=k,\dots,0$, with $\E_l$ the expectation in $z_l$ \cite[(5.4)--(5.5)]{Tal07}. By the
Cole--Hopf formula \eqref{eq:pde-colehopf} and
Lemma~\ref{lem:pde-measurable}\ref{lem:pde-measurable-iv}, $X_l=\Phi_\alpha(q_l,z_0+\dots+z_{l-1})$
for every $l$, so $X_0=\Phi_\alpha(0,0)$. Since $\vth'(t)=t\xi''(t)$, the sum equals
$\int_0^1t\xi''(t)\alpha([0,t])\,dt$. The right-hand side is therefore $P_\alpha$, which proves
\eqref{eq:moments-guerra} for this $\alpha$.

Now let $\alpha\in\Pspace([0,1])$ with $\alpha(\{0\})\ge\theta$, and put $c(t)=\alpha([0,t])$, so
that $c\ge\theta$ on $[0,1)$. For $n>1/(1-\theta)$ let $c_n=\theta$ on $[0,1/n)$ and
\[
 c_n(t)=\min\Bigl\{1-\tfrac1n,\ \max\bigl\{\theta,\ n^{-1}\lfloor nc(t)\rfloor\bigr\}\Bigr\}
 \qquad\text{for } t\in[1/n,1).
\]
Then $c_n$ is a nondecreasing step function with finitely many values in $[\theta,1-1/n]$, equal to
$\theta$ on $[0,1/n)$. Its distinct values and its jump points are parameters as above, so
\eqref{eq:moments-guerra} holds for the measure $\alpha_n$ with distribution function $c_n$ on
$[0,1)$. Since $c\ge\theta$, we have $|c_n-c|\le1/n$ on $[1/n,1)$, hence
$\|c_n-c\|_{L^1}\le2/n$. By the Lipschitz bound
$|P(c_1)-P(c_2)|\le\frac12\xi''(1)\|c_1-c_2\|_{L^1}$ of
Corollary~\ref{cor:pde-step}\ref{cor:pde-step-c}, which holds for measurable coefficients by
Lemma~\ref{lem:pde-measurable}\ref{lem:pde-measurable-iii}, and since $P(\cdf\zeta)=P_\zeta$
(Lemma~\ref{lem:pde-measurable}\ref{lem:pde-measurable-iv}), $P_{\alpha_n}\to P_\alpha$. This
proves \eqref{eq:moments-guerra}. The last claim follows by taking the infimum over $\alpha$,
because $\Lambda(\theta)=\theta\inf\{P_\alpha:\alpha(\{0\})\ge\theta\}$.
\end{proof}

\begin{remark}[Convexity and $\alpha(\{0\})\ge\theta$]\label{rem:moments-interpolation}
For $\alpha$ with finitely many atoms, the time derivative of Guerra's interpolation equals
$-\frac12\int_0^1t\xi''(t)\alpha([0,t])\,dt$ minus averages of $\xi(R)-\xi(q)-\xi'(q)(R-q)$ over
overlaps $R$, with weight $\frac12(\alpha(\{0\})-\theta)\ge0$ at $q=0$ and $\frac12\alpha(\{q\})$ at
the atoms $q>0$; convexity on $[-1,1]$ makes these
remainders nonnegative. For non-convex $\xi$ a remainder at a negative overlap can be negative, and
the upper half of \condFM\ needs another proof (for $\xi$ that is not convex on $[-1,1]$ this step
is sketched in \cite[proof of Theorem~2.1]{CGKLM}; we do not rely on it).
\end{remark}

\paragraph{The lower half for even covariances.} For even $\xi$, Talagrand shows that Guerra's
bound, optimized over $\alpha$, is asymptotically sharp. When only even degrees are present, he
explains in \cite[Section~7]{Tal07} how to modify the proof of the Parisi
formula \cite{PF} to compute $\lim_N(N\theta)^{-1}\log\E\ZN^\theta$ for $\theta\in(0,1)$
\cite[Section~1]{Tal07}. The limit is the infimum of the right-hand side of Guerra's bound over
all step parameters with first exponent $\theta$ \cite[Sections~5 and~7]{Tal07}, argued there
by stated modifications of the proof of \cite{PF}. The admissible step parameters
are those of the proof of Lemma~\ref{lem:moments-guerra}, now with $q_1=0$ allowed; they describe
exactly the step measures $\alpha$ with $\alpha(\{0\})\ge\theta$ and $\alpha(\{1\})>0$. By the
approximation in that proof, their infimum is
$\inf\{P_\alpha:\alpha(\{0\})\ge\theta\}$. In our normalization the result is therefore
\begin{equation}\label{eq:moments-talagrand}
 \lim_{N\to\infty}\frac1N\log\E\ZN^\theta=\Lambda(\theta)\qquad\text{for every }\theta\in(0,1),
\end{equation}
for even $\xi$. We use it only for even polynomials and reach the other even $\xi\in\classC$ by
truncation.

\begin{proof}[Proof of Proposition~\ref{prop:max-fm}]
Recall the setting: $\xi$ satisfies \condH, $\xi\in\classC$, and $\theta\in(0,1)$.

\ref{prop:max-fm-a} This is Lemma~\ref{lem:moments-guerra}.

\ref{prop:max-fm-b} Let $\theta\le\thstar$. The Parisi measure satisfies
$\zP(\{0\})=\thstar\ge\theta$, so it is admissible in the constrained problem and
$\Lambda(\theta)=\theta P_{\zP}=\theta\feq$. By \eqref{eq:max-crude-M}, $\MN\ge N\feq$ for every
realization. On the event $B_t=\{\MN\le\log\ZN+tN\}$ we therefore have
$\ZN^\theta\ge e^{\theta(\MN-tN)}\ge e^{\theta N(\feq-t)}$, so that
$\E\ZN^\theta\ge e^{\theta N(\feq-t)}\Prob(B_t)$. By Proposition~\ref{prop:max-upper},
$\Prob(B_t)\to1$. Hence $\liminf_NN^{-1}\log\E\ZN^\theta\ge\theta(\feq-t)$ for every $t>0$. This
part uses only \condR.

\ref{prop:max-fm-c} Let $\xi$ be even. If $\xi$ is a polynomial, \eqref{eq:moments-talagrand} is
the claim. In general, truncate $\xi$ and $\HN$ at degree $D$ as in
Appendix~\ref{app:gauss-bounds}, with $D$ at least the smallest active degree. Then
$\xi_{\le D}$ is a nonzero even polynomial, so \eqref{eq:moments-talagrand} holds for the
truncated model:
$\lim_NN^{-1}\log\E(\ZN^{(D)})^\theta=\Lambda^{(D)}(\theta)$. By
Lemma~\ref{lem:gauss-bounds}\ref{lem:gauss-bounds-d3} and~\ref{lem:gauss-bounds-d1},
\[
 \liminf_{N\to\infty}\frac1N\log\E\ZN^\theta
 \ge\Lambda^{(D)}(\theta)-\theta\varepsilon_D-\frac{\theta^2}2\xi_{>D}(1)
 \ge\Lambda(\theta)-\frac32\xi'_{>D}(1)-\theta\varepsilon_D-\frac{\theta^2}2\xi_{>D}(1),
\]
and the right-hand side tends to $\Lambda(\theta)$ as $D\to\infty$. No Kac--Rice argument enters.

\ref{prop:max-fm-d} Let $\xi$ be convex on $[-1,1]$ but not even. Then $\xi$ has two active
degrees (see the beginning of Appendix~\ref{app:kr}), and Proposition~\ref{prop:kr-main} applies.
If $\thstar=1$, part~\ref{prop:max-fm-b} covers every
$\theta\in(0,1)$. Otherwise, by Theorem~\ref{thm:exist}\ref{thm:exist-b}, every
$\theta\in(\thstar,1)$ equals $\thf$ for the level $f=\Lambda'(\theta)\in(\feq,\fone)$, and
$\Sigma(f)=\Lambda(\theta)-\theta f$. Fix such $\theta$ and $f$, and let $\varepsilon>0$. Every
witness $(\zeta,x)$ of $\Eex{\varepsilon}{f}$ has $\FTAP(x)\ge N(f-\varepsilon)$, and
$\MN\ge\FTAP(x)$, so $\Eex{\varepsilon}{f}\subset\{\MN\ge N(f-\varepsilon)\}$ and
\[
 \E e^{\theta\MN}\ge e^{\theta N(f-\varepsilon)}\,\Prob\bigl(\Eex{\varepsilon}{f}\bigr).
\]
Proposition~\ref{prop:kr-main} gives, at this $\varepsilon$,
\[
 \liminf_{N\to\infty}\frac1N\log\Prob\bigl(\Eex{\varepsilon}{f}\bigr)\ge\Sigma(f).
\]
Hence
\[
 \liminf_{N\to\infty}\frac1N\log\E e^{\theta\MN}\ge\theta(f-\varepsilon)+\Sigma(f)
 =\Lambda(\theta)-\theta\varepsilon .
\]
Let $\varepsilon\downarrow0$; Theorem~\ref{thm:max}\ref{thm:max-b} turns the bound for
$\E e^{\theta\MN}$ into the same bound for $\E\ZN^\theta$. Part~\ref{prop:max-fm-b} covers
$\theta\le\thstar$. The lower bound of Proposition~\ref{prop:kr-main} is proved in
Appendix~\ref{app:kr} by a first-moment count of critical points and a deterministic
spherical-code bound; it uses neither \condFM\ nor
Theorem~\ref{thm:exist}\ref{thm:exist-c}--\ref{thm:exist-e}, so the argument is not circular.

Finally, part~\ref{prop:max-fm-a} together with~\ref{prop:max-fm-c} for even $\xi$, and
with~\ref{prop:max-fm-b} and~\ref{prop:max-fm-d} otherwise, gives \condFM.
\end{proof}

\subsection{Pure even covariances: a second route}\label{app:moments-pure}

Here $\xi(t)=\beta^2t^p$, $\beta>0$, $p\ge2$ even (SK and the pure even $p$-spin models); \condH\
holds and $\xi\in\classC$. Proposition~\ref{prop:max-fm}\ref{prop:max-fm-a}, \ref{prop:max-fm-b}
leave the lower half of \condFM\ for $\theta\in(\thstar,1)$, which~\ref{prop:max-fm-c} takes from
\cite[Section~7]{Tal07}. The Kac--Rice route of~\ref{prop:max-fm-d} needs $\HN(x)$ to keep a
positive conditional variance given $\nabla\HN(x)$, which Euler's identity
$x\cdot\nabla\HN(x)=p\HN(x)$ removes. We perturb $\xi$ by an independent second degree,
as Boursier does \cite[Section~3.2]{Bou26b}, to which
Proposition~\ref{prop:kr-main} applies; the TAP free energies differ by the perturbing field plus
$O(\gamma^2N)$, so we transfer the maximum, not the critical points
(Remark~\ref{rem:moments-pure-scope}).

\begin{proposition}[Pure even covariances]\label{prop:moments-pure}
Let $\xi(t)=\beta^2t^p$ with $\beta>0$ and $p\ge2$ even, and let $f\in(\feq,\fone)$. For every
$\varepsilon>0$,
\[
 \liminf_{N\to\infty}\frac1N\log\Prob\bigl(|N^{-1}\MN-f|\le\varepsilon\bigr)\ge\Sigma(f).
\]
\end{proposition}

The proof uses Proposition~\ref{prop:kr-main} and Lemma~\ref{lem:moments-guerra} for the
covariances $\beta^2t^p+\gamma^2t^l$, and Theorems~\ref{thm:max} and~\ref{thm:exist}\ref{thm:exist-b},
\ref{thm:exist-c}; it does not use \cite[Section~7]{Tal07}.

\begin{corollary}[A second route to the fractional-moment formula]\label{cor:moments-pure-all}
Let $\xi$ be as in Proposition~\ref{prop:moments-pure}. Then \condFM\ holds. Consequently
$\lim_NN^{-1}\log\E e^{\theta\MN}=\Lambda(\theta)$ for every $\theta\in(0,1)$, and the conclusions
of Theorem~\ref{thm:exist}\ref{thm:exist-d} hold at every $f\in(\feq,\fone)$. These statements
are thereby proved without \cite[Section~7]{Tal07}.
\end{corollary}

\begin{proof}
The upper half of \condFM\ is Proposition~\ref{prop:max-fm}\ref{prop:max-fm-a}, and the lower half
for $\theta\le\thstar$ is Proposition~\ref{prop:max-fm}\ref{prop:max-fm-b}. Let
$\theta\in(\thstar,1)$. By Theorem~\ref{thm:exist}\ref{thm:exist-b}, $\theta=\thf$ for
$f=\Lambda'(\theta)\in(\feq,\fone)$, and $\Sigma(f)=\Lambda(\theta)-\theta f$. For $\varepsilon>0$,
$\E e^{\theta\MN}\ge e^{\theta N(f-\varepsilon)}\Prob(|N^{-1}\MN-f|\le\varepsilon)$, so
Proposition~\ref{prop:moments-pure} gives
\[
 \liminf_{N\to\infty}\frac1N\log\E e^{\theta\MN}\ge\theta(f-\varepsilon)+\Sigma(f)
 =\Lambda(\theta)-\theta\varepsilon .
\]
Let $\varepsilon\downarrow0$ and apply Theorem~\ref{thm:max}\ref{thm:max-b}. This proves \condFM;
the other statements are Theorem~\ref{thm:exist}\ref{thm:exist-c}, \ref{thm:exist-d}.
\end{proof}

\paragraph{The perturbation.} Fix an even $l\ge2$ with $l\ne p$. Let $H^{(l)}_N$ be a copy of the
degree-$l$ field $H_{N,l}$ of Appendix~\ref{app:gauss-bounds}, independent of $\HN$; it is a
centered Gaussian field on $\cube$ with covariance $N\ov xy^l$. For $\gamma>0$ put
\begin{equation}\label{eq:moments-perturbation}
 \xi_\gamma(t)=\xi(t)+\gamma^2t^l,\qquad
 H^\gamma_N=\HN+\gamma H^{(l)}_N,\qquad
 d_\gamma=l\gamma^2,\qquad
 S_N=\max_{\cube}|H^{(l)}_N| .
\end{equation}
The field $H^\gamma_N$ is centered Gaussian with covariance $N\xi_\gamma(\ov xy)$. The polynomial
$\xi_\gamma$ is even and has two active degrees, so it satisfies \condH, lies in $\classC$, and
Proposition~\ref{prop:kr-main} applies to it; the degree $l$ is even so that $\xi_\gamma$ stays
convex on $[-1,1]$. Objects built from $(\xi_\gamma,H^\gamma_N)$ carry a superscript or subscript
$\gamma$: $F^\gamma_\zeta$, $F^\gamma_{\mathrm{TAP}}$, $M^\gamma_N$, $Z^\gamma_N$,
$\mathcal E^\gamma_{N,\varepsilon}(f)$, $P^\gamma_\zeta$, $\feq^\gamma$, $\Lambda_\gamma$,
$\Lambar_\gamma$, $\Sigma_\gamma$, $\fone^\gamma$. Every event below lives on the product space of
$(\HN,H^{(l)}_N)$; an event determined by $\HN$ alone has the same probability as in the original
model.

The comparison of the two models rests on writing the TAP functions through a correction term. For
a covariance $\xi$ satisfying \condR, $\alpha\in\Pspace([0,1])$, $r\in[0,1]$ and
$\mu\in\Pspace([-1,1])$, let $\Phi^\xi_\alpha$ be the solution of \eqref{eq:intro-parisi-pde}
for $\xi$, let $\Phi^{\xi,*}_\alpha$ be its concave conjugate, and put
\[
 V^\xi_\alpha(r,\mu)=\int\Phi^{\xi,*}_\alpha(r,m)\,\mu(dm)
 -\frac12\int_r^1t\xi''(t)\alpha([0,t])\,dt,\qquad
 V^\xi_{\mathrm{TAP}}(r,\mu)=\inf_{\alpha\in\Pspace([0,1])}V^\xi_\alpha(r,\mu) .
\]
These are $V_\alpha$ of \eqref{eq:hess-V} and $V_{\mathrm{TAP}}$ of
Appendix~\ref{app:matrix-hessian}, with the covariance shown. By \eqref{eq:intro-tap-function}
and \eqref{eq:intro-ftap}, for every $x\in\cube$,
\begin{equation}\label{eq:moments-identities}
 F_\zeta(x)=\HN(x)+N\,V^\xi_\zeta(q_x,\mu_x),\qquad
 \FTAP(x)=\HN(x)+N\,V^\xi_{\mathrm{TAP}}(q_x,\mu_x).
\end{equation}

\begin{lemma}[Comparison of covariances]\label{lem:moments-cov}
Let $\xi_a$ and $\xi_b$ satisfy \condR, and assume that $\kappa=\xi_b-\xi_a$ has nonnegative
coefficients. For all $\alpha\in\Pspace([0,1])$, $r\in[0,1]$, $y\in\R$ and
$\mu\in\Pspace([-1,1])$,
\begin{equation}\label{eq:moments-pde-compare}
 0\le\Phi^{\xi_b}_\alpha(r,y)-\Phi^{\xi_a}_\alpha(r,y)\le\kappa'(1)-\kappa'(r),
\end{equation}
and each of $|V^{\xi_b}_\alpha(r,\mu)-V^{\xi_a}_\alpha(r,\mu)|$,
$|V^{\xi_b}_{\mathrm{TAP}}(r,\mu)-V^{\xi_a}_{\mathrm{TAP}}(r,\mu)|$ and
$|P^{\xi_b}_\alpha-P^{\xi_a}_\alpha|$ is
at most $\kappa'(1)$.
\end{lemma}

\begin{proof}
The proof of Lemma~\ref{lem:gauss-bounds}\ref{lem:gauss-bounds-d1} uses only that $\xi$ and
$\xi_{\le D}$ satisfy \condR\ and that $\xi_{>D}=\xi-\xi_{\le D}$ has nonnegative coefficients.
With $(\xi,\xi_{\le D},\xi_{>D})$ replaced by $(\xi_b,\xi_a,\kappa)$, its representation
\eqref{eq:gauss-delta} reads, for $\alpha$ with finitely many atoms and $u=\Phi^{\xi_b}_\alpha$,
\[
 \Phi^{\xi_b}_\alpha(r,y)-\Phi^{\xi_a}_\alpha(r,y)
 =\frac12\,\E\int_r^1\kappa''(s)\bigl(u_{xx}+\cdf\alpha\,u_x^2\bigr)(s,X_s)\,ds ,
\]
where $X$ is a diffusion started at $(r,y)$. By Lemma~\ref{lem:pde-class} with $G=1$,
$0\le u_{xx}+\cdf\alpha\,u_x^2\le2$, which gives \eqref{eq:moments-pde-compare} for such $\alpha$.
For general $\alpha$, approximate $\cdf\alpha$ in $L^1$ by distribution functions of measures with
finitely many atoms; by Lemma~\ref{lem:pde-measurable}\ref{lem:pde-measurable-i},
\ref{lem:pde-measurable-iv}, applied to each covariance, the solutions converge uniformly, and
\eqref{eq:moments-pde-compare} passes to the limit.

If two functions satisfy $\psi_a\le\psi_b\le\psi_a+K$, their concave conjugates
$m\mapsto\inf_y\{\psi(y)-my\}$ satisfy the same inequalities. So
$0\le\Phi^{\xi_b,*}_\alpha-\Phi^{\xi_a,*}_\alpha\le\kappa'(1)$ on $[0,1]\times(-1,1)$, and at
$m=\pm1$ by continuity. The integral terms of $V$ differ by
$\frac12\int_r^1t\kappa''(t)\alpha([0,t])\,dt\in[0,\frac12(\kappa'(1)-\kappa(1))]$. Hence
$V^{\xi_b}_\alpha-V^{\xi_a}_\alpha\in[-\frac12\kappa'(1),\kappa'(1)]$, which gives the bound for $V$;
the bound for $V_{\mathrm{TAP}}$ follows by taking infima over the common set $\Pspace([0,1])$.
Finally,
$\Phi^\xi_\alpha(0,\cdot)$ is even and convex, so $\Phi^{\xi,*}_\alpha(0,0)=\Phi^\xi_\alpha(0,0)$ and
$P^\xi_\alpha=V^\xi_\alpha(0,\delta_0)$.
\end{proof}

\begin{lemma}[The perturbed models]\label{lem:moments-perturb}
In the setting \eqref{eq:moments-perturbation}:
\begin{romanlist}
\item\label{lem:moments-perturb-i} For every realization and every $x\in\cube$,
 $|F^\gamma_{\mathrm{TAP}}(x)-\FTAP(x)-\gamma H^{(l)}_N(x)|\le d_\gamma N$. Hence
 $|M^\gamma_N-\MN|\le\gamma S_N+d_\gamma N$.
\item\label{lem:moments-perturb-ii} $|P^\gamma_\zeta-P_\zeta|\le d_\gamma$ for every
 $\zeta\in\Pspace([0,1])$. Moreover $|\feq^\gamma-\feq|\le d_\gamma$,
 $|\Lambar_\gamma-\Lambar|\le d_\gamma$ on $[0,1]$, and $|\Sigma_\gamma-\Sigma|\le d_\gamma$
 on $\R$.
\item\label{lem:moments-perturb-iii} $\liminf_{\gamma\downarrow0}\fone^\gamma\ge\fone$.
 Consequently, for every $f\in(\feq,\fone)$ there is $\gamma_0(f)>0$ such that
 $f\in(\feq^\gamma,\fone^\gamma)$ for all $\gamma\in(0,\gamma_0(f))$.
\end{romanlist}
\end{lemma}

\begin{proof}
\ref{lem:moments-perturb-i} By \eqref{eq:moments-identities} for both models,
$F^\gamma_{\mathrm{TAP}}(x)-\FTAP(x)-\gamma H^{(l)}_N(x)
=N[V^{\xi_\gamma}_{\mathrm{TAP}}-V^\xi_{\mathrm{TAP}}](q_x,\mu_x)$. Apply
Lemma~\ref{lem:moments-cov} with $(\xi_a,\xi_b)=(\xi,\xi_\gamma)$, for which $\kappa(t)=\gamma^2t^l$
and $\kappa'(1)=d_\gamma$.
Taking maxima over $\cube$ gives the bound on $|M^\gamma_N-\MN|$.

\ref{lem:moments-perturb-ii} The first bound is Lemma~\ref{lem:moments-cov}. The quantities $\feq$
and $\Lambda(\theta)/\theta$, $0<\theta<1$, are infima of $P_\zeta$ over sets of order
parameters that do not depend on the covariance, so $|\feq^\gamma-\feq|\le d_\gamma$ and
$|\Lambda_\gamma-\Lambda|\le\theta d_\gamma\le d_\gamma$ on $[0,1)$. By
Theorem~\ref{thm:exist}\ref{thm:exist-c}, $\Lambar(1)=\xi(1)/2$ and
$\Lambar_\gamma(1)=\xi_\gamma(1)/2$, which differ by $\gamma^2/2\le d_\gamma$. Taking infima of
$\Lambar_\gamma(\theta)-\theta f$ and $\Lambar(\theta)-\theta f$ over $\theta\in[0,1]$ gives the
bound for $\Sigma$.

\ref{lem:moments-perturb-iii} By Theorem~\ref{thm:exist}\ref{thm:exist-b}, \ref{thm:exist-c},
applied to $\xi$ and to $\xi_\gamma$, the functions $\Lambar$ and $\Lambar_\gamma$ are convex and
finite on $[0,1]$: $\Lambda$ is continuously differentiable on $(0,1)$ with nondecreasing
derivative, $\Lambda(\theta)\to0$ as $\theta\downarrow0$ because
$\feq\le\Lambda(\theta)/\theta\le P_{\delta_0}$, and $\Lambar(1)$ is finite. For a convex function
$g$ that is finite on $[0,1]$, the chord slope $(g(1)-g(1-h))/h$ is nonincreasing in
$h\in(0,1]$, and $g'_-(1)=\sup_{0<h\le1}(g(1)-g(1-h))/h$. By~\ref{lem:moments-perturb-ii}, for
every $h\in(0,1]$,
\[
 \fone^\gamma\ge\frac{\Lambar_\gamma(1)-\Lambar_\gamma(1-h)}h
 \ge\frac{\Lambar(1)-\Lambar(1-h)}h-\frac{2d_\gamma}h .
\]
Let $\gamma\downarrow0$, then take the supremum over $h$. For the last claim,
$\feq^\gamma\le\feq+d_\gamma<f$ once $d_\gamma<f-\feq$, and $\fone^\gamma>f$ for all
small $\gamma$ because $\liminf_{\gamma\downarrow0}\fone^\gamma\ge\fone>f$.
\end{proof}

\begin{lemma}[The perturbing field]\label{lem:moments-field}
Let $c_l=2l\sqrt{2/\pi}$. For every $N$ and every $u\ge0$,
\[
 \Prob\bigl(S_N\ge(c_l+u)N\bigr)\le e^{-u^2N/2}.
\]
\end{lemma}

\begin{proof}
By \eqref{eq:gauss-sup-p} with $p=l$ and $\rho=1$, $\E S_N\le c_lN$. The variable $S_N$ is the
supremum of the centered Gaussian process $(x,s)\mapsto sH^{(l)}_N(x)$ on $\cube\times\{-1,1\}$,
whose variance $Nq_x^l$ is at most $N$, so \eqref{eq:gauss-conc} gives
$\Prob(S_N-\E S_N\ge uN)\le e^{-u^2N/2}$.
\end{proof}

\begin{proof}[Proof of Proposition~\ref{prop:moments-pure}]
Fix $f\in(\feq,\fone)$ and $\varepsilon>0$, and put $\theta=\thf$. By
Theorem~\ref{thm:exist}\ref{thm:exist-b}, $\theta\in(0,1)$ and $\Sigma(f)=\Lambda(\theta)-\theta f$;
also $-\infty<\Sigma(f)\le\Lambar(0)=0$. Choose $\gamma\in(0,\gamma_0(f))$, with $\gamma_0(f)$ from
Lemma~\ref{lem:moments-perturb}\ref{lem:moments-perturb-iii}, so small that
\begin{equation}\label{eq:moments-gamma}
 d_\gamma\le\min\Bigl\{1,\frac{\theta\varepsilon}8\Bigr\},\qquad
 \frac\varepsilon{4\gamma}>c_l,\qquad
 \frac12\Bigl(\frac\varepsilon{4\gamma}-c_l\Bigr)^2>1-\Sigma(f).
\end{equation}
Put
\[
 a_N=\Prob\bigl(M^\gamma_N\ge N(f-\tfrac\varepsilon2)\bigr),\qquad
 b_N=\Prob\bigl(M^\gamma_N>N(f+\tfrac\varepsilon2)\bigr),\qquad
 p_N=\Prob\bigl(\gamma S_N>\tfrac{\varepsilon N}4\bigr).
\]
Inside the first event and outside the other two, $|M^\gamma_N-Nf|\le\varepsilon N/2$ and
$\gamma S_N\le\varepsilon N/4$, and Lemma~\ref{lem:moments-perturb}\ref{lem:moments-perturb-i}
with $d_\gamma<\varepsilon/4$ gives $|\MN-Nf|\le\varepsilon N$. Hence
\[
 \Prob\bigl(|N^{-1}\MN-f|\le\varepsilon\bigr)\ge a_N-b_N-p_N .
\]
We bound the three terms.
\begin{itemize}
\item \emph{The existence term.} Every witness $(\zeta,x)$ of
 $\mathcal E^\gamma_{N,\varepsilon/2}(f)$ has $F^\gamma_{\mathrm{TAP}}(x)\ge N(f-\varepsilon/2)$, so
 $a_N\ge\Prob(\mathcal E^\gamma_{N,\varepsilon/2}(f))$. Since $f\in(\feq^\gamma,\fone^\gamma)$,
 Proposition~\ref{prop:kr-main} applies to $\xi_\gamma$ at the level $f$, and it gives
 $\liminf_NN^{-1}\log\Prob(\mathcal E^\gamma_{N,\varepsilon/2}(f))\ge\Sigma_\gamma(f)$. With
 Lemma~\ref{lem:moments-perturb}\ref{lem:moments-perturb-ii},
 $\liminf_NN^{-1}\log a_N\ge\Sigma(f)-d_\gamma$.
\item \emph{The overshoot term.} Markov's inequality gives
 $b_N\le e^{-\theta N(f+\varepsilon/2)}\E e^{\theta M^\gamma_N}$. By
 Theorem~\ref{thm:max}\ref{thm:max-b} and Lemma~\ref{lem:moments-guerra} for $\xi_\gamma$, and by
 Lemma~\ref{lem:moments-perturb}\ref{lem:moments-perturb-ii},
 \[
  \limsup_{N\to\infty}\frac1N\log\E e^{\theta M^\gamma_N}
  =\limsup_{N\to\infty}\frac1N\log\E(Z^\gamma_N)^\theta
  \le\Lambda_\gamma(\theta)\le\Lambda(\theta)+d_\gamma .
 \]
 Hence, by~\eqref{eq:moments-gamma},
 \[
  \limsup_{N\to\infty}\frac1N\log b_N\le\Lambda(\theta)+d_\gamma-\theta f-\frac{\theta\varepsilon}2
  =\Sigma(f)+d_\gamma-\frac{\theta\varepsilon}2\le\Sigma(f)-d_\gamma-\frac{\theta\varepsilon}4 .
 \]
\item \emph{The field term.} By Lemma~\ref{lem:moments-field} with
 $u=\varepsilon/(4\gamma)-c_l>0$ and by~\eqref{eq:moments-gamma},
 $\limsup_NN^{-1}\log p_N\le-\frac12(\varepsilon/(4\gamma)-c_l)^2<\Sigma(f)-1\le\Sigma(f)-d_\gamma$.
\end{itemize}
Put $s_\gamma=\Sigma(f)-d_\gamma>-\infty$. Then $\liminf_NN^{-1}\log a_N\ge s_\gamma$ and
$\limsup_NN^{-1}\log(b_N+p_N)<s_\gamma$. So for every sufficiently small $\eta>0$ and all large
$N$, $a_N\ge e^{N(s_\gamma-\eta)}$ and $b_N+p_N\le e^{N(s_\gamma-2\eta)}$, whence
$a_N-b_N-p_N\ge e^{N(s_\gamma-\eta)}(1-e^{-\eta N})$ and
$\liminf_NN^{-1}\log(a_N-b_N-p_N)\ge s_\gamma-\eta$.
Therefore
\[
 \liminf_{N\to\infty}\frac1N\log\Prob\bigl(|N^{-1}\MN-f|\le\varepsilon\bigr)\ge\Sigma(f)-d_\gamma .
\]
The left-hand side involves $\HN$ alone and does not depend on $\gamma$. Letting
$\gamma\downarrow0$ proves the proposition.
\end{proof}

The order of the limits is: fix $f$ and $\varepsilon$, choose $\gamma$, let $N\to\infty$, and only
then let $\gamma\downarrow0$. No estimate needs to be uniform in $\gamma$: the perturbation enters
only through the deterministic constant $d_\gamma$ and the explicit Gaussian tail of
Lemma~\ref{lem:moments-field}. The inclusions used above hold for every realization, so if the
existence event is not known to be measurable, the bounds hold for its inner and outer
probabilities alike.

\begin{remark}[What the transfer gives]\label{rem:moments-pure-scope}
Corollary~\ref{cor:moments-pure-all} and Theorem~\ref{thm:exist}\ref{thm:exist-a} yield
generalized TAP critical points at rate $\Sigma(f)$ at every level, of the kind in
Remark~\ref{rem:exist-points}. Boursier uses this perturbation \cite[Section~3.2]{Bou26b} but
transfers retained critical points, by \cite[Proposition~B.1]{Bou26b} under the strict Plefka
condition, false at marginal first contacts; our transfer of $\MN$ needs no
Hessian estimate, only the deterministic comparison of
Lemma~\ref{lem:moments-perturb}\ref{lem:moments-perturb-i} and the Gaussian tail of
Lemma~\ref{lem:moments-field}.
\end{remark}